\documentclass[11pt]{article}
\usepackage{graphicx} % Required for inserting images
\usepackage{format}

\title{From Good Starts to Optimal Inference: Generalized Latent Factor Models with Missingness and Implicit Regularization}
\begin{document}
\author{Chengzhu Huang \and Yuqi Gu}
\date{Department of Statistics, Columbia University}
\maketitle
\begin{abstract}
Generalized latent factor models provide a flexible framework for analyzing high-dimensional non-Gaussian data, but principled estimation and uncertainty quantification under missingness remain substantially less developed. We develop a theory that connects a computationally tractable nonconvex procedure directly to statistical inference for nonlinear latent factor models with exponential-family links and partially observed entries. Our procedure combines a link-aware double-SVD initialization, a unilateral refinement that achieves rowwise consistency, and vanilla gradient descent. We show that the refined initializer enters a region of incoherence and contraction and that gradient descent remains in this region through implicit regularization, contracting rapidly down to the statistical estimation error without explicit incoherence or balancing regularization. Our central result is a uniform rowwise linear approximation for the actual output of gradient descent that isolates the leading score fluctuations from higher-order estimation and optimization errors. These expansions yield asymptotically valid individual and Gaussian multiplier-bootstrap simultaneous inference for latent factors, together with simultaneous confidence bands for missing-entry means, without requiring an additional debiasing step. The resulting estimation rate matches a restricted-class minimax lower bound up to logarithmic factors, while the theory accommodates severe missingness, weak low-rank signals, and diminishing local curvature. Simulations support the theoretical findings, and an application to large language model evaluation illustrates uncertainty-aware estimation and ranking of latent model capabilities.
\end{abstract}     


\section{Introduction}
Latent factor models \citep{bartholomew2008analysis} have found widespread use in the social sciences, psychology, recommendation systems, and many other fields. When linear models are inadequate for heterogeneous discrete observations, generalized latent factor models \citep{bartholomew2008analysis,skrondal2004generalized,huber2004estimation,kidzinski2022generalized,chen2024note} provide a flexible framework for binary, ordinal, and count data. More importantly, amid the rapid development of large language models, there is a pressing need for trustworthy evaluation procedures that scale to many models and many questions. Given the binary correct/incorrect responses of a collection of language models to a common set of questions, latent factor models offer a powerful way to uncover their underlying capabilities and distinctive strengths. 

In many applications, however, only a subset of the data is observed. Survey nonresponse is often unavoidable, while evaluating every language model on every available question can be prohibitively expensive. Statistical procedures must therefore accommodate missingness in both latent factor estimation and the subsequent uncertainty quantification.

For linear low-rank models with missing data, a substantial body of work has developed estimators together with rigorous theoretical guarantees \citep{candes2012exact,chen2019noisy,xia2021normal}. Their extension to generalized latent factor models remains considerably less developed, even though latent factor models for non-Gaussian data have wide-ranging applications in practice.
More broadly, estimation and inference for nonlinear latent factor models under missingness are not yet fully understood, nor is the extent to which missingness and sparsity compromise the reliability of the resulting estimates. In this work, we develop statistical foundations for estimation and uncertainty quantification in generalized latent factor models with partially observed data, with an emphasis on statistical inference, optimality, and computational tractability.

To formulate the problem, let $\mb R \in \bb R^{n \times p}$ denote the underlying observation matrix, let $\Omega\subseteq[n]\times[p]$ denote the set of observed locations, and let $\bo\zeta^* \in \bb R^p$, $\mb X^* \in \bb R^{n \times r}$, and $\mb Y^* \in \bb R^{p \times r}$ denote the intercept vector and the left and right factor matrices, respectively. We consider an exponential-family model under which, conditional on the latent parameters and the observation set $\Omega$, the likelihood is given by
\begin{equation}
     l(\mb R \mid [\bo\zeta^*, \mb X^*, \mb Y^*], \Omega) = \prod_{(i,j)\in \Omega}\exp\big[\gamma(R_{i,j}) + R_{i,j}M^*_{i,j} - \Psi(M^*_{i,j})\big]
     \label{eq: likelihood function}
\end{equation}
with $M^*_{i,j}\coloneqq \zeta_j^* + {\mb X_i^*}\t \mb Y_j^*$ for every $(i,j)\in \Omega$. Throughout, we adopt a uniform sampling model in which each entry belongs to $\Omega$ independently with probability $\pi$, which may \emph{decay} with the dimensions $n$ and $p$.
As a concrete example, consider the logistic-link model for binary entries, under which the responses are conditionally independent and satisfy
    \begin{align}\label{eq: bernoulli model example}
        & R_{i,j} \mid M_{i,j}^* \sim \sf{Bern}\Big(\frac{\exp\big(M_{i,j}^*\big)}{1 + \exp\big(M_{i,j}^*\big)}\Big), \qquad (i,j)\in\Omega,
    \end{align}
    and only entries with $(i,j) \in \Omega$ are observed. This model includes the one-bit matrix completion problem as a special case \citep{davenport20141}. 
As is standard in the latent variable literature, we impose a low-rank structure through the natural-parameter matrix $\mb M^*$. In contrast to the linear setting, however, this low-rank natural-parameter matrix is observed \emph{only indirectly} through the nonlinear conditional-mean map $\psi\coloneqq\Psi'$, which constitutes one main challenge of the problem. In particular, this nonlinear formulation undermines two main strategies commonly used in prior work and makes them inapplicable: (i) obtaining the maximizer of \eqref{eq: likelihood function} in closed form, and (ii) conducting spectral analysis of the data matrix under a “low-rank signal plus noise” decomposition. Therefore, this problem calls for a new computational and theoretical framework.

Existing theory addresses important pieces of this problem separately: inference under missingness is substantially better developed for linear low-rank models, whereas recent theory for nonlinear latent factor models often characterizes likelihood-based estimators without establishing that a computationally tractable algorithm reaches an inferentially regular solution. Our goal is to bridge these statistical and computational gaps. We show that a fully data-driven procedure can be initialized in the appropriate local region, that vanilla gradient descent remains there through implicit regularization and, most importantly, that gradient descent's actual output admits uniform rowwise linear expansions that directly enable optimal statistical inference.

    \subsection{Local Landscape and Downstream Inference}
    In view of the objective \eqref{eq: likelihood function}, a natural approach is to characterize its global maximizer \citep{wang2022maximum,li2023statistical,ouyang2024statistical}. This perspective, however, leaves fundamental statistical and computational issues unresolved. Under severe missingness, the global likelihood maximizer may lie in a spurious region with poor or even ill-posed statistical properties, while computing it is generally intractable. We therefore focus on the output of a concrete optimization procedure:gradient descent, and seek to establish statistical guarantees through an analysis of its optimization dynamics.

The success of gradient descent hinges on the local geometry of the empirical objective. Under full observation, suitably regularized likelihoods often inherit the geometry of their population counterparts. After accounting for the nonidentifiable directions, restricted strong convexity and smoothness hold throughout an appropriate $\ell_2$ neighborhood \citep{ma2020universal,li2023statistical}, with empirical curvature comparable to its population counterpart along every admissible direction. 
Missingness makes this transfer substantially more delicate. The sampling operator can distort the empirical curvature outside the region of incoherence and introduce spurious local minima even in linear models. Consequently, proximity in a global norm alone does not guarantee that the gradient-descent trajectory lies in a region of contraction.

    Prior work on simpler low-rank signal-plus-noise problems has nevertheless revealed an \emph{implicit-regularization} phenomenon \citep{ma2018implicit,chen2019noisy,chen2021bridging,chi2019nonconvex}. Once initialized in a \emph{region of incoherence and contraction}, gradient-based methods remain sufficiently incoherent, encounter population-like curvature, and contract toward the optimal error rate. 
    This motivates our first question:
    \begin{quote}
        \emph{Question 1: Under both nonlinearity and missingness, can one characterize a region in which the empirical loss retains population-like curvature and gradient descent contracts?} 
    \end{quote}

    Characterizing such a region is only part of the challenge, since the result becomes operational only if a data-driven procedure can reach it. This leads to a second question:
    \begin{quote}
         \emph{Question 2: Can a computationally tractable initialization scheme reach this region and furnish a good start for gradient descent?}
    \end{quote}
    The initialization problem is nontrivial. Spectral estimators used for linear low-rank models \citep{ma2018implicit} are no longer suitable because the nonlinear mean matrix generally exhibits rank inflation and is no longer low-rank. Although link-aware spectral procedures can provide global recovery guarantees \citep{ma2020universal,zhang2020note}, their existing analyses do not furnish the rowwise and leave-one-out controls required to study the subsequent gradient-descent trajectory.

    Finally, reaching the \emph{region of incoherence and contraction} matters not only for optimization. The stationarity equations at a sufficiently regular solution may yield a tractable stochastic expansion, opening the door to uncertainty quantification for the recovered latent structure. This raises our third question:
    \begin{quote}
        \emph{Question 3: Can one derive a distributional approximation for the output of gradient descent and use it to conduct valid statistical inference?}
    \end{quote}

    \subsection{Linear Approximations for Gradient Descent} 
    We expand on \emph{Question 3} in the following. Classical M-estimation provides a useful starting point. Consider a \emph{fixed-dimensional} smooth convex problem with random objective $L_n$, let $\bo\theta^*$ denote its population target, and let $\hat{\bo\theta}_n$ be a stationary point. Expanding the score along the line segment from $\bo\theta^*$ to $\hat{\bo\theta}_n$ gives
    \begin{align}
        \hat{\bo\theta}_n-\bo\theta^*
        &=-\Big\{\int_0^1\nabla^2L_n\big(\bo\theta^*+s(\hat{\bo\theta}_n-\bo\theta^*)\big)\,\mathrm ds\Big\}^{-1}\nabla L_n(\bo\theta^*)\approx-\Big\{\nabla^2\bb E[L_n](\bo\theta^*)\Big\}^{-1}\nabla L_n(\bo\theta^*),
        \label{eq: linear approximation heuristic}
    \end{align}
    provided that the averaged empirical Hessian is close to its population counterpart. The distributional behavior of $\hat{\bo\theta}_n-\bo\theta^*$ is then governed by the score $\nabla L_n(\bo\theta^*)$.
    
    At first sight, one might expect the same reasoning to apply once an approximate stationary point enters the \emph{region of incoherence and contraction}, with $-\big\{\nabla^2\bb E[L_n](\bo\theta^*)\big\}^{-1}\nabla L_n(\bo\theta^*)$ serving as a proxy for $\hat{\bo\theta}_n-\bo\theta^*$.\footnote{Here, $\hat{\bo\theta}_n$ and $\bo\theta^*$ generically represent estimated and population latent parameters, including intercepts and factors.} Extending this heuristic to high dimensions is considerably more delicate, particularly because factor parametrizations are rotationally nonidentifiable. Several recent works \citep{wang2022maximum,li2023statistical,ouyang2024statistical} have characterized global minimizers of suitably regularized maximum likelihood estimators. Our goal is instead to connect tractable computation directly with valid inference: we aim to derive rowwise linear expansions and distributional approximations for the output of a computationally tractable procedure, while proving that the algorithm reaches the regime in which these conclusions hold. The resulting conclusions are fundamentally statistical, but they are accompanied by end-to-end algorithmic guarantees.

    A linear approximation is useful not only for characterizing the fluctuations of individual estimates, but also for quantifying uncertainty in the many comparisons that arise in latent factor analysis. Examples include ranking subjects by latent ability, comparing items by difficulty or discriminatory power, and assessing unobserved entries \citep{fan2023spectral,fan2025ranking}. Such tasks involve large collections of dependent contrasts and therefore cannot be justified by pointwise consistency or marginal normal approximations alone. 

    \subsection{Main Contributions}
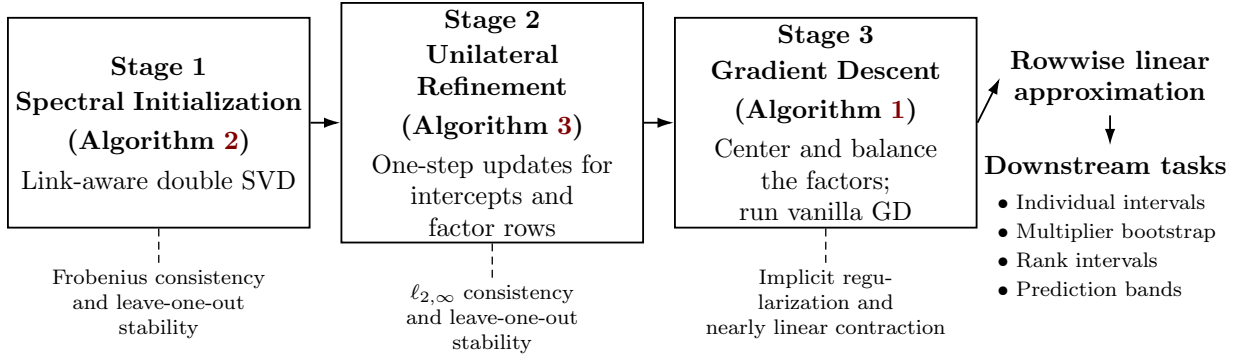
\begin{figure}[tb]
\centering
\resizebox{\linewidth}{!}{%
\begin{tikzpicture}[
    node distance=4mm,
    box/.style={
        draw=black,
        line width=0.7pt,
        text width=3.85cm,
        minimum height=2.8cm,
        align=center,
        inner sep=2.5pt,
        font=\small
    },
    description/.style={
        text width=3.45cm,
        align=center,
        inner sep=1pt,
        font=\scriptsize
    },
    flow/.style={
        -{Latex[length=2.2mm,width=1.4mm]},
        line width=0.7pt,
        draw=black
    }
]

\node[box] (spectral) {
    \textbf{Stage 1}\\[0.5mm]
    \textbf{Spectral Initialization}\\[1mm]
    \textbf{(Algorithm~\ref{alg: USVT})}\\[1mm]
    Link-aware double SVD
};

\node[box, right=of spectral] (refinement) {
    \textbf{Stage 2}\\[0.5mm]
    \textbf{Unilateral Refinement}\\[1mm]
    \textbf{(Algorithm~\ref{alg: OS})}\\[1mm]
    One-step updates for\\
    intercepts and factor rows
};

\node[box, right=of refinement] (gd) {
    \textbf{Stage 3}\\[0.5mm]
    \textbf{Gradient Descent}\\[1mm]
    \textbf{(Algorithm~\ref{alg: GD})}\\[1mm]
    Center and balance\\
    the factors; run vanilla GD
};

\node[align=center, anchor=west, font=\small] (linear)
    at ([xshift=3mm,yshift=6mm]gd.east) {
    \textbf{Rowwise linear}\\
    \textbf{approximation}
};

\node[align=center, below=4mm of linear, font=\small] (inference) {
    \textbf{Downstream tasks}
};

\node[description, text width=3cm, align=left, below=1mm of inference] (tasks) {
    $\bullet$ Individual intervals\\[0.5mm]
    $\bullet$ Multiplier bootstrap\\[0.5mm]
    $\bullet$ Rank intervals\\[0.5mm]
    $\bullet$ Prediction bands
};

\draw[flow] (spectral) -- (refinement);
\draw[flow] (refinement) -- (gd);
\draw[flow] (gd.east) -- (linear.west);
\draw[flow] (linear.south) -- (inference.north);

\node[description, below=5mm of spectral] (spectral-result) {
    Frobenius consistency\\
    and \mbox{leave-one-out} stability
};

\node[description, below=5mm of refinement] (refinement-result) {
    $\ell_{2,\infty}$ consistency\\
    and \mbox{leave-one-out} stability
};

\node[description, below=5mm of gd] (gd-result) {
    Implicit regularization and\\
    nearly linear contraction
};

\draw[densely dashed, line width=0.5pt] (spectral.south) -- (spectral-result.north);
\draw[densely dashed, line width=0.5pt] (refinement.south) -- (refinement-result.north);
\draw[densely dashed, line width=0.5pt] (gd.south) -- (gd-result.north);

\end{tikzpicture}
}
\caption{Overview of the estimation and inference pipeline. Each estimation stage builds on the preceding stage, with its theoretical guarantees summarized below the corresponding box. Rowwise linear approximations of the gradient descent output enable the downstream inference procedures shown on the right.}

\label{fig: estimation and inference pipeline}
\end{figure}

Building on a nonconvex optimization framework, this work advances the understanding of gradient descent dynamics under both nonlinear observations and missing data, while addressing the associated initialization challenge through a carefully designed multi-stage pipeline. Beyond estimation, we show that the output of the computational pipeline itself supports optimal individual and simultaneous statistical inference. Our affirmative answers to \emph{Questions~1--3} and technical contributions are summarized in the following principal contributions. 
\begin{enumerate}
\item We characterize a \emph{region of incoherence and contraction} for generalized latent factor models with missing data, thereby answering \emph{Question~1}. 
It is further revealed that the gradient descent iterate remains in the \emph{region of incoherence and contraction} and converge at a nearly linear rate to a statistically accurate stationary point, without explicit incoherence or balancing regularization. 
Our guarantees accommodate severe missingness, weak low-rank signals, and diminishing local curvature of the link. We show that the resulting aggregate estimation rate matches the restricted-class minimax lower bound up to logarithmic factors.

\item To warm-start gradient descent, we analyze a nonlinear spectral procedure that incorporates knowledge of the link function. We establish both a Frobenius-norm recovery guarantee and leave-one-out error control, yielding globally consistent estimators of $\bo\zeta^*$, $\mb X^*$, and $\mb Y^*$. Building on this spectral estimator, we show that a unilateral refinement step upgrades Frobenius-norm consistency to entrywise $\ell_{2,\infty}$ consistency. These guarantees place the refined estimator within the \emph{region of incoherence and contraction}, thereby providing the computationally tractable initialization sought in \emph{Question~2}.

\item We establish \emph{rowwise linear approximations} for the output of gradient descent in nonlinear latent factor models with missing observations. These expansions isolate the leading score fluctuations from higher-order estimation and optimization errors, thereby making the output of a nonconvex algorithm amenable to statistical inference. They yield asymptotically valid individual inference, with limiting covariance attaining the corresponding oracle Cram\'er--Rao lower bounds, without requiring an additional debiasing step. 
Building on the same expansions, we develop a Gaussian multiplier bootstrap for simultaneous inference on latent factors and missing entries with theoretical guarantees. 

\item On the technical side, we develop a leave-one-out analysis spanning \emph{the entire algorithmic pipeline}, which yields fine-grained control of the iterative estimators despite the dependencies induced by missing observations. Building on this control, we establish the rowwise linear approximation through a delicate analysis of the global Hessian structure, complemented by an auxiliary-sequence argument.

\end{enumerate}
Figure~\ref{fig: estimation and inference pipeline} provides a schematic summary of the three-stage computational pipeline, their corresponding theoretical guarantees, and downstream statistical inference results. 

\paragraph{Paper organization.} The remainder of the paper is organized as follows. Section~\ref{sec: gd theory} presents the main estimation and linear approximation results for warm-started gradient descent, together with a discussion of the underlying assumptions. Section~\ref{sec: prior algorithmic design} describes the prior two-stage algorithmic design for warm start. Section~\ref{sec: inference} develops the downstream statistical inference procedures.  Section~\ref{sec: numerical simulations} provides numerical experiments that illustrate the practical performance of our methods. Section~\ref{sec: real data application} applies our methods to evaluate large language models with uncertainty quantification. Section~\ref{sec: discussion} concludes with a discussion of future directions. All proofs are deferred to the Supplementary Material.

\section{Warm-Started Gradient Descent: Estimation and Linear Approximations}
\label{sec: gd theory}
As described in the diagram (Figure~\ref{fig: estimation and inference pipeline}), the overall procedure consists of three stages, each serving a distinct role.
For clarity, however, we present the algorithmic analyses in an order reverse to that of the algorithm itself.
We first present the estimation guarantees and linear approximation properties for the gradient descent algorithm, from which we develop insights into the region of incoherence and contraction needed for the warm start. Then, we turn to the earlier stages of the pipeline that lead to the warm start. Notably, the output of each stage enjoys a progressively stronger form of consistency, ultimately meeting the requirements for initializing gradient descent.

\subsection{Notation}
\label{subsec: notation}

For any positive integer $q$, let $[q]\coloneqq\{1,\ldots,q\}$. We use $\mb 1_q$, $\mb 0_q$, and $\mb I_q$ for the all-ones vector, the zero vector, and the identity matrix of dimension $q$, respectively. For $\mb A\in\bb R^{q_1\times q_2}$, $\mb A_{i,\cdot}\in\bb R^{1\times q_2}$ denotes its $i$th row in row-vector form, whereas $\mb A_i\coloneqq\mb A_{i,\cdot}^{\top}\in\bb R^{q_2}$ denotes the same row in column-vector form; similarly, $\mb A_{\cdot,j}$ denotes its $j$th column. For a vector $\mb a$, $\|\mb a\|_2$ denotes its Euclidean norm. For a matrix $\mb A$, we write $\|\mb A\|$, $\|\mb A\|\fb$, and $\|\mb A\|_{2,\infty}\coloneqq\max_{i\in[q_1]}\|\mb A_i\|_2$ for its spectral, Frobenius, and maximum rowwise Euclidean norms, respectively, and denote its $k$th largest singular value by $\sigma_k(\mb A)$. The symbols $(\cdot)^\top$ and $\circ$ represent the transpose and Hadamard product. For $\Omega\subseteq[n]\times[p]$, the sampling operator $\mc P_\Omega$ is defined by $[\mc P_\Omega(\mb A)]_{ij}\coloneqq A_{ij}\ind\{(i,j)\in\Omega\}$, or equivalently $\mc P_\Omega(\mb A)=\bo\Omega\circ\mb A$, where $\bo\Omega$ is the associated binary mask. The operator $\mc P^\infty_a$ denotes coordinatewise clipping to
$[-a,a]$. The map $\mc P^{2,\infty}_a(\mb A)$ rescales each row of $\mb A$
whose $\ell_2$ norm exceeds $a$ so that its norm equals $a$.
For $\bo\theta=[\bo\zeta,\mb X,\mb Y]$, define $\mc M(\bo\theta)\coloneqq\mb 1_n\bo\zeta^\top+\mb X\mb Y^\top$. We write $\mc O(r)\coloneqq\{\mb O\in\bb R^{r\times r}:\mb O^\top\mb O=\mb I_r\}$ for the orthogonal group, and set $d\coloneqq n\vee p$ and $\rho\coloneqq(n\vee p)/(n\wedge p)$. Finally, $a_n\lesssim b_n$ means $a_n\leq Cb_n$ for some constant $C>0$, $a_n\asymp b_n$ means that both $a_n\lesssim b_n$ and $b_n\lesssim a_n$, and $a_n \ll b_n$ means that there exists a sufficiently small $c>0$ such that $a_n \leq c b_n$.

\subsection{Model Assumptions}
    We first introduce the basic model setup, which accommodates a wide class of continuous and discrete distributions. 
    Let $\mb U^*\bo \Sigma^*\mb V^{*\top}$ denote the rank-$r$ SVD of $\mb X^*{\mb Y^*}\t$, with singular values $\sigma_1^*\geq\cdots\geq\sigma_r^*$, and define $\kappa\coloneqq\sigma_1^*/\sigma_r^*$. We introduce the incoherence degree $\mu$ of the singular vectors as follows:
    \begin{definition}
        The parameters $[\bo \zeta^*, \mb X^*, \mb Y^*]$ are said to be $\mu$-incoherent if $
            \frac{n}{r}\norm{\mb U^*}\ti^2  \vee \frac{p}{r} \norm{\mb V^*}\ti^2 \leq \mu$ and $\norm{\bo \zeta^* - \zeta_0 \cdot \mb 1_p}_\infty \leq \sqrt{\mu r \sigma_r^* / np}$ for some reference level $\zeta_0$.  
    \end{definition}
    We impose the following assumptions on sampling, the response distributions, the link function, and the signal.

    \begin{assumption}[Sampling scheme and noise assumptions]
    \label{assumption: main noise and sampling}
    \begin{itemize}
        \item[(a)] Each entry is independently observed with probability $\pi \in (0,1]$. 
        \item[(b)] The observation follows the formulation in \eqref{eq: likelihood function}. For some constant $c>1$, the entries $R_{i,j}$ are mutually independent and satisfy
        $\max_{i,j}\mathrm{Var}(R_{i,j})^{1/2}\vee d^{-c-10}\leq \sigma$ and $|\log\sigma|\lesssim\log d$. Moreover, with probability at least $1-O(d^{-c-10})$, $\max_{i,j}|R_{i,j}|\leq B$, where, without loss of generality, $B\geq \max_{x\in\mc D_\psi}|\psi(x)|$.
    \end{itemize}
    \end{assumption}

    \begin{assumption}[Link function conditions]
    \label{assumption: main link function}
    Let $\psi\coloneqq\Psi'$.
        Set $\mc D_\psi\coloneqq[\zeta_0-c_{\mc D_\psi}\mu r\sigma_1^*/\sqrt{np},\allowbreak\zeta_0+c_{\mc D_\psi}\mu r\sigma_1^*/\sqrt{np}]$ for a sufficiently large fixed constant $c_{\mc D_\psi} \geq 1$ as the domains of both $\Psi$ and $\psi$.  A possibly diverging or vanishing $\underline c_\psi$ and a constant-order smoothness parameter $L_\psi$ satisfy, uniformly over $x\in\mc D_\psi$,
        \begin{equation*}
            \underline c_\psi\leq\psi'(x)\leq\sigma^2,
            \qquad
            |\psi^{(k)}(x)|\leq\sigma^2 k!L_\psi^{k-1},
            \quad k\in\bb N^+.
        \end{equation*}
        Define $\kappa_\psi\coloneqq\sigma^2/\underline c_\psi$.
    \end{assumption}
     \begin{remark}
        Assumption~\ref{assumption: main noise and sampling} encompasses the regular one-parameter exponential-family models commonly used in generalized latent factor analysis
\citep{moustaki2000generalized,chen2020structured}, including Bernoulli, exponential, and Poisson. For unbounded distributions, $B$ depends on the tail behavior, with its relation to $\sigma$ specified later in Assumption~\ref{assumption: main signal strength and incoherence degree}. In the Bernoulli case, for example, $B = 1$, $\psi'(x) = e^x / (1 + e^x)^2$, and $\underline c_\psi = \min_{x\in \mc D_\psi} e^x / (1 + e^x)^2$, with the corresponding smoothness parameters given in Lemma~\ref{lemma: bernoulli smoothness}.
    \end{remark}
    \begin{remark}
        The identity $\mathrm{Var}(R_{i,j})=\psi'(M_{i,j}^*)$ explains our use of $\sigma^2$ as a common scale for the response variances at the truth and the derivatives of $\psi$ over $\mc D_\psi$. The level of $\psi$ is controlled separately through $\bar\sigma$ in Assumption~\ref{assumption: main signal strength and incoherence degree}.
    \end{remark}

    We impose the following simplified conditions on the low-rank signal in the natural-parameter scale by letting the less critical parameters $r,\rho, \kappa, \kappa_\psi$ be absorbed into the constant factors.  The parameter-explicit assumptions are collected in Appendix~\ref{sec: supplementary preliminary setups}. 
    \begin{assumption}[Signal conditions]
    \label{assumption: main signal strength and incoherence degree}
    We assume $r,\rho, \kappa, \kappa_\psi = O(1)$ and $\xi=O((\log d)^{r/2})$ and that the following conditions hold: 
    \begin{itemize}
        \item[(a)] \emph{Signal Strength: } Define $\bar \sigma \coloneqq \sigma \vee \max_{x \in \mc D_\psi} |\psi(x)|$. Assume $\sigma_r^*\gtrsim1$ and
        \begin{align}
            \sigma_r^*
            &\gg \frac{\bar \sigma
            \mu^{\frac52}(\log d)^{\frac{5r}{2}+3}}{\underline c_\psi}
            \sqrt{\frac{d}{\pi}}, \qquad
            \pi(n\wedge p)
            \gg \mu^3(\log d)^{4r}.
            \label{eq: consolidated sampling condition}
        \end{align}
        \item[(b)] \emph{Tail condition: }
        The relation between $\sigma$ and $B$ is given by $\mu B \sqrt{\log d} \ll \sigma \sqrt{\pi(n \wedge p)}$.
    \end{itemize}
    \end{assumption}

    \begin{remark}
        The conditions on the incoherence degree are standard in the vast literature of linear low-rank estimation and are generally necessary for entrywise estimation and inference; to name a few, \cite{candes2012exact,chen2015incoherence,abbe2020entrywise}. 
    \end{remark}

    \begin{remark}
    If, in addition, $\mu$ is of constant order, then for Gaussian matrix completion with $\Psi(x)=\sigma^2 x^2/2$, our signal-strength condition has the same signal-to-noise scaling, up to polylogarithmic factors, as standard noisy matrix-completion guarantees \citep{chen2019noisy}, namely, $\sigma_r(\bb E[\mb R]) \stackrel{\text{up to polylogs}}{\gtrsim} \sigma \sqrt{d / \pi}.$
        More broadly, 
        our assumptions are designed to cover several weak-signal regimes that frequently arise in incomplete data with nonlinear links. In particular, the theory accounts for: 
        \begin{itemize}
            \item weak low-rank signal strength, with $\sigma_r^*$ approaching $\bar \sigma\sqrt{d / \pi} / \underline c_\psi$ up to logarithmic factors;
            \item severe missingness, with $\pi$ as small as $(n\wedge p)^{-1}$ up to logarithmic factors; 
            \item weak local curvature of the link, with $\underline c_\psi=o(1)$. 
        \end{itemize}
        These regimes may overlap, and together they capture settings in which both statistical recovery and algorithmic stability are delicate.  
    \end{remark}

    \paragraph{Identifiability conditions.}
    \label{subsec: identifiability}
    Before proceeding, we note that the parameters in the generalized latent factor model are not uniquely identifiable without further conditions, just as their linear factor model counterparts. To disentangle the intercepts from the latent factors, we treat $\boldsymbol\zeta^*$ as the average response level and impose the following assumptions on the intercept and the factors: 
    \begin{assumption}[Basic Identifiability]
    \label{assumption: main identifiability}
    The factors $\mb X^*$, $\mb Y^*$ and the intercept vector $\bo \zeta^*$ obey: 
        \begin{align}
      \mb 1\t \mb X^* &= \mb 0,     \label{eq: identifiability condition 1} \\ 
    {\mb X^*}\t \mb X^* &=  \omega {\mb Y^*}\t \mb Y^* \label{eq: identifiability condition 2}. 
\end{align} 
    \end{assumption}
Following the discussions in \cite{bai2012statistical} and \cite{fan2021robust}, the parameter triple $[\bo \zeta^* , \mb X^*, \mb Y^*]$ is identifiable only up to an orthogonal rotation; that is, $(\mb X^*, \mb Y^*)$ and $(\mb X^* \mb R, \mb Y^* \mb R)$ produce the same distribution for any $\mb R \in \mc O(r)$. In particular, the first condition in \eqref{eq: identifiability condition 1} centers the columns of $\mb X^*$, which renders $\bo \zeta^*$ identifiable. The second condition in \eqref{eq: identifiability condition 2} fixes the relative scaling between $\mb X^*$ and $\mb Y^*$. The scaling factor $\omega>0$ is set \emph{a priori} to accommodate Bayesian settings in which the latent factors naturally operate on different scales. For instance, if the entries of $\mb X^*$ and $\mb Y^*$ are independently drawn from their respective priors, then one typically has $\omega \asymp n/p$ since $\sigma_k(\mb X^*) \asymp \sqrt{n}$ and $\sigma_k(\mb Y^*) \asymp \sqrt{p}$ for $k\in[r]$. One may further strengthen the identifiability by removing the rotational ambiguity, as we shall elaborate in Section~\ref{subsec: gd theory}.

    \subsection{Gradient Descent and Theoretical Guarantees}
    \label{subsec: gd theory}
    As all preceding steps (cf.~Algorithms~\ref{alg: USVT}~and~\ref{alg: OS}) 
    are designed to guarantee the success of gradient descent, we begin by identifying the initialization conditions it requires and the consequences that follow. For a centering parameter $c_\perp>0$, consider
    \begin{equation}
        \label{eq: loss function}
     L([\bo \zeta, \mb X, \mb Y]) = \sum_{(i,j) \in \Omega}\big[-R_{i,j}(\zeta_j + \mb X_i\t \mb Y_j) + \Psi(\zeta_j + \mb X_i\t \mb Y_j)\big] + c_\perp\frac{\pi\sigma_r^*}{n\sqrt{\omega}}\bignorm{\mb 1_n\t\mb X}_2^2.
    \end{equation} 
    This objective augments the negative log-likelihood by a quadratic penalty that enforces the centering condition in \eqref{eq: identifiability condition 1}. Algorithm~\ref{alg: GD} summarizes the resulting gradient-descent procedure.
    \begin{algorithm}[tb]
\caption{Vanilla Gradient Descent (GD)}
\label{alg: GD}
\KwIn{Refined estimates $[\widetilde{\bo\zeta}^{0},\widetilde{\mb X}^{0},\widetilde{\mb Y}^{0}]$, scale $\omega>0$, centering parameter $c_\perp>0$, positive block step sizes $(\eta_\zeta,\eta_{\mb X},\eta_{\mb Y})$, and iteration count $t_0\in\bb N$}
\textbf{Centering and Balancing: }
Set $\mb J_n\coloneqq\mb I_n-n^{-1}\mb 1_n\mb 1_n\t$, 
$
    \widetilde{\mb X}^{0,\sf c}\coloneqq\mb J_n\widetilde{\mb X}^{0}$, and $
    \bo\zeta^0\coloneqq\widetilde{\bo\zeta}^{0}$. 
Compute the top-$r$ SVD
$\widetilde{\mb X}^{0,\sf c}\widetilde{\mb Y}^{0\t}
=\mb U^0\bo\Sigma^0\mb V^{0\t}$, and set
\[
    \mb X^0\coloneqq\omega^{1/4}\mb U^0(\bo\Sigma^0)^{1/2},\qquad
    \mb Y^0\coloneqq\omega^{-1/4}\mb V^0(\bo\Sigma^0)^{1/2}.
\]
\textbf{Gradient Descent Update: }\For{$t=1,2,\ldots,t_0$}{
    \begin{align*}
        \bo \zeta^{t} &= \bo \zeta^{t-1} - \eta_\zeta
        \nabla_{\bo \zeta}L([\bo \zeta^{t-1},\mb X^{t-1},\mb Y^{t-1}]),\\
        \mb X^{t} &= \mb X^{t-1} - \eta_{\mb X} \nabla_{\mb X}L([\bo \zeta^{t-1},\mb X^{t-1},\mb Y^{t-1}]),\\
        \mb Y^{t} &= \mb Y^{t-1} - \eta_{\mb Y}\nabla_{\mb Y}L([\bo \zeta^{t-1},\mb X^{t-1},\mb Y^{t-1}]).
    \end{align*}
}

\KwOut{The factor estimates $\bo \zeta^{t_0}$, $\mb X^{t_0}$, $\mb Y^{t_0}$, and the augmented matrix estimate
$
\hat{\mb M}= \mb 1_n \bo\zeta^{t_0\top} + \mb X^{t_0} \mb Y^{t_0\top}
$. 
}
\end{algorithm}

The centering-and-balancing step in Algorithm~\ref{alg: GD}
leaves the intercept estimate unchanged and replaces the
low-rank component by its centered and balanced factorization:
\begin{align*}
    \mb 1_n\bo\zeta^{0\t}+\mb X^0\mb Y^{0\t}
    &=\mb 1_n\widetilde{\bo\zeta}^{0\t}
      +\mb J_n\widetilde{\mb X}^{0}\widetilde{\mb Y}^{0\t},\qquad 
    \mb 1_n\t\mb X^0=\mb 0,\qquad
    {\mb X^0}\t\mb X^0
    =\omega\,{\mb Y^0}\t\mb Y^0.
\end{align*}

\emph{Regularization Choice. } We briefly explain why the objective contains neither an incoherence penalty nor a balancing penalty. As noted in the introduction, rowwise consistency prevents the iterates from entering the unfavorable parts of the landscape created by missingness. One could enforce this property by adding $\norm{\mb X}_{2,\infty}^2+\norm{\mb Y}_{2,\infty}^2$, but such a penalty would complicate the optimization. Inspired by \citet{ma2018implicit}, we instead prove that vanilla gradient descent enjoys \emph{implicit regularization}, in that it implicitly preserves rowwise consistency when initialized from a rowwise-consistent estimate.
Prior work also commonly adds the balancing penalty $\norm{\mb X^\top \mb X-\mb Y^\top \mb Y}\fb$ when $\omega=1$ to control the relative scales of the two factors. Algorithm~\ref{alg: GD} imposes balance only once, through an exact refactorization of the warm start. No balancing penalty is used thereafter. We show that the scaled learning rates preserve approximate balance along the ensuing trajectory. Section~\ref{subsec: mechanism} describes these two forms of implicit regularization in detail.

\bigskip

We are now positioned to present our estimation error control for the gradient
descent iterates.

\begin{theorem}
\label{thm: GD}
Suppose that Assumptions~\ref{assumption: main noise and sampling},
\ref{assumption: main link function},
\ref{assumption: main signal strength and incoherence degree}, and
\ref{assumption: main identifiability} hold. We follow the initialization procedures depicted in Figure~\ref{fig: estimation and inference pipeline} (with the tuning in Theorems~\ref{thm: USVT} and~\ref{thm: OS}), let the penalty parameter $c_\perp$ satisfy
$\sigma^2
    \vee(\log d)^4
    \ll c_\perp
    \ll\frac{\underline c_\psi d^{10}}{\log d},
$
and determine the learning rates as 
$$
    \qquad \eta \coloneqq\frac{1}{8c_\perp\pi\sigma_1^*},
    \quad \eta_\zeta\coloneqq\eta\frac{\sigma_r^*}{n},\quad
    \eta_{\mb X}\coloneqq\eta\sqrt{\omega},\quad
    \eta_{\mb Y}\coloneqq\frac{\eta}{\sqrt{\omega}}. 
$$
 For each
$0\leq t\leq t_0 \coloneqq\Big\lfloor
    \frac{c_t\log d}{\underline c_\psi\pi\sigma_r^*\eta}
    \Big\rfloor$ with a sufficiently large $c_t>0$, let $\mb O^t$ be a measurable choice satisfying
$
    \mb O^t\in\argmin_{\mb O\in\mc O(r)}\Big\{
    \omega^{-1/2}\bignorm{\mb X^t\mb O-\mb X^*}\fb^2
    +\omega^{1/2}\bignorm{\mb Y^t\mb O-\mb Y^*}\fb^2
    \Big\}$.
Then, with probability at least $1-O(d^{-c})$, the final iterate satisfies
\begin{align*}
    &\bignorm{\bo\zeta^{t_0}-\bo\zeta^*}_2
    \lesssim
    \frac{\sigma\sqrt{d\log d}}
    {\underline c_\psi\sqrt{\pi n}},\quad \omega^{-\frac14}
    \bignorm{\mb X^{t_0}\mb O^{t_0}-\mb X^*}\fb
    \vee
    \omega^{\frac14}
    \bignorm{\mb Y^{t_0}\mb O^{t_0}-\mb Y^*}\fb
    \lesssim
    \frac{\sigma\sqrt{d\log d}}
    {\underline c_\psi\sqrt{\pi\sigma_r^*}}.
\end{align*}
\end{theorem}
    
    \begin{remark}
        The error bounds follow from a nearly linear one-step contraction.
        Define the jointly aligned and scaled error
        $\bo\Delta^t\coloneqq[\sqrt{n/\sigma_r^*}(\bo\zeta^t-\bo\zeta^*),
        \omega^{-1/4}(\mb X^t\mb O^t-\mb X^*),
        \omega^{1/4}(\mb Y^t\mb O^t-\mb Y^*)]$. For some constants $c,C>0$, $
            \bignorm{\bo\Delta^t}_2
            \leq
            \big(1-c\eta\underline c_\psi\pi\sigma_r^*\big)
            \bignorm{\bo\Delta^{t-1}}_2
            +C\eta\sigma\sqrt{\pi dr\sigma_1^*\log d}$.
        Iterating this relation gives the geometric decay and the
        stochastic-error floor in Theorem~\ref{thm: GD}; the proof also
        establishes the corresponding rowwise recursion.
        \label{remark:linear-convergence-of-errors}
    \end{remark}
    \begin{remark}
    In line with the implicit-regularization phenomenon identified by
    \citet{ma2018implicit}, our analysis permits a learning rate close to the
    most aggressive rate allowed by population gradient descent. Specifically,
    Theorem~\ref{thm: GD} permits
    $\eta\asymp(c_\perp\pi\sigma_1^*)^{-1}$, whereas the largest stable learning
    rate for the population negative log-likelihood is of order
    $(\sigma^2\pi\sigma_1^*)^{-1}$. When $\sigma$ is also of constant order,
    the two rates differ only by a polylogarithmic factor
    in $d$.
    \end{remark}

    More importantly, beyond the estimation error control described above, we provide fine-grained linear approximations by quantifying their residual errors from the gradient descent iterates. 
    \begin{theorem}[Linear Approximation for Output]
    \label{thm: linear approximations}
        Suppose that the assumptions of Theorem~\ref{thm: GD} hold. Suppose
        further that
        \begin{align}
            & \sigma^2\pi(n\wedge p)
            \gg \mu^7(\log d)^{5r+6},\quad \sigma^2\mu^{\frac94}(\log d)^{\frac{5r}{4}+3}
            \sqrt{\bar\sigma\sigma_1^*}
            (\pi /d)^{\frac14}\lesssim c_\perp.
        \end{align}
        Let $\mb E\coloneqq\mb R-\bb E[\mb R]=\mb R-\psi(\mb M^*)$ denote the response-noise matrix.
        Set $(\eta_\zeta,\eta_{\mb X},\eta_{\mb Y})$ and $t_0$ according to Theorem~\ref{thm: GD}. Then, with probability at
        least $1-O(d^{-c})$, the final output admits the following uniform
        rowwise linear approximations in terms of the score noise:
            \begin{align}
                & \sqrt{\frac n{\sigma_r^*}}\max_{j\in[p]}\Big| \big(\zeta_j^{t_0} - \zeta_j^*\big) - \sqrt{\frac{\sigma_r^*}{n}}\,\mc P_{\Omega}(\mb E)_{\cdot, j}\t \mb X^{\sf{app},*} \big(\mb H^{*-1}_{\mb Y^{\sf{app}}_{j}}\big)_{\cdot,1} \Big|  \\
                & \vee\ \omega^{-\frac14}\max_{i\in[n]}\Bignorm{\big( \mb X^{t_0}_{i,\cdot } \mb O^{t_0} - \mb X^*_{i,\cdot} \big) -\mc P_{\Omega}(\mb E)_{i,\cdot} \mb Y^* \mb H^{*-1}_{\mb X_i}}_2 \\
                & \vee\ \omega^{\frac14}\max_{j\in[p]}\Bignorm{\big( \mb Y^{t_0}_{j,\cdot } \mb O^{t_0} - \mb Y^*_{j,\cdot} \big) - \mc P_{\Omega}(\mb E)_{\cdot, j}\t \mb X^{\sf{app},*} \big(\mb H^{*-1}_{\mb Y^{\sf{app}}_{j}}\big)_{\cdot, 2:r+1}}_2\lesssim\frac{1}{\sigma\sqrt{\pi\sigma_1^*}(\log d)^2},
                \label{eq: main linear approximation error}
            \end{align}
            where we define the augmented left factor $\mb X^{\sf{app},*} \coloneqq \big(\sqrt{\sigma_r^*/n} \mb 1_n, \mb X^* \big) \in \bb R^{n \times (r+1)}$ and the observed information matrices 
    \begin{align}
        & \mb H^*_{\mb X_i} \coloneqq \sum_{j\in[p],(i,j)\in\Omega} \psi'(M^*_{i,j}) \mb Y_j^* \mb Y_j^{*\top}\in\bb R^{r\times r}, \qquad \mb H^*_{\mb Y^{\sf{app}}_{j}} \coloneqq \sum_{i\in[n], (i,j)\in \Omega} \psi'(M_{i,j}^*) \mb X^{\sf{app},*}_{i} \mb X^{\sf{app},*\top}_{i}.  \label{eq: definition of fisher matrices}
    \end{align}
    Consequently, the final iterate obeys the sharp rowwise error bound
    \begin{align}
        &\sqrt{\frac{n}{\sigma_r^*}}
        \bignorm{\bo\zeta^{t_0}-\bo\zeta^*}_\infty
        \vee \omega^{-\frac14}
        \bignorm{\mb X^{t_0}\mb O^{t_0}-\mb X^*}\ti
        \vee \omega^{\frac14}
        \bignorm{\mb Y^{t_0}\mb O^{t_0}-\mb Y^*}\ti
        \lesssim
        \frac{\sigma\sqrt{\mu\log d}}
        {\underline c_\psi\sqrt{\pi\sigma_r^*}}.
        \label{eq: rowwise error from linear approximation}
    \end{align}
    \end{theorem}
    The parameter-explicit counterpart of Theorems~\ref{thm: GD}~and~\ref{thm: linear approximations}, which allows $\rho$, $\xi$, $r$, $\kappa$, and
$\kappa_\psi$ to vary, is given in 
    Section~\ref{sec: parameter-explicit gradient results}.
    \begin{remark}
        The bound in \eqref{eq: main linear approximation error} is stated at the level needed for the subsequent inference results. Lemma~\ref{lemma: gradient replacement 3} in the Supplementary Material provides the corresponding explicit rate.
    \end{remark}

    We highlight the following implications. 
    \begin{itemize}
        \item \emph{Sample complexity. } When the remaining structural quantities are bounded, our assumptions require $\pi(n\wedge p)$ to grow only polylogarithmically in $d$. Thus, the sampling requirement lies within logarithmic factors of the $1/(n\wedge p)$ degrees-of-freedom scale familiar from fixed-rank matrix completion \citep{candes2009power,candes2012exact}. In this low-sampling regime, missingness destroys local $\ell_2$ convexity outside the region analyzed in Section~\ref{subsec: mechanism}.
        \item \emph{Sparse scenarios. } Consider the Bernoulli model with the logistic link in the balanced regime
$n\asymp p$, and write $q_n\coloneqq\psi(\zeta_0)$. When
$\mu$ is fixed, $\sigma_r^*\asymp\sigma_1^*\asymp n$,
and the factor contribution is uniformly bounded, all success probabilities
are of order $q_n$. Under the sparse-link regularity condition,
$\sigma\asymp q_n^{1/2}$, $\underline c_\psi\asymp\sigma^2\asymp q_n$,
$L_\psi\asymp1$. The fixed-rank approximation factor contributes only a
polylogarithmic term. The threshold and approximation conditions are
then satisfied if $\pi n q_n\gg(\log n)^{5r+4}$. Under full observation model without missing data, the
estimation and rowwise linearization theory therefore permits
   $q_n=n^{-1}(\log n)^c$ for any fixed $c>5r+4$. This response regime is strictly sparser than
   the $q_n\asymp n^{-1+\epsilon}$ regime
considered by \citet{li2023statistical} and \cite{cui2026convexity}, whereas
\citet{ouyang2024statistical} assumes curvature bounded away from zero
and therefore excludes vanishing Bernoulli probabilities. Polylogarithmic
expected degrees have also been accommodated for Frobenius-norm estimation
by \citet{ma2020universal} and for spectral inference under a linear
random-dot-product link by \citet{rubin2022statistical}. The
distinction of our result is the combination of a nonlinear generalized
latent factor model, missing observations, uniform rowwise linearization,
and inferential guarantees.
        \item \emph{No need for debiasing. } 
        The linear approximation is readily seen to be mean zero. In contrast to convex programming approaches \citep{chen2019inference} for the linear matrix completion, our method does not require an additional debiasing step, since the underlying optimization procedure avoids bias-inducing regularization. 
    \end{itemize}

\paragraph{Rotation adjustment.}
The Procrustes matrix $\mb O^{t_0}$ in Theorem~\ref{thm: linear approximations} is used only for analysis and is neither required by the algorithm nor observable from the data. Certain downstream tasks in Section~\ref{sec: inference}, however, require an estimable rotation to identify individual factor coordinates. We present two structural conditions under which $\mb O^{t_0}$ can be replaced by a data-driven rotation. 
\begin{itemize}
    \item First, distinct singular values supply a canonical orientation. 
    Suppose that the population orientation is fixed by
\eq{
    \omega^{-1/2}{\mb X^*}\t\mb X^*
    =
    \omega^{1/2}{\mb Y^*}\t\mb Y^*
    =
    \bo\Sigma^*
    \coloneqq
    \diag(\sigma_1^*,\ldots,\sigma_r^*),
    \qquad
    \sigma_1^*>\cdots>\sigma_r^*>0.
    \label{eq: distinct singular value conditions}
}
Let
$
    \Delta_{\sf{eig}}
    \coloneqq
    \min_{k\neq\ell}|\sigma_k^*-\sigma_\ell^*|$,
and let $\hat{\mb O}^{\sf{eig}}$ contain the ordered eigenvectors of
$\omega^{-1/2}{\mb X^{t_0}}\t\mb X^{t_0}$. When
$\Delta_{\sf{eig}}$ dominates the Gram-matrix perturbation,
$\hat{\mb O}^{\sf{eig}}$ consistently estimates $\mb O^{t_0}$ up to
columnwise sign changes.
\item 
Alternatively, the rotation can be identified semiparametrically through fourth-order structure in the factor distribution by the Varimax criterion. Varimax is widely used in factor analysis to obtain a simple and interpretable loading structure \citep{bartholomew2008analysis,skrondal2004generalized,rohe2023vintage}. It is particularly natural for question loadings in LLM evaluation: each question is expected to probe only a few latent capabilities, suggesting an approximately sparse loading structure. 
In particular, for
$\mb M\in\bb R^{n\times r}$, define
\begin{align*}
    \mc V_n(\mb Q;\mb M)
    &\coloneqq
    \sum_{k=1}^r
    \Big\{
        \frac{1}{n}\sum_{i=1}^n(\mb M\mb Q)_{i,k}^4
        -
        \Big[\frac{1}{n}\sum_{i=1}^n(\mb M\mb Q)_{i,k}^2\Big]^2
    \Big\}. 
\end{align*}
When the entries of $\mb M$ are i.i.d. from a standardized leptokurtic
distribution, the population criterion is maximized precisely by signed
permutation matrices. 
Exploiting this identification, 
\citet{rohe2023vintage} established consistency of varimax applied to
the PCA estimator $\hat{\mb Z}$ in the linear latent model
$\bb E[\mb A\mid\mb Z,\mb Y]=\mb Z\mb B\mb Y^\top$. 

We apply the same principle directly to the gradient-descent output in
our nonlinear setting. Suppose that $\mb X^*=\mb J_n\mb Z$, where the
rows of $\mb Z$ are i.i.d. and their centered coordinates are
independent, standardized, sub-Gaussian, and leptokurtic. The population
criterion identifies the coordinate axes up to signed permutations,
while Theorem~\ref{thm: linear approximations} controls the discrepancy
of the sample criterion. Consequently, $\hat{\mb O}^{\sf var} \in\argmax_{\mb Q\in\mc O(r)}
    \mc V_n(\mb Q;\mb X^{t_0})$
consistently estimates $\mb O^{t_0}$ up to a signed permutation. The
same construction applies to $\mb Y^*$.
\end{itemize}
Formally, Proposition~\ref{prop:rotation-recovery} in Section~\ref{subsec: rotation-free factor estimation} of the Supplementary Material
quantifies the resulting rotation-estimation errors, which provide
admissible choices of $\delta_{\sf{rotation}}$ in
Theorems~\ref{thm: individual asymptotic normality}
and~\ref{thm: simultaneous inference}.

    \paragraph{Minimax lower bound.}
    We next give minimax lower bounds for the parameter estimation errors.
    The result applies to the canonical one-parameter exponential family in
    \eqref{eq: likelihood function}, provided that the Fisher information is
    uniformly comparable over the local parameter interval specified below.

    \begin{theorem}
    \label{thm: minimax lower bound}
    Let $1\leq r\leq r_0$ for a fixed positive integer $r_0$, and fix
    $s,\omega>0$. For $c_\zeta\in\bb R$, $c_1>0$, and a sufficiently large
    constant $c_2>0$ depending only on $r_0$, define
    \begin{align}
        &\mc S_{\sf E}(r,s,\omega)
        \coloneqq\Big\{[\bo\zeta,\mb X,\mb Y]:{}
        \ \bignorm{\bo\zeta-c_\zeta\mb 1_p}_\infty\leq c_1,
        \notag\\
        &\quad \mb 1_n\t\mb X=\mb 0_r\t,
        \quad \mb X\t\mb X=\omega\mb Y\t\mb Y,\quad \bignorm{\mb X}\ti\leq
        c_2\omega^{1/4}\sqrt{\frac{s}{n}},
        \quad \bignorm{\mb Y}\ti\leq
        c_2\omega^{-1/4}\sqrt{\frac{s}{p}},\notag\\
        &\quad \operatorname{rank}(\mb X\mb Y\t)=r,
        \quad 0.9s\leq\sigma_r(\mb X\mb Y\t)
        \leq\sigma_1(\mb X\mb Y\t)\leq1.1s\Big\}.
        \label{eq: minimax lower-bound class}
    \end{align}
    Suppose $c_1+c_2^2s/\sqrt{np}\leq B_0$ for a fixed $B_0>0$, and let
    $\mc D_{\sf E}\coloneqq
    [c_\zeta-c_1-c_2^2s/\sqrt{np},
    c_\zeta+c_1+c_2^2s/\sqrt{np}]$. Assume that $\mc D_{\sf E}$ lies
    in the interior of the natural-parameter space and that, uniformly in
    $n$ and $p$,
    $
        0<\underline c_{\psi,{\sf E}}
        \coloneqq\inf_{u\in\mc D_{\sf E}}\psi'(u)
        \asymp\sup_{u\in\mc D_{\sf E}}\psi'(u)
        \eqqcolon\sigma_{\sf E}^2$.
    Write $
        d_\omega\big((\widetilde{\mb X},\widetilde{\mb Y}),(\mb X,\mb Y)\big)
        \coloneqq\min_{\mb O\in\mc O(r)}\Big\{
        \omega^{-1/2}\bignorm{\widetilde{\mb X}\mb O-\mb X}\fb^2
        +\omega^{1/2}\bignorm{\widetilde{\mb Y}\mb O-\mb Y}\fb^2
        \Big\}^{1/2}$.
    There exists $c_{\sf{lb}}>0$
    such that, for all sufficiently large $n\wedge p$,
    \begin{equation}
    \begin{aligned}
        \inf_{\widehat{\bo\zeta}}
        \sup_{\bo\theta\in\mc S_{\sf E}(r,s,\omega)}
        \bb E_{\bo\theta}\bignorm{\widehat{\bo\zeta}-\bo\zeta}_2
        &\geq c_{\sf{lb}}\sqrt p\Big\{
        c_1\wedge
        \frac{\sigma_{\sf E}}
        {\underline c_{\psi,{\sf E}}\sqrt{\pi n}}\Big\},\\
        \inf_{(\widehat{\mb X},\widehat{\mb Y})}
        \sup_{\bo\theta\in\mc S_{\sf E}(r,s,\omega)}
        \bb E_{\bo\theta}
        d_\omega\big((\widehat{\mb X},\widehat{\mb Y}),(\mb X,\mb Y)\big)
        &\geq c_{\sf{lb}}\sqrt s\Big\{
        1\wedge\frac{\sigma_{\sf E}}{\underline c_{\psi,{\sf E}}}
        \sqrt{\frac{n\vee p}{\pi s^2}}\Big\}.
    \end{aligned}
    \label{eq: minimax lower bound}
    \end{equation}
    \end{theorem}

    If $c_\zeta=\zeta_0$, $\mc D_{\sf E}\subseteq\mc D_\psi$,
    $\sigma_{\sf E}\asymp\sigma$, and
    $\underline c_{\psi,{\sf E}}\asymp\underline c_\psi$, then the
    high-probability aggregate common-rotation upper bound implied by
    Theorem~\ref{thm: GD} has the same parameter scaling as the second
    restricted-class minimax bound above, up to logarithmic factors. 

    \subsection{Proof architecture for the gradient stage}
    Several highlights of our technical contributions of the gradient descent analysis are in order. We begin by characterizing the optimization landscape, and then discuss the key technical ingredients underlying the linear approximations. 
    \label{subsec: mechanism}
    \subsubsection{Local curvature and implicit regularization }
    \label{subsubsec: implicit regularization mechanism}
    We now shed light on the behavior of gradient descent iterates and on why a warm start is necessary. As a point of reference, consider the linear setting of low-rank matrix completion. The seminal work \cite{ma2018implicit} offered a new perspective through the lens of implicit regularization, showing that gradient descent can automatically constrain its trajectory to remain within the \emph{region of incoherence and contraction (RIC)}. This property makes it possible to obtain precise control of the Hessian along tangent directions at any point in the RIC. In this paper, we show that this implicit regularization phenomenon extends to generalized latent factor models with missingness. Moreover, in our analysis, implicit regularization manifests itself in two distinct senses: (i) locational regularization within the \emph{RIC}, and (ii) approximate balancedness between $\mb X^*$ and $\mb Y^*$.
\begin{itemize}
\item \emph{Implicit confinement to the RIC. }
In a nutshell, the \emph{region of incoherence and contraction (RIC)} comprises parameters that are close to $[\bo \zeta^*, \mb X^*, \mb Y^*]$ in both global and rowwise norms. 
As it turns out, restricted strong convexity and smoothness continue to hold after augmenting the empirical loss with a balancing penalty. Specifically, define 
\begin{equation}
L^{\sf{aug}}(\bo \theta) = L(\bo \theta) + c_{\sf{aug}}\pi \bignorm{\mb X\t \mb X - \omega \mb Y\t \mb Y}\fb^2.  
\label{eq: augmented function definition}
\end{equation}
with $c_{\sf{aug}} \asymp \underline c_\psi$. The curvature of $L^{\sf{aug}}$, along directions that approximately lie in the tangent space at $[\bo \zeta, \mb X, \mb Y]$, admits upper and lower bounds of nearly the same order as those of its population counterpart, obtained by taking expectation with respect to both $\mb R$ and $\Omega$. Meanwhile, the behavior of the actual iterates hinges on the curvature of $L^{\sf{aug}}$, enabling us to obtain near linear convergence rates (cf.~Remark~\ref{remark:linear-convergence-of-errors}). However, this invariance of the trajectory does not follow from local curvature condition alone; it depends on finer structural properties of the loss and the factorized parametrization. We build upon the implicit-regularization perspective of \citet{ma2018implicit}, \cite{chen2019noisy}, \cite{chen2021bridging} , and \cite{fan2025covariates} and show that, when suitably initialized, vanilla gradient descent remains in the \emph{RIC}, \emph{without} an explicit incoherence or balancing penalty.

We emphasize that our implicit regularization analysis relies crucially on a leave-one-out argument, which decouples the dependence between the iterate $[\bo \zeta^t, \mb X^t, \mb Y^t]$ and a fixed row or column of $\mc P_{\Omega}(\mb R - \psi(\mb M^*))$. Because gradient descent is initialized by the first two stages, the coupled construction
must begin at the spectral estimator and continue through one-step refinement; see Section~\ref{sec: prior algorithmic design}.

\item \emph{Implicit balancedness. }
We also prove that throughout the iterates, the gradient descent algorithm is able to maintain the balancedness conditions appearing in the identifiability condition \eqref{eq: identifiability condition 1}, without explicitly imposing regularization such as $\norm{{\mb X}\t \mb X - \omega \mb Y\t \mb Y }\fb$. 
This form of implicit regularization is noteworthy for two reasons. First, unlike \cite{tu2016low}, \cite{zheng2016convergence}, and \cite{park2018finding}, we do not impose any explicit regularization term to guarantee them. Second, they resonate with the discussion in \cite{soltanolkotabi2023implicit} for matrix sensing, but here it arises in a subtler notion of balancedness that is essential for the downstream inferential procedures. 
\end{itemize}

\subsubsection{Key Ingredients of Linear Approximations}
\label{subsubsec: key ingredients in the proof}
Having established that the gradient descent iterates provably surf through a nice landscape, we now elaborate on how the linear approximation in Theorem~\ref{thm: linear approximations} is derived from the final gradient descent iterate. This derivation constitutes the main technical contribution of our gradient descent analysis.

Let $\hat{\bo \theta}^{t_0}$ denote the concatenation of
$[\bo \zeta^{t_0}, \mb X^{t_0}, \mb Y^{t_0}]$. At a high level, for a \emph{convex} objective, the desired linear approximation typically follows from two ingredients: (i) the gradient at the current iterate is sufficiently small \emph{in the rowwise sense}, and (ii) the estimation error admits a heuristic representation like \eqref{eq: linear approximation heuristic}.

A critical complication, however, is that, leaving rotational ambiguity alone, the actual objective function $L$ in \eqref{eq: loss function}, unlike the augmented one in Section~\ref{subsubsec: implicit regularization mechanism}, is \emph{nonconvex} and fails to be strongly convex in all directions, even locally. This degeneracy stems from the factorized parametrization: the latent factors are invariant under simultaneous rescalings, rendering the Hessian matrix singular. Consequently, controlling the gradient norm rowwisely becomes challenging, for instance, $\norm{\nabla_{\mb X} L(\bo \theta^{t_0})}\ti$, particularly when the iteration time $t_0$ is chosen to be rather small, while the singularity of the Hessian also rules out directly applying the inverse appearing in \eqref{eq: linear approximation heuristic}.

To address these challenges, 
 we combine three complementary technical ingredients in the theoretical analysis, summarized below. 
\begin{itemize}
    \item \emph{Rowwise score decomposition. } Our analysis is guided by the following decomposition. For
$\bo\theta=[\bo\zeta,\mb X,\mb Y]$, write $\nabla_{\mb X_{i,\cdot}} L(\bo \theta) = \mc P_{\Omega}\big( \psi(\mc M(\bo \theta)) - \mb R)_{i,\cdot} \mb Y$. Given $\tilde{\bo \theta}^{t_0} = [\bo \zeta^{t_0}, \tilde{\mb X}^{t_0}, \tilde{\mb Y}^{t_0}]$ with $\tilde{\mb X}^{t_0} = \mb X^{t_0} \mb O^{t_0}$ and $\tilde{\mb Y}^{t_0} = \mb Y^{t_0} \mb O^{t_0}$, we have 
    \begin{align}
        & \underbrace{\nabla_{\mb X_{i,\cdot}} L(\tilde{\bo \theta}^{t_0}) }_{\text{controlled by stationarity}} - \underbrace{\nabla_{\mb X_{i,\cdot}} L(\bo \theta^*)}_{\text{linear form of noise}} \\ 
        = &\underbrace{\nabla_{\mb X_{i,\cdot}} L([\bo \zeta^*, \tilde{\mb X}^{t_0}, \mb Y^*]) - \nabla_{\mb X_{i,\cdot}} L(\bo \theta^*)}_{\approx ( \tilde{\mb X}^{t_0}_{i,\cdot} - \mb X^*_{i,\cdot}) \nabla^2_{\mb X_{i,\cdot}} L(\bo \theta^*)} + \underbrace{\nabla_{\mb X_{i,\cdot}} L(\tilde{\bo \theta}^{t_0})  - \nabla_{\mb X_{i,\cdot}} L([\bo \zeta^*, \tilde{\mb X}^{t_0}, \mb Y^*])}_{\approx \tilde{\mb X}_{i,\cdot}^{t_0} \mb r_{\mb X_i,\mb Y} + \mb r_{\mb X_i, \bo \zeta} \text{ (nuisance term)}}, 
        \label{eq: rowwise analysis heuristic}
    \end{align}
    where we let $\mb D_i \coloneqq \diag( \mc P_\Omega(\psi'(\mc M(\bo \theta^*)))_{i,\cdot})$ and define
    $$\mb r_{\mb X_i,\mb Y} \coloneqq (\tilde{\mb Y}^{t_0} - \mb Y^*)\t \mb D_{i} \mb Y^*, \qquad \mb r_{\mb X_i, \bo \zeta} \coloneqq (\bo \zeta^{t_0} - \bo \zeta^*)\t \mb D_i \mb Y^*. $$ 
    The first approximation follows by treating the corresponding gradient difference as an $r$-dimensional problem evaluated at $\tilde{\mb X}_{i,\cdot}^{t_0}$ and $\mb X_{i,\cdot}^{t_0}$. The second follows from a Taylor expansion of $\nabla_{\mb X_{i,\cdot}}L$ with respect to $\bo\zeta$ and $\mb Y$. Once both $\nabla_{\mb X_{i,\cdot}}L(\bo\theta^{t_0})$ and the nuisance term are shown to be negligible, we obtain
$$
\tilde{\mb X}^{t_0}_{i,\cdot} - \mb X^*_{i,\cdot} = -\nabla_{\mb X_{i,\cdot}} L(\bo \theta^*)\big[\nabla^2_{\mb X_{i,\cdot}} L(\bo \theta^*)\big]^{-1} + \text{negligible error}. 
$$

    \item \emph{Global control of nuisance term. } Controlling the term $\tilde{\mb X}_{i,\cdot}^{t_0} \mb r_{\mb X_i,\mb Y} + \mb r_{\mb X_i, \bo \zeta}$  requires delicate control of $
    \norm{\mb r_{\mb X_i,\mb Y}}\fb$ and $\norm{\mb r_{\mb X_i,\bo \zeta}}_2$. 
    For example, the desired bound for
$\norm{\mb r_{\mb X_i,\mb Y}}\fb$ does not follow directly from the
rowwise error control in Theorem~\ref{thm: GD}. Nor can we directly apply
the argument in \citet{chen2019inference} for the (linear) matrix completion,
which concerns
the unweighted matrix
$(\tilde{\mb Y}^{t_0}-\mb Y^*)\t\mb Y^*$. The additional difficulty stems from the additional random weight matrix $\diag\!\big(\mc P_\Omega(\psi'(\bo\theta^*))_{i,\cdot}\big)$ due to the nonlinearity.
To tackle this, we derive a global approximation for the aggregated estimation error $\bo \Delta \coloneqq [\bo \zeta^{t_0} - \bo \zeta^*, \tilde{\mb X}^{t_0} - \mb X^*, \tilde{\mb Y}^{t_0} - \mb Y^*]$: 
    \begin{align}
        & \mb C \bo \Delta= -\nabla_{\bo \theta} L(\bo \theta^*) + \bo \Psi,
    \end{align}
    where $\mb C$ denotes the effective Fisher-information component that resembles the Hessian of $L$ and $\bo \Psi$ collects the residual effects. Owing to both the nonlinear structure and the nonidentifiability of the factorized parametrization, $\mb C$ has a highly intricate and singular structure. Consequently, this relation does not directly yield a linear approximation for $\tilde{\mb X}^{t_0}_{i,\cdot}-\mb X^*_{i,\cdot}$. Nevertheless, it provides a route to obtaining sharper bounds for $\mb v^\top\bo\Delta$ with $\mb v\t \bo \Delta = (\mb r_{\mb X_i, \mb Y})_{k_1,k_2}$, {as opposed to the row-wise characterization}. Precisely, denote the range of $\mb C$ by 
    \begin{equation}
    \mc V^{\sf{range}} \coloneqq \big(\underbrace{ \mc V^{\sf{scaling}} \oplus \mc V^{\sf{rotation}}}_{\text{the null space of $\mb C$}}\big)^{\perp},
    \end{equation}
    where $\mc V^{\sf{scaling}}$ and $\mc V^{\sf{rotation}}$ correspond to the scaling and rotational ambiguities, respectively. We then decompose
    \begin{align}
        & \mb v\t \bo \Delta =  \mb v\t \mc P_{\mc V^{\sf{range}}}(\bo \Delta) + \mb v\t \mc P_{\mc V^{\sf{scaling}} \oplus \mc V^{\sf{rotation}}}(\bo \Delta)\\ 
        = & - \mb v\t \mb C^{\dagger} \nabla_{\bo \theta}L(\bo \theta^*)  + \mb v\t \mb C^{\dagger} \bo \Psi + \mb v\t \mc P_{\mc V^{\sf{scaling}} \oplus \mc V^{\sf{rotation}}}(\bo \Delta). 
        \label{eq: vtDelta decomposition}
    \end{align}
    The three terms on the right-hand side of \eqref{eq: vtDelta decomposition} are controlled through, respectively: (i) a decoupling argument among $\mb v$, $\mb C^\dagger$, and $\nabla_{\bo\theta}L(\bo\theta^*)$; (ii) perturbation analysis of $\mb C$ and control of the residual $\bo\Psi$; and (iii) the approximate balancedness of $\mb X^{t_0}$ and $\mb Y^{t_0}$.

\item \emph{Auxiliary iterates for rowwise stationarity. }
An auxiliary sequence, initialized at the same point as Algorithm~\ref{alg: GD}, is analyzed under the augmented objective $L^{\sf{aug}}$ in \eqref{eq: augmented function definition}. The proof compares this sequence with the original iterates and combines decay of the auxiliary gradient with the approximate-balancedness bound. This comparison supplies the bound on $\norm{\nabla_{\mb X_{i,\cdot}} L(\tilde{\bo \theta}^{t_0})}_2$ required in \eqref{eq: rowwise analysis heuristic}. 
\end{itemize}

\section{Constructing Good Starts for the Gradient Descent}
\label{sec: prior algorithmic design}
As noted earlier, a fundamental prerequisite for analyzing the gradient-descent trajectory is that the initializer lies within the \emph{region of incoherence and contraction}. For models with a low-rank mean structure, a spectral initializer constructed from the observed matrix $\mc P_\Omega(\mb R)$ is often sufficient \citep{ma2018implicit}, owing to recent advances in entrywise eigenspace and singular-subspace perturbation theory \citep{abbe2020entrywise,agterberg2022entrywise,cai2021subspace}. In the generalized model, however, directly applying a truncated SVD to the observation matrix, whose expectation is $\pi\psi\bigl(\mb 1_n{\bo \zeta^*}\t+\mb X^*{\mb Y^*}\t\bigr)$, typically yields little direct information about the latent factors. The entrywise nonlinear transformation $\psi$ distorts the underlying low-rank structure and, in general, does not preserve the rank of the latent signal matrix. 

Nonetheless, we resolve this mismatch through a two-step SVD procedure, similar to the spirit of \citet{ma2020universal} and \citet{zhang2020note}. For this initial double-SVD stage, we develop a new analysis tailored to the present model and establish Frobenius-norm consistency and leave-one-out stability of the resulting factor estimates. 
The remaining challenge is then to upgrade such global consistency to row-wise consistency, as required for the initializer to enter the \emph{region of incoherence and contraction} introduced in Section~\ref{subsec: mechanism}. To this end, we incorporate a one-step refinement procedure, inspired by \citet{chen2024note}, to sharpen the spectral initializer.

    \subsection{First Step: Double SVD with Link Inversion}
    \label{subsec: USVT}
    The first step of the optimization pipeline is to construct an initializer that is sufficiently close to the true parameters under a suitable metric. Specifically, we first apply an SVD to the zero-filled observation matrix $\mc P_\Omega(\mb R)$, then apply $\psi^{-1}$ entrywise to the resulting estimate (e.g., $\psi^{-1}$ is the logit function in the Bernoulli case), and finally perform a second rank-$r$ truncated SVD, following the general spirit of \cite{ma2020universal,zhang2020note}. To facilitate the leave-one-out analysis, we replace \emph{universal singular value thresholding} \citep{chatterjee2015matrix} with soft singular value thresholding \citep{koltchinskii2011nuclear}. The complete procedure is presented in Algorithm~\ref{alg: USVT}.

\begin{algorithm}[tb]
\caption{Double SVD with link inversion}
\label{alg: USVT}
\KwIn{$\mc P_\Omega(\mb R)\in\bb R^{n\times p}$, threshold $\lambda^{\sp}$, truncation domain $\psi(\mc D_\psi)$, sampling rate $\pi$}
\KwOut{$\hat{\bo \zeta}^{\sp}$, $\hat{\mb X}^{\sp}$, $\hat{\mb Y}^{\sp}$}

Compute the SVD $\mc P_\Omega(\mb R)=\bar{\mb U}\bar{\bo\Sigma}\bar{\mb V}^\top$, where $\bar{\bo\Sigma}=\diag(\bar\sigma_i)$\;

Compute the soft-thresholded approximation
$$
\hat{\mb R}^{\sp}
=
\pi^{-1}\bar{\mb U}\diag\big((\bar\sigma_i-\lambda^{\sp})_+\big)\bar{\mb V}^\top .
$$

Clip $\hat{\mb R}^{\sp}$ entrywise into $\psi(\mc D_\psi)$ and invert:
$$
\hat{\mb M}^{\sp}
=
\psi^{-1}\!\Big(\operatorname{clip}_{\psi(\mc D_\psi)}\big(\hat{\mb R}^{\sp}\big)\Big).
$$

Set
$
\hat{\bo\zeta}^{\sp}= \hat{\mb M}^{\sp\top} \mb 1_n/n$.
Compute the rank-$r$ truncated SVD
$$
\text{Top-$r$ SVD}\big(\hat{\mb M}^{\sp}-\mb 1_n(\hat{\bo\zeta}^{\sp})^\top\big)
= \hat{\mb U}^{\sp} \hat{\bo \Sigma}^{\sp} \hat{\mb V}^{\sp\top};$$
and output $[\hat{\bo \zeta}^{\sp}, \hat{\mb X}^{\sp}, \hat{\mb Y}^{\sp}]$ with
$
\hat{\mb X}^{\sp}= \omega^{\frac14} \hat{\mb U}^{\sp} \hat{\bo \Sigma}^{\sp\frac12}$ and $
\hat{\mb Y}^{\sp}= \omega^{-\frac14} \hat{\mb V}^{\sp} \hat{\bo \Sigma}^{\sp\frac12}$.
\end{algorithm} 

    While the spectrum of $\pi \psi(\mb 1_n \bo \zeta^{*\top} + \mb X^*\mb Y^{*\top})$ is not fully understood, this approach does not require the transformed mean matrix to be rank-$r$. The threshold specified below balances the normalized global estimation error against the leave-one-out discrepancy. Under the stated sampling and signal-strength conditions, it also dominates the base noise scale $\bar\sigma\sqrt{\pi d}$ and the nonlinear approximation error. This feature is particularly useful for nonlinear transformations, which generally inflate the rank of the mean matrix.

     Formally, we have the following theorem, whose proof is deferred to Section~\ref{sec: proof of usvt}. 
     \begin{theorem}[Spectral Initialization]
        \label{thm: USVT}
        Suppose that Assumptions~\ref{assumption: main noise and sampling},
        \ref{assumption: main link function},
        \ref{assumption: main signal strength and incoherence degree}, and
        \ref{assumption: main identifiability} hold. Choose
        \begin{equation*}
            \lambda^{\sp}
            \asymp
            \sigma\sqrt{\bar\sigma\sigma_1^*}\,
            \pi^{\frac34}d^{\frac14}\mu^{\frac14}
            (\log d)^{\frac{3r}{4}},
        \end{equation*}
        with a sufficiently large leading constant.
        Then, with
        probability at least $1-d^{-c}$, there exists a common rotation
        $\mb O^{\sp}\in\mc O(r)$ satisfying
        \begin{align*}
                 &\bignorm{\hat{\bo \zeta}^{\sf{sp}} - \bo \zeta^*}_2 \lesssim \frac{\lambda^{\sp}(\log d)^{\frac r2}}{\underline c_\psi\pi\sqrt{n}}  , \qquad \omega^{-\frac14}\bignorm{\hat{\mb X}^{\sf{sp}}\mb O^{\sp}- \mb X^*}\fb \vee \omega^{\frac14}\bignorm{\hat{\mb Y}^{\sf{sp}}\mb O^{\sp}- \mb Y^*}\fb \lesssim \frac{\lambda^{\sp}(\log d)^{\frac r2}}{\underline c_\psi\pi\sqrt{\sigma_r^*}} .
        \end{align*}
     \end{theorem}

     Here, $\lambda^{\sp}$ governs the tradeoff between the global estimation error and the leave-one-out discrepancy, and the displayed choice balances these two quantities; Theorem~\ref{thm: detailed USVT} in the supplement provides the general characterization. 

     \begin{remark}
        The enlarged scale $\bar\sigma$ accounts for the additional fluctuation
        induced by zero-filling the missing entries during spectral
        initialization. An alternative imputation scheme reduces the
        contribution of $\max_{x\in\mc D_\psi}|\psi(x)|$ to order
        $\sigma^2\mu r\sigma_1^*/\sqrt{np}$; see
        Remark~\ref{remark: alternative imputation scheme} for details.
    \end{remark}

    \subsection{Second Step: Unilateral Refinement}
    \label{subsec: OS}
    
    The \emph{unilateral refinement (UR)} procedure in Algorithm~\ref{alg: OS} first clip the spectral estimates to control their $\ell_\infty$ or $\ell\ti$ norm. It then performs two independent likelihood refinements: one estimates $\mb X$ with $(\bo\zeta^0,\mb Y^0)$ fixed, while the other estimates $(\bo\zeta,\mb Y)$ with $\mb X^0$ fixed. Both subproblems decompose into low-dimensional convex programs and can be solved efficiently using first- or second-order methods.

\begin{algorithm}[tb]
\caption{Unilateral Refinement (UR)}
\label{alg: OS}
\KwIn{Initial estimates $\hat{\bo \zeta}$, $\hat{\mb X}$, $\hat{\mb Y}$; midpoint $\zeta_0$ of $\mc D_\psi$; trimming radii $\iota_{\bo\zeta}$, $\iota_{\mb X}^{(\omega)}$, $\iota_{\mb Y}^{(\omega)}$}
\KwOut{$\hat{\bo \zeta}^{\ur}$, $\hat{\mb X}^{\ur}$, $\hat{\mb Y}^{\ur}$}

Initialize by trimming the input
\begin{equation*}
\bo \zeta^{0}=\zeta_0 \mb 1_p + \mc P^\infty_{\iota_{\bo\zeta}}(\hat{\bo \zeta} - \zeta_0 \mb 1_p),
\qquad
\mb X^{0}=\mc P^{2,\infty}_{\iota_{\mb X}^{(\omega)}}(\hat{\mb X}),
\qquad
\mb Y^{0}=\mc P^{2,\infty}_{\iota_{\mb Y}^{(\omega)}}(\hat{\mb Y}).
\end{equation*}
$\mc P^\infty_a$ denotes coordinatewise clipping, and $\mc P^{2,\infty}_a(\mb A)$ shrinks the rows of $\mb A$ to have $\ell_2$ norm at most $a$. 
The factor-truncation radii are given by $
    \iota_{\mb X}^{(\omega)}=\omega^{1/4}\iota_{\mb X}$, $
    \iota_{\mb Y}^{(\omega)}=\omega^{-1/4}\iota_{\mb Y}$.

Update the estimates by
\begin{align}
& \hat{\mb X}^{\ur}
=
\arg\min_{\mb X}
\sum_{(i,j)\in\Omega}
\Bigl[
-R_{ij}\bigl(\zeta_j^0+\mb X_{i,\cdot}\mb Y_j^0\bigr)
+\Psi\bigl(\zeta_j^0+\mb X_{i,\cdot}\mb Y_j^0\bigr)
\Bigr], \\ 
& (\hat{\bo \zeta}^{\ur},\hat{\mb Y}^{\ur})
=
\arg\min_{\bo \zeta,\mb Y}
\sum_{(i,j)\in\Omega}
\Bigl[
-R_{ij}\bigl(\zeta_j+\mb X^0_{i,\cdot}\mb Y_j\bigr)
+\Psi\bigl(\zeta_j+\mb X^0_{i,\cdot}\mb Y_j\bigr)
\Bigr].
\end{align}
\end{algorithm} 

    The intuition behind this refinement is that, when estimating the $i$th row of the latent factor matrix, the optimizer $\hat{\mb X}_i^{\ur}$ isolates the influence of the $i$th-row noise from that of the remaining rows, thereby enabling sharper row-specific error control. This improvement is further corroborated by the numerical experiments in Section~\ref{sec: numerical simulations}. Formally, we have the following theorem whose proof is presented in Section~\ref{sec: unilateral refinement analysis}. 
    \begin{theorem}[Unilateral Refinement]
        \label{thm: OS}
        Suppose that Assumptions~\ref{assumption: main noise and sampling},
        \ref{assumption: main link function},
        \ref{assumption: main signal strength and incoherence degree}, and
        \ref{assumption: main identifiability} hold. Choose the
        truncation radii so that
            $\iota_{\bo\zeta}
            = c_{\iota} \sqrt{\frac{\mu\sigma_1^*\sigma_r^*}{np}}$,$
            \iota_{\mb X}^{(\omega)}
            = c_{\iota}  \omega^{\frac14}\sqrt{\frac{\mu\sigma_1^*}{n}}$, and $
            \iota_{\mb Y}^{(\omega)}
            = c_{\iota}  \omega^{-\frac14}\sqrt{\frac{\mu\sigma_1^*}{p}}$ for some sufficiently large constant $c_{\iota}$.
        Let $[\bo\zeta^0,\mb X^0,\mb Y^0]$ denote the trimmed input in
        Algorithm~\ref{alg: OS}. Define $\mb O^0 \coloneqq \argmin_{\mb O \in \mc O(r)} \omega^{-\frac12}\bignorm{\mb X^0\mb O-\mb X^*}\fb^2
        +\omega^{\frac12}\bignorm{\mb Y^0\mb O-\mb Y^*}\fb^2$.
        Then, with probability at least $1-O(d^{-c})$, the output of
        Algorithm~\ref{alg: OS}, initialized by Algorithm~\ref{alg: USVT},
        satisfies
        \begin{align}
            \max\Big\{
            &\sqrt{\frac{n}{\sigma_r^*}}
            \bignorm{\hat{\bo\zeta}^{\ur}-\bo\zeta^*}_\infty,
            \omega^{-\frac14}
            \bignorm{\hat{\mb X}^{\ur}\mb O^0-\mb X^*}\ti,
            \omega^{\frac14}
            \bignorm{\hat{\mb Y}^{\ur}\mb O^0-\mb Y^*}\ti
            \Big\}\\
            \lesssim{}&
            \frac{\bar\sigma\sqrt{\log d}}
            {\underline c_\psi\sqrt{\pi\sigma_r^*}}
            +\frac{\sigma\sqrt{\bar\sigma}}
            {\underline c_\psi}
            \Big(\frac{\mu^3}{\pi d}\Big)^{\frac14}
            (\log d)^{\frac{5r}{4}}
            +\frac{\sigma^2\sqrt{\sigma_r^*}
            \mu^{\frac32}\log d}
            {\underline c_\psi\pi d^{\frac32}}.
            \label{eq: plugged-in row-wise error for unilateral refinement}
        \end{align}
    \end{theorem}

    The three terms in the error bound arise, respectively, from the intrinsic
    refinement noise, the balanced contribution of the global and
    leave-one-out errors from the spectral stage, and the nonlinear
    missingness remainder. Thus, the refinement turns
    the spectral initializer's global consistency and stability into uniform
    rowwise control. This, in turn, allows us to identify a rather mild
    condition under which the output of the first two steps enters the
    \emph{region of incoherence and contraction}.
    
    \begin{remark}
    \citet{chen2024note} analyze one-step refinements with and without data splitting. Their bound for the no-splitting case contains an additional $\pi^{-1/2}$ factor, although this procedure performs better in their simulations.
    In our setting, the analogous difficulty is the dependence in terms such as $\mc P_\Omega(\mb E)_{i,\cdot}\hat{\mb Y}^{\sf{sp}}$. We control this dependence by leave-one-out coupling throughout the pipeline, without further sample splitting.
    \end{remark}

    Together with the centering-and-balancing step in Algorithm~\ref{alg: GD}, these guarantees furnish the warm start needed for Theorem~\ref{thm: GD}. Under the additional conditions of Theorem~\ref{thm: linear approximations}, the resulting gradient-descent output admits rowwise linear expansions, which underpin the inference procedures developed in the next section.

\section{Individual and Simultaneous Uncertainty Quantification}
\label{sec: inference}
    The linear expansions in Theorem~\ref{thm: linear approximations} reduce inference to the distribution of rowwise score sums. We shall start from individual inference, as a direct consequence of asymptotic normality. Going beyond, we extend the scenario to the simultaneous inference, leveraging the Gaussian multiplier bootstrap technique. To streamline the notation, we introduce the appended factor estimates $\mb X^{\sf{app},t_0}\coloneqq(\sqrt{\sigma_r^* /n}\mb 1_n,\mb X^{t_0})$ and $\mb Y^{\sf{app},t_0}\coloneqq(\sqrt{n/\sigma_r^*}\bo\zeta^{t_0},\mb Y^{t_0})$ along with their population counterparts $\mb X^{\sf{app},*}\coloneqq (\sqrt{\sigma_r^* / n} \mb 1_n, \mb X^*)$ and $\mb Y^{\sf{app},*} \coloneqq (\sqrt{n/\sigma_r^*} \bo \zeta^*, \mb Y^*)$.

    \subsection{Individual Inference for single entries and single rows } 
    \label{subsec:individual-inference}
    Since the leading score terms in Theorem~\ref{thm: linear approximations} have mean zero, their asymptotic normality follows under mild conditions on the noise, as formalized below.
    \begin{theorem}
    \label{thm: individual inference}
        Suppose that the conditions of Theorem~\ref{thm: linear approximations} hold.
        Define
        $\mb O^{\sf{app},t_0}\coloneqq\diag(1,\mb O^{t_0})$. Then, for every fixed $i\in[n]$ and $j\in[p]$,
        \begin{align*}
            \sup_{\substack{A\subseteq\bb R^r\\ A\ \mathrm{convex}}}
            &\Big|\bb P\Big[(\mb H^*_{\mb X_i})^{\frac12}
            \big(\mb O^{t_0\top}\mb X^{t_0\top}_{i,\cdot}-\mb X^{*\top}_{i,\cdot}\big)\in A\Big]
            -\bb P_{\mb Z\sim\mc N(\mb 0,\mb I_r)}[\mb Z\in A]\Big|=o(1),\\
            \sup_{\substack{A\subseteq\bb R^{r+1}\\ A\ \mathrm{convex}}}
            &\Big|\bb P\Big[(\mb H^*_{\mb Y^{\sf{app}}_{j}})^{\frac12}
            \big(\mb O^{\sf{app},t_0\top}\mb Y^{\sf{app},t_0\top}_{j,\cdot}
            -\mb Y^{\sf{app},*\top}_{j,\cdot}\big)\in A\Big]
            -\bb P_{\mb Z\sim\mc N(\mb 0,\mb I_{r+1})}[\mb Z\in A]\Big|=o(1).
        \end{align*}
        Here, $\mb H^*_{\mb X_i}$ and $\mb H^*_{\mb Y^{\sf{app}}_{j}}$ are defined in \eqref{eq: definition of fisher matrices}.
    \end{theorem}
    A natural way to estimate the covariance matrices that appear in these asymptotic distributions is to plug the gradient descent estimate into the information matrices. Let
$\hat{\mb X}^{t_0}\coloneqq\mb X^{t_0}\hat{\mb O}^{t_0}$ and
$\hat{\mb Y}^{t_0}\coloneqq\mb Y^{t_0}\hat{\mb O}^{t_0}$, where
$\hat{\mb O}^{t_0}$ estimates the population Procrustes rotation $\mb O^{t_0}$. Set
$\hat{\mb M}\coloneqq\mb 1_n\bo\zeta^{t_0\top}+\hat{\mb X}^{t_0}\hat{\mb Y}^{t_0\top}$.
The plug-in information estimators are
\begin{align}
    \hat{\mb H}_{\mb X_i} \coloneqq \sum_{j\in[p], (i,j)\in \Omega} \psi'(\hat M_{i,j}) \hat{\mb Y}_j^{t_0} \hat{\mb Y}_j^{t_0\top} , \qquad \hat{\mb H}_{\mb Y^{\sf{app}}_{j}} \coloneqq \sum_{i\in[n], (i,j)\in \Omega} \psi'(\hat M_{i,j}) \hat{\mb X}^{\sf{app},t_0}_{i} \hat{\mb X}^{\sf{app},t_0\top}_{i}.
    \label{eq:plug-in-cov-estimators}
\end{align}
The resulting plug-in confidence regions satisfy the following coverage guarantee. 
    \begin{theorem}
    \label{thm: individual asymptotic normality}
        Suppose that the conditions of Theorem~\ref{thm: individual inference} hold and, in addition,
        \begin{align*}
            &\sigma\sqrt{\pi\sigma_1^*}\,\delta_{\sf{rotation}}
            \big(\norm{\mb X^*}\ti\vee\norm{\mb Y^*}\ti\big)=o(1),
        \end{align*}
        where $\delta_{\sf{rotation}}\coloneqq
        \bignorm{\hat{\mb O}^{t_0}-\mb O^{t_0}}$.
        For fixed $i\in[n]$ and $j\in[p]$, define
        \begin{align}
         & \mc I_{\bo \zeta_j } \coloneqq \Big\{x \in \bb R: \Big| \Big[\frac{\sigma_r^*}{n}(\hat{\mb H}_{\mb Y^{\sf{app}}_{j}}^{-1})_{1,1}\Big]^{-\frac12} (x - \zeta^{t_0}_{j}) \Big| \leq z_{1-\frac{\alpha}{2}}\Big\},
        \\
        & \mc I_{\mb X_i} \coloneqq \Big\{\mb x \in \bb R^r: \Big\| \hat{\mb H}_{{\mb X}_i}^{\frac12}\big(\mb x - \hat{\mb X}^{t_0\top}_{i,\cdot}\big)\Big\|_2^2 \leq q_{\chi,r,1-\alpha} \Big\}, \\
        & \mc I_{\mb Y_j} \coloneqq \Big\{\mb y\in \bb R^r: \Big\|\big[(\hat{\mb H}_{{\mb Y^{\sf{app}}_{j}}}^{-1})_{2:r+1,2:r+1} \big]^{-\frac12} \big(\mb y - \hat{\mb Y}^{t_0\top}_{j,\cdot}\big)\Big\|_2^2 \leq q_{\chi,r,1-\alpha} \Big\}.
        \end{align}
        Here, $z_{1-\alpha/2}$ is the $(1-\alpha/2)$-quantile of the standard normal distribution, and $q_{\chi,r,1-\alpha}$ is the $(1-\alpha)$-quantile of the $\chi^2$ distribution with $r$ degrees of freedom. Then
        \begin{align}
            &\sup_{0 < \alpha < 1} \Big| \bb P\Big[\zeta_{j}^* \in \mc I_{\bo \zeta_j} \Big] - (1 - \alpha) \Big| = o(1), \qquad
            \sup_{0 < \alpha < 1} \Big| \bb P\Big[\mb X_{i,\cdot}^{*\top} \in \mc I_{\mb X_{i}} \Big] - (1 - \alpha) \Big| = o(1),\\
            &
            \sup_{0 < \alpha < 1} \Big| \bb P\Big[\mb Y_{j,\cdot}^{*\top} \in \mc I_{\mb Y_j} \Big] - (1 - \alpha) \Big| = o(1).
        \end{align}
    \end{theorem}

    \begin{remark}
    To assess the optimality of our inference procedure, we observe that the inverses of the normalizing terms match the Cram\'er–Rao lower bounds for estimating an individual entry or row of $\bo \zeta^*$, $\mb X^*$, or $\mb Y^*$, under the oracle setting in which all remaining parameters are known. Formally, suppose that we have access to the intercept $\bo \zeta^*$ as well as the right factor $\mb Y^*$. Then for any unbiased estimator $\hat{\mb x}$ for $\mb X_{i}^*$, we have 
    \begin{equation}
         \mathrm{Cov}(\hat{\mb x}\mid \Omega) \succeq \mb H^{*-1}_{\mb X_i} = \mathrm{CRLB}(\mb X_{i}^*\mid \Omega). 
    \end{equation}
    This lower bound coincides with the limiting covariance in Theorem~\ref{thm: individual inference}, thereby establishing oracle efficiency. The same argument applies to $\bo \zeta^*$ and $\mb Y^*$. 
        \end{remark}

    \subsection{Simultaneous Confidence Bands via Gaussian Multiplier Bootstrap}
    \label{subsec:bootstrap}
The preceding results quantify uncertainty for a fixed number of parameters. Many applications, however, require simultaneous comparisons over a collection whose size may grow with the dimensions. Fortunately, our new uniform rowwise linear approximations results make it possible to approximate the distributions of the resulting maximum statistics by a Gaussian multiplier bootstrap.

We begin with a unified framework covering a broad class of maxima over linear functionals of the latent factors. Specifically, consider targets of the form
\begin{equation}
\begin{aligned}
\mathcal T_{\mb X}
&\coloneqq
\max_{l\in[L]}
\Big\{
\sum_{(i,k)\in\mc S_{\mb X,l}}
\hat a_{\mb X,l,i,k}
\big(
\hat{\mb X}^{t_0}-\mb X^*
\big)_{i,k}
\Big\},
\qquad
\mc S_{\mb X,l}\subseteq[n]\times[r],
\nonumber\\
\mathcal T_{\mb Y}
&\coloneqq
\max_{l\in[L]}
\Big\{
\sum_{(j,k)\in\mc S_{\mb Y,l}}
\hat a_{\mb Y,l,j,k}
\big(
\widehat{\mb Y}^{\sf{app},t_0}-\mb Y^{\sf{app},*}
\big)_{j,k}
\Big\},
\qquad
\mc S_{\mb Y,l}\subseteq[p]\times[r+1].
\end{aligned}
\label{eq:factor-focused-statistics}
\end{equation}
Here the index sets $\mc S_{\mb X,l}$ and $\mc S_{\mb Y,l}$, together with a sequence of
estimated weights $\hat a_{\mb X,l,i,k}$ and $\hat a_{\mb Y,l,j,k}$, specify the particular contrasts of
interest, including simultaneous confidence bands, and
pairwise comparisons.
We emphasize that the number $L$ of contrasts may diverge, while the sizes of $\mc S_{\mb X,l}$ and $\mc S_{\mb Y,l}$ are fixed, as needed in most applications and as assumed in our theory; for example, in the case of simultaneous inference for all of the $n$ entries in the column $\mb X^*_{\cdot,1}$, the index set $\mc S_{\mb X, l}$ would be a singleton $\{(l,1)\}$ for $l\in[n]$, and $\hat a_{\mb X,l,l,1} = 1$. 
 
To approximate the distributions of the maximum statistics introduced above, we construct bootstrap analogues of their linear expansions using the Gaussian multiplier bootstrap. 
\begin{enumerate}
    \item For each $b\in[N_{\sf boot}]$, generate independent standard Gaussian multipliers
    $
    \{\xi_{\mb X,i}^{[b]}: i\in[n]\}$ and $
    \{\xi_{\mb Y,j}^{[b]}: j\in[p]\}$.

    \item Let $\hat{\mb E}\coloneqq\mc P_\Omega\{\mb R-\psi(\hat{\mb M})\}$. For each bootstrap replication $b\in[N_{\sf boot}]$, construct the multiplier-score rows
    \begin{align}
    (\mb S_{\mb X}^{[b]})_{i,\cdot}
    &\coloneqq \hat{\mb E}_{i,\cdot}
    \diag\big(\xi_{\mb Y,j}^{[b]}\big)_{j\in[p]}
    \hat{\mb Y}^{t_0},\quad i\in[n];\quad 
    (\mb S_{\mb Y^{\sf{app}}}^{[b]})_{j,\cdot}\coloneqq \hat{\mb E}_{\cdot,j}\t
    \diag\big(\xi_{\mb X,i}^{[b]}\big)_{i\in[n]}
    \hat{\mb X}^{\sf{app},t_0},\quad j\in[p].
    \end{align}
    Based on these components, define the linearized bootstrap quantities
    \begin{equation}
    \begin{aligned}
        (\bo \Delta^{[b]}_{\mb X})_{i,\cdot}
        &\coloneqq(\mb S_{\mb X}^{[b]})_{i,\cdot}\hat{\mb H}_{\mb X_i}^{-1},\qquad
        (\bo \Delta^{[b]}_{\mb Y^{\sf{app}}})_{j,\cdot}
        \coloneqq(\mb S_{\mb Y^{\sf{app}}}^{[b]})_{j,\cdot}\hat{\mb H}_{\mb Y^{\sf{app}}_{j}}^{-1}.
        \end{aligned}
        \label{eq: bootstrap Delta}
    \end{equation}

    \item Approximate the distributions of $\mathcal T_{\mb X}$ and $\mc T_{\mb Y}$ by the empirical distributions, over $b\in[N_{\sf boot}]$, of
    \begin{align} 
    & \mathcal T_{\mb X}^{[b]}
\coloneqq
\max_{l\in[L]}
\Big\{
\sum_{(i,k)\in\mc S_{\mb X,l}}
\hat a_{\mb X,l,i,k} (\bo \Delta_{\mb X}^{[b]})_{i,k}
\Big\},\qquad 
\mathcal T_{\mb Y}^{[b]}
\coloneqq
\max_{l\in[L]}
\Big\{
\sum_{(j,k)\in\mc S_{\mb Y,l}}
\hat a_{\mb Y,l,j,k}
(\bo \Delta_{\mb Y^{\sf{app}}}^{[b]})_{j,k}
\Big\}, 
\label{eq:factor-focused-bootstrap-statistics}
    \end{align}

\end{enumerate} 

The validity of this procedure follows from combining our rowwise linear approximations with high-dimensional central limit theorems \citep{chernozhukov2013gaussian,chernozhukov2022improved,chernozhukov2023nearly}. Its score-perturbation construction is related to existing bootstrap approaches for simultaneous inference \citep{belloni2018uniformly,zhang2017simultaneous,chernozhukov2013gaussian}. 
To our knowledge, our result is among the first to justify simultaneous inference for the intercepts and both latent factor matrices in generalized latent factor models with missing observations. A central difficulty is that both factor matrices are unknown: the score for $\mb X$, for example, depends on the estimated $\mb Y$. Establishing bootstrap validity therefore requires controlling how estimation errors in the two factors propagate through the score statistics.

\begin{theorem}
\label{thm: simultaneous inference}
    Suppose that the conditions of Theorem~\ref{thm: linear approximations} hold and that
    $\max_{l\in[L]}\{|\mc S_{\mb X,l}|\vee|\mc S_{\mb Y,l}|\}$ is upper bounded by a constant. Let $a_{\mb X,l,i,k}$ and $a_{\mb Y,l,j,k}$ denote the oracle counterparts
of the estimated weights. 
    Let $\bar a_{\mb X}$ and $\underline a_{\mb X}>0$ denote the maximum and
    minimum, over $l\in[L]$, of
    $\{\sum_{(i,k)\in\mc S_{\mb X,l}}a_{\mb X,l,i,k}^2\}^{1/2}$,
    and define $\bar a_{\mb Y}$ and $\underline a_{\mb Y}>0$ analogously.
    Assume that
    $\bar a_{\mb X}/\underline a_{\mb X}\vee
    \bar a_{\mb Y}/\underline a_{\mb Y}=O(1)$, and set
    \begin{equation*}
        \delta_a\coloneqq\max_{l\in[L], (i,k) \in \mc S_{\mb X,l}}\big\{|\hat a_{\mb X,l,i,k}-a_{\mb X,l,i,k}|\big\} \vee \max_{l\in[L], (j,k) \in \mc S_{\mb Y,l}}\big\{ 
        |\hat a_{\mb Y,l,j,k}-a_{\mb Y,l,j,k}|\big\}.
    \end{equation*}
    Let $\delta_{\sf{rotation}}\coloneqq\bignorm{\hat{\mb O}^{t_0}-\mb O^{t_0}}$.
    In addition to $\log L\lesssim\log d$, assume that 
    \begin{align}
        &\sqrt{\mu}(\log d)^2\delta_{\sf{rotation}}
        \Big\{1\vee\frac{\sigma\sqrt{\pi}\sigma_1^*}{\sqrt{n\wedge p}}\Big\}=o(1),\qquad
        \delta_a\sqrt{\mu}\log d
        =o(\underline a_{\mb X}\wedge\underline a_{\mb Y}). 
        \label{eq: conditions for simultaneous inference}
    \end{align}
    The oracle weights $\{a_{\mb X, l,i,k}\}$ and $\{a_{\mb Y, l, j,k}\}$
    may depend on $\bo\Omega$ but not on $\mb E$.
    Let $\mc G_{\mb X,1-\alpha}$ and $\mc G_{\mb Y,1-\alpha}$ be the exact $(1-\alpha)$ conditional quantiles of a fresh multiplier draw from \eqref{eq:factor-focused-bootstrap-statistics}, given the data. Then, for every fixed $\alpha\in(0,1)$,
    \begin{align}
        & \big| \bb P\big[ \mc T_{\mb X}\leq \mc G_{\mb X,1-\alpha} \big] - (1-\alpha) \big| = o(1),\qquad
        \big| \bb P\big[ \mc T_{\mb Y}\leq \mc G_{\mb Y,1-\alpha} \big] - (1-\alpha) \big| = o(1).
    \end{align}
    If $\hat{\mc G}_{\mb X,1-\alpha}$ and $\hat{\mc G}_{\mb Y,1-\alpha}$ are the empirical quantiles based on $N_{\sf boot}$ independent draws and $N_{\sf boot}\to\infty$, the same conclusions hold with $\mc G$ replaced by $\hat{\mc G}$; the additional conditional distribution error is $O_{\bb P}(N_{\sf boot}^{-1/2})$.

\end{theorem}
In words, provided that the rotation and weight estimation errors are sufficiently small, the multiplier bootstrap yields asymptotically valid simultaneous inference for a broad class of maximum statistics over $L$ functionals, where \(L\) may grow polynomially with \(d\). The requisite control of the rotation error may follow from Proposition~\ref{prop:rotation-recovery} of the Supplementary Material.

We next consider two applications of simultaneous inference.
\paragraph{Ranking intervals for factor entries.}
Factor ranking is central to LLM evaluation, with public leaderboards ordering models
by aggregate benchmark scores or pairwise human preferences
\citep{chiang2024chatbot,alzahrani2024benchmarks}. Such rankings can be
sensitive to minor evaluation choices \citep{alzahrani2024benchmarks}. Our
latent factors instead yield capability-specific rankings of models and
questions; since each rank is determined by pairwise differences,
Theorem~\ref{thm: simultaneous inference} yields simultaneous confidence
intervals for these ranks.

Specifically, we fix $k\in[r]$ and define the descending rank by $\operatorname{rank}_k(i)\coloneqq
1+|\{l\in[n]:X_{l,k}^*>X_{i,k}^*\}|$. Similar to \citet{fan2025ranking}, we can
construct simultaneous rank intervals from pairwise comparisons. For
$i\neq j$, let $\sigma_{X_{i,k}}^2\coloneqq
(\mb H_{\mb X_i}^{*-1})_{k,k}$ and $\hat\sigma_{X_{i,k}}^2\coloneqq
(\hat{\mb H}_{\mb X_i}^{-1})_{k,k}$, and set
$s_{i,j,k}^2\coloneqq\sigma_{X_{i,k}}^2+\sigma_{X_{j,k}}^2$ and
$\hat s_{i,j,k}^2\coloneqq\hat\sigma_{X_{i,k}}^2+
\hat\sigma_{X_{j,k}}^2$. Conditional independence of the score rows makes
$s_{i,j,k}^2$ the conditional variance of the leading pairwise score
difference. Writing $\Delta_{i,k}^{[b]}$ for the $(i,k)$th entry of
$\bo\Delta_{\mb X}^{[b]}$, define
\begin{align}
    \mc T_k^{\sf{rank}}
    &\coloneqq\max_{i\neq j}\hat s_{i,j,k}^{-1}
    \big|\hat X_{i,k}^{t_0}-\hat X_{j,k}^{t_0}-X_{i,k}^*+X_{j,k}^*\big|,
    \label{eq: ranking stat}\\
    \mc T_k^{\sf{rank},[b]}
    &\coloneqq\max_{i\neq j}\hat s_{i,j,k}^{-1}
    \big|\Delta_{i,k}^{[b]}-\Delta_{j,k}^{[b]}\big|.
    \label{eq: bootstrap ranking stat}
\end{align}
Let $\hat q_{1-\alpha}^{\sf{rank}}$ be the empirical $(1-\alpha)$-quantile of
$\mc T_k^{\sf{rank},[b]}$. We then define
\begin{equation}
    \begin{aligned}
    \mc I_i&\coloneqq[L_i,U_i]\cap[n],\\
    L_i&\coloneqq1+\big|\{l\in[n]:\hat X_{l,k}^{t_0}-\hat X_{i,k}^{t_0}>\hat q_{1-\alpha}^{\sf{rank}}\hat s_{l,i,k}\}\big|,\\
    U_i&\coloneqq n-\big|\{l\in[n]:\hat X_{l,k}^{t_0}-\hat X_{i,k}^{t_0}<-\hat q_{1-\alpha}^{\sf{rank}}\hat s_{l,i,k}\}\big|.
    \end{aligned}
    \label{eq: ranking CI construction}
\end{equation}

The following result is a direct consequence of
Theorem~\ref{thm: simultaneous inference}; $\underline a_{\sf{rank}}$ below
is proportional to the minimum Euclidean norm of the pairwise contrast
weights.
\begin{proposition}[Simultaneous rank coverage]
\label{thm:left-factor-ranking-consistency}
Suppose that the conditions of Theorem~\ref{thm: linear approximations} hold.
Set
$\underline a_{\sf{rank}}\coloneqq\min_{i\neq j}s_{i,j,k}^{-1}$,
$\bar a_{\sf{rank}}\coloneqq\max_{i\neq j}s_{i,j,k}^{-1}$, and
$\delta_{a,{\sf{rank}}}\coloneqq\max_{i\neq j}
|\hat s_{i,j,k}^{-1}-s_{i,j,k}^{-1}|$,
and assume $\underline a_{\sf{rank}}>0$,
$\bar a_{\sf{rank}}/\underline a_{\sf{rank}}=O(1)$, and
$N_{\sf boot}\to\infty$. With
$\delta_{\sf{rotation}}\coloneqq
\bignorm{\hat{\mb O}^{t_0}-\mb O^{t_0}}$, assume further that
\begin{align*}
    &\sqrt{\mu}(\log d)^2\delta_{\sf{rotation}}
    \Big\{1\vee\frac{\sigma\sqrt{\pi}\sigma_1^*}{\sqrt{n\wedge p}}\Big\}=o(1),\qquad
    \delta_{a,{\sf{rank}}}\sqrt{\mu}\log d
    =o(\underline a_{\sf{rank}}).
\end{align*}
Then, for every fixed $\alpha\in(0,1)$,
\begin{equation*}
    \bb P\big[\operatorname{rank}_k(i)\in\mc I_i
    \text{ for every }i\in[n]\big]\geq1-\alpha-o(1).
\end{equation*}
\end{proposition}
The same construction applies to the rows of $\mb Y^*$, yielding simultaneous
intervals for the capability-specific ranks of the questions.

\paragraph{Simultaneous missing-entry prediction.}
We finally consider simultaneous prediction of the unobserved means. This
task is invariant under a common rotation of the estimated factors and hence
requires no rotation estimator; equivalently, one may take
$\hat{\mb O}^{t_0}=\mb I_r$ in the definitions below. For
$(i,j)\in\Omega^{\complement}$ and each bootstrap draw, define
\begin{align*}
    & A_{i,j}^{[b]}\coloneqq(\bo\Delta_{\mb X}^{[b]})_{i,\cdot}
    \hat{\mb Y}_j^{t_0}, \quad
    C_{i,j}^{[b]}\coloneqq(\bo\Delta_{\mb Y^{\sf{app}}}^{[b]})_{j,\cdot}
    \hat{\mb X}^{\sf{app},t_0}_{i},\quad \\
    & \hat s_{\mb X,i,j}^2\coloneqq(\hat{\mb Y}_j^{t_0})^\top
    \hat{\mb H}_{\mb X_i}^{-1}\hat{\mb Y}_j^{t_0}, \quad
    \hat s_{\mb Y^{\sf{app}},i,j}^2\coloneqq
    (\hat{\mb X}^{\sf{app},t_0}_{i})^\top
    \hat{\mb H}_{\mb Y^{\sf{app}}_{j}}^{-1}\hat{\mb X}^{\sf{app},t_0}_{i}.
\end{align*}
Consider the bootstrap maximum statistics
\begin{align*}
    \mc T_{\mb X}^{\sf{ent},[b]}\coloneqq
    \max_{(i,j)\in\Omega^{\complement}}
    \frac{|A_{i,j}^{[b]}|}{\hat s_{\mb X,i,j}}, \quad
    \mc T_{\mb Y^{\sf{app}}}^{\sf{ent},[b]}\coloneqq
    \max_{(i,j)\in\Omega^{\complement}}
    \frac{|C_{i,j}^{[b]}|}{\hat s_{\mb Y^{\sf{app}},i,j}}.
\end{align*}
Let $\hat q^{\sf{ent}}_{\mb X,1-\frac{\alpha}{2}}$ and
$\hat q^{\sf{ent}}_{\mb Y^{\sf{app}},1-\frac{\alpha}{2}}$ denote the empirical
$(1-\frac{\alpha}{2})$-quantiles of
$\mc T_{\mb X}^{\sf{ent},[b]}$ and
$\mc T_{\mb Y^{\sf{app}}}^{\sf{ent},[b]}$, respectively. A ratio with a zero
numerator and denominator is understood to be zero.
The resulting simultaneous band is
\begin{equation}
\begin{aligned}
    \mc I_{i,j}^{\sf ent}\coloneqq\Big[&
    \psi\big(\hat M_{i,j}-\hat q^{\sf{ent}}_{\mb X,1-\frac{\alpha}{2}}
    \hat s_{\mb X,i,j}-\hat q^{\sf{ent}}_{\mb Y^{\sf{app}},1-\frac{\alpha}{2}}
    \hat s_{\mb Y^{\sf{app}},i,j}\big),\\
    &\psi\big(\hat M_{i,j}+\hat q^{\sf{ent}}_{\mb X,1-\frac{\alpha}{2}}
    \hat s_{\mb X,i,j}+\hat q^{\sf{ent}}_{\mb Y^{\sf{app}},1-\frac{\alpha}{2}}
    \hat s_{\mb Y^{\sf{app}},i,j}\big)\Big],\qquad
    (i,j)\in\Omega^{\complement}.
\end{aligned}
    \label{eq: prediction interval}
\end{equation}
The following result establishes its simultaneous coverage.
\begin{proposition}[Simultaneous coverage for missing-entry means]
\label{prop:missing-entry-prediction}
Suppose that the conditions of Theorem~\ref{thm: linear approximations} hold
and $N_{\sf boot}\to\infty$. Then, for every fixed $\alpha\in(0,1)$,
\begin{equation*}
    \bb P\Big[\psi(M_{i,j}^*)\in\mc I_{i,j}^{\sf ent}
    \text{ for all }(i,j)\in\Omega^{\complement}\Big]
    \geq1-\alpha-o(1).
\end{equation*}
\end{proposition}

 \section{Related Literature}
 
    The literature on latent factor models is vast; see \citet{bartholomew2008analysis} for an overview. Low-rank estimation under linear observation models is studied by, among others, \citet{candes2012exact,chen2019noisy,chen2020nonconvex,chen2019inference}. Nonlinear observation models and one-bit matrix completion are considered in \citet{davenport20141,cai2013max,chen2024note}.

Inference for linear low-rank models, including incomplete panels, has been developed in \citet{chen2019inference,xie2023efficient,chernozhukov2023inference,xiong2023large,choi2024inference,su2025inference}. Related likelihood-based theory for nonlinear panels, comparison, and networks is studied in \citet{wang2022maximum,chen2021nonlinear,gao2023binary,li2023statistical,chen2023statistical,ouyang2024statistical,tian2026bridging,li2026low,gao2022community}. By contrast, we characterize the actual output of a computationally tractable, warm-started nonconvex pipeline under severe missingness and derive uniform rowwise linear approximations that support simultaneous inference for the intercepts, latent factors, and missing-entry means. 

    Our optimization dynamics analysis is related to the implicit-regularization literature for linear latent factor models \citep{chen2020nonconvex,chen2019inference,ma2018implicit} and nonlinear low-rank estimation \citep{cui2026convexity}. Beyond rowwise confinement, we quantify the implicit balancedness of the factor iterates, a phenomenon related to \citet{du2018algorithmic,ma2021beyond,soltanolkotabi2023implicit}, with refined bounds tailored to downstream statistical inference.

\section{Simulation Studies}
\label{sec: numerical simulations}

This section evaluates the empirical performance of the proposed computational pipeline and statistical inference methods.
We first examine the estimation errors of the final gradient descent estimator (Algorithm~\ref{alg: GD}) under varying signal strength and sampling rates, in comparison with the preceding stages (Algorithms~\ref{alg: USVT}~and~\ref{alg: OS}). 
Then, we conduct experiments to verify the inferential guarantees provided in Section~\ref{sec: inference}.

Consider the logistic latent factor model
$$
R_{ij} \mid M^*_{ij} \sim \mathrm{Bernoulli}\big( 1 / (1 + \exp(-M^*_{ij})) \big),
\qquad
M^*_{i,j} = \zeta^*_j + \mb X_i^{*\top} \mb Y_j^*,
$$
where each entry is observed independently with probability $\pi$.  The left latent factor is generated from $ \sqrt{\lambda n}\cdot \bar{\mb U} \mb E_U \cdot \diag(\sqrt{3/2},1)$, where $\bar{\mb U} \in \mc O(n,n-1)$ is an arbitrary orthonormal matrix with $\mb 1_n\t \bar{\mb U} = \mb 0$ and  $\mb E_U$ is sampled from $\mathsf{Unif}(\mathrm{Stiefel}(r,n-1))$, 
while the right latent factor follows $ \sqrt{\lambda p}\cdot \mb E_V \cdot \diag(\sqrt{3/2},1)$ with $\mb E_V \sim \mathsf{Unif}(\mathrm{Stiefel}(r,p))$. The resulting matrix $\mb X^*{\mb Y^*}\t$ is then rescaled through the top-$r$ singular value decomposition so that the identifiability constraints with $\omega=1$ are satisfied. We sample each entry of $\bo\zeta^*$ independently from $\sf{Unif}([-0.1,0.1])$ and use $\pi=0.5$ as the baseline sampling rate. The factor learning rate is $\eta=0.5/(\pi\sigma_1^*)$. In the numerical implementation, we omit the input trimming in Algorithm~\ref{alg: OS} and apply no clipping during unilateral refinement. Throughout, we align the two factor estimates at each stage to $(\mb X^*,\mb Y^*)$ using their common orthogonal Procrustes rotation, defined as in Theorem~\ref{thm: GD}; the rotation for the final iterate is denoted by $\mb O^{t_0}$. We find that the algorithm is relatively insensitive to the choice of hyperparameters; the remaining implementation details are reported in Section~\ref{sec: experimental implementation details} of the Supplementary Material.

\subsection{Experiments with Various Signal Strength and Missing Levels}

\begin{figure}[tb]
    \centering
        \includegraphics[width=\linewidth]{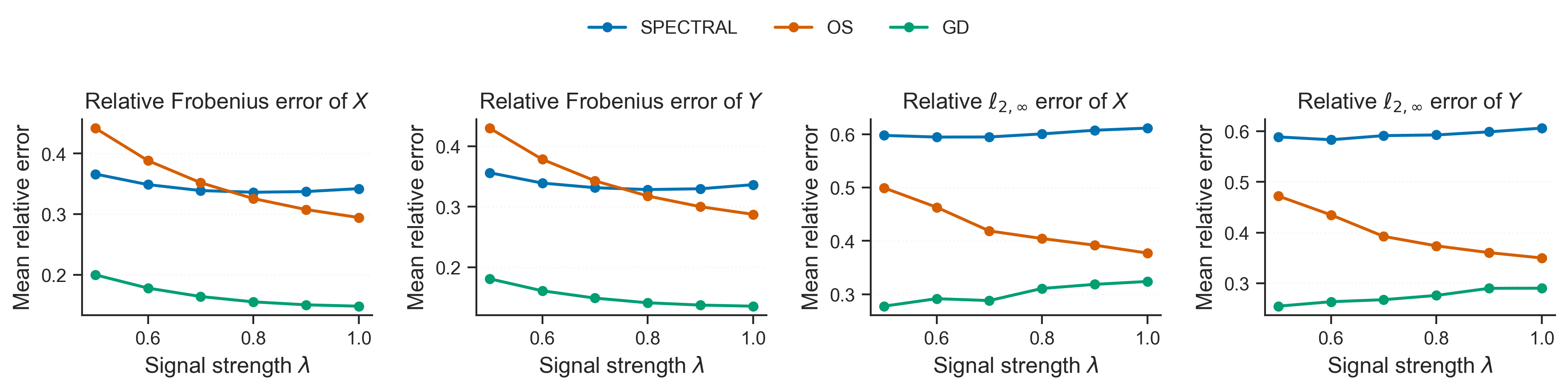}

    \caption{Estimation errors across algorithmic stages under varying signal strength $\lambda$. Here, the relative Frobenius error is defined as $\norm{\mb X^{t_0}\mb O^{t_0}-\mb X^*}\fb/\norm{\mb X^*}\fb$, whereas the relative $\ell\ti$ error is $\norm{\mb X^{t_0}\mb O^{t_0}-\mb X^*}\ti/\norm{\mb X^*}\ti$. The corresponding errors for $\mb Y^{t_0}$ are defined analogously.}
    \label{fig:rho_rho_sweep}
    
\end{figure}

\begin{figure}[tb]
        \centering
        \includegraphics[width=\linewidth]{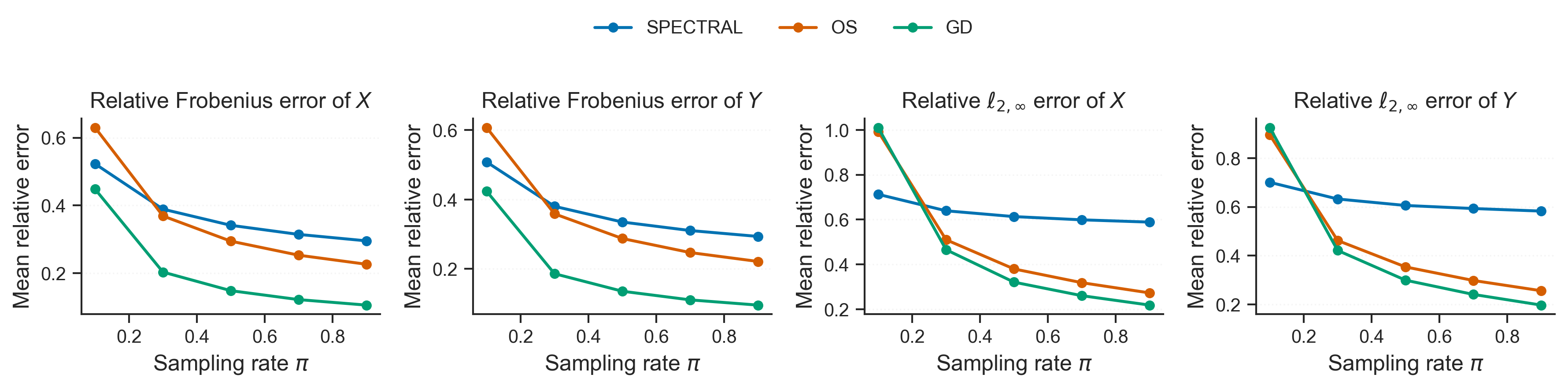}
    \caption{Estimation errors across algorithmic stages under varying sampling rate $\pi$; the errors are calculated using the same convention as Figure~\ref{fig:rho_rho_sweep}. }
    \label{fig:rho_pi_sweep}
\end{figure}

We begin with the full estimator. The goal of this experiment is to identify the practical stability region of the pipeline after all three stages have been applied. We set $n=1000$, $p=500$, and $r=2$, and vary the signal strength $\lambda \in \{0.5,0.6,\cdots, 1.0\}$ and the sampling rate $\pi \in \{0.2,0.4,\cdots, 1.0\}$. Figures~\ref{fig:rho_rho_sweep}~and~\ref{fig:rho_pi_sweep} plot the mean relative errors against the signal strength $\lambda$ and the sampling rate $\pi$, with separate panels for $\pi$ and separate curves for $\lambda$. 


Across 200 Monte Carlo replications, the relative Frobenius errors generally decrease as either the signal strength $\lambda$ or the observation probability $\pi$ increases. The spectral estimator provides an informative initialization, the one-step refinement substantially improves rowwise accuracy, and the final gradient-descent stage typically produces a further reduction in estimation error. These gains become particularly pronounced as $\pi$ increases, indicating that denser observations improve the stability of both model- and task-side factor recovery.\footnote{Given weak signals, the one-step refinement can occasionally increase the relative Frobenius error compared with the spectral initializer. Nevertheless, its improvement in rowwise error mostly persists in these regimes, consistent with the theoretical motivation for this refinement step.} Overall, the results demonstrate the benefits of the proposed multistage procedure and show that its performance improves systematically under stronger signals and more complete observation.

\subsection{Uncertainty Quantification}
\begin{figure}[tb]
    \centering

    \begin{subfigure}[t]{0.242\textwidth}
        \centering
        \includegraphics[width=\linewidth]
        {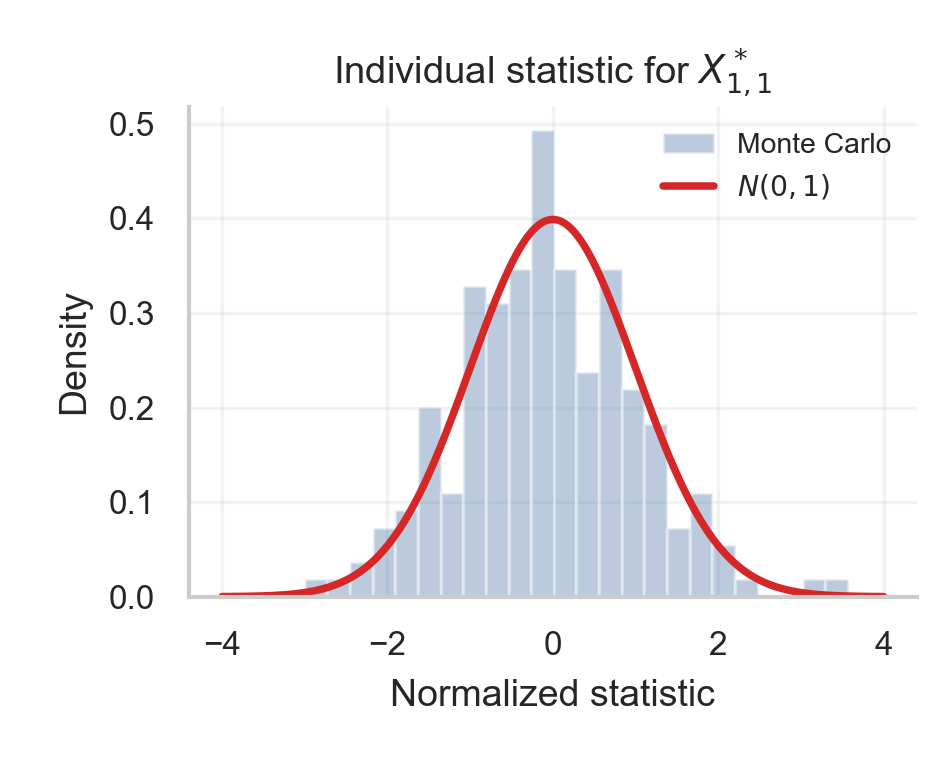}
        \caption{Individual studentized statistic.}
        \label{fig:x11-histogram}
    \end{subfigure}
    \hfill
    \begin{subfigure}[t]{0.242\textwidth}
        \centering
        \includegraphics[width=\linewidth]
        {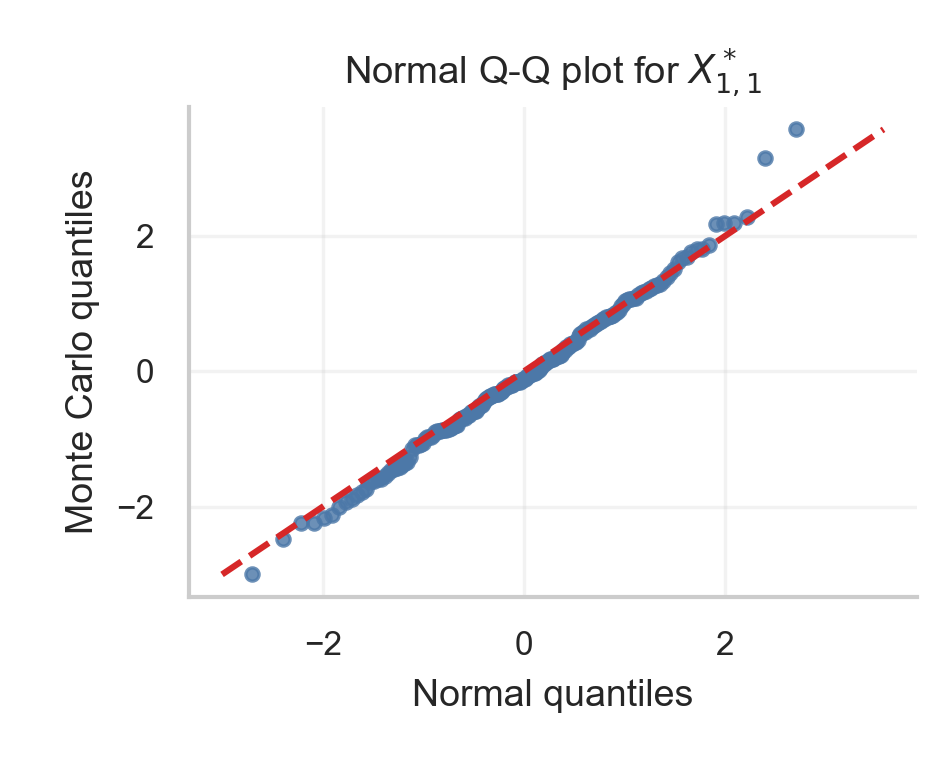}
        \caption{Normal Q--Q diagnostic.}
        \label{fig:x11-normal-qq}
    \end{subfigure}
    \hfill
    \begin{subfigure}[t]{0.242\textwidth}
        \centering
        \includegraphics[width=\linewidth]
        {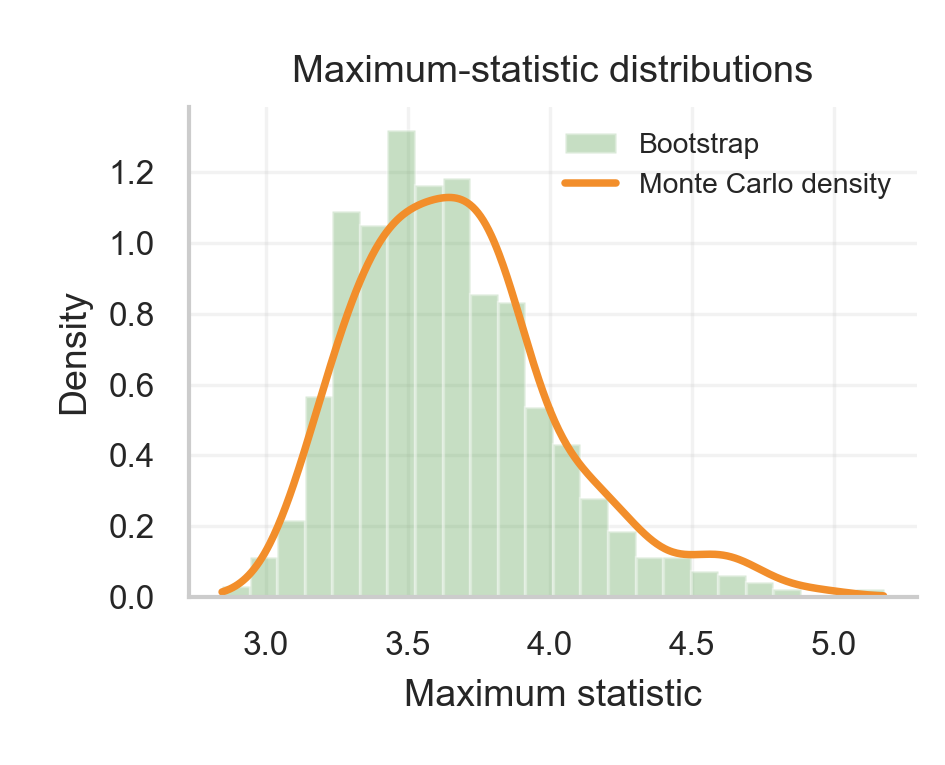}
        \caption{Maximum-statistic distributions.}
        \label{fig:x1-max-histogram}
    \end{subfigure}
    \hfill
    \begin{subfigure}[t]{0.242\textwidth}
        \centering
        \includegraphics[width=\linewidth]
        {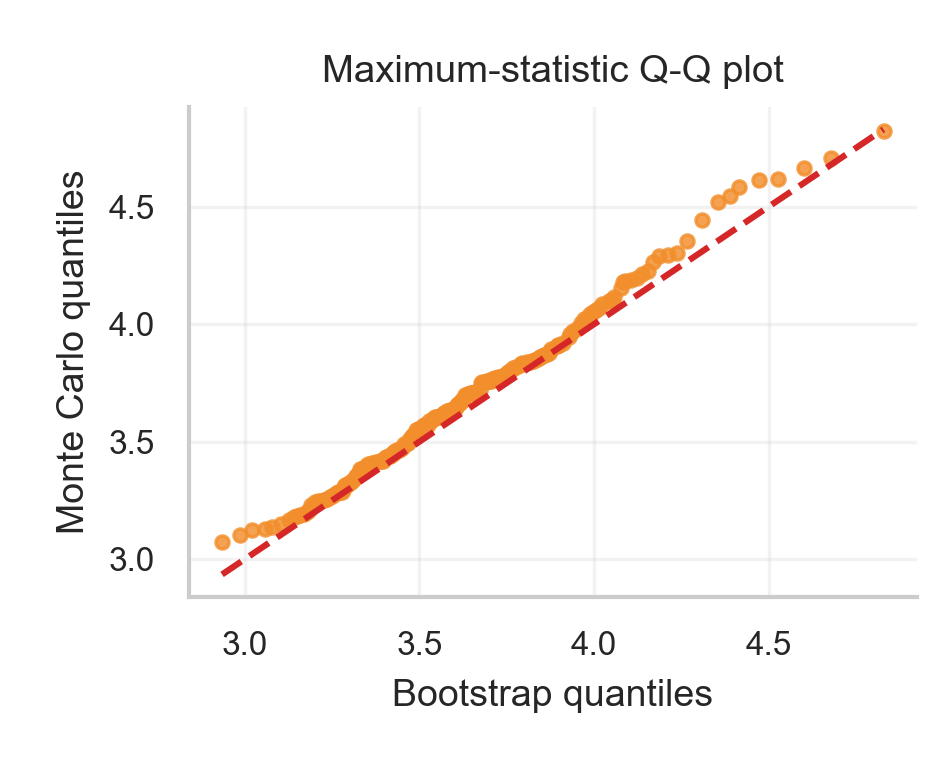}
        \caption{Bootstrap Q--Q diagnostic.}
        \label{fig:x1-max-qq}
    \end{subfigure}

\caption[Distributional diagnostics for inference on the first model factor.]{
Distributional diagnostics for inference on the first model factor, based on
$n=2000$, $p=1500$, $r=2$, $\lambda=1$, and $\pi=0.5$. The responses are regenerated over 200
Monte Carlo replications; the multiplier bootstrap uses 1000 draws from one
response realization. \textbf{(a)} Monte Carlo distribution of the
studentized error for $X^*_{1,1}$ with the standard normal density.
\textbf{(b)} Corresponding normal Q--Q plot. \textbf{(c)} Bootstrap histogram
and Monte Carlo density of the maximum studentized error.
\textbf{(d)} Q--Q comparison of the bootstrap and Monte Carlo maximum
statistics. 
}
    \label{fig:x1-inference-diagnostics}
\end{figure}

We first conduct numerical experiments to verify the asymptotic normality of the entries as well as the consistency of the bootstrap procedures. We start with the asymptotic normality of $(\mb X^*)_{1,1}$, as predicted by Theorem~\ref{thm: individual asymptotic normality}. We then turn to the distributions of the maximum statistic $\max_{i\in[n]} \big|(\hat{\mb X}^{t_0}-\mb X^*)_{i,1}\big|/\sqrt{(\mb H^{*-1}_{\mb X_i})_{1,1}}$ and its bootstrap counterpart $\max_{i\in[n]} \big|(\bo\Delta^{[b]}_{\mb X})_{i,1}\big|/\sqrt{(\hat{\mb H}^{-1}_{\mb X_i})_{1,1}}$. Here, $\hat{\mb X}^{t_0}\coloneqq\mb X^{t_0}\mb O^{t_0}$, $\bo\Delta_{\mb X}^{[b]}$ is defined in \eqref{eq: bootstrap Delta}, and $\mb H_{\mb X_i}^*$ and $\hat{\mb H}_{\mb X_i}$ are defined in \eqref{eq: definition of fisher matrices} and \eqref{eq:plug-in-cov-estimators}, respectively. The corresponding histograms and Q--Q plots are presented in Figure~\ref{fig:x1-inference-diagnostics}.

Figures~\ref{fig:x11-histogram}~and~\ref{fig:x11-normal-qq} show that the
studentized error for $X^*_{1,1}$ is well approximated by a standard normal
distribution. Figures~\ref{fig:x1-max-histogram}
and~\ref{fig:x1-max-qq} further show that the Gaussian multiplier bootstrap
accurately approximates the maximum statistic used for simultaneous
inference.

Next, we consider two simultaneous confidence-band construction tasks based on
the Gaussian multiplier bootstrap procedures developed in
Section~\ref{sec: inference}: confidence intervals for the rankings
constructed from the pairwise maximum statistic (cf.~\eqref{eq: ranking stat},
\eqref{eq: bootstrap ranking stat}, and~\eqref{eq: ranking CI construction})
and studentized prediction intervals for the missing entries
(cf.~\eqref{eq: prediction interval}).
Table~\ref{tab:nested-bootstrap-coverage} reports their empirical coverage
rates and average interval widths at the nominal $95\%$ level. The rank intervals attain $100\%$ simultaneous coverage in all settings, although their widths reflect some conservativeness. This is a familiar cost of simultaneously inverting a large collection of pairwise comparisons \citep{bhattacharya2026inferences,fan2025ranking}. 
The prediction bands exceed the nominal level except in the sparsest setting,
$\pi=0.1$.

\begin{table}[tb]
    \centering
    \caption{Empirical coverage of nominal 95\% multiplier-bootstrap
    confidence bands under one-at-a-time variation of signal strength and
    observation probability. Each setting uses $n=2000$, $p=1500$, rank
    $r=2$, and 200 Monte Carlo replications. We use
    $B=500$ multiplier draws. }
    \label{tab:nested-bootstrap-coverage}
    \small
    \begin{tabular}{@{}lcccc@{}}
        \toprule
        Setting
        & Rank coverage (\%)
        & Mean rank width
        & Entry coverage (\%)
        & Mean entry CI width \\
        \midrule
        \multicolumn{5}{@{}l}{\textit{Signal-strength sweep
        ($\pi=0.5$)}} \\
        $\lambda=0.5$ & 100.0 & 910.1 & 99.5 & 2.206 \\
        $\lambda=0.6$ & 100.0 & 801.6 & 98.5 & 2.315 \\
        $\lambda=0.7$ & 100.0 & 723.4 & 100.0 & 2.439 \\
        $\lambda=0.8$ & 100.0 & 666.6 & 98.5 & 2.578 \\
        $\lambda=0.9$ & 100.0 & 622.2 & 100.0 & 2.725 \\
        $\lambda=1.0$ & 100.0 & 589.9 & 99.5 & 2.886 \\
        \addlinespace
        \multicolumn{5}{@{}l}{\textit{Observation-probability sweep
        ($\lambda=1$)}} \\
        $\pi=0.1$ & 100.0 & 1443.9 & 79.5 & 7.513 \\
        $\pi=0.3$ & 100.0 & 767.9 & 99.5 & 3.815 \\
        $\pi=0.5$ & 100.0 & 589.9 & 99.5 & 2.886 \\
        $\pi=0.7$ & 100.0 & 497.5 & 99.5 & 2.415 \\
        $\pi=0.9$ & 100.0 & 437.7 & 100.0 & 2.113 \\
        \bottomrule
    \end{tabular}
    \par\smallskip
    \begin{minipage}{0.96\linewidth}
        \footnotesize
        \textit{Notes:} Rank coverage requires the intervals for all
        entries of $\mb X^*_{\cdot,1}$ to contain their descending ranks
        simultaneously. Mean rank width averages the number of integer
        ranks in these intervals over all entries and replications.
        Entry coverage is simultaneous over all entries in
        $\Omega^{\complement}$. The reported entry-CI width is measured
        on the linear-predictor scale and averages the full width
        $
        2\left\{
        \hat q^{\sf ent}_{\mb X,0.975}\hat s_{\mb X,i,j}
        +
        \hat q^{\sf ent}_{\mb Y^{\sf app},0.975}
        \hat s_{\mb Y^{\sf app},i,j}
        \right\}$.
    \end{minipage}
\end{table}

\section{Real Data Application to LLM Evaluation}
\label{sec: real data application}
We analyze the binary model--question response matrix constructed from the primary split of the \emph{Metabench} dataset \citep{kipnis2025metabench}. The rows of the data matrix are indexed by the LLMs and the columns are indexed by the questions, while the binary value indicates if the corresponding model answers the question correctly. We fit the data by the Bernoulli generalized latent factor model \eqref{eq: bernoulli model example}. 
\emph{Metabench} aggregates item-wise evaluation records from the Open LLM Leaderboard and selects a compact set of diagnostic questions across ARC, GSM8K, HellaSwag, MMLU, TruthfulQA, and WinoGrande. 
For the present analysis we use the 858 primary Metabench items, consisting of 145 ARC questions, 237 GSM8K questions, 93 HellaSwag questions, 96 MMLU questions, 154 TruthfulQA questions, and 133 WinoGrande questions. The processed response matrix contains 5,055 models with complete responses on all 858 items. We set the factor number to $r=3$, under which the estimated structure remains stable and interpretable across random data splits, and apply the Varimax rotation described in Section~\ref{subsec: gd theory}. 

\begin{figure}[ht]
    \centering
    \begin{subfigure}[t]{0.33\linewidth}
        \centering
        \includegraphics[width=\linewidth]
            {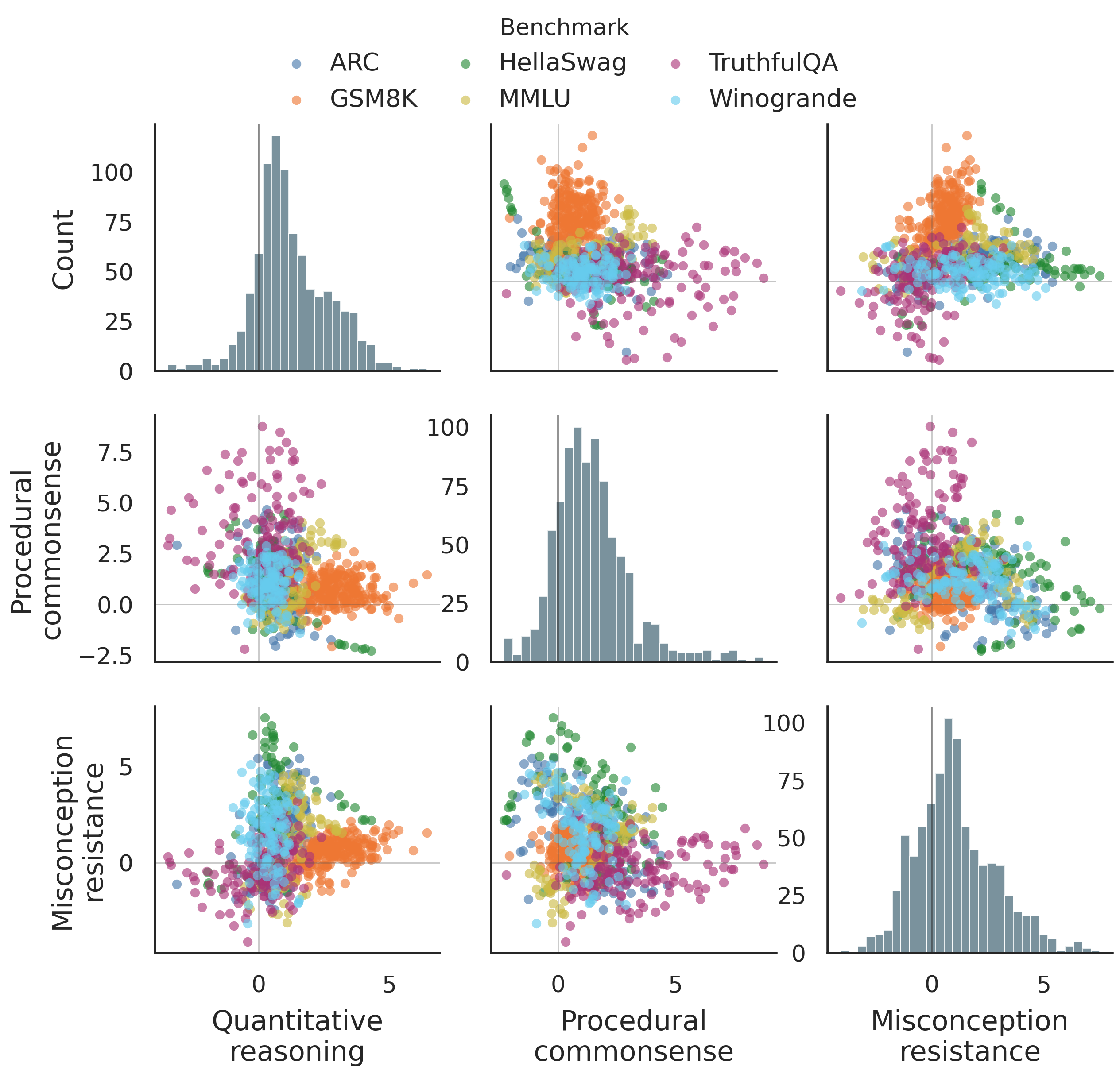}
        \caption{Pairwise relationships among the rotated question factors.}
        \label{fig:pairplot}
    \end{subfigure}
    \hfill
    \begin{subfigure}[t]{0.65\linewidth}
        \centering
        \includegraphics[width=\linewidth]
            {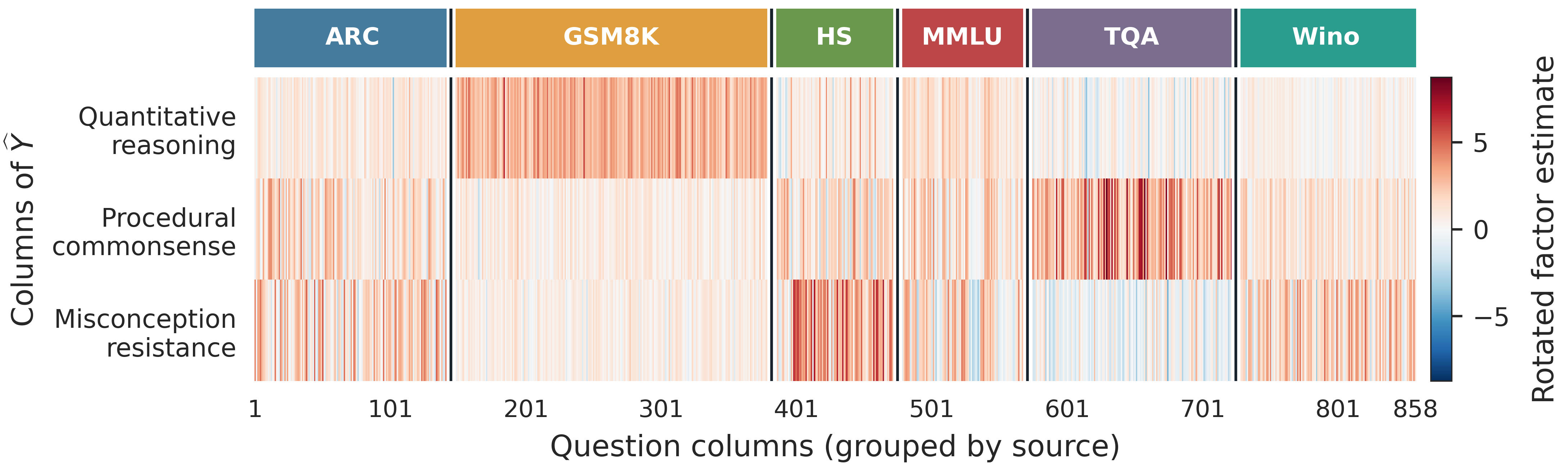}
        \caption{Rotated question-factor estimates grouped by evaluation source.}
        \label{fig:heatmap}
    \end{subfigure}
    \caption{Structure of the Varimax-rotated question-factor estimates.}
    \label{fig:metabench-y-factors}
\end{figure}

We begin by examining the estimated right factors, which characterize variation across questions. Figure~\ref{fig:pairplot} displays their pairwise relationships, while Figure~\ref{fig:heatmap} reveals a pronounced block-sparse structure across benchmark sources. Questions from \emph{GSM8K} load most strongly on the first factor; because these questions emphasize multistep arithmetic and mathematical problem solving, we call this factor \emph{quantitative reasoning}. The second factor is driven primarily by \emph{HellaSwag} and \emph{WinoGrande}, which require models to apply commonsense knowledge to identify plausible events and resolve contextual relationships, motivating the label \emph{procedural commonsense}. Finally, \emph{TruthfulQA} loads most strongly on the third factor; because this benchmark tests whether models resist common misconceptions and produce truthful responses, we call this factor \emph{misconception resistance}. Together, these loading patterns provide a coherent substantive interpretation of the estimated latent structure. 

\begin{figure}[htbp]

\centering

    \includegraphics[width=0.95\linewidth]{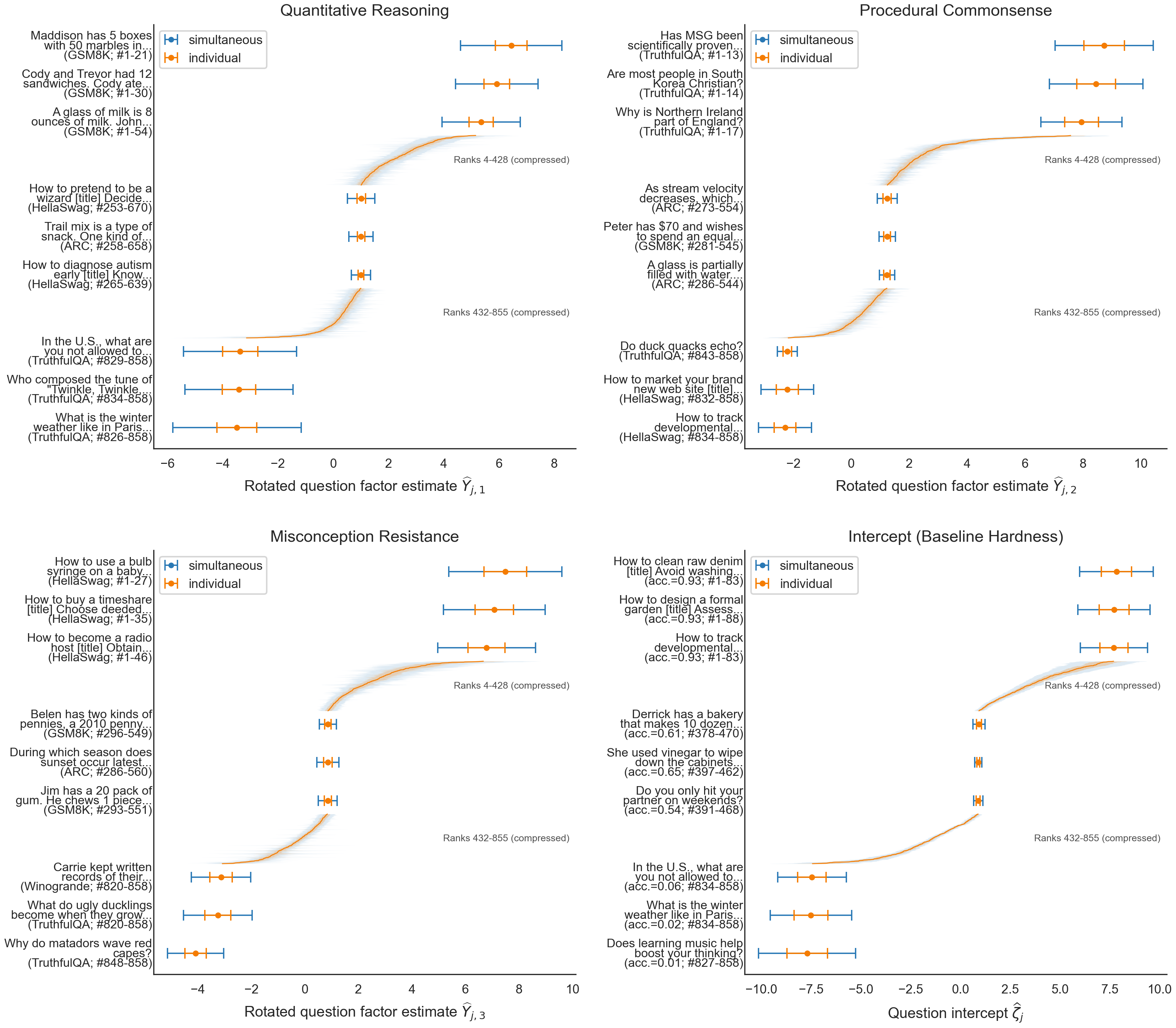}

\caption{Ranking uncertainty for the rotated question factors and question intercepts. Each panel highlights the three highest-ranked, three middle-ranked, and three lowest-ranked questions, with blue and orange segments denoting simultaneous and pointwise confidence intervals, respectively. Questions between these groups are displayed in compressed rank order: blue and orange shaded bands represent their simultaneous and pointwise confidence intervals, while the orange curve traces their point estimates. Parenthetical values report simultaneous rank confidence intervals; labels in the intercept panel additionally report each question’s average accuracy.}
\label{fig:question_ci}
\end{figure}

\begin{figure}[htbp]

\centering

\includegraphics[width=0.95\linewidth]{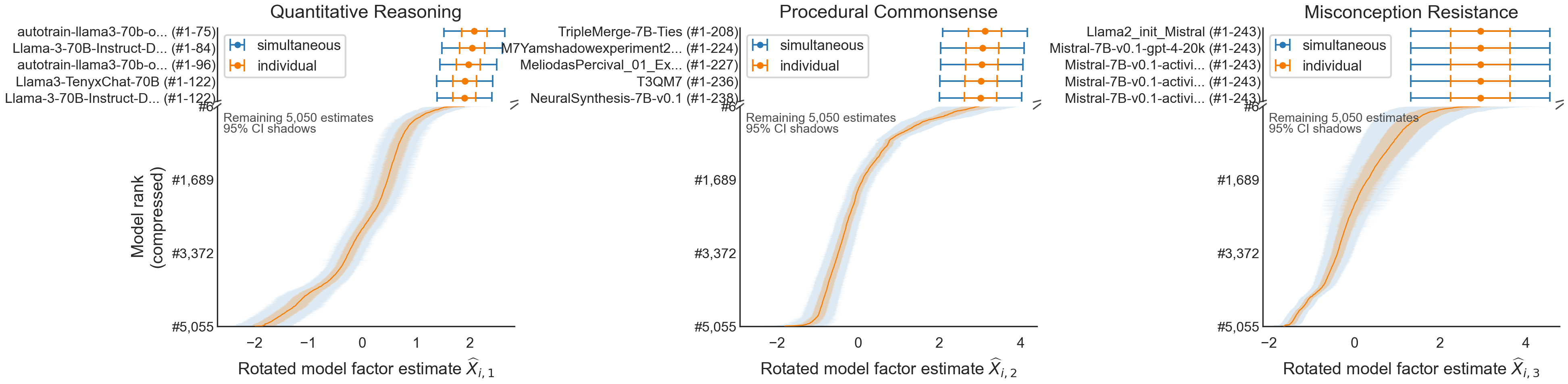}
\caption{Ranking uncertainty for the rotated model factors. Point estimates and confidence intervals for the factor values and model rankings follow the conventions of Figure~\ref{fig:question_ci}; only the five highest-ranked models under each factor are shown. }
\label{fig:model_ci}
\end{figure}

Figure~\ref{fig:question_ci} reports pointwise and simultaneous confidence intervals, together with simultaneous rank confidence intervals, for the question-level estimates, following the procedures developed in Sections~\ref{subsec:individual-inference} and~\ref{subsec:bootstrap}. The benchmark composition of the highest-ranked questions reinforces the factor interpretations: questions from \emph{GSM8K}, \emph{HellaSwag}, and \emph{TruthfulQA} are prominent under the quantitative-reasoning, procedural-commonsense, and misconception-resistance factors, respectively. The accompanying rank intervals quantify the uncertainty in these orderings rather than treating the estimated rankings as exact. In addition, the question-intercept estimate $\bo \zeta^{t_0}$ closely tracks the corresponding average accuracies, supporting their interpretation as baseline easiness, or equivalently as negative difficulty under the present sign convention.
Turning to the model-level results, Figure~\ref{fig:model_ci} displays confidence intervals and simultaneous rank intervals for the five models with the largest estimated scores along each latent capability. Although the point estimates suggest an ordering, every displayed simultaneous rank interval includes rank one, and many span dozens or even hundreds of positions. The data therefore do not support a decisive ranking among the leading models, cautioning against treating small differences in aggregate benchmark performance as conclusive evidence of model superiority.

\section{Conclusion and Discussion}
\label{sec: discussion}
We have introduced and analyzed a three-step algorithm for uncertainty quantification in generalized latent factor models with missing data. First, we obtain a Frobenius-consistent initialization via a double-SVD procedure. Next, we refine these estimates through a unilateral refinement procedure to achieve rowwise consistency. Finally, we apply vanilla gradient descent, benefiting from its implicit regularization, to locate an approximate stationary point of the nonconvex likelihood.

Our theoretical and empirical results confirm that this workflow delivers both accurate point estimates and valid confidence bands. The linear approximation of the gradient-descent iterates underpins our inferential framework, allowing us to characterize their limiting distributions and to construct individual and simultaneous confidence sets via the Gaussian multiplier bootstrap. Moreover, by analyzing the loss landscape within the \emph{region of incoherence and contraction}, we show why gradient descent remains robust to local minima and permits optimal estimation errors. The optimality manifests in three complementary senses. First, the aggregate error matches the restricted-class minimax lower bound in Theorem~\ref{thm: minimax lower bound} up to logarithmic factors and the limiting variance coincides with the corresponding Cram\'er-Rao lower bound.  Second, for balanced logistic Bernoulli models, where the scale and sparsity of the average success probability is of interest, the theory requires only $\pi n\psi(\zeta_0)\gg(\log n)^{5r+4}$; under full observation, it therefore permits success probabilities of order $n^{-1}(\log n)^c$ for any fixed $c>5r+4$. Third, the sampling rate $\pi$ attains the information-theoretic lower bound up to logarithmic factors.

Our simulations demonstrate that the proposed method achieves high estimation precision and reliable coverage. 
This work thus provides a practical and theoretically grounded toolkit for inference in high-dimensional latent factor settings with missingness. Its application to the \emph{Metabench} dataset validates its effectiveness in extracting hidden factors from LLM responses and in uncertainty quantification of their latent capabilities. 

Several paths are left open for further work: (i) such nonlinearity-encoded factor structure is also natural for symmetric network data as well as tensor data, thus it would be interesting to see whether we can have computationally tractable and statistically optimal procedures for such models; (ii) While our analysis reaches the fixed-rank degrees-of-freedom sampling scale up to logarithmic factors in the balanced regime $(n \approx p)$, the singular subspace perturbation theory for the linear models \citep{cai2018rate,zhang2022heteroskedastic,cai2021subspace} suggests that, in the unbalanced regime ($n \gg p$ or $n \ll p$),
the error rates for the left and right factors should scale according to their respective dimensions. Characterizing the sharp dependence on these dimensions in the present nonlinear setting is an important question for future study.

\bibliographystyle{apalike}
\bibliography{ref}
\newpage
\begin{appendices}
\setcounter{tocdepth}{3}
\tableofcontents{}
\clearpage

\section{Preliminary Setups}
\label{sec: supplementary preliminary setups}
\section*{Roadmap of the Supplementary Material}

The theoretical analysis in this Supplementary Material follows the computational pipeline, from
initialization to statistical inference.
Section~\ref{sec: supplementary preliminary setups} records the
parameter-explicit assumptions and notation, allowing the parameters
treated as constants in the main text to vary.

Section~\ref{sec: proof of usvt} analyzes the double-SVD initializer.
We first establish global estimation accuracy by tracking the errors
through singular value thresholding, entrywise link inversion, and
the second low-rank factorization. We then construct leave-one-row-out
and leave-one-column-out counterparts and bound their discrepancies
from the original estimates. These stability bounds provide the
decoupling needed for the subsequent refinement analysis.

Section~\ref{sec: unilateral refinement analysis} analyzes the two
unilateral likelihood refinements, each holding the preliminary
estimates on the opposite side fixed. Leave-one-out decoupling and
control of the sampled Hessians yield rowwise error bounds despite
the dependence between the preliminary estimates and the observations.
We further establish leave-one-out stability of the refined estimates.
Together, these results supply the accuracy and stability requirements
for the warm-started gradient-descent analysis, without additional
sample splitting.

Section~\ref{sec: supplementary gradient descent} presents the
parameter-explicit gradient-descent and linear-approximation guarantees,
whose proofs proceed through three dependent steps:
\begin{enumerate}
    \item \emph{Control of the gradient-descent trajectory.}
    We establish the local geometry and use leave-one-out induction
    to control estimation errors, incoherence, centering, and
    balancedness under a general warm start. Verifying that the
    first two stages supply the required initialization then yields
    the gradient-descent guarantee for the full pipeline.

    \item \emph{From trajectory control to linear approximation.}
    An auxiliary sequence under the augmented objective provides
    rowwise approximate stationarity. We then control the weighted
    nuisance terms through a global error decomposition that
    accounts for the nonidentifiable directions. These ingredients,
    combined with a rowwise score expansion, yield the linear
    representations in
    Section~\ref{subsec: proof of inference results}.

    \item \emph{From linear approximation to inference.}
    Section~\ref{subsec:proof-of-inference-thm} combines these
    representations with Gaussian and multiplier-bootstrap
    approximations to establish the individual and simultaneous
    inferential guarantees and their downstream applications.
    Section~\ref{subsec: rotation-free factor estimation} studies
    data-driven recovery of the rotations needed for interpretable
    factor inference.
\end{enumerate}

Section~\ref{sec: proof of minimax lower bound} establishes the minimax
lower bound. The remaining sections collect auxiliary results on
sampling, concentration, matrix perturbation, and distributional
approximation, followed by the numerical implementation details in
Section~\ref{sec: experimental implementation details}. 

\subsection{Parameter-Explicit Assumptions}
For completeness, we record the parameter-explicit assumptions used throughout
the supplementary analysis. These conditions allow $r$, $\kappa$, and
$\kappa_\psi$ to vary, and reduce to Assumptions~\ref{assumption: main noise and sampling}--\ref{assumption: main identifiability}
under the fixed-parameter regime adopted in the main text.

Let $\mb U^*\bo\Sigma^*\mb V^{*\top}$ denote the rank-$r$ SVD of
$\mb X^*\mb Y^{*\top}$, with singular values
$\sigma_1^*\geq\cdots\geq\sigma_r^*$, and define
$\kappa\coloneqq\sigma_1^*/\sigma_r^*$. We say that
$[\bo\zeta^*,\mb X^*,\mb Y^*]$ is $\mu$-incoherent if
$\frac nr\norm{\mb U^*}\ti^2\vee\frac pr\norm{\mb V^*}\ti^2\leq\mu$ and
$\norm{\bo\zeta^*-\zeta_0\mb 1_p}_\infty
\leq\sqrt{\mu r\sigma_r^*/(np)}$ for some reference level $\zeta_0$.

\begin{assumption}[Sampling scheme and noise assumptions]
\label{assumption: noise and sampling}
\begin{itemize}
    \item[(a)] Each entry is independently observed with probability
    $\pi\in(0,1]$.
    \item[(b)] The observation follows the formulation in
    \eqref{eq: likelihood function}. For some constant $c>1$, the entries
    $R_{i,j}$ are mutually independent and satisfy
    $\max_{i,j}\mathrm{Var}(R_{i,j})^{1/2}\vee d^{-c-10}\leq\sigma$ and
    $|\log\sigma|\lesssim\log d$. Moreover, with probability at least
    $1-O(d^{-c-10})$, $\max_{i,j}|R_{i,j}|\leq B$, where, without loss of
    generality, $B\geq\max_{x\in\mc D_\psi}|\psi(x)|$.
\end{itemize}
\end{assumption}

\begin{assumption}[Link function conditions]
\label{assumption: link function}
{\emergencystretch=1.5em
Let $\psi\coloneqq\Psi'$. Set
$\mc D_\psi\coloneqq
[\zeta_0-c_{\mc D_\psi}\mu r\sigma_1^*/\sqrt{np},\allowbreak
\zeta_0+c_{\mc D_\psi}\mu r\sigma_1^*/\sqrt{np}]$
for a sufficiently large fixed constant $c_{\mc D_\psi}\geq1$ as the domains of both $\Psi$ and
$\psi$. A possibly diverging or vanishing $\underline c_\psi$ and a
constant-order smoothness parameter $L_\psi$ satisfy, uniformly over
$x\in\mc D_\psi$,\par}

\begin{equation*}
    \underline c_\psi\leq\psi'(x)\leq\sigma^2,
    \qquad
    |\psi^{(k)}(x)|\leq\sigma^2k!L_\psi^{k-1},
    \quad k\in\bb N^+.
\end{equation*}
Define $\kappa_\psi\coloneqq\sigma^2/\underline c_\psi$.
\end{assumption}

\begin{assumption}[Signal conditions]
\label{assumption: signal strength and incoherence degree}
Assume the following conditions hold for $\mb X^*$ and $\mb Y^*$:
\begin{itemize}
    \item[(a)] \emph{Signal strength:} Define
    $\bar\sigma\coloneqq\sigma\vee\max_{x\in\mc D_\psi}|\psi(x)|$.
    Assume $\sigma_r^*\gtrsim1$ and
    \begin{align}
        \sigma_r^*
        &\gg
        \frac{\bar\sigma\kappa^{11}\kappa_\psi^5
        \mu^{\frac52}r^{\frac92}\rho^{\frac32}
        \xi^5(\log d)^3}{\underline c_\psi}
        \sqrt{\frac d\pi},
        \qquad
        \pi(n\wedge p)
        \gg\rho\mu^3r^2\kappa^3\kappa_\psi^4\xi^8,
        \label{eq: general consolidated sampling condition}
    \end{align}
    where
    $\xi\coloneqq
    [c_\xi r\log d\,(L_\psi\mu r^2\sigma_1^*/\sqrt{np}\vee1)]^{r/2}$
    for some sufficiently large constant $c_\xi$.
    \item[(c)] \emph{Tail condition:}
    The relation between $\sigma$ and $B$ is given by
    $\kappa^4\kappa_\psi^4\mu r^2B\sqrt{\log d}
    \ll\sigma\sqrt{\pi(n\wedge p)}$.
\end{itemize}
\end{assumption}

\begin{assumption}[Basic Identifiability]
\label{assumption: identifiability}
The factors $\mb X^*,\mb Y^*$ and the intercept vector $\bo\zeta^*$ obey
\begin{align}
    \mb 1_n^\top\mb X^*&=\mb 0_r^\top,\label{eq: general identifiability condition 1}\\
    {\mb X^*}^\top\mb X^*&=\omega{\mb Y^*}^\top\mb Y^*.
    \label{eq: general identifiability condition 2}
\end{align}
\end{assumption}

\paragraph{Notations.} 
We collectively state (or restate) the notations used in the following context. 
We write $\diag(\mb A_1,\cdots,\mb A_N)$ for the (block) diagonal matrix $\left( \begin{matrix}
    \mb A_1 & & \\ 
    & \ddots & \\ 
    & & \mb A_N
\end{matrix}\right)$ where each $\mb A_i$ is either a scalar or a square matrix. And we let $\sf{vec}(\mb A)$ for $\mb A\in \bb R^{n_1 \times n_2}$ be the flattening vector $(\mb A_{1,\cdot}, \cdot, \mb A_{n_1,\cdot})\t \in \bb R^{n_1n_2}$. We denote by the boldface matrix $\bo \Omega\in \bb R^{n\times p}$ the binary matrix realization of $\Omega$, and recall that $\mc P_{\Omega}(\mb R) = \mb R \circ \bo \Omega$. For an arbitrary matrix $\mb A\in \bb R^{n_1\times n_2}$, we denote by $\sigma_{m}(\mb A)$ the $m$-th largest singular value of $\mb A$. 
We define $[-a]$ as $\{-a ,-a+1,\cdots,-1\}$ for $a \in \bb N^+$. 

Let $\mb R$ be the sample matrix $n$-by-$p$ with all entries following 
$$
p(\mb R \mid [\bo \zeta^*, \mb X^*, \mb Y^*]) = \prod_{(i,j)\in[n]\times [p]} \exp\big[ \gamma(R_{i,j}) + R_{i,j} M^*_{i,j} - \Psi(M^*_{i,j}) \big]. 
$$
We write $\mb R^* \coloneqq \bb E[\mb R] = \psi(\mb M^*) = \psi(\mb 1_n {\bo \zeta^*}\t + \mb X^*{\mb Y^*}\t)$ and distinguish the response noise from the additional fluctuation induced by missingness:
\begin{equation}
    \mb E \coloneqq \mb R-\mb R^*,
    \qquad
    \bar{\mb E} \coloneqq \mc P_\Omega(\mb R)-\pi\mb R^*
    =\mc P_\Omega(\mb E)+(\bo\Omega-\pi\mb 1_n\mb 1_p\t)\circ\mb R^*.
    \label{eq: response and spectral noise matrices}
\end{equation}
The entries of both matrices are centered. In particular, independence of $\bo\Omega$ and $\mb R$ gives
\begin{equation}
    \bb E\big[(\bar{\mb E}_{i,j})^2\big]
    =\pi\mathrm{Var}(R_{i,j})+\pi(1-\pi)(R_{i,j}^*)^2
    \leq 2\pi\bar\sigma^2.
    \label{eq: variance of spectral noise}
\end{equation}

Regarding the pipeline, the output of each stage is listed in Table~\ref{tab:pipeline-notation}.
\begin{table}[tb]
\centering
\caption{Notation for pipeline inputs and outputs.}
\label{tab:pipeline-notation}
\begin{tabular}{lll}
\toprule
Notation & Description \\
\midrule
$[\mb X^*, \mb Y^*, \bo \zeta^*]$ & Ground-truth latent factors, intercept, and signal matrix \\
$(\mb R, \Omega)$ & Data matrix and observed locations \\
$[\widehat {\mb X}^{\sf{sp}}, \widehat {\mb Y}^{\sf{sp}}, \widehat {\bo \zeta}^{\sf{sp}}]$ & double-SVD output \\
 $[\widehat {\mb X}^{\ur}, \widehat{\mb Y}^{\ur}, \widehat {\bo \zeta}^{\ur}]$ & Unilateral refinement output \\
$[\mb X^{t}, \mb Y^{t}, \bo \zeta^{t} ]$ & the $t$-th step of the gradient descent \\
\bottomrule
\end{tabular}
\end{table}

    We collect the parameters of interest into a single vector
$\bo \theta \coloneqq (\sqrt{\frac{n}{\sigma_r^*}}\bo \zeta\t,\allowbreak \sf{vec}(\mb X)\t,\allowbreak \sf{vec}(\mb Y)\t)\t \in \bb R^{p + nr + pr}$, where the intercept component is suitably rescaled. For notational convenience, we occasionally denote $\bo \theta$ by $[\bo\zeta, \mb X, \mb Y]$ hereinafter. Under the theoretical calibration in Theorem~\ref{thm: parameter-explicit GD}, the objective function in \eqref{eq: loss function} becomes
\begin{align}
    & L([\bo \zeta, \mb X, \mb Y]) \coloneqq \sum_{(i,j) \in \Omega}\big[-R_{i,j}( \zeta_j + {\mb X_i}\t \mb Y_j) + \Psi( \zeta_j + {\mb X_i}\t \mb Y_j)\big] + c_{\perp}\frac{\pi \sigma_r^*}{n\sqrt{\omega}}\bignorm{\mb 1_n\t \mb X}_2^2.
\end{align}  
For convenience, we introduce the following shorthands for the derivatives of $L([\bo \zeta, \mb X, \mb Y])$: 
\begin{align} 
    &  \nabla_{\bo \zeta} L([\bo \zeta, \mb X, \mb Y]) \coloneqq \mc P_\Omega(\psi(\mb 1_n {\bo\zeta}\t+\mb X{\mb Y}\t) - \mb R)\t \mb 1_n  \in \bb R^{p}, \\ 
    & \nabla_{\mb X} L([\bo \zeta, \mb X, \mb Y]) \coloneqq \mc P_\Omega(\psi(\mb 1_n {\bo\zeta}\t+\mb X{\mb Y}\t) - \mb R) \mb Y + 2 c_\perp\frac{\pi\sigma_r^*}{n\sqrt{\omega}} \mb 1_n \mb 1_n\t \mb X \in \bb R^{n\times r}, \\
    & \nabla_{\mb Y} L ([\bo \zeta, \mb X, \mb Y]) \coloneqq \mc P_\Omega(\psi(\mb 1_n {\bo\zeta}\t+\mb X{\mb Y}\t) - \mb R)\t \mb X \in \bb R^{p \times r}.
    \label{eq: forms of derivatives}
\end{align}

\paragraph{\texorpdfstring{Scaling relation between $\mb X$ and $\mb Y$}{Scaling relation between X and Y}.}
As discussed earlier, we introduce a prescribed scaling parameter $\omega>0$ to accommodate settings in which the two latent factors operate on different scales. Nevertheless, it is sufficient to restrict attention to the balanced case ${\mb X^*}^{\top}\mb X^* = {\mb Y^*}^{\top}\mb Y^*$ without loss of generality, as shown by the following lemma.

\begin{lemma}
\label{lemma: scaling equivalence}
Let the superscript $\omega$ denote the output of the algorithm under scaling parameter $\omega$. Suppose that the refinement radii satisfy $\iota_{\mb X}^{(\omega)}=\omega^{1/4}\iota_{\mb X}^{(1)}$ and $\iota_{\mb Y}^{(\omega)}=\omega^{-1/4}\iota_{\mb Y}^{(1)}$, and that the gradient stage uses the same centering parameter, $c_\perp^{(\omega)}=c_\perp^{(1)}$, with step sizes satisfying
$\eta_\zeta^{(\omega)}=\eta_\zeta^{(1)}$,
$\eta_{\mb X}^{(\omega)}=\omega^{1/2}\eta_{\mb X}^{(1)}$, and
$\eta_{\mb Y}^{(\omega)}=\omega^{-1/2}\eta_{\mb Y}^{(1)}$.
Then the outputs of the algorithmic pipeline with scaling parameter $\omega$ are related to those with $\omega = 1$ through the following identities:
\begin{align}
    & \hat{\mb X}^{\sp,\omega}\omega^{-\frac14} = \hat{\mb X}^{\sp,1}, 
    \qquad
    \hat{\mb Y}^{\sp,\omega}\omega^{\frac14} = \hat{\mb Y}^{\sp,1}, \\
    & \hat{\mb X}^{\ur,\omega}\omega^{-\frac14} = \hat{\mb X}^{\ur,1}, 
    \qquad
    \hat{\mb Y}^{\ur,\omega}\omega^{\frac14} = \hat{\mb Y}^{\ur,1}, \\
    & \hat{\mb X}^{t,\omega}\omega^{-\frac14} = \hat{\mb X}^{t,1},
    \qquad
    \hat{\mb Y}^{t,\omega}\omega^{\frac14} = \hat{\mb Y}^{t,1},
    \qquad \forall\, t \in [t_0].
\end{align}
\end{lemma}

\begin{proof}
The first two identities follow from a change of variables. The centering step
commutes with this rescaling, while the intercept estimate is left unchanged.
The SVD refactorization then introduces the factors $\omega^{1/4}$ and
$\omega^{-1/4}$. The factor $\omega^{-1/2}$ in the centering penalty makes the
entire objective invariant under this rescaling. Consequently, the scaled
gradient updates preserve the relation inductively.
\end{proof}

As a result, the scaling parameter $\omega$ is assumed to be $1$ throughout the remaining analysis. 

\paragraph{Example of Parameters in Assumption~\ref{assumption: link function}.}
\label{subsec: link discussion}

We illustrate Assumption~\ref{assumption: link function} using the sparse
logistic link.
\begin{lemma}
\label{lemma: bernoulli smoothness}
Let $\psi(x)=1/(1+e^{-x})$ and
$\mc D_\psi=[\zeta_0-b,\zeta_0+b]$, where $b=O(1)$ and
$\zeta_0+b\leq-1$. For every $k\geq1$,
\begin{equation*}
    \sup_{x\in\mc D_\psi}|\psi^{(k)}(x)|
    \leq C e^{\zeta_0+b}k!2^k.
\end{equation*}
Moreover, $\psi'(x)\asymp e^x$ uniformly on $\mc D_\psi$. Consequently,
Assumption~\ref{assumption: link function} holds with
$\sigma^2\asymp e^{\zeta_0+b}$, $L_\psi\asymp1$,
$\underline c_\psi\asymp e^{\zeta_0-b}$, and
$\kappa_\psi\asymp e^{2b}$.
\end{lemma}
\begin{proof}
For $x\in\mc D_\psi$ and $|z-x|=1/2$, one has
$\operatorname{Re}(z)\leq\zeta_0+b+1/2\leq-1/2$. Hence
$|\psi(z)|=|e^z/(1+e^z)|\leq C e^{\zeta_0+b}$. Cauchy's integral formula
therefore gives
\begin{equation*}
    |\psi^{(k)}(x)|
    \leq k!2^k\sup_{|z-x|=1/2}|\psi(z)|
    \leq C e^{\zeta_0+b}k!2^k.
\end{equation*}
The curvature claim follows from
$\psi'(x)=e^x/(1+e^x)^2\asymp e^x$ on the stated interval.
\end{proof}

\section{Theoretical Analysis of Double-SVD Initialization}
\label{sec: proof of usvt}
\subsection{A General Version of Theorem~\ref{thm: USVT}}
\label{sec: proof of USVT}
We provide a more general version that includes the leave-one-out error characterization. Specifically, we aim to establish the stability of the output of Algorithm~\ref{alg: USVT} when the noise in a single row or column is removed. 

We first define the SVD of $\mb R^*$ as $\bar{\mb U}^*\bar{\bo \Sigma}^* \bar{\mb V}^{*\top}$ with $\bar{\bo \Sigma}^* = \diag(\bar \sigma_l^*)$. For arbitrary $i\in[n]$ and $j\in[p]$, define the coupled observation matrices as follows:
\begin{align*}
    & (\mb R^{(i)})_{m_1, m_2} = \begin{cases}
        \pi R_{i,m_2}^* & \text{if } m_1 = i , \\ 
        0  & \text{if } m_1 \neq i, (m_1,m_2) \notin \Omega, \\ 
        R_{m_1,m_2} & \text{otherwise, }
    \end{cases}; \\
    & (\mb R^{(-j)})_{m_1,m_2} = \begin{cases}
        \pi R_{m_1,j}^* & \text{if } m_2 = j , \\ 
        0  & \text{if } m_2 \neq j, (m_1,m_2) \notin \Omega, \\ 
        R_{m_1,m_2} & \text{otherwise. }
    \end{cases}
\end{align*}

Further, we denote the number of the non-vanishing singular values after the screening of Algorithm~\ref{alg: USVT} by $\hat r^{\sf{sp}} \coloneqq \big|\big\{ 
k\in [n\wedge p]: \bar \sigma_k \geq  \lambda^{\sp}
\big\} \big|$.  
Given the leave-one-out versions of the observation matrix, we define the corresponding low-rank approximations as 
\begin{equation}
    \begin{aligned}
    \hat{\mb R}^{\sp,(i)} &\coloneqq \pi^{-1} \bar{\mb U}^{(i)} \operatorname{diag}\big( ( \bar\sigma_l^{(i)}-\lambda^{\sp})_+ \big){\bar{\mb V}^{(i)\top}}, \\
    \hat{\mb R}^{\sp,(-j)} &\coloneqq \pi^{-1}\bar{\mb U}^{(-j)} \operatorname{diag}\big( ( \bar\sigma_l^{(-j)}-\lambda^{\sp})_+ \big) {\bar{\mb V}^{(-j)\top}},
    \end{aligned}
\end{equation}
where $\bar{\mb U}^{(i)} \bar{\bo \Sigma}^{(i)} \bar{\mb V}^{(i)}$ is the SVD of $\mb R^{(i)} $, and $\bar{\bo \Sigma}^{(i)} = \diag(\bar \sigma_l^{(i)})$ collects the singular values of $\mb R^{(i)}$. 
The same definitions apply for $\bar{\mb U}^{(-j)}$, $\bar{\mb V}^{(-j)}$, and $\bar{\bo \Sigma}^{(-j)}$ corresponding to $\mb R^{(-j)}$.  
Based on these low-rank approximations, we denote the resulting final outputs by $[\hat{\bo\zeta}^{\sp,(i)}, \hat{\mb X}^{\sp,(i)}, \hat{\mb Y}^{\sp,(i)}]$ and $[\hat{\bo\zeta}^{\sp,(-j)}, \hat{\mb X}^{\sp,(-j)}, \hat{\mb Y}^{\sp,(-j)}]$, respectively.

\begin{theorem}[A Detailed Characterization of Algorithm~\ref{alg: USVT}]
\label{thm: detailed USVT}
Suppose that Assumptions~\ref{assumption: noise and sampling},
\ref{assumption: link function},
\ref{assumption: signal strength and incoherence degree}, and
\ref{assumption: identifiability} hold, and choose
\begin{equation*}
    \lambda^{\sp}
    \asymp
    \sigma\sqrt{\bar\sigma\sigma_1^*}\,
    \pi^{\frac34}d^{\frac14}(\mu r)^{\frac14}\xi^{\frac32}
\end{equation*}
with a sufficiently large leading constant. Write
\begin{equation*}
    \mb F^*
    \coloneqq
    \begin{pmatrix}\mb X^*\\ \mb Y^*\end{pmatrix},
    \qquad
    \widehat{\mb F}^{\sp}
    \coloneqq
    \begin{pmatrix}\widehat{\mb X}^{\sp}\\ \widehat{\mb Y}^{\sp}\end{pmatrix},
    \qquad
    \widehat{\mb F}^{\sp,(l)}
    \coloneqq
    \begin{pmatrix}
        \widehat{\mb X}^{\sp,(l)}\\
        \widehat{\mb Y}^{\sp,(l)}
    \end{pmatrix},
    \qquad l\in[n]\cup[-p].
\end{equation*}
Define the common Procrustes rotations
\begin{align*}
    \mb O^{\sp}
    &\in\argmin_{\mb O\in\mc O(r)}
    \bignorm{\widehat{\mb F}^{\sp}\mb O-\mb F^*}\fb,
    &
    \mb O^{\sp,(l)}
    &\in\argmin_{\mb O\in\mc O(r)}
    \bignorm{\widehat{\mb F}^{\sp,(l)}\mb O
    -\widehat{\mb F}^{\sp}}\fb,
    \quad l\in[n]\cup[-p].
\end{align*}
Then, with probability at least $1-O(d^{-c})$,
        \begin{align}
        &\bignorm{\hat{\bo \zeta}^{\sf{sp}} - \bo \zeta^*}_2 \lesssim \frac{\lambda^{\sp} \xi}{\underline c_\psi\pi\sqrt{n}}, \quad \bignorm{\hat{\mb X}^{\sf{sp}}\mb O^{\sp}- \mb X^*}\fb \vee \bignorm{\hat{\mb Y}^{\sf{sp}}\mb O^{\sp}- \mb Y^*}\fb  \lesssim \frac{\kappa \lambda^{\sp} \xi}{\underline c_\psi\pi\sqrt{\sigma_r^*}},
        \label{eq: USVT zeta X and Y bounds}
        \\ 
        &\begin{aligned}
        & \underline c_\psi  \! \! \max_{l \in[n]\cup [-p]}\Big\{ \sqrt{ \pi n}\bignorm{\hat{\bo \zeta}^{\sp} - \hat{\bo \zeta}^{\sp, (l)}}_2, \! \sqrt{\pi\sigma_r^*}\bignorm{\hat{\mb X}^{\sf{sp}, (l)} \mb O^{\sp,(l)}-\hat{\mb X}^{\sf{sp}}}\fb, \\ 
        & \qquad \sqrt{\pi \sigma_r^*}\bignorm{\hat{\mb Y}^{\sf{sp}, (l)} \mb O^{\sp,(l)}-\hat{\mb Y}^{\sf{sp}}}\fb \Big\}  
        \lesssim    \frac{\kappa\sigma^2( r^*)^2 \bar\sigma \pi \sigma_1^* \sqrt{ \mu r \rho  }}{ \lambda^{\sp}}  \eqqcolon \xi^{\sf{loo}},
        \end{aligned}
         \label{eq: USVT zeta X Y loo error}
    \end{align}
    where $\xi \coloneqq \big[c_\xi r\log d\,(L_\psi \mu r^2 \sigma_1^* / \sqrt{np}\vee1)\big]^{\frac r2}$ for some sufficiently large constant $c_\xi$.
\end{theorem}

\begin{remark}[Baseline-centered imputation]
\label{remark: alternative imputation scheme}
One may apply the first singular value thresholding step to
$\mc P_\Omega\big(\mb R-\psi(\zeta_0)\mb 1_n\mb 1_p\t\big)$ and add
$\psi(\zeta_0)\mb 1_n\mb 1_p\t$ back before link inversion. The corresponding
spectral noise is
$\mc P_\Omega(\mb E)+(\bo\Omega-\pi\mb 1_n\mb 1_p\t)\circ
\{\mb R^*-\psi(\zeta_0)\mb 1_n\mb 1_p\t\}$. The mean-value theorem and
the stated incoherence bounds give
$\max_{i,j}|R_{i,j}^*-\psi(\zeta_0)|
\lesssim\sigma^2\mu r\sigma_1^*/\sqrt{np}$. Thus, the proof below applies
with $\bar\sigma$ replaced by
$\sigma\vee\{\sigma^2\mu r\sigma_1^*/\sqrt{np}\}$.
\end{remark}

In the sequel, we separately establish the global consistency bound \eqref{eq: USVT zeta X and Y bounds} in Section~\ref{subsec: proof of USVT part 1} and the leave-one-out discrepancy bound \eqref{eq: USVT zeta X Y loo error} in Section~\ref{subsec: proof of USVT part 2}.

\subsection{Consistency of the Spectral Estimator}
\label{subsec: proof of USVT part 1}

\begin{proof}[Proof of Eq.~\eqref{eq: USVT zeta X and Y bounds}]

We begin with the analysis of the soft singular value thresholding, which is largely aligned with the proof of Theorem~4 in \cite{xu2018rates}. This in turn allows us to have a characterization of the intercept and factor estimation. The prescribed threshold and Assumption~\ref{assumption: signal strength and incoherence degree}(a) imply $\lambda^{\sp}\gg\bar\sigma\sqrt{\pi d}\,\xi^4$ and $\lambda^{\sp}\xi/(\pi\underline c_\psi)\ll\sigma_r^*$, which are the two threshold consequences used below.

    \emph{Analysis of the soft singular value thresholding. } 
    We justify the error rate of $\hat{\mb R}^{\sp}$ with respect to $\mb R^*$ in Algorithm~\ref{alg: USVT}. Applying \citet[Lemma~1]{xu2018rates} gives that
    \begin{equation}
        \bignorm{\hat{\mb R}^{\sf{sp}} - \mb R^*}\fb^2 \leq 16 \pi^{-2} \min_{0 \leq s \leq n\wedge p} \Big[(\lambda^{\sp})^2 s + \big(1 + \frac{\norm{\mc P_\Omega(\mb R) - \pi \mb R^*}}{\lambda^{\sp} - \norm{\mc P_\Omega(\mb R) - \pi \mb R^*}}\big)^2 \pi^2 \sum_{i \geq s + 1}^{n \wedge p} (\bar \sigma_i^*)^2 \Big],
        \label{eq: bound in Lemma 1 xu2018}
    \end{equation}
    where $\bar \sigma_i^*$ is defined as the $i$-th largest singular value of $\mb R^*$. 

    Combining \eqref{eq: bound in Lemma 1 xu2018} together with Lemma~\ref{lemma: noise spectral norm concentration} immediately implies that 
    \begin{equation}
        \bignorm{\hat{\mb R}^{\sf{sp}} -\mb R^*}\fb^2 \leq 16\pi^{-2} \min_{0 \leq s \leq n\wedge p} \Big[( \lambda^{\sp})^2 s + 4\pi^2 \sum_{i \geq s + 1}^{n \wedge p}  (\bar \sigma_i^* )^2 \Big] \label{eq: error rate of USVT 1}
    \end{equation}
    holds with probability at least $1-O(d^{-c-2})$ since $\lambda^{\sp} \gg \bar\sigma\sqrt{\pi d}$. 

    Then we consider a uniform $l$-partition $\mc E_l 
    $ over the domain of $\mb X_i^*$. Specifically, let $M_X\coloneqq \|\mb X^*\|_\infty$. For each $q\in[l]$, define
$$
I_q
\coloneqq
\Big[
-M_X+\frac{2(q-1)M_X}{l},
-M_X+\frac{2qM_X}{l}
\Big),
$$
where the right endpoint of $I_l$ is included. For every
$\bo\nu=(\nu_1,\ldots,\nu_r)\in[l]^r$, let
$$
    E_{\bo\nu}
    \coloneqq
    I_{\nu_1}\times\cdots\times I_{\nu_r}
    =\prod_{a=1}^r I_{\nu_a},
$$
and define the occupied-cell index set
$$
    \mc I_l\coloneqq
    \big\{\bo\nu\in[l]^r:\mb X_i^*\in E_{\bo\nu}
    \text{ for some }i\in[n]\big\}.
$$
For each $\bo\nu\in\mc I_l$, choose an index $i_{\bo\nu}$ such that
$\mb X_{i_{\bo\nu}}^*\in E_{\bo\nu}$ and set
$\mb z_{\bo\nu}\coloneqq\mb X_{i_{\bo\nu}}^*$. Then
$$
    \mc E_l\coloneqq\{E_{\bo\nu}:\bo\nu\in[l]^r\}
$$
partitions $[-M_X,M_X]^r$ into $l^r$ hyperrectangles.

For each $\bo\nu\in\mc I_l$ and $j\in[p]$, define the order-$(k-1)$
Taylor polynomial centered at $\mb z_{\bo\nu}$ by
\begin{align}
    \mc P_{\bo\nu,k-1}^{(j)}(\mb x)
    \coloneqq
    \sum_{a=0}^{k-1}
    \frac{\psi^{(a)}(\zeta_j^*+\mb z_{\bo\nu}^{\top}\mb Y_j^*)}{a!}
    \big\{(\mb x-\mb z_{\bo\nu})^{\top}\mb Y_j^*\big\}^{a}.
\end{align}
Here and below, $\psi^{(0)}\coloneqq\psi$. Assemble these cellwise
polynomials into the piecewise-polynomial approximation
\begin{align}
    \mc P_{k-1}([\zeta_j^*,\mb x,\mb Y_j^*])
    \coloneqq
    \sum_{\bo\nu\in\mc I_l}
    \mc P_{\bo\nu,k-1}^{(j)}(\mb x)
    \ind\{\mb x\in E_{\bo\nu}\}.
\end{align}
If $\mb X_i^*\in E_{\bo\nu}$, then
$\|\mb X_i^*-\mb z_{\bo\nu}\|_\infty\leq 2M_X/l$. Moreover, both
$\zeta_j^*+{\mb X_i^*}^{\top}\mb Y_j^*$ and
$\zeta_j^*+\mb z_{\bo\nu}^{\top}\mb Y_j^*$ are entries of the true
predictor matrix and hence belong to $\mc D_\psi$; their connecting segment
also lies in $\mc D_\psi$. Taylor's theorem and
Assumption~\ref{assumption: link function} therefore give
\begin{align}
    &\Big|
        \psi(\zeta_j^*+{\mb X_i^*}^{\top}\mb Y_j^*)
        -\mc P_{k-1}([\zeta_j^*,\mb X_i^*,\mb Y_j^*])
    \Big|
    \leq
    \sigma^2L_\psi^{k-1}
    \Big(
        \frac{2r\|\mb X^*\|_\infty\|\mb Y^*\|_\infty}{l}
    \Big)^k.
    \label{eq: piecewise polynomial approximation}
\end{align} 

    Define $\mb N\in\bb R^{n\times p}$ by
$
    N_{ij}
    \coloneqq
    \mc P_{k-1}([\zeta_j^*,\mb X_i^*,\mb Y_j^*])$, 
which satisfies  
    \begin{equation}
    \operatorname{rank}(\mb N)
    \leq
    l^{r} \sum_{i=0}^{k-1} \binom{i+r-1}{r-1},
    \label{eq: upper bound for rank}
\end{equation}
since $|\mc I_l|\leq l^r$ and, within each occupied cell, the degree-$a$
multinomial expansion consists of $\binom{a+r-1}{r-1}$ rank-one terms.
Define $l \coloneqq  \lceil  e \big(2 r \bignorm{\mb X^*}_\infty\bignorm{\mb Y^*}_\infty L_\psi \vee 1 \big) \rceil$ and $r_0\coloneqq \min\{s \in \bb N^+:
    s \geq l, \quad
    s^{\frac{1}{r}} \geq 2l \log s\} \wedge (n\wedge p)$.

    For every $s \geq r_0$, there is an integer $k_0$ such that $l^{r} \sum_{i=0}^{k_0-1}\binom{i+r-1}{r-1} \leq s \leq l^{r} \sum_{i=0}^{k_0}\binom{i+r-1}{r-1}$. Substituting this choice of $k_0$ into the definition of $\mb N$, we obtain that
    \begin{equation}
        \begin{aligned}
        \sum_{i \geq s +1}^{n\wedge p} (\bar \sigma_i^*)^2
        &\leq \bignorm{\mb R^* - \mb N}\fb^2 \\
        &\leq d^2 \Big[\sigma^2L_\psi^{k_0-1}\big(2r\bignorm{\mb X^*}_\infty\bignorm{\mb Y^*}_\infty/l\big)^{k_0}\Big]^2 \\
        &\lesssim \sigma^4d^2 e^{-2 k_0 - 2r}  s^2
        \lesssim \sigma^4d^2 e^{- s^{\frac{1}{r}} l^{-1}}.
        \end{aligned}
        \label{eq: low-rank approximation residual error control}
    \end{equation}
    Here, the last inequality is reasoned by 
    \begin{equation}
        e^{-2r-2k_0} s^2  \leq e^{-2(r+k_0)} e^{s^{\frac{1}{r}}l^{-1}} \leq e^{-2 s^{\frac{1}{r}} l^{-1}} e^{s^{\frac{1}{r}}l^{-1}} = e^{- s^{\frac{1}{r}} l^{-1}} 
    \end{equation}
    according to the definition of $r_0$ and the fact $ s\leq l^r \sum_{i=0}^{k_0} \binom{i+r-1 }{ r-1} \leq l^{r}(k_0+1)(k_0 + r-1)^{r-1} \leq l^{r}(k_0+r)^{r}$.

    Set $s_\psi\coloneqq\lceil c_1rl\log d\rceil^r\wedge(n\wedge p)$ for a
    sufficiently large numerical constant $c_1$. Since $l$ is bounded by a
    constant multiple of
    $(L_\psi\mu r^2\sigma_1^*/\sqrt{np})\vee1$, the consolidated
    sampling condition implies $\log(c_1rl\log d)\lesssim\log d$.
    Enlarging $c_1$ therefore gives
    $s_\psi^{1/r}\geq2l\log s_\psi$, and hence $s_\psi\geq r_0$, whenever
    $s_\psi<n\wedge p$. In that case,
    \eqref{eq: low-rank approximation residual error control} and the
    bounds $\bar\sigma\geq\sigma$ and $|\log\sigma|\lesssim\log d$ give
    $\sigma^4d^2\exp(-s_\psi^{1/r}/l)\lesssim
    \bar\sigma^2/(\pi d)$; if $s_\psi=n\wedge p$, the spectral tail is zero.
    Consequently, \eqref{eq: upper bound for rank},
    \eqref{eq: error rate of USVT 1}, and the definition of $\lambda^{\sp}$
    imply
    \begin{align}
    & \sum_{i \geq s_\psi +1}^{n\wedge p} (\bar \sigma_i^*)^2 \lesssim \bar\sigma^2 / (\pi d)
    \label{eq: upper bound for F approximation error} \\ 
        & \bignorm{\hat{\mb R}^{\sf{sp}} - \mb R^*}\fb^2
        \lesssim \pi^{-2}\Big[(\lambda^{\sp})^2s_\psi
        +4\pi^2\sigma^4d^2
        \exp(-s_\psi^{\frac1r}l^{-1})\Big]
        \lesssim (\lambda^{\sp})^2\xi^2/\pi^2
        \label{eq: upper bound for hat R - pi R^*}, 
    \end{align}
    where the last inequality follows from $s_\psi\leq\xi^2$ after taking
    $c_\xi$ sufficiently large. Thus, the polynomial-scale condition on
    $\sigma$ absorbs its contribution to the approximation rank into the
    $\log d$ factor.

    \emph{Analysis of inversion and truncation. } It follows immediately from \eqref{eq: upper bound for hat R - pi R^*} and the validity of the truncation range that, with probability at least $1-O(d^{-c-2})$,
    \begin{align}
        & \bignorm{\operatorname{clip}_{\psi(\mc D_\psi)}(\hat{\mb R}^{\sf{sp}}) - \mb R^*}\fb \lesssim \lambda^{\sp} \xi / \pi. 
    \end{align}
    Applying the inverse mapping to obtain $\check{\bo\Theta}$ and using the lower bound on $\psi'$, we further have
    \begin{align}
        & \bignorm{\hat{\mb M}^{\sp} - \mb M^*}\fb \lesssim \frac{\lambda^{\sp} \xi}{\pi \underline c_\psi}
    \end{align}
    with probability at least $1-O(d^{-c-2})$. 

    Recalling that $\hat{\bo \zeta}^{\sf{sp}} = \big(\mb 1_n\t\hat{\mb M}^{\sf{sp}}\big)\t /n$, we directly have 
    \begin{equation}
        \bignorm{\hat{\bo \zeta}^{\sf{sp}} - \bo \zeta^*}_2 \lesssim \frac{\lambda^{\sp} \xi}{\underline c_\psi \pi \sqrt{n}}. 
        \label{eq: USVT zeta error bound}
    \end{equation}
    In addition, we notice $(\mb I - \mb 1_n \mb 1_n\t  /n ) \mb M^* = \mb X^* \mb Y^{*\top}$ and derive that
    \begin{equation}
        \norm{(\mb I - \mb 1_n\mb 1_n\t /n) \hat{\mb M}^{\sf{sp}} - \mb X^*\mb Y^{*\top}} \lesssim  \frac{\lambda^{\sp} \xi}{\pi \underline c_\psi} \label{eq: estimation error of XYt for spectral method}. 
    \end{equation}
    Since the centering matrix is an orthogonal projection, the Frobenius-norm counterpart also obeys
    \begin{equation}
        \bignorm{(\mb I - \mb 1_n\mb 1_n\t/n)\hat{\mb M}^{\sp}-\mb X^*\mb Y^{*\top}}\fb
        \leq \bignorm{\hat{\mb M}^{\sp}-\mb M^*}\fb
        \lesssim \frac{\lambda^{\sp}\xi}{\pi\underline c_\psi}.
        \label{eq: frobenius estimation error of XYt for spectral method}
    \end{equation}

    It remains to control the two factor matrices under a common rotation.  Let
    \begin{equation*}
        \hat{\mb A}^{\sp}
        \coloneqq
        \big(\mb I_n-\mb 1_n\mb 1_n\t/n\big)\hat{\mb M}^{\sp},
        \qquad
        \mb A^*\coloneqq\mb X^*\mb Y^{*\top},
    \end{equation*}
    and write their top-$r$ SVDs as
    $\hat{\mb U}^{\sp}\hat{\bo\Sigma}^{\sp}\hat{\mb V}^{\sp\top}$
    and $\mb U^*\bo\Sigma^*\mb V^{*\top}$, respectively. Under the balanced
    parametrization adopted in this supplement, there exists a deterministic
    $\mb Q^*\in\mc O(r)$ such that
    \begin{equation*}
        \mb X^*=\mb U^*(\bo\Sigma^*)^{1/2}\mb Q^*,
        \qquad
        \mb Y^*=\mb V^*(\bo\Sigma^*)^{1/2}\mb Q^*.
    \end{equation*}
    Define the positive-eigenspace bases of the corresponding Hermitian
    dilations by
    \begin{align*}
        \hat{\mb W}^{\sp}
        \coloneqq \frac{1}{\sqrt 2}
        \begin{pmatrix}\hat{\mb U}^{\sp}\\ \hat{\mb V}^{\sp}\end{pmatrix},
        \qquad
        \mb W^*
        \coloneqq \frac{1}{\sqrt 2}
        \begin{pmatrix}\mb U^*\\ \mb V^*\end{pmatrix},\\
        \bar{\mb O}^{\sp}
        \in\arg\min_{\mb O\in\mc O(r)}
        \bignorm{\hat{\mb W}^{\sp}\mb O-\mb W^*}\fb,
        \qquad
        \mb O_{\sf H}^{\sp}\coloneqq\bar{\mb O}^{\sp}\mb Q^*.
    \end{align*}
    Under the prescribed threshold, Assumption~\ref{assumption: signal strength and incoherence degree}(a) and \eqref{eq: estimation error of XYt for spectral method} give
    \begin{equation}
        \norm{\hat{\mb A}^{\sp}-\mb A^*}
        \lesssim \frac{\lambda^{\sp}\xi}{\pi\underline c_\psi}
        \leq c_0\sigma_r^*,
        \label{eq: spectral estimator perturbation neighborhood}
    \end{equation}
    where $c_0$ is the numerical constant in
    Lemma~\ref{lemma: perturbation theory for joint factors}; the last
    inequality follows from
    Assumption~\ref{assumption: signal strength and incoherence degree}(a).
    Applying Lemma~\ref{lemma: perturbation theory for joint factors} with
    $\mb M^{(1)}=\hat{\mb A}^{\sp}$ and $\mb M^{(2)}=\mb A^*$ shows that both
    $\norm{\hat{\mb X}^{\sp}\mb O_{\sf H}^{\sp}-\mb X^*}\fb$ and
    $\norm{\hat{\mb Y}^{\sp}\mb O_{\sf H}^{\sp}-\mb Y^*}\fb$ are bounded by
    $\kappa\lambda^{\sp}\xi/
    (\underline c_\psi\pi\sqrt{\sigma_r^*})$ up to a universal constant.
    By the Procrustes optimality of $\mb O^{\sp}$ defined in the theorem, its
    joint factor error is no larger than that under $\mb O_{\sf H}^{\sp}$.
    Consequently,
    \begin{equation}
        \norm{\hat{\mb X}^{\sp}\mb O^{\sp}-\mb X^*}\fb
        \vee
        \norm{\hat{\mb Y}^{\sp}\mb O^{\sp}-\mb Y^*}\fb
        \lesssim
        \frac{\kappa\lambda^{\sp}\xi}
        {\underline c_\psi\pi\sqrt{\sigma_r^*}}.
        \label{eq: joint spectral norm error of spectral factor estimates}
    \end{equation}

\end{proof}

\subsection{Leave-One-Out Analysis for Algorithm~\ref{alg: USVT}}
\label{subsec: proof of USVT part 2}
This subsection develops the leave-one-out perturbation bound \eqref{eq: USVT zeta X Y loo error} needed for the subsequent unilateral refinement and gradient descent analysis.

The following lemma is concerned with the leave-one-out error with respect to $\hat{\mb R}^{\sp}$, whose proof is deferred to Section~\ref{subsubsec: proof of leave-one-out error for USVT}. 
\begin{lemma}
    \label{lemma: leave-one-out for the first USVT}
    Suppose that Assumptions~\ref{assumption: noise and sampling},
    \ref{assumption: link function},
    \ref{assumption: signal strength and incoherence degree}, and
    \ref{assumption: identifiability} hold, and take $\lambda^{\sp}$ as in
    Theorem~\ref{thm: detailed USVT}.
    Then, with probability at least $1-O(d^{-c})$, the following bound holds uniformly over all leave-one-row-out and leave-one-column-out versions of the spectral algorithm:
    \begin{align}
        & \max_{i\in[n]}\norm{\hat{\mb R}^{\sp,(i)} - \hat{\mb R}^{\sp}}\fb
        \vee \max_{j\in[p]}\norm{\hat{\mb R}^{\sp,(-j)} - \hat{\mb R}^{\sp}}\fb
        \lesssim \frac{\sigma^2(r^*)^2\bar\sigma\sigma_1^*\sqrt{\mu r\rho\pi}}{\lambda^{\sp}} .
    \end{align} 
\end{lemma}

\bigskip

\begin{proof}[Proof of Eq.~\eqref{eq: USVT zeta X Y loo error}. ]
    By applying $\psi^{-1}$ to the truncated matrices, it follows that, with probability at least $1 - O(d^{-c})$,
\begin{align}
    \max_{i\in[n]} \bignorm{\hat{\mb M}^{\sf{sp},(i)} - \hat{\mb M}^{\sp}}\fb
    \vee
    \max_{j\in[p]} \bignorm{\hat{\mb M}^{\sf{sp},(-j)} - \hat{\mb M}^{\sp}}\fb
    \lesssim
    \frac{\sigma^2( r^*)^2 \bar\sigma  \sigma_1^* \sqrt{ \mu r \rho \pi }}{\underline c_\psi \lambda^{\sp}}.
\end{align}

    For $l\in[n]\cup[-p]$, set
    \begin{equation*}
        \hat{\mb A}^{\sp,(l)}
        \coloneqq
        \big(\mb I_n-\mb 1_n\mb 1_n\t/n\big)\hat{\mb M}^{\sp,(l)},
    \end{equation*}
    and let $\hat{\mb U}^{\sp,(l)}\hat{\bo\Sigma}^{\sp,(l)}
    \hat{\mb V}^{\sp,(l)\top}$ be its top-$r$ SVD.  Define
    \begin{align*}
        \hat{\mb W}^{\sp,(l)}
        \coloneqq \frac{1}{\sqrt 2}
        \begin{pmatrix}\hat{\mb U}^{\sp,(l)}\\
        \hat{\mb V}^{\sp,(l)}\end{pmatrix},\\
        \mb O_{\sf H}^{\sp,(l)}
        \in\arg\min_{\mb O\in\mc O(r)}
        \bignorm{\hat{\mb W}^{\sp,(l)}\mb O-\hat{\mb W}^{\sp}}\fb.
    \end{align*}
    By \eqref{eq: spectral estimator perturbation neighborhood} and Weyl's
    inequality, the singular values of $\hat{\mb A}^{\sp}$ satisfy
    \begin{equation*}
        \hat\sigma_r^{\sp}-\hat\sigma_{r+1}^{\sp}
        \gtrsim\sigma_r^*,
        \qquad
        \frac{\hat\sigma_1^{\sp}}{\hat\sigma_r^{\sp}}
        \lesssim\kappa.
    \end{equation*}
    In addition, Lemma~\ref{lemma: leave-one-out for the first USVT}, the
    nonexpansiveness of left multiplication by
    $\mb I_n-\mb 1_n\mb 1_n\t/n$, and the bound $r^*\lesssim\xi^2$ in
    \eqref{eq: upper bound on r*} give
    \begin{align}
        \max_{l\in[n]\cup[-p]}
        \norm{\hat{\mb A}^{\sp,(l)}-\hat{\mb A}^{\sp}}
        &\leq
        \max_{l\in[n]\cup[-p]}
        \bignorm{\hat{\mb M}^{\sp,(l)}-\hat{\mb M}^{\sp}}\fb\nonumber\\
        &\lesssim
        \frac{\sigma^2(r^*)^2\bar\sigma\sigma_1^*\sqrt{\mu r\rho\pi}}
        {\underline c_\psi\lambda^{\sp}}
        \leq c_0\hat\sigma_r^{\sp}.
        \label{eq: leave-one-out perturbation neighborhood}
    \end{align}
    Here, the last inequality follows from the threshold choice and the
    signal-strength and sampling conditions. Apply
    Lemma~\ref{lemma: perturbation theory for joint factors} with
    $\mb M=\hat{\mb A}^{\sp,(l)}$ and
    $\mb M^*=\hat{\mb A}^{\sp}$. Its reference factors are
    $(\hat{\mb X}^{\sp},\hat{\mb Y}^{\sp})$, with $\mb Q^*=\mb I_r$.
    The preceding two displays imply the factor rate below with
    $\mb O_{\sf H}^{\sp,(l)}$ in place of $\mb O^{\sp,(l)}$, uniformly over
    $l\in[n]\cup[-p]$. By the Procrustes optimality of $\mb O^{\sp,(l)}$
    defined in the theorem, its joint factor error is no larger. The intercept
    bound is unaffected by the rotation. Consequently,
    \begin{align}
        &\bignorm{\hat{\mb X}^{\sp,(l)}\mb O^{\sp,(l)}
            -\hat{\mb X}^{\sp}}\fb
        \vee
        \bignorm{\hat{\mb Y}^{\sp,(l)}\mb O^{\sp,(l)}
            -\hat{\mb Y}^{\sp}}\fb
        \lesssim
        \frac{\xi^{\sf{loo}}}
        {\underline c_\psi\sqrt{\pi\sigma_r^*}},
        \nonumber\\
        &\bignorm{\hat{\bo\zeta}^{\sp}
            -\hat{\bo\zeta}^{\sp,(l)}}_2
        \lesssim
        \frac{\xi^{\sf{loo}}}
        {\underline c_\psi\sqrt{\pi n}}.
        \label{eq: leave-one-out perturbation for USVT output estimations}
    \end{align}

    For later use, we pass these estimates through the trimming step at the
    beginning of Algorithm~\ref{alg: OS}. Let the deterministic radii contain
    the population intercept and factor rows in their respective projection
    sets and satisfy
    $\iota_{\bo\zeta}+\iota_{\mb X}^{(1)}\iota_{\mb Y}^{(1)}
    \leq c_{\mc D_\psi}\mu r\sigma_1^*/\sqrt{np}$. Define
    \begin{align*}
        \bo\zeta^0
        &\coloneqq\mc P^\infty_{\iota_{\bo\zeta},\zeta_0}
        (\hat{\bo\zeta}^{\sp}),\quad 
        \mb X^0 \coloneqq\mc P^{2,\infty}_{\iota_{\mb X}^{(1)}}
        (\hat{\mb X}^{\sp}),\quad 
        \mb Y^0
        \coloneqq\mc P^{2,\infty}_{\iota_{\mb Y}^{(1)}}
        (\hat{\mb Y}^{\sp}),\\
        \bo\zeta^{0,(l)}
        &\coloneqq\mc P^\infty_{\iota_{\bo\zeta},\zeta_0}
        (\hat{\bo\zeta}^{\sp,(l)}),\\
        &\quad
        \mb X^{0,(l)}
        \coloneqq\mc P^{2,\infty}_{\iota_{\mb X}^{(1)}}
        (\hat{\mb X}^{\sp,(l)}),\quad
        \mb Y^{0,(l)}
        \coloneqq\mc P^{2,\infty}_{\iota_{\mb Y}^{(1)}}
        (\hat{\mb Y}^{\sp,(l)}),
        \qquad l\in[n]\cup[-p].
    \end{align*}
    The rowwise projection is right-orthogonally equivariant and
    nonexpansive in the Frobenius norm:
    \begin{equation*}
        \mc P_a^{2,\infty}(\mb A\mb O)
        =\mc P_a^{2,\infty}(\mb A)\mb O,
        \qquad
        \bignorm{\mc P_a^{2,\infty}(\mb A)
            -\mc P_a^{2,\infty}(\mb B)}\fb
        \leq\bignorm{\mb A-\mb B}\fb.
    \end{equation*}
    Centered entrywise clipping is likewise nonexpansive. The chosen
    projection sets contain the population intercept and factor rows, so
    \begin{align*}
        \bignorm{\bo\zeta^0-\bo\zeta^*}_2
        &\leq\bignorm{\hat{\bo\zeta}^{\sp}-\bo\zeta^*}_2,\\
        \bignorm{\mb X^0\mb O^{\sp}-\mb X^*}\fb
        \vee\bignorm{\mb Y^0\mb O^{\sp}-\mb Y^*}\fb
        &\leq
        \bignorm{\hat{\mb X}^{\sp}\mb O^{\sp}-\mb X^*}\fb
        \vee\bignorm{\hat{\mb Y}^{\sp}\mb O^{\sp}-\mb Y^*}\fb.
    \end{align*}
    Hence the global bounds in
    \eqref{eq: USVT zeta X and Y bounds} pass unchanged to the trimmed input.
    The same contraction argument and rotational equivariance give, uniformly
    over $l\in[n]\cup[-p]$,
    \begin{align}
        &\bignorm{\bo\zeta^{0,(l)}-\bo\zeta^0}_2
        \lesssim
        \frac{\xi^{\sf{loo}}}{\underline c_\psi\sqrt{\pi n}},\nonumber\\
        &\bignorm{\mb X^{0,(l)}\mb O^{\sp,(l)}-\mb X^0}\fb
        \vee
        \bignorm{\mb Y^{0,(l)}\mb O^{\sp,(l)}-\mb Y^0}\fb
        \lesssim
        \frac{\xi^{\sf{loo}}}
        {\underline c_\psi\sqrt{\pi\sigma_r^*}}.
        \label{eq: leave-one-out perturbation for truncated estimates in step 2}
    \end{align}
    Finally, for every $i\in[n]$ and $j\in[p]$, the original initializer
    satisfies
    \begin{equation*}
        |\zeta_j^0-\zeta_0|
        +\bignorm{\mb X_i^0}_2\bignorm{\mb Y_j^0}_2
        \leq\iota_{\bo\zeta}
        +\iota_{\mb X}^{(1)}\iota_{\mb Y}^{(1)},
    \end{equation*}
    and, uniformly over $l\in[n]\cup[-p]$, every leave-one-out initializer
    satisfies
    \begin{equation*}
        |\zeta_j^{0,(l)}-\zeta_0|
        +\bignorm{\mb X_i^{0,(l)}}_2\bignorm{\mb Y_j^{0,(l)}}_2
        \leq\iota_{\bo\zeta}
        +\iota_{\mb X}^{(1)}\iota_{\mb Y}^{(1)}.
    \end{equation*}
    The radius condition therefore places these predictors, as well as the
    interpolation paths used in the refinement proof, inside $\mc D_\psi$.
    \end{proof}

\subsubsection{Proof of Lemma~\ref{lemma: leave-one-out for the first USVT}}
\label{subsubsec: proof of leave-one-out error for USVT}
We first remark on a useful byproduct of the preceding arguments that the number of singular values of $\mb R^*$ exceeding the threshold $\lambda^{\sp}$ is relatively small. Consequently, there must exist at least one singular value of non-negligible magnitude, and there must also be a nontrivial gap among these singular values. This observation suggests that the previous results are sufficient to support a leave-one-out perturbation analysis developed in \cite{abbe2020entrywise,chen2021spectral}.

Specifically, we group the eigenvalues into blocks separated by nontrivial eigengaps, enabling the application of established leave-one-out error control argument within each block. A potential concern is that the soft-thresholding step partially truncates some singular values, which could complicate the analysis. This issue is mitigated by observing that the soft-thresholding operator is Lipschitz, so controlling the unthresholded quantity is sufficient to derive the final bound.
    
    Let $
\mc E^{\sp}\coloneqq
\left\{\big\|\mc P_\Omega(\mb R)-\pi\mb R^*\big\|
\leq\lambda^{\sp}/2\right\}$.
By Lemma~\ref{lemma: noise spectral norm concentration} and Weyl's inequality,
$\mathbb P(\mc E^{\sp})\geq1-O(d^{-c-2})$, and, on $\mc E^{\sp}$,
$$
\hat r^{\sp}\leq r^*
\coloneqq
\left|\left\{i\in[n\wedge p]:
\pi\sigma_i(\mb R^*)\geq\lambda^{\sp}/2\right\}\right|.
$$
    We also introduce a slightly larger quantity $\bar r^*$ as follows 
    \begin{equation}
    \bar r^* \coloneqq  \left|\left\{i\in[n\wedge p]:
\pi\sigma_i(\mb R^*)\geq\lambda^{\sp}/8\right\}\right| \leq \big[c_r r\log d\,(L_\psi \mu r^2 \sigma_1^* / \sqrt{np}\vee1)\big]^{r}
    \label{eq: upper bound on r*}
    \end{equation}
    with some constant $c_r$ from the Taylor approximation in the proof of Theorem~\ref{thm: detailed USVT}. In fact, taking Eq. \eqref{eq: low-rank approximation residual error control} and the reasoning in \eqref{eq: upper bound for hat R - pi R^*} with $s=s_\psi$ tells us that
    $ \pi^2 \sum_{i \geq s_\psi+1}^{n\wedge p} (\bar \sigma_i^*)^2 \leq (\lambda^{\sp} / 2)^2$,
    indicating that there are at most $s_\psi$ singular values of $\pi \mb R^*$ that are greater than $ \lambda^{\sp} / 8$.
    
   \paragraph{A useful singular-value partition.}
We partition the leading singular values of the expectation matrix
$\mb R^*$. Write
\begin{equation*}
    \tau_s\coloneqq\bar\sigma_s^*,
    \qquad s\in[n\wedge p],
    \qquad\text{and}\qquad
    \tau_{n\wedge p+1}\coloneqq0.
\end{equation*}
Starting from $r_0=0$, whenever $r_{k-1}<r^*$, define
\begin{align}
    \mc I_k
    &\coloneqq
    \Big\{
        q\in\{r_{k-1}+1,\ldots,\bar r^*\}:
        \frac{\tau_{r_{k-1}+1}}{\tau_q}\leq4,\quad
        \tau_q-\tau_{q+1}\geq\frac{\tau_q}{2\bar r^*}
    \Big\}, \\
    r_k&\coloneqq\max\mc I_k,
    \qquad
    \mc S_k\coloneqq\{r_{k-1}+1,\ldots,r_k\}.
\end{align}
This construction is similar in spirit to that of
\citet{zhou2025deflated}.

We next verify that $\mc I_k$ is nonempty. Let
\begin{equation*}
    q_k\coloneqq
    \min\Big(
        \Big\{
            q\in\{r_{k-1}+1,\ldots,\bar r^*\}:
            \tau_q-\tau_{q+1}\geq\frac{\tau_q}{2\bar r^*}
        \Big\}
        \cup\{\bar r^*+1\}
    \Big).
\end{equation*}
Suppose, to the contrary, that $q_k=\bar r^*+1$. By the minimality of
$q_k$,
\begin{equation*}
    \tau_{\bar r^*+1}
    >
    \Big(1-\frac{1}{2\bar r^*}\Big)^{\bar r^*-r_{k-1}}
    \tau_{r_{k-1}+1}
    \geq
    \frac{1}{2}\tau_{r^*}
    \geq\frac{\lambda^{\sp}}{4\pi}.
\end{equation*}
This contradicts the definition of $\bar r^*$, which implies
$\tau_{\bar r^*+1}<\lambda^{\sp}/(8\pi)$. Hence
$q_k\leq\bar r^*$. Another application of the preceding recursion gives
\begin{equation*}
    \tau_{q_k}
    \geq
    \Big(1-\frac{1}{2\bar r^*}\Big)^{q_k-r_{k-1}-1}
    \tau_{r_{k-1}+1}
    \geq\frac{1}{2}\tau_{r_{k-1}+1}.
\end{equation*}
Thus $q_k\in\mc I_k$, proving that $\mc I_k$ is nonempty.

The recursion therefore terminates at
$\bar k\coloneqq\min\{k:r_k\geq r^*\}$. Its blocks are contiguous and
satisfy
\begin{equation*}
    \bigcup_{k=1}^{\bar k-1}\mc S_k
    \subseteq[r^*]
    \subseteq
    \bigcup_{k=1}^{\bar k}\mc S_k
    \subseteq[\bar r^*].
\end{equation*}
Moreover, when $\bar k\geq2$, the gaps between distinct blocks satisfy
\begin{equation}
    \min_{\substack{k_1,k_2\in[\bar k],\,k_1\neq k_2\\
                    i_1\in\mc S_{k_1},\,i_2\in\mc S_{k_2}}}
    \big|\bar\sigma_{i_1}^*-\bar\sigma_{i_2}^*\big|
    \geq
    \frac{\bar\sigma_{r^*}^*}{2\bar r^*}
    \geq
    \frac{\lambda^{\sp}}{4\pi\bar r^*}
    \gg\bar\sigma\sqrt{\frac{d}{\pi}}.
    \label{eq: lower bound on the gap among groups}
\end{equation}
Finally, the last block is separated from the lower tail since
$\tau_{r_{\bar k}}-\tau_{r_{\bar k}+1}
\geq\lambda^{\sp}/(16\pi\bar r^*)$.

    \paragraph{Applying Wedin's theorem.}
    The task of controlling the leave-one-out fluctuation of $\bar{\mb U}_{\cdot, 1:\hat r^{\sf{sp}}}$ and $\bar{\mb V}_{\cdot, 1:\hat r^{\sf{sp}}}$ can be decomposed as upper bounding smaller singular subspaces corresponding to the partition. Recall the spectral noise matrix $\bar{\mb E}$ defined in \eqref{eq: response and spectral noise matrices}. 
    A direct application of Wedin's theorem \citet[Theorem~4.1]{stewart1990matrix} together with \eqref{eq: upper bound for F approximation error}, \eqref{eq: lower bound on the gap among groups}, and Lemma~\ref{lemma: noise spectral norm concentration} yields 
    \begin{align}
        & \distf(\bar {\mb U}_{\cdot, \mc S_k}, \bar {\mb U}^{(i)}_{\cdot, \mc S_k}) \vee \distf(\bar{\mb V}_{\cdot, \mc S_k}, \bar{\mb V}^{(i)}_{\cdot, \mc S_k}) \lesssim \frac{r^* \bignorm{\bar{\mb E}_{i,\cdot} \bar{\mb V}^{(i)}_{\cdot, \mc S_k}}_2 + r^* \bignorm{\bar{\mb E}_{i}}_2 \bignorm{\bar{\mb U}^{(i)}_{\cdot, \mc S_k}}\ti  }{\pi \min_{i \in \mc S_k} \bar \sigma_i^*}, \label{eq: wedin theorem for leave-one-out step 1} 
    \end{align}

    For illustration, we examine the left singular vectors contained in $\bar{\mb U}_{\cdot,\mc S_k}^{(i)}$. 
    To upper bound the right-hand side of \eqref{eq: wedin theorem for leave-one-out step 1}, we invoke the Bernstein inequality to derive that, for each fixed $(i,k)$, with conditional probability at least $1-O(d^{-c-2})$ given $\mc P_{-i,\cdot}(\bar{\mb E})$, one has 
    \begin{align}
        & \bignorm{\bar{\mb E}_{i,\cdot} \bar{\mb V}^{(i)}_{\cdot, \mc S_k}} \lesssim \bar\sigma \sqrt{\pi r^* \log d } + B \bignorm{\bar{\mb V}^{(i)}_{\cdot, \mc S_k}}\ti \log d. 
        \label{eq: bernstein for wedin + leave-one-out}
    \end{align}
    Since $\bar k\leq r^*\leq n\wedge p\leq d$ and $n+p\leq2d$, a union bound over the row or column index and the block index shows that the row and column analogues of \eqref{eq: bernstein for wedin + leave-one-out} hold simultaneously with probability at least $1-O(d^{-c})$.
    By \eqref{eq: bernstein for wedin + leave-one-out} and Lemma~\ref{lemma: noise spectral norm concentration}, one has  
    \begin{align} 
    & \distf(\bar {\mb U}_{\cdot, \mc S_k}, \bar {\mb U}^{(i)}_{\cdot, \mc S_k}) \vee \distf(\bar{\mb V}_{\cdot, \mc S_k}, \bar{\mb V}^{(i)}_{\cdot, \mc S_k})  \\ 
        \lesssim &\frac{\bar\sigma r^* \sqrt{\pi r^* \log d } + r^* \big( B \bignorm{\bar{\mb V}^{(i)}_{\cdot, \mc S_k}}\ti \log d + \bar\sigma \sqrt{\pi d}\bignorm{\bar{\mb U}^{(i)}_{\cdot, \mc S_k}}\ti \big) }{\pi \min_{i \in \mc S_k} \bar \sigma_i^*}  \\ 
        \lesssim &  \frac{\bar\sigma r^* \sqrt{\pi r^* \log d }}{\pi \min_{i \in \mc S_k} \bar \sigma_i^*} + \frac{ r^* \big( B \bignorm{\bar{\mb V}_{\cdot, \mc S_k}}\ti \log d + \bar\sigma \sqrt{\pi d}\bignorm{\bar{\mb U}_{\cdot, \mc S_k}}\ti \big) }{\pi \min_{i \in \mc S_k} \bar \sigma_i^*} \\
         & + \frac{ r^* \big( B  \log d + \bar\sigma \sqrt{\pi d} \big) }{\pi \min_{i \in \mc S_k} \bar \sigma_i^*}  \big[ \distf(\bar {\mb U}_{\cdot, \mc S_k}, \bar {\mb U}^{(i)}_{\cdot, \mc S_k}) \vee \distf(\bar{\mb V}_{\cdot, \mc S_k}, \bar{\mb V}^{(i)}_{\cdot, \mc S_k}) \big] .
    \end{align}
    This allows us to rearrange the terms regarding $\distf(\bar{\mb U}^{(i)}_{\cdot, \mc S_k} ,\bar{\mb U}_{\cdot, \mc S_k}) $ to arrive at 
    \begin{align}
        & \distf(\bar{\mb U}_{\cdot, \mc S_k}, \bar{\mb U}^{(i)}_{\cdot, \mc S_k}) \lesssim \frac{\bar\sigma r^* \sqrt{\pi r^* \log d }}{\pi \min_{i \in \mc S_k} \bar \sigma_i^*} + \frac{ r^* \big( B \bignorm{\bar{\mb V}_{\cdot, \mc S_k}}\ti \log d + \bar\sigma \sqrt{\pi d}\bignorm{\bar{\mb U}_{\cdot, \mc S_k}}\ti \big) }{\pi \min_{i \in \mc S_k} \bar \sigma_i^*} 
        \label{eq: leave-one-out error after rearrangement}
    \end{align}
    simultaneously for all $i\in[n]$ and $k\in[\bar k]$ with probability at least $1-O(d^{-c})$, with the help of the condition 
    $$
     B \sqrt{\log d} \ll \bar\sigma \sqrt{\pi d} \ll \pi \min_{i\in\mc S_k} \bar \sigma_i^* / r^*.
    $$

        \paragraph{Singular subspace perturbation analysis.} In view of \eqref{eq: leave-one-out error after rearrangement}, it requires a row-wise control of $\bar{\mb U}_{\cdot, \mc S_k}$ as well as $\bar{\mb V}_{\cdot, \mc S_k}$. 
        
        To this end, we first evaluate its population counterpart $\norm{\bar{\mb U}^*_{\cdot, 1:r^*}}\ti$. Notice the relation 
        \begin{align}
            & \norm{\bar {\mb U}^*_{\cdot,1:r^*}}\ti = \norm{\psi(\mb 1_n \bo \zeta^{*\top} + \mb X^*\mb Y^{*\top})\bar{\mb V}_{\cdot,1:r^*}^{*}(\bar{\bo \Sigma}^*_{1:r^*, 1:r^*})^{-1}}\ti \\ 
            \lesssim & \bignorm{\psi(\mb 1 {\bo \zeta^*}\t + \mb X^*{\mb Y^*}\t) - \psi(\mb 1 {\bo \zeta^*}\t) }\ti / \lambda^{\sp} + \norm{\psi(\bo \zeta^*)\t \bar{\mb V}_{\cdot,1:r^*}^{*} (\bar{\bo \Sigma}^*_{1:r^*, 1:r^*})^{-1}}_2  \\ 
            = & \bignorm{\psi(\mb 1 {\bo \zeta^*}\t + \mb X^*{\mb Y^*}\t) - \psi(\mb 1 {\bo \zeta^*}\t) }\ti / \lambda^{\sp} + \norm{\psi(\mb 1_n\bo \zeta^{*\top}) \bar{\mb V}_{\cdot,1:r^*}^{*} (\bar{\bo \Sigma}^*_{1:r^*, 1:r^*})^{-1}}\fb / \sqrt{n} \\ 
            \leq & 2\pi \bignorm{\psi(\mb 1 {\bo \zeta^*}\t + \mb X^*{\mb Y^*}\t) - \psi(\mb 1 {\bo \zeta^*}\t) }\ti / \lambda^{\sp} + \norm{\bar{\mb U}_{\cdot,1:r^*}^{*}}\fb / \sqrt{n}\\ 
            \lesssim & \frac{\sigma^2\pi \sigma_1^* \sqrt{\mu r / n}}{\lambda^{\sp}} + \sqrt{\frac{r^*}{n}} \lesssim  \frac{\sigma^2\pi \sigma_1^* \sqrt{\mu r / n}}{\lambda^{\sp}}.
            \label{eq: upper bound on bar U}
        \end{align}
        by the choice of $\lambda^{\sp}$. 
        We next control $\bar{\mb V}^*_{\cdot,1:r^*}$ without requiring the
        entries of $\psi(\bo\zeta^*)$ to be bounded away from zero. Set
        \begin{equation*}
            \mb a^*\coloneqq\psi(\bo\zeta^*),\qquad
            a_0\coloneqq\psi(\zeta_0),\qquad
            \delta_\zeta^*\coloneqq
            \norm{\bo\zeta^*-\zeta_0\mb 1_p}_\infty,
        \end{equation*}
        {\emergencystretch=3em
        and define
        $\mb D_\psi^*\coloneqq\mb R^*-\mb 1_n\mb a^{*\top}$.
        The mean-value theorem, factor incoherence, and
        $\|({\bar{\bo\Sigma}^*_{1:r^*,1:r^*}})^{-1}\|\allowbreak\leq
        2\pi/\lambda^{\sp}$ give\par}
        \begin{align}
            \bignorm{\mb D_\psi^{*\top}}\ti
            &\lesssim \sigma^2\sigma_1^*\sqrt{\frac{\mu r}{p}},
            \nonumber\\
            \Bignorm{\mb D_\psi^{*\top}\bar{\mb U}^*_{\cdot,1:r^*}
            ({\bar{\bo\Sigma}^*_{1:r^*,1:r^*}})^{-1}}\ti
            &\lesssim
            \frac{\sigma^2\pi\sigma_1^*\sqrt{\mu r/p}}
            {\lambda^{\sp}}
            \eqqcolon\gamma_{\mb V}.
            \label{eq: varying component of bar V}
        \end{align}
        Write the remaining intercept component as
        \begin{equation*}
            \bo\Gamma_{\mb V}^*
            \coloneqq
            \mb a^*\mb 1_n\t\bar{\mb U}^*_{\cdot,1:r^*}
            ({\bar{\bo\Sigma}^*_{1:r^*,1:r^*}})^{-1},
        \end{equation*}
        so that
        $\bar{\mb V}^*_{\cdot,1:r^*}
        =\mb D_\psi^{*\top}\bar{\mb U}^*_{\cdot,1:r^*}
        ({\bar{\bo\Sigma}^*_{1:r^*,1:r^*}})^{-1}
        +\bo\Gamma_{\mb V}^*$.

        We distinguish two cases. If
        $|a_0|\geq2\sigma^2\delta_\zeta^*$ and $a_0=0$, then this
        condition forces $\delta_\zeta^*=0$; the mean-value bound gives
        $\mb a^*=\mb 0$, and hence $\bo\Gamma_{\mb V}^*=\mb 0$. If instead
        $|a_0|\geq2\sigma^2\delta_\zeta^*$ and $a_0\neq0$, then
        $|a_j^*-a_0|\leq\sigma^2\delta_\zeta^*$ for every $j\in[p]$,
        and hence
        $\max_j|a_j^*|/\|\mb a^*\|_2\leq3/\sqrt p$. Therefore,
        \begin{align*}
            \bignorm{\bo\Gamma_{\mb V}^*}\ti
            &\leq \frac{3}{\sqrt p}\bignorm{\bo\Gamma_{\mb V}^*}\fb\leq \frac{3}{\sqrt p}
            \Big(\sqrt{r^*}+\sqrt p\,\gamma_{\mb V}\Big)
            \lesssim\gamma_{\mb V}.
        \end{align*}
        The second inequality uses the preceding decomposition,
        $\|\bar{\mb V}^*_{\cdot,1:r^*}\|_{\mathrm F}=\sqrt{r^*}$, and
        \eqref{eq: varying component of bar V}. The last inequality follows
        from $r^*\lesssim\xi^2$, the definition of $\lambda^{\sp}$, and the
        signal-strength condition, which together imply
        $\sqrt{r^*/p}\lesssim\gamma_{\mb V}$.
        If instead $|a_0|<2\sigma^2\delta_\zeta^*$, then
        $\max_j|a_j^*|\leq3\sigma^2\delta_\zeta^*$, and thus
        \begin{align*}
            \bignorm{\bo\Gamma_{\mb V}^*}\ti
            &\leq
            \frac{2\pi\sqrt n}{\lambda^{\sp}}\max_{j\in[p]}|a_j^*|\lesssim
            \frac{\sigma^2\pi\sqrt{\mu r\sigma_r^*/p}}
            {\lambda^{\sp}}
            \lesssim\gamma_{\mb V},
        \end{align*}
        where the last two inequalities use
        $\delta_\zeta^*\leq\sqrt{\mu r\sigma_r^*/(np)}$ and the
        explicit lower bound $\sigma_r^*\gtrsim1$ in
        Assumption~\ref{assumption: signal strength and incoherence degree},
        which gives $\sqrt{\sigma_r^*}\lesssim\sigma_r^*\leq\sigma_1^*$.
        Combining
        the two cases with \eqref{eq: varying component of bar V} yields
        \begin{equation}
            \bignorm{\bar{\mb V}^*_{\cdot,1:r^*}}\ti
            \lesssim
            \frac{\sigma^2\pi\sigma_1^*\sqrt{\mu r/p}}
            {\lambda^{\sp}}.
            \label{eq: upper bound on bar V}
        \end{equation}

        Define $\mc S \coloneqq \cup_{l\in[\bar k]} \mc S_l$. It suffices to control $\norm{\bar{\mb U}_{\cdot, \mc S}}\ti$, which requires us to bound
        $\dist\ti(\bar{\mb U}_{\cdot,1:r^*}, \bar{\mb U}^*_{\cdot, 1:r^*})$. We leverage the general eigen-subspace perturbation theory developed in \cite{lei2019unified} together with the Hermitian dilation trick: 
        
        Let $\mb A \coloneqq \left[\begin{matrix}
            \mb 0 & \mc P_{\Omega}(\mb R) / \pi \\ 
            \mc P_{\Omega}(\mb R)\t /\pi & \mb 0
        \end{matrix} \right] \in \bb R^{(n+p) \times (n+p)}$ and $\mb A^* \coloneqq \left[\begin{matrix} 
        \mb 0 & \mb R^{*,\mc S} \\ 
        \mb R^{*,\mc S\top} & \mb 0 \end{matrix}
        \right]$ where $\mb R^{*,\mc S} \coloneqq \bar{\mb U}^*_{\cdot, \mc S} \bar{\bo \Sigma}^*_{\mc S, \mc S} \bar{\mb V}^{*\top}_{\cdot, \mc S}$. Note that the SVD of $\mb A^*$ is written as $\tilde{\mb U}^* \tilde{\bo \Sigma}^* \tilde{\mb U}^{*\top}$ where $\tilde{\mb U}^* \coloneqq \frac{1}{\sqrt{2}} \left[
        \begin{matrix}
            \bar{\mb U}^*_{\cdot, \mc S} & \bar{\mb U}^*_{\cdot, \mc S} \\ 
            \bar{\mb V}^*_{\cdot, \mc S} & - \bar{\mb V}^*_{\cdot, \mc S}
        \end{matrix}
        \right] \in \mc O(|\mc S|)$ and $\tilde{\bo \Sigma}^* \coloneqq \left[ \begin{matrix}
            \bar{\bo \Sigma}^*_{\mc S, \mc S} & \mb 0 \\ 
            \mb 0 & -\bar{\bo \Sigma}^*_{\mc S, \mc S}
        \end{matrix}\right]$. Similarly, the top-$(|\mc S|)$ SVD (in magnitude) of $\mb A$ is written as $\tilde{\mb U}\tilde{\bo \Sigma} \tilde{\mb U}\t$. 
        It therefore remains to control $\dist\ti(\tilde{\mb U}, \tilde{\mb U}^*)$ using \citet[Theorem~2.3]{lei2019unified}. We first specify the quantities involved in this result:
        \begin{align}
            & \delta = c_\delta d^{-c-2},\quad L_1 (\delta) = c_{L_1}\bar\sigma \sqrt{d / \pi}, \quad L_2(\delta) = 0, \\
            & L_3(\delta) = c_{L_3} r^* \big[\bar\sigma \sqrt{\pi d r \log d} + B \log d \big] / \lambda^{\sp} , \\ 
            & \lambda_-(\delta) = E_+(\delta) = \bar E_+(\delta) = E_\infty(\delta) = L_1(\delta), \\
            & b_\infty(\delta) = c_{b_\infty} B  \log d / \pi, \quad b_2(\delta) = c_{b_2} \bar\sigma \sqrt{r \log d / \pi}, \\ 
            & \eta(\delta) = c_{\eta} \bar\sigma \sqrt{d / \pi}, \quad \sigma(\delta) = c_{\sigma} \bar\sigma \sqrt{d / \pi}. 
        \end{align}
        Here, we use the approximation error bound \eqref{eq: upper bound for F approximation error} and Lemma~\ref{lemma: noise spectral norm concentration}. The condition $\lambda^{\sp}/\pi \geq 4\{\sigma(\delta)+L_1(\delta)+\lambda_-(\delta)\}$ holds by the definition of $\lambda^{\sp}$. Therefore, applying \citet[Theorem~2.3]{lei2019unified} together with \eqref{eq: upper bound on bar U} and \eqref{eq: upper bound on bar V} yields that with probability at least $1-O(d^{-c-2})$
        \begin{align}
            & \dist\ti\big(\tilde{\mb U}, \mb A \tilde{\mb U}^* (\tilde{\bo \Lambda}^*)^{-1} \big) \\
            \lesssim & \frac{\pi r^* }{\lambda^{\sp}} \Big\{ \sigma(\delta) \Big[ \norm{\tilde{\mb U}^*}\ti + \frac{\pi \bignorm{(\mb A - \mb A^*) \tilde{\mb U}^*}\ti}{\lambda^{\sp}} \Big] + \frac{ \pi E_+(\delta) b_2(\delta)}{\lambda^{\sp}} + r^* \bar E_+(\delta)\bignorm{\tilde{\mb U}^*}\ti \Big\} \\
            \lesssim & \frac{\pi r^* }{\lambda^{\sp}} \Big\{ \bar\sigma \sqrt{d / \pi} \Big[\frac{\sigma^2\pi \sigma_1^* \sqrt{\mu r / (n\wedge p)}}{\lambda^{\sp}}  + \frac{\pi \bar\sigma \sqrt{r \log d/\pi}}{\lambda^{\sp}}\Big]
             + \frac{\bar\sigma \sqrt{d } \bar\sigma \sqrt{r \log d}}{\lambda^{\sp}} \\ 
             & + \frac{r^* \pi \bar\sigma \sqrt{d / \pi}}{\lambda^{\sp}}\frac{\sigma^2\pi \sigma_1^* \sqrt{\mu r / (n\wedge p)}}{\lambda^{\sp}}\Big\}\\
            \lesssim & \frac{\sigma^2\pi \sigma_1^* \sqrt{\mu r / (n\wedge p)}}{\lambda^{\sp}} .
        \end{align}

        On the other hand, we have from the Bernstein inequality that with probability at least $1-O(d^{-c-2})$ 
        \begin{align}
            & \dist\ti(\mb A \tilde{\mb U}^* (\tilde{\bo \Lambda}^*)^{-1}, \tilde{\mb U}^*) \lesssim \frac{ r^* \bar\sigma \sqrt{\pi r \log d}}{\lambda^{\sp}} + \frac{B \log d}{\lambda^{\sp}} \norm{\tilde{\mb U}^*}\ti \lesssim \frac{\sigma^2 \pi \sigma_1^* \sqrt{\mu r / (n\wedge p)}}{\lambda^{\sp}} .
        \end{align}

        Putting these together implies that with probability at least $1-O(d^{-c-2})$ 
        \begin{align}
            & \max_{k \in[\bar k]}\big\{ \norm{\bar{\mb U}_{\cdot, \mc S_k}}\ti \vee \norm{\bar{\mb V}_{\cdot, \mc S_k}}\ti \big\} \lesssim \norm{\bar{\mb U}^*_{\cdot, 1: r^*}}\ti \vee  \norm{\bar{\mb V}^*_{\cdot, 1: r^*}}\ti +   \dist\ti\big(\tilde{\mb U}, \tilde{\mb U}^*)\\ 
            \leq & \norm{\bar{\mb U}^*_{\cdot, 1: r^*}}\ti \vee  \norm{\bar{\mb V}^*_{\cdot, 1: r^*}}\ti + \dist\ti\big(\tilde{\mb U}, \tilde{\mb U}^*)\\ 
            \leq & \norm{\bar{\mb U}^*_{\cdot, 1: r^*}}\ti \vee  \norm{\bar{\mb V}^*_{\cdot, 1: r^*}}\ti + \dist\ti\big(\tilde{\mb U}, \mb A \tilde{\mb U}^* (\tilde{\bo \Lambda}^*)^{-1} \big)  + \dist\ti(\mb A  \tilde{\mb U}^* (\tilde{\bo \Lambda}^*)^{-1}, \tilde{\mb U}^*) \\
            \lesssim & \frac{\sigma^2\pi \sigma_1^* \sqrt{\mu r / (n\wedge p)}}{\lambda^{\sp}}.
            \label{eq: upper bound on U_r}
        \end{align}

    To proceed, we substitute \eqref{eq: upper bound on U_r} into \eqref{eq: leave-one-out error after rearrangement} and obtain, simultaneously for all $i\in[n]$ and $k\in[\bar k]$, with probability at least $1-O(d^{-c})$ that 
    \begin{align}
        & \distf(\bar{\mb U}_{\cdot, \mc S_k}, \bar{\mb U}^{(i)}_{\cdot, \mc S_k}) \vee \distf(\bar{\mb V}_{\cdot, \mc S_k}, \bar{\mb V}^{(i)}_{\cdot, \mc S_k}) \\ 
        \lesssim & \frac{\bar\sigma r^* \sqrt{\pi r^* \log d }}{\pi \min_{i \in \mc S_k} \bar \sigma_i^*} + \frac{r^*(B  \log d + \bar\sigma \sqrt{\pi d} ) }{\pi \min_{i \in \mc S_k} \bar \sigma_i^*} \frac{\sigma^2\pi \sigma_1^* \sqrt{\mu r / (n\wedge p)}}{\lambda^{\sp}} \\
        \lesssim &   \frac{\sigma^2 r^* \bar\sigma  \sigma_1^* \sqrt{ \mu r \rho \pi }}{\min_{i \in \mc S_k} \bar{\sigma}_i^* \lambda^{\sp}}   \eqqcolon \bar \zeta^{\sf{loo}}_k, \qquad \forall k \in [\bar k] .
        \label{eq: leave-one-out error for bar U}
        \end{align}

    The same fixed-index bound and union argument give the leave-one-column-out control below with probability at least $1-O(d^{-c})$:
    \begin{align}
        & \max_{k \in [\bar k]} \max_{j \in[p]}  \Big\{\distf(\bar{\mb U}_{\cdot, \mc S_k},\bar{\mb U}^{(-j)}_{\cdot, \mc S_k})\vee \distf(\bar{\mb V}_{\cdot, \mc S_k},\bar{\mb V}^{(-j)}_{\cdot, \mc S_k}) \Big\}  \lesssim \bar \zeta^{\sf{loo}}_k . 
    \end{align}

    \paragraph{Controlling Matrix Estimation Errors.}
    Finally, with the above singular subspace analysis in place, we are able to establish the guarantee for the leave-one-out coupling. 
    We first distinguish scalar soft thresholding from its matrix-valued
    counterpart.  For a matrix $\mb A=\mb P\diag(\sigma_l(\mb A))\mb Q\t$,
    define
    \begin{equation}
        \operatorname{SVT}_{s}(\mb A)
        \coloneqq
        \mb P\diag\big((\sigma_l(\mb A)-s)_+\big)\mb Q\t.
        \label{eq: matrix soft singular value thresholding}
    \end{equation}
    This map is the proximal operator of $s\norm{\cdot}_*$ and is therefore
    nonexpansive in the Frobenius norm. Recall that an operator $\mathcal T$ is nonexpansive under the
Frobenius norm if
$$
    \norm{\mathcal T(\mb A)-\mathcal T(\mb B)}_{\mathrm F}
    \leq \norm{\mb A-\mb B}_{\mathrm F}
$$
for all conformable matrices $\mb A$ and $\mb B$; see also \citet[Sections~2.3 and~6.7.3]{parikh2014proximal}.
Importantly, this is a property of the full matrix-valued SVT operator. It is also bi-orthogonally equivariant:
    for arbitrary orthogonal matrices $\mb Q_1$ and $\mb Q_2$,
    \begin{equation}
        \operatorname{SVT}_{s}(\mb Q_1\mb A\mb Q_2\t)
        =\mb Q_1\operatorname{SVT}_{s}(\mb A)\mb Q_2\t.
        \label{eq: equivariance of matrix SVT}
    \end{equation}
    In particular, when $\mb A$ is a nonnegative diagonal matrix,
    $\operatorname{SVT}_{s}(\mb A)=\mc T_s(\mb A)$. 
    
    On $\mc E^{\sp}$, zeroing the $i$th row cannot increase the spectral
    norm, and hence
    $\norm{\mb R^{(i)}-\pi\mb R^*}
    \leq\norm{\mc P_\Omega(\mb R)-\pi\mb R^*}
    \leq\lambda^{\sp}/2$. Weyl's inequality and the definition of $r^*$
    therefore give $\bar\sigma_s\vee\bar\sigma_s^{(i)}<\lambda^{\sp}$
    for every $s>r^*$. Since $[r^*]\subseteq
    \bigcup_{k\in[\bar k]}\mc S_k$, all thresholded components outside
    these blocks vanish for both matrices. The following telescoping identity
    is therefore exact:
    \begin{align}
        &\pi  \hat{\mb R}^{\sp,(i)} - \pi \hat{\mb R}^{\sp} \\
        = & \sum_{k \in [\bar k]}\Big( \bar{\mb U}^{(i)}_{\cdot, \mc S_k} \operatorname{SVT}_{\lambda^{\sp}}(\bar{\bo \Sigma}^{(i)}_{\mc S_k}) {\bar{\mb V}^{(i)\top}_{\cdot, \mc S_k} }  - \bar{\mb U}_{\cdot, \mc S_k} \operatorname{SVT}_{\lambda^{\sp}}(\bar{\bo \Sigma}_{\mc S_k}) {\bar{\mb V}_{\cdot, \mc S_k} }\t\Big) \\ 
        = & \sum_{k \in [\bar k]}\Big[ \big(\bar{\mb U}^{(i)}_{\cdot, \mc S_k}  - \bar{\mb U}_{\cdot, \mc S_k} \mb O_{\mb U,k}^{(i)} \big)\operatorname{SVT}_{\lambda^{\sp}}(\bar{\bo \Sigma}^{(i)}_{\mc S_k}) {\bar{\mb V}^{(i)\top}_{\cdot, \mc S_k} } \\
        & \qquad +  \bar{\mb U}_{\cdot, \mc S_k}\mb O_{\mb U,k}^{(i)}  \operatorname{SVT}_{\lambda^{\sp}}(\bar{\bo \Sigma}_{\mc S_k}^{(i)} ) \big(\bar{\mb V}^{(i)}_{\cdot, \mc S_k} - \bar{\mb V}_{\cdot, \mc S_k} {\mb O_{\mb V,k}^{(i)}}\big)\t \\ 
        & \qquad +\bar{\mb U}_{\cdot,\mc S_k}
\big(
    \mb O_{\mb U,k}^{(i)}
    \operatorname{SVT}_{\lambda^{\sp}}
    (\bar{\bo\Sigma}_{\mc S_k}^{(i)})
    \mb O_{\mb V,k}^{(i)\top}
    -
    \operatorname{SVT}_{\lambda^{\sp}}
    (\bar{\bo\Sigma}_{\mc S_k})
\big)
\bar{\mb V}_{\cdot,\mc S_k}^{\top} \Big] \\ 
        = & \sum_{k \in [\bar k]}\Big[ \big(\bar{\mb U}^{(i)}_{\cdot, \mc S_k}  - \bar{\mb U}_{\cdot, \mc S_k} \mb O_{\mb U,k}^{(i)} \big)\operatorname{SVT}_{\lambda^{\sp}}(\bo \Sigma^{(i)}_{\mc S_k}) {\bar{\mb V}^{(i)\top}_{\cdot, \mc S_k} } \\
        & \qquad +  \bar{\mb U}_{\cdot, \mc S_k} \mb O_{\mb U,k}^{(i)} \operatorname{SVT}_{\lambda^{\sp}}(\bar{\bo \Sigma}_{\mc S_k}^{(i)} ) \big(\bar{\mb V}^{(i)}_{\cdot, \mc S_k} - \bar{\mb V}_{\cdot, \mc S_k} {\mb O_{\mb V,k}^{(i)}}\big)\t \\ 
        & \qquad + \bar{\mb U}_{\cdot, \mc S_k} \big(\operatorname{SVT}_{\lambda^{\sp}}(  \mb O^{(i)}_{\mb U,k} \bar{\bo \Sigma}_{\mc S_k}^{(i)} {\mb O^{(i)\top}_{\mb V,k}}) -  \operatorname{SVT}_{\lambda^{\sp}}( \bar{\bo \Sigma}_{\mc S_k}  )  \big)\bar{\mb V}_{\cdot, \mc S_k}\t \Big],
        \label{eq: leave-one-out perturbation of USVT matrix estimation}
    \end{align}
    where we introduce the rotation matrices
     $\mb O^{(i)}_{\mb U,k}$ and $\mb O_{\mb V,k}^{(i)}$ as follows: 
    \begin{align}
        & \mb O^{(i)}_{\mb U,k} \coloneqq \arg\min_{\mb O \in \mc O(|\mc S_k|)} \norm{\bar{\mb U}^{(i)}_{\cdot, \mc S_k} - \bar{\mb U}_{\cdot, \mc S_k} \mb O}\fb, \\
        & \mb O^{(i)}_{\mb V,k} \coloneqq \arg\min_{\mb O \in \mc O(|\mc S_k|)} \norm{\bar{\mb V}^{(i)}_{\cdot, \mc S_k} - \bar{\mb V}_{\cdot, \mc S_k} \mb O}\fb. 
    \end{align}

    By \eqref{eq: equivariance of matrix SVT} and the nonexpansiveness of
    $\operatorname{SVT}_{\lambda^{\sp}}$, we have 
    \eq{
    \Bignorm{\operatorname{SVT}_{\lambda^{\sp}}(  \mb O^{(i)}_{\mb U,k} \bar{\bo \Sigma}_{\mc S_k}^{(i)} {\mb O^{(i)\top}_{\mb V,k}}) -  \operatorname{SVT}_{\lambda^{\sp}}( \bar{\bo \Sigma}_{\mc S_k}) }\fb \leq  \Bignorm{\mb O^{(i)}_{\mb U,k} \bar{\bo \Sigma}_{\mc S_k}^{(i)} {\mb O^{(i)\top}_{\mb V,k}} -  \bar{\bo \Sigma}_{\mc S_k}}\fb.
    }
    Consequently, the triangle inequality yields that 
    \longeq{
    & \Bignorm{\mb O^{(i)}_{\mb U,k} \bar{\bo \Sigma}_{\mc S_k}^{(i)} {\mb O^{(i)\top}_{\mb V,k}} -  \bar{\bo \Sigma}_{\mc S_k}}\fb  =  \Bignorm{\mb O^{(i)}_{\mb U,k} \bar{\bo \Sigma}_{\mc S_k}^{(i)}  -  \bar{\bo \Sigma}_{\mc S_k} {\mb O^{(i)}_{\mb V,k}}}\fb\\ 
    \leq & \Bignorm{\bar{\mb U}_k\t {\bar{\mb U}^{(i)}_k} \bar{\bo \Sigma}_{\mc S_k}^{(i)} -  \bar{\bo \Sigma}_{\mc S_k}\bar{\mb V}_k\t \bar{\mb V}_k^{(i)}}\fb \\
    &\quad + \Bignorm{{\bar{\mb U}^{\top}_k} \bar{\mb U}^{(i)}_k - \mb O_{\mb U,k}^{(i)}}\fb \norm{\bar{\bo \Sigma}_{\mc S_k}^{(i)}} + \Bignorm{{\bar{\mb V}^{\top}_k} \bar{\mb V}^{(i)}_k - \mb O_{\mb V,k}^{(i)}}\fb \norm{\bar{\bo \Sigma}_{\mc S_k}}  \\
    \leq & \Bignorm{ {\bar{\mb U}^{\top}_k} (\mc P_\Omega(\mb R) - \mb R^{(i)}) \bar{\mb V}_k^{(i)}}\fb \\
    &\quad + \Bignorm{{\bar{\mb U}^{\top}_k} \bar{\mb U}^{(i)}_k - \mb O_{\mb U,k}^{(i)}}\fb \norm{\bar{\bo \Sigma}_{\mc S_k}^{(i)}} + \Bignorm{{\bar{\mb V}^{\top}_k} \bar{\mb V}^{(i)}_k - \mb O_{\mb V,k}^{(i)}}\fb \norm{\bar{\bo \Sigma}_{\mc S_k}} \\
    \lesssim  & \bignorm{\bar{\mb E}_{i,\cdot} \bar{\mb V}^{(i)}_{\cdot, \mc S_k}} + (\bar \zeta^{\sf{loo}}_{k})^2 \pi \max_{i\in \mc S_k} \bar \sigma_i^*  \\ 
    \lesssim& \bar \zeta_k^{\sf{loo}}\pi  \max_{i\in \mc S_k} \bar \sigma_i^*
    \label{eq: leave-one-out spectral core}
    }
    with probability at least $1- O(d^{-c})$, where we used \eqref{eq: bernstein for wedin + leave-one-out}, \eqref{eq: upper bound on U_r}, and \eqref{eq: leave-one-out error for bar U} to derive the intermediate upper bound: 
    \begin{align}
        & \bignorm{\bar{\mb E}_{i,\cdot} \bar{\mb V}^{(i)}_{\cdot, \mc S_k}} \lesssim \bar\sigma \sqrt{\pi r^* \log d} + B\norm{\bar{\mb V}_{\cdot,\mc S_k}^{(i)}}\ti \log d\\ 
        \lesssim & \bar\sigma \sqrt{\pi r^* \log d} + B  \frac{\sigma^2\pi \sigma_1^* \sqrt{\mu r / (n\wedge p)}}{\lambda^{\sp}} \log d
        \lesssim  \bar \zeta_k^{\sf{loo}} \pi \max_{i\in \mc S_k} \bar \sigma_i^*. 
    \end{align}

Combining \eqref{eq: leave-one-out spectral core} with \eqref{eq: leave-one-out error for bar U} and \eqref{eq: leave-one-out perturbation of USVT matrix estimation}, we bound the difference $\hat{\mb R}^{\sp} - \hat{\mb R}^{\sp,(i)}$ in terms of the Frobenius norm as follows: 
\begin{align}
    \max_{i\in[n]}\norm{\hat{\mb R}^{\sp,(i)} - \hat{\mb R}^{\sp} }\fb 
    \lesssim{}\;&
    \sum_{k \in [\bar k]} \bar \zeta^{\sf{loo}}_{k} \max_{i \in \mc S_k} \bar \sigma_i^* \lesssim  \frac{\sigma^2( r^*)^2 \bar\sigma  \sigma_1^* \sqrt{ \mu r \rho \pi }}{\lambda^{\sp}},
\end{align}
with probability at least $1 - O(d^{-c})$, which concludes the proof for the leave-one-row-out error. The proof for the leave-one-column-out error similarly follows, which is omitted here.

\section{Unilateral Refinement Analysis}
\label{sec: unilateral refinement analysis}

This section studies the unilateral refinement procedure (Algorithm~\ref{alg: OS}) initialized by the spectral estimator (Algorithm~\ref{alg: USVT}). The main message is that a single refinement step is enough to strengthen a spectral/Frobenius-norm consistent estimator into a row-wise consistent one. Moreover, extending the leave-one-out analysis developed for the spectral estimator, we establish a corresponding leave-one-out analysis for the refinement step, which in turn prepares the ground for the decoupling analysis of the final gradient descent stage.

\label{subsec: proof of OS}
\subsection{A General Version of Theorem~\ref{thm: OS}}

We begin with a general formulation of the unilateral refinement algorithm by imposing conditions on the initial estimators and their leave-one-out counterparts, without yet specifying their explicit forms. 

\begin{theorem}\label{thm: general OS}
Consider the model and the assumptions in Theorem~\ref{thm: USVT}.  Write
\begin{equation*}
    \mb F^0\coloneqq
    \begin{pmatrix}\mb X^0\\ \mb Y^0\end{pmatrix},
    \qquad
    \mb F^*\coloneqq
    \begin{pmatrix}\mb X^*\\ \mb Y^*\end{pmatrix},
\end{equation*}
and define the common population alignment
\begin{equation}
    \mb O^0
    \coloneqq
    \widehat{\mb H}\big((\mb F^0)\t\mb F^*\big)
    \in\argmin_{\mb O\in\mc O(r)}
    \bignorm{\mb F^0\mb O-\mb F^*}\fb.
    \label{eq: common population alignment for UR}
\end{equation}
Suppose that the initialization $[{\bo \zeta^0}, \mb X^0, \mb Y^0]$ satisfies
\begin{align}
    & \underline c_\psi \max\Big\{ \sqrt{n \pi}\bignorm{\bo \zeta^0 - \bo \zeta^*}_2,
    \sqrt{\pi \sigma_r^*}\bignorm{\mb X^0\mb O^0-\mb X^*}\fb,
    \sqrt{\pi \sigma_r^*}\bignorm{\mb Y^0\mb O^0-\mb Y^*}\fb \Big\} \leq  \xi^0
\end{align}
with $\frac{\xi^0}{\underline c_\psi \sqrt{\pi \sigma_r^*}}\ll \sqrt{\sigma_r^*}$, and suppose that
\begin{equation*}
    \bignorm{\mb X^0}\ti
    \lesssim\sqrt{\frac{\mu r\sigma_1^*}{n}},
    \qquad
    \sqrt{\frac{n}{\sigma_r^*}}
    \bignorm{\bo\zeta^0-\zeta_0\mb 1_p}_\infty
    \vee\bignorm{\mb Y^0}\ti
    \lesssim\sqrt{\frac{\mu r\sigma_1^*}{p}}.
\end{equation*}
{\emergencystretch=3em
Moreover, there exist leave-one-out versions $[{\bo \zeta^{0, (i)}},\allowbreak  \mb X^{0, (i)},\allowbreak  \mb Y^{0, (i)}]$
 and $[{\bo \zeta^{0, (-j)}},\allowbreak  \mb X^{0, (-j)},\allowbreak  \mb Y^{0, (-j)}]$ which are functions of the observation such that 
(i) $[{\bo \zeta^{0, (i)}},\allowbreak  \mb X^{0, (i)},\allowbreak  \mb Y^{0, (i)}]$ (resp. $[{\bo \zeta^{0, (-j)}},\allowbreak  \mb X^{0, (-j)},\allowbreak  \mb Y^{0, (-j)}]$) is independent of $\mb E_{i,\cdot}$ and $\bo \Omega_{i,\cdot}$ (resp. $\mb E_{\cdot, j}$ and $\bo \Omega_{\cdot, j}$) for every $i \in[n]$ (resp. $j \in[p]$); (ii) the following conditions hold for every $i \in[n]$ and $j \in[p]$:
letting\par}

\begin{equation*}
\begin{aligned}
    & \mb F^{0,(l)}\coloneqq
    \begin{pmatrix}\mb X^{0,(l)}\\ \mb Y^{0,(l)}\end{pmatrix},
    \qquad
    \mb O^{0,(l)}
    \coloneqq
    \widehat{\mb H}\big((\mb F^{0,(l)})\t\mb F^0\big)
    \in\argmin_{\mb O\in\mc O(r)}
    \bignorm{\mb F^{0,(l)}\mb O-\mb F^0}\fb,\\
    &\quad \qquad l\in[n]\cup[-p],
\end{aligned}
\end{equation*}
one has 
\begin{align*}
    & \underline c_\psi\max\Big\{\sqrt{n \pi}\bignorm{\bo \zeta^{0, (i)} - \bo \zeta^0}_2,
    \sqrt{\pi \sigma_r^*}\bignorm{\mb X^{0,(i)}\mb O^{0,(i)}-\mb X^0}\fb,\\
    &\quad \sqrt{\pi \sigma_r^*}\bignorm{\mb Y^{0,(i)}\mb O^{0,(i)}-\mb Y^0}\fb\Big\}
    \lesssim \xi^{0,\sf{loo}}, \\ 
    & \underline c_\psi\max\Big\{\sqrt{n \pi}\bignorm{\bo \zeta^{0, (-j)} - \bo \zeta^0}_2,
    \sqrt{\pi \sigma_r^*}\bignorm{\mb X^{0,(-j)}\mb O^{0,(-j)}-\mb X^0}\fb,\\
    &\quad \sqrt{\pi \sigma_r^*}\bignorm{\mb Y^{0,(-j)}\mb O^{0,(-j)}-\mb Y^0}\fb\Big\}
    \lesssim \xi^{0,\sf{loo}},
\end{align*}
Assume that the same rowwise bounds hold with the initializer replaced by each
leave-one-out version, and that
$\frac{\sqrt{\mu r}\kappa\xi^{0,\sf{loo}}}
{\underline c_\psi \sqrt{\pi \sigma_r^*}}
\ll \sqrt{\frac{\mu r\sigma_1^*}{d}}$.
Define
\begin{align}
    \xi^{1}\coloneqq{}& \bar\sigma \sqrt{\kappa r \log d}
    +\kappa_{\psi}\kappa \sqrt{\mu r}\big[\xi^0 / \sqrt{n\wedge p} + \kappa \xi^{0,\sf{loo}} \big]  + \frac{\sigma^2\kappa^{\frac32}\sigma_r^*}{\sqrt\pi}
    \sqrt{\frac{\mu^3 r^3}{n p (n\wedge p)}}\log d.
    \label{eq: general unilateral refinement rate}
\end{align}
Assume in addition that
\begin{align}
    & \kappa_\psi^2\kappa^2\mu r\log d
    \ll \pi(n\wedge p), \qquad  \frac{\xi^1}{\underline c_\psi\sqrt{\pi\sigma_r^*}}
    \ll \sqrt{\frac{\mu r\sigma_1^*}{d}}, \nonumber\\
    & \kappa_\psi\frac{\kappa\xi^{0,\sf{loo}}}
    {\underline c_\psi\sqrt{\pi\sigma_r^*}}
    \Big[\sqrt{\frac{\kappa}{\pi\sigma_r^*}}
    +\kappa\sqrt{\frac{\mu r\sigma_1^*}{n\wedge p}}\Big]
    \ll 1, \qquad  \kappa_\psi^2\kappa^3\mu r\sqrt{\rho}\,\xi^1
    \ll \underline c_\psi\sqrt{\pi(n\wedge p)}.
    \label{eq: general unilateral refinement smallness conditions}
\end{align}
Then, the unilateral refinement algorithm (Algorithm~\ref{alg: OS}) satisfies the following guarantees with probability at least $1 - O(d^{-c})$:  for every $i\in[n]$ and $j\in[p]$, 
\begin{align}
    & \underline c_\psi \max\Big\{ \sqrt{n\pi}\bignorm{\hat{\bo \zeta}^{\ur} - \bo\zeta^*}_\infty,
    \sqrt{\pi \sigma_r^*}\bignorm{\hat{\mb X}^{\ur}\mb O^0-\mb X^*}\ti,\\
    &\qquad \sqrt{\pi \sigma_r^*}\bignorm{\hat{\mb Y}^{\ur}\mb O^0-\mb Y^*}\ti \Big\}
    \lesssim \xi^1,
    \label{eq: row-wise error for unilateral refinement}
    \\
    & 
    \begin{aligned}
    & \underline c_\psi\max\Big\{\sqrt{n \pi}\bignorm{\hat{\bo \zeta}^{\ur, (i)} - \hat{\bo \zeta}^{\ur}}_2,
    \sqrt{\pi \sigma_r^*}\bignorm{\hat{\mb X}^{\ur,(i)}\mb O^{0,(i)}-\hat{\mb X}^{\ur}}\fb,\\
    &\qquad \sqrt{\pi \sigma_r^*}\bignorm{\hat{\mb Y}^{\ur,(i)}\mb O^{0,(i)}-\hat{\mb Y}^{\ur}}\fb\Big\} \lesssim \kappa_\psi \kappa^{\frac32} \sqrt{\mu r\rho}\,\xi^1, \\
    & \underline c_\psi\max\Big\{\sqrt{n \pi}\bignorm{\hat{\bo \zeta}^{\ur, (-j)} - \hat{\bo \zeta}^{\ur}}_2,
    \sqrt{\pi \sigma_r^*}\bignorm{\hat{\mb X}^{\ur,(-j)}\mb O^{0,(-j)}-\hat{\mb X}^{\ur}}\fb,\\
    &\qquad \sqrt{\pi \sigma_r^*}\bignorm{\hat{\mb Y}^{\ur,(-j)}\mb O^{0,(-j)}-\hat{\mb Y}^{\ur}}\fb\Big\}\lesssim \kappa_\psi \kappa^{\frac32} \sqrt{\mu r \rho }\xi^1. 
    \end{aligned}
    \label{eq: leave-one-out error for unilateral refinement}
\end{align}
\end{theorem}

For the population-level analysis, it is useful to introduce the auxiliary
rotation
\begin{equation*}
    \widetilde{\mb O}^{0,(l)}
    \coloneqq 
    \widehat{\mb H}\big((\mb F^{0,(l)})\t\mb F^*\big),
    \qquad l\in[n]\cup[-p].
\end{equation*}
where $\hat{\mb H}(\mb A) = \mb A(\mb A\t \mb A)^{-\frac12}$. 
Unlike the relative rotation $\mb O^{0,(l)}$ appearing in the theorem,
$\widetilde{\mb O}^{0,(l)}$ depends only on the $l$th leave-one-out
initializer and the deterministic population factor. The two alignments are
linked by the stability of the Procrustes map. Indeed, setting
$\bar{\mb F}^{0,(l)}\coloneqq\mb F^{0,(l)}\mb O^{0,(l)}$, the orthogonal
equivariance of the polar factor gives
\begin{equation*}
    \widehat{\mb H}\big((\mb F^0\mb O^0)\t\mb F^*\big)=\mb I_r,
    \qquad
    \widehat{\mb H}\big((\bar{\mb F}^{0,(l)}\mb O^0)\t\mb F^*\big)
    ={\mb O^0}\t{\mb O^{0,(l)}}\t\widetilde{\mb O}^{0,(l)}.
\end{equation*}
Notice that, by the assumptions on $\xi^0$ and $\xi^{0,\sf{loo}}$, 
\begin{equation*}
    \bignorm{\mb F^0\mb O^0-\mb F^*}\norm{\mb F^*}
    \leq\frac{1}{2}\sigma_r^2(\mb F^*),
    \qquad
    \bignorm{\mb F^{0,(l)}\mb O^{0,(l)}-\mb F^0}\norm{\mb F^*}
    \leq\frac{1}{4}\sigma_r^2(\mb F^*).
\end{equation*}
Applying Lemma~\ref{lemma: perturbation theory of optimal rotation under F norm}
to $\mb F^0\mb O^0$ and $\bar{\mb F}^{0,(l)}\mb O^0$ then gives
\begin{equation}
    \bignorm{\mb F^{0,(l)}\widetilde{\mb O}^{0,(l)}-\mb F^0\mb O^0}\fb
    \leq 5\kappa_{\mb F^*}
    \bignorm{\mb F^{0,(l)}\mb O^{0,(l)}-\mb F^0}\fb,
    \qquad
    \kappa_{\mb F^*}
    \coloneqq\frac{\norm{\mb F^*}^2}{\sigma_r^2(\mb F^*)}=\kappa.
    \label{eq: transfer from relative to population-aligned LOO rotations}
\end{equation}
Here the last identity follows from the balanced population parametrization.
Thus, for the spectral initialization, one may take $\mb O^0$ to be the
stacked-factor Procrustes representative of the common rotation in
Theorem~\ref{thm: USVT}. In particular, a relative leave-one-out error of
order $\xi^{0,\sf{loo}}$ yields a population-aligned error of order
$\kappa \xi^{0,\sf{loo}}$. 

\begin{remark}
    For an actual initialization constructed from the observations, the information-theoretic lower bound suggests that $\xi^0 / \sqrt{n\wedge p} \vee \xi^{0,\sf{loo}}$ should dominate $\sigma\sqrt{r\log d}$. We therefore focus on the terms involving $\xi^0$ and $\xi^{0,\sf{loo}}$, and interpret the error cost of the refinement as an additional multiplicative factor depending on $\kappa_{\psi}$, $\kappa$, $\mu$, and $r$, in exchange for achieving row-wise consistency in a uniformly spread manner. 
\end{remark}

\subsection{Proof of Theorem~\ref{thm: OS}}

By Lemma~\ref{lemma: scaling equivalence}, it suffices first to work in the
balanced parametrization. On the event in
Theorem~\ref{thm: detailed USVT}, the sharp global bound in
\eqref{eq: USVT zeta X and Y bounds} and the choice of $\lambda^{\sp}$ in
Theorem~\ref{thm: USVT} give
\begin{equation*}
    \xi^0
    \lesssim
    \kappa\sigma\sqrt{\bar\sigma\sigma_1^*}
    (\pi\mu r d)^{\frac14}\xi^{\frac52}.
\end{equation*}
Moreover, \eqref{eq: USVT zeta X Y loo error} and
\eqref{eq: upper bound on r*} yield
\begin{equation*}
    \xi^{0,\sf{loo}}
    \lesssim
    \kappa\sigma\sqrt{\bar\sigma\sigma_1^*}
    \Big(\frac{\pi\mu r\rho^2}{d}\Big)^{\frac14}
    \xi^{\frac52}.
\end{equation*}
Since $d/(n\wedge p)=\rho$, these two estimates have the common normalized
order
\begin{equation}
    \frac{\xi^0}{\sqrt{n\wedge p}}\vee\xi^{0,\sf{loo}}
    \lesssim
    \kappa\sigma\sqrt{\bar\sigma\sigma_1^*}
    \Big(\frac{\pi\mu r\rho^2}{d}\Big)^{\frac14}
    \xi^{\frac52}.
    \label{eq: balanced spectral inputs for unilateral refinement}
\end{equation}

The centered entrywise projection and the rowwise projections are
nonexpansive, and the latter are right-orthogonally equivariant. The choices
of the radii in Theorem~\ref{thm: OS} ensure that these projection sets contain
the population intercept and factor rows. Consequently, the global spectral
bounds pass unchanged, up to numerical constants, to the trimmed initializer.
By the optimality of the common Procrustes rotation $\mb O^0$, the same bounds
hold under $\mb O^0$. The centered intercept projection also gives
\begin{equation*}
    \sqrt{\frac{n}{\sigma_r^*}}
    \bignorm{\bo\zeta^0-\zeta_0\mb 1_p}_\infty
    \lesssim\sqrt{\frac{\mu r\sigma_1^*}{p}},
\end{equation*}
and the factor projections give the remaining rowwise input bounds in
Theorem~\ref{thm: general OS}.

The leave-one-out spectral estimators remain independent of the omitted row
or column after deterministic trimming. The same projection properties and
the optimality of the relative common rotations give the required
leave-one-out bounds with rate $\xi^{0,\sf{loo}}$. Equation
\eqref{eq: transfer from relative to population-aligned LOO rotations}
accounts for the additional factor $\kappa$ when these relative errors are
transferred to population-aligned errors.

Substituting \eqref{eq: balanced spectral inputs for unilateral refinement}
into \eqref{eq: general unilateral refinement rate} gives
\begin{align*}
    \kappa_\psi\kappa\sqrt{\mu r}
    \Big[\frac{\xi^0}{\sqrt{n\wedge p}}
    +\kappa\xi^{0,\sf{loo}}\Big]
    \lesssim{}&
    \kappa_\psi\kappa^3
    \sigma\sqrt{\bar\sigma\sigma_1^*}
    \Big(\frac{\pi\mu^3r^3\rho^2}{d}\Big)^{\frac14}
    \xi^{\frac52}.
\end{align*}
Using $np(n\wedge p)=d^3/\rho^2$, we conclude that the rate $\xi^1$ in
Theorem~\ref{thm: general OS} satisfies
\begin{align*}
    \xi^1\lesssim{}&
    \bar\sigma\sqrt{\kappa r\log d}
    +\kappa_\psi\kappa^3
    \sigma\sqrt{\bar\sigma\sigma_1^*}
    \Big(\frac{\pi\mu^3r^3\rho^2}{d}\Big)^{\frac14}
    \xi^{\frac52}+\frac{\sigma^2\kappa^{\frac32}\sigma_r^*
    \mu^{\frac32}r^{\frac32}\rho\log d}
    {\sqrt\pi\,d^{\frac32}}.
\end{align*}

It remains to verify the smallness requirements of the general theorem. The
sampling and signal-strength bounds in
Assumption~\ref{assumption: signal strength and incoherence degree}(a),
applied term by term to the three terms in the preceding display, show that
its right-hand side is of smaller order than
\begin{equation*}
    \min\Big\{
    \underline c_\psi\sigma_r^*\sqrt{\frac{\pi\kappa\mu r}{d}},
    \frac{\underline c_\psi\sqrt{\pi(n\wedge p)}}
    {\kappa_\psi^2\kappa^3\mu r\sqrt\rho},
    \frac{\underline c_\psi\sigma_r^*\sqrt\pi}
    {\kappa_\psi^4\kappa^7\mu r\rho\sqrt d\log d}
    \Big\}.
\end{equation*}
The same comparison, together with
\eqref{eq: balanced spectral inputs for unilateral refinement}, verifies the
initialization requirements on $\xi^0$ and $\xi^{0,\sf{loo}}$. The first bound
in the preceding minimum gives the required rowwise smallness of $\xi^1$ and
the corresponding condition on $\xi^{0,\sf{loo}}$, while the second gives
$\kappa_\psi^2\kappa^3\mu r\sqrt\rho\,\xi^1
\ll\underline c_\psi\sqrt{\pi(n\wedge p)}$. Thus all assumptions of
Theorem~\ref{thm: general OS} hold. Applying that theorem proves the asserted
rowwise bound in the balanced parametrization. Scaling back through
Lemma~\ref{lemma: scaling equivalence} introduces the factors
$\omega^{-1/4}$ and $\omega^{1/4}$ in the $\mb X$- and $\mb Y$-errors,
respectively. This completes the proof. \qed

\subsection{Proof of Theorem~\ref{thm: general OS}}

To distinguish the refinement objective from the penalized objective $L$ used in the gradient-descent analysis, throughout this section we write
\begin{equation}
    \widetilde L([\bo \zeta,\mb X,\mb Y])
    \coloneqq
    \sum_{(i,j)\in\Omega}
    \left[
    -R_{i,j}\bigl(\zeta_j+\mb X_{i,\cdot}\mb Y_j\bigr)
    +\Psi\bigl(\zeta_j+\mb X_{i,\cdot}\mb Y_j\bigr)
    \right]
    \label{eq: unilateral refinement objective}
\end{equation}
for the non-penalized negative log-likelihood underlying Algorithm~\ref{alg: OS}. To facilitate the decoupling argument, we introduce its leave-one-out variants: for every $i\in[n]$ and $j \in[p]$,
\begin{align}
    & \widetilde L^{(i)}([\bo \zeta,\mb X,\mb Y])
    \coloneqq
    \sum_{(i',j)\in\Omega, i' \neq i}
    \left[
    -R_{i',j}\bigl(\zeta_j+\mb X_{i',\cdot}\mb Y_j\bigr)
    +\Psi\bigl(\zeta_j+\mb X_{i',\cdot}\mb Y_j\bigr)
    \right] \\ 
    & \qquad + \pi \sum_{j \in [p]}\left[
    - \psi(M_{i,j}^*)\bigl(\zeta_j+\mb X_{i,\cdot}\mb Y_j\bigr)
    +\Psi\bigl(\zeta_j+\mb X_{i,\cdot}\mb Y_j\bigr)
    \right]  , \\ 
    & \widetilde L^{(-j)}([\bo \zeta,\mb X,\mb Y])
    \coloneqq
    \sum_{(i,j')\in\Omega, j'\neq j}
    \left[
    -R_{i,j'}\bigl(\zeta_{j'}+\mb X_{i,\cdot}\mb Y_{j'}\bigr)
    +\Psi\bigl(\zeta_{j'}+\mb X_{i,\cdot}\mb Y_{j'}\bigr)
    \right] \\ 
    & \qquad + \pi \sum_{i \in [n]}\left[
    - \psi(M_{i,j}^*)\bigl(\zeta_j+\mb X_{i,\cdot}\mb Y_j\bigr)
    +\Psi\bigl(\zeta_j+\mb X_{i,\cdot}\mb Y_j\bigr)
    \right]. 
\end{align}

We use the following probability bookkeeping throughout this proof. A concentration bound for one fixed row or column is invoked with failure probability $O(d^{-c-1})$, whereas a bound for one fixed pair of row or column indices is invoked with failure probability $O(d^{-c-2})$. Since $n+p\leq2d$ and $np\leq d^2$, the corresponding uniform events have failure probability $O(d^{-c})$. The numerical constants in the Bernstein thresholds are enlarged accordingly; this does not alter any displayed rate because the thresholds already scale with $\log d$. The truncation event in Assumption~\ref{assumption: noise and sampling}, whose failure probability is $O(d^{-c-10})$, is absorbed into all these events.

Fix the constant $c_{\mc D_\psi}$ in Assumption~\ref{assumption: link function}
sufficiently large relative to the constants in the centered intercept and
factor row bounds of Theorem~\ref{thm: general OS}, including the fixed-point
neighborhoods used below. For each parameter triple in these neighborhoods,
$|\mc M(\bo\theta)_{i,j}-\zeta_0|
\leq|\zeta_j-\zeta_0|+\norm{\mb X_i}_2\norm{\mb Y_j}_2
\lesssim\mu r\sigma_1^*/\sqrt{np}$.
Thus the global and leave-one-out predictors, and the mixed predictors used
in the unilateral updates, belong to $\mc D_\psi$. Since this domain is an
interval, the natural-parameter interpolation segments used in the
mean-value and Taylor expansions remain in the same domain.

\subsubsection{Rowwise Error Control for Algorithm~\ref{alg: OS}. }
This part aims to establish \eqref{eq: row-wise error for unilateral refinement}. 
The core of our argument combines a fixed-point theorem with the convexity of the loss function when one factor is held fixed, which establishes the existence of a stationary point. In addition, the uniqueness of the solution follows from the strong convexity of the objective function.

\paragraph{\texorpdfstring{Update the left factor $\mb X$}{Update the left factor X}.}
We first write the gradient and Hessian forms of $\widetilde L$ with respect to $\mb X$ as: 
\begin{align}
    & \nabla_{\mb X} \widetilde L([\bo \zeta, \mb X, \mb Y]) = \Big[\big(\psi(\mb 1_n\bo \zeta\t + \mb X {\mb Y}\t ) - \mb R \big) \circ \bo \Omega\Big] \mb Y \in \bb R^{n \times r},  
    \label{eq: gradient of loss function with respect to X} 
    \\ 
    & \nabla^2_{\mb X} \widetilde L([\bo \zeta, \mb X, \mb Y])= \sf{diag}\bigg[\mb Y\t \diag\big(\psi'(\zeta_j + \mb x_i\t \mb y_j) \cdot \ind\{ (i,j) \in \Omega \} \big)_{j \in [p]}\mb Y \bigg]_{ i \in [n]} \in \bb R^{nr \times nr}. 
    \label{eq: hessian of loss function with respect to X}
\end{align}
{\emergencystretch=1.5em
In addition, we denote
\begin{align*}
\nabla_{\mb X_i} \widetilde L([\bo \zeta, \mb X, \mb Y])
&\coloneqq [\nabla_{\mb X} \widetilde L([\bo \zeta, \mb X, \mb Y])]_{i,\cdot}\in\bb R^{1\times r}
\end{align*}
and
\begin{align*}
\nabla^2_{\mb X_i} \widetilde L([\bo \zeta, \mb X, \mb Y])
&\coloneqq \nabla^2_{\mb X} \widetilde L([\bo \zeta, \mb X, \mb Y])_{(i-1)r+1: ir, (i-1)r +1: ir}.
\end{align*}
\par}

Expanding the loss function $\widetilde L$ with respect to the left factor $\mb X$, one directly obtains from the stationary condition that
\begin{align}
    & \mb 0 =  \sf{vec}(\nabla_{\mb X} \widetilde L([\bo \zeta^0, \mb X, \mb Y^0 \mb O^0]))\\
    & =   \sf{vec}(\nabla_{\mb X} \widetilde L([\bo \zeta^0, \mb X^*, \mb Y^0 \mb O^0])) + \nabla^2_{\mb X} \widetilde L([\bo \zeta^0, \mb X^*, \mb Y^0 \mb O^0]) \sf{vec}(\mb X - \mb X^*)\\
    & + \int_0^1 \Big[\nabla^2_{\mb X} \widetilde L\big([\bo \zeta^0, (1-t)\mb X^* + t \mb X, \mb Y^0 \mb O^0]\big) - \nabla^2_{\mb X} \widetilde L([\bo \zeta^0, \mb X^*, \mb Y^0 \mb O^0]) \Big]\mathrm d t \cdot \sf{vec}(\mb X\\
    &\quad - \mb X^*). 
    \label{eq: taylor expansion of loss function with respect to X (matrix form)}
\end{align}
To characterize the difference $\mb X - \mb X^*$, it amounts to understanding the gradient and the Hessian term presented in the above expansion. 

It is worth mentioning that the matrices $\mb E$ and $\bo \Omega$ are both dependent on $\mb Y^0$, resulting in extra complication. 
This challenge is addressed by introducing the leave-one-row-out version of $\bo \zeta^0$ and $\mb Y^0 \mb O^0$; in particular, the leave-one-out errors for the double SVD stage output have been controlled by \eqref{eq: leave-one-out perturbation for truncated estimates in step 2} and their general forms are assumed in Theorem~\ref{thm: general OS}. 
Replacing the original forms by the leave-one-out counterparts, we obtain for the $i$-th row of $\mb X$ that 
\begin{align}
    & \mb 0=   \nabla_{\mb X_i} \widetilde L([\bo \zeta^{0, (i)}, \mb X^*, \mb Y^{0, (i)}\tilde{\mb O}^{0,(i)}])\t +\nabla_{\mb X_i} \widetilde L([\bo \zeta^0, \mb X^*, \mb Y^0 \mb O^0])\t - \nabla_{\mb X_i} \widetilde L([\bo \zeta^{0,(i)}, \mb X^*,\\
    &\quad \mb Y^{0, (i)}\tilde{\mb O}^{0,(i)}])\t\\
    & + \big[\nabla^2_{\mb X_i} \widetilde L([\bo \zeta^0, \mb X^*, \mb Y^0 \mb O^0])\big]\cdot (\mb X - \mb X^*)_{i} \\ 
    & + \int_0^1 \big[\nabla^2_{\mb X_i} \widetilde L\big([\bo \zeta^0, (1-t)\mb X^* + t \mb X, \mb Y^0 \mb O^0] \big) - \nabla^2_{\mb X_i} \widetilde L([\bo \zeta^0, \mb X^*, \mb Y^0 \mb O^0]) \big]\mathrm d t \cdot (\mb X - \mb X^*)_{i }. 
    \label{eq: taylor expansion of loss function with respect to X}
\end{align}
Next, multiplying \eqref{eq: taylor expansion of loss function with respect to X} by $\big(\nabla^2_{\mb X_i} \widetilde L([\bo \zeta^0, \mb X^*, \mb Y^0 \mb O^0])\big)^{-1}$ and rearranging the terms yields that 
\begin{align}
    & \mb X_{i}  - {\mb X^*_{i}}   
    =  -\big[\nabla_{\mb X_i}^2 \widetilde L([\bo \zeta^0, \mb X^*, \mb Y^0 \mb O^0])\big]^{-1}   \Big\{ \nabla_{\mb X_i} \widetilde L([\bo \zeta^{0, (i)}, \mb X^*, \mb Y^{0, (i)}\tilde {\mb O}^{0,(i)}])\t \\
    & \qquad  + \Big[\nabla_{\mb X_i} \widetilde L([\bo \zeta^0, \mb X^*, \mb Y^0 \mb O^0]) - \nabla_{\mb X_i} \widetilde L([\bo \zeta^{0,(i)}, \mb X^*,  \mb Y^{0, (i)}\tilde {\mb O}^{0,(i)}])\Big]\t
        \Big\}  - \big[\nabla_{\mb X_i}^2 \widetilde L([\bo \zeta^0,\\
    &\quad \mb X^*, \mb Y^0 \mb O^0])\big]^{-1}\\
        & \cdot \int_0^1 \Big[\nabla_{\mb X_i}^2 \widetilde L\big([\bo \zeta^0, (1-t)\mb X^* + t \mb X, \mb Y^0 \mb O^0] \big) - \nabla^2_{\mb X_i} \widetilde L([\bo \zeta^0, \mb X^*, \mb Y^0 \mb O^0]) \Big]\mathrm d t \cdot  (\mb X - \mb X^*)_{i} 
    \label{eq: taylor expansion of loss function with respect to X after rearranging}.
\end{align}
The terms on the right-hand side will be parsed separately in what follows.  

\bigskip
 
\emph{Upper bounding the gradient terms in Eq.~\eqref{eq: taylor expansion of loss function with respect to X after rearranging}. }
Regarding the gradient at the true parameter $[\bo \zeta^*, \mb X^*, \mb Y^*]$, it follows from \eqref{eq: gradient of loss function with respect to X} that 
\begin{align}
    & \norm{\nabla_{\mb X_i} \widetilde L([\bo \zeta^{0, (i)}, \mb X^*, \mb Y^{0, (i)}\tilde{\mb O}^{0,(i)}]) }_2
    \leq  \norm{\mc P_\Omega(\mb E)_{i,\cdot}\mb Y^{0, (i)}}_2 \\ 
    & \qquad\\
    &\quad + \Bignorm{\Big\{ \big[\psi(\mb 1_n{ \bo \zeta^{0,(i)}}\t + \mb X^* {{}\tilde{\mb O}^{0,(i)}}\t {\mb Y^{0, (i)}}\t ) - \psi(\mb 1_n {\bo \zeta^*}\t + \mb X^*{\mb Y^*}\t  )\big] \circ \bo \Omega \Big\}_{i,\cdot}  \mb Y^{0, (i)} \tilde{\mb O}^{0, (i)} }_2\\
    & \lesssim   \underbrace{\norm{\mc P_\Omega(\mb E)_{i,\cdot}\mb Y^{0, (i)} \mb O^{0, (i)}}_2}_{\alpha_1} + \underbrace{  \bignorm{\big[(\bo \zeta^{0,(i)} - \bo \zeta^*)\t \circ \bo\Omega_{i,\cdot} \big] \diag(\psi'(m_{i,j}))_{j\in[p]} \mb Y^{0,(i)}} }_{\alpha_2}\\
    & + \underbrace{\bignorm{ \big[\big(\mb X^*_{i,\cdot} (\mb Y^{0,(i)} \tilde{\mb O}^{0,(i)} - \mb Y^*)\t\big)\circ \bo \Omega_{i,\cdot} \big] \diag(\psi'(m'_{i,j}))_{j\in[p]}\mb Y^{0,(i)}} }_{\alpha_3},
    \label{eq: decomposition of gradient term in os analysis}
    \end{align}
    where we used the mean value theorem \emph{twice} for $\psi(\zeta_j^{0,(i)} + (\mb X^* {{}\tilde{\mb O}^{0,(i)\top}} {\mb Y^{0, (i)}\top} )_{i,j}) - \psi(\zeta_j^* + (\mb X^*{\mb Y^*}\t)_{i,j})$ and noticed that the intermediate values, denoted by $m_{i,j}$ and $m_{i,j}'$ for $j\in[p]$, always fall within $\mc D_\psi$. Hence, $\psi'(m_{i,j})$ and $\psi'(m'_{i,j})$ are bounded above by $\sigma^2$. We analyze each term in \eqref{eq: decomposition of gradient term in os analysis} in what follows.
    \begin{itemize}
        \item For $\alpha_1$, the upper bound is a direct consequence of the matrix Bernstein inequality together with the mutual independence between $\mc P_\Omega(\mb E)_{i,\cdot}$ and $\mb Y^{0,(i)}$: 
        \begin{align}
            & \norm{\mc P_\Omega(\mb E)_{i,\cdot}\mb Y^{0, (i)} }_2 \lesssim \sigma \sqrt{\pi r \log d} \bignorm{\mb Y^{0,(i)}} +  B \bignorm{\mb Y^{0,(i)}}\ti \log d \\ 
            \lesssim & \sigma \sqrt{\pi r \log d} (\sigma_1^*)^{\frac12} + B\sqrt{\frac{\mu r}{p}} (\sigma_1^*)^{\frac12} \log d \lesssim \sigma \sqrt{\pi r \log d} (\sigma_1^*)^{\frac12} 
        \end{align}
        with probability at least $1-O(d^{-c-1})$, where the last line follows from the condition $B \sqrt{\mu r \log d} \ll \sigma \sqrt{\pi p}$. 
        \item Regarding $\alpha_2$, it first follows from the matrix Bernstein inequality conditional on $(\bo \zeta^{0,(i)} - \bo \zeta^*)$ that 
        \begin{align}
            & \alpha_2 = \bignorm{\big((\bo \zeta^{0,(i)} - \bo \zeta^*)\t \circ \bo\Omega_{i,\cdot} \big)  \diag(\psi'(m_{i,j}))_{j\in[p]} \mb Y^{0,(i)}}  \\ 
            \lesssim & \bignorm{\pi \big(\bo \zeta^{0,(i)} - \bo \zeta^*\big)\t  \diag(\psi'(m_{i,j}))_{j\in[p]} \mb Y^{0,(i)}}   + \sqrt{V \log d} + L \log d \\ 
            \lesssim & \pi \sigma^2 \bignorm{\bo \zeta^{0,(i)} - \bo \zeta^*} (\sigma_1^*)^{\frac12} + \sigma^2 \sqrt{\pi} \bignorm{\bo \zeta^{0,(i)} - \bo \zeta^*} \sqrt{\frac{\mu r}{p}} (\sigma_1^*)^{\frac12} \sqrt{\log d} + \frac{\sigma^2\mu r(\sigma_1^*)^{\frac32}\log d}{\sqrt{np^2}}  \\
            \lesssim & \pi \sigma^2 \sqrt{\mu r }\bignorm{\bo \zeta^{0,(i)} - \bo \zeta^*} (\sigma_1^*)^{\frac12}  + \frac{\sigma^2\mu r(\sigma_1^*)^{\frac32}\log d}{\sqrt{np^2}} \\
            \lesssim &  \kappa_\psi \sqrt{\pi \mu r \kappa / n }(\xi^0 + \xi^{0,\sf{loo}}) (\sigma_1^*)^{\frac12}  +  \sigma^2 \mu r(\sigma_1^*)^{\frac32} \log d  / \sqrt{ (np^2) }
        \end{align}
        with probability at least $1-O(d^{-c-1})$, where 
        we used the facts that
        $\bignorm{\bo \zeta^{0,(i)} - \bo \zeta^*}_\infty
        \vee\bignorm{\bo \zeta^{0,(i)}-\zeta_0\mb 1_p}_\infty
        \lesssim\sigma_1^*\sqrt{\mu r/(np)}$, and
        the quantities $V$ and $L$ arising from the matrix Bernstein inequality are evaluated as 
        \begin{align*}
            & V \coloneqq  \sum_{j\in[p]} \pi ( 1 - \pi) (\zeta^{0,(i)}_j - \zeta^*_{j})^2 \bignorm{\mb Y^{0,(i)}_j}_2^2 \psi'(m_{i,j})^2 ,\\ 
            & L\coloneqq \max_{j\in[p]} \big|\zeta^{0,(i)}_j - \zeta^*_{j}\big| \bignorm{\mb Y^{0,(i)}_j}_2 \psi'(m_{i,j}) \lesssim \frac{\sigma^2 \mu r (\sigma_1^*)^{\frac32}}{ \sqrt{np^2} }  .
        \end{align*}
        
        \item As for the final term $\alpha_3$, its upper bound can be derived in essentially the same way as for $\alpha_2$. In particular, with probability at least $1-O(d^{-c-1})$, 
\begin{align}
    & \alpha_3 =\bignorm{ \big[\big(\mb X^*_{i,\cdot} (\mb Y^{0,(i)} \tilde{\mb O}^{0,(i)} - \mb Y^*)\t\big)\circ \bo \Omega_{i,\cdot} \big] \diag(\psi'(m'_{i,j}))_{j\in[p]}\mb Y^{0,(i)}} \\ 
    \lesssim & \bignorm{ \pi\big(\mb X^*_{i,\cdot} (\mb Y^{0,(i)} \tilde{\mb O}^{0,(i)} - \mb Y^*)\t\big) \diag(\psi'(m'_{i,j}))_{j\in[p]}\mb Y^{0,(i)}}  \\ 
    & + \sigma^2 \sqrt{\frac{\mu r}{n}} (\sigma_1^*)^{\frac12} \sqrt{\pi } \bignorm{\mb Y^{0,(i)} \tilde{\mb O}^{0,(i)} - \mb Y^*}\fb \sqrt{\frac{\mu r}{p}} (\sigma_1^*)^{\frac12} \sqrt{\log d}\\
    & + \sigma^2 \sqrt{\frac{\mu r}{n}} (\sigma_1^*)^{\frac12} \bignorm{\mb Y^{0,(i)} \tilde{\mb O}^{0,(i)} - \mb Y^*}\ti \sqrt{\frac{\mu r}{p}} (\sigma_1^*) \log d \\
    \lesssim & \kappa_\psi \sqrt{\pi \mu r\kappa / n}  (\xi^0 + \kappa\xi^{0,\sf{loo}}) (\sigma_1^*)^{\frac12} + \sigma^2 \sqrt{\frac{\mu^3 r^3}{np^2}}  (\sigma_1^*)^{\frac32} \log d .
\end{align}
  
    \end{itemize}

    Substituting the above pieces into \eqref{eq: decomposition of gradient term in os analysis} immediately gives, for the fixed row $i$, with probability at least $1-O(d^{-c-1})$ 
    \begin{align}
        & \norm{\nabla_{\mb X_i} \widetilde L([\bo \zeta^{0, (i)}, \mb X^*, \mb Y^{0, (i)}\tilde{\mb O}^{0,(i)}]) }_2 \\
        \lesssim &  \sigma \sqrt{\pi r \log d} (\sigma_1^*)^{\frac12}  +   \kappa_\psi \sqrt{\kappa \pi \mu r / n }(\xi^0 + \kappa\xi^{0,\sf{loo}}) (\sigma_1^*)^{\frac12} + \sigma^2 \sqrt{\frac{\mu^3 r^3}{np^2}}  (\sigma_1^*)^{\frac32} \log d.
        \label{eq: gradient term evaluated at leave-one-out location}
    \end{align}

We next control the gradient difference term. To this end, observe that
\begin{align}
& \nabla_{\mb X_i} \widetilde L([\bo \zeta^0, \mb X^*, \mb Y^0 \mb O^0])
- \nabla_{\mb X_i} \widetilde L([\bo \zeta^{0,(i)}, \mb X^*, \mb Y^{0,(i)}\tilde{\mb O}^{0,(i)}]) \nonumber\\
    & = {} 
\mc P_\Omega\big(\psi(\mb 1_n \bo \zeta^{0\top} + \mb X^*(\mb Y^0 \mb O^0)\t)-\mb R\big)_{i,\cdot}\mb Y^0 \mb O^0
-
\mc P_\Omega\big(\psi(\mb 1_n \bo \zeta^{0,(i)\top}\\
    &\quad + \mb X^*(\mb Y^{0,(i)}\tilde{\mb O}^{0,(i)})\t)-\mb R\big)_{i,\cdot}
\mb Y^{0,(i)}\tilde{\mb O}^{0,(i)} \nonumber\\
    & = {} 
\mc P_\Omega\big(
\psi(\mb 1_n \bo \zeta^{0,(i)\top} + \mb X^*(\mb Y^{0,(i)}\tilde{\mb O}^{0,(i)})\t)
-\psi(\mb 1_n \bo \zeta^{*\top} + \mb X^*{\mb Y^*}\t)
-\mb E
\big)_{i,\cdot}
\big(\mb Y^0 \mb O^0\\
    &\quad -\mb Y^{0,(i)}\tilde{\mb O}^{0,(i)}\big) \nonumber\\
& +
\mc P_\Omega\big(
\psi(\mb 1_n \bo \zeta^{0\top} + \mb X^*(\mb Y^0\mb O^0)\t)
-
\psi(\mb 1_n \bo \zeta^{0,(i)\top} + \mb X^*(\mb Y^{0,(i)}\tilde{\mb O}^{0,(i)})\t)\big)_{i,\cdot}\\ 
& \quad \cdot 
\big[ \mb Y^{0,(i)}\tilde{\mb O}^{0,(i)} +\big(\mb Y^0 \mb O^0 - \mb Y^{0,(i)}\tilde{\mb O}^{0,(i)} \big) \big].
\end{align}
Consequently, for the $i$-th row, we have
\begin{align}
& \norm{
\Big[
\nabla_{\mb X} \widetilde L([\bo \zeta^0,\mb X^*,\mb Y^0\mb O^0])
-
\nabla_{\mb X} \widetilde L([\bo \zeta^{0,(i)},\mb X^*,\mb Y^{0,(i)}\tilde{\mb O}^{0,(i)}])
\Big]_{i,\cdot}
}_2 \nonumber\\
    & \lesssim {} 
\Big[
\bignorm{
\mc P_\Omega\big(
\psi(\mb 1_n \bo \zeta^{0,(i)\top} + \mb X^*(\mb Y^{0,(i)}\tilde{\mb O}^{0,(i)})\t)
-
\psi(\mb 1_n \bo\zeta^{*\top} + \mb X^*{\mb Y^*}\t)
\big)_{i,\cdot}
}_2\\
    &\quad +
\norm{\mc P_\Omega(\mb E_{i,\cdot})}_2
\Big]
\bignorm{\mb Y^0\mb O^0-\mb Y^{0,(i)}\tilde{\mb O}^{0,(i)}} \nonumber\\
& +
\norm{
\psi(\mb 1_n \bo \zeta^{0\top} + \mb X^*(\mb Y^0\mb O^0)\t)_{i,\cdot}
-
\psi(\mb 1_n \bo \zeta^{0,(i)\top} + \mb X^*(\mb Y^{0,(i)}\tilde{\mb O}^{0,(i)})\t)_{i,\cdot}
}_2\\ 
& \qquad \cdot 
\Big[ \bignorm{\diag(\Omega_{i,j})_{j\in[p]}\mb Y^{0,(i)}\tilde{\mb O}^{0,(i)}} + \bignorm{\mb Y^0\mb O^0-\mb Y^{0,(i)}\tilde{\mb O}^{0,(i)}} \Big] \nonumber\\
    & \lesssim {} 
\Big[
\sigma^2
\sqrt{\frac{\pi \mu r \sigma_1^* \log d}{n}}
\frac{\xi^0+\kappa\xi^{0,\sf{loo}}}{\underline c_\psi \sqrt{\pi\sigma_r^*}} + \sigma^2 \mu r \sigma_1^* \log d / \sqrt{np}
+
\sigma\sqrt{\pi p}
\Big]
\frac{\kappa\xi^{0,\sf{loo}}}{\underline c_\psi \sqrt{\pi\sigma_r^*}}\\
    &\quad +
\sigma^2
\sqrt{\frac{\mu r \sigma_1^*}{n}}
\frac{\kappa\xi^{0,\sf{loo}}}{\underline c_\psi \sqrt{\pi\sigma_r^*}}
\sqrt{\pi\sigma_1^*} \nonumber\\
    & \ll  
\underline c_\psi \pi\sigma_r^*
\frac{\kappa\xi^{0,\sf{loo}}}{\underline c_\psi \sqrt{\pi\sigma_r^*}},
\label{eq: leave-one-out difference for gradient wrt X}
\end{align}
with probability at least $1-O(d^{-c-1})$. Here, the last inequality follows from the assumptions $\frac{\xi^0 + \kappa\xi^{0,\sf{loo}}}{\underline c_\psi \sqrt{\pi \sigma_r^*}} \ll \sqrt{\sigma_r^*}$,
$\kappa_\psi^2\kappa^2\mu r\log d\ll\pi(n\wedge p)$ in
\eqref{eq: general unilateral refinement smallness conditions}, and
$\underline c_\psi \sigma_r^* \gg \sigma\sqrt{p/\pi}$. We used the following bound from the Bernstein inequality: with probability at least $1-O(d^{-c-1})$, 
\begin{align}
    & \bignorm{
\mc P_\Omega\big(
\psi(\mb 1_n \bo \zeta^{0,(i)\top} + \mb X^*(\mb Y^{0,(i)}\tilde{\mb O}^{0,(i)})\t)
-
\psi(\mb 1_n \bo\zeta^{*\top} + \mb X^*{\mb Y^*}\t)
\big)_{i,\cdot}
}_2 
\\ 
\lesssim & \sigma^2 \sqrt{\frac{\mu r \sigma_1^*}{n}} \Big[\frac{\sqrt{\pi}(\xi^0 + \kappa\xi^{0,\sf{loo}})}{\underline c_\psi \sqrt{\pi \sigma_r^*}} \sqrt{\log d} + \sqrt{\frac{\mu r \sigma_1^*}{p}} \log d\Big] .
\end{align} 
By \citet[Corollary~2.15]{brailovskaya2022universality} and the independence between
$\{\Omega_{i,j}\}_{j\in[p]}$ and $\mb Y^{0,(i)}$, we have
\begin{align}
&\norm{\diag(\Omega_{i,j})_{j\in[p]}\mb Y^{0,(i)}\tilde{\mb O}^{0,(i)}}
\lesssim \sqrt{\pi \sigma_1^*} + (\frac{\pi \mu r \sigma_1^*}{p})^{\frac14}(\pi \sigma_1^*)^{\frac14}(\log d)^{\frac34} +  (\frac{\pi \mu r \sigma_1^*}{p})^{\frac12}\sqrt{\log d} \\ 
& \hspace{3cm} +  (\frac{ \mu r \sigma_1^*}{p})^{\frac16} (\pi \sigma_1^*)^{\frac13} (\log d)^{\frac23} + (\frac{ \mu r \sigma_1^*}{p})^{\frac12} \log d \lesssim \sqrt{\pi \sigma_1^*} , 
\label{eq: upper bound on Omega_i Y^{0,(i)}}
\end{align}
with probability at least $1-O(d^{-c-1})$, where we used
$\mu r(\log d)^4 \ll \pi(n\wedge p)$. 

\bigskip

\emph{Bounding the Hessian term in Eq.~\eqref{eq: taylor expansion of loss function with respect to X after rearranging}. }
We next turn to dealing with the Hessian $\nabla^2_{\mb X_i}\widetilde L([\bo \zeta^0, \mb X^*, \mb Y^0\mb O^0])$ with respect to the $i$-th row of $\mb X$. Substituting the Hessian matrix with the leave-one-out counterparts, it follows that
\begin{align}
    & \norm{\Big(\nabla^2_{\mb X_i} \widetilde L([\bo \zeta^0, \mb X^*, \mb Y^0 \mb O^0]) - \nabla^2_{\mb X_i} \widetilde L([\bo \zeta^{0, (i)}, \mb X^*, \mb Y^{0, (i)}\widetilde{\mb O}^{0,(i)}])\Big)}  \\ 
    & \leq   \bigg\|{\mb O^0}\t {\mb Y^0}\t \diag\Big(\psi'(\zeta_j^0 + \mb X_{i,\cdot}^*(\mb Y^0 \mb O^0)_{j}) \ind\{(i,j)\in \bo \Omega\}\Big)_{j\in[p]} \mb Y^0 {\mb O^0}\\
    & - {{}\tilde{\mb O}^{0, (i)}}\t  {\mb Y^{0, (i)}}\t \diag\Big(\psi'(\zeta_j^0 + \mb X_{i,\cdot}^*(\mb Y^0 \mb O^0)_{j}) \ind\{(i,j)\in \bo \Omega\}\Big)_{j\in[p]} \mb Y^{0, (i)} \tilde{\mb O}^{0, (i)} \bigg\| \\
    & +    \bigg\|   {\mb Y^{0, (i)}}\t   \sf{diag}\Big(\big(\psi'(\zeta_j^0  + {\mb x_i^*}\t {\mb O^{0}}\t\mb Y_j^0) - \psi'(\zeta_j^{0, (i)}  + {\mb x_i^*}\t {{}\tilde{\mb O}^{0, (i)}}\t \mb Y_j^{0,(i)}) \big) \ind\{ (i,\\
    &\quad j) \in \Omega \} \Big)_{j \in [p]} \mb Y^{0, (i)}\bigg\|\\
    & \lesssim   \sigma^2 \frac{\kappa\xi^{0,\sf{loo}}}{\underline c_\psi\sqrt{\pi \sigma_r^*}} \sqrt{\pi \sigma_1^*} + \sigma^2 \Big(\frac{\xi^{0,\sf{loo}}}{\underline c_\psi\sqrt{\pi n}} + \sqrt{\frac{\mu r}{n}} (\sigma_1^*)^{\frac12} \frac{\kappa\xi^{0,\sf{loo}}}{\underline c_\psi\sqrt{\pi \sigma_r^*}} \Big) \pi \sigma_1^* \ll \underline c_\psi \pi \sigma_r^*
    \label{eq: leave-one-out perturbation for Hessian of loss function with respect to X}
\end{align}
with probability at least $1-O(d^{-c-1})$. Here, we reason the second inequality using the triangle inequality as follows: 
\begin{align}
    & \bigg\|{\mb O^0}\t {\mb Y^0}\t \diag\Big(\psi'(\zeta_j^0 + \mb X_{i,\cdot}^*(\mb Y^0 \mb O^0)_{j}) \ind\{(i,j)\in \bo \Omega\}\Big)_{j\in[p]} \mb Y^0 {\mb O^0} \\ 
    & - {{}\tilde{\mb O}^{0, (i)}}\t  {\mb Y^{0, (i)}}\t \diag\Big(\psi'(\zeta_j^0 + \mb X_{i,\cdot}^*(\mb Y^0 \mb O^0)_{j}) \ind\{(i,j)\in \bo \Omega\}\Big)_{j\in[p]} \mb Y^{0, (i)} \tilde{\mb O}^{0, (i)} \bigg\| \\ 
    & \lesssim   \max_{i\in[n], j\in[p]}\psi'(\zeta_j^0 + \mb X_{i,\cdot}^*(\mb Y^0 \mb O^0)_{j})  \cdot \Big[2 \bignorm{\mb Y^0\mb O^0 - \mb Y^{0,(i)} \widetilde{\mb O}^{0,(i)}} \norm{\diag(\Omega_{i,j})_{j\in[p]}\mb Y^0}\\
    & \qquad + \bignorm{\mb Y^0\mb O^0 - \mb Y^{0,(i)} \widetilde{\mb O}^{0,(i)}}^2 
    \Big]\\ 
    & \lesssim   \sigma^2 \frac{\kappa\xi^{0,\sf{loo}}}{\underline c_\psi\sqrt{\pi \sigma_r^*}} \sqrt{\pi \sigma_1^*}  +\sigma^2\big( \frac{\kappa\xi^{0,\sf{loo}}}{\underline c_\psi \sqrt{\pi \sigma_r^*}} \big)^2  \lesssim \sigma^2 \frac{\kappa\xi^{0,\sf{loo}}}{\underline c_\psi \sqrt{\pi \sigma_r^*}} (\pi \sigma_1^*)^{\frac12} \ll \underline c_\psi \pi \sigma_r^*,
    \label{eq: Hessian bound 1}\\
    & \Big\|   {\mb Y^{0, (i)}}\t   \sf{diag}\Big(\big(\psi'(\zeta_j^0  + {\mb X_i^*}\t {\mb O^{0}}\t\mb Y_j^0) - \psi'(\zeta_j^{0, (i)}  + {\mb X_i^*}\t {{}\tilde{\mb O}^{0, (i)}}\t \mb Y_j^{0,(i)}) \big) \ind\{ (i,\\
    &\quad j) \in \Omega \} \Big)_{j \in [p]} \mb Y^{0, (i)}\Big\|\\
    & \lesssim   \bignorm{\diag(\Omega_{i,j})_{j\in[p]} \mb Y^{0,(i)}}^2 \max_{j\in[p]} \big|\psi'(\zeta_j^0  + {\mb X_i^*}\t {\mb O^{0}}\t\mb Y_j^0) - \psi'(\zeta_j^{0, (i)}\\
    &\quad + {\mb X_i^*}\t {{}\tilde{\mb O}^{0, (i)}}\t \mb Y_j^{0,(i)}) \big|\\
    & \lesssim   \bignorm{\diag(\Omega_{i,j})_{j\in[p]} \mb Y^{0,(i)}}^2 \sigma^2 \big(\norm{\bo\zeta^{0} - \bo \zeta^{0,(i)}}_\infty + \norm{\mb X^*}\ti \bignorm{\mb Y^0\mb O^0 - \mb Y^{0,(i)} \widetilde{\mb O}^{0,(i)}}\ti \big)\\
    & \lesssim    \sigma^2 \Big[\frac{\xi^{0,\sf{loo}}}{\underline c_\psi \sqrt{\pi n}} + \sqrt{\frac{\mu r}{n}} (\sigma_1^*)^{\frac12} \frac{\kappa\xi^{0,\sf{loo}}}{\underline c_\psi \sqrt{\pi \sigma_r^*}} \Big] \pi \sigma_1^* \ll \underline c_\psi\pi \sigma_r^*
    \label{eq: Hessian bound 2}
\end{align}
with probability at least $1-O(d^{-c-1})$ by the assumption $\frac{\kappa\xi^{0,\sf{loo}}}{\underline c_\psi \sqrt{\pi \sigma_r^*}} \ll \sqrt{\frac{\mu r \sigma_1^*}{n \wedge p}} \ll \sqrt{\pi \sigma_r^*} / (\kappa \kappa_\psi^2)$ along with \eqref{eq: upper bound on Omega_i Y^{0,(i)}}. 

With the above in place, we are able to control the Hessian matrix. Invoking the matrix Bernstein inequality together with the independence between $\bo \Omega_{i,\cdot}$ and the leave-one-out coupling yields that 
\begin{align}
    & \norm{\nabla^2_{\mb X_i} \widetilde L([\bo \zeta^{0, (i)}, \mb X^*, \mb Y^{0, (i)}\tilde{\mb O}^{0,(i)}]) - \bb E_{\Omega}\big[\nabla^2_{\mb X_i} \widetilde L([\bo \zeta^{0, (i)}, \mb X^*, \mb Y^{0, (i)}\tilde{\mb O}^{0,(i)}]) \big]}\\  
    \lesssim &\sigma^2 \bignorm{\mb Y^{0, (i)}}\ti \bignorm{\mb Y^{0, (i)}} \sqrt{\pi \log d}  + \sigma^2 \bignorm{\mb Y^{0, (i)}}\ti^2 \log d \\
    \lesssim & \sigma^2 \sqrt{\frac{\mu r}{p}} \sigma_1^* \sqrt{\pi \log d} + \sigma^2 \frac{\mu r}{p} \sigma_1^* \log d \ll \underline c_\psi \pi \sigma_r^*
\end{align}
with probability at least $1-O(d^{-c-1})$ ($\bb E_{\bo \Omega_{i,\cdot}}$ represents the expectation over the randomness of $\bo \Omega_{i,\cdot}$). Indeed, the leading term in the preceding display, divided by $\underline c_\psi\pi\sigma_r^*$, is of order
$\kappa_\psi\kappa\sqrt{\mu r\log d/(\pi p)}$; hence both terms are $o(\underline c_\psi\pi\sigma_r^*)$ under
$\kappa_\psi^2\kappa^2\mu r\log d\ll\pi(n\wedge p)$ from
\eqref{eq: general unilateral refinement smallness conditions}. Moreover, the following population bounds hold on the same event:
\begin{align}
    & \norm{\bb E_{\Omega}\big[\nabla^2_{\mb X_i} \widetilde L([\bo \zeta^{0, (i)}, \mb X^*, \mb Y^{0, (i)}\tilde{\mb O}^{0,(i)}]) \big] }\lesssim \sigma^2 \pi \sigma_1^*, \\
    & \sigma_{\min}(\bb E_{\Omega}\big[\nabla^2_{\mb X_i} \widetilde L([\bo \zeta^{0, (i)}, \mb X^*, \mb Y^{0, (i)}\tilde{\mb O}^{0,(i)}]) \big]) \gtrsim \underline c_\psi \pi \sigma_r^*. 
\end{align}
This in turn implies that
\begin{align}
    & \norm{\nabla^2_{\mb X_i} \widetilde L([\bo \zeta^{0, (i)}, \mb X^*, \mb Y^{0, (i)}\widetilde{\mb O}^{0,(i)}])} \lesssim \sigma^2 \pi  \sigma_1^*, \\
    & \sigma_{\min}\big(\nabla^2_{\mb X_i} \widetilde L([\bo \zeta^{0, (i)}, \mb X^*, \mb Y^{0, (i)}\widetilde{\mb O}^{0,(i)}])\big) \gtrsim \underline c_\psi \pi  \sigma_r^* .
\end{align}

Plugging the above bounds together with \eqref{eq: leave-one-out perturbation for Hessian of loss function with respect to X} into the Hessian matrix evaluated at $[\bo \zeta^*, \mb X^*, \mb Y^0 \mb O^0]$ gives, for the fixed row $i$, with probability at least $1-O(d^{-c-1})$ 
\begin{align}
    & \norm{\nabla^2_{\mb X_i} \widetilde L([\bo \zeta^0, \mb X^*, \mb Y^0 \mb O^0])} \lesssim \pi \sigma^2 \sigma_1^*, \quad
      \sigma_{\min}\big(\nabla^2_{\mb X_i} \widetilde L([\bo \zeta^0, \mb X^*, \mb Y^0 \mb O^0])\big) \gtrsim \pi \underline c_\psi \sigma_r^*. 
        \label{eq: upper bound and lower bound on singular values of Hessian wrt X}
\end{align}

\bigskip 

\emph{Bounding the terms in \eqref{eq: taylor expansion of loss function with respect to X after rearranging}. }
With the above preparations in place, we proceed with controlling the terms in \eqref{eq: taylor expansion of loss function with respect to X after rearranging}. 
\begin{itemize}
\item
Combining \eqref{eq: gradient term evaluated at leave-one-out location}, \eqref{eq: leave-one-out difference for gradient wrt X}, and \eqref{eq: upper bound and lower bound on singular values of Hessian wrt X} yields an upper bound on the first term of \eqref{eq: taylor expansion of loss function with respect to X after rearranging}: 
\begin{align}
    & \Big\|\big[\nabla_{\mb X_i}^2 \widetilde L([\bo \zeta^0, \mb X^*, \mb Y^0 \mb O^0])\big]^{-1}   \Big\{ \nabla_{\mb X_i} \widetilde L([\bo \zeta^{0, (i)}, \mb X^*, \mb Y^{0, (i)}\tilde{\mb O}^{0,(i)}])\t \\
    & \qquad  + \Big[\nabla_{\mb X_i} \widetilde L([\bo \zeta^0, \mb X^*, \mb Y^0 \mb O^0]) - \nabla_{\mb X_i} \widetilde L([\bo \zeta^{0,(i)}, \mb X^*,  \mb Y^{0, (i)}\tilde{\mb O}^{0,(i)}])\Big]\t \Big\}
        \Big\| \\
        \lesssim & \frac{\sigma \sqrt{\kappa r \log d}}{\underline c_\psi  \sqrt{\pi \sigma_r^*}}   +  \kappa_{\psi}\kappa \sqrt{ \mu r } \frac{\xi^0 / \sqrt{n} + \kappa \xi^{0,\sf{loo}}  }{\underline c_\psi \sqrt{\pi \sigma_r^*}} + \kappa_\psi \kappa \sqrt{\frac{\mu^3 r^3}{np^2}}  (\sigma_1^*)^{\frac12}\log d / \pi \eqqcolon L_1
    \label{eq: upper bound on the first term wrt X}
\end{align}
holds with probability at least $1-O(d^{-c-1})$. By the definition of $\xi^1$, one has $L_1\lesssim\xi^1/(\underline c_\psi\sqrt{\pi\sigma_r^*})$. Consequently, the second condition in \eqref{eq: general unilateral refinement smallness conditions} gives $L_1 \ll \sqrt{\mu r\sigma_1^*/n}$.

\item
In view of \eqref{eq: taylor expansion of loss function with respect to X after rearranging}, we are left with controlling the remainder term in \eqref{eq: taylor expansion of loss function with respect to X after rearranging}. As before, we invoke the leave-one-out variants to decouple the inherent dependencies. To this end, we decompose the quantity of interest into three components as follows:
\begin{align}
    & \Bignorm{\nabla^2_{\mb X_i} \widetilde L\big([\bo \zeta^0, (1-t)\mb X^* + t \mb X, \mb Y^0 \mb O^0] \big)\\
    &\qquad - \nabla^2_{\mb X_i} \widetilde L([\bo \zeta^0, \mb X^*, \mb Y^0 \mb O^0])}\\
    & =   \Bignorm{\nabla^2_{\mb X_i} \widetilde L\big([\bo \zeta^{0,(i)}, (1-t)\mb X^* + t \mb X, \mb Y^{0,(i)}\tilde{\mb O}^{0,(i)}]\big)\\
    &\qquad - \nabla^2_{\mb X_i} \widetilde L([\bo \zeta^{0,(i)}, \mb X^*, \mb Y^{0, (i)}\tilde{\mb O}^{0,(i)}]) }\\
    & + \Bignorm{
        \nabla^2_{\mb X_i} \widetilde L\big([\bo \zeta^{0,(i)}, (1-t)\mb X^* + t \mb X, \mb Y^{0,(i)}\tilde{\mb O}^{0,(i)}] \big)\\
    &\qquad - \nabla^2_{\mb X_i} \widetilde L\big([\bo \zeta^0, (1-t)\mb X^* + t \mb X, \mb Y^0 \mb O^0] \big)}\\
    & + \norm{
    \nabla^2_{\mb X_i} \widetilde L\big([\bo \zeta^{0,(i)}, \mb X^*, \mb Y^{0,(i)} \tilde{\mb O}^{0,(i)}]\big)  - \nabla^2_{\mb X_i} \widetilde L\big([\bo \zeta^0, \mb X^*, \mb Y^0 \mb O^0] \big)
    }\\
    & \eqqcolon   \delta_1 + \delta_2 + \delta_3.
    \label{eq: decomposition for the remainder term of X}
\end{align}

For the first term $\delta_1$, the mean-value theorem implies that there exists $t_j \in [0,1]$ for $j \in [p]$ such that
\begin{align*}
     \delta_1 \leq {}& \max_{j\in[p]}\Big|\psi''\big(\zeta_j^{0,(i)} + [(1-t_j)\mb X^*_{i,\cdot}+t_j\mb X_{i,\cdot}]
     (\mb Y^{0,(i)}\tilde{\mb O}^{0,(i)})_j\big)\Big| \\
     & \quad \times \Big|(\mb X-\mb X^*)_{i,\cdot}(\mb Y^{0,(i)}\tilde{\mb O}^{0,(i)})_j\Big|
     \bignorm{\diag(\Omega_{i,j})_{j\in[p]}\mb Y^{0,(i)}}^2.
\end{align*}
Then we follow a similar argument to the one used for \eqref{eq: Hessian bound 2}. Assumption~\ref{assumption: link function} gives $|\psi''(x)|\lesssim\sigma^2$, and hence
\begin{align}
    & \delta_1 \leq \bignorm{\diag(\Omega_{i,j})_{j\in[p]} \mb Y^{0,(i)}}^2 \sigma^2 \norm{(\mb X- \mb X^*)_i}_2 \sqrt{\frac{\mu r}{p}}(\sigma_1^*)^{\frac12}\\
    &\quad \lesssim  \sigma^2 \sqrt{\frac{\mu r}{p}} \pi (\sigma_1^*)^{3/2}  \norm{\mb X_{i}  - {\mb X^*_{i}} }_2
    \label{eq: bound for first Hessian remainder}
\end{align}
with probability at least $1-O(d^{-c-1})$. 

For the second and third terms, assuming $\norm{\mb X- \mb X^*}_2 \ll \sqrt{\frac{\mu r}{n}} (\sigma_1^*)^{\frac12}$, 
one follows a similar argument to \eqref{eq: leave-one-out perturbation for Hessian of loss function with respect to X} to obtain that
\begin{align}
    & \delta_2 \vee \delta_3  \lesssim \sigma^2 \frac{\kappa \xi^{0,\sf{loo}}}{\underline c_\psi \sqrt{\pi \sigma_r^*}} \sqrt{\pi \sigma_1^*}
    + \sigma^2 \sqrt{\frac{\mu r}{n}}\pi (\sigma_1^*)^{\frac32} \frac{\kappa \xi^{0,\sf{loo}}}{\underline c_\psi \sqrt{\pi \sigma_r^*}}
    \label{eq: bound for coupled Hessian remainders}
\end{align}
holds with probability at least $1-O(d^{-c-1})$. 

Putting the above components together, the bound for the second term in \eqref{eq: taylor expansion of loss function with respect to X after rearranging} is then completed by
\begin{align}
    & \Big\| - \big[\nabla_{\mb X_i}^2 \widetilde L([\bo \zeta^0, \mb X^*, \mb Y^0 \mb O^0])\big]^{-1} \\
    & \cdot \int_0^1 \Big(\nabla^2_{\mb X_i} \widetilde L\big([\bo \zeta^0, (1-t)\mb X^* + t \mb X, \mb Y^0 \mb O^0] \big) - \nabla^2_{\mb X_i} \widetilde L([\bo \zeta^0, \mb X^*, \mb Y^0 \mb O^0]) \Big)\mathrm d t \cdot  (\mb X\\
    &\quad - \mb X^*)_{i}  \Big\|\\
    & \lesssim   \frac{1}{\underline c_\psi  \pi\sigma_r^*}  \Big[\sigma^2 \sqrt{\frac{\mu r}{p}} \pi (\sigma_1^*)^{3/2}  \norm{\mb X_{i}  - {\mb X^*_{i}} }_2^2\\
        & \qquad + \big( \sigma^2 \frac{\kappa \xi^{0,\sf{loo}}}{\underline c_\psi \sqrt{\pi \sigma_r^*}} \sqrt{\pi \sigma_1^*}
    + \sigma^2 \sqrt{\frac{\mu r}{n}} (\sigma_1^*)^{\frac12} \frac{\kappa \xi^{0,\sf{loo}}}{\underline c_\psi\sqrt{\pi \sigma_r^*}} \pi \sigma_1^* \big)\norm{\mb X_i - \mb X^*_i}_2
        \Big] \eqqcolon L_2
        \label{eq: upper bound on the remainder term wrt X}
\end{align}
with probability at least $1-O(d^{-c-1})$. 
\end{itemize}

\bigskip 

\emph{Applying Brouwer's fixed point theorem. } 
Equipped with the characterization of the terms in \eqref{eq: taylor expansion of loss function with respect to X after rearranging}, the last step amounts to applying Brouwer's fixed point theorem in conjunction with strong convexity. For $\mb x\in\bb R^r$, let $\mb X^{[i]}(\mb x)$ equal $\mb X^*$ with its $i$th row replaced by $\mb x\t$, and write $\mb X_t^{[i]}(\mb x)\coloneqq\mb X^{[i]}((1-t)\mb X_i^*+t\mb x)$. Set $\mb H_i^*\coloneqq\nabla_{\mb X_i}^2\widetilde L([\bo\zeta^0,\mb X^*,\mb Y^0\mb O^0])$. On the $r$-dimensional ball of radius $\rho_1\coloneqq C_1L_1$, define
\begin{align*}
    F_{1,i}(\mb x)
    \coloneqq{}& \mb X_i^*-(\mb H_i^*)^{-1}
    \big\{\nabla_{\mb X_i}\widetilde L([\bo\zeta^0,\mb X^*,\mb Y^0\mb O^0])\big\}\t\\
    &-(\mb H_i^*)^{-1}\int_0^1
    \Big\{\nabla_{\mb X_i}^2\widetilde L([\bo\zeta^0,\mb X_t^{[i]}(\mb x),\mb Y^0\mb O^0])
    -\mb H_i^*\Big\}\,\mathrm dt\,(\mb x-\mb X_i^*).
    \label{eq: definition of F_1,i}
\end{align*}

Invoking the upper bounds \eqref{eq: upper bound on the first term wrt X} and \eqref{eq: upper bound on the remainder term wrt X}, we conclude that, for a fixed $i$ and an appropriately chosen $C_1$, the image of $F_{1,i}(\cdot)$ is contained within its domain 
\eq{
\mc B_i \coloneqq \Big\{\mb x: \norm{\mb x - \mb X_i^*}_2 \leq C_1 L_1 
 \Big\}\label{eq: domain of F1i}
}
with probability at least $1-O(d^{-c-1})$. Indeed, under this domain,
\begin{align}
    & L_2 \lesssim \Big[\kappa_{\psi} \kappa  \sqrt{\frac{\mu r}{p}} (\sigma_1^*)^{1/2}  L_1  + \big( \kappa_\psi \frac{\kappa \xi^{0,\sf{loo}}}{\underline c_\psi\sqrt{\pi \sigma_r^*}} \sqrt{\frac{\kappa}{\pi \sigma_r^*}}
    + \kappa_\psi \kappa \sqrt{\frac{\mu r}{n}} (\sigma_1^*)^{\frac12} \frac{\kappa \xi^{0,\sf{loo}}}{\underline c_\psi \sqrt{\pi \sigma_r^*}}  \big)
        \Big] L_1  \ll  L_1
\end{align}
holds with probability at least $1-O(d^{-c-1})$. The third condition in \eqref{eq: general unilateral refinement smallness conditions} controls the two terms in parentheses in the preceding display. Moreover, the definition of $\xi^1$ and the last condition in \eqref{eq: general unilateral refinement smallness conditions} imply
\begin{align}
& \kappa_\psi \kappa \sqrt{\frac{\mu r}{p}} (\sigma_1^*)^{\frac12}L_1
\lesssim \frac{\kappa_\psi\kappa^{\frac32}\sqrt{\mu r}\,\xi^1}{\underline c_\psi\sqrt{\pi p}}
\ll 1.
\end{align}

A union bound over $i\in[n]$ now shows that, with probability at least $1-O(d^{-c})$,
\[
    \sup_{\mb x\in\mathcal B_i}
    \bignorm{F_{1,i}(\mb x)-\mb X_i^*}_2
    \leq \frac{C_1L_1}{2},
\]
where $\mathcal B_i$ denotes the domain in
\eqref{eq: domain of F1i}. Hence,
$F_{1,i}(\mathcal B_i)\subseteq\mathcal B_i$, and Brouwer's
fixed-point theorem guarantees the existence of a fixed point in
$\mathcal B_i$. Since the corresponding rowwise objective is convex
and is strongly convex on $\mathcal B_i$, this fixed point is its
unique minimizer. By rotational equivariance, the fixed point is
$(\hat{\mb X}^{\ur}\mb O^0)_i$, where $\hat{\mb X}^{\ur}$ solves
the optimization problem
\[
    \min_{\mb X\in\bb R^{n\times r}}
    \sum_{(i,j)\in\Omega}
    \left[
        -R_{ij}\bigl(\zeta_j^0+\mb X_i^\top\mb Y_j^0\bigr)
        +\Psi\bigl(\zeta_j^0+\mb X_i^\top\mb Y_j^0\bigr)
    \right]
\]
in Algorithm~\ref{alg: OS}. We therefore obtain
\begin{align}
    \dist\ti(\hat{\mb X}^{\ur},\mb X^*)
    &\leq
    \bignorm{\hat{\mb X}^{\ur}\mb O^0-\mb X^*}\ti 
    =
    \max_{i\in[n]}
    \bignorm{
        (\hat{\mb X}^{\ur}\mb O^0)_{i,\cdot}
        -\mb X^*_{i,\cdot}
    }_2
    \lesssim L_1
\end{align}
with probability at least $1-O(d^{-c})$.

\paragraph{\texorpdfstring{Update the intercept $\bo \zeta$ and the right factor $\mb Y$}{Update the intercept and the right factor Y}.}
The analysis of $\hat{\bo \zeta}^{\ur}$ and $\hat{\mb Y}^{\ur}$ follows largely the same arguments used to bound
$\dist_{\sf ti}(\hat{\mb X}^{\ur},\mb X^*)$, with the only additional ingredient being a change-of-variables argument that incorporates the intercept parameter $\bo \zeta^*$. For brevity, we focus on explaining how this change of variables reduces the problem to a setting analogous to the preceding bound for $\mb X$. 

We first rewrite the loss function in an augmented form without placing the
possibly diverging baseline $\zeta_0$ in an augmented factor. For
$\bo\zeta^{\sf c}\coloneqq\bo\zeta-\zeta_0\mb 1_p$, define
\begin{align}
    \check L_{\zeta_0}(\mb X^{\sf{app}},\mb Y^{\sf{app},\sf c})
    \coloneqq
    \sum_{(i,j)\in\Omega}
    \Big[
    -R_{i,j}\big(\zeta_0+\mb X^{\sf{app}\top}_{i}\mb Y^{\sf{app},\sf c}_{j}\big)
    +\Psi\big(\zeta_0+\mb X^{\sf{app}\top}_{i}\mb Y^{\sf{app},\sf c}_{j}\big)
    \Big],
\end{align}
where, for $\bo \zeta \in \bb R^p$, $\mb X\in \bb R^{n\times r}$, and $\mb Y\in \bb R^{p\times r}$,
\begin{align}
    \mb X^{\sf{app}}
    \coloneqq
    \Big(
    \sqrt{\frac{\sigma_r^*}{ n}}\mb 1_n,\,
    \mb X
    \Big)
    \in \bb R^{n\times (r+1)},
    \qquad
    \mb Y^{\sf{app},\sf c}
    \coloneqq
    \Big(
    \sqrt{\frac{n}{\sigma_r^*}}\bo \zeta^{\sf c},\,
    \mb Y
    \Big)
    \in \bb R^{p\times (r+1)} .
    \label{eq: tilde L,X,Y}
\end{align}
Analogously to \eqref{eq: gradient of loss function with respect to X} and
\eqref{eq: hessian of loss function with respect to X}, the gradient and Hessian with respect to $\mb Y^{\sf{app},\sf c}$ take the forms
\begin{align}
    & \nabla_{\mb Y^{\sf{app},\sf c}}\check L_{\zeta_0}
     = 
    \big[
    \big(
    \psi(\zeta_0\mb 1_n\mb 1_p^\top
    +\mb X^{\sf{app}}\mb Y^{\sf{app},\sf c\top})-\mb R
    \big)\circ \bo \Omega
    \big]^\top
    \mb X^{\sf{app}},
    \label{eq: forms of derivatives (tilde L)}\\
    & \nabla^2_{\mb Y^{\sf{app},\sf c}}
    \check L_{\zeta_0}(\mb X^{\sf{app}},\mb Y^{\sf{app},\sf c})
     =
    \diag\Big[
    \mb X^{\sf{app}\top}
    \diag\big(
    \psi'(\zeta_0+\mb X^{\sf{app}\top}_{i}\mb Y^{\sf{app},\sf c}_{j})
    \ind\{(i,\\
    &\quad j)\in\Omega\}
    \big)_{i\in[n]}
    \mb X^{\sf{app}}
    \Big]_{j\in[p]}
    \in \bb R^{p(r+1)\times p(r+1)} .
    \label{eq: forms of Hessian (tilde L)}
\end{align}

Next, the rotations for the augmented variables are defined by embedding the usual $r$-dimensional rotations into the lower-right block:
\begin{align*}
     \mb O^{\sf{app},0}
    &\coloneqq
    \begin{pmatrix}
        1 & \mb 0 \\
        \mb 0 & \mb O^0
    \end{pmatrix},
    \qquad 
     \mb O^{\sf{app},0,(-j)}
    \coloneqq
    \begin{pmatrix}
        1 & \mb 0 \\
        \mb 0 &  \mb O^{0,(-j)}
    \end{pmatrix} \quad \text{for}~j\in[p]. 
\end{align*}

With this augmentation, the analysis of $\hat{\bo \zeta}^{\ur}$ and $\hat{\mb Y}^{\ur}$ can be viewed as the analogue of Step~1 after switching the roles of $\mb X$ and $\mb Y$. The main difference in the leave-one-out argument is that we now remove the contribution of a single column noise vector, rather than that of a single row noise vector. Apart from this change, the same perturbation and concentration arguments apply to the augmented variables $(\mb X^{\sf{app}},\mb Y^{\sf{app},\sf c})$. Define
\begin{align*}
    \hat{\mb Y}^{\sf{app},\ur,\sf c}
    &\coloneqq\Big(
    \sqrt{\frac{n}{\sigma_r^*}}
    (\hat{\bo\zeta}^{\ur}-\zeta_0\mb 1_p),\,
    \hat{\mb Y}^{\ur}\Big),\\
    \mb Y^{\sf{app},*,\sf c}
    &\coloneqq\Big(
    \sqrt{\frac{n}{\sigma_r^*}}
    (\bo\zeta^*-\zeta_0\mb 1_p),\,
    \mb Y^*\Big).
\end{align*}
For each fixed $j$, the preceding arguments are invoked with failure probability $O(d^{-c-1})$; a union bound over $j\in[p]$ therefore yields the following uniform bound with failure probability $O(d^{-c})$:
\begin{align}
    &\sqrt{\frac{n}{\sigma_r^*}}
    \bignorm{\hat{\bo\zeta}^{\ur}-\bo\zeta^*}_\infty
    \vee\bignorm{\hat{\mb Y}^{\ur}\mb O^0-\mb Y^*}\ti\nonumber\\
    &\qquad\leq
    \bignorm{\hat{\mb Y}^{\sf{app},\ur,\sf c}\mb O^{\sf{app},0}
    -\mb Y^{\sf{app},*,\sf c}}\ti
    \lesssim\frac{\xi^1}
    {\underline c_\psi\sqrt{\pi\sigma_r^*}}
\end{align}
with probability at least $1- O(d^{-c})$. 

\subsubsection{Leave-one-out Error Controls for Algorithm~\ref{alg: OS}.}
We now control the difference between the outputs of Algorithm~\ref{alg: OS} when applied to the original initialization and to its leave-one-out counterpart (cf.~\eqref{eq: leave-one-out error for unilateral refinement}). We focus primarily on the leave-one-row-out variant $\hat{\mb X}^{\ur,(i)}$ as a canonical example, while shedding light on the similarity and distinction for treating $\hat{\mb X}^{\ur,(-j)}$.  
In addition, the leave-one-out analysis with respect to $\bo\zeta$ and $\mb Y$ follows similarly from the change of variables discussed in \eqref{eq: tilde L,X,Y}. Recall that $[\hat{\bo \zeta}^{\ur,(i)},\hat{\mb X}^{\ur,(i)}, \hat{\mb Y}^{\ur,(i)}]$ is obtained via solving 
\begin{align}
    & \hat{\mb X}^{\ur,(i)} \coloneqq \arg\min_{\mb X \in \bb R^{n \times r}} \widetilde L^{(i)}([\hat{\bo \zeta}^{0,(i)}, \mb X, \hat{\mb Y}^{0,(i)}]), \\ 
    & (\hat{\bo \zeta}^{\ur,(i)}, \hat{\mb Y}^{\ur,(i)}) \coloneqq  \arg\min_{(\bo \zeta, \mb Y) \in \bb R^p \times \bb R^{p \times r}} \widetilde L^{(i)}([\bo \zeta, \hat{\mb X}^{0,(i)}, \mb Y]), 
\end{align}
where the objective function $\widetilde L^{(i)}$ is defined as 
\begin{align}
    & \widetilde L^{(i)}([\bo \zeta, \mb X, \mb Y]) \coloneqq \sum_{(i',j)\in \Omega, i'\neq i} \Big[- R_{i',j}(\zeta_{j} + \mb X_{i'}\t \mb Y_j) + \Psi( \zeta_j + \mb X_{i'}\t \mb Y_j)\Big]\\ 
    & \qquad + \pi \sum_{j\in[p]} \Big[- \psi(\zeta_j^* +{ \mb X^*_i}\t \mb Y^*_j)(\zeta_j + \mb X_i\t \mb Y_j) + \Psi(\zeta_j + \mb X_i\t \mb Y_j)\Big] . 
\end{align}

The Brouwer fixed point theorem again plays an important role in proving the proximity between
$[\hat{\bo \zeta}^{\ur}, \hat{\mb X}^{\ur}, \hat{\mb Y}^{\ur}]$ and
$[\hat{\bo \zeta}^{\ur,(i)}, \hat{\mb X}^{\ur,(i)}, \hat{\mb Y}^{\ur,(i)}]$.
The key difference is that the Taylor expansion is now taken w.r.t. $\widetilde L^{(i)}$ around
$[\hat{\bo \zeta}^{\ur}, \hat{\mb X}^{\ur}, \hat{\mb Y}^{\ur}]$ rather than w.r.t. $\widetilde L$ around the ground truth.

We express the stationary condition for $\hat{\mb X}^{\ur,(i)}$ analogous to
\eqref{eq: taylor expansion of loss function with respect to X after rearranging}, but in a joint manner. 
Specifically, the condition can be expressed as 
\begin{align}
    & \hat{\mb X}^{\ur, (i)} \mb O^{0,(i)}  = \mb G_i(\hat{\mb X}^{\ur, (i)} \mb O^{0,(i)}) 
\end{align}
where, for each $i'\in[n]$, the $i'$-th slice of length $r$ of the matrix-valued function $\mb G_i(\mb X): \bb R^{n \times r} \rightarrow \bb R^{n\times r}$ is defined below. Write $\mb X_t(\mb X)\coloneqq(1-t)\hat{\mb X}^{\ur}+t\mb X$.
\begin{align}
    \mb H_{i,i'}^0
    \coloneqq{}&\nabla_{\mb X_{i'}}^2\widetilde L^{(i)}
    ([\hat{\bo\zeta}^{0,(i)},\hat{\mb X}^{\ur},
    \hat{\mb Y}^{0,(i)}\mb O^{0,(i)}]),
    \notag\\
    \big[\mb G_i(\mb X)\big]_{i'}
    \coloneqq{}&\hat{\mb X}^{\ur}_{i'}-(\mb H_{i,i'}^0)^{-1}
    \Big\{
    \nabla_{\mb X_{i'}}
    \widetilde L^{(i)}([\hat{\bo\zeta}^{0,(i)},\hat{\mb X}^{\ur},
    \hat{\mb Y}^{0,(i)}\mb O^{0,(i)}])\t
    \notag\\
    &\quad+\int_0^1
    \Big[
    \nabla_{\mb X_{i'}}^2
    \widetilde L^{(i)}([\hat{\bo\zeta}^{0,(i)},\mb X_t(\mb X),
    \hat{\mb Y}^{0,(i)}\mb O^{0,(i)}])
    -\mb H_{i,i'}^0
    \Big] \,\mathrm dt\,
    (\mb X_{i'}-\hat{\mb X}_{i'}^{\ur})
    \Big\}.
    \label{eq: definition of G}
\end{align}
The leave-one-column-out counterpart $\mb G_{-j}(\mb X)$ is defined analogously by replacing $\widetilde L^{(i)}$, $\hat{\bo \zeta}^{0,(i)}$, and $\hat{\mb Y}^{0,(i)}\mb O^{0,(i)}$ with $\widetilde L^{(-j)}$, $\hat{\bo \zeta}^{0,(-j)}$, and $\hat{\mb Y}^{0,(-j)}\mb O^{0,(-j)}$, respectively.

Then it suffices to control the magnitude of the terms in the function within a certain domain, which are expanded in the following: 
\begin{itemize}
    \item Invoking the stationary condition of $\hat {\mb X}^{\ur}$ implies that 
    \begin{align}
        & \norm{ \nabla_{\mb X}
    \widetilde L^{(i)}([\hat{\bo \zeta}^{0,(i)}, \hat{\mb X}^{\ur}, \hat{\mb Y}^{0,(i)}\mb O^{0,(i)}])}\fb \\ 
    \leq &\norm{ \nabla_{\mb X}
    \widetilde L([\hat{\bo \zeta}^{0,(i)}, \hat{\mb X}^{\ur}, \hat{\mb Y}^{0,(i)} \mb O^{0,(i)}]) - \nabla_{\mb X}
    \widetilde L([\hat{\bo \zeta}^{0}, \hat{\mb X}^{\ur}, \hat{\mb Y}^{0}])}\fb \\ 
    + & \norm{\nabla_{\mb X}(\widetilde L - \widetilde L^{(i)}) ([\hat{\bo \zeta}^{0,(i)}, \hat{\mb X}^{\ur}, \hat{\mb Y}^{0,(i)} \mb O^{0,(i)}])}\fb. 
    \label{eq: decomposition of nabla L^(i)}
    \end{align}
    Regarding the first term on the right-hand side of \eqref{eq: decomposition of nabla L^(i)}, it is seen that 
    \begin{align}
        & \norm{ \nabla_{\mb X}
    \widetilde L([\hat{\bo \zeta}^{0,(i)}, \hat{\mb X}^{\ur}, \hat{\mb Y}^{0,(i)}\mb O^{0,(i)}]) - \nabla_{\mb X}
    \widetilde L([\hat{\bo \zeta}^{0}, \hat{\mb X}^{\ur}, \hat{\mb Y}^{0}])}\fb \\ 
    \leq & \norm{\big[\big(\psi(\mb 1_n \hat{\bo \zeta}^{0,(i)\top} + \hat{\mb X}^{\ur} \mb O^{0,(i)\top}\hat{\mb Y}^{0,(i)\top}) -  \psi(\mb 1_n \hat{\bo \zeta}^{0\top} + \hat{\mb X}^{\ur} \hat{\mb Y}^{0\top})\big) \circ \bo \Omega \big] \hat{\mb Y}^{0,(i)}}\fb \\ 
    & + \norm{\mc P_{\Omega}\big( \mb R - \psi(\mb 1_n\hat {\bo \zeta}^{0\top} + \hat{\mb X}^{\ur} \hat{\mb Y}^{0\top}) \big) \big(\hat{\mb Y}^{0,(i)}\mb O^{0,(i)} - \hat{\mb Y}^0\big)}\fb  \\ 
    \leq & \underbrace{ \norm{\mc P_{\Omega}\big(\psi(\mb 1_n \hat{\bo \zeta}^{0,(i)\top} + \hat{\mb X}^{\ur}\mb O^{0,(i)\top}\hat{\mb Y}^{0,(i)\top}) -  \psi(\mb 1_n \hat{\bo \zeta}^{0\top} + \hat{\mb X}^{\ur} \hat{\mb Y}^{0\top})\big) \hat{\mb Y}^{0,(i)}}\fb }_{\beta_1} \\ 
    & + \underbrace{\norm{\mc P_{\Omega}\big(\psi(\mb 1_n \bo \zeta^{*\top} + \mb X^*{\mb Y^*}\t)  - \psi(\mb 1_n \hat{\bo \zeta}^{0\top} + \hat{\mb X}^{\ur} \hat{\mb Y}^{0\top}) \big) \big(\hat{\mb Y}^{0,(i)} \mb O^{0,(i)}  - \hat{\mb Y}^0 \big) }\fb}_{\beta_2}  \\ 
    & + \underbrace{\norm{\mc P_{\Omega}(\mb E)} \bignorm{\hat{\mb Y}^{0,(i)} \mb O^{0,(i)} - \hat{\mb Y}^0}\fb}_{\beta_3}. 
    \label{eq: decomposition of leave-one-out gradient difference}
    \end{align}
    The terms $\beta_i$, $i=1,2,3$ above can be reasoned as follows: 
    \begin{itemize}
    \item For the first term $\beta_1$, we invoke Lemma~\ref{lemma: generalized version of lemma 4.4 in candes and recht} together with Cauchy's inequality to derive that 
    \begin{align*}
    & \beta_1\\
    &\quad \leq \norm{\mc P_{ \Omega}\big(\psi(\mb 1_n \hat{\bo \zeta}^{0,(i)\top} + \hat{\mb X}^{\ur} \mb O^{0,(i)\top}\hat{\mb Y}^{0,(i)\top}) -  \psi(\mb 1_n \hat{\bo \zeta}^{0\top} + \hat{\mb X}^{\ur} \hat{\mb Y}^{0\top})\big) }\fb\\
    &\qquad \max_{i'\in[n]}\bignorm{ \diag(\bo \Omega_{i',\cdot}) \hat{\mb Y}^{0,(i)}}\\
    & \lesssim   \sigma^2 \sqrt{\pi \sigma_1^*} \frac{\xi^{0,\sf{loo}}}{\underline c_\psi \sqrt{\pi \sigma_r^*}} \sqrt{\pi \sigma_1^*} = \sigma^2 \pi \sigma_1^* \frac{\xi^{0,\sf{loo}}}{\underline c_\psi \sqrt{\pi \sigma_r^*}}
    \end{align*}
    holds, for each fixed leave-one-out index $i$, with probability at least $1-O(d^{-c-1})$. Here, the bound for $\bignorm{ \diag(\bo \Omega_{i',\cdot}) \hat{\mb Y}^{0,(i)}}$ follows from \eqref{eq: upper bound on Omega_i Y^{0,(i)}} if $i' = i$ and from 
    \begin{align}
        & \bignorm{ \diag(\bo \Omega_{i',\cdot}) \hat{\mb Y}^{0,(i)}}_2 \leq \bignorm{ \diag(\bo \Omega_{i',\cdot}) \hat{\mb Y}^{0,(i')}} + 2 \max_{i_1\in[n]} \bignorm{\hat{\mb Y}^{0,(i_1)} \mb O^{0,(i_1)} - \hat{\mb Y}^0} \\
        \lesssim &\sqrt{\pi \sigma_1^*} + \frac{\xi^{0,\sf{loo}}}{\underline c_\psi \sqrt{\pi \sigma_r^*}} \lesssim  \sqrt{\pi \sigma_1^*}
        \label{eq: upper bound for beta1 (i)}
    \end{align}
    with probability at least $1-O(d^{-c-2})$ for each fixed pair $(i,i')$ by \eqref{eq: upper bound on Omega_i Y^{0,(i)}} and the condition on the leave-one-out error if $i' \neq i$. 
        \item To upper bound $\beta_2$, we resort to Lemma~\ref{lemma: lemma 4.4 chen2019} together with a variational form. By the mean-value theorem, one has 
        \begin{align}
        &  \psi(\mb 1_n \bo \zeta^{*\top} + \mb X^*{\mb Y^*}\t)  - \psi(\mb 1_n \hat{\bo \zeta}^{0\top} + \hat{\mb X}^{\ur} \hat{\mb Y}^{0\top})\\ 
        = & \underbrace{\int_0^1 \psi'((1 - t)(\mb 1_n \bo \zeta^{*\top} + \mb X^*{\mb Y^*}\t) + t(\mb 1_n \hat{\bo \zeta}^{0\top} + \hat{\mb X}^{\ur} \hat{\mb Y}^{0\top} )) \mathrm d t}_{\eqqcolon \mb S_\psi} \\ 
        & \circ \underbrace{\Big[\mb 1_n(\bo\zeta^*-\hat{\bo\zeta}^{0})^\top
+\mb X^*
    \big(\mb Y^*-\hat{\mb Y}^{0}\mb O^0\big)^\top +
\big(\mb X^*-\hat{\mb X}^{\ur}\mb O^0\big)
    \big(\hat{\mb Y}^{0}\mb O^0\big)^\top \Big] }_{\eqqcolon \mb D_\psi} 
        \eqqcolon  \bo \Delta_\psi.  
        \end{align}
        For $t\in[0,1]$, let $\mb S_{\psi,t}$ denote the integrand in
        $\mb S_\psi$. The row bounds above and the domain verification at
        the start of this proof apply to both $\bo\theta^*$ and
        $[\hat{\bo\zeta}^0,\hat{\mb X}^{\ur},\hat{\mb Y}^0]$.
        Their endpoint predictors and every natural-parameter interpolant
        therefore lie in $\mc D_\psi$, so
        $\mb S_{\psi,t}\in\mc M_\psi$ for every $t\in[0,1]$.
        Lemma~\ref{lemma: injectivity of P_Omega under nonlinearity},
        linearity, and the triangle inequality give
        \begin{align}
            \bignorm{\mc P_\Omega(\mb S_\psi)-\pi\mb S_\psi}
            &\leq \int_0^1
            \bignorm{\mc P_\Omega(\mb S_{\psi,t})-\pi\mb S_{\psi,t}}\,\mathrm dt
            \lesssim \sigma^2\xi^4\sqrt{\pi d}.
            \label{eq: uniform nonlinear sampling operator bound}
        \end{align}
        And applying the triangle inequality yields that 
        $$
        \beta_2 \leq  \norm{\big(\mc P_{\Omega}(\bo \Delta_\psi) - \pi \bo \Delta_\psi\big) \big(\hat{\mb Y}^{0,(i)} \mb O^{0,(i)}  - \hat{\mb Y}^0 \big) }\fb + \pi \bignorm{\bo \Delta_\psi  \big(\hat{\mb Y}^{0,(i)} \mb O^{0,(i)}  - \hat{\mb Y}^0 \big)}\fb . 
        $$
        Regarding the first term on the right-hand side above, we invoke the
        variational characterization of the Frobenius norm,
        Lemma~\ref{lemma: lemma 4.4 chen2019}, and the
        Lemma~\ref{lemma: injectivity of P_Omega under nonlinearity} bound
        \eqref{eq: uniform nonlinear sampling operator bound} to obtain
        \begin{align*}
    & \norm{\big(\mc P_{\Omega}(\bo \Delta_\psi) - \pi \bo \Delta_\psi\big) \big(\hat{\mb Y}^{0,(i)} \mb O^{0,(i)}  - \hat{\mb Y}^0 \big) }\fb\\
    & =   \max_{\mb A \in \bb R^{n \times r}: \norm{\mb A}\fb = 1} \ip{\bo \Delta_\psi, \mc P_{\Omega}\big(\mb A \big(\hat{\mb Y}^{0,(i)} \mb O^{0,(i)}  - \hat{\mb Y}^0 \big)\t \big) - \pi\mb A \big(\hat{\mb Y}^{0,(i)} \mb O^{0,(i)}  - \hat{\mb Y}^0 \big)\t}\\
    & \stackrel{\text{by Lemma~\ref{lemma: lemma 4.4 chen2019}}}{\leq}   \norm{\mc P_\Omega(\mb S_\psi) - \pi \mb S_\psi} \Big[\bignorm{\bo \zeta^*  - \hat{\bo \zeta}^0 }_\infty  + \norm{\mb X^*}\ti \bignorm{\mb Y^*  - \hat{\mb Y}^0 \mb O^0 }\ti\\
            & \qquad + \bignorm{\hat{\mb X}^{\ur} \mb O^0 - \mb X^*}\ti \bignorm{\hat{\mb Y}^0}\ti  \Big]   \bignorm{\hat{\mb Y}^{0,(i)} \mb O^{0,(i)}  - \hat{\mb Y}^0  }\fb \\ 
    & \lesssim   \sigma^2 \sqrt{\pi d} \xi^4 \Big[\sqrt{\frac{\mu^2 r^2}{np}} \sigma_1^*  \Big] \frac{\xi^{0,\sf{loo}}}{\underline c_\psi \sqrt{\pi \sigma_r^*}}\\
    & \lesssim   \sigma^2 \pi \sigma_1^* \frac{\xi^{0,\sf{loo}}}{\underline c_\psi \sqrt{\pi \sigma_r^*}}
        \end{align*}
        holds, for each fixed $i$, with probability at least $1-O(d^{-c-1})$, since $\mu^2 r^2 \xi^8 \ll \pi (n\wedge p)$. 

        On the other hand, for each fixed $i$, with probability at least $1-O(d^{-c-1})$, 
        \begin{align}
            & \pi \norm{\bo \Delta_\psi  \big(\hat{\mb Y}^{0,(i)} \mb O^{0,(i)}  - \hat{\mb Y}^0 \big)}\fb  \lesssim \sigma^2 \pi \sigma_1^* \frac{\xi^{0,\sf{loo}}}{\underline c_\psi \sqrt{\pi \sigma_r^*}},
        \end{align}
        which leads to 
        \begin{align}
            & \beta_2 \lesssim   \sigma^2 \pi \sigma_1^* \frac{\xi^{0,\sf{loo}}}{\underline c_\psi \sqrt{\pi \sigma_r^*}}.
            \label{eq: upper bound for beta2}
        \end{align}

        \item Lastly, we have for $\beta_3$, for each fixed $i$, with probability at least $1-O(d^{-c-1})$, 
        \begin{align}
            & \beta_3 \lesssim \sigma \sqrt{\pi d} \frac{\xi^{0,\sf{loo}}}{\underline c_\psi \sqrt{\pi \sigma_r^*}} \ll \sigma^2 \pi \sigma_1^* \frac{\xi^{0,\sf{loo}}}{\underline c_\psi \sqrt{\pi \sigma_r^*}},
        \end{align}
        by Lemma~\ref{lemma: noise spectral norm concentration}.
    \end{itemize}  
    Substituting the above into \eqref{eq: decomposition of leave-one-out gradient difference} yields, for each fixed $i$, with probability at least $1-O(d^{-c-1})$, 
    \begin{align}
        &  \norm{ \nabla_{\mb X}
    \widetilde L([\hat{\bo \zeta}^{0,(i)}, \hat{\mb X}^{\ur}, \hat{\mb Y}^{0,(i)}\mb O^{0,(i)} ]) - \nabla_{\mb X}
    \widetilde L([\hat{\bo \zeta}^{0}, \hat{\mb X}^{\ur}, \hat{\mb Y}^{0}])}\fb 
    \lesssim  \sigma^2 \pi \sigma_1^* \frac{\xi^{0,\sf{loo}}}{\underline c_\psi \sqrt{\pi \sigma_r^*}} .
    \label{eq: bound for nabla L difference}
    \end{align}
    A similar argument also leads to the same upper bound
    \begin{align}
    \norm{ \nabla_{\mb X}
    \widetilde L([\hat{\bo \zeta}^{0,(-j)}, \hat{\mb X}^{\ur}, \hat{\mb Y}^{0,(-j)} \mb O^{0,(-j)}]) - \nabla_{\mb X}
    \widetilde L([\hat{\bo \zeta}^{0}, \hat{\mb X}^{\ur}, \hat{\mb Y}^{0}])}\fb  
    \lesssim  \sigma^2 \pi \sigma_1^* \frac{\xi^{0,\sf{loo}}}{\underline c_\psi \sqrt{\pi \sigma_r^*}} .
    \label{eq: eq: bound for nabla L difference (column)}
    \end{align}
    
    On the other hand, by construction, the objectives $\widetilde L$ and
$\widetilde L^{(i)}$ differ only in their contributions from the $i$th
row. Consequently,
\[
    \nabla_{\mb X_{i'}}
    \big(\widetilde L-\widetilde L^{(i)}\big)
    ([\hat{\bo \zeta}^{0,(i)}, \hat{\mb X}^{\ur}, \hat{\mb Y}^{0,(i)}\mb O^{0,(i)}])
    =\mb 0
    \qquad\text{for every }i'\neq i.
\]
When $i'=i$, adding and subtracting the corresponding leave-one-out
quantities yields
\begin{equation}
\begin{split}
    &\bignorm{
        \nabla_{\mb X_i}
        \big(\widetilde L-\widetilde L^{(i)}\big)
        ([\hat{\bo\zeta}^{0,(i)},\hat{\mb X}^{\ur},\hat{\mb Y}^{0,(i)}\mb O^{0,(i)}])
    }_2
    \lesssim
    \gamma_1+\gamma_2+\gamma_3+\gamma_4 + \gamma_5 .
\end{split}
\label{eq: decomposition of nabla L - L(i)}
\end{equation}

Here,
\begin{align}
    \gamma_1
    &\coloneqq
    \Bignorm{
        \Big[
            \psi\big(
                \hat{\bo\zeta}^{0,(i)\top}
                +\mb X^*_{i,\cdot} \tilde{\mb O}^{0,(i)\top}
                \hat{\mb Y}^{0,(i)\top}
            \big)
            \circ\bo\Omega_{i,\cdot}
            \\ 
            & \qquad -
            \pi\psi\big(
                \hat{\bo\zeta}^{0,(i)\top}
                +\mb X^*_{i,\cdot} \tilde{\mb O}^{0,(i)\top}
                \hat{\mb Y}^{0,(i)\top}
            \big)
        \Big]
        \hat{\mb Y}^{0,(i)}
    }_2,
    \label{eq: definition of gamma1}\\
    \gamma_2
    &\coloneqq \Bignorm{\big[\psi\big(
            \hat{\bo\zeta}^{0,(i)\top}
            +\hat{\mb X}_{i,\cdot}^{\ur} \mb O^{0,(i)\top}\hat{\mb Y}^{0,(i)\top}
        \big)
        -
        \psi\big(
            \hat{\bo\zeta}^{0,(i)\top}
            +\mb X^*_{i,\cdot}\tilde{\mb O}^{0,(i)\top}
            \hat{\mb Y}^{0,(i)\top}
        \big) \big] \circ \bo \Omega_{i,\cdot}}_2 \\ 
        & \qquad \cdot 
        \Big( \bignorm{
            \diag(\bo\Omega_{i,\cdot})
            \hat{\mb Y}^{0,(i)}
        }
        +
        \bignorm{
            \hat{\mb Y}^{0,(i)}\mb O^{0,(i)}
            -\hat{\mb Y}^{0}
        }  \Big) \\ 
        & \qquad + \pi \Bignorm{\psi\big(
             \hat{\bo\zeta}^{0,(i)\top}
            +\hat{\mb X}_{i,\cdot}^{\ur} \mb O^{0,(i)\top}\hat{\mb Y}^{0,(i)\top}
        \big)
        -
        \psi\big(
            \hat{\bo\zeta}^{0,(i)\top}
            +\mb X^*_{i,\cdot}\tilde{\mb O}^{0,(i)\top}
            \hat{\mb Y}^{0,(i)\top}
        \big)}_2\\
        &\qquad \cdot\bignorm{\hat{\mb Y}^{0}} , 
    \label{eq: definition of gamma2}\\
    \gamma_3
    &\coloneqq
    \Bignorm{
        \psi\big(
            \hat{\bo\zeta}^{0,(i)\top}
            +\mb X^*_{i,\cdot} \tilde{\mb O}^{0,(i)\top}
            \hat{\mb Y}^{0,(i)\top}
        \big)
        \circ\bo\Omega_{i,\cdot}
        -
        \pi\psi\big(
            \hat{\bo\zeta}^{0,(i)\top}
            +\mb X^*_{i,\cdot} \tilde{\mb O}^{0,(i)\top}
            \hat{\mb Y}^{0,(i)\top}
        \big)
    }_2
    \nonumber\\
    &\hspace{1cm}\cdot
    \bignorm{
        \hat{\mb Y}^{0,(i)}\mb O^{0,(i)}
        -\hat{\mb Y}^{0}
    },
    \label{eq: definition of gamma3}\\
    \gamma_4
    &\coloneqq
    \bignorm{
        \mb E_{i,\cdot}
        \diag(\bo\Omega_{i,\cdot})
        \hat{\mb Y}^{0}
    }_2,\\
    \gamma_5&\coloneqq \bignorm{\big(\mc P_{\Omega}(\psi(\mb M^*)) - \pi \psi(\mb M^*) \big)_{i,\cdot} \mb Y^{0}}_2. 
    \label{eq: definition of gamma4}
\end{align}

    Invoking the matrix Bernstein inequality, we have for $\gamma_1$, for each fixed $i$, with probability at least $1-O(d^{-c-1})$,  
    \begin{align}
        & \gamma_1 \lesssim \bar \sigma \sqrt{\pi \sigma_1^* \log d} + B \sqrt{\frac{\mu r\sigma_1^*}{p}} \log d\lesssim \bar\sigma \sqrt{\pi \sigma_1^* \log d}.
    \end{align}

    With respect to $\gamma_2$, the Taylor expansion along with \eqref{eq: row-wise error for unilateral refinement}, \eqref{eq: transfer from relative to population-aligned LOO rotations}, and \eqref{eq: upper bound on Omega_i Y^{0,(i)}} first implies, for each fixed $i$, that with probability at least $1-O(d^{-c-1})$ 
    \begin{align}
        & \Bignorm{\big[\psi\big(
            \hat{\bo\zeta}^{0,(i)\top}
            +\hat{\mb X}_{i,\cdot}^{\ur} \mb O^{0,(i)\top}\hat{\mb Y}^{0,(i)\top}
        \big)
        -
        \psi\big(
            \hat{\bo\zeta}^{0,(i)\top}
            +\mb X^*_{i,\cdot}\tilde{\mb O}^{0,(i)\top}
            \hat{\mb Y}^{0,(i)\top}
        \big) \big] \circ \bo \Omega_{i,\cdot}}_2 \\ 
    & \leq   \sigma^2 \norm{\mb X^*_{i,\cdot}\tilde{\mb O}^{0,(i)\top} - \hat{\mb X}^{\ur}_{i,\cdot} \mb O^{0,(i)\top}}_2 \bignorm{\diag(\bo\Omega_{i,\cdot}) \hat{\mb Y}^{0,(i)}}
        \lesssim  \sigma^2  \sqrt{\pi \sigma_1^*} \frac{\sqrt{\kappa \mu r}\xi^1}{\underline c_\psi \sqrt{\pi \sigma_r^*}} , \\
        & \Bignorm{\psi\big(
            \hat{\bo\zeta}^{0,(i)\top}
            +\hat{\mb X}_{i,\cdot}^{\ur} \mb O^{0,(i)\top}\hat{\mb Y}^{0,(i)\top}
        \big)
        -
        \psi\big(
            \hat{\bo\zeta}^{0,(i)\top}
            +\mb X^*_{i,\cdot}\tilde{\mb O}^{0,(i)\top}
            \hat{\mb Y}^{0,(i)\top}
        \big)}_2\\
    &\quad \lesssim \sigma^2 \sqrt{\sigma_1^*} \frac{\sqrt{\mu r \kappa} \xi^1}{\underline c_\psi \sqrt{\pi \sigma_r^*}},
    \end{align}
    where we used the relation from Lemma~\ref{lemma: perturbation theory of optimal rotation under F norm} that 
    \begin{align}
        & \norm{\mb X^*_{i,\cdot}\tilde{\mb O}^{0,(i)\top} - \hat{\mb X}^{\ur}_{i,\cdot} \mb O^{0,(i)\top}}_2 \leq \norm{\mb X^*}\ti \bignorm{\mb O^{0,(i) }\mb O^0 -  \tilde{\mb O}^{0,(i)}}_2 + \bignorm{\mb X^{0} \mb O^{0} - \mb X^*}\ti \\ 
        \leq & \sqrt{\frac{\mu r \sigma_1^*}{n}} \frac{\sqrt{n} \xi^1 }{\underline c_\psi \sqrt{\pi }\sigma_r^* } +  \frac{\xi^1}{\underline c_\psi \sqrt{\pi \sigma_r^*}} \lesssim  \frac{\sqrt{\kappa \mu r}\xi^1}{\underline c_\psi \sqrt{\pi \sigma_r^*}}. 
    \end{align}

    In addition, one has from the initial leave-one-out error condition and \eqref{eq: upper bound on Omega_i Y^{0,(i)}} that, for each fixed $i$, with probability at least $1-O(d^{-c-1})$ 
    \begin{align}
        & \bignorm{
            \diag(\bo\Omega_{i,\cdot})
            \hat{\mb Y}^{0,(i)}
        }
        \vee 
        \bignorm{
            \hat{\mb Y}^{0,(i)}\mb O^{0,(i)}
            -\hat{\mb Y}^{0} } \lesssim \sqrt{\pi \sigma_1^*}. 
    \end{align}
    
    Taken together, these bounds yield, for each fixed $i$, with probability at least $1-O(d^{-c-1})$,
    \begin{align}
        & \gamma_2 \lesssim \sigma^2 \sqrt{\pi \sigma_1^*} \frac{\sqrt{\kappa \mu r}\xi^1}{\underline c_\psi \sqrt{\pi \sigma_r^*}} \sqrt{\pi \sigma_1^*} + \pi \sigma^2 \sqrt{\pi \sigma_1^*} \frac{\sqrt{\kappa \mu r}\xi^1}{\underline c_\psi \sqrt{ \sigma_r^*}} \sqrt{\sigma_1^*}\lesssim  \sigma^2 \pi \sigma_1^* \frac{\sqrt{\kappa \mu r}\xi^1}{\underline c_\psi \sqrt{\pi \sigma_r^*}},
    \end{align}
 where $\xi^1$ was introduced in \eqref{eq: general unilateral refinement rate}.

    Further, applying the matrix Bernstein inequality with respect to $\bo \Omega_{i,\cdot}$ implies that 
    \begin{align}
        & \gamma_3 \lesssim
        \bar \sigma \sqrt{\pi d} \frac{\xi^{0,\sf{loo}}}{\underline c_\psi \sqrt{\pi \sigma_r^*}} 
        \lesssim  \sigma^2 \pi \sigma_1^* \frac{\xi^{0,\sf{loo}}}{\underline c_\psi \sqrt{\pi \sigma_r^*}}
    \end{align}
    holds, for each fixed $i$, with probability at least $1-O(d^{-c-1})$.

    To proceed, we bound $\gamma_4$ using the leave-one-out trick: for each fixed $i$, with probability at least $1-O(d^{-c-1})$, it holds that
    \begin{align*}
    & \gamma_4 \leq \norm{\mb E_{i,\cdot} \diag(\bo\Omega_{i,\cdot}) \mb Y^{0,(i)}} + \norm{\mb E_{i,\cdot} \diag(\bo\Omega_{i,\cdot}) } \bignorm{\mb Y^{0,(i)} \mb O^{0,(i)} - \mb Y^0} \lesssim \sigma \sqrt{\pi r \sigma_1^* \log d}\\
    &\quad + \sigma \sqrt{\pi p} \frac{\xi^{0,\sf{loo}}}{\underline c_\psi \sqrt{\pi \sigma_r^*}}.
    \end{align*}

    Lastly, we leverage the leave-one-out coupling along with the matrix Bernstein inequality to obtain that 
    \begin{align}
        & \gamma_5 \leq \bignorm{\big(\mc P_{\Omega}(\psi(\mb M^*)) - \pi \psi(\mb M^*) \big)_{i,\cdot} \mb Y^{0,(i)}}_2 \\ 
        & \qquad + \bignorm{\big(\mc P_{\Omega}(\psi(\mb M^*)) - \pi \psi(\mb M^*) \big)_{i,\cdot}}_2 \bignorm{\mb Y^{0,(i)} \mb O^{0,(i)} - \mb Y^0} \\ 
        \lesssim & \bar \sigma \sqrt{\pi \sigma_1^* \log d} + B \sqrt{\frac{\mu r \sigma_1^*}{p}} \log d + \bar \sigma \sqrt{\pi p} \frac{\xi^{0,\sf{loo}}}{\underline c_\psi \sqrt{\pi \sigma_r^*}}\lesssim \bar \sigma \sqrt{\pi \sigma_1^* \log d}. 
    \end{align}

    Plugging these terms into \eqref{eq: decomposition of nabla L - L(i)} yields, for each fixed $i$, with probability at least $1-O(d^{-c-1})$,
    \begin{align}
        & \norm{\nabla_{\mb X}(\widetilde L - \widetilde L^{(i)}) ([\hat{\bo \zeta}^{0}, \hat{\mb X}^{\ur}, \hat{\mb Y}^{0}])}\fb \lesssim  \sigma^2 \pi \sigma_1^* \frac{\xi^{1}}{\underline c_\psi \sqrt{\pi \sigma_r^*}}.
    \end{align}

    For the leave-one-column-out counterpart, direct algebra gives
    \begin{align*}
        & \nabla_{\mb X}(\widetilde L-\widetilde L^{(-j)})
        ([\hat{\bo\zeta}^{0,(-j)},\hat{\mb X}^{\ur},
        \hat{\mb Y}^{0,(-j)}\mb O^{0,(-j)}]) \\
        ={}& \Big[(\bo\Omega_{\cdot,j}-\pi\mb 1_n)\circ
        \psi\big(\mb 1_n\hat\zeta_j^{0,(-j)}
        +\hat{\mb X}^{\ur}\mb O^{0,(-j)\top}
        \hat{\mb Y}_j^{0,(-j)}\big) \\
        &\qquad-\big\{\mc P_\Omega(\mb R)_{\cdot,j}
        -\pi\psi(\mb M^*)_{\cdot,j}\big\}\Big]
        (\hat{\mb Y}^{0,(-j)}\mb O^{0,(-j)})_{j,\cdot}.
    \end{align*}
    The rowwise bounds ensure that the predictors above lie in $\mc D_\psi$. Hence, for each fixed $j$, with probability at least $1-O(d^{-c-1})$,
    \begin{align*}
        & \norm{\nabla_{\mb X}(\widetilde L-\widetilde L^{(-j)})
        ([\hat{\bo\zeta}^{0,(-j)},\hat{\mb X}^{\ur},
        \hat{\mb Y}^{0,(-j)}\mb O^{0,(-j)}])}\fb \\
        \leq{}& \Big[
        \bignorm{\mc P_\Omega(\mb R)_{\cdot,j}
        -\pi\psi(\mb M^*)_{\cdot,j}}_2
        +\sigma^2\bignorm{\bo\Omega_{\cdot,j}-\pi\mb 1_n}_2
        \Big]\bignorm{\hat{\mb Y}^{0,(-j)}}\ti \\
        \lesssim{}& \bar\sigma
        \sqrt{\frac{\pi n\mu r\sigma_1^*}{p}}
        \lesssim \sigma^2\pi\sigma_1^*
        \frac{\sqrt{\kappa\mu r\rho}\,\xi^1}
        {\underline c_\psi\sqrt{\pi\sigma_r^*}}.
    \end{align*}

    A union bound over the leave-one-out index, together with \eqref{eq: bound for nabla L difference} and \eqref{eq: eq: bound for nabla L difference (column)}, implies that, with probability at least $1-O(d^{-c})$, 
    \begin{align}
        &  \norm{ \nabla_{\mb X}
    \widetilde L^{(i)}([\hat{\bo \zeta}^{0,(i)}, \hat{\mb X}^{\ur}, \hat{\mb Y}^{0,(i)}\mb O^{0,(i)}])}\fb \lesssim \sigma^2 \pi \sigma_1^* \frac{\sqrt{\kappa \mu r \rho}\xi^{1}}{\underline c_\psi \sqrt{\pi \sigma_r^*}}.
    \label{eq: leave-one-out gradient norm upper bound}, \\ 
    & \norm{ \nabla_{\mb X}
    \widetilde L^{(-j)}([\hat{\bo \zeta}^{0,(-j)}, \hat{\mb X}^{\ur}, \hat{\mb Y}^{0,(-j)}\mb O^{0,(-j)}])}\fb \lesssim \sigma^2 \pi \sigma_1^* \frac{\sqrt{\kappa \mu r \rho }\xi^{1}}{\underline c_\psi \sqrt{\pi \sigma_r^*}}.
    \label{eq: leave-one-out gradient norm upper bound (column)}
    \end{align}

    Regarding the Hessian $\nabla_{\mb X_{i'}}^2
    \widetilde L^{(i)}([\hat{\bo \zeta}^{0,(i)}, \hat{\mb X}^{\ur}, \hat{\mb Y}^{0,(i)}\mb O^{0,(i)}])$, in the case $i' = i$, it is immediate from the bound for $\hat{\mb X}^{\ur}_{i}$ that 
    \begin{align}
        &c_1 \underline c_\psi \pi \sigma_r^*\mb I_r\preceq  \nabla_{\mb X_{i'}}^2
    \widetilde L^{(i)}([\hat{\bo \zeta}^{0,(i)}, \hat{\mb X}^{\ur}, \hat{\mb Y}^{0,(i)}\mb O^{0,(i)}]) \preceq c_2 \sigma^2 \pi \sigma_1^* \mb  I_r \label{eq: error control on leave-one-out Hessian}
    \end{align}
    for some constants $c_1$ and $c_2$ with probability at least $1-O(d^{-c-1})$ for each fixed $i$. If $i' \neq i$, we proceed by substituting $(\hat{\bo \zeta}^{0,(i)}, \hat{\mb Y}^{0,(i)})$ with $(\hat{\bo \zeta}^{0,(i')}, \hat{\mb Y}^{0,(i')})$ in the Hessian. In particular, one has 
    \begin{align}
        & \Bignorm{\nabla_{\mb X_{i'}}^2
    \widetilde L^{(i)}([\hat{\bo \zeta}^{0,(i)}, \hat{\mb X}^{\ur}, \hat{\mb Y}^{0,(i)}\mb O^{0,(i)}])\\ 
    & \qquad - 
   \mb O^{0,(i')\top} \hat{\mb Y}^{0,(i')\top} \diag(\psi'(\hat\zeta_j^{0,(i')} + \mb X^*_{i',\cdot} \tilde{\mb O}^{0,(i')\top} \hat{\mb Y}_j^{0,(i')} )\ind\{\Omega_{i',j} = 1\})_{j\in[p]}\\
    &\qquad \hat{\mb Y}^{0,(i')} \mb O^{0,(i')}} \\ 
    \lesssim & \Big\|\nabla_{\mb X_{i'}}^2
    \widetilde L^{(i)}([\hat{\bo \zeta}^{0,(i)}, \hat{\mb X}^{\ur}, \hat{\mb Y}^{0,(i)}\mb O^{0,(i)}])\\ 
    & - 
    \mb O^{0,(i')\top} \hat{\mb Y}^{0,(i')\top} \diag(\psi'(\hat\zeta_j^{0,(i)} + \hat{\mb X}^{\ur}_{i',\cdot} \mb O^{0,(i)\top} \hat{\mb Y}_j^{0,(i)}) \ind\{\Omega_{i',j} = 1\})_{j\in[p]}\\
    &\qquad \hat{\mb Y}^{0,(i')} \mb O^{0,(i')}\Big\| \\ 
    & +  \Big\|
    \hat{\mb Y}^{0,(i')\top} \diag\big([\psi'(\hat\zeta_j^{0,(i)} + \hat{\mb X}^{\ur}_{i',\cdot} \mb O^{0,(i)\top} \hat{\mb Y}_j^{0,(i)} ) \\ 
    & -\psi'( \hat\zeta_j^{0,(i')} + \mb X^*_{i',\cdot} \tilde{\mb O}^{0,(i')\top} \hat{\mb Y}_j^{0,(i')} ) ] \ind\{\Omega_{i',j} = 1\}  \big)_{j\in[p]} \hat{\mb Y}^{0,(i')}\| \\ 
    \lesssim & \sigma^2 \cdot \Big[\big( \bignorm{\hat{\mb Y}^{0,(i')} \mb O^{0,(i')} - \hat{\mb Y}^0} + \bignorm{\hat{\mb Y}^{0,(i)} \mb O^{0,(i)} - \hat{\mb Y}^0}\big)\\
    &\qquad \bignorm{\diag(\bo \Omega_{i',\cdot}) \hat{\mb Y}^{0,(i')}}\\
    &  + \big( \bignorm{\hat{\mb Y}^{0,(i')} \mb O^{0,(i')} - \hat{\mb Y}^0} + \bignorm{\hat{\mb Y}^{0,(i)} \mb O^{0,(i)} - \hat{\mb Y}^0}\big)^2\Big] \\ 
    &  + \max_{j\in[p]} \big|\psi'(\hat\zeta_j^{0,(i')} + \mb X^*_{i',\cdot} \tilde{\mb O}^{0,(i')\top} \hat{\mb Y}_j^{0,(i')}) - \psi'(\hat\zeta_j^{0,(i)} + \hat{\mb X}^{\ur}_{i',\cdot} \mb O^{0,(i)\top} \hat{\mb Y}_j^{0,(i)}) \big|\\
    &\qquad \bignorm{\diag(\bo \Omega_{i',\cdot}) \hat{\mb Y}^{0,(i')}}^2 . 
    \label{eq: leave-one-out hessian control}
    \end{align}
    The last term in this display is controlled using
    $|\psi'(a)-\psi'(b)|\lesssim\sigma^2|a-b|$,
    uniformly for $a,b\in\mc D_\psi$.

    For each fixed pair $(i,i')$, all concentration inputs to \eqref{eq: leave-one-out hessian control} are invoked with failure probability $O(d^{-c-2})$. A union bound over $(i,i')\in[n]^2$ therefore makes \eqref{eq: error control on leave-one-out Hessian} valid uniformly with failure probability $O(d^{-c})$.

    The leave-one-column-out Hessian is handled similarly. Fix any $i'\in[n]$ and $j\in[p]$. Its $k$th summand is weighted by $\ind\{\Omega_{i',k}=1\}$ when $k\neq j$ and by $\pi$ when $k=j$. To recover the independence needed for concentration, we replace $(\hat{\bo\zeta}^{0,(-j)},\hat{\mb Y}^{0,(-j)})$ by $(\hat{\bo\zeta}^{0,(i')},\hat{\mb Y}^{0,(i')})$ and replace $\hat{\mb X}_{i',\cdot}^{\ur}$ by $\mb X^*_{i',\cdot}$, with the corresponding rotations applied as above. The resulting comparison Hessian is independent of $\bo\Omega_{i',\cdot}$, and the expected weight of each summand remains $\pi$, including that corresponding to $k=j$. Hence, the same matrix Bernstein argument and population curvature bounds used above imply that its eigenvalues lie between constant multiples of $\underline c_\psi\pi\sigma_r^*$ and $\sigma^2\pi\sigma_1^*$. Moreover, the leave-one-out error bounds and the same perturbation decomposition as in \eqref{eq: leave-one-out hessian control} show that the operator-norm difference between the original and comparison Hessians is $o(\underline c_\psi\pi\sigma_r^*)$. Invoking the fixed-pair concentration bound at level $1-O(d^{-c-2})$ and taking a union bound over $(i',j)\in[n]\times[p]$, Weyl's inequality therefore yields, uniformly over $i'\in[n]$ and $j\in[p]$,
\begin{equation}
    c_1\underline c_\psi\pi\sigma_r^*\mb I_r
    \preceq
    \nabla_{\mb X_{i'}}^2
    \widetilde L^{(-j)}
    ([\hat{\bo\zeta}^{0,(-j)},\hat{\mb X}^{\ur},
      \hat{\mb Y}^{0,(-j)}\mb O^{0,(-j)}])
    \preceq
    c_2\sigma^2\pi\sigma_1^*\mb I_r.
    \label{eq: error control on leave-one-column-out Hessian}
\end{equation}
    
    Combining the above decomposition with the upper bounds on the leave-one-out errors and \eqref{eq: upper bound on Omega_i Y^{0,(i)}} implies that 
    \begin{align*}
        &  \Bignorm{\nabla_{\mb X_{i'}}^2
    \widetilde L^{(i)}([\hat{\bo \zeta}^{0,(i)}, \hat{\mb X}^{\ur}, \hat{\mb Y}^{0,(i)}\mb O^{0,(i)}])  \\ 
    & - 
    \mb O^{0,(i')\top} \hat{\mb Y}^{0,(i')\top} \diag(\psi'(\hat\zeta_j^{0,(i')} + \mb X^*_{i',\cdot} \tilde{\mb O}^{0,(i')\top} \hat{\mb Y}_j^{0,(i')} )\ind\{\Omega_{i',j} = 1\})_{j\in[p]}\\
    &\qquad \hat{\mb Y}^{0,(i')} \mb O^{0,(i')} } \ll \underline c_\psi \pi \sigma_r^*
    \end{align*}
    holds with probability at least $1 - O(d^{-c})$, which also finally leads to \eqref{eq: error control on leave-one-out Hessian}.

    Combining \eqref{eq: leave-one-out gradient norm upper bound} with \eqref{eq: error control on leave-one-out Hessian}, we arrive at the conclusion that with probability at least $1- O(d^{-c})$ 
    \begin{align} 
    & \bignorm{-\big[\nabla_{\mb X}^2
    \widetilde L^{(i)}([\hat{\bo \zeta}^{0,(i)}, \hat{\mb X}^{\ur}, \hat{\mb Y}^{0,(i)}\mb O^{0,(i)}])\big]^{-1}
    \sf{vec}\big(\nabla_{\mb X}
    \widetilde L^{(i)}([\hat{\bo \zeta}^{0,(i)}, \hat{\mb X}^{\ur}, \hat{\mb Y}^{0,(i)}\mb O^{0,(i)}])\big) }_2 \\
    \leq & C_1^{\sf{loo}} \frac{\kappa_\psi\kappa^{\frac32} \sqrt{\mu r \rho }\xi^{1}}{\underline c_\psi \sqrt{\pi \sigma_r^*}} \coloneqq L_1^{\sf{loo}}. 
    \label{eq: leave-one-out analysis for Hessian inverse times gradient} 
    \end{align}
    for some constant $C_1^{\sf{loo}}>0$.

    For the leave-one-column-out counterpart, putting \eqref{eq: leave-one-out gradient norm upper bound (column)} together with \eqref{eq: error control on leave-one-column-out Hessian} also yields that with probability at least $1- O(d^{-c})$  
    \begin{align}
        &  \bignorm{-\big[\nabla_{\mb X}^2
    \widetilde L^{(-j)}([\hat{\bo \zeta}^{0,(-j)}, \hat{\mb X}^{\ur}, \hat{\mb Y}^{0,(-j)}\mb O^{0,(-j)}])\big]^{-1}
    \\ &\qquad \sf{vec}\big(\nabla_{\mb X}
    \widetilde L^{(-j)}([\hat{\bo \zeta}^{0,(-j)}, \hat{\mb X}^{\ur}, \hat{\mb Y}^{0,(-j)} \mb O^{0,(-j)}])\big) }_2  \\
    & \leq   C_1^{\sf{loo}}\frac{\kappa_\psi\kappa^{\frac32} \sqrt{\mu r \rho}\xi^{1}}{\underline c_\psi \sqrt{\pi \sigma_r^*}}
    =L_1^{\sf{loo}}
    \label{eq: leave-one-out analysis for Hessian inverse times gradient (column)} 
    \end{align}
    Since $\bignorm{\sf{vec}(\mb A)}_2=\bignorm{\mb A}\fb$, these two displays give the Frobenius bounds for the full Newton steps used below.

    \item As for the Hessian-related integral part of \eqref{eq: definition of G}, define the domain $\big\{\mb X_{i'}: \bignorm{\mb X_{i'} - \hat{\mb X}^{\ur}_{i'}}_2 \leq C_2^{\sf{loo}} L_1^{\sf{loo}} \big\}$. The needed arguments highly resemble the proof for \eqref{eq: upper bound on the remainder term wrt X}, whose details are omitted for brevity: 
    \begin{align}
        & \Bignorm{
    \big[\nabla_{\mb X_{i'}}^2
    \widetilde L^{(i)}([\hat{\bo \zeta}^{0,(i)}, \hat{\mb X}^{\ur}, \hat{\mb Y}^{0,(i)}\mb O^{0,(i)}])\big]^{-1}
    \\ 
    &  \cdot 
    \int_0^1
    \Big[
    \nabla_{\mb X_{i'}}^2
    \widetilde L^{(i)}([\hat{\bo \zeta}^{0,(i)}, (1 - t)\hat{\mb X}^{\ur} + t \mb X, \hat{\mb Y}^{0,(i)} \mb O^{0,(i)}])\\
    &\qquad -
   \nabla_{\mb X_{i'}}^2
    \widetilde L^{(i)}([\hat{\bo \zeta}^{0,(i)}, \hat{\mb X}^{\ur}, \hat{\mb Y}^{0,(i)}\mb O^{0,(i)}])
    \Big]
    \mathrm d t \\ 
    & 
    \big(
    \mb X_{i'} - \hat{\mb X}_{i'}^{\ur}
    \big)} \leq  \frac{C_3^{\sf{loo}}}{\underline c_\psi  \pi\sigma_r^*}  \Big[\sigma^2 \sqrt{\frac{\mu r}{p}} \pi (\sigma_1^*)^{3/2} \bignorm{\mb X_{i'}  - \hat{\mb X}^{\ur}_{i'} }_2^2
        \Big] 
        \label{eq: leave-one-out analysis for Hessian residual}
    \end{align}
    uniformly over the leave-one-out index and $i'\in[n]$, with probability at least $1-O(d^{-c})$, for some constant $C_3^{\sf{loo}}>0$. The leave-one-column-out bound follows from the same fixed-pair argument and union bound. 
    \item Now it boils down to applying Brouwer's fixed point theorem. Concretely, we deduce from \eqref{eq: leave-one-out analysis for Hessian inverse times gradient} as well as \eqref{eq: leave-one-out analysis for Hessian residual} that, under the domain $\big\{\mb X: \bignorm{\mb X - \hat{\mb X}^{\ur}}\fb \leq C_2^{\sf{loo}} L_1^{\sf{loo}} \big\}$,
    \begin{align}
         \norm{\mb G_i(\mb X) - \hat{\mb X}^{\ur}}\fb
         &\leq L_1^{\sf{loo}}
         +C_3^{\sf{loo}}\kappa_\psi\kappa
         \sqrt{\frac{\mu r\sigma_1^*}{p}}
         \big(C_2^{\sf{loo}}L_1^{\sf{loo}}\big)^2 \\
         &\leq 2L_1^{\sf{loo}},
        \label{eq: G(X) contraction}
    \end{align}
    Here, the last condition in \eqref{eq: general unilateral refinement smallness conditions} gives
    $\kappa_\psi\kappa\sqrt{\mu r\sigma_1^*/p}\,L_1^{\sf{loo}}\ll1$.

    The same argument for the leave-one-column-out version gives
        \begin{align}
         \norm{\mb G_{-j}(\mb X) - \hat{\mb X}^{\ur}}\fb
         &\leq L_1^{\sf{loo}}
         +C_3^{\sf{loo}}\kappa_\psi\kappa
         \sqrt{\frac{\mu r\sigma_1^*}{p}}
         \big(C_2^{\sf{loo}}L_1^{\sf{loo}}\big)^2
         \leq 2L_1^{\sf{loo}}.
        \label{eq: G(X) contraction (column)}
    \end{align}
    Since $C_2^{\sf{loo}}>2$ can be chosen as a fixed constant, Eqs.~\eqref{eq: G(X) contraction} and \eqref{eq: G(X) contraction (column)} show that the corresponding maps send their domains into themselves. Brouwer's fixed-point theorem and strong convexity therefore give unique solutions within these domains with probability at least $1-O(d^{-c})$. In particular,
    \begin{align}
        &  \max_{i\in[n]}\bignorm{\hat{\mb X}^{\ur,(i)} \mb O^{0,(i)} - \hat{\mb X}^{\ur}}\fb\lesssim \frac{\kappa_\psi\kappa^{\frac32} \sqrt{\mu r\rho}\,\xi^{1}}{\underline c_\psi \sqrt{\pi \sigma_r^*}},
    \end{align}
    and the same bound holds for the leave-one-column-out solution. The nonlinear remainder is absorbed by the last condition in \eqref{eq: general unilateral refinement smallness conditions}.

\end{itemize}

\section{Gradient Descent from an \texorpdfstring{$\ell_{2,\infty}$}{l2infty}-Consistent Initialization}
\label{sec: supplementary gradient descent}

The final stage of the proposed pipeline applies Algorithm~\ref{alg: GD} to the one-step estimator. The preceding refinement provides a jointly aligned initialization that is accurate in both Frobenius and rowwise norms. Our goal is to show that these guarantees persist along the gradient-descent trajectory and, ultimately, yield the linear representations used for inference. The argument follows an implicit-regularization strategy: rather than imposing explicit incoherence or balancing constraints, we prove that a suitably initialized trajectory remains in a region where the empirical loss has stable local geometry.

\subsection{Parameter-Explicit Pipeline Guarantees}
\label{sec: parameter-explicit gradient results}
We first state the parameter-explicit counterparts of
Theorems~\ref{thm: GD} and~\ref{thm: linear approximations}. Unlike their
main-text specializations, the following results allow $r$, $\kappa$, and
$\kappa_\psi$ to vary.

\begin{theorem}[Parameter-Explicit Gradient-Descent Guarantee]
\label{thm: parameter-explicit GD}
Suppose that Assumptions~\ref{assumption: noise and sampling},
\ref{assumption: link function},
\ref{assumption: signal strength and incoherence degree}, and
\ref{assumption: identifiability} hold. Initialize
Algorithm~\ref{alg: GD} with the output of Algorithm~\ref{alg: OS}, itself
initialized by Algorithm~\ref{alg: USVT}, using the parameter-explicit tuning
choices in Theorem~\ref{thm: detailed USVT} and
Lemma~\ref{lemma: pipeline satisfies GD initialization}. For the theoretical calibration,
choose
\begin{align*}
    &\kappa r\sigma^2
    \vee\kappa^2\kappa_\psi^2(\log d)^4
    \ll c_\perp
    \ll\frac{\underline c_\psi d^{10}}{\kappa\log d},
    \qquad
    \eta\coloneqq\frac{1}{8c_\perp\pi\sigma_1^*},\\
    &\eta_\zeta\coloneqq\eta\frac{\sigma_r^*}{n},
    \qquad
    \eta_{\mb X}\coloneqq\eta\sqrt{\omega},
    \qquad
    \eta_{\mb Y}\coloneqq\frac{\eta}{\sqrt{\omega}},
    \qquad
    t_0\coloneqq
    \Big\lfloor\frac{c_t\log d}
    {\underline c_\psi\pi\sigma_r^*\eta}\Big\rfloor,
\end{align*}
where $c_t>0$ is a sufficiently large constant. The sampling condition in
Assumption~\ref{assumption: signal strength and incoherence degree} ensures
that this interval is nonempty and that $t_0\leq d^{10}$ for all sufficiently
large $d$. For each $0\leq t\leq t_0$, let $\mb O^t$ be a measurable choice
satisfying
\begin{equation*}
    \mb O^t\in\argmin_{\mb O\in\mc O(r)}
    \Big\{\omega^{-1/2}\bignorm{\mb X^t\mb O-\mb X^*}\fb^2
    +\omega^{1/2}\bignorm{\mb Y^t\mb O-\mb Y^*}\fb^2\Big\}.
\end{equation*}
Then, with probability at least $1-O(d^{-c})$, the final iterate satisfies
\begin{align*}
    &\bignorm{\bo\zeta^{t_0}-\bo\zeta^*}_2
    \lesssim
    \frac{\sigma\sqrt{\kappa d r\log d}}
    {\underline c_\psi\sqrt{\pi n}},\quad
    \omega^{-\frac14}
    \bignorm{\mb X^{t_0}\mb O^{t_0}-\mb X^*}\fb
    \vee
    \omega^{\frac14}
    \bignorm{\mb Y^{t_0}\mb O^{t_0}-\mb Y^*}\fb\\
    &\quad \lesssim
    \frac{\sigma\sqrt{\kappa d r\log d}}
    {\underline c_\psi\sqrt{\pi\sigma_r^*}}.
\end{align*}
\end{theorem}

\begin{theorem}[Parameter-Explicit Linear Approximations]
\label{thm: parameter-explicit linear approximations}
Suppose that Assumptions~\ref{assumption: noise and sampling},
\ref{assumption: link function},
\ref{assumption: signal strength and incoherence degree}, and
\ref{assumption: identifiability} hold. Run the full pipeline with the tuning
choices in Theorem~\ref{thm: detailed USVT} and
Lemma~\ref{lemma: pipeline satisfies GD initialization}, and choose $c_\perp$
in the range specified in Theorem~\ref{thm: parameter-explicit GD}. Suppose
further that
\begin{align}
    \sigma^2\pi(n\wedge p)
    &\gg\kappa^{24}\kappa_\psi^{10}\mu^7r^{11}\rho^3
    \xi^{10}(\log d)^6,
    \label{eq: effective information condition for linear approximation}\\
    \sigma^2\kappa^{\frac{23}{2}}\kappa_\psi^7
    \mu^{\frac94}r^{\frac{13}{4}}\rho^2(\log d)^3
    \sqrt{\bar\sigma\sigma_1^*}
    \Big(\frac{\pi}{d}\Big)^{\frac14}\xi^{\frac52}
    &\lesssim c_\perp.
    \label{eq: cperp condition for linear approximation}
\end{align}
Let $\mb E\coloneqq\mb R-\bb E[\mb R]=\mb R-\psi(\mb M^*)$ denote the
response-noise matrix. Set $(\eta_\zeta,\eta_{\mb X},\eta_{\mb Y})$ and
$t_0$ according to
Theorem~\ref{thm: parameter-explicit GD}. Then, with probability at least
$1-O(d^{-c})$, the final output admits the following uniform rowwise linear
approximations in terms of the score noise:
\begin{align}
    &\sqrt{\frac n{\sigma_r^*}}\max_{j\in[p]}
    \Big|(\zeta_j^{t_0}-\zeta_j^*)
    -\sqrt{\frac{\sigma_r^*}{n}}\,
    \mc P_\Omega(\mb E)_{\cdot,j}^\top\mb X^{\sf{app},*}
    \big(\mb H_{\mb Y^{\sf{app}}_{j}}^{*-1}\big)_{\cdot,1}\Big|
    \nonumber\\
    &\vee\ \omega^{-\frac14}\max_{i\in[n]}
    \Bignorm{(\mb X_{i,\cdot}^{t_0}\mb O^{t_0}-\mb X_{i,\cdot}^*)
    -\mc P_\Omega(\mb E)_{i,\cdot}\mb Y^*
    \mb H_{\mb X_i}^{*-1}}_2
    \nonumber\\
    &\vee\ \omega^{\frac14}\max_{j\in[p]}
    \Bignorm{(\mb Y_{j,\cdot}^{t_0}\mb O^{t_0}-\mb Y_{j,\cdot}^*)
    -\mc P_\Omega(\mb E)_{\cdot,j}^\top\mb X^{\sf{app},*}
    \big(\mb H_{\mb Y^{\sf{app}}_{j}}^{*-1}\big)_{\cdot,2:r+1}}_2
    \lesssim
    \frac{1}{\sigma\sqrt{\pi\sigma_1^*}(\log d)^2},
    \label{eq: parameter-explicit linear approximation error}
\end{align}
where
$\mb X^{\sf{app},*}\coloneqq
(\sqrt{\sigma_r^*/n}\mb 1_n,\mb X^*)\in\bb R^{n\times(r+1)}$ and
\begin{align*}
    \mb H_{\mb X_i}^*
    &\coloneqq
    \sum_{\substack{j\in[p]\\(i,j)\in\Omega}}
    \psi'(M_{i,j}^*)\mb Y_j^*\mb Y_j^{*\top},
    \qquad
    \mb H_{\mb Y^{\sf{app}}_{j}}^*
    \coloneqq
    \sum_{\substack{i\in[n]\\(i,j)\in\Omega}}
    \psi'(M_{i,j}^*)\mb X^{\sf{app},*}_{i}\mb X^{\sf{app},*\top}_{i}.
\end{align*}
Consequently, the final iterate obeys
\begin{align}
    &\sqrt{\frac n{\sigma_r^*}}
    \bignorm{\bo\zeta^{t_0}-\bo\zeta^*}_\infty
    \vee\omega^{-\frac14}
    \bignorm{\mb X^{t_0}\mb O^{t_0}-\mb X^*}\ti
    \vee\omega^{\frac14}
    \bignorm{\mb Y^{t_0}\mb O^{t_0}-\mb Y^*}\ti
    \lesssim
    \frac{\sigma\sqrt{\kappa\mu r\log d}}
    {\underline c_\psi\sqrt{\pi\sigma_r^*}}.
    \label{eq: parameter-explicit rowwise error}
\end{align}
\end{theorem}

\paragraph{Supplementary gradient-descent setup.}
We adopt the notation and balanced parameterization from the preliminary setups; in particular, Lemma~\ref{lemma: scaling equivalence} allows us to set $\omega=1$ throughout this section. Recall the shorthand introduced in the preliminary setup: for any conformable triple $(\bo a,\mb A,\mb B)$,
$
    [\bo a,\mb A,\mb B]
    \coloneqq
    \big(
        \sqrt{\frac{n}{\sigma_r^*}}\bo a^\top,\,
        \sf{vec}(\mb A)^\top,\,
        \sf{vec}(\mb B)^\top
    \big)^\top
$, 
so that the intercept block is automatically rescaled whenever the bracket notation is used. 
Accordingly, we write the full gradient with respect to this rescaled parameter vector as
\begin{align*}
    \nabla L(\bo\theta)
    \coloneqq
    \Big(
    \sqrt{\frac{\sigma_r^*}{n}}\nabla_{\bo\zeta}L(\bo\theta)\t,
    \sf{vec}(\nabla_{\mb X}L(\bo\theta))\t,
    \sf{vec}(\nabla_{\mb Y}L(\bo\theta))\t
    \Big)\t.
\end{align*}
The same convention applies to the full gradients of $L^{\sf{aug}}$, $L^{\sf{diff}}$, and $L^{(l)}$, while block gradients denote derivatives with respect to the corresponding unscaled blocks.
Recall that $\mc M(\bo\theta)=\mb 1_n\bo\zeta\t+\mb X\mb Y\t$ for $\bo\theta=[\bo\zeta,\mb X,\mb Y]$.

Write the augmented loss as
\begin{align}
    L^{\sf{aug}}([\bo\zeta,\mb X,\mb Y])
    &\coloneqq L([\bo\zeta,\mb X,\mb Y])
      +L^{\sf{diff}}(\mb X,\mb Y),
      \label{eq: definition of Laug}\\
    L^{\sf{diff}}(\mb X,\mb Y)
    &\coloneqq c_{\sf{aug}}\pi
      \bignorm{\mb X\t\mb X-\mb Y\t\mb Y}\fb^2.
      \label{eq: definition of Ldiff}
\end{align}
Writing $\mb B(\mb X,\mb Y)\coloneqq\mb X\t\mb X-\mb Y\t\mb Y$, the derivatives of the added term are
\begin{equation}
    \nabla_{\mb X}L^{\sf{diff}}
    =4c_{\sf{aug}}\pi\mb X\mb B(\mb X,\mb Y),
    \qquad
    \nabla_{\mb Y}L^{\sf{diff}}
    =-4c_{\sf{aug}}\pi\mb Y\mb B(\mb X,\mb Y).
    \label{eq: derivatives of Ldiff}
\end{equation}
The balancing term strengthens the curvature along the otherwise problematic scaling directions, but it does not remove the rotational invariance of the factorization. Accordingly, the local curvature result below is restricted to suitably aligned directions. We will show separately that the original iterates remain approximately centered, incoherent, and balanced, so that the additional gradient in \eqref{eq: derivatives of Ldiff} is negligible along the trajectory. This comparison transfers the contraction properties of the augmented dynamics to the actual iterates of Algorithm~\ref{alg: GD}.

\subsection{Auxiliary Results and Proofs of Theorems~\ref{thm: parameter-explicit GD} and~\ref{thm: parameter-explicit linear approximations}}
\label{subsec: auxiliary results for GD}

We first collect the auxiliary results needed to prove
Theorem~\ref{thm: parameter-explicit GD}.
The first result controls the local geometry of the augmented loss. The second
shows that the first two stages of the proposed pipeline deliver the required
initialization, and the third establishes the gradient-descent guarantees from
an arbitrary initialization satisfying these requirements. Their proofs are
deferred to the subsequent subsections. The general rowwise score-replacement
lemmas used to deduce
Theorem~\ref{thm: parameter-explicit linear approximations} are developed in
Section~\ref{subsec: proof of inference results}.

We begin by stating precisely the initialization assumption and the inductive
condition used below. For the actual and leave-one-out sequences, write
\begin{equation*}
    \mb F^t \coloneqq
    \begin{pmatrix}\mb X^t\\ \mb Y^t\end{pmatrix},
    \qquad
    \mb F^* \coloneqq
    \begin{pmatrix}\mb X^*\\ \mb Y^*\end{pmatrix},
    \qquad
    \mb F^{t,(l)} \coloneqq
    \begin{pmatrix}\mb X^{t,(l)}\\ \mb Y^{t,(l)}\end{pmatrix},
    \quad t\geq 0,\quad l\in[n]\cup[-p].
\end{equation*}
For $\bo\theta_k=[\bo\zeta_k,\mb X_k,\mb Y_k]$ and
$\mb F_k=(\mb X_k\t,\mb Y_k\t)\t$, $k\in\{1,2\}$, fix a measurable common
Procrustes rotation $\mb O(\bo\theta_1,\bo\theta_2)
\in\argmin_{\mb O\in\mc O(r)}
\bignorm{\mb F_1\mb O-\mb F_2}\fb$.
For the actual and leave-one-out sequences, abbreviate
\begin{align*}
    \mb O^t
    &\coloneqq\mb O(\bo\theta^t,\bo\theta^*),
    &
    \mb O^{t,(l)}
    &\coloneqq\mb O(\bo\theta^{t,(l)},\bo\theta^t),
    &
    \widetilde{\mb O}^{t,(l)}
    &\coloneqq\mb O(\bo\theta^{t,(l)},\bo\theta^*),
    \quad l\in[n]\cup[-p].
\end{align*}
Using this same rotation, define
\begin{align*}
    \distf(\bo\theta_1,\bo\theta_2)^2
    &\coloneqq
    \frac{n}{\sigma_r^*}\norm{\bo\zeta_1-\bo\zeta_2}_2^2
    +\bignorm{\mb F_1\mb O(\bo\theta_1,\bo\theta_2)-\mb F_2}\fb^2,\\
    \dist\ti(\bo\theta_1,\bo\theta_2)
    &\coloneqq
    \sqrt{\frac{n}{\sigma_r^*}}\norm{\bo\zeta_1-\bo\zeta_2}_\infty
    +\bignorm{\mb F_1\mb O(\bo\theta_1,\bo\theta_2)-\mb F_2}\ti.
\end{align*}

\subsubsection{Local Geometry}
We establish the following lemma to control the curvature of the augmented loss function $L^{\sf{aug}}$ in \eqref{eq: definition of Laug}. The proof is deferred to Section~\ref{subsec: local geometry proof}.
\begin{lemma}\label{lemma: bounds for Hessian}
       Suppose that Assumptions~\ref{assumption: noise and sampling},~\ref{assumption: link function},~\ref{assumption: signal strength and incoherence degree},~and~\ref{assumption: identifiability} hold and $c_\perp \gg \kappa r\sigma^2$.
       Set $\epsilon\coloneqq\epsilon_0/\kappa_\psi^2$ for a sufficiently small absolute constant $\epsilon_0>0$. With the prescribed choice $c_{\sf{aug}}\asymp\underline c_\psi$, take $c_0\coloneqq1\wedge(4c_{\sf{aug}}/\underline c_\psi)$, $\underline c_{\sf{Hess}}\coloneqq c_0\underline c_\psi/64$, and $\bar c_{\sf{Hess}}\coloneqq4$. The following inequalities hold uniformly with probability at least $1- O(d^{-c-1})$:
       \begin{align}
        \bo \Delta\t \nabla_{\bo \theta}^2 L^{\sf{aug}}(\bo \theta)\bo \Delta
        &\geq \underline c_{\sf{Hess}}\pi\sigma_r^*\norm{\bo \Delta}_2^2,
        \label{eq: local Hessian lower bound}\\
        \bignorm{\nabla_{\bo\theta}^2L^{\sf{aug}}(\bo\theta)}
        &\leq \bar c_{\sf{Hess}}c_\perp\pi\sigma_1^*.
        \label{eq: local Hessian upper bound}
       \end{align}
       These inequalities hold for $\bo \Delta = [\bo \Delta_{\bo \zeta}, \bo \Delta_{\mb X}, \bo \Delta_{\mb Y}]$ satisfying
       \begin{equation}
\begin{aligned}
    & \bo \Delta_{\mb F} \coloneqq \begin{bmatrix}
            \bo \Delta_{\mb X}\\
            \bo \Delta_{\mb Y}
        \end{bmatrix}
        \in \Big\{\mb F_2 \hat{\mb O} - \mb F_1:
        \hat{\mb O} \in \arg\min_{\mb O \in \mc O(r)}\norm{\mb F_2 \mb O - \mb F_1 }\fb,~
        \norm{\mb F_1 - \mb F^*}\\
    &\quad \leq \sqrt{\epsilon \sigma_r^*} / \kappa \Big\},
        \qquad
        \mb F \coloneqq \begin{pmatrix}\mb X \\ \mb Y \end{pmatrix},
        \label{eq: direction set}
\end{aligned}
       \end{equation}
       and for $\bo \theta = [\bo \zeta, \mb X, \mb Y]$ staying in the \emph{region of incoherence and contraction} (\emph{RIC}), namely, the set of parameters $[\bo \zeta, \mb X, \mb Y]$ satisfying
       \begin{align}
        &  \norm{\mb X - \mb X^*}\ti \leq \frac{\sqrt{\epsilon}\norm{\mb X^*}}{\kappa \sqrt{n}},~\norm{\mb Y - \mb Y^*}\ti \leq \frac{ \sqrt{\epsilon}\norm{\mb Y^*}}{\kappa \sqrt{p}}  ,~ \norm{\bo \zeta - \bo \zeta^*}_\infty \leq \frac{\sqrt{\epsilon}\sigma_1^* }{\kappa \sqrt{np}}. 
        \label{eq: rowwise consistency of parameters in geometry lemma}
       \end{align}
       The same conclusions hold for tangent directions satisfying
       $\mb F_1\t\bo\Delta_{\mb F}=\bo\Delta_{\mb F}\t\mb F_1$.
    \end{lemma}
    \begin{remark}
{\emergencystretch=1.5em
    The Hessian is necessarily degenerate along rotational directions. Lemma~\ref{lemma: bounds for Hessian} therefore controls aligned differences and their horizontal tangent directions, which are precisely the directions used in the contraction arguments below.\par}

\end{remark} 

\subsubsection{Gradient Descent from a General Warm Start}

We first formalize the centering-and-balancing step in
Algorithm~\ref{alg: GD}. For any preliminary estimate
$[\widetilde{\bo\zeta},\widetilde{\mb X},\widetilde{\mb Y}]$, let
\begin{equation}
    \widetilde{\mb X}^{\sf c}
    \coloneqq\mb J_n\widetilde{\mb X},
    \label{eq: centering map for GD initialization}
\end{equation}
where $\mb J_n\coloneqq\mb I_n-n^{-1}\mb 1_n\mb 1_n\t$. Given a top-$r$
SVD $\widetilde{\mb X}^{\sf c}\widetilde{\mb Y}\t
    =\mb U\bo\Sigma\mb V\t$,
define
\begin{equation}
    \mathfrak C\big(
    [\widetilde{\bo\zeta},\widetilde{\mb X},\widetilde{\mb Y}]
    \big)
    \coloneqq
    [\widetilde{\bo\zeta},\mb U\bo\Sigma^{1/2},
    \mb V\bo\Sigma^{1/2}].
    \label{eq: centering and balancing map}
\end{equation}
We fix an arbitrary measurable convention for the SVD; its orthogonal
nonuniqueness is immaterial because all factor errors are defined after a
common Procrustes alignment.
This section works under $\omega=1$; the factors in
Algorithm~\ref{alg: GD} recover the prescribed general scaling by multiplying
the last two components in \eqref{eq: centering and balancing map} by
$\omega^{1/4}$ and $\omega^{-1/4}$, respectively.

We now formalize the gradient-descent initialization requirements. In
particular, we use the radius
$\epsilon=\epsilon_0/\kappa_\psi^2$ from
Lemma~\ref{lemma: bounds for Hessian} to make the required RIC neighborhood
explicit.
\begin{assumption}[Warm-Start Initialization]
\label{assumption: GD initialization}
There exist leave-one-out initializations
$[\bo\zeta^{0,(i)},\mb X^{0,(i)},\mb Y^{0,(i)}]$, $i\in[n]$, and
$[\bo\zeta^{0,(-j)},\mb X^{0,(-j)},\mb Y^{0,(-j)}]$, $j\in[p]$, such that the
former is independent of $(\mb E_{i,\cdot},\bo\Omega_{i,\cdot})$ and the latter
is independent of $(\mb E_{\cdot,j},\bo\Omega_{\cdot,j})$. Moreover, the
actual and leave-one-out initializations satisfy
\begin{equation}
\begin{aligned}
    & \mb 1_n\t\mb X^0=\mb 0,\qquad
    {\mb X^0}\t\mb X^0={\mb Y^0}\t\mb Y^0,
    \qquad
    \mb 1_n\t\mb X^{0,(l)}=\mb 0,\qquad
    {\mb X^{0,(l)}}\t\mb X^{0,(l)}\\
    &\quad ={\mb Y^{0,(l)}}\t\mb Y^{0,(l)}
    \quad\text{for all }l\in[n]\cup[-p],
    \label{eq: exact GD initialization constraints}
\end{aligned}
\end{equation}
and
\begin{subequations}
\label{eq: GD initialization requirement}
\begin{align}
   & \underline c_\psi \sqrt{\pi \sigma_r^*}\max\Big\{
    \sqrt{\frac{n}{\sigma_r^*}}\norm{\bo\zeta^0-\bo\zeta^*}_\infty,
    \bignorm{\mb X^0\mb O^0-\mb X^*}\ti,
    \bignorm{\mb Y^0\mb O^0-\mb Y^*}\ti
    \Big\}
    \lesssim \xi^{\sf{gd}},
    \label{eq: GD initialization global}\\
   & \underline c_\psi \sqrt{\pi \sigma_r^*}\max_{i\in[n]}\max\Big\{
    \sqrt{\frac{n}{\sigma_r^*}}\bignorm{\bo\zeta^{0,(i)}-\bo\zeta^0}_2,
    \bignorm{\mb X^{0,(i)}\mb O^{0,(i)}-\mb X^0}\fb,\\
    &\quad \bignorm{\mb Y^{0,(i)}\mb O^{0,(i)}-\mb Y^0}\fb
    \Big\}
    \lesssim \xi^{\sf{gd},\sf{loo}},
    \label{eq: GD initialization row LOO}\\
    &\underline c_\psi \sqrt{\pi \sigma_r^*}\max_{j\in[p]}\max\Big\{
    \sqrt{\frac{n}{\sigma_r^*}}\bignorm{\bo\zeta^{0,(-j)}-\bo\zeta^0}_2,
    \bignorm{\mb X^{0,(-j)}\mb O^{0,(-j)}-\mb X^0}\fb,\\
    &\quad \bignorm{\mb Y^{0,(-j)}\mb O^{0,(-j)}-\mb Y^0}\fb
    \Big\}
    \lesssim \xi^{\sf{gd},\sf{loo}}.
    \label{eq: GD initialization column LOO}
\end{align}
\end{subequations}
They satisfy the following relation:
\begin{equation}
\begin{aligned}
    &\kappa_\psi^3\kappa^3\log d
    \Big\{\kappa\kappa_\psi\sqrt{\mu\rho r}\,
    \xi^{\sf{gd}}+\xi^{\sf{gd},\sf{loo}}\Big\}\ll \underline c_\psi\sigma_r^*
    \sqrt{\frac{\pi}{\kappa d}}.
\end{aligned}
    \label{eq: epsilon GD initialization requirement}
\end{equation}

For notational convenience, define
\begin{align}
    \delta_{\sf{loo}}^0
    &\coloneqq
    \frac{\xi^{\sf{gd},\sf{loo}}}
    {\underline c_\psi\sqrt{\pi\sigma_r^*}},\qquad
    \delta_{\sf{row}}^0
    \coloneqq
    \frac{\xi^{\sf{gd}}}
    {\underline c_\psi\sqrt{\pi\sigma_r^*}}.
    \label{eq: GD initialization error envelopes}
\end{align}
The corresponding Frobenius initialization scale is
$\sqrt d\,\delta_{\sf{row}}^0$.
\end{assumption}

Note that the actual input to Algorithm~\ref{alg: GD} may not fulfill \eqref{eq: exact GD initialization constraints} exactly, but the centering-and-balancing map $\mathfrak C$ in \eqref{eq: centering and balancing map} ensures that the one after preprocessing in Algorithm~\ref{alg: GD} does. The following lemma establishes the proximity between the input and the output of the centering-and-balancing map, which is used to verify Assumption~\ref{assumption: GD initialization}. 
\begin{lemma}[Initialization delivered by the warm-start pipeline]
\label{lemma: pipeline satisfies GD initialization}
Suppose that Assumptions~\ref{assumption: noise and sampling},
\ref{assumption: link function},
\ref{assumption: signal strength and incoherence degree}, and~\ref{assumption: identifiability} hold.
Let Algorithm~\ref{alg: USVT} use the threshold in
Theorem~\ref{thm: detailed USVT}, and let Algorithm~\ref{alg: OS} use
$\iota_{\bo\zeta}=c_\iota\sqrt{\mu r\sigma_1^*\sigma_r^*/(np)}$,
$\iota_{\mb X}^{(\omega)}
=c_\iota\omega^{1/4}\sqrt{\mu r\sigma_1^*/n}$, and
$\iota_{\mb Y}^{(\omega)}
=c_\iota\omega^{-1/4}\sqrt{\mu r\sigma_1^*/p}$ for a sufficiently large
constant $c_\iota$. Write the output of
Algorithm~\ref{alg: OS} as
$[\widetilde{\bo\zeta}^{0},\widetilde{\mb X}^{0},
\widetilde{\mb Y}^{0}]
\coloneqq[\hat{\bo\zeta}^{\ur},\hat{\mb X}^{\ur},\hat{\mb Y}^{\ur}]$, and set
$[\bo\zeta^0,\mb X^0,\mb Y^0]
\coloneqq\mathfrak C(
[\widetilde{\bo\zeta}^{0},\widetilde{\mb X}^{0},
\widetilde{\mb Y}^{0}])$.
Write the corresponding unprocessed leave-one-out refinements as
$[\widetilde{\bo\zeta}^{0,(l)},\widetilde{\mb X}^{0,(l)},
\widetilde{\mb Y}^{0,(l)}]$, and apply the same map to define
$[\bo\zeta^{0,(l)},\mb X^{0,(l)},\mb Y^{0,(l)}]$ for
$l\in[n]\cup[-p]$. Then the independence properties, the identities in
\eqref{eq: exact GD initialization constraints}, and the bounds in
\eqref{eq: GD initialization requirement} and
\eqref{eq: epsilon GD initialization requirement}, as well as
Condition~\ref{condition: inductive contraction for GD} at $t=0$, hold with
probability at least $1-O(d^{-c})$, with
\begin{equation}
    \xi^{\sf{gd}}
    \coloneqq
    \kappa^{5/2}\sqrt{\mu r\rho}\,\xi^1,
    \qquad
    \xi^{\sf{gd},\sf{loo}}
    \coloneqq
    \kappa_\psi\kappa^3\sqrt{\mu r\rho}\,\xi^1,
    \label{eq: GD initialization rates from pipeline}
\end{equation}
where $\xi^1$ is defined in Theorem~\ref{thm: general OS} and obeys the
explicit concrete-pipeline bound established in the proof of
Theorem~\ref{thm: OS}.
\end{lemma}

The following theorem establishes the gradient-descent guarantees from a general warm start. The proof is deferred to Section~\ref{subsec: GD proof}. 
\begin{theorem}[Gradient Descent from a General Warm Start]
\label{thm: general GD}
Suppose that Assumptions~\ref{assumption: noise and sampling},
\ref{assumption: link function}, \ref{assumption: signal strength and incoherence degree}, and
\ref{assumption: identifiability} hold. Suppose further that
Assumption~\ref{assumption: GD initialization} holds for some
$\xi^{\sf{gd}}$ and $\xi^{\sf{gd},\sf{loo}}$, and that
Condition~\ref{condition: inductive contraction for GD} holds at $t=0$. Let
$c_\perp\gg\kappa r\sigma^2$ and
$c_{\perp}\gtrsim\kappa^2\kappa_\psi^2(\log d)^4$,
$\eta\leq 1 / (8 c_\perp \pi\sigma_1^*)$, and
$t_0=\lfloor c_t\log d/(\eta \underline c_\psi \pi\sigma_r^*) \rfloor\leq d^{10}$ for a sufficiently large constant $c_t$.
Then, with
probability at least $1-O(d^{-c})$, the actual and leave-one-out iterates remain
in the RIC for every $0\leq t\leq t_0$. Moreover, let
$\mb O^{t_0}\coloneqq\mb O(\bo\theta^{t_0},\bo\theta^*)$ be the common
Procrustes rotation used in the definition of $\dist\ti$. The final iterate
in the balanced parametrization satisfies
\begin{align}
    \bignorm{\bo\zeta^{t_0}-\bo\zeta^*}_2
    &\lesssim
    \frac{\sigma\sqrt{\kappa d r\log d}}
    {\underline c_\psi\sqrt{\pi n}},
    \label{eq: conditional GD zeta rate}\\
    \bignorm{\mb X^{t_0}\mb O^{t_0}-\mb X^*}\fb
    &\lesssim
    \frac{\sigma\sqrt{\kappa dr\log d}}
    {\underline c_\psi\sqrt{\pi\sigma_r^*}},
    \qquad
    \bignorm{\mb Y^{t_0}\mb O^{t_0}-\mb Y^*}\fb
    \lesssim
    \frac{\sigma\sqrt{\kappa dr\log d}}
    {\underline c_\psi\sqrt{\pi\sigma_r^*}}.
    \label{eq: conditional GD factor rates}
\end{align}
Furthermore, the initialization-dependent terms in
Condition~\ref{condition: inductive contraction for GD} are absorbed into the
stochastic-error terms at $t=t_0$, and hence
\begin{subequations}
\label{eq: final GD inductive bounds}
\begin{align}
&\max_{l\in[n]\cup[-p]}
\distf(\bo\theta^{t_0},\bo\theta^{t_0,(l)})
\lesssim
\frac{\sigma\kappa\kappa_\psi
\sqrt{\kappa\mu\rho r^2\log d}}
{\underline c_\psi\sqrt{\pi\sigma_r^*}},
\label{eq: final GD LOO coupling}\\
&\distf(\bo\theta^{t_0},\bo\theta^*)
\lesssim
\frac{\sigma\sqrt{d\kappa r\log d}}
{\underline c_\psi\sqrt{\pi\sigma_r^*}},
\label{eq: final GD global error}\\
    &\dist\ti(\bo\theta^{t_0},\bo\theta^*)\vee\max_{l\in[n]\cup[-p]}
\dist\ti(\bo\theta^{t_0,(l)},\bo\theta^*)
\lesssim
\frac{\sigma\kappa\kappa_\psi
\sqrt{\kappa\mu\rho r^2\log d}}
{\underline c_\psi\sqrt{\pi\sigma_r^*}},
\label{eq: final GD rowwise errors}\\
&\bignorm{\mb 1_n\t\mb X^{t_0}}_2
\vee\max_{l\in[n]\cup[-p]}
\bignorm{\mb 1_n\t\mb X^{t_0,(l)}}_2
\lesssim
\sqrt n\frac{\sigma\kappa^{\frac32}\kappa_\psi\sqrt{dr\log d}
+\kappa\kappa_\psi\sqrt d\big(\xi^{\sf{gd}}+\xi^{\sf{gd},\sf{loo}}\big)}
{c_\perp\sqrt{\pi\sigma_r^*}},
\label{eq: final GD centering}\\
&\bignorm{{\mb X^{t_0}}\t\mb X^{t_0}
-{\mb Y^{t_0}}\t\mb Y^{t_0}}\fb
\vee\max_{l\in[n]\cup[-p]}
\bignorm{{\mb X^{t_0,(l)}}\t\mb X^{t_0,(l)}
-{\mb Y^{t_0,(l)}}\t\mb Y^{t_0,(l)}}\fb\\
&\qquad\lesssim
\sqrt{\sigma_1^*}
\frac{\sigma\kappa^{\frac32}\kappa_\psi\sqrt{dr\log d}
+\kappa\kappa_\psi\sqrt d\big(\xi^{\sf{gd}}+\xi^{\sf{gd},\sf{loo}}\big)}
{c_\perp\sqrt{\pi\sigma_r^*}}.
\label{eq: final GD balancedness}
\end{align}
\end{subequations}
The same conclusions hold with exceptional probability $O(d^{-c-q})$ for
any fixed $q\in\{1,2,3\}$, after increasing the fixed constants in the
concentration bounds. Assumption~\ref{assumption: noise and sampling}, whose
truncation event has exceptional probability $O(d^{-c-10})$, provides the
required probability slack.
\end{theorem}

\begin{proof}[Proof of Theorem~\ref{thm: parameter-explicit GD}]
By Lemma~\ref{lemma: scaling equivalence}, it suffices first to consider the
balanced parametrization $\omega=1$. The detailed spectral guarantee in
Theorem~\ref{thm: detailed USVT} verifies the global and leave-one-out
conditions required by Theorem~\ref{thm: general OS}. Consequently,
Lemma~\ref{lemma: pipeline satisfies GD initialization} shows that the
centered-and-balanced refactorization of the first two stages, together with
its leave-one-out counterparts, satisfies
Assumption~\ref{assumption: GD initialization} and
Condition~\ref{condition: inductive contraction for GD} at $t=0$. In
particular, the initialization
invoked in the gradient-descent analysis is a consequence of the preceding
pipeline rather than an additional assumption.

We may therefore apply Theorem~\ref{thm: general GD}. Equations
\eqref{eq: conditional GD zeta rate},
\eqref{eq: conditional GD factor rates}, and
\eqref{eq: final GD rowwise errors} give the claimed bounds when $\omega=1$.
The weighted common-alignment criterion in
Theorem~\ref{thm: parameter-explicit GD} is invariant
under the scaling equivalence. Thus the same common rotation applies for
general $\omega$, while Lemma~\ref{lemma: scaling equivalence} multiplies the
$\mb X$- and $\mb Y$-errors by $\omega^{1/4}$ and $\omega^{-1/4}$,
respectively, and leaves the intercept unchanged. This proves the theorem.
\end{proof}

\subsection{Proof of Lemma~\ref{lemma: bounds for Hessian}}
\label{subsec: local geometry proof}
        Firstly, it follows from direct algebra that
        \begin{equation}
        \begin{split}
             &\bo \Delta\t \nabla^2 L^{\sf{aug}}(\bo \theta) \bo \Delta  
            =  \alpha_1 + \alpha_2,\label{eq: quadratic form of the augmented function} 
            \end{split}
        \end{equation}
        where $\alpha_1$ and $\alpha_2$ are defined as
        \begin{align}
            & \alpha_1 \coloneqq \norm{(\psi')^{\frac12}(\mc M(\bo \theta)) \circ \mc P_{\Omega}(\mb 1_n\bo \Delta_{\bo \zeta}\t + \bo \Delta_{\mb X}\mb Y\t + \mb X \bo \Delta_{\mb Y}\t)}\fb^2 \\ 
            & + 2 \pi \Big[c_{\perp} \frac{\sigma_r^*}{n}\bignorm{\mb 1_n\t \bo \Delta_{\mb X}}_2^2   + c_{\sf{aug}}\bignorm{\bo \Delta_{\mb X}\t \mb X + \mb X\t {\bo \Delta_{\mb X}} - \bo \Delta_{\mb Y}\t \mb Y - \mb Y\t {\bo \Delta_{\mb Y}}}\fb^2\Big], \\ 
            & \alpha_2 \coloneqq  2 \big \langle \mc P_{\Omega}(\psi(\mc M(\bo \theta)) - \mb R), \bo \Delta_{\mb X}{\bo \Delta_{\mb Y}}\t \big \rangle  + 4 c_{\sf{aug}} \pi \bigip{\mb X\t \mb X - \mb Y\t \mb Y, \bo \Delta_{\mb X}\t \bo \Delta_{\mb X} -\bo \Delta_{\mb Y}\t \bo \Delta_{\mb Y}}. 
        \end{align}
        Here and throughout, $\mc P_\Omega$ denotes the sampling operator introduced in Section~\ref{subsec: notation}.

        We start with parsing the first term $\alpha_1$ in \eqref{eq: quadratic form of the augmented function}. Fix $c_{\mc D_\psi}$ in Assumption~\ref{assumption: link function} sufficiently large relative to the fixed constants in the population centered intercept and factor row bounds and in \eqref{eq: rowwise consistency of parameters in geometry lemma}. Throughout the RIC, these bounds give $|\mc M(\bo\theta)_{i,j}-\zeta_0|\leq|\zeta_j-\zeta_0|+\norm{\mb X_i}_2\norm{\mb Y_j}_2\lesssim\mu r\sigma_1^*/\sqrt{np}$. Consequently, the predictor entries lie in $\mc D_\psi$, as do all natural-parameter interpolants between two such predictors. We can therefore lower bound every entry of $\psi'(\mc M(\bo \theta))$ and obtain the following decomposition for $c_0= 1\wedge \frac{4c_{\sf{aug}}}{\underline c_\psi}$, which separates out the randomness of $\Omega$ from the deterministic part:
        \begin{align*}
    & \alpha_1 \geq   \underline c_\psi \bignorm{\mc P_{\Omega}(\mb 1_n \bo \Delta_{\bo \zeta}\t + \bo \Delta_{\mb X} \mb Y\t + \mb X \bo \Delta_{\mb Y}\t)}\fb^2\\
    &\quad +   c_{\sf{aug}}\pi \bignorm{\bo \Delta_{\mb X}\t \mb X + \mb X\t {\bo \Delta_{\mb X}} - \bo \Delta_{\mb Y}\t \mb Y - \mb Y\t {\bo \Delta_{\mb Y}} }\fb^2\\
    & \qquad + c_{\perp} \pi \frac{\sigma_r^*}{n}\bignorm{\mb 1_n\t \bo \Delta_{\mb X}}_2^2\\
    & \geq   c_0  \pi  \underline c_\psi \Big[ \bignorm{\mb 1_n\bo \Delta_{\bo\zeta}\t  + \bo \Delta_{\mb X}\mb Y\t + \mb X \bo \Delta_{\mb Y}\t }\fb^2\\
    &\quad + \frac{1}{4} \bignorm{\bo \Delta_{\mb X}\t \mb X + \mb X\t {\bo \Delta_{\mb X}} - \bo \Delta_{\mb Y}\t \mb Y - \mb Y\t {\bo \Delta_{\mb Y}} }\fb^2    \Big]\\
                & + c_\perp \frac{ \pi \sigma_r^*}{n}\bignorm{\mb 1_n\t \bo \Delta_{\mb X}}_2^2  - 
                c_0 \underline c_\psi \Big| \bignorm{\mc P_{\Omega}(\mb 1_n\bo \Delta_{\bo \zeta}\t + \bo \Delta_{\mb X}\mb Y\t + \mb X \bo \Delta_{\mb Y}\t)}\fb^2\\
    &\quad -  \pi  \bignorm{\mb 1_n\bo \Delta_{\bo\zeta}\t  + \bo \Delta_{\mb X}\mb Y\t + \mb X \bo \Delta_{\mb Y}\t }\fb^2 \Big|\\
    & \geq   c_0 \pi  \underline c_\psi \Big[
                 \underbrace{ \bignorm{\bo \Delta_{\mb X}\mb Y\t + \mb X \bo \Delta_{\mb Y}\t }\fb^2 + \frac{1}{4} \norm{\bo \Delta_{\mb X}\t \mb X + \mb X\t {\bo \Delta_{\mb X}} - \bo \Delta_{\mb Y}\t \mb Y - \mb Y\t {\bo \Delta_{\mb Y}}}\fb^2}_{\alpha_{1,1}}
                \Big]  \\ 
                & + c_0 \pi  \underline c_\psi \Big[\underbrace{ \bignorm{\mb 1_n\bo \Delta_{\bo \zeta}\t }\fb^2+  \frac{c_\perp \sigma_r^*}{\underline c_\psi n} \bignorm{\mb 1_n\t \bo \Delta_{\mb X}}_2^2+ 2 \big\langle \mb 1_n\bo \Delta_{\bo \zeta}\t  , \bo \Delta_{\mb X}\mb Y\t \big\rangle + 2 \big\langle \mb 1_n\bo \Delta_{\bo \zeta}\t, \mb X \bo \Delta_{\mb Y}\t \big\rangle  }_{\alpha_{1,2}} \Big] \\
                &  - c_0 \underline c_\psi \underbrace{  \Big|\bignorm{\mc P_{\Omega}(\mb 1_n\bo \Delta_{\bo \zeta}\t + \bo \Delta_{\mb X}\mb Y\t + \mb X \bo \Delta_{\mb Y}\t)}\fb^2 -  \pi  \bignorm{\mb 1_n\bo \Delta_{\bo\zeta}\t  + \bo \Delta_{\mb X}\mb Y\t + \mb X \bo \Delta_{\mb Y}\t }\fb^2 \Big|}_{\alpha_{1,3}}. 
        \end{align*}
        
        In the sequel, we separately bound each component above. 

    \bigskip 
    \begin{itemize}
        \item 
        \emph{Lower bounding $\alpha_{1,1}$. }
        We examine the first term $\alpha_{1,1}$ first. Expanding the quadratic term gives
        \begin{align}
    & \alpha_{1,1} = \bignorm{\bo \Delta_{\mb X} \mb Y\t }\fb^2 + \bignorm{\mb X \bo \Delta_{\mb Y}\t }\fb^2 + \frac{1}{2}\bignorm{\bo\Delta_{\mb X}\t \mb X}\fb^2 + \frac12 \bignorm{\bo \Delta_{\mb Y}\t \mb Y}\fb^2\\
    & +\frac{1}{2}\bigip{\bo \Delta_{\mb X}\t \mb X, \mb X\t \bo \Delta_{\mb X}} + \frac{1}{2}\bigip{\bo \Delta_{\mb Y}\t \mb Y, \mb Y\t \bo \Delta_{\mb Y}}  +\bigip{\mb X\t \bo \Delta_{\mb X}, \bo \Delta_{\mb Y}\t \mb Y} -  \bigip{\mb X\t \bo \Delta_{\mb X}, \mb Y\t \bo \Delta_{\mb Y}}\\
    & =   \bignorm{\bo \Delta_{\mb X} \mb Y\t }\fb^2 + \bignorm{\mb X \bo \Delta_{\mb Y}\t }\fb^2 + \frac{1}{2}\bignorm{\bo\Delta_{\mb X}\t \mb X- \bo \Delta_{\mb Y}\t \mb Y}\fb^2 + \frac12 \big\langle\bo \Delta_{\mb X}\t \mb X + \bo \Delta_{\mb Y}\t \mb Y, \mb X\t \bo\Delta_{\mb X}\\
    &\quad + \mb Y\t \bo \Delta_{\mb Y} \big\rangle .
\end{align}

A useful observation to deal with the last inner-product term is that, if $\bo \Delta_{\mb X}\t \mb X + \bo \Delta_{\mb Y}\t \mb Y$ can be replaced by a symmetric version, then the inner product becomes the Frobenius norm of the surrogate quantity. Recall that, by the definition of $\bo \Delta$, there exist $\mb F_1$ and $\mb F_2$ such that $\bo \Delta_{\mb F} = \mb F_2 \hat{\mb O} - \mb F_1$ where $\hat{\mb O}\coloneqq \argmin_{\mb O \in \mc O(r)} \norm{\mb F_2 \mb O - \mb F_1}\fb$. 
We thus invoke the optimality condition to give that
\begin{align}
    &\hat{\mb O}\t \mb F_2\t (\mb F_2 \hat{\mb O}- \mb F_1) = (\mb F_2 \hat{\mb O}- \mb F_1)\t \mb F_2 \hat{\mb O},
\end{align}
which implies that $\mb F_1\t (\mb F_2 \hat{\mb O} - \mb F_1)$ is symmetric.  As a consequence, we replace $\mb F$ with $\mb F_1$ in the inner product and derive that 
\begin{align}
    & \big\langle \bo \Delta_{\mb X}\t \mb X + \bo \Delta_{\mb Y}\t \mb Y, \mb X\t \bo\Delta_{\mb X} + \mb Y\t \bo \Delta_{\mb Y} \big\rangle \\  
    = & \ip{\bo \Delta_{\mb F}\t \mb F, \mb F\t \bo \Delta_{\mb F}} = \ip{\bo \Delta_{\mb F}\t (\mb F_1 + (\mb F - \mb F_1)), (\mb F_1 + (\mb F - \mb F_1))\t \bo \Delta_{\mb F} } \\  
    = & \ip{\bo \Delta_{\mb F}\t \mb F_1, \bo \Delta_{\mb F}\t \mb F_1} + \ip{\bo \Delta_{\mb F}\t (\mb F - \mb F_1), (\mb F - \mb F_1)\t \bo \Delta_{\mb F}} + 2 \ip{\bo \Delta_{\mb F}\t (\mb F - \mb F_1), \mb F_1\t \bo \Delta_{\mb F}}\\
    \geq&  \bignorm{\mb F_1\t \bo\Delta_{\mb F}}\fb^2 - \norm{(\mb F_1   - \mb F)\t\bo \Delta_{\mb F}}\fb^2 - 2\norm{(\mb F_1- \mb F)\t \bo \Delta_{\mb F}}\fb\norm{\mb F_1\t \bo \Delta_{\mb F}}\fb \\  
    \geq & \big[  - \norm{\mb F_1 - \mb F}^2 - 2 \norm{\mb F_1 - \mb F} \norm{\mb F_1}\big] \norm{\bo \Delta_{\mb F}}\fb^2  \\ 
    \geq & -\big[   \epsilon + 2 \sqrt{\epsilon} \big] \sigma_r^* \norm{\bo \Delta_{\mb F}}\fb^2. 
\end{align}

Therefore, given a sufficiently small $\epsilon$, we obtain the following lower bound for $\alpha_{1,1}$:
\begin{align} 
    & \alpha_{1,1} \geq \frac{1}{8} \sigma_r^* \big(\norm{\bo \Delta_{\mb X}}\fb^2 + \norm{\bo \Delta_{\mb Y}}\fb^2\big). 
\end{align}

\item \emph{Lower Bounding $\alpha_{1,2}$. }
To lower bound the second term $\alpha_{1,2}$, we notice that
\begin{align*}
    & \bignorm{\mb 1_n\bo \Delta_{\bo \zeta}\t }\fb^2 +  \frac{c_{\perp}\sigma_r^*}{\underline c_\psi n} \bignorm{\mb 1_n\t \bo \Delta_{\mb X}}_2^2 + 2 \bigip{\mb 1_n\bo\Delta_{\bo \zeta}\t, \bo \Delta_{\mb X}\mb Y\t }  +  2 \bigip{\mb 1_n\bo \Delta_{\bo \zeta}\t, \mb X \bo \Delta_{\mb Y}\t }\\ 
    & \geq   \bignorm{\mb 1_n\bo \Delta_{\bo \zeta}\t }\fb^2 +   \frac{c_{\perp}\sigma_r^*}{\underline c_\psi n}\bignorm{\mb 1_n\t \bo \Delta_{\mb X}}_2^2 + 2 \bigip{\mb 1_n\t  \bo \Delta_{\mb X}, \bo \Delta_{\bo \zeta}\t \mb Y } - 2 \bignorm{\mb 1_n\t \mb X}_2 \norm{\bo \Delta_{\bo \zeta}}_2 \norm{\bo \Delta_{\mb Y}}\\
    & \geq   \frac12 n \bignorm{\bo \Delta_{\bo \zeta}}_2^2 + \frac{n}{2\norm{\mb Y}\fb^2} \bignorm{\bo \Delta_{\bo \zeta}\t \mb Y}\fb^2 +   \frac{c_{\perp}\sigma_r^*}{\underline c_\psi n}\bignorm{\mb 1_n\t \bo \Delta_{\mb X}}_2^2 + 2\bigip{\mb 1_n\t \bo \Delta_{\mb X}, \bo \Delta_{\bo \zeta}\t \mb Y}\\
    &\quad -  2 \bignorm{\mb 1_n\t \mb X}_2 \norm{\bo \Delta_{\bo \zeta}}_2 \norm{\bo \Delta_{\mb Y}}\\
    & \geq   \frac{1}{2} n\norm{\bo \Delta_{\bo \zeta}}_2^2 + \norm{\frac{\sqrt{n}}{\sqrt{2}\norm{\mb Y}\fb} \bo \Delta_{\bo \zeta}\t \mb Y + \frac{\sqrt{2}\norm{\mb Y}\fb}{\sqrt{n}} \mb 1_n\t \bo \Delta_{\mb X}}_2^2 + \big( \frac{c_{\perp}\sigma_r^*}{\underline c_\psi n}\\
    &\quad - \frac{2\norm{\mb Y}\fb^2}{n}\big) \norm{\mb 1_n\t \bo \Delta_{\mb X}}_2^2\\
    &  -2 \bignorm{\mb 1_n\t \mb X}_2 \norm{\bo \Delta_{\bo \zeta}}_2 \norm{\bo \Delta_{\mb Y}}\\ 
    & \geq   \frac{1}{2} n \norm{\bo \Delta_{\bo \zeta}}_2^2 -  2 \bignorm{\mb 1_n\t \mb X}_2 \norm{\bo \Delta_{\bo \zeta}}_2 \norm{\bo \Delta_{\mb Y}},
\end{align*}
where we invoked the fact $\frac{2\norm{\mb Y}\fb^2}{n} \leq \frac{4\kappa r \sigma_r^*}{n} \leq \frac{c_{\perp} \sigma_r^*}{\underline c_\psi n}$ from the condition on $\mb Y$ and $c_\perp$. 

\item 
\emph{Upper bounding $\alpha_{1,3}$. }
Now it remains to upper bound the third term $\alpha_{1,3}$ defined in \eqref{eq: quadratic form of the augmented function}, i.e., establishing the approximate injectivity along the direction in interest. We begin by the following decomposition: 
\begin{align*}
    &  \Big|\norm{\mc P_{\Omega}(\mb 1_n\bo \Delta_{\bo \zeta}\t + \bo \Delta_{\mb X}\mb Y\t + \mb X \bo \Delta_{\mb Y}\t)}\fb^2  -  \pi  \norm{\mb 1_n\bo \Delta_{\bo\zeta}\t  + \bo \Delta_{\mb X}\mb Y\t + \mb X \bo \Delta_{\mb Y}\t }\fb^2 \Big| \\ 
    \leq & \underbrace{
     \Big|\norm{\mc P_{\Omega}(\mb 1_n\bo \Delta_{\bo \zeta}\t + \bo \Delta_{\mb X}\mb Y\t + \mb X \bo \Delta_{\mb Y}\t)}\fb^2 - 
                \norm{\mc P_{\Omega}(\mb 1_n\bo \Delta_{\bo \zeta}\t + \bo \Delta_{\mb X}{\mb Y^*}\t + \mb X^* \bo \Delta_{\mb Y}\t)}\fb^2}_{\alpha_{1,3,1}} \\ 
    & +\pi  \underbrace{\Big| 
    \norm{\mb 1_n\bo \Delta_{\bo\zeta}\t  + \bo \Delta_{\mb X}\mb Y\t + \mb X \bo \Delta_{\mb Y}\t }\fb^2 - \bignorm{\mb 1_n\bo \Delta_{\bo\zeta}\t  + \bo \Delta_{\mb X}{\mb Y^*}\t + \mb X^* \bo \Delta_{\mb Y}\t }\fb^2
    \Big|}_{\alpha_{1,3,2}} \\ 
    & + \underbrace{\Big| \norm{\mc P_\Omega\big(\mb 1_n\bo \Delta_{\bo\zeta}\t  + \bo \Delta_{\mb X}{\mb Y^*}\t + \mb X^* \bo \Delta_{\mb Y}\t \big)}\fb^2- \pi \bignorm{\mb 1_n\bo \Delta_{\bo\zeta}\t  + \bo \Delta_{\mb X}{\mb Y^*}\t + \mb X^* \bo \Delta_{\mb Y}\t }\fb^2\Big|}_{\alpha_{1,3,3}}.
    \label{eq: decomposition of sampling deviation}
\end{align*}

The terms $\alpha_{1,3,i}$, $i=1,2,3$ are controlled below. 
\begin{itemize}
    \item For $\alpha_{1,3,1}$, we first consider the following identity: 
    \begin{align*}
        & \sum_{(i,j)\in \Omega}\big((\bo \Delta_{\bo \zeta})_j + 
                 \bigip{\bo \Delta_{\mb X} \mb Y\t + \mb X {\bo \Delta_{\mb Y}}\t , e_i^n{e_j^p}\t }    
                \big)^2 - 
                \big((\bo \Delta_{\bo \zeta})_j\\
    &\quad + 
                 \bigip{\bo \Delta_{\mb X} {\mb Y^*}\t + \mb X^* {\bo \Delta_{\mb Y}}\t , e_i^n{e_j^p}\t }    
                \big)^2   \\ 
    & =   \sum_{(i,j) \in \Omega}\bigip{\bo \Delta_{\mb X}(\mb Y - \mb Y^*)\t + (\mb X - \mb X^*)\bo \Delta_{\mb Y}\t , e_i^n {e_j^p}\t }^2\\
                & + 2 \sum_{(i,j)\in \Omega} \Big[(\bo \Delta_{\bo \zeta})_j + 
                 \bigip{\bo \Delta_{\mb X} {\mb Y^*}\t + \mb X^* {\bo \Delta_{\mb Y}}\t , e_i^n{e_j^p}\t }    
                \Big]\\
    &\quad \cdot \bigip{\bo \Delta_{\mb X}(\mb Y - \mb Y^*)\t + (\mb X - \mb X^*)\bo \Delta_{\mb Y}\t , e_i^n {e_j^p}\t }.
    \end{align*}
    Invoking this identity together with the Cauchy-Schwarz inequality yields that 
    \begin{align}
         & \alpha_{1,3,1} \leq 3\Big[\sum_{(i,j)\in \Omega}\big( (\bo \Delta_{\mb X})_{i, \cdot} (\mb Y - \mb Y^*)_j \big)^2\Big]  + 3 \Big[\sum_{(i,j) \in \Omega}\big( (\mb X - \mb  X^*)_{i, \cdot} (\bo \Delta_{\mb Y})_j \big)^2 \Big] \\ 
    & + 2 \bignorm{\mc P_{\Omega}\big(\mb 1_n{\bo \Delta}_{\bo \zeta}\t + \bo \Delta_{\mb X} {\mb Y^*}\t + \mb X^*{\bo \Delta_{\mb Y}}\t\big) }\fb\\
    &\quad \cdot \bignorm{\mc P_{\Omega}\big(\bo \Delta_{\mb X}(\mb Y - \mb Y^*)\t + (\mb X - \mb X^*)\bo \Delta_{\mb Y}\t \big)}\fb.
         \label{eq: decomposition of alpha131}
    \end{align}
    For the first term in \eqref{eq: decomposition of alpha131}, one has 
    \begin{align}
        & \sum_{(i,j) \in \Omega} \big( (\bo \Delta_{\mb X})_{i, \cdot} (\mb Y - \mb Y^*)_{j, \cdot}\t \big)^2 
        \leq  \sum_{i \in[n]} \norm{(\bo \Delta_{\mb X})_{i }}_2^2 \cdot \big(\sum_{j\in[p]} \ind\{(i,j)\in \Omega\} \big) \cdot \norm{\mb Y - \mb Y^*}\ti ^2 \\ 
    & \leq   \max_{i \in[n]} \Big(\sum_{j\in[p]} \ind\{(i,j)\in \Omega\} \Big) \norm{\bo \Delta_{\mb X}}\fb^2 \norm{\mb Y - \mb Y^*}\ti^2  
        \leq  2\pi p\norm{\bo \Delta_{\mb X}}\fb^2 \norm{\mb Y - \mb Y^*}\ti^2
        \label{eq: first component of alpha131}
    \end{align}
    with probability at least $1-  O(d^{-c-1})$. 

    Similarly, we have for the second term in \eqref{eq: decomposition of alpha131} that 
    \begin{equation}
        \sum_{(i,j)\in \Omega}\big( (\mb X - \mb  X^*)_{i, \cdot} (\bo \Delta_{\mb Y})_{j} \big)^2   \leq 2\pi n \norm{\bo \Delta_{\mb Y}}\fb^2 \norm{\mb X - \mb X^*}\ti^2
        \label{eq: second component of alpha131}
    \end{equation}
    holds with probability at least $1- O(d^{-c-1})$. 

    To control the remaining factor, define
    \begin{align*}
        \widetilde{\mb U}^*
        &\coloneqq \big[n^{-1/2}\mb 1_n,\mb U^*\big],\\
        \mc T_{\sf aug}^*
        &\coloneqq
        \Big\{
        \widetilde{\mb U}^*\mb A\t+\mb B{\mb V^*}\t:
        \mb A\in\bb R^{p\times(r+1)},
        \mb B\in\bb R^{n\times r}
        \Big\}.
    \end{align*}
    The centering condition ${\mb 1_n}\t\mb X^*=\mb 0$ ensures that $\widetilde{\mb U}^*$ has orthonormal columns. Moreover, the incoherence of $\mb U^*$ and $\mb V^*$ implies that the augmented tangent space $\mc T_{\sf aug}^*$ has incoherence of order $\mu$. In particular, $\mb 1_n{\bo \Delta}_{\bo \zeta}\t
    +\bo \Delta_{\mb X}{\mb Y^*}\t
    +\mb X^*{\bo \Delta_{\mb Y}}\t
    \in \mc T_{\sf aug}^*$.
    The standard Bernoulli-sampling concentration bound for a fixed tangent space, using the rectangular extension of \citet[Theorem~4.1 and Corollary~4.3]{candes2012exact}, therefore gives
    \begin{align}
    & \bignorm{\mc P_{\Omega}\big(\mb 1_n{\bo \Delta}_{\bo \zeta}\t + \bo \Delta_{\mb X} {\mb Y^*}\t + \mb X^*{\bo \Delta_{\mb Y}}\t\big) }\fb \lesssim \sqrt{\pi } \Big[ \sqrt{n}\norm{\bo \Delta_{\bo \zeta}}_2 + \sqrt{\sigma_1^*}\big(\norm{\bo \Delta_{\mb X}}\fb\\
    &\quad + \norm{\bo \Delta_{\mb Y}}\fb \big) \Big]
    \end{align}
    with probability at least $1-O(d^{-c-1})$, provided that
    $\mu r\log d\ll \pi(n\wedge p)$.
    Besides, it follows directly from the triangle inequality, \eqref{eq: first component of alpha131}, and \eqref{eq: second component of alpha131} that 
    \begin{align}
        & \bignorm{\mc P_{\Omega}\big(\bo \Delta_{\mb X}(\mb Y - \mb Y^*)\t + (\mb X - \mb X^*)\bo \Delta_{\mb Y}\t \big)}\fb \\ 
        \lesssim & \sqrt{\pi}\Big[\sqrt{p}\norm{\bo \Delta_{\mb X}}\fb \norm{\mb Y - \mb Y^*}\ti + \sqrt{n}\norm{\bo \Delta_{\mb Y}}\fb \norm{\mb X - \mb X^*}\ti\Big]
    \end{align}
    with probability exceeding $1 - O(d^{-c-1})$. 

    Taking these bounds collectively, we finally conclude that 
    \begin{align*}
        & \alpha_{1,3,1} \lesssim \pi p\norm{\bo \Delta_{\mb X}}\fb^2 \norm{\mb Y - \mb Y^*}\ti^2 + \pi n \norm{\bo \Delta_{\mb Y}}\fb^2 \norm{\mb X - \mb X^*}\ti^2  \\ 
        & \quad + \sqrt{\pi}\Big[\sqrt{p}\norm{\bo \Delta_{\mb X}}\fb \norm{\mb Y - \mb Y^*}\ti + \sqrt{n}\norm{\bo \Delta_{\mb Y}}\fb \norm{\mb X - \mb X^*}\ti\Big] \\ 
    & \quad \cdot \sqrt{\pi } \Big[ \sqrt{n}\norm{\bo \Delta_{\bo \zeta}}_2 + \sqrt{\sigma_1^*}\big(\norm{\bo \Delta_{\mb X}}\fb + \norm{\bo \Delta_{\mb Y}}\fb \big) \Big]\\
    & \ll   \pi \big[ \norm{\bo \Delta_{\mb X}}\fb^2 +  \norm{\bo \Delta_{\mb Y}}\fb^2 + \bignorm{ (n / \sigma_r^*)^{\frac12}\bo \Delta_{\bo \zeta}}\fb^2 \big]\big[\sqrt{p}\norm{\mb Y - \mb Y^*}\ti\\
    &\quad + \sqrt{n} \norm{\mb X - \mb X^*}\ti \big]\sqrt{\sigma_1^* }\\
    & \lesssim   \pi \norm{\bo \Delta}_2^2 \sigma_r^*.
    \end{align*}
    with probability at least $1- O(d^{-c-1})$, provided the condition 
    \eq{
    \sqrt{n}\norm{\mb X - \mb X^*}\ti \vee \sqrt{p}\norm{\mb Y - \mb Y^*}\ti \leq \sqrt{\epsilon \sigma_1^* } / \kappa. 
    \label{eq: condition on the rowwise consistency of X and Y}
    }

    \item For the second term $\alpha_{1,3,2}$, applying the triangle inequality and the Cauchy-Schwarz inequality, we deduce that 
    \eq{
\begin{aligned}
    & \alpha_{1,3,2} \lesssim  \big( \norm{\bo \Delta_{\mb X}}\fb^2 +  \norm{\bo \Delta_{\mb Y}}\fb^2 + \bignorm{(n / \sigma_r^*)^{\frac12}\bo \Delta_{\bo \zeta}}\fb^2 \big)\big(\sqrt{p}\norm{\mb Y - \mb Y^*}\ti\\
    &\quad + \sqrt{n} \norm{\mb X - \mb X^*}\ti \big)\sqrt{\sigma_1^*} 
         \ll  \norm{\bo \Delta}_2^2 \sigma_r^* 
\end{aligned}
    }
    with probability at least $1- O(d^{-c-1})$. 

    \item Regarding the last term $\alpha_{1,3,3}$, applying the same restricted-isometry bound on $\mc T_{\sf aug}^*$ yields with probability at least $1 - O(d^{-c-1})$ that
    \begin{align}
        & \alpha_{1,3,3} = \Big| \bignorm{\mc P_{ \Omega}\big(\mb 1_n\bo \Delta_{\bo \zeta}\t + \mb X^* \bo \Delta_{\mb Y}\t + \bo \Delta_{\mb X}{ \mb Y^*}\t \big)}\fb^2 - \pi \bignorm{\mb 1_n\bo \Delta_{\bo \zeta}\t + \mb X^* \bo \Delta_{\mb Y}\t + \bo \Delta_{\mb X}{ \mb Y^*}\t }\fb^2 \Big| \\
    & \lesssim   \sqrt{\frac{\pi \mu d r  \log d}{ np }} \Big( \sigma_1^* \big( \norm{\bo \Delta_{\mb X}}\fb^2 +\norm{\bo \Delta_{\mb Y}}\fb^2\big) + n \norm{\bo \Delta_{\bo \zeta}}_2^2 \Big).
    \end{align}
\end{itemize}

Taken together, these bounds yield
\begin{equation}
    \alpha_{1,3} / \norm{\bo \Delta}_2^2 \ll \pi \sigma_r^*,
\end{equation}
with probability at least $1-O(d^{-c-1})$, provided that
$\mu r\kappa^2\log d\ll \pi(n\wedge p)$.

    \end{itemize}

Putting the pieces together, we deduce that, with probability at least $1- O(d^{-c-1})$, 
\begin{align}
    & \alpha_1 / (c_0 \underline c_\psi) \geq    \frac{1}{8} \pi \sigma_r^* \big(\norm{\bo \Delta_{\mb X}}\fb^2 + \norm{\bo \Delta_{\mb Y}}\fb^2\big) +  \frac{1}{2} \pi n \norm{\bo \Delta_{\bo \zeta}}_2^2\\
    & \qquad - 2 \pi \bignorm{\mb 1_n\t \mb X}_2\norm{\bo \Delta_{\bo \zeta}}_2 \norm{\bo \Delta_{\mb Y}}\fb - \frac{1}{32} \pi \sigma_r^* \norm{\bo \Delta}\fb^2 \\ 
    & \geq   \frac{1}{8} \pi \sigma_r^* \big(\norm{\bo \Delta_{\mb X}}\fb^2 + \norm{\bo \Delta_{\mb Y}}\fb^2\big) +  \frac{1}{2} \pi n \norm{\bo \Delta_{\bo \zeta}}_2^2 - \frac{1}{32} \pi \sqrt{n\sigma_r^*}\norm{\bo \Delta_{\bo \zeta}}_2 \norm{\bo \Delta_{\mb Y}}\fb\\
    & \qquad - \frac{1}{32} \pi \sigma_r^* \big( \norm{\bo \Delta_{\mb X}}\fb^2 + \norm{\bo \Delta_{\mb Y}}_2^2 + \frac{n}{\sigma_r^*}\norm{\bo \Delta_{\bo \zeta}}_2^2\big) \\ 
    & \geq   \frac{1}{32} \pi \sigma_r^* \big(\norm{\bo \Delta_{\mb X}}\fb^2 + \norm{\bo \Delta_{\mb Y}}\fb^2 + \frac{n}{\sigma_r^*} \norm{\bo \Delta_{\bo \zeta}}_2^2\big) + \frac{1}{64} \pi \big(\sqrt{n} \norm{\bo \Delta_{\bo \zeta}}_2 - \sqrt{\sigma_r^*}\norm{\bo \Delta_{\mb Y}}\fb\big)^2\\
    & \geq   \frac{1}{32} \pi \sigma_r^*\norm{\bo \Delta}\fb^2,
    \label{eq: bound for alpha1}
\end{align}
where we used the condition $\bignorm{\mb 1_n\t \mb X}_2 \lesssim \sqrt{\epsilon n \sigma_r^*}$ for some sufficiently small $\epsilon$ in the first inequality. 

\bigskip

We are left to control the second term $\alpha_2$ in \eqref{eq: quadratic form of the augmented function}. 
We examine the first term in $\alpha_2$ in \eqref{eq: quadratic form of the augmented function} by the following decomposition: 
\begin{align}
    & \Big|  \sum_{(i,j)\in \Omega }\big(\psi(\mc M(\bo \theta)_{i,j}) - R_{i,j } \big) \bigip{e_i^n{e_j^p}\t, \bo \Delta_{\mb X}{\bo \Delta_{\mb Y}}\t } \Big|\\
    &\quad = \big| \bigip{\mc P_\Omega(\psi(\mc M(\bo \theta)) - \mb R), \mc P_\Omega(\bo \Delta_{\mb X} \bo \Delta_{\mb Y}\t)} \big|\\
    & \leq   \big| \bigip{\mc P_\Omega(\mb E), \mc P_\Omega(\bo \Delta_{\mb X} \bo \Delta_{\mb Y}\t)} \big| + \big| \bigip{\mc P_\Omega(\psi(\mc M(\bo \theta)) - \psi(\mc M(\bo \theta^*))), \mc P_\Omega(\bo \Delta_{\mb X} \bo \Delta_{\mb Y}\t)} \big|
    \label{eq: gradient term in hessian analysis}
\end{align}

Regarding the first term on the RHS of \eqref{eq: gradient term in hessian analysis}, we invoke Lemma~\ref{lemma: noise spectral norm concentration} to obtain with probability at least $1 - O(d^{-c-1})$ that
\begin{align}
    &  \big| \bigip{\mc P_\Omega(\mb E), \mc P_\Omega(\bo \Delta_{\mb X} \bo \Delta_{\mb Y}\t)} \big| \leq \norm{\mc P_{\Omega}(\mb E)} \norm{\bo \Delta_{\mb F}}\fb^2 \lesssim \sigma \sqrt{\pi d} \norm{\bo \Delta_{\mb F}}\fb^2 . 
\end{align}
Regarding the second term therein, we make use of the following decomposition
\begin{align}
    & \Big|\bigip{\mc P_\Omega(\psi(\mc M(\bo \theta)) - \psi(\mc M(\bo \theta^*))), \mc P_\Omega(\bo \Delta_{\mb X} \bo \Delta_{\mb Y}\t)}\Big| \\ 
    \leq & \Big|\bigip{(\bo \Omega - \pi \mb 1_n \mb 1_p\t )\circ (\psi(\mc M(\bo \theta)) - \psi(\mc M(\bo \theta^*))), \bo \Delta_{\mb X} \bo \Delta_{\mb Y}\t}\Big| \\ 
    &  + \pi \Big|\bigip{\psi(\mc M(\bo \theta)) - \psi(\mc M(\bo \theta^*)), \bo \Delta_{\mb X} \bo \Delta_{\mb Y}\t}| . 
\end{align}
For $s\in[0,1]$, define
$\mb H_s\coloneqq\psi'\big((1-s)\mc M(\bo\theta^*)
+s\mc M(\bo\theta)\big)$ and
$\overline{\mb H}\coloneqq\int_0^1\mb H_s\,\mathrm ds$.
The RIC row bounds and the domain verification above imply that
$\mb H_s\in\mc M_\psi$ for every $s\in[0,1]$.
The Lemma~\ref{lemma: injectivity of P_Omega under nonlinearity} bound
\eqref{eq: uniform nonlinear sampling operator bound}, applied with endpoints
$\bo\theta$ and $\bo\theta^*$, gives
$\bignorm{\mc P_\Omega(\overline{\mb H})-\pi\overline{\mb H}}
\lesssim\sigma^2\xi^4\sqrt{\pi d}$.
The mean-value theorem yields
$\psi(\mc M(\bo\theta))-\psi(\mc M(\bo\theta^*))
=\overline{\mb H}\circ\big[\mc M(\bo\theta)-\mc M(\bo\theta^*)\big]$.
Consequently, Lemma~\ref{lemma: lemma 4.4 chen2019} and
\eqref{eq: uniform nonlinear sampling operator bound} imply
\begin{align*}
    & \Big|\bigip{(\bo \Omega - \pi \mb 1_n \mb 1_p\t )
    \circ \big[\psi(\mc M(\bo \theta))-\psi(\mc M(\bo \theta^*))\big],
    \bo \Delta_{\mb X} \bo \Delta_{\mb Y}\t}\Big|  \\ 
    ={}& \Big|\bigip{
    \mc P_\Omega(\overline{\mb H})-\pi\overline{\mb H},
    \big[\mc M(\bo \theta)-\mc M(\bo \theta^*)\big]
    \circ (\bo \Delta_{\mb X} \bo \Delta_{\mb Y}\t)} \Big| \\
    \lesssim & \sigma^2 \xi^4  \sqrt{\pi d} \Big[\norm{\bo \zeta - \bo \zeta^*}_\infty  + \sqrt{\frac{\mu r\sigma_1^* }{n}} \norm{\mb Y - \mb Y^*}\ti+ \sqrt{\frac{\mu r \sigma_1^*}{p}}\norm{\mb X - \mb X^*}\ti \Big] \norm{\bo \Delta_{\mb F}}\fb^2,
\end{align*}
The second term can be directly bounded by
\begin{align}
    & \pi \Big|\bigip{\psi(\mc M(\bo \theta)) - \psi(\mc M(\bo \theta^*)), \bo \Delta_{\mb X} \bo \Delta_{\mb Y}\t}\Big|
    \lesssim \sqrt{\epsilon}\,\sigma^2\pi\sigma_r^*
    \norm{\bo\Delta_{\mb F}}\fb^2
    \ll \underline c_\psi\pi\sigma_r^*
    \norm{\bo\Delta_{\mb F}}\fb^2,
\end{align}
where the last inequality follows from
$\sqrt{\epsilon}\,\kappa_\psi=\sqrt{\epsilon_0}\ll1$.

Consequently, combining these two bounds concludes that, with probability at least $1- O(d^{-c-1})$, 
\begin{align*}
    & \Big|  \sum_{(i,j)\in \Omega }\big(\psi(\mc M(\bo \theta)_{i,j}) - R_{i,j } \big) \bigip{e_i^n{e_j^p}\t, \bo \Delta_{\mb X}{\bo \Delta_{\mb Y}}\t } \Big|
    \lesssim   \Big[\sigma \sqrt{\pi d} + \sigma^2\xi^4 \sqrt{\pi d} \big( \norm{\bo \zeta - \bo \zeta^*}_\infty  \\
    & \qquad   +  \sqrt{\frac{\mu r\sigma_1^* }{n}}\norm{\mb Y - \mb Y^*}\ti + \sqrt{\frac{\mu r\sigma_1^* }{p}} \norm{\mb X - \mb X^*}\ti\big) \Big] \norm{\bo \Delta_{\mb F}}\fb^2 + \frac{c_0 \underline c_\psi \pi \sigma_r^*}{512} \norm{\bo \Delta_{\mb F}}\fb^2 \\
    \leq & \frac{c_0 \underline c_\psi  \pi \sigma_r^*}{256}  \norm{\bo \Delta_{\mb F}}\fb^2. 
\end{align*}
due to the conditions \eqref{eq: rowwise consistency of parameters in geometry lemma}, $\mu r\kappa_\psi \xi^8 \ll \pi (n\wedge p)$, and $\sigma \sqrt{\pi d} \leq \bar \sigma \sqrt{\pi d}\ll \underline c_\psi \pi \sigma_r^*$. 

On the other hand, the second term can be directly controlled as follows: 
\begin{align*}
    & \Big| c_{\sf{aug}} \pi \bigip{\mb X\t \mb X - \mb Y\t \mb Y, \bo \Delta_{\mb X}\t \bo \Delta_{\mb X} -\bo \Delta_{\mb Y}\t \bo \Delta_{\mb Y}} \Big|\\ 
    & \leq   c_{\sf{aug}} \pi \Big(\norm{\mb X - \mb X^*}^2 + 2 \norm{\mb X - \mb X^*} \norm{\mb X^*} + \norm{\mb Y - \mb Y^*}^2 + 2 \norm{\mb Y - \mb Y^*} \norm{\mb X^*}\Big) \norm{\bo \Delta}\fb^2\\
    &\quad \leq \frac{c_0\underline c_\psi \pi \sigma_r^*}{512}  \norm{\bo \Delta}\fb^2
\end{align*}
by taking $\epsilon_0$ sufficiently small, where we used the fact that ${\mb X^*}\t \mb X^* = {\mb Y^*}\t \mb Y^*$ and the conditions on $\norm{\mb X - \mb X^*}$ and $\norm{\mb Y - \mb Y^*}$. Accounting for the factors $2$ and $4$ in the definition of $\alpha_2$ then gives
\begin{align}
    \big|\alpha_2\big| \leq \frac{c_0\underline c_\psi\pi\sigma_r^*}{64}\norm{\bo\Delta}_2^2.
    \label{eq: bound for alpha2}
\end{align}
This bound holds uniformly over arbitrary directions $\bo\Delta$ with probability at least $1-O(d^{-c-1})$. Since $\alpha_1\geq0$, it also gives $\nabla^2L^{\sf{aug}}(\bo\theta)\succeq-\underline c_{\sf{Hess}}\pi\sigma_r^*\mb I$ throughout the RIC.

Putting \eqref{eq: bound for alpha1} and \eqref{eq: bound for alpha2} together, we finally arrive at the conclusion that, with probability exceeding $1 -O(d^{-c-1})$, 
\begin{align}
    & \bo \Delta\t  \nabla^2 L^{\sf{aug}}(\bo \theta) \bo \Delta \geq \frac{c_0 \underline c_\psi \pi \sigma_r^*}{64} \norm{\bo \Delta}_2^2. 
\end{align} 

For the upper bound, apply the sampling estimates used for $\alpha_{1,3}$
with $\psi'\leq\sigma^2$, and use
$c_{\sf{aug}}\asymp\underline c_\psi\leq\sigma^2$ and
\eqref{eq: bound for alpha2}. Separating the centering penalty gives,
uniformly over arbitrary vectors $\mb v$,
\begin{align*}
    \big|\mb v\t\nabla^2L^{\sf{aug}}(\bo\theta)\mb v\big|
    &\leq\big(2c_\perp\pi\sigma_r^*
    +C\kappa r\sigma^2\pi\sigma_1^*\big)\norm{\mb v}_2^2\\
    &\leq4c_\perp\pi\sigma_1^*\norm{\mb v}_2^2,
\end{align*}
where the last inequality uses $c_\perp\gg\kappa r\sigma^2$.
This proves \eqref{eq: local Hessian upper bound} with
$\bar c_{\sf{Hess}}=4$.

Finally, suppose that $\mb F_1\t\bo\Delta_{\mb F}$ is symmetric. For all sufficiently small $u>0$, the identity matrix is a Procrustes alignment of $\mb F_1+u\bo\Delta_{\mb F}$ with $\mb F_1$, because $(\mb F_1+u\bo\Delta_{\mb F})\t\mb F_1$ is symmetric positive definite. Applying the preceding bounds to $u\bo\Delta$ and using homogeneity proves the tangent-direction claim.

\qed

\subsubsection{Proof of Lemma~\ref{lemma: pipeline satisfies GD initialization}}
Theorem~\ref{thm: detailed USVT} supplies the global and leave-one-out inputs
required by Theorem~\ref{thm: general OS}; in particular, one may take
$\xi^0\lesssim\kappa\lambda^{\sp}\xi/\sqrt\pi$ and
$\xi^{0,\sf{loo}}\lesssim\xi^{\sf{loo}}$. The choice of $\lambda^{\sp}$ in
Theorem~\ref{thm: USVT}, together with the sampling and signal-strength bounds
in Assumption~\ref{assumption: signal strength and incoherence degree}(a),
guarantees the initialization and smallness requirements in
Theorem~\ref{thm: general OS}, as verified in its specialization to the
concrete pipeline. Applying that
theorem first controls the unprocessed refinement. We next verify that the map
in \eqref{eq: centering and balancing map} preserves these guarantees at the
rates claimed in \eqref{eq: GD initialization rates from pipeline}.

We use the following deterministic observation. Consider factors
$\widetilde{\mb X}$ and $\widetilde{\mb Y}$ that admit a common alignment
$\mb O$ with $(\mb X^*,\mb Y^*)$, and set
\begin{equation*}
    \varepsilon_{2,\infty}
    \coloneqq
    \bignorm{\mb J_n\widetilde{\mb X}\mb O-\mb X^*}\ti
    \vee
    \bignorm{\widetilde{\mb Y}\mb O-\mb Y^*}\ti,
    \qquad
    \varepsilon_F
    \coloneqq
    \bignorm{\mb J_n\widetilde{\mb X}\mb O-\mb X^*}\fb
    \vee
    \bignorm{\widetilde{\mb Y}\mb O-\mb Y^*}\fb.
\end{equation*}
If $\varepsilon_F\ll\sqrt{\sigma_r^*}$, the balanced factors
$[\mb X^{\sf c},\mb Y^{\sf c}]$ obtained from the SVD in
\eqref{eq: centering and balancing map} satisfy
\begin{align}
    &\bignorm{\mb X^{\sf c}\mb O^{\sf c}-\mb X^*}\ti
    \vee
    \bignorm{\mb Y^{\sf c}\mb O^{\sf c}-\mb Y^*}\ti \nonumber\\
    &\qquad\lesssim
    \varepsilon_{2,\infty}
    +\norm{\mb F^*}\ti
    \frac{\kappa^{3/2}}{\sigma_r^*}
    \bigl(\sqrt{\sigma_1^*}\,\varepsilon_F+\varepsilon_F^2\bigr)
    \label{eq: rowwise stability of centering and balancing}
\end{align}
for the common Procrustes rotation $\mb O^{\sf c}$.
To verify this claim, replace the preliminary factors by their commonly
aligned versions and, with a slight abuse of notation, continue to denote them
by $\widetilde{\mb X}$ and $\widetilde{\mb Y}$. Write
\begin{equation*}
    \mb G_{\mb X}
    \coloneqq
    (\mb J_n\widetilde{\mb X})\t\mb J_n\widetilde{\mb X},
    \qquad
    \mb G_{\mb Y}
    \coloneqq
    \widetilde{\mb Y}\t\widetilde{\mb Y},
    \qquad
    \mb A
    \coloneqq
    \mb G_{\mb X}^{-1/2}
    \big(
        \mb G_{\mb X}^{1/2}\mb G_{\mb Y}\mb G_{\mb X}^{1/2}
    \big)^{1/4}.
\end{equation*}
The pair $[\mb J_n\widetilde{\mb X}\mb A,\,
\widetilde{\mb Y}\mb A^{-\top}]$ is a balanced factorization of
$\mb J_n\widetilde{\mb X}\widetilde{\mb Y}\t$ and hence agrees, up to a
common orthogonal transformation, with its SVD factors. Let
$\mb G^*\coloneqq{\mb X^*}\t\mb X^*={\mb Y^*}\t\mb Y^*$. Direct expansion gives
$\norm{\mb G_{\mb X}-\mb G^*}
\vee\norm{\mb G_{\mb Y}-\mb G^*}
\lesssim\sqrt{\sigma_1^*}\,\varepsilon_F+\varepsilon_F^2$.
The assumed smallness of $\varepsilon_F$ places the spectra of both Gram
matrices in a fixed enlargement of $[\sigma_r^*,\sigma_1^*]$. The resolvent
representations for positive-definite matrix powers therefore imply
\begin{align*}
    \norm{\mb G_{\mb X}^{-1/2}-{\mb G^*}^{-1/2}}
    &\lesssim
    \frac{\sqrt{\sigma_1^*}\,\varepsilon_F+\varepsilon_F^2}
    {(\sigma_r^*)^{3/2}},\\
    \norm{
      (\mb G_{\mb X}^{1/2}\mb G_{\mb Y}\mb G_{\mb X}^{1/2})^{1/4}
      -(\mb G^*)^{1/2}
    }
    &\lesssim
    \frac{\kappa^{3/2}}
    {\sqrt{\sigma_r^*}}
    \bigl(\sqrt{\sigma_1^*}\,\varepsilon_F+\varepsilon_F^2\bigr).
\end{align*}
Expanding $\mb A-\mb I_r$, applying these two bounds, and using
$\norm{(\mb G^*)^{1/2}}\leq\sqrt{\sigma_1^*}$ gives
$\norm{\mb A-\mb I_r}
\vee\norm{\mb A^{-\top}-\mb I_r}
\lesssim\frac{\kappa^{3/2}}{\sigma_r^*}
\bigl(\sqrt{\sigma_1^*}\,\varepsilon_F+\varepsilon_F^2\bigr)$.
The bound for $\mb A^{-\top}-\mb I_r$ follows from the same display and the
inverse perturbation identity.
Taking row norms proves
\eqref{eq: rowwise stability of centering and balancing} for one common
rotation. Lemma~\ref{lemma: perturbation theory of optimal rotation under F norm}
then transfers the bound to the Procrustes rotation $\mb O^{\sf c}$, with the
same order.

For the actual refinement, Theorem~\ref{thm: general OS} and
$\mb 1_n\t\mb X^*=0$ yield
$\varepsilon_{2,\infty}
\lesssim\xi^1/(\underline c_\psi\sqrt{\pi\sigma_r^*})$ and
$\varepsilon_F\lesssim\sqrt d\,\varepsilon_{2,\infty}$.
Together with
$\norm{\mb F^*}\ti\lesssim
\sqrt{\mu r\sigma_1^*/(n\wedge p)}$, the signal and sampling conditions absorb
the quadratic term in \eqref{eq: rowwise stability of centering and balancing}
and give
\begin{equation}
    \underline c_\psi\sqrt{\pi\sigma_r^*}
    \Bigl\{
    \bignorm{\mb X^0\mb O^0-\mb X^*}\ti
    \vee
    \bignorm{\mb Y^0\mb O^0-\mb Y^*}\ti
    \Bigr\}
    \lesssim
    \kappa^{5/2}\sqrt{\mu r\rho}\,\xi^1.
    \label{eq: processed GD initialization rowwise rate}
\end{equation}
Since the preprocessing leaves the intercept unchanged, its contribution to
\eqref{eq: GD initialization global} follows directly from
Theorem~\ref{thm: general OS}. Together with
\eqref{eq: processed GD initialization rowwise rate}, this proves the global
initialization bound.

It remains to compare the processed actual and leave-one-out estimates. After
aligning each unprocessed leave-one-out pair with the unprocessed actual pair,
the contraction of $\mb J_n$ and the factor-norm bounds give
\begin{align*}
    &\bignorm{
       \mb J_n\widetilde{\mb X}^{0,(l)}
       \widetilde{\mb Y}^{0,(l)\t}
       -
       \mb J_n\widetilde{\mb X}^{0}
       \widetilde{\mb Y}^{0\t}
     }\fb\\
    &\qquad\lesssim
    \sqrt{\sigma_1^*}
    \Bigl\{
        \bignorm{
            \widetilde{\mb X}^{0,(l)}\mb O_{\ur}^{0,(l)}
            -\widetilde{\mb X}^{0}
        }\fb
        +
        \bignorm{
            \widetilde{\mb Y}^{0,(l)}\mb O_{\ur}^{0,(l)}
            -\widetilde{\mb Y}^{0}
        }\fb
    \Bigr\},
    \qquad l\in[n]\cup[-p],
\end{align*}
where $\mb O_{\ur}^{0,(l)}$ denotes the common relative rotation supplied by
Theorem~\ref{thm: general OS}. Applying
Lemma~\ref{lemma: perturbation theory for joint factors}, now with the
processed actual low-rank matrix as the reference, yields
\begin{align*}
    &\bignorm{\mb X^{0,(l)}\mb O^{0,(l)}-\mb X^0}\fb
    \vee
    \bignorm{\mb Y^{0,(l)}\mb O^{0,(l)}-\mb Y^0}\fb\\
    &\qquad\lesssim
    \kappa^{3/2}
    \Bigl\{
        \bignorm{
            \widetilde{\mb X}^{0,(l)}\mb O_{\ur}^{0,(l)}
            -\widetilde{\mb X}^{0}
        }\fb
        +
        \bignorm{
            \widetilde{\mb Y}^{0,(l)}\mb O_{\ur}^{0,(l)}
            -\widetilde{\mb Y}^{0}
        }\fb
    \Bigr\}.
\end{align*}
The intercept estimates are unchanged by the preprocessing, so their
leave-one-out bounds follow directly from Theorem~\ref{thm: general OS}.
The same theorem therefore implies
\eqref{eq: GD initialization row LOO} and
\eqref{eq: GD initialization column LOO} with
$\xi^{\sf{gd},\sf{loo}}$ in
\eqref{eq: GD initialization rates from pipeline}.
Together with the strengthened signal-strength condition and the sampling
condition, these rates also imply
\eqref{eq: epsilon GD initialization requirement}. Indeed, after substituting
the spectral rates into $\xi^1$, its dominant global-error contribution is of
order $\bar\sigma\kappa_\psi\kappa^4\sqrt{\mu r\rho}\,\xi^5$.
Combining this bound with \eqref{eq: GD initialization rates from pipeline}
gives precisely the parameter factors imposed by the signal-strength
condition. Here $\epsilon_0$ is fixed, while the explicit powers of
$\kappa_\psi$ in the same condition ensure the required radius
$\epsilon=\epsilon_0/\kappa_\psi^2$.

Finally, \eqref{eq: centering and balancing map} gives
\eqref{eq: exact GD initialization constraints} by construction. It also gives
$\bo\zeta^0=\widetilde{\bo\zeta}^0$ and
$\mb X^0\mb Y^{0\t}
=\mb J_n\widetilde{\mb X}^{0}\widetilde{\mb Y}^{0\t}$.
The same identities hold for every leave-one-out estimate. Since the
preprocessing is deterministic, it also preserves the row-wise or
column-wise independence furnished by the coupled spectral and refinement
procedures. The global and leave-one-out bounds above, together with the
triangle inequality, verify all proximity bounds in
Condition~\ref{condition: inductive contraction for GD} at $t=0$ after
enlarging the constants $C_1,\ldots,C_8$ if necessary. The exact centering and
balancing identities verify the remaining two bounds.
\qed

\subsection{Proof of Theorem~\ref{thm: general GD}}
\label{subsec: GD proof}
We now introduce the leave-one-out gradient-descent sequences used to prove
Theorem~\ref{thm: general GD}.

    Starting from these leave-one-out initializations, we further define the leave-one-out iterates recursively. For each $l\in [n] \cup [-p]$ and each iteration $t\ge 0$, let $[\bo \zeta^{t,(l)}, \mb X^{t,(l)}, \mb Y^{t,(l)}]$ evolve according to
    \begin{align*}
    & \bo \zeta^{t+1, (l)}  =  \bo \zeta^{t, (l)} +  \eta \frac{\sigma_r^*}{n} \bo \Delta_{\mb R}^{t,(l)\top} \mb 1_n, \qquad \mb X^{t+1, (l)} \coloneqq \mb X^{t, (l)} + \eta \bo \Delta_{\mb R}^{t,(l)} \mb Y^{t, (l)}\\
    &\quad - 2 \eta c_{\perp} \frac{\pi \sigma_r^*}{n} \mb 1_n\mb 1_n\t \mb X^{t,(l)},\\
        &  \mb Y^{t+1, (l)} \coloneqq \mb Y^{t, (l)} + \eta \bo \Delta_{\mb R}^{t,(l)\top} \mb X^{t, (l)}, \\ 
    & \text{where~} \bo \Delta_{\mb R}^{t, (l)} \coloneqq  \pi \mc P^{(l)}\big(\mb R^* - \psi(\mb 1_n \bo\zeta^{t, (l)\top}+ \mb X^{t, (l)} {\mb Y^{t, (l)} }\t  )\big) + \mc P^{(l),\complement}_\Omega\big(\mb R - \psi(\mb 1_n \bo \zeta^{t, (l)\top}\\
    &\quad + \mb X^{t, (l)} {\mb Y^{t, (l)} }\t  )\big).
    \end{align*}
    Here, the projection operator $\mc P^{(l)}$ acts on an observation by preserving either its $l$-th row (for $l\in[n]$) or its $(-l)$-th column (for $l \in \{-p,\cdots, -1\}$), while $\mc P^{(l),\complement}_\Omega\coloneqq \mc P^{(l),\complement}\circ\mc P_\Omega$ preserves the observed entries in the complement of the $l$-th row or the $(-l)$-th column. These iterates can be viewed as the gradient descent updates for the leave-one-out objective function
    \begin{align}
         & L^{(l)}([\bo \zeta, \mb X, \mb Y]) 
         \\ &\coloneqq   
	     \begin{cases}
            \begin{aligned}
            &\sum_{(i,j)\in\Omega,\,i\neq l}
            \Big\{-R_{i,j}(\zeta_j+\mb X_i\t\mb Y_j)
            +\Psi(\zeta_j+\mb X_i\t\mb Y_j)\Big\}\\
            &\quad+\pi\sum_{j\in[p]}
            \Big\{-R^*_{l,j}(\zeta_j+\mb X_l\t\mb Y_j)
            +\Psi(\zeta_j+\mb X_l\t\mb Y_j)\Big\}
            \\ &\quad+c_{\perp}\frac{\pi\sigma_r^*}{n}
            \norm{\mb 1_n\t\mb X}_2^2,
            \end{aligned}
            & l\in[n],\\[2ex]
            \begin{aligned}
            &\sum_{(i,j)\in\Omega,\,-j\neq l}
            \Big\{-R_{i,j}(\zeta_j+\mb X_i\t\mb Y_j)
            +\Psi(\zeta_j+\mb X_i\t\mb Y_j)\Big\}\\
            &\quad+\pi\sum_{i\in[n]}
            \Big\{-R^*_{i,-l}(\zeta_{-l}+\mb X_i\t\mb Y_{-l})
            +\Psi(\zeta_{-l}+\mb X_i\t\mb Y_{-l})\Big\}
            \\ &\quad+c_{\perp}\frac{\pi\sigma_r^*}{n}
            \norm{\mb 1_n\t\mb X}_2^2,
            \end{aligned}
            & l\in[-p].
        \end{cases}
        \label{eq: leave-one-out objective function}
	    \end{align}

    We use the following probability bookkeeping throughout this stage. Work
    first on the single truncation event in
    Assumption~\ref{assumption: noise and sampling}, whose complement has
    probability $O(d^{-c-10})$, and fix any $q\in\{0,1,2,3\}$. Conditional on
    this event, the concentration bounds for a fixed iteration and a fixed
    leave-one-out index are invoked with failure probability
    $O(d^{-c-q-12})$. A union bound over $l\in[n]\cup[-p]$ therefore gives
    failure probability $O(d^{-c-q-11})$ for each one-step result. The
    local-geometry and sampling-operator events are uniform over the RIC and
    are selected at the same or a stronger probability level, so they require
    no further union over the iterates. Since $t_0\leq d^{10}$, a final union
    bound over the gradient-descent steps retains failure probability
    $O(d^{-c-q})$. Adding the probability of the shared truncation event does
    not change this order. Taking $q=0$ yields the stated probability in
    Theorem~\ref{thm: general GD}; the stronger levels will be used below.

    Our analysis is based on an inductive argument. We assume that the bounds in Condition~\ref{condition: inductive contraction for GD} hold at step $t$ and show that they continue to hold at step $t+1$ with high probability.

\begin{condition}[Inductive Bounds at Step $t$]
\label{condition: inductive contraction for GD}
For a given $t\in\{0,\ldots,t_0\}$, the following bounds hold for some
appropriately chosen constants $C_k$, $k\in[10]$:
\begin{subequations}
\begin{align}
    & \max_{l\in[n]\cup[-p]}\distf(\bo\theta^t,\bo\theta^{t,(l)})
    \leq C_1\nu^t \Big\{\kappa\kappa_\psi\sqrt{\mu\rho r}\,\delta_{\sf{row}}^0+\delta_{\sf{loo}}^0 \Big\}
    +C_2\frac{\sigma\kappa\kappa_\psi
    \sqrt{\kappa\mu\rho r^2\log d}}
    {\underline c_\psi\sqrt{\pi\sigma_r^*}},   
    \label{eq: condition -- F-norm consistency between thetat and thetatl}\\
    & \distf(\bo\theta^t,\bo\theta^*)
    \leq C_3\nu^t\sqrt d\,\delta_{\sf{row}}^0
    +C_4\frac{\sigma\sqrt{d\kappa r\log d}}
    {\underline c_\psi\sqrt{\pi\sigma_r^*}},   
    \label{eq: condition -- F-norm consistency between theta and thetat}\\
    & \dist\ti(\bo\theta^t,\bo\theta^*)
    \leq C_5\nu^t\Big\{
        \kappa\kappa_\psi\sqrt{\mu\rho r}\,\delta_{\sf{row}}^0+\delta_{\sf{loo}}^0
    \Big\}
    +C_6\frac{\sigma\kappa\kappa_\psi
    \sqrt{\kappa\mu\rho r^2\log d}}
    {\underline c_\psi\sqrt{\pi\sigma_r^*}},   
    \label{eq: condition -- two-to-infinity-norm consistency between theta and thetat}\\
    & \sqrt{\frac{n}{\sigma_r^*}}\max_{j\in[p]} \big|(\bo \zeta^{t,(-j)} - \bo \zeta^*)_j\big| \vee \max_{i\in[n]} \bignorm{(\mb X^{t,(i)} \tilde{\mb O}^{t,(i)} - \mb X^*)_i}_2\\
    &\quad \vee  \max_{j\in[p]} \bignorm{(\mb Y^{t,(-j)} \tilde{\mb O}^{t,(-j)} - \mb Y^*)_j}_2\\
    & \leq    C_7\nu^{t}\Big\{
\kappa\kappa_\psi\sqrt{\mu\rho r}\,\delta_{\sf{row}}^0+\delta_{\sf{loo}}^0
\Big\}
+C_8\frac{\sigma\kappa\kappa_\psi
\sqrt{\kappa\mu\rho r^2\log d}}
{\underline c_\psi\sqrt{\pi\sigma_r^*}},   
    \label{eq: condition -- rowwise consistency of LOO iterates}\\
        & \bignorm{\mb 1_n\t\mb X^t}_2
    \vee\max_{l\in[n]\cup[-p]}
    \bignorm{\mb 1_n\t\mb X^{t,(l)}}_2
    \leq C_9 \sqrt{n} \frac{\sigma \kappa^{\frac32} \kappa_\psi\sqrt{dr\log d} + \kappa \kappa_\psi \sqrt{d} \big(\xi^{\sf{gd}} + \xi^{\sf{gd},\sf{loo}}\big) }{c_{\perp} \sqrt{\pi \sigma_r^*}},   
    \label{eq: condition -- orthogonality between 1 and Xt} \\ 
    & \bignorm{{\mb X^t}\t\mb X^t-{\mb Y^t}\t\mb Y^t}\fb
    \vee\max_{l\in[n]\cup[-p]}
    \bignorm{{\mb X^{t,(l)}}\t\mb X^{t,(l)}
    -{\mb Y^{t,(l)}}\t\mb Y^{t,(l)}}\fb\\
    & \leq   C_{10} \frac{t}{t_0}\sqrt{\sigma_1^*} \frac{\sigma \kappa^{\frac32} \kappa_\psi \sqrt{dr \log d} + \kappa \kappa_\psi \sqrt{d} \big(\xi^{\sf{gd}} + \xi^{\sf{gd},\sf{loo}} \big) }{c_{\perp} \sqrt{\pi \sigma_r^*}},   
    \label{eq: condition -- balancedness between X and Y}
\end{align}
\end{subequations}
where $\nu\coloneqq1-\eta\underline c_{\sf{Hess}}\pi\sigma_r^*/4$.
\end{condition}

\begin{remark}
    The distance bounds consist of a geometrically decaying initialization
    term and a stochastic-error term, while the final two bounds keep the
    iterates approximately balanced and centered. In particular, the rowwise bound
    \eqref{eq: condition -- two-to-infinity-norm consistency between theta and thetat}
    is a coarse trajectory-level control used to close the induction and keep
    the iterates in the RIC. It is not intended as the final rowwise
    statistical rate. At $t=t_0$, the linear expansion in
    Theorem~\ref{thm: linear approximations} sharpens this bound to
    \eqref{eq: rowwise error from linear approximation}.
\end{remark}

We note that the polar-factor perturbation bound (Lemma~\ref{lemma: perturbation theory of optimal rotation under F norm}) and the RIC row-norm bounds imply
\begin{align}
\max_{l\in[n]\cup[-p]}\dist\ti(\bo\theta^{t,(l)},\bo\theta^*)
&\lesssim \dist\ti(\bo\theta^t,\bo\theta^*)
 +\Big\{1+\kappa\sqrt{\frac{\mu r}{n\wedge p}}\Big\}
 \max_{l\in[n]\cup[-p]}
 \distf(\bo\theta^t,\bo\theta^{t,(l)})
\nonumber\\
&\lesssim \dist\ti(\bo\theta^t,\bo\theta^*)
 +\max_{l\in[n]\cup[-p]}
 \distf(\bo\theta^t,\bo\theta^{t,(l)}). 
\label{eq: rowwise transfer from actual to LOO}
\end{align}
Thus, the same
uniform leave-one-out bound as in \eqref{eq: condition -- two-to-infinity-norm consistency between theta and thetat} follows from \eqref{eq: condition -- F-norm consistency between thetat and thetatl} and \eqref{eq: condition -- two-to-infinity-norm consistency between theta and thetat}. 

    In the remainder of the analysis, we aim to establish the following for the $(t+1)$-th step: 
    \begin{enumerate}
        \item Prove the \emph{inductive contraction conditions} for $t+1$ with high probability. 
        \item Prove that the gradient sequence of $L$ vanishes as the iteration proceeds with high probability. 
    \end{enumerate}

\bigskip

The following lemma controls the discrepancy between the actual and
leave-one-out trajectories for the $(t+1)$-th step (cf. \eqref{eq: condition -- F-norm consistency between thetat and thetatl} for the $(t+1)$-th step). 
\begin{lemma}[Leave-one-out coupling]
\label{lemma: leave-one-out error}
Suppose that the assumptions of Theorem~\ref{thm: general GD} hold and that
Condition~\ref{condition: inductive contraction for GD} is valid at step $t$.
Assume, in addition, that
$\pi(n\wedge p)\gg
\rho\kappa^2\mu^2r^2(\sigma^2/\underline c_{\sf{Hess}})^2\log d$.
Then, for sufficiently large constants $C_1$ and $C_2$, with probability at
least $1-O(d^{-c-1})$,
\begin{align}
&\max_{l\in[n]\cup[-p]}
\distf(\bo\theta^{t+1},\bo\theta^{t+1,(l)})
\leq C_1\nu^{t+1}\Big\{\kappa\kappa_\psi\sqrt{\mu\rho r}\,\delta_{\sf{row}}^0+\delta_{\sf{loo}}^0\Big\}
+\frac{C_2\sigma\kappa\kappa_\psi
\sqrt{\kappa\mu\rho r^2\log d}}
{\underline c_\psi\sqrt{\pi\sigma_r^*}}.
\label{eq: LOO coupling bound}
\end{align}
\end{lemma}

Lemma~\ref{lemma: leave-one-out error} handles the coupling component of the
induction. We next propagate the global and rowwise errors that correspond to the inductive hypotheses \eqref{eq: condition -- F-norm consistency between theta and thetat},~\eqref{eq: condition -- two-to-infinity-norm consistency between theta and thetat}, and \eqref{eq: condition -- rowwise consistency of LOO iterates} to the next step.

\begin{lemma}[Global and rowwise error propagation]
\label{lemma: gd error controls}
Suppose that the assumptions of Theorem~\ref{thm: general GD} hold and that
Condition~\ref{condition: inductive contraction for GD} is valid at step $t$.
Then, for sufficiently large constants $C_3$ and $C_4$, with probability at
least $1-O(d^{-c-1})$,
\begin{align}
\distf(\bo\theta^{t+1},\bo\theta^*)
&\leq C_3\nu^{t+1}\sqrt d\,\delta_{\sf{row}}^0
+C_4\frac{\sigma\sqrt{\kappa rd\log d}}
 {\underline c_\psi\sqrt{\pi\sigma_r^*}}.
\label{eq: propagated global GD error}
\end{align}
For sufficiently large constants $C_5,C_6,C_7,C_8$, the same event satisfies
\begin{align}
&\dist\ti(\bo\theta^{t+1},\bo\theta^*)
\leq C_5\nu^{t+1}\Big\{
\kappa\kappa_\psi\sqrt{\mu\rho r}\,\delta_{\sf{row}}^0+\delta_{\sf{loo}}^0
\Big\}
+C_6\frac{\sigma\kappa\kappa_\psi
\sqrt{\kappa\mu\rho r^2\log d}}
{\underline c_\psi\sqrt{\pi\sigma_r^*}},
\label{eq: propagated rowwise GD error}\\
    & \sqrt{\frac{n}{\sigma_r^*}}\max_{j\in[p]} \big| (\bo \zeta^{t+1,(-j)} - \bo \zeta^*)_j \big| \vee \max_{i\in[n]} \bignorm{(\mb X^{t+1,(i)} \mb O^{t+1,(i)} - \mb X^*)_i}_2\\
    &\quad \vee  \max_{j\in[p]} \bignorm{(\mb Y^{t+1,(-j)} \mb O^{t+1,(-j)} - \mb Y^*)_j}_2\\
& \quad \leq   C_7\nu^{t+1}\Big\{
\kappa\kappa_\psi\sqrt{\mu\rho r}\,\delta_{\sf{row}}^0+\delta_{\sf{loo}}^0
\Big\}
+C_8\frac{\sigma\kappa\kappa_\psi
\sqrt{\kappa\mu\rho r^2\log d}}
{\underline c_\psi\sqrt{\pi\sigma_r^*}}.
\label{eq: propagated LOO rowwise error}
\end{align}
\end{lemma}

    \begin{lemma}
	    \label{lemma: balancedness and orthogonality}
	        Suppose that the assumptions of Theorem~\ref{thm: general GD} hold. Let $\eta \leq \frac{1}{8 c_{\perp}\pi \sigma_1^*}$ and $t_0 = \lfloor \frac{c_t \log d}{\eta \underline c_\psi \pi \sigma_r^* }\rfloor\leq d^{10}$ for some $c_t > 0$. Fix $0\leq t<t_0$ and suppose that Condition~\ref{condition: inductive contraction for GD} holds for every $0\leq s\leq t$. Then the following holds with probability at least $1 - O(d^{-c-1})$:
	        \begin{align}
	        & \bignorm{\mb 1_n\t \mb X^{t+1}}_2
            \vee\max_{l\in[n]\cup[-p]}
            \bignorm{\mb 1_n\t\mb X^{t+1,(l)}}_2\\
    &\quad \leq C_9  \sqrt{n} \frac{\sigma  \kappa^{\frac32} \kappa_\psi \sqrt{ d r \log d} + \kappa \kappa_\psi \sqrt{d} ( \xi^{\sf{gd}} + \xi^{\sf{gd},\sf{loo}} )  }{c_{\perp} \sqrt{\pi \sigma_r^*}},\\
	            & \bignorm{{\mb X^{t+1}}\t \mb X^{t+1} - {\mb Y^{t+1}}\t \mb Y^{t+1}}\fb
            \vee\max_{l\in[n]\cup[-p]}
            \bignorm{{\mb X^{t+1,(l)}}\t\mb X^{t+1,(l)}
            -{\mb Y^{t+1,(l)}}\t\mb Y^{t+1,(l)}}\fb\\
    & \leq   C_{10}\frac{t+1}{t_0} \sqrt{\sigma_1^*} \frac{\sigma  \kappa^{\frac32} \kappa_\psi \sqrt{ d r \log d} + \kappa \kappa_\psi \sqrt{d} ( \xi^{\sf{gd}} + \xi^{\sf{gd},\sf{loo}} ) }{c_{\perp} \sqrt{\pi \sigma_r^*}}.
	        \end{align}
    \end{lemma}
    \begin{remark}
            Several prior statistical analyses introduce explicit Frobenius-norm or balancing regularization into the actual or auxiliary nonconvex objective \citep{chen2019noisy,fan2025covariates}. By contrast, we do not penalize the imbalance ${\mb X^t}\t\mb X^t-{\mb Y^t}\t\mb Y^t$, thereby avoiding the associated regularization bias. Related optimization work shows that gradient methods can maintain or induce approximate balancedness without such a penalty in asymmetric low-rank matrix factorization and matrix sensing \citep{du2018algorithmic,ma2021beyond,ye2021global,soltanolkotabi2023implicit}. Our analysis requires a more stringent, inference-oriented form of this phenomenon: the Gram-matrix discrepancy must remain sufficiently small throughout the trajectory that its contribution along the scaling directions is negligible in the subsequent rowwise linear expansion and hence does not affect the limiting distribution.
    \end{remark}

We argue by induction over $t$. The initialization bounds
\eqref{eq: GD initialization requirement}, the explicit smallness conditions
\eqref{eq: epsilon GD initialization requirement}, and
Condition~\ref{condition: inductive contraction for GD} at $t=0$ establish the
base case. Lemma~\ref{lemma: leave-one-out error}
propagates the coupling between the actual and leave-one-out iterates, while
Lemma~\ref{lemma: gd error controls} propagates the Frobenius and rowwise
estimation bounds. Lemma~\ref{lemma: balancedness and orthogonality}, in turn,
keeps the iterates approximately balanced and centered.

We make the RIC implication explicit. Since $\nu^t\leq1$,
\eqref{eq: epsilon GD initialization requirement} makes the
initialization-dependent parts of the proximity bounds in
Condition~\ref{condition: inductive contraction for GD} smaller than half of
the corresponding RIC radii with
$\epsilon=\epsilon_0/\kappa_\psi^2$. The
signal-strength condition gives the same conclusion for the stochastic-error
parts.
Consequently, uniformly over $0\leq t\leq t_0$,
\begin{align*}
    &\distf(\bo\theta^t,\bo\theta^*)
    \vee
    \max_{l\in[n]\cup[-p]}
    \distf(\bo\theta^t,\bo\theta^{t,(l)})
    \leq \frac{\sqrt{\epsilon\sigma_r^*}}{2\kappa},\\
    &\dist\ti(\bo\theta^t,\bo\theta^*)
    \vee
    \max_{l\in[n]\cup[-p]}
    \dist\ti(\bo\theta^{t,(l)},\bo\theta^*)
    \vee
    \max_{l\in[n]\cup[-p]}
    \distf(\bo\theta^t,\bo\theta^{t,(l)})
    \leq \frac12\sqrt{\frac{\epsilon\sigma_r^*}{\kappa d}}.
\end{align*}
In the balanced parametrization,
$\bignorm{\mb X^*}=\bignorm{\mb Y^*}=\sqrt{\sigma_1^*}$ and
$d=n\vee p$. Hence the second display is stronger than each of the three
rowwise inequalities in
\eqref{eq: rowwise consistency of parameters in geometry lemma}. The first
display also verifies the global proximity requirement in
\eqref{eq: direction set}. The additional coupling bound and the triangle
inequality give the same conclusion when a leave-one-out iterate is aligned
relative to the actual iterate. Thus every aligned actual and leave-one-out
iterate, as well as each line segment used in the one-step argument, lies in
the RIC.
Lemma~\ref{lemma: bounds for Hessian} therefore applies at the next iteration.
Conditional on the shared truncation event, the proofs of
Lemmas~\ref{lemma: leave-one-out error},
\ref{lemma: gd error controls}, and
\ref{lemma: balancedness and orthogonality} may be run uniformly over all
leave-one-out indices with failure probability $O(d^{-c-11})$ at each fixed
iteration. Since $t_0\leq d^{10}$, a union bound over $0\leq t<t_0$ closes the
induction with conditional failure probability $O(d^{-c-1})$. Adding the
$O(d^{-c-10})$ failure probability of the truncation event yields the stated
probability $1-O(d^{-c})$.

Recall that $\nu=1-\eta\underline c_{\sf{Hess}}\pi\sigma_r^*/4$. Since
$t_0\geq c_t\log d/(\eta\underline c_\psi\pi\sigma_r^*)-1$, the inequality
$1-x\leq\exp(-x)$ gives
$\nu^{t_0}\leq
\exp\{\eta\underline c_{\sf{Hess}}\pi\sigma_r^*/4\}
d^{-c_t\underline c_{\sf{Hess}}/(4\underline c_\psi)}
\lesssim d^{-c_t\underline c_{\sf{Hess}}/(4\underline c_\psi)}$.
Together with $\sigma\geq d^{-c-10}$ from
    Assumption~\ref{assumption: noise and sampling}, this ensures, for a
sufficiently large $c_t$, that the initialization-dependent terms in the
contraction bounds are absorbed by their stochastic-error terms.
Substituting $t=t_0$ into
Condition~\ref{condition: inductive contraction for GD} and using
\eqref{eq: rowwise transfer from actual to LOO} therefore gives
\eqref{eq: final GD LOO coupling}--\eqref{eq: final GD rowwise errors}. The
centering and balancedness conclusions
\eqref{eq: final GD centering}--\eqref{eq: final GD balancedness} follow by
setting $t=t_0$ in the last two inductive bounds.

\qed

\subsection{Proofs of the Auxiliary Lemmas for Theorem~\ref{thm: general GD}}
\label{subsec: proofs of auxiliary lemmas for general GD}

\subsubsection{Proof of Lemma~\ref{lemma: leave-one-out error}}
We first fix $l\in[n]$ and align the leave-one-out iterate with the actual
iterate by $\mb O^{t,(l)}$. Define
\begin{align*}
    & \bo\Delta^{t,(l)}
 \coloneqq[\bo\zeta^t-\bo\zeta^{t,(l)},
\mb X^t-\mb X^{t,(l)}\mb O^{t,(l)},
\mb Y^t-\mb Y^{t,(l)}\mb O^{t,(l)}],\\
    & \gamma_1
 \coloneqq\bo\Delta^{t,(l)}-\eta\Big\{
\nabla L^{\sf{aug}}([\bo\zeta^t,\mb X^t,\mb Y^t])
-\nabla L^{\sf{aug}}([\bo\zeta^{t,(l)},
\mb X^{t,(l)}\mb O^{t,(l)},\mb Y^{t,(l)}\mb O^{t,(l)}])\Big\},\\
    & \gamma_2
 \coloneqq\eta\Big\{
\nabla L^{\sf{diff}}(\mb X^t,\mb Y^t)
-\nabla L^{\sf{diff}}(\mb X^{t,(l)}\mb O^{t,(l)},
\mb Y^{t,(l)}\mb O^{t,(l)})\Big\},\\
    & \gamma_3
 \coloneqq\eta\Big\{
\nabla L^{(l)}([\bo\zeta^{t,(l)},
\mb X^{t,(l)}\mb O^{t,(l)},\mb Y^{t,(l)}\mb O^{t,(l)}])
-\nabla L([\bo\zeta^{t,(l)},
\mb X^{t,(l)}\mb O^{t,(l)},\\
    &\quad \mb Y^{t,(l)}\mb O^{t,(l)}])\Big\}.
\end{align*}
Rotational equivariance and the optimality of the Procrustes rotation at step
$t+1$ then give
\begin{equation}
\distf(\bo\theta^{t+1},\bo\theta^{t+1,(l)})
\leq\norm{\gamma_1}_2+\norm{\gamma_2}_2+\norm{\gamma_3}_2.
\label{eq: leave-one-out gradient descent error analysis}
\end{equation}

We examine the three terms in turn.
\begin{itemize}
    \item The inductive hypotheses ensure that the segment joining the two
    aligned iterates remains in the RIC. Let $\overline{\mb H}^{t,(l)}$ be the
    average Hessian of $L^{\sf{aug}}$ along this segment.
    The bound \eqref{eq: bound for alpha2} holds for arbitrary directions.
Since $\alpha_1\geq0$, every Hessian $\mb H$ along the segment satisfies
$\mb H\succeq-\underline c_{\sf{Hess}}\pi\sigma_r^*\mb I$.
Together with \eqref{eq: local Hessian upper bound}, this places the
spectrum of $\overline{\mb H}^{t,(l)}$ in
$[-\underline c_{\sf{Hess}}\pi\sigma_r^*,
\bar c_{\sf{Hess}}c_\perp\pi\sigma_1^*]$.
Spectral calculus and \eqref{eq: local Hessian lower bound} therefore give
    \begin{align}
        \bignorm{\overline{\mb H}^{t,(l)}\bo\Delta^{t,(l)}}_2^2
    &\leq
    \big(\bar c_{\sf{Hess}}c_\perp\pi\sigma_1^*
    -\underline c_{\sf{Hess}}\pi\sigma_r^*\big)
    {\bo\Delta^{t,(l)}}\t
    \overline{\mb H}^{t,(l)}\bo\Delta^{t,(l)}
    \nonumber\\
    &\quad+
    \underline c_{\sf{Hess}}\bar c_{\sf{Hess}}c_\perp
    \pi^2\sigma_r^*\sigma_1^*
    \bignorm{\bo\Delta^{t,(l)}}_2^2
    \nonumber\\
    &\leq
    2\bar c_{\sf{Hess}}c_\perp\pi\sigma_1^*
    {\bo\Delta^{t,(l)}}\t
    \overline{\mb H}^{t,(l)}\bo\Delta^{t,(l)}.
        \label{eq: averaged Hessian multiplication bound LOO}
    \end{align}
    Since $\bar c_{\sf{Hess}}=4$, the stated step-size condition gives
    $\eta\leq1/(2\bar c_{\sf{Hess}}c_\perp\pi\sigma_1^*)$.
    Thus, \eqref{eq: local Hessian lower bound} and
    \eqref{eq: averaged Hessian multiplication bound LOO} yield
    \begin{align*}
        \norm{\gamma_1}_2^2
        &=\bignorm{\big(\mb I-\eta\overline{\mb H}^{t,(l)}\big)
        \bo\Delta^{t,(l)}}_2^2\\
        &=\bignorm{\bo\Delta^{t,(l)}}_2^2
        -2\eta{\bo\Delta^{t,(l)}}\t\overline{\mb H}^{t,(l)}\bo\Delta^{t,(l)}
        +\eta^2\bignorm{\overline{\mb H}^{t,(l)}\bo\Delta^{t,(l)}}_2^2\\
        &\leq\bignorm{\bo\Delta^{t,(l)}}_2^2
        -\eta{\bo\Delta^{t,(l)}}\t\overline{\mb H}^{t,(l)}\bo\Delta^{t,(l)}\\
        &\leq\big(1-\eta\underline c_{\sf{Hess}}\pi\sigma_r^*\big)
        \bignorm{\bo\Delta^{t,(l)}}_2^2.
    \end{align*}
    Taking square roots and using $\sqrt{1-x}\leq1-x/2$ for $x\in[0,1]$ gives
    \begin{align}
    \bignorm{\gamma_1}_2
    &\leq\Big(1-\frac{\eta\underline c_{\sf{Hess}}\pi\sigma_r^*}{2}\Big)
    \distf(\bo\theta^t,\bo\theta^{t,(l)}).
    \label{eq: contraction of gamma1}
    \end{align}
    Recall that $\nu=1-\eta\underline c_{\sf{Hess}}\pi\sigma_r^*/4$; the
    slower bookkeeping rate leaves enough slack to absorb the transient
    sampling fluctuation below.

    \item Regarding $\gamma_2$, the derivative formula
    \eqref{eq: derivatives of Ldiff} and the RIC factor-norm bounds give
    \begin{align}
    \norm{\gamma_2}_2
    \lesssim{}&\eta c_{\sf{aug}}\pi\sqrt{\sigma_1^*}
    \Big\{
    \bignorm{{\mb X^t}\t\mb X^t-{\mb Y^t}\t\mb Y^t}\fb
    \nonumber\\
    &\qquad\qquad
    +\bignorm{{\mb X^{t,(l)}}\t\mb X^{t,(l)}
    -{\mb Y^{t,(l)}}\t\mb Y^{t,(l)}}\fb
    \Big\}.
    \label{eq: balancing-gradient contribution to LOO coupling}
    \end{align}
    The balancedness hypothesis in
    Condition~\ref{condition: inductive contraction for GD} therefore yields
    \begin{align}
    \norm{\gamma_2}_2
    &\lesssim \eta\sigma
    \sqrt{\pi\mu\rho r\sigma_1^*\log d}.
    \label{eq: completed gamma2 bound}
    \end{align}
    \item Regarding the term $\gamma_3$, direct calculation gives, for $l \in[n]$,
    \begin{align*}
        & \gamma_3 = \eta\cdot \mathrm{vec}\Big(
        \Big\{\mc P_{l,\cdot}\Big(
        \big[\psi(\mc M([\bo \zeta^{t,(l)}, \mb X^{t,(l)} \mb O^{t,(l)}, \mb Y^{t,(l)} \mb O^{t,(l)}]))-\psi(\mb M^*)\big]
        \circ (\pi\mb 1_n\mb 1_p\t\\
    &\quad -\bo\Omega)\Big)\Big\}\t
        \mb 1_n\sqrt{\frac{\sigma_r^*}{n}}, \\
        & \quad \mc P_{l,\cdot}\Big(
        \big[\psi(\mc M([\bo \zeta^{t,(l)}, \mb X^{t,(l)} \mb O^{t,(l)}, \mb Y^{t,(l)} \mb O^{t,(l)}]))-\psi(\mb M^*)\big]
        \circ (\pi\mb 1_n\mb 1_p\t\\
    &\quad -\bo\Omega)\Big)
        \mb Y^{t,(l)}\mb O^{t,(l)}, \\
        & \quad \Big\{\mc P_{l,\cdot}\Big(
        \big[\psi(\mc M([\bo \zeta^{t,(l)}, \mb X^{t,(l)} \mb O^{t,(l)}, \mb Y^{t,(l)} \mb O^{t,(l)}]))-\psi(\mb M^*)\big]
        \circ (\pi\mb 1_n\mb 1_p\t\\
    &\quad -\bo\Omega)\Big)\Big\}\t
        \mb X^{t,(l)}\mb O^{t,(l)}\Big) \\
        & \quad + \eta \cdot \mathrm{vec}\Big(
        \mc P_{l,\cdot}(\mb E \circ \bo \Omega)\t \mb 1_n\sqrt{\frac{\sigma_r^*}{n}},
        \mc P_{l,\cdot}(\mb E \circ \bo \Omega) \mb Y^{t,(l)} \mb O^{t,(l)},
        \mc P_{l,\cdot}(\mb E \circ \bo \Omega)\t \mb X^{t,(l)} \mb O^{t,(l)} \Big).
    \end{align*}
    With respect to the terms containing $\mb E \circ \bo \Omega$, we invoke the matrix Bernstein inequality to derive that 
    \begin{align}
        & \norm{\mc P_{l,\cdot}(\mb E \circ \bo \Omega)\t \mb 1_n}\fb \leq \norm{\mc P_{l,\cdot}(\mb E) \circ \bo \Omega }\fb  \lesssim \sigma \sqrt{\pi p}   , \\
        & \norm{\mc P_{l,\cdot}(\mb E \circ \bo \Omega)\t {\mb X}^{t,(l)}}\fb \leq \norm{\mc P_{l,\cdot}(\mb E \circ \bo \Omega)}\fb \bignorm{{\mb X}^{t,(l)}}\ti  \lesssim \sigma \sqrt{\pi r p \sigma_1^* }  \cdot \sqrt{\frac{\mu r }{n}} , \\ 
        & \norm{\mc P_{l,\cdot}(\mb E \circ \bo \Omega) \mb Y^{t,(l)}}\fb \lesssim \sigma \sqrt{\sigma_1^*\pi  r \log d}, 
        \label{eq: leave-one-out gradient descent analysis -- single noise vector}
    \end{align}
    with probability at least $1- O(d^{-c-2})$ for the fixed index $l$. Moreover, conditionally on the leave-one-out iterate, the matrix Bernstein inequality applied to the centered missingness fluctuation gives
    \begin{align*}
        & \norm{\mc P_{l,\cdot}\Big(
        \big[\psi(\mc M([\bo \zeta^{t,(l)}, \mb X^{t,(l)} \mb O^{t,(l)}, \mb Y^{t,(l)} \mb O^{t,(l)}]))-\psi(\mb M^*)\big]
        \circ(\bo\Omega-\pi\mb 1_n\mb 1_p\t)\Big)\mb Y^{t,(l)}}_2 \\
        \lesssim{} & \sqrt{\pi}\,
        \norm{\big[\psi(\mc M([\bo \zeta^{t,(l)}, \mb X^{t,(l)} \mb O^{t,(l)}, \mb Y^{t,(l)} \mb O^{t,(l)}]))-\psi(\mb M^*)\big]_{l,\cdot}}_\infty
        \bignorm{\mb Y^{t,(l)}}\fb\sqrt{\log d} \\
        & \quad +
        \norm{\big[\psi(\mc M([\bo \zeta^{t,(l)}, \mb X^{t,(l)} \mb O^{t,(l)}, \mb Y^{t,(l)} \mb O^{t,(l)}]))-\psi(\mb M^*)\big]_{l,\cdot}}_\infty
        \bignorm{\mb Y^{t,(l)}}\ti\log d \\
        \lesssim{} & \sigma^2
        \norm{\big[\mc M([\bo \zeta^{t,(l)}, \mb X^{t,(l)} \mb O^{t,(l)}, \mb Y^{t,(l)} \mb O^{t,(l)}])-\mb M^*\big]_{l,\cdot}}_\infty
        \Big(\sqrt{\pi\sigma_1^*r\log d}
        +\sqrt{\frac{\mu r\sigma_1^*}{p}}\log d\Big)
    \end{align*}
    with probability at least $1-O(d^{-c-2})$ for the fixed index $l$. The corresponding intercept
    and right-factor blocks are bounded, up to the same predictor-error
    factor, by
    $\sigma^2\sqrt{\pi p\sigma_r^*/n}$ and
    $\sigma^2\bignorm{\mb X^{t,(l)}}\ti\sqrt{\pi p}$,
    respectively. Both are dominated by the second term in
    \eqref{eq: gamma3 preliminary row bound}.

    Putting these inequalities together and using
    $\pi(n\wedge p)\gg\log d$ yields
    \begin{align}
    \norm{\gamma_3}_2
    \lesssim{}&\eta\sigma
    \sqrt{\pi\mu \rho r^2 \sigma_1^*\log d}
    \nonumber\\
    &+\eta\sigma^2
    \norm{\big[\mc M([\bo\zeta^{t,(l)},
    \mb X^{t,(l)}\mb O^{t,(l)},
    \mb Y^{t,(l)}\mb O^{t,(l)}])-\mb M^*\big]_{l,\cdot}}_\infty
    \sqrt{\pi\mu r\rho\sigma_1^*\log d}.
    \label{eq: gamma3 preliminary row bound}
    \end{align}
    Since $\mc M$ is invariant under a common factor rotation, we may use the
    population alignment $\widetilde{\mb O}^{t,(l)}$ when bounding the predictor
    error in the second line. In particular,
    \begin{align}
    &\norm{\big[\mc M([\bo\zeta^{t,(l)},
    \mb X^{t,(l)}\widetilde{\mb O}^{t,(l)},
    \mb Y^{t,(l)}\widetilde{\mb O}^{t,(l)}])-\mb M^*\big]_{l,\cdot}}_\infty
    \nonumber\\
    &\quad\leq\norm{\bo\zeta^{t,(l)}-\bo\zeta^*}_\infty
    +\bignorm{\mb X^{t,(l)}\widetilde{\mb O}^{t,(l)}-\mb X^*}\ti
    \bignorm{\mb Y^*}\ti\\
    &\quad +\bignorm{\mb X^{t,(l)}\widetilde{\mb O}^{t,(l)}}\ti
    \bignorm{\mb Y^{t,(l)}\widetilde{\mb O}^{t,(l)}-\mb Y^*}\ti
    \nonumber\\
    &\quad\lesssim
    \sqrt{\frac{\mu r\sigma_1^*}{n\wedge p}}
    \dist\ti(\bo\theta^{t,(l)},\bo\theta^*)\\
    & \quad  \lesssim
    \nu^t\sqrt{\frac{\mu r\sigma_1^*}{n\wedge p}}
    \Big\{\kappa\kappa_\psi\sqrt{\mu\rho r}\,\delta_{\sf{row}}^0+\delta_{\sf{loo}}^0\Big\}
    +\sqrt{\frac{\mu r\sigma_1^*}{n\wedge p}}
    \frac{\sigma\kappa\kappa_\psi
    \sqrt{\kappa\mu\rho r^2\log d}}
    {\underline c_\psi\sqrt{\pi\sigma_r^*}}.
    \label{eq: row predictor error for gamma3}
    \end{align}
    Here, the last inequality follows directly from \eqref{eq: rowwise transfer from actual to LOO}.
    Combining the last two displays and using the stated
    sampling condition gives
    \begin{align}
    \norm{\gamma_3}_2
    \lesssim{}&\eta\sigma
    \sqrt{\pi\mu \rho r^2\sigma_1^*\log d}
    +\eta\sigma^2\mu r\sigma_1^*
    \sqrt{\frac{\pi\rho\log d}{n\wedge p}}\,\nu^t
    \Big(\kappa\kappa_\psi\sqrt{\mu\rho r}\,\delta_{\sf{row}}^0+\delta_{\sf{loo}}^0\Big)
    \nonumber\\
    &+\eta\sigma^2\mu r\sigma_1^*
    \sqrt{\frac{\pi\rho\log d}{n\wedge p}}
    \frac{\sigma\kappa\kappa_\psi
    \sqrt{\kappa\mu\rho r^2\log d}}
    {\underline c_\psi\sqrt{\pi\sigma_r^*}} \\ 
    \lesssim & \eta\sigma
    \sqrt{\pi\mu \rho r^2 \sigma_1^*\log d} + \frac{\mu r \sqrt{\rho \log d}}{\sqrt{\pi (n\wedge p)}}\eta \sigma^2 \pi \sigma_1^* \nu^t \Big\{\kappa\kappa_\psi\sqrt{\mu\rho r}\,\delta_{\sf{row}}^0+\delta_{\sf{loo}}^0 \Big\}.
    \label{eq: completed gamma3 row bound}
    \end{align}
    Indeed, the coefficient in the second term relative to the local curvature
    is
    $\frac{\sigma^2}{\underline c_{\sf{Hess}}}\kappa\mu r
    \allowbreak \sqrt{\frac{\rho\log d}{\pi(n\wedge p)}} \allowbreak \ll 1$.
    Hence, after fixing the constants in the inductive bounds, the transient
    term is absorbed by the slack in \eqref{eq: contraction of gamma1}, while
    the last term is absorbed into the stochastic error bound.
\end{itemize}

For completeness, consider next $l=\allowbreak -j\in[-p]$. The removed column is
independent of
$[\bo\zeta^{t,(-j)},\allowbreak \mb X^{t,(-j)},\allowbreak \mb Y^{t,(-j)}]$.
Conditional Bernstein bounds for the intercept, left-factor, and right-factor
blocks are, respectively, of orders
$\sigma\allowbreak \sqrt{\pi\sigma_r^*\log d}$,
$\sigma\allowbreak \sqrt{\pi n\mu r\sigma_1^*/p}$, and
$\sigma\allowbreak \sqrt{\pi\sigma_1^*r\log d}$.
The centered-missingness terms obey the same bounds after replacing $\sigma$
by
$\sigma^2
\allowbreak \norm{[\mc M(\bo\theta^{t,(-j)})-\mb M^*]_{\cdot,j}}_\infty$.
For each fixed $j$, these conditional concentration bounds hold with probability
at least $1-O(d^{-c-2})$.
Consequently,
\begin{align*}
\norm{\gamma_3}_2
\lesssim{}&\eta\sigma
\sqrt{\pi\mu r\rho\sigma_1^*\log d}
+\eta\sigma^2
\norm{[\mc M(\bo\theta^{t,(-j)})-\mb M^*]_{\cdot,j}}_\infty
\sqrt{\pi\mu r\rho\sigma_1^*\log d}.
\end{align*}
The columnwise predictor error is controlled by
\begin{align*}
\norm{[\mc M(\bo\theta^{t,(-j)})-\mb M^*]_{\cdot,j}}_\infty
&\lesssim
\sqrt{\frac{\mu r\sigma_1^*}{n\wedge p}}
\dist\ti(\bo\theta^{t,(-j)},\bo\theta^*),
\end{align*}
so \eqref{eq: completed gamma3 row bound} holds for $l=\allowbreak -j$ as well.

Substituting the preceding bounds into
\eqref{eq: leave-one-out gradient descent error analysis} gives
\begin{align}
& \distf(\bo\theta^{t+1},\bo\theta^{t+1,(l)})
\leq{}\Big(1-\frac{\eta\underline c_{\sf{Hess}}\pi\sigma_r^*}{2}\Big)
\distf(\bo\theta^t,\bo\theta^{t,(l)})
\\ 
& \qquad  + C \Big[\frac{\mu r \sqrt{\rho \log d}}{\sqrt{\pi (n\wedge p)}}\eta \sigma^2 \pi \sigma_1^* \nu^t \Big\{\kappa\kappa_\psi\sqrt{\mu\rho r}\,\delta_{\sf{row}}^0+\delta_{\sf{loo}}^0 \Big\}+  \eta\sigma
    \sqrt{\pi\mu \rho r^2 \sigma_1^*\log d} \Big] \\ 
    & \lesssim  C_1 \nu^{t+1}  \Big\{\kappa\kappa_\psi\sqrt{\mu\rho r}\,\delta_{\sf{row}}^0+\delta_{\sf{loo}}^0 \Big\}+  C_2 \frac{\sigma \kappa \kappa_\psi
    \sqrt{\kappa \mu \rho r^2 \log d}}{\underline c_\psi \sqrt{\pi \sigma_r^*}} ,
\label{eq: leave-one-out one-step recursion}
\end{align}
where the transient term in \eqref{eq: completed gamma3 row bound} is absorbed
using the slack between the contraction factor in
\eqref{eq: contraction of gamma1} and $\nu$ by choosing appropriate constants. 
A union bound over $l\in[n]\cup[-p]$ gives failure probability
$O(d^{-c-1})$ and completes the proof.
\qed

\subsubsection{Proof of Lemma~\ref{lemma: gd error controls}}
\emph{Frobenius-norm contraction.}
Let $\mb O^t=\allowbreak \mb O(\bo\theta^t,\allowbreak \bo\theta^*)$ and set
$\bo\Delta^t=\allowbreak [\bo\zeta^t-\bo\zeta^*,\allowbreak 
\mb X^t\mb O^t-\mb X^*,\allowbreak \mb Y^t\mb O^t-\mb Y^*]$.
Rotational equivariance and the optimality of the Procrustes rotation at step
$t+1$ give
\begin{align*}
\distf(\bo\theta^{t+1},\bo\theta^*)
\leq{}&\bignorm{\bo\Delta^t
-\eta\Big\{\nabla L^{\sf{aug}}(
[\bo\zeta^t,\mb X^t\mb O^t,\mb Y^t\mb O^t])
-\nabla L^{\sf{aug}}(\bo\theta^*)\Big\}}_2\\
&+\eta\bignorm{\nabla L^{\sf{aug}}(\bo\theta^*)}_2
+\eta\bignorm{\nabla L^{\sf{diff}}(\mb X^t,\mb Y^t)}_2.
\end{align*}
{\emergencystretch=1.5em
Let $\overline{\mb H}^t$ denote the average Hessian of $L^{\sf{aug}}$ along
the segment joining $\bo\theta^*$ and
$[\bo\zeta^t,\allowbreak \mb X^t\mb O^t,\allowbreak \mb Y^t\mb O^t]$. The inductive hypotheses place
this segment in the RIC and make $\bo\Delta^t$ an admissible direction in
Lemma~\ref{lemma: bounds for Hessian}. The argument establishing
\eqref{eq: averaged Hessian multiplication bound LOO} gives
$\bignorm{\overline{\mb H}^t\bo\Delta^t}_2^2
\leq 2\bar c_{\sf{Hess}}c_\perp\pi\sigma_1^*
{\bo\Delta^t}\t\overline{\mb H}^t\bo\Delta^t$.
Consequently, for
$\eta\leq1/(2\bar c_{\sf{Hess}}c_\perp\pi\sigma_1^*)$,
expanding the squared norm and applying
\eqref{eq: local Hessian lower bound} yield\par}

\begin{align}
\bignorm{(\mb I-\eta\overline{\mb H}^t)\bo\Delta^t}_2^2
&\leq\big(1-\eta\underline c_{\sf{Hess}}\pi\sigma_r^*\big)
\norm{\bo\Delta^t}_2^2.
\label{eq: one-step analysis under F norm}
\end{align}
It follows that the first term is at most
$\nu\distf(\bo\theta^t,\bo\theta^*)$.

At the population parameter, both the balancing and centering gradients
vanish. Matrix Bernstein's inequality, applied blockwise to the score, gives
\begin{align}
\bignorm{\nabla L^{\sf{aug}}(\bo\theta^*)}_2
=\bignorm{\nabla L(\bo\theta^*)}_2
\lesssim\sigma\sqrt{\pi rd\sigma_1^*\log d}
\label{eq: score norm at population for GD}
\end{align}
uniformly with probability at least $1-O(d^{-c-1})$. Moreover, the gradient form~\eqref{eq: derivatives of Ldiff}, the RIC factor-norm bounds, and the balancedness hypothesis~\eqref{eq: condition -- balancedness between X and Y} imply
$\bignorm{\nabla L^{\sf{diff}}(\mb X^t,\mb Y^t)}_2
\lesssim\sigma\sqrt{\pi\mu\rho r\sigma_1^*\log d}$. Since
$\mu r\leq n\wedge p$, this term is dominated by
\eqref{eq: score norm at population for GD}. We therefore obtain
\begin{align}
\distf(\bo\theta^{t+1},\bo\theta^*)
&\leq\nu\distf(\bo\theta^t,\bo\theta^*)
+C\eta\sigma\sqrt{\pi rd\sigma_1^*\log d}.
\label{eq: global GD one-step recursion}
\end{align}
Iterating this recursion and using
$1-\nu=\eta\underline c_{\sf{Hess}}\pi\sigma_r^*/4$, together with
$c_{\sf{aug}}\asymp\underline c_\psi$ and hence
$\underline c_{\sf{Hess}}\asymp\underline c_\psi$, proves
\eqref{eq: propagated global GD error}.

\medskip
\emph{Rowwise contraction.}
We give the argument for a left-factor row. Fix $i\in[n]$, and recall that
$\widetilde{\mb O}^{t,(i)}$ aligns $\mb F^{t,(i)}$ with $\mb F^*$. 
The triangle inequality gives
\begin{align}
& \bignorm{(\mb X^{t+1}\mb O^{t+1}-\mb X^*)_{i,\cdot}}_2
\leq{}
\bignorm{(\mb X^{t+1,(i)}\widetilde{\mb O}^{t+1,(i)}
-\mb X^*)_{i,\cdot}}_2\\
    &\quad +\bignorm{\big(\mb X^{t+1}\mb O^{t+1}
-\mb X^{t+1,(i)}\widetilde{\mb O}^{t+1,(i)}\big)_{i,\cdot}}_2 \\ 
    & \leq   \bignorm{(\mb X^{t+1,(i)}\widetilde{\mb O}^{t+1,(i)}
-\mb X^*)_{i,\cdot}}_2 + \bignorm{\mb X^{t+1,(i)}}\ti
\bignorm{\widetilde{\mb O}^{t+1,(i)} - \mb O^{t+1,(i)}\mb O^{t+1}} \\
& + \norm{ \mb X^{t+1,(i)} \mb O^{t+1,(i)} - \mb X^{t+1} }\fb.
\label{eq: leave-one-out decomposition}
\end{align}
Define $
\mb D_i^t
\coloneqq\sf{diag}\Big(
\int_0^1\psi'\Big(
(1-s)M_{i,j}^*+s\mc M(\bo\theta^{t,(i)})_{i,j}
\Big)\,\mathrm ds
\Big)_{j\in[p]}$.
The mean-value theorem and the population update in the omitted row show that
the first term on the right-hand side of
\eqref{eq: leave-one-out decomposition} is bounded as follows: 
\begin{align}
    & \bignorm{(\mb X^{t+1,(i)}\widetilde{\mb O}^{t+1,(i)}
-\mb X^*)_{i,\cdot}}_2 \\ 
    & \leq  \bignorm{
(\mb X^{t,(i)}\widetilde{\mb O}^{t,(i)}-\mb X^*)_{i,\cdot}
\big(\mb I_r-\eta\pi{\mb Y^*}\t\mb D_i^t\mb Y^*\big)}_2\\
    &\quad +C\eta\pi\sigma^2\sqrt{\sigma_1^*}
\bignorm{\mb Y^{t,(i)}\widetilde{\mb O}^{t,(i)}-\mb Y^*}
\bignorm{(\mb X^{t,(i)}\widetilde{\mb O}^{t,(i)}-\mb X^*)_{i,\cdot}}_2
\nonumber\\
&\quad+C\eta\pi\sigma^2\Big\{
\bignorm{\bo\zeta^{t,(i)}-\bo\zeta^*}_2
\bignorm{\mb Y^{t,(i)}\widetilde{\mb O}^{t,(i)}-\mb Y^*}
+\norm{\mb X^*_{i,\cdot}}_2
\bignorm{\mb Y^{t,(i)}\widetilde{\mb O}^{t,(i)}-\mb Y^*}^2
\Big\}
\nonumber\\
&\quad+C\eta\pi\Bignorm{\big\{
\bo\zeta^{t,(i)}-\bo\zeta^*
+\big(\mb Y^{t,(i)}\widetilde{\mb O}^{t,(i)}-\mb Y^*\big)
\mb X_i^*\big\}\t\mb D_i^t\mb Y^*}_2
\nonumber\\
&\quad+\frac{2\eta c_\perp\pi\sigma_r^*}{n}
\bignorm{\mb 1_n\t\mb X^{t,(i)}}_2
+\bignorm{\mb X^{t+1,(i)}}\ti
\bignorm{\widetilde{\mb O}^{t,(i)\top}
\widetilde{\mb O}^{t+1,(i)}-\mb I_r}.
\label{eq: corrected held-out row decomposition}
\end{align}

The same squared-norm argument used in
\eqref{eq: one-step analysis under F norm} gives
\begin{align}
\bignorm{(\mb X^{t,(i)}\widetilde{\mb O}^{t,(i)}-\mb X^*)_{i,\cdot}
\big(\mb I_r-\eta\pi{\mb Y^*}\t\mb D_i^t\mb Y^*\big)}_2
&\leq\Big(1-\frac{\eta\underline c_\psi\pi\sigma_r^*}{2}\Big)
\bignorm{(\mb X^{t,(i)}\widetilde{\mb O}^{t,(i)}-\mb X^*)_{i,\cdot}}_2.
\label{eq: held-out row contraction}
\end{align}

By the triangle inequality,
$\distf(\bo\theta^{t,(i)},\bo\theta^*)$ is at most
$\distf(\bo\theta^t,\bo\theta^{t,(i)})+
\distf(\bo\theta^t,\bo\theta^*)$, where the latter terms are controlled by \eqref{eq: condition -- F-norm consistency between thetat and thetatl} and \eqref{eq: condition -- F-norm consistency between theta and thetat}. This leads to the following bound due to the condition $\kappa^2 \kappa_\psi^2 \mu r \ll n\wedge p$: 
\begin{align}
    & \distf(\bo \theta^{t,(i)}, \bo \theta^*) \lesssim \nu^t \big[ \sqrt d\,\delta_{\sf{row}}^0+\delta_{\sf{loo}}^0 \big]
+ \frac{\sigma
\sqrt{d \kappa r \log d}}
{\underline c_\psi\sqrt{\pi\sigma_r^*}}. 
\label{eq: LOO global error used in row propagation}
\end{align}

Consequently, the second and the third terms of \eqref{eq: corrected held-out row decomposition} are
absorbed into the contraction by the inductive hypotheses: 
\begin{align}
    & \eta\pi\sigma^2\sqrt{\sigma_1^*}
\bignorm{\mb Y^{t,(i)}\widetilde{\mb O}^{t,(i)}-\mb Y^*}
\bignorm{(\mb X^{t,(i)}\widetilde{\mb O}^{t,(i)}-\mb X^*)_{i,\cdot}}_2\\
    &\quad \ll  \eta \pi \sigma^2 \sigma_1^* \bignorm{(\mb X^{t,(i)}\widetilde{\mb O}^{t,(i)}-\mb X^*)_{i,\cdot}}_2,\\
& \eta\pi\sigma^2\Big\{
\bignorm{\bo\zeta^{t,(i)}-\bo\zeta^*}_2
\bignorm{\mb Y^{t,(i)}\widetilde{\mb O}^{t,(i)}-\mb Y^*}
+\norm{\mb X^*_{i,\cdot}}_2
\bignorm{\mb Y^{t,(i)}\widetilde{\mb O}^{t,(i)}-\mb Y^*}^2
\Big\} \\ 
    & \ll   \eta \pi \sigma^2 \sigma_1^* \Big[\nu^{t}
\Big[\sqrt{\mu\rho r}\,\delta_{\sf{row}}^0
+\sqrt{\frac{\mu\rho r}{d}}\,\delta_{\sf{loo}}^0\Big]
+ \frac{\sigma 
\sqrt{\kappa\mu\rho r^2\log d}}
{\underline c_\psi\sqrt{\pi\sigma_r^*}} \Big]
\label{eq: second term in the LOO rowwise}
\end{align}
Moreover, the remaining term satisfies
\begin{align}
&
\eta \pi \Bignorm{\Big\{
\bo\zeta^{t,(i)}-\bo\zeta^*
+\big(\mb Y^{t,(i)}\widetilde{\mb O}^{t,(i)}-\mb Y^*\big)
\mb X_i^*\Big\}\t\mb D_i^t\mb Y^*}_2
\nonumber\\
\leq & \eta \pi \sigma^2\sqrt{\sigma_1^*}
\Big\{\bignorm{\bo\zeta^{t,(i)}-\bo\zeta^*}_2
+\Bignorm{\big(\mb Y^{t,(i)}\widetilde{\mb O}^{t,(i)}-\mb Y^*\big)
\mb X_i^*}_2\Big\}
\nonumber\\
\leq & \eta \pi
\sigma^2 \Big\{ \sigma_1^* / {\sqrt{\kappa n}}
+ \sqrt{\frac{\mu r}{n}} \sigma_1^* \Big\}
\distf(\bo\theta^{t,(i)},\bo\theta^*)
\nonumber\\
\lesssim & \eta \pi
\sigma^2  \sigma_1^* \Big[ \nu^t \Big\{\sqrt{\mu\rho r}\,\delta_{\sf{row}}^0
+\sqrt{\frac{\mu\rho r}{d}}\,\delta_{\sf{loo}}^0\Big\}
+ \frac{\sigma
\sqrt{\kappa\mu\rho r^2\log d}}
{\underline c_\psi\sqrt{\pi\sigma_r^*}}\Big].
\label{eq: weighted cross-factor bound from Frobenius error}
\end{align}

It remains to control the change in the population alignment. For this purpose,
let
\begin{align*}
\mb G^{t,(i)}
\coloneqq{}&\mc P_\Omega^{(i),\complement}\big[
\psi(\mc M(\bo\theta^{t,(i)}))-\mb R
\big]
+\pi\mc P^{(i)}\Big(
\psi(\mc M(\bo\theta^{t,(i)}))-\mb R^*
\Big),
\end{align*}
and define
\begin{align*}
\mb F^{t+1,(i),\sf{srg}}
\coloneqq
\begin{bmatrix}
\mb X^{t,(i)}\widetilde{\mb O}^{t,(i)}
-\eta\mb G^{t,(i)}\mb Y^*
-2\eta c_\perp\frac{\pi\sigma_r^*}{n}
\mb 1_n\mb 1_n\t\mb X^{t,(i)}\widetilde{\mb O}^{t,(i)}\\
\mb Y^{t,(i)}\widetilde{\mb O}^{t,(i)}
-\eta{\mb G^{t,(i)}}\t\mb X^*
\end{bmatrix}.
\end{align*}
This surrogate retains the response noise, the missingness fluctuation, and
the centering update. Since
${\mb F^*}\t\mb F^{t,(i)}\widetilde{\mb O}^{t,(i)}$ is symmetric and
${\mb X^*}\t\mb G^{t,(i)}\mb Y^*
+{\mb Y^*}\t{\mb G^{t,(i)}}\t\mb X^*$ is symmetric, the identity is the
Procrustes alignment of $\mb F^{t+1,(i),\sf{srg}}$ with $\mb F^*$. Lemma~\ref{lemma: perturbation theory of optimal rotation under F norm} now
gives
\begin{align}
\bignorm{\widetilde{\mb O}^{t,(i)\top}
\widetilde{\mb O}^{t+1,(i)}-\mb I_r}
&\lesssim\frac{\sqrt{\sigma_1^*}}{\sigma_r^*}
\bignorm{\mb F^{t+1,(i)}\widetilde{\mb O}^{t,(i)}
-\mb F^{t+1,(i),\sf{srg}}}.
\label{eq: corrected rotation deviation}
\end{align}

Separating the sampled-complement and held-out-row contributions in the two
factor updates yields the following decomposition:
\begin{align}
&\bignorm{\mb F^{t+1,(i)}\widetilde{\mb O}^{t,(i)}
-\mb F^{t+1,(i),\sf{srg}}}
\nonumber\\
&\quad\leq \eta\Big\{
\underbrace{\Bignorm{
\mc P_\Omega^{(i),\complement}\big[
\psi(\mc M(\bo\theta^{t,(i)}))-\mb R
\big]
\big(\mb Y^{t,(i)}\widetilde{\mb O}^{t,(i)}-\mb Y^*\big)
}\fb}_{\gamma_1}
\nonumber\\
&\qquad+\pi\underbrace{\Bignorm{
\mc P^{(i)}\Big(
\psi(\mc M(\bo\theta^{t,(i)}))-\mb R^*
\Big)
\big(\mb Y^{t,(i)}\widetilde{\mb O}^{t,(i)}-\mb Y^*\big)
}\fb}_{\gamma_2}
\nonumber\\
&\qquad+\underbrace{\Bignorm{
\Big\{\mc P_\Omega^{(i),\complement}\big[
\psi(\mc M(\bo\theta^{t,(i)}))-\mb R
\big]\Big\}\t
\big(\mb X^{t,(i)}\widetilde{\mb O}^{t,(i)}-\mb X^*\big)
}\fb}_{\gamma_3}
\nonumber\\
&\qquad+\pi\underbrace{\Bignorm{
\Big\{\mc P^{(i)}\Big(
\psi(\mc M(\bo\theta^{t,(i)}))-\mb R^*
\Big)\Big\}\t
\big(\mb X^{t,(i)}\widetilde{\mb O}^{t,(i)}-\mb X^*\big)
}\fb}_{\gamma_4}\Big\}. 
\label{eq: corrected rotation-surrogate remainder}
\end{align}
In the sequel, we shall bound each of the four terms $\gamma_1,\ldots,\gamma_4$ in \eqref{eq: corrected rotation-surrogate remainder} separately. 
\begin{itemize}
    \item For $\gamma_1$, it first follows from the triangle inequality that 
    \begin{align}
        & \gamma_1 \leq \Bignorm{
\mc P_\Omega^{(i),\complement}\big[
\psi(\mc M(\bo\theta^{t,(i)}))-\psi(\mb M^*)
\big]\big(\mb Y^{t,(i)}\widetilde{\mb O}^{t,(i)}-\mb Y^*\big)
}\fb\\
    &\quad + \Bignorm{
\mc P_\Omega^{(i),\complement}\big[
\mb E
\big]\big(\mb Y^{t,(i)}\widetilde{\mb O}^{t,(i)}-\mb Y^*\big)
}\fb \\ 
    & \leq   \Bignorm{\mc P_\Omega\big[
\psi(\mc M(\bo\theta^{t,(i)}))-\psi(\mb M^*)
\big]\big(\mb Y^{t,(i)}\widetilde{\mb O}^{t,(i)}-\mb Y^*\big)
}\fb\\
    &\quad + \Bignorm{
\mc P_\Omega^{(i),\complement}\big[
\mb E
\big]\big(\mb Y^{t,(i)}\widetilde{\mb O}^{t,(i)}-\mb Y^*\big)
}\fb . 
\label{eq: decomposition of gamma1}
    \end{align}

    The inductive row bounds and
    \eqref{eq: rowwise transfer from actual to LOO} place
    $\bo\theta^{t,(i)}$ in the RIC after a common factor rotation.
    The domain verification in the local-geometry proof therefore places
    its predictor, the true predictor, and their natural-parameter
    interpolants in $\mc D_\psi$. To control the first term above, we apply the
    Lemma~\ref{lemma: injectivity of P_Omega under nonlinearity} bound
    \eqref{eq: uniform nonlinear sampling operator bound} with these
    admissible endpoints $\bo\theta^{t,(i)}$ and $\bo\theta^*$, and then invoke
    Lemma~\ref{lemma: lemma 4.4 chen2019} and
    \eqref{eq: LOO global error used in row propagation}. This gives, with
    probability at least $1-O(d^{-c-2})$,
    \begin{align}
        & \Bignorm{
\mc P_\Omega\big[
\psi(\mc M(\bo\theta^{t,(i)}))-\psi(\mb M^*)
\big]\big(\mb Y^{t,(i)}\widetilde{\mb O}^{t,(i)}-\mb Y^*\big)
}\fb  \\ 
\leq & \Bignorm{
\big(\bo \Omega-  \pi \mb 1_n \mb 1_p\t \big) \circ \big[
\psi(\mc M(\bo\theta^{t,(i)}))-\psi(\mb M^*)
\big]\big(\mb Y^{t,(i)}\widetilde{\mb O}^{t,(i)}-\mb Y^*\big)
}\fb \\ 
& + \pi \norm{\psi(\mc M(\bo\theta^{t,(i)}))-\psi(\mb M^*)
\big(\mb Y^{t,(i)}\widetilde{\mb O}^{t,(i)}-\mb Y^*\big)}\fb \\ 
\lesssim & \sigma^2 \sqrt{\pi d}\xi^4 \sqrt{\frac{\mu \rho r}{d}} \sigma_1^* \bignorm{\mb Y^{t,(i)}\widetilde{\mb O}^{t,(i)} - \mb Y^*}\ti + \sigma^2 \pi \sqrt{\sigma_1^*} \distf(\bo \theta^{t,(i)}, \bo \theta^*)^2 \\
\ll & \underline c_\psi \pi \sigma_r^* \sqrt{\frac{n\wedge p}{\mu r \kappa^2}}\Big[ \nu^{t}\Big\{
\kappa\kappa_\psi\sqrt{\mu\rho r}\,\delta_{\sf{row}}^0+\delta_{\sf{loo}}^0
\Big\}
+\frac{\sigma\kappa\kappa_\psi
\sqrt{\kappa\mu\rho r^2\log d}}
{\underline c_\psi\sqrt{\pi\sigma_r^*}} \Big]. 
    \end{align}

On the other hand, the second term can be directly bounded by 
\begin{align}
    & \Bignorm{
\mc P_\Omega^{(i),\complement}\big[
\mb E
\big]\big(\mb Y^{t,(i)}\widetilde{\mb O}^{t,(i)}-\mb Y^*\big)
}\fb \\ 
\leq & \norm{\mc P_\Omega(\mb E)} \norm{\mb Y^{t,(i)}\widetilde{\mb O}^{t,(i)}-\mb Y^*}\fb  \\
\ll & \underline c_\psi \pi \sigma_r^* \sqrt{\frac{n\wedge p}{\mu r \kappa^2}}\Big[ \nu^{t}\Big\{
\kappa\kappa_\psi\sqrt{\mu\rho r}\,\delta_{\sf{row}}^0+\delta_{\sf{loo}}^0
\Big\}
+\frac{\sigma\kappa\kappa_\psi
\sqrt{\kappa\mu\rho r^2\log d}}
{\underline c_\psi\sqrt{\pi\sigma_r^*}} \Big]. 
\end{align}
This holds with probability at least $1-O(d^{-c-2})$ by
Lemma~\ref{lemma: noise spectral norm concentration}.

\item The second term $\gamma_2$ can be directly handled by \eqref{eq: condition -- rowwise consistency of LOO iterates} and the Taylor expansion to yield
\begin{align}
    & \gamma_2 \lesssim \sigma^2 \pi \sqrt{\sigma_1^*} \bignorm{(\mb X^{t,(i)}\widetilde{\mb O}^{t,(i)}-\mb X^*)_{i,\cdot}}_2 \distf(\bo\theta^{t,(i)},\bo\theta^*) + \sigma^2 \pi \sqrt{\frac{\mu \rho r \sigma_1^*}{d}} \distf(\bo\theta^{t,(i)},\\
    &\quad \bo\theta^*)^2\\
    & \ll   \underline c_\psi \pi \sigma_r^* \sqrt{\frac{n\wedge p}{\mu r \kappa^2}}\Big[ \nu^{t}\Big\{
\kappa\kappa_\psi\sqrt{\mu\rho r}\,\delta_{\sf{row}}^0+\delta_{\sf{loo}}^0
\Big\}
+\frac{\sigma\kappa\kappa_\psi
\sqrt{\kappa\mu\rho r^2\log d}}
{\underline c_\psi\sqrt{\pi\sigma_r^*}} \Big]. 
\end{align}

\item A decomposition similar to \eqref{eq: decomposition of gamma1}  gives
\begin{equation}
\begin{aligned}
    & \gamma_3 \leq \Bignorm{\mc P_\Omega^{(i),\complement}\big[
\psi(\mc M(\bo\theta^{t,(i)}))-\psi(\mb M^*)
\big]\t\big(\mb X^{t,(i)}\widetilde{\mb O}^{t,(i)}-\mb X^*\big)
}\fb\\
    &\quad + \Bignorm{
\mc P_\Omega^{(i),\complement}\big[
\mb E
\big]\t\big(\mb X^{t,(i)}\widetilde{\mb O}^{t,(i)}-\mb X^*\big)
}\fb, 
\end{aligned}
\end{equation}
where both terms can be controlled by the same argument as for $\gamma_1$.  

\item The fourth term $\gamma_4$ can be handled as follows:
\begin{align}
    & \gamma_4 \leq \pi \Bignorm{\big[
\psi(\mc M(\bo\theta^{t,(i)}))-\psi(\mb M^*)
\big]_i
}_2  \norm{\mb X^{t,(i)}\widetilde{\mb O}^{t,(i)}-\mb X^*}\ti \\ 
\lesssim & \sigma^2 \sqrt{n} \sqrt{\frac{\mu \rho r \sigma_1^*}{n}} \dist\ti(\bo\theta^{t,(i)},\bo\theta^*)^2 \\
\ll & \underline c_\psi \pi \sigma_r^* \sqrt{\frac{n\wedge p}{\mu r \kappa^2}}\Big[ \nu^{t}\Big\{
\kappa\kappa_\psi\sqrt{\mu\rho r}\,\delta_{\sf{row}}^0+\delta_{\sf{loo}}^0
\Big\}
+\frac{\sigma\kappa\kappa_\psi
\sqrt{\kappa\mu\rho r^2\log d}}
{\underline c_\psi\sqrt{\pi\sigma_r^*}} \Big]. 
\end{align}
\end{itemize}

Combining the preceding bounds for $\gamma_1,\ldots,\gamma_4$ and substituting into \eqref{eq: corrected rotation-surrogate remainder} gives
\begin{align}
    & \bignorm{\mb F^{t+1,(i)}\widetilde{\mb O}^{t,(i)}
-\mb F^{t+1,(i),\sf{srg}}}
\ll \eta \underline c_\psi \pi \sigma_r^* \sqrt{\frac{n\wedge p}{\mu r \kappa^2}}\Big[ \nu^{t}\Big\{
\kappa\kappa_\psi\sqrt{\mu\rho r}\,\delta_{\sf{row}}^0+\delta_{\sf{loo}}^0
\Big\}\\
    &\quad +\frac{\sigma\kappa\kappa_\psi
\sqrt{\kappa\mu\rho r^2\log d}}
{\underline c_\psi\sqrt{\pi\sigma_r^*}} \Big].
\label{eq: corrected rotation-surrogate remainder bound}
\end{align} 

In view of \eqref{eq: corrected rotation deviation}, Assumption~\ref{assumption: signal strength and incoherence degree}~and~\eqref{eq: epsilon GD initialization requirement} ensure that 
\begin{align}
 & \bignorm{\mb X^{t+1,(i)}}\ti
\bignorm{\widetilde{\mb O}^{t,(i)\top}
\widetilde{\mb O}^{t+1,(i)}-\mb I_r}\\
\ll &\eta \kappa \sqrt{\frac{\mu r}{n}} \Big\{\frac{\sigma^2 \sqrt{\pi }\sigma_1^*}{\kappa \kappa_\psi}   \Big[ \nu^t \big[ \sqrt d\,\delta_{\sf{row}}^0+\delta_{\sf{loo}}^0 \big]
+ \frac{\sigma
\sqrt{d \kappa r \log d}}
{\underline c_\psi\sqrt{\pi\sigma_r^*}} \Big] + \sigma \sqrt{\pi d } \Big\} \\
 \ll & \eta \underline c_\psi \pi \sigma_r^* \Big[ \nu^{t}\Big\{
\kappa\kappa_\psi\sqrt{\mu\rho r}\,\delta_{\sf{row}}^0+\delta_{\sf{loo}}^0
\Big\}
+\frac{\sigma\kappa\kappa_\psi
\sqrt{\kappa\mu\rho r^2\log d}}
{\underline c_\psi\sqrt{\pi\sigma_r^*}} \Big]. 
\label{eq: rotation error contration}
\end{align}

Collecting the preceding bounds \eqref{eq: held-out row contraction}, \eqref{eq: second term in the LOO rowwise}, \eqref{eq: weighted cross-factor bound from Frobenius error}, \eqref{eq: rotation error contration}, and~\eqref{eq: condition -- orthogonality between 1 and Xt} into \eqref{eq: corrected held-out row decomposition} yields
\begin{align}
& \bignorm{(\mb X^{t+1,(i)}\tilde{\mb O}^{t+1,(i)}-\mb X^*)_{i,\cdot}}_2
\leq{}\Big(1-\frac{\eta\underline c_\psi\pi\sigma_r^*}{4}\Big)
\bignorm{(\mb X^{t,(i)}\widetilde{\mb O}^{t,(i)}-\mb X^*)_{i,\cdot}}_2  \\ 
& \qquad + C\big[ \eta \underline c_\psi \pi \sigma_r^*\big] \Big[ \nu^{t}\Big\{
\kappa\kappa_\psi\sqrt{\mu\rho r}\,\delta_{\sf{row}}^0+\delta_{\sf{loo}}^0
\Big\}
+\frac{\sigma\kappa\kappa_\psi
\sqrt{\kappa\mu\rho r^2\log d}}
{\underline c_\psi\sqrt{\pi\sigma_r^*}} \Big]. 
\label{eq: held-out row recursion under general link}
\end{align}
for some constant $C$. By taking sufficiently large constants $C_7$ and $C_8$ in \eqref{eq: condition -- rowwise consistency of LOO iterates}, the first term on the right-hand side of \eqref{eq: held-out row recursion under general link} is absorbed into the contraction, yielding \eqref{eq: propagated LOO rowwise error}.

Lastly, noticing the relation \eqref{eq: leave-one-out decomposition}
\begin{align}
    & \norm{(\mb X^{t+1} \tilde{\mb O}^{t+1} - \mb X^*)_i}_2 
 \leq \norm{(\mb X^{t+1,(i)} \mb O^{t+1,(i)} - \mb X^*)_i}_2
 +\Big\{1+\kappa\sqrt{\frac{\mu r}{n\wedge p}}\Big\}
 \distf(\bo\theta^{t+1},\\
    &\quad \bo\theta^{t+1,(i)})
\nonumber\\
&\leq  \norm{(\mb X^{t+1,(i)} \mb O^{t+1,(i)} - \mb X^*)_i}_2
 + 2
 \distf(\bo\theta^{t+1},\bo\theta^{t+1,(i)})
\end{align}
we conclude \eqref{eq: propagated rowwise GD error} from \eqref{eq: propagated LOO rowwise error} and Lemma~\ref{lemma: leave-one-out error} for $C_5 = 2C_1 + C_7$ and $C_6 = 2C_2 + C_8$.

For the columnwise argument, fix $j\in[p]$ and use the augmented factors
\begin{align*}
\mb X^{\sf{app},t,(-j)}
&=\big(\sqrt{\sigma_r^*/n}\mb 1_n,\mb X^{t,(-j)}\big),
&
\mb Y^{\sf{app},t,(-j)}
&=\big(\sqrt{n/\sigma_r^*}\bo\zeta^{t,(-j)},\mb Y^{t,(-j)}\big),
\end{align*}
{\emergencystretch=1.5em
together with their population counterparts. The associated common rotation
is $\operatorname{diag}(1,\widetilde{\mb O}^{t,(-j)})$.
Under this scaling, the proof of \eqref{eq: propagated rowwise GD error} is identical to the rowwise case, with the roles of $\mb X$ and $\mb Y$ interchanged. We omit the details for brevity.
\qed\par}

\subsubsection{Proof of Lemma~\ref{lemma: balancedness and orthogonality}}
We first establish four auxiliary bounds used in the proof.

    To arrive at this subtle balancedness, we first need the following lemma to sort out the relation between the balancedness and the learning rate as well as the last-step gradient norm. 
    \begin{lemma}
    \label{lemma: bound for Xt+1t Xt+1 - Yt+1t Yt+1}
    Work in the balanced parametrization $\omega=1$. Fix $0\leq t<t_0$ and
    suppose that the current iterate satisfies
    $\bignorm{\mb X^t}\vee\bignorm{\mb Y^t}
    \leq 2\sqrt{\sigma_1^*}$, and
    let $\bo\theta^{t+1}$ be obtained from $\bo\theta^t$ by one update of
    Algorithm~\ref{alg: GD} with any step size $\eta>0$. Then for some constant $C>0$,
    \begin{align*}
        &\bignorm{{\mb X^{t+1}}\t\mb X^{t+1}-{\mb Y^{t+1}}\t\mb Y^{t+1}}\fb \\
    & \leq{} \bignorm{{\mb X^t}\t\mb X^t-{\mb Y^t}\t\mb Y^t}\fb
        +C\Big\{\frac{\eta c_{\perp}\pi\sigma_r^*}{n}
        \bignorm{\mb 1_n\t\mb X^t}_2^2 +\eta^2\bignorm{\nabla_{\bo\theta}L^{\sf{aug}}(\bo\theta^t)}_2^2\\
    &\quad +c_{\sf{aug}}^2\eta^2\pi^2\sigma_1^*
        \bignorm{{\mb X^t}\t\mb X^t-{\mb Y^t}\t\mb Y^t}\fb^2\Big\}. 
    \end{align*}
    \end{lemma}
    \begin{proof}
        Firstly, for the $(t+1)$-th step, we notice that the cross term $\mb X^{t\top} (\nabla_{\mb X} L(\bo \theta^t) - 2 c_{\perp} \frac{\pi \sigma_r^*}{n}\mb 1_n\mb 1_n\t \mb X^t) + (\nabla_{\mb X}L(\bo \theta^t) - 2 c_{\perp} \frac{\pi \sigma_r^*}{n}\mb 1_n\mb 1_n\t \mb X^t)\t \mb X^t - \mb Y^{t\top} \nabla_{\mb Y} L(\bo \theta^t) - \nabla_{\mb Y}L(\bo \theta^t)\t \mb Y^t$ vanishes, which implies the following decomposition
        \begin{align}
            & {\mb X^{t+1}}\t\mb X^{t+1} - {\mb Y^{t+1}}\t \mb Y^{t+1}\\
            = &{\mb X^{t}}\t \mb X^t - {\mb Y^t}\t \mb Y^t - 4 c_{\perp} \frac{\pi \sigma_r^*}{n} \eta {\mb X^t}\t \mb 1_n\mb 1_n\t \mb X^t\\
            &  +  \underbrace{\eta^2 \nabla_{\mb X}L(\bo \theta^t)\t \nabla_{\mb X}L(\bo \theta^t)}_{\beta_1} -  \underbrace{\eta^2 \nabla_{\mb Y} L(\bo \theta^t)\t \nabla_{\mb Y}L(\bo \theta^t) }_{\beta_2},
        \end{align}
     based on which we delve into the decompositions of $\beta_1$ and $\beta_2$ in the sequel. 

    \begin{itemize}
        \item With respect to $\beta_1$, we recall the definition of $L^{\sf{aug}}$ and $L^{\sf{diff}}$ in \eqref{eq: definition of Laug} and \eqref{eq: definition of Ldiff} and write 
    \begin{align}
        &\beta_1 = \eta^2 \Big[\nabla_{\mb X} L^{\sf{aug}}(\bo \theta^t) - \nabla_{\mb X} L^{\sf{diff}}(\bo \theta^t)\Big]\t  \Big[\nabla_{\mb X} L^{\sf{aug}}(\bo \theta^t) - \nabla_{\mb X} L^{\sf{diff}}(\bo \theta^t)\Big] \\ 
    & =   \eta^2 \Big[\nabla_{\mb X} L^{\sf{aug}}(\bo \theta^t)\t \nabla_{\mb X} L^{\sf{aug}}(\bo \theta^t) -\nabla_{\mb X} L^{\sf{aug}}(\bo \theta^t)\t  \nabla_{\mb X} L^{\sf{diff}}(\bo \theta^t)\\
    &\quad - \nabla_{\mb X} L^{\sf{diff}}(\bo \theta^t)\t \nabla_{\mb X} L^{\sf{aug}}(\bo \theta^t)\Big]\\
        & + \eta^2 \nabla_{\mb X} L^{\sf{diff}}(\bo \theta^t)\t \nabla_{\mb X} L^{\sf{diff}}(\bo \theta^t)\\ 
    & =   \eta^2 \Big[\nabla_{\mb X} L^{\sf{aug}}(\bo \theta^t)\t \nabla_{\mb X} L^{\sf{aug}}(\bo \theta^t) -\nabla_{\mb X} L^{\sf{aug}}(\bo \theta^t)\t  \nabla_{\mb X} L^{\sf{diff}}(\bo \theta^t)\\
    &\quad - \nabla_{\mb X} L^{\sf{diff}}(\bo \theta^t)\t \nabla_{\mb X} L^{\sf{aug}}(\bo \theta^t)\Big]\\
        & + 16 c_{\sf{aug}}^2\eta^2\pi^2  \Big[{\mb X^t}\t \mb X^t -{\mb Y^t}\t \mb Y^t\Big] {\mb X^t}\t \mb X^t \Big[{\mb X^t}\t \mb X^t -{\mb Y^t}\t \mb Y^t\Big]. 
    \end{align}
    It immediately follows from the triangle inequality that
    \begin{align*}
        & \norm{\beta_1}\fb \lesssim \eta^2 \norm{\nabla_{\mb X} L^{\sf{aug}}(\bo \theta^t)}\fb^2 + \pi \eta^2 \sqrt{\sigma_1^*} \bignorm{{\mb X^t}\t \mb X^t - {\mb Y^t}\t \mb Y^t}\fb \norm{\nabla_{\mb X} L^{\sf{aug}}(\bo \theta^t)}\fb \\ 
        &  \qquad + c_{\sf{aug}}^2 \eta^2 \pi^2 \bignorm{{\mb X^t}\t \mb X^t - {\mb Y^t}\t \mb Y^t}\fb^2 \sigma_1^* 
    \end{align*}
    \item Analogously, one has for $\beta_2$ that 
    \begin{align}
        & \norm{\beta_2}\fb \lesssim \eta^2 \norm{\nabla_{\mb Y} L^{\sf{aug}}(\bo \theta^t)}\fb^2 +  \eta^2 \pi \sqrt{\sigma_1^*} \bignorm{{\mb X^t}\t \mb X^t - {\mb Y^t}\t \mb Y^t}\fb \bignorm{\nabla_{\mb Y} L^{\sf{aug}}(\bo \theta^t)}\fb \\ 
        &\qquad  +  c_{\sf{aug}}^2\eta^2 \pi^2  \bignorm{{\mb X^t}\t \mb X^t - {\mb Y^t}\t \mb Y^t}\fb^2 \sigma_1^*. 
    \end{align}
    \end{itemize}
    
    Combining these two parts yields that 
    \begin{align*}
        & \bignorm{{\mb X^{t+1}}\t\mb X^{t+1}  - {\mb Y^{t+1}}\t \mb Y^{t+1} }\fb  
        \leq   \bignorm{{\mb X^{t}}\t \mb X^t - {\mb Y^t}\t \mb Y^t}\fb \\
        &  + C\Big[ \frac{\eta c_{\perp} \pi \sigma_r^*}{n}\bignorm{\mb 1_n\t \mb X^t}_2^2 + \eta^2 \norm{\nabla L^{\sf{aug}}(\bo \theta^t)}\fb^2\\ 
    & \qquad  +  \eta^2 \pi \sqrt{\sigma_1^*} \bignorm{{\mb X^t}\t \mb X^t - {\mb Y^t}\t \mb Y^t}\fb \norm{\nabla L^{\sf{aug}}(\bo \theta^t)}\fb\\
    &\quad +   
         c_{\sf{aug}}^2\eta^2 \pi^2 \bignorm{{\mb X^t}\t \mb X^t - {\mb Y^t}\t \mb Y^t}\fb^2 \sigma_1^* \Big] \\ 
    & \leq    \bignorm{{\mb X^{t}}\t \mb X^t - {\mb Y^t}\t \mb Y^t}\fb   + C\Big[ \frac{\eta c_{\perp} \pi \sigma_r^*}{n}\bignorm{\mb 1_n\t \mb X^t}_2^2 + \eta^2 \norm{\nabla L^{\sf{aug}}(\bo \theta^t)}\fb^2\\
    &\quad + 
        c_{\sf{aug}}^2 \eta^2 \pi^2 \bignorm{{\mb X^t}\t \mb X^t - {\mb Y^t}\t \mb Y^t}\fb^2 \sigma_1^* \Big] 
    \end{align*}
    for some constant $C$, which proves the claim.

    \end{proof}

    Regarding the orthogonality of $\mb X^t$ to $\mb 1_n$, the penalty $c_{\perp}\norm{\mb 1_n\t\mb X^t}_2^2$ maintains approximate centering along the gradient-descent trajectory.
    \begin{lemma}
    \label{lemma: bound for 1t Xt}
    Suppose that the assumptions of Theorem~\ref{thm: general GD} hold. Fix $0\leq t<t_0$ and suppose that Condition~\ref{condition: inductive contraction for GD} holds at step $t$. Then, with probability at least $1-O(d^{-c-1})$,
    \begin{align*}
    & \bignorm{\mb 1_n\t\mb X^{t+1}}_2
        \leq{}  \big(1-2c_{\perp}\eta\pi\sigma_r^*\big)
        \bignorm{\mb 1_n\t\mb X^t}_2
        +C\eta\sigma^2\pi\sqrt n\,\sigma_1^*
        \distf(\bo\theta^t,\bo\theta^*)\\
    &\quad +C\eta\sigma n\sqrt{\pi\sigma_1^*\log d},
    \end{align*}
    for some universal constant $C>0$.
    \end{lemma}
    \begin{proof}
    The step-size condition ensures that
    $1-2c_{\perp}\eta\pi\sigma_r^*\geq0$. Since
    $\mb R=\psi(\mb M^*)+\mb E$, the gradient-descent update gives
    \begin{align*}
        \bignorm{\mb 1_n\t\mb X^{t+1}}_2
        \leq{}& \big(1-2c_{\perp}\eta\pi\sigma_r^*\big)
        \bignorm{\mb 1_n\t\mb X^t}_2+\eta\beta_1+\eta\beta_2,\\
        \beta_1
        &\coloneqq
        \bignorm{\mb 1_n\t\mc P_\Omega\big[
        \psi(\mc M(\bo\theta^t))-\psi(\mb M^*)
        \big]\mb Y^t}_2,\\
        \beta_2
        &\coloneqq
        \bignorm{\mb 1_n\t\mc P_\Omega(\mb E)\mb Y^t}_2.
    \end{align*}

    We first bound $\beta_1$. Set
    $H_{i,j}^t\coloneqq\int_0^1\psi'\big((1-s)M_{i,j}^*
    +s\mc M(\bo\theta^t)_{i,j}\big)\,\mathrm ds$,
    $\mb B^t\coloneqq\mc P_\Omega(\mb H^t)-\pi\mb H^t$, and
    $\mb D^t\coloneqq\mc M(\bo\theta^t)-\mb M^*$. The mean-value theorem gives
    $\psi(\mc M(\bo\theta^t))-\psi(\mb M^*)=\mb H^t\circ\mb D^t$.
    The inductive RIC bounds place the endpoint predictors and their
    natural-parameter interpolants in $\mc D_\psi$ and supply the factor
    row bounds defining $\mc M_\psi$. The
    Lemma~\ref{lemma: injectivity of P_Omega under nonlinearity} bound
    \eqref{eq: uniform nonlinear sampling operator bound}, applied with
    endpoints $\bo\theta^t$ and $\bo\theta^*$, yields
    $\norm{\mb B^t}\lesssim\sigma^2\xi^4\sqrt{\pi d}$.
    Then one has 
    \begin{align*}
    & \mb D^t
        ={} \mb 1_n(\bo\zeta^t-\bo\zeta^*)\t
        +(\mb X^t\mb O^t-\mb X^*)(\mb Y^t\mb O^t)\t
        +\mb X^*(\mb Y^t\mb O^t-\mb Y^*)\t,\\
    & \beta_1
        \leq{} \pi\sqrt n\,\bignorm{\mb H^t\circ\mb D^t}\fb
        \norm{\mb Y^t}
        +
        \norm{\mb B^t}\bignorm{\mb Y^t\mb O^t}\ti
        \Big\{
        \sqrt n\,\bignorm{\bo\zeta^t-\bo\zeta^*}_2\\
    &\quad +\norm{\mb Y^t}\bignorm{\mb X^t\mb O^t-\mb X^*}\fb
        +\bignorm{\mb X^*}\fb
        \bignorm{\mb Y^t\mb O^t-\mb Y^*}\fb
        \Big\}\\
    & \lesssim{} 
        \sigma^2\pi\sqrt n\,\sigma_1^*
        \distf(\bo\theta^t,\bo\theta^*).
    \end{align*}
    The last step uses
    $\xi^4r\sqrt{\mu d/(\pi np)}\ll1$, which follows from
    Assumption~\ref{assumption: signal strength and incoherence degree}(a).

    For $\beta_2$, insert the leave-one-row-out factor and condition on
    $\mb Y^{t,(i)}$. For each fixed $i$, the matrix Bernstein bounds below
    hold with failure probability $O(d^{-c-2})$; a union bound over $i\in[n]$
    yields failure probability $O(d^{-c-1})$. Thus,
    \begin{align*}
        \beta_2
        \leq{}&
        n\max_{i\in[n]}\Big\{
        \bignorm{\mc P_\Omega(\mb E)_{i,\cdot}\mb Y^{t,(i)}}_2
        +
        \bignorm{\mc P_\Omega(\mb E)_{i,\cdot}}_2
        \bignorm{\mb Y^{t,(i)}\mb O^{t,(i)}-\mb Y^t}\fb
        \Big\}\\
        \lesssim{}&n\Big\{
        \sigma\sqrt{\pi\sigma_1^*r\log d}
        +B\sqrt{\frac{\mu r\sigma_1^*}{p}}\,r\log d
        +\big(\sigma\sqrt{\pi p}+B\sqrt{\log d}\big)
        \max_{i\in[n]}\distf(\bo\theta^t,\bo\theta^{t,(i)})
        \Big\}\\
        \lesssim{} &  \sigma n\sqrt{\pi\sigma_1^*\log d}.
    \end{align*}
    Here the last inequality follows from
    \eqref{eq: condition -- F-norm consistency between thetat and thetatl},
    \eqref{eq: epsilon GD initialization requirement}, and the stated
    sampling, signal-strength, and tail conditions. Combining the bounds for
    $\beta_1$ and $\beta_2$ with the update recursion proves the claim.
    \end{proof}
     The next lemma characterizes the vanishing rate of the gradient under the augmented function as $t$ grows. 
    \begin{lemma}
        \label{lemma: bound on nabla Laug}
        Suppose that the assumptions of Theorem~\ref{thm: general GD} hold.
Fix $0\leq t<t_0$ and suppose that
Condition~\ref{condition: inductive contraction for GD} holds for all $0 \leq s \leq t$. Then we have for some constant $C>0$ that
        \begin{align*}
            & \norm{\nabla L^{\sf{aug}}(\bo \theta^{t+1})}_2 \leq \nu\norm{\nabla L^{\sf{aug}}(\bo \theta^{t})}_2 + C \eta c_\perp \sigma^2 \pi^2(\sigma_1^*)^{\frac32} \bignorm{{\mb X^t}\t \mb X^t - {\mb Y^t}\t \mb Y^t }\fb
            \\
            \lesssim & \nu^{t+1} \norm{\nabla L^{\sf{aug}}(\bo \theta^0)}_2 + C c_\perp \kappa_\psi \kappa \pi(\sigma_1^*)^{\frac12}\max_{0 \leq s \leq t} \bignorm{{\mb X^s}\t \mb X^s-{\mb Y^s}\t \mb Y^s}\fb 
        \end{align*}
        with probability at least $1- O(d^{-c-1})$. 
    \end{lemma}
    \begin{proof}
        Define $\overline{\mb H}_{\sf g}^t\coloneqq\int_0^1
        \nabla^2L^{\sf{aug}}\bigl((1-s)\bo\theta^t+s\bo\theta^{t+1}\bigr)\,\mathrm ds$.
        Taylor's theorem and $\nabla L=\nabla L^{\sf{aug}}-\nabla L^{\sf{diff}}$ give
        \begin{align*}
            \nabla L^{\sf{aug}}(\bo\theta^{t+1})
            &=\big(\mb I-\eta\overline{\mb H}_{\sf g}^t\big)
            \nabla L^{\sf{aug}}(\bo\theta^t)
            +\eta\overline{\mb H}_{\sf g}^t\nabla L^{\sf{diff}}(\bo\theta^t).
        \end{align*}
        Rotational invariance implies that the factor block $\mb G^t$ of $\nabla L^{\sf{aug}}(\bo\theta^t)$ is a horizontal tangent direction at $\mb F^t$, namely, ${\mb F^t}\t\mb G^t$ is symmetric. Hence the tangent-direction version of Lemma~\ref{lemma: bounds for Hessian} applies to $\nabla L^{\sf{aug}}(\bo\theta^t)$ at every point of the interpolation segment. The inductive hypotheses ensure separately that this segment remains in the RIC. The argument establishing \eqref{eq: averaged Hessian multiplication bound LOO}, applied to this tangent direction, therefore yields
        \begin{equation*}
            \bignorm{\overline{\mb H}_{\sf g}^t\nabla L^{\sf{aug}}(\bo\theta^t)}_2^2
            \leq2\bar c_{\sf{Hess}}c_\perp\pi\sigma_1^*
            \nabla L^{\sf{aug}}(\bo\theta^t)\t
            \overline{\mb H}_{\sf g}^t\nabla L^{\sf{aug}}(\bo\theta^t).
        \end{equation*}
        Using $\eta\leq1/(2\bar c_{\sf{Hess}}c_\perp\pi\sigma_1^*)$ and
        squaring the first term in the preceding decomposition gives
        \begin{align*}
            &\bignorm{\big(\mb I-\eta\overline{\mb H}_{\sf g}^t\big)
            \nabla L^{\sf{aug}}(\bo\theta^t)}_2^2\\
            &\quad=\bignorm{\nabla L^{\sf{aug}}(\bo\theta^t)}_2^2
            -2\eta\nabla L^{\sf{aug}}(\bo\theta^t)\t
            \overline{\mb H}_{\sf g}^t\nabla L^{\sf{aug}}(\bo\theta^t)
            +\eta^2\bignorm{\overline{\mb H}_{\sf g}^t
            \nabla L^{\sf{aug}}(\bo\theta^t)}_2^2\\
            &\quad\leq\big(1-\eta\underline c_{\sf{Hess}}\pi\sigma_r^*\big)
            \bignorm{\nabla L^{\sf{aug}}(\bo\theta^t)}_2^2.
        \end{align*}
        Finally, the upper-Hessian bound and \eqref{eq: derivatives of Ldiff} give
        \begin{align*}
            \eta\bignorm{\overline{\mb H}_{\sf g}^t\nabla L^{\sf{diff}}(\bo\theta^t)}_2
            &\leq  \sigma^2 c_\perp \eta\pi^2(\sigma_1^*)^{3/2}
            \bignorm{{\mb X^t}\t\mb X^t-{\mb Y^t}\t\mb Y^t}\fb.
        \end{align*}
        Combining the last two displays proves the one-step recursion. Summing
        the resulting geometric series and using
        $1-\nu=\eta\underline c_{\sf{Hess}}\pi\sigma_r^*/4$ together with
        $\underline c_{\sf{Hess}}\asymp\underline c_\psi$ proves the second bound.
    \end{proof} 

    \begin{lemma}
        \label{lemma: gradient norm at theta0}
        Suppose that the assumptions of Theorem~\ref{thm: general GD} hold. Then, with probability at least $1-O(d^{-c-1})$,
        \begin{align}
            & \norm{\nabla L^{\sf{aug}}(\bo \theta^0)}_2 = \norm{\nabla L(\bo \theta^0)}_2 \lesssim \sigma \sqrt{d \pi r\sigma_1^*} + \kappa_\psi
            \sqrt{\kappa d \pi\sigma_1^*} \xi^{\sf{gd}}. 
        \end{align}
    \end{lemma}
    \begin{proof}
        The identities $\mb 1_n\t\mb X^0=\mb 0$ and
        ${\mb X^0}\t\mb X^0={\mb Y^0}\t\mb Y^0$ imply that the centering
        and balancing penalties have vanishing gradients at $\bo\theta^0$.
        Hence
        \begin{align}
            \norm{\nabla L^{\sf{aug}}(\bo\theta^0)}_2
            ={}&\norm{\nabla L(\bo\theta^0)}_2\\
            \leq{}&\bignorm{\mc P_\Omega\big[\psi(\mc M(\bo\theta^0))-\mb R\big]\mb Y^0}\fb
            +\bignorm{\mc P_\Omega\big[\psi(\mc M(\bo\theta^0))-\mb R\big]\t\mb X^0}\fb\\
            &+\sqrt{\frac{\sigma_r^*}{n}}
            \bignorm{\mc P_\Omega\big[\psi(\mc M(\bo\theta^0))-\mb R\big]\t\mb 1_n}_2.
        \end{align} 

        We treat the first term; the other two follow by the same argument.
        Set $\mb D^0\coloneqq\mc M(\bo\theta^0)-\mb M^*$,
        $H_{i,j}^0\coloneqq\int_0^1\psi'\big((1-s)M_{i,j}^*
        +s\mc M(\bo\theta^0)_{i,j}\big)\,\mathrm ds$, and
        $\mb B^0\coloneqq\mc P_\Omega(\mb H^0)-\pi\mb H^0$. The mean-value
        theorem gives $\psi(\mc M(\bo\theta^0))-\psi(\mb M^*)
        =\mb H^0\circ\mb D^0$. The initialization requirements place
        $\bo\theta^0$ in the RIC after a common factor rotation. Hence both
        endpoints satisfy the factor row bounds, and their predictors and
        natural-parameter interpolants belong to $\mc D_\psi$. The
        Lemma~\ref{lemma: injectivity of P_Omega under nonlinearity} bound
        \eqref{eq: uniform nonlinear sampling operator bound}, applied with
        endpoints $\bo\theta^0$ and $\bo\theta^*$, yields
        $\norm{\mb B^0}\lesssim\sigma^2\xi^4\sqrt{\pi d}$.
        Therefore, the variational characterization of the Frobenius norm and
        Lemma~\ref{lemma: lemma 4.4 chen2019} give
        \begin{align*}
            &\bignorm{\mc P_\Omega\big[\psi(\mc M(\bo\theta^0))
            -\psi(\mb M^*)\big]\mb Y^0}\fb\\
            \leq{}&\pi\bignorm{(\mb H^0\circ\mb D^0)\mb Y^0}\fb
            +\bignorm{(\mb B^0\circ\mb D^0)\mb Y^0}\fb\\
            \lesssim{}&\pi\bignorm{(\mb H^0\circ\mb D^0)\mb Y^0}\fb
            +\norm{\mb B^0}\bignorm{\mb Y^0}\fb
            \Big\{
            \norm{\bo\zeta^0-\bo\zeta^*}_\infty
            \\ 
            & \quad +\norm{\mb X^*}\ti
            \bignorm{\mb Y^0\mb O^0-\mb Y^*}\ti
            +\bignorm{\mb X^0\mb O^0-\mb X^*}\ti
            \bignorm{\mb Y^0}\ti        
            \Big\}\\
            \lesssim{}&\kappa_\psi
            \sqrt{\kappa d \pi\sigma_1^*}\,
            \xi^{\sf{gd}}.
        \end{align*}
        Here, the last inequality follows from
        \eqref{eq: GD initialization global}, the incoherence bounds, and the
        stated sampling condition; in particular, the latter absorbs the
        fluctuation term involving $\norm{\mb B^0}$.

        For the noise component, \eqref{eq: response noise spectral concentration}
        and $\bignorm{\mb Y^0}\fb\lesssim\sqrt{r\sigma_1^*}$ give directly
        \begin{align*}
            \bignorm{\mc P_\Omega(\mb E)\mb Y^0}\fb
            \leq \norm{\mc P_\Omega(\mb E)}\bignorm{\mb Y^0}\fb
            \lesssim \sigma\sqrt{\pi d r\sigma_1^*}.
        \end{align*}
        The same argument applies to $\mc P_\Omega(\mb E)\t\mb X^0$. Finally,
        $\sqrt{\sigma_r^*/n}\,\bignorm{\mc P_\Omega(\mb E)\t\mb 1_n}_2
        \leq\sqrt{\sigma_r^*}\norm{\mc P_\Omega(\mb E)}$ and is therefore
        absorbed by the preceding bound. Combining the nonlinear and noise
        bounds proves the claim.
    \end{proof}

Having established these bounds, we now prove Lemma~\ref{lemma: balancedness and orthogonality}.

We invoke Lemma~\ref{lemma: bound for 1t Xt}, \eqref{eq: condition -- F-norm consistency between theta and thetat}, and \eqref{eq: condition -- orthogonality between 1 and Xt} to control $\norm{\mb 1_n\t \mb X^{t+1}}_2$ as follows:
        \begin{align*}
            & \bignorm{\mb 1_n\t \mb X^{t+1}}_2 \leq \big(1-2c_{\perp}\eta\pi\sigma_r^*\big)
        \bignorm{\mb 1_n\t\mb X^t}_2
        +C\eta\kappa^{\frac12}\kappa_\psi \sqrt{\pi d n \sigma_1^*}
        \big(\xi^{\sf{gd}}+\xi^{\sf{gd},\sf{loo}}\big)\\
    &\quad +C\eta\sigma \kappa \kappa_\psi \sqrt{ nd r \pi \sigma_1^*\log d}\\
    & \leq  (1 - 2 c_{\perp} \eta \pi \sigma_r^*) \frac{\sigma  \kappa^{\frac32} \kappa_\psi \sqrt{nd r \log d} + \sqrt{n} \kappa\kappa_{\psi} \sqrt{d} \big(\xi^{\sf{gd}}+\xi^{\sf{gd},\sf{loo}}\big)}{c_{\perp} \sqrt{\pi \sigma_r^*}}\\
    &\quad + C\eta\kappa^{\frac12}\kappa_\psi \sqrt{\pi d n \sigma_1^*} \big(\xi^{\sf{gd}}+\xi^{\sf{gd},\sf{loo}}\big)
        +C\eta\sigma \kappa \kappa_\psi \sqrt{ nd r \pi \sigma_1^*\log d} \\ 
    & \lesssim    \sqrt{n}\frac{\sigma  \kappa^{\frac32} \kappa_\psi \sqrt{ d r \log d} + \kappa\kappa_{\psi} \sqrt{d} \big(\xi^{\sf{gd}}+\xi^{\sf{gd},\sf{loo}}\big)}{c_{\perp} \sqrt{\pi \sigma_r^*}}
        \end{align*}
        with probability at least $1 -O(d^{-c-1})$.

        On the other hand, to arrive at the conclusion that, with probability at least $1-O(d^{-c-1})$,
        \begin{equation}
        \bignorm{{\mb X^{t+1}}\t \mb X^{t+1} - {\mb Y^{t+1}}\t \mb Y^{t+1}}\fb \leq \frac{t+1}{t_0}\sqrt{\sigma_1^*}\frac{\sigma  \kappa^{\frac32} \kappa_\psi \sqrt{ d r \log d} + \kappa\kappa_{\psi} \sqrt{d} \big(\xi^{\sf{gd}}+\xi^{\sf{gd},\sf{loo}}\big)}{c_{\perp} \sqrt{\pi \sigma_r^*}}
        \end{equation}
        holds for the $(t+1)$-th step, it suffices to establish the following in view of Lemma~\ref{lemma: bound for Xt+1t Xt+1 - Yt+1t Yt+1}: 
        \begin{align*}
    & t_0 \underline c_\psi^2 \eta^2\pi^2 \sigma_1^*\bignorm{{\mb X^t}\t \mb X^t - {\mb Y^t}\t \mb Y^t}\fb \ll  1,\\
    & \max\Big\{ \frac{\eta c_{\perp} \pi \sigma_r^*}{n} \bignorm{\mb 1_n\t \mb X^t}_2^2,~ \eta^2 \norm{\nabla L^{\sf{aug}}(\bo \theta^t)}\fb^2
            \Big\}\\
    &\quad \lesssim  \frac{1}{t_0}\sqrt{\sigma_1^*}\frac{\sigma  \kappa^{\frac32} \kappa_\psi \sqrt{ d r \log d} + \kappa\kappa_{\psi} \sqrt{d} \big(\xi^{\sf{gd}}+\xi^{\sf{gd},\sf{loo}}\big)}{c_{\perp} \sqrt{\pi \sigma_r^*}}.
        \end{align*}
     These conditions are verified as follows. 
     
     First, we have from \eqref{eq: condition -- balancedness between X and Y} that 
     \begin{align}
        &  t_0 \underline c_\psi^2 \eta^2\pi^2 \sigma_1^*\bignorm{{\mb X^t}\t \mb X^t - {\mb Y^t}\t \mb Y^t}\fb \\
    & \lesssim    \kappa^{\frac12} \log d \frac{\sigma  \kappa^{\frac32} \kappa_\psi \sqrt{ d r \log d} + \kappa \kappa_\psi \sqrt{d} \big(\xi^{\sf{gd}}+\xi^{\sf{gd},\sf{loo}}\big) }{c_{\perp} \sigma_r^*}\\
    &\quad \lesssim   \kappa^{\frac12} \frac{\sigma  \kappa^{\frac32} \kappa_\psi \sqrt{ d r \log d} + \kappa\kappa_{\psi} \sqrt{d} \big(\xi^{\sf{gd}}+\xi^{\sf{gd},\sf{loo}}\big)}{\underline c_\psi \sqrt{\pi }\sigma_r^*}
        \ll 1, 
        \label{eq: upper bound for balance term in the induction}
     \end{align}
     due to the conditions $\eta t_0 \pi \underline c_\psi \sigma_r^* \asymp \log d$, $ \eta c_{\perp} \pi \sigma_1^* \ll 1$, $\sigma \kappa^3 \kappa_\psi \sqrt{dr \log d / \pi} \ll \underline c_\psi \sigma_r^*$, and $\kappa^{\frac32} \kappa_\psi \sqrt{d}\big(\xi^{\sf{gd}}+\xi^{\sf{gd},\sf{loo}}\big) \ll \underline c_\psi \sqrt{\pi} \sigma_r^*$. 
     Moreover,  
        noticing the condition from \eqref{eq: condition -- orthogonality between 1 and Xt}, one has 
        \begin{align}
            & \frac{\eta c_{\perp} \pi \sigma_r^*}{n} \bignorm{\mb 1_n\t \mb X^t}_2^2 \lesssim \Big[ \eta c_\perp \pi \sigma_r^*\frac{\sigma  \kappa^{\frac32} \kappa_\psi \sqrt{ d r \log d} + \kappa\kappa_{\psi} \sqrt{d} \big(\xi^{\sf{gd}}+\xi^{\sf{gd},\sf{loo}}\big)}{c_{\perp} \sqrt{\pi \sigma_r^*}}\Big] \\ 
    & \qquad \cdot  \Big[ \frac{\sigma  \kappa^{\frac32} \kappa_\psi \sqrt{ d r \log d} + \kappa\kappa_{\psi} \sqrt{d} \big(\xi^{\sf{gd}}+\xi^{\sf{gd},\sf{loo}}\big)}{c_{\perp} \sqrt{\pi \sigma_r^*}} \Big]\\
    &\quad \ll \frac{1}{t_0}\sqrt{\sigma_1^*} \frac{\sigma  \kappa^{\frac32} \kappa_\psi \sqrt{ d r \log d} + \kappa \kappa_\psi \sqrt{d} \big(\xi^{\sf{gd}}+\xi^{\sf{gd},\sf{loo}}\big) }{c_{\perp} \sqrt{\pi \sigma_r^*}},
        \end{align}
        {\emergencystretch=1.5em
since $\eta t_0 \pi \underline c_\psi \sigma_r^* \asymp \log d$, $\sigma \kappa \kappa_\psi \sqrt{d r \log d / \pi}(\log d)^{\frac32} \ll \sigma_r^*$, and $\kappa \kappa_\psi \sqrt{d}\log d \big(\xi^{\sf{gd}}+\xi^{\sf{gd},\sf{loo}}\big)  \ll \underline c_\psi \sqrt{\pi} \sigma_r^*$. Moreover, it is seen from Lemmas~\ref{lemma: bound on nabla Laug}~and~\ref{lemma: gradient norm at theta0} that\par}
 
        \begin{align*}
	            & \eta^2  \norm{\nabla L^{\sf{aug}}(\bo \theta^t)}\fb^2 \lesssim \eta^2 \norm{\nabla L^{\sf{aug}}(\bo \theta^0)}_2^2 +  \eta^2 c_\perp^2 \kappa_\psi^2 \kappa^2 \pi^2 \sigma_1^*\max_{0\leq s\leq t}\bignorm{{\mb X^s}\t \mb X^s  - {\mb Y^s}\t \mb Y^s}\fb^2  \\ 
    & \lesssim   \eta^2 \big[\sigma^2 d \pi r\sigma_1^* + \kappa_\psi^2  \kappa d \pi \sigma_1^* \big(\xi^{\sf{gd}}+\xi^{\sf{gd},\sf{loo}}\big)^2 \big]\\
    &\quad + \eta^2 c_\perp^2 \kappa_\psi^2 \kappa^2 \pi^2 \sigma_1^* \max_{0\leq s\leq t}\bignorm{{\mb X^s}\t \mb X^s  - {\mb Y^s}\t \mb Y^s}\fb^2\\
    & \lesssim   \Big[\frac{1}{t_0} \sqrt{\sigma_1^*} \frac{\sigma  \kappa^{\frac32} \kappa_\psi \sqrt{ d r \log d} + \kappa \kappa_\psi \sqrt{d} \big(\xi^{\sf{gd}}+\xi^{\sf{gd},\sf{loo}}\big) }{c_{\perp} \sqrt{\pi \sigma_r^*}}\Big]\\
    &\quad \cdot \underbrace{\Big[\frac{\sigma r\sqrt{ d \log d / \pi }}{\underline c_\psi \sigma_r^* }  + \frac{\kappa_\psi  \sqrt{d}\big(\xi^{\sf{gd}}+\xi^{\sf{gd},\sf{loo}}\big) \log d }{\underline c_\psi \sqrt{\pi} \sigma_r^* } \Big] }_{\ll 1}\\
	            & \quad + \eta^2 c_\perp^2 \kappa_\psi^2 \kappa^2 \pi^2 \sigma_1^*\max_{0\leq s\leq t}\bignorm{{\mb X^s}\t \mb X^s  - {\mb Y^s}\t \mb Y^s}\fb^2 . 
        \end{align*}
        where we used the conditions $\frac{\sigma r\sqrt{d\log d/\pi}}{\underline c_\psi\sigma_r^*}\ll1$ and $\frac{\kappa_\psi\sqrt{d}\log d\big(\xi^{\sf{gd}}+\xi^{\sf{gd},\sf{loo}}\big)}{\underline c_\psi\sqrt{\pi}\sigma_r^*}\ll1$. Thus, it remains to show that 
        \begin{align}
    & \eta^2 c_\perp^2 \kappa_\psi^2 \kappa^2 \pi^2 \sigma_1^*\max_{0\leq s\leq t}\bignorm{{\mb X^s}\t \mb X^s  - {\mb Y^s}\t \mb Y^s}\fb^2\\
    &\quad \ll \frac{1}{t_0}\sqrt{\sigma_1^*} \frac{\sigma  \kappa^{\frac32} \kappa_\psi \sqrt{ d r \log d} + \kappa \kappa_\psi \sqrt{d} \big(\xi^{\sf{gd}}+\xi^{\sf{gd},\sf{loo}}\big) }{c_{\perp} \sqrt{\pi \sigma_r^*}}.
        \end{align}
        Recalling the relations $\eta c_\perp \pi \sigma_1^* \ll 1$, $t_0 \eta \underline c_\psi \pi \sigma_r^* \asymp \log d$ and the inductive hypothesis \eqref{eq: condition -- balancedness between X and Y}, we have 
        \begin{align*}
	            & \eta^2 c_\perp^2 \kappa_\psi^2 \kappa^2 \pi^2 \sigma_1^* \max_{0\leq s\leq t}\bignorm{{\mb X^s}\t \mb X^s  - {\mb Y^s}\t \mb Y^s}\fb^2  \\ 
    & \lesssim   \Big[t_0 \eta \kappa_\psi^2 \kappa^{\frac52} \pi \sigma_r^* \frac{\sigma  \kappa^{\frac32} \kappa_\psi \sqrt{ d r \log d} +\kappa \kappa_\psi \sqrt{d} \big(\xi^{\sf{gd}}+\xi^{\sf{gd},\sf{loo}}\big)}{\sqrt{\pi} \sigma_r^*} \Big]\\
    &\qquad \Big[ \frac{1}{t_0}\sqrt{\sigma_1^*} \frac{\sigma  \kappa^{\frac32} \kappa_\psi \sqrt{ d r \log d} + \kappa \kappa_\psi \sqrt{d} \big(\xi^{\sf{gd}}+\xi^{\sf{gd},\sf{loo}}\big) }{c_{\perp} \sqrt{\pi \sigma_r^*}} \Big]\\
    & \lesssim    \Big[ \kappa_\psi^2 \kappa^{\frac52} \log d \frac{\sigma  \kappa^{\frac32} \kappa_\psi \sqrt{ d r \log d} + \kappa \kappa_\psi \sqrt{d} \big(\xi^{\sf{gd}}+\xi^{\sf{gd},\sf{loo}}\big)}{\underline c_\psi \sqrt{\pi} \sigma_r^*} \Big]\\
    &\qquad \Big[\frac{1}{t_0}\sqrt{\sigma_1^*} \frac{\sigma  \kappa^{\frac32} \kappa_\psi \sqrt{ d r \log d} +\kappa \kappa_\psi \sqrt{d} \big(\xi^{\sf{gd}}+\xi^{\sf{gd},\sf{loo}}\big) }{c_{\perp} \sqrt{\pi \sigma_r^*}}\Big]\\
    & \ll   \frac{1}{t_0}\sqrt{\sigma_1^*} \frac{\sigma  \kappa^{\frac32} \kappa_\psi \sqrt{ d r \log d} +\kappa \kappa_\psi \sqrt{d} \big(\xi^{\sf{gd}}+\xi^{\sf{gd},\sf{loo}}\big) }{c_{\perp} \sqrt{\pi \sigma_r^*}},
        \end{align*}
        {\emergencystretch=1.5em
where the last inequality follows from the conditions $ \sigma\kappa_\psi^3 \kappa^4 \allowbreak \sqrt{dr}(\log d)^{\frac32} \vee\allowbreak  \kappa_\psi^3 \kappa^{\frac72} \log d \allowbreak \sqrt{d}\big(\xi^{\sf{gd}}+\xi^{\sf{gd},\sf{loo}}\big) \allowbreak \ll \underline c_\psi \allowbreak \sqrt{\pi} \sigma_r^*$.\par}

        The same algebra applies to each leave-one-out objective
        $L^{(l)}$. Replacing one row or column by its population counterpart
        cannot increase the concentration bounds used above, and the
        initialization is exactly centered and balanced for every $l$.
        For each fixed $l$, the conditional concentration argument has failure
        probability $O(d^{-c-2})$. Taking a union bound over
        $l\in[n]\cup[-p]$ proves the two uniform leave-one-out bounds with
        failure probability $O(d^{-c-1})$.
\qed

    \subsection{Analysis of the Auxiliary Sequence under the Augmented Objective} 
    In addition to analyzing the actual GD iterates and their leave-one-out counterparts, we make a detour to examine the vanilla gradient descent sequence under the augmented objective $L^{\sf{aug}}$, as outlined in Section~\ref{subsubsec: key ingredients in the proof}. 
    Denote this sequence by $\bo \theta^{\sf{aug}, t} \coloneqq [\bo \zeta^{\sf{aug}, t}, \mb X^{\sf{aug},t}, \mb Y^{\sf{aug},t}]$. It plays a crucial role in establishing the approximate stationary condition in Section~\ref{subsec: proof of inference results}. Formally, it is initialized at $\bo \theta^{\sf{aug}, 0} = \bo \theta^0$ and updated according to $\bo \theta^{\sf{aug},t} = \bo \theta^{\sf{aug},t-1} - \eta \nabla_{\bo \theta} L^{\sf{aug}}(\bo \theta^{\sf{aug},t-1})$. Correspondingly, we introduce its leave-one-out variant $\{\bo \theta^{\sf{aug},t,(l)}\}$ with $l\in[n] \cup [-p]$ by $\bo \theta^{\sf{aug},t,(l)} = \bo \theta^{\sf{aug},t-1,(l)} - \eta \nabla_{\bo \theta} L^{(l)}(\bo \theta^{\sf{aug},t-1,(l)}) - \eta \nabla_{\bo \theta} L^{\sf{diff}}(\bo \theta^{\sf{aug},t-1,(l)})$ with $L^{(l)}$ defined in \eqref{eq: leave-one-out objective function} and $\bo \theta^{\sf{aug},0,(l)} \coloneqq \bo \theta^{0,(l)}$.
    For its comparison with the actual sequence, let
    \begin{align*}
        \mb F^{\sf{aug},t}
        &\coloneqq
        \begin{pmatrix}
            \mb X^{\sf{aug},t}\\
            \mb Y^{\sf{aug},t}
        \end{pmatrix},
        \quad 
        \mb O_{\sf{aug}}^t
        \coloneqq
        \mb O(\bo\theta^{\sf{aug},t},\bo\theta^t)
        \in\argmin_{\mb O\in\mc O(r)}
        \bignorm{\mb F^{\sf{aug},t}\mb O-\mb F^t}\fb,
        \\
        \bo\Delta_{\sf{aug}}^t
        &\coloneqq
        [\bo\zeta^{\sf{aug},t}-\bo\zeta^t,
        \mb X^{\sf{aug},t}\mb O_{\sf{aug}}^t-\mb X^t,
        \mb Y^{\sf{aug},t}\mb O_{\sf{aug}}^t-\mb Y^t].
    \end{align*}
    Thus $\norm{\bo\Delta_{\sf{aug}}^t}_2
    =\distf(\bo\theta^{\sf{aug},t},\bo\theta^t)$.
    For simplicity, introduce the static comparison radius
    \begin{align}
        \delta^{\sf{aug}}
        \coloneqq{}&
        \Big\{\frac{\kappa\kappa_\psi\underline c_\psi}{c_\perp}
        +\frac{1}{\sqrt{\pi(n\wedge p)}}\Big\}
        \frac{\sigma\kappa\kappa_\psi
        \sqrt{\kappa\mu r d\log d}}
        {\underline c_\psi\sqrt{\pi\sigma_r^*}}
        +\frac{\kappa^2\kappa_\psi^2\underline c_\psi}{c_\perp}
        \sqrt d\,(\delta_{\sf{row}}^0+\delta_{\sf{loo}}^0),
        \label{eq: definition of delta^aug}
    \end{align}
    and the decaying initialization radius
    \begin{align}
        \tau_{\sf aug}^t
        \coloneqq
        \frac{\nu^t}{\underline c_\psi\sqrt{\sigma_r^*}}
        \Big\{\frac{\xi^{\sf gd}+\xi^{\sf gd,\sf loo}}{\sqrt{\mu r}}
        +\kappa\kappa_\psi\sqrt{\rho}\,\xi^{\sf gd}\Big\},
        \qquad 0\leq t\leq t_0.
        \label{eq: definition of auxiliary initialization transient}
    \end{align}
    The sampling condition in Lemma~\ref{lemma: leave-one-out error} gives
    \begin{align}
        \frac{\sigma^2}{\underline c_{\sf{Hess}}}\kappa\mu r
        \sqrt{\frac{\rho\log d}{\pi(n\wedge p)}}\ll1.
        \label{eq: auxiliary LOO sampling ratio}
    \end{align}
    In particular, the factor $\{\pi(n\wedge p)\}^{-1/2}$ in
    \eqref{eq: definition of delta^aug} multiplies only the stochastic term.
    \begin{lemma}
    \label{lemma: auxiliary sequence analysis}
        Under the assumptions of Theorem~\ref{thm: general GD}, there is a
        sufficiently large numerical constant $D_2$ such that, with
        probability at least $1-O(d^{-c})$, the following holds for every
        $0\leq t\leq t_0$:
        \begin{align}
            & \distf(\bo \theta^{\sf{aug},t},\bo \theta^t) \lesssim \delta^{\sf{aug}},\quad \dist\ti(\bo \theta^{\sf{aug},t} , \bo \theta^t)  \leq D_2 \kappa\kappa_\psi \sqrt{\frac{\mu \rho r}{d}} (\delta^{\sf{aug}}+\tau_{\sf aug}^t) \log d, \\
            & \norm{\mb 1_n\t \big(\mb X^{\sf{aug},t}\mb O_{\sf{aug}}^t - \mb X^t \big)}_2 \lesssim  \kappa\sqrt{n \mu r} \frac{\sigma^2}{c_{\perp}}  (\delta^{\sf{aug}}+\tau_{\sf aug}^t) \log d.
        \end{align}
        Moreover, $\tau_{\sf aug}^{t_0}\lesssim\delta^{\sf{aug}}$, so the
        rowwise and centering bounds at $t=t_0$ hold with
        $\delta^{\sf{aug}}$ in place of
        $\delta^{\sf{aug}}+\tau_{\sf aug}^{t_0}$.
    \end{lemma}
    \begin{proof}
        Following the probability bookkeeping above
        Condition~\ref{condition: inductive contraction for GD}, invoke each
        fixed-iteration, fixed-index concentration bound with failure probability
        $O(d^{-c-12})$. Since $n+p\leq2d$ and $t_0\leq d^{10}$, a union bound,
        together with the uniform RIC event and the $q=1$ conclusion of
        Theorem~\ref{thm: general GD}, gives a joint event of probability at
        least $1-O(d^{-c})$ on which all estimates below hold simultaneously.

        We justify the closeness between $\bo \theta^{\sf{aug},t}$ and $\bo \theta^t$ by induction. Assume that the following four hypotheses hold at time $t$:
        \begin{align}
            & \distf(\bo \theta^{\sf{aug},t},\bo \theta^t) \leq D_1 \delta^{\sf{aug}},
            \label{eq: auxiliary induction Frobenius}\\
            & \dist\ti(\bo \theta^{\sf{aug},t} , \bo \theta^t)  \leq  D_2 \kappa \kappa_\psi \sqrt{\frac{\mu \rho r}{d}} (\delta^{\sf{aug}}+\tau_{\sf aug}^t) \log d,
            \label{eq: auxiliary induction rowwise}\\
            & \norm{\mb 1_n\t \big(\mb X^{\sf{aug},t} \mb O^t_{\sf{aug}} - \mb X^t \big)}_2 \leq D_3\kappa\sqrt{n \mu r} \frac{\sigma^2}{c_{\perp}}  (\delta^{\sf{aug}}+\tau_{\sf aug}^t) \log d,
            \label{eq: auxiliary induction centering}\\
            & \max_{l\in[n]\cup[-p]}
            \distf(\bo\theta^{\sf aug,t},\bo\theta^{\sf aug,t,(l)})\leq D_4\Bigg[\nu^t\Big\{\kappa\kappa_\psi\sqrt{\mu\rho r}\,\delta_{\sf{row}}^0+\delta_{\sf{loo}}^0\Big\}\nonumber\\
            &\qquad+\frac{\sigma\kappa\kappa_\psi\sqrt{\kappa\mu\rho r^2\log d}}
            {\underline c_\psi\sqrt{\pi\sigma_r^*}}+\frac{\sigma^2}{\underline c_{\sf{Hess}}}
            \kappa^2\kappa_\psi\mu r\sqrt{\frac{\rho\log d}{\pi(n\wedge p)}}
            \sqrt{\frac{\mu\rho r}{d}}(\delta^{\sf aug}+\tau_{\sf aug}^t)\log d\Bigg].
            \label{eq: auxiliary induction LOO}
        \end{align}
        
        Next, fix an arbitrary $t \in \{0,\ldots,t_0-1\}$. Rotational
        equivariance allows us to rotate the auxiliary update by
        $\mb O_{\sf{aug}}^t$, while the optimality of
        $\mb O_{\sf{aug}}^{t+1}$ makes the resulting iterate a valid comparison
        candidate at step $t+1$. Moreover, the factor block of
        $\bo\Delta_{\sf{aug}}^t$ is Procrustes aligned and hence belongs to the
        direction set in \eqref{eq: direction set}. The squared-norm argument
        proving \eqref{eq: one-step analysis under F norm}, together with
        \eqref{eq: condition -- balancedness between X and Y}, therefore gives
        \begin{align*}
            & \distf(\bo \theta^{\sf{aug},t+1},\bo \theta^{t+1})\\
            \leq & \nu \distf(\bo \theta^{\sf{aug},t},\bo \theta^t) + \eta \bignorm{\nabla L^{\sf{diff}}(\bo \theta^t)}_2 \\ 
            \leq & \nu\distf(\bo \theta^{\sf{aug},t},\bo \theta^t) + C \eta \underline c_\psi \pi \sigma_1^*  \frac{\sigma \kappa^{\frac32} \kappa_\psi \sqrt{d r\log d} + \kappa \kappa_\psi \sqrt{d} (\xi^{\sf{gd}}+\xi^{\sf{gd},\sf{loo}} )}{c_{\perp} \sqrt{\pi \sigma_r^*}} \\
            \leq & \nu D_1\delta^{\sf{aug}}+(1-\nu)D_1\delta^{\sf{aug}}
            =D_1\delta^{\sf{aug}}.
        \end{align*}
        Here, the last inequality follows from the definition of
        $\delta^{\sf{aug}}$, $1-\nu=\eta\underline c_{\sf{Hess}}
        \pi\sigma_r^*/4$, and $\underline c_{\sf{Hess}}\asymp
        \underline c_\psi$: for a sufficiently large constant $D_1$, the
        one-step forcing term is at most $(1-\nu)D_1\delta^{\sf{aug}}$.
        Thus \eqref{eq: auxiliary induction Frobenius} holds at time $t+1$ on
        the joint event.

        \bigskip

        Turning our attention to the row-wise error for the $(t+1)$-th step,
        rotational equivariance allows us to right-multiply the auxiliary
        update by $\mb O_{\sf{aug}}^t$. The corresponding linear predictors,
        residuals, row norms, and Gram-matrix discrepancy are unchanged. We
        take the $i$-th row of $\mb X$ as an example and obtain
        \begin{align*}
            & \bignorm{\mb X^{\sf{aug},t+1}_i\mb O_{\sf{aug}}^t - \mb X^{t+1}_i}_2 \\ 
    & \leq   \Bignorm{\mb X^{\sf{aug},t}_i\mb O_{\sf{aug}}^t - \mb X^{t}_i - \eta \big[\mc P_\Omega\big(\psi(\mc M(\bo \theta^{\sf{aug},t})) - \mb R\big)\big]_{i,\cdot}\mb Y^{\sf{aug},t}\mb O_{\sf{aug}}^t\\
    &\qquad + \eta \big[\mc P_\Omega\big(\psi(\mc M(\bo \theta^{t})) - \mb R\big)\big]_{i,\cdot}\mb Y^{t}}\\
            & + 2 \eta c_{\perp}\frac{\pi \sigma_r^*}{n}\bignorm{\mb 1_n\t \big(\mb X^{t} - \mb X^{\sf{aug},t}\mb O_{\sf{aug}}^t\big)}_2  \\ 
            & + 4 \eta c_{\sf{aug}} \pi
            \norm{\mb X^{\sf{aug},t}\mb O_{\sf{aug}}^t}\ti \\
            &\qquad \cdot \norm{(\mb X^{\sf{aug},t}\mb O_{\sf{aug}}^t)\t
            (\mb X^{\sf{aug},t}\mb O_{\sf{aug}}^t)
            -(\mb Y^{\sf{aug},t}\mb O_{\sf{aug}}^t)\t
            (\mb Y^{\sf{aug},t}\mb O_{\sf{aug}}^t)}\fb  \\ 
    & \leq   \Bignorm{\mb X^{\sf{aug},t}_{i,\cdot}\mb O_{\sf{aug}}^t - \mb X^{t}_{i,\cdot}\\
    &\qquad - \eta \big[\big( \psi(\mc M([\bo\zeta^t,\mb X^{\sf{aug},t}\mb O_{\sf{aug}}^t,\mb Y^t])) - \psi(\mc M(\bo\theta^t)) \big) \circ \bo \Omega \big]_{i,\cdot} \mb Y^t}\\
    & + \eta \Bignorm{\big[\big(\psi(\mc M([\bo\zeta^{\sf{aug},t},
            \mb X^{\sf{aug},t}\mb O_{\sf{aug}}^t,
            \mb Y^{\sf{aug},t}\mb O_{\sf{aug}}^t]))\\
    &\qquad - \psi(\mc M([\bo\zeta^t,
            \mb X^{\sf{aug},t}\mb O_{\sf{aug}}^t,\mb Y^t])) \big)
            \circ \bo \Omega \big]_{i,\cdot} \mb Y^t}_2 \\ 
    & + \eta \norm{\big[\mc P_\Omega\big(\psi(\mc M(\bo \theta^{\sf{aug},t})) - \mb R\big)\big]_{i,\cdot} (\mb Y^{\sf{aug},t}\mb O_{\sf{aug}}^t - \mb Y^t)}\\
    &\quad + 2 \eta c_{\perp}\frac{\pi \sigma_r^*}{n}\bignorm{\mb 1_n\t\big(\mb X^{\sf{aug},t}\mb O_{\sf{aug}}^t- \mb X^{t}\big)}_2\\
            & + 4 \eta c_{\sf{aug}} \pi
            \norm{\mb X^{\sf{aug},t}\mb O_{\sf{aug}}^t}\ti\\
    &\qquad \cdot\norm{(\mb X^{\sf{aug},t}\mb O_{\sf{aug}}^t)\t
            (\mb X^{\sf{aug},t}\mb O_{\sf{aug}}^t)
            -(\mb Y^{\sf{aug},t}\mb O_{\sf{aug}}^t)\t
            (\mb Y^{\sf{aug},t}\mb O_{\sf{aug}}^t)}\fb
            \label{eq: auxiliary sequence decomposition}
        \end{align*}

        We justify the bounds for the terms in the above decomposition as follows.
        \begin{itemize}
        \item To begin with, the mean-value theorem and the leave-one-out argument used in \eqref{eq: upper bound and lower bound on singular values of Hessian wrt X} produce an $r\times r$ positive-semidefinite Hessian whose eigenvalues lie between constant multiples of $\underline c_\psi\pi\sigma_r^*$ and $\sigma^2\pi\sigma_1^*$. Squaring the norm and using the corresponding multiplication inequality, as in \eqref{eq: held-out row contraction}, yields on the joint event that
        \begin{align}
            &  \norm{\mb X^{\sf{aug},t}_{i,\cdot}\mb O_{\sf{aug}}^t - \mb X^{t}_{i,\cdot} - \eta \big[\big( \psi(\mc M([\bo\zeta^t,\mb X^{\sf{aug},t}\mb O_{\sf{aug}}^t,\mb Y^t])) - \psi(\mc M(\bo\theta^t)) \big) \circ \bo \Omega \big]_{i,\cdot} \mb Y^t}   \\ 
            \leq & \Big(1 - \frac12\eta \underline c_{\sf{Hess}} \pi \sigma_r^*\Big) \bignorm{\mb X^{\sf{aug},t}_{i}\mb O_{\sf{aug}}^t - \mb X^{t}_{i}}_2, 
        \end{align} 
        
        \item To use the inductive bounds from the preceding GD
        analysis without replacing the decaying initialization terms by a
        time-uniform bound, observe that
        \begin{align}
            &\max_{l\in[n]\cup[-p]}\distf(\bo\theta^t,\bo\theta^{t,(l)})\lesssim\nu^t\Big\{\kappa\kappa_\psi\sqrt{\mu\rho r}\,\delta_{\sf{row}}^0+\delta_{\sf{loo}}^0\Big\}
            +\frac{\sigma\kappa\kappa_\psi\sqrt{\kappa\mu\rho r^2\log d}}
            {\underline c_\psi\sqrt{\pi\sigma_r^*}}\nonumber\\
            &\quad\lesssim\sqrt{\frac{\mu r}{\pi}}
            (\delta^{\sf aug}+\tau_{\sf aug}^t).
            \label{eq: explicit time-dependent actual LOO bound}
        \end{align}
        The last inequality follows from
        \eqref{eq: definition of delta^aug}--\eqref{eq: definition of auxiliary initialization transient},
        $d=\rho(n\wedge p)$, $\pi\leq1$, and $\mu\geq1$.
        To propagate the augmented leave-one-out coupling jointly with the
        rowwise comparison, set
        \begin{align*}
            \mathfrak A_t
            &\coloneqq\max_{l\in[n]\cup[-p]}
            \distf(\bo\theta^{\sf aug,t},\bo\theta^{\sf aug,t,(l)}),\qquad
            \mathfrak R_t
            \coloneqq\kappa\kappa_\psi\sqrt{\frac{\mu\rho r}{d}}
            (\delta^{\sf aug}+\tau_{\sf aug}^t)\log d.
        \end{align*}
        A standard chaining argument verifies the RIC conditions:
        \eqref{eq: auxiliary induction Frobenius}--\eqref{eq: auxiliary induction rowwise}
        and the actual-GD bounds control the augmented full iterate;
        \eqref{eq: auxiliary induction LOO} and the triangle inequality control
        every augmented leave-one-out iterate; and convexity controls the
        aligned interpolation segments.
        The required slack follows from the initialization, signal-strength,
        and sampling conditions.

        We now repeat the decomposition in the proof of
        Lemma~\ref{lemma: leave-one-out error}, comparing the augmented full
        update under $L^{\sf aug}$ with the augmented leave-one-out update
        under $L^{(l)}+L^{\sf diff}$. The common balancing gradient is included
        in the contracting $L^{\sf aug}$ term, so the only forcing is the
        held-out score. The local Hessian bound and conditional Bernstein
        inequality, uniformly over $l\in[n]\cup[-p]$, give
        \begin{align}
            \mathfrak A_{t+1}
            \leq{}&\Big(1-\frac{\eta\underline c_{\sf{Hess}}
            \pi\sigma_r^*}{2}\Big)\mathfrak A_t
            +C\eta\sigma\sqrt{\pi\mu\rho r^2\sigma_1^*\log d}\nonumber\\
            &+C\eta\sigma^2\mu r\sigma_1^*
            \sqrt{\frac{\pi\rho\log d}{n\wedge p}}
            (\mathfrak A_t+\mathfrak R_t)\nonumber\\
            &+C\eta\sigma^2\mu r\sigma_1^*
            \sqrt{\frac{\pi\rho\log d}{n\wedge p}}\,
            \nu^t\Big\{\kappa\kappa_\psi\sqrt{\mu\rho r}\,\delta_{\sf{row}}^0+\delta_{\sf{loo}}^0\Big\}\nonumber\\
            &+C\eta\sigma^2\mu r\sigma_1^*
            \sqrt{\frac{\pi\rho\log d}{n\wedge p}}
            \frac{\sigma\kappa\kappa_\psi\sqrt{\kappa\mu\rho r^2\log d}}
            {\underline c_\psi\sqrt{\pi\sigma_r^*}}.
            \label{eq: augmented LOO one-step recursion}
        \end{align}
        Substituting the inductive bounds into
        \eqref{eq: augmented LOO one-step recursion}, the stochastic forcing
        is $O(1-\nu)$ times the stochastic floor, while
        \eqref{eq: auxiliary LOO sampling ratio} absorbs the remaining
        predictor error into the local curvature. Since
        $\tau_{\sf aug}^{t+1}=\nu\tau_{\sf aug}^t$,
        \eqref{eq: auxiliary induction LOO} holds at time $t+1$ for
        sufficiently large $D_4$:
        \begin{align}
            \mathfrak A_{t+1}
            \leq{}&D_4\Bigg[\nu^{t+1}\Big\{\kappa\kappa_\psi\sqrt{\mu\rho r}\,\delta_{\sf{row}}^0+\delta_{\sf{loo}}^0\Big\}
            +\frac{\sigma\kappa\kappa_\psi\sqrt{\kappa\mu\rho r^2\log d}}
            {\underline c_\psi\sqrt{\pi\sigma_r^*}}\nonumber\\
            &\quad+\frac{\sigma^2}{\underline c_{\sf{Hess}}}\kappa^2\kappa_\psi\mu r
            \sqrt{\frac{\rho\log d}{\pi(n\wedge p)}}
            \sqrt{\frac{\mu\rho r}{d}}(\delta^{\sf aug}+\tau_{\sf aug}^{t+1})\log d\Bigg].
            \label{eq: time-dependent augmented LOO comparison}
        \end{align}
        The base case follows from
        $\bo\theta^{\sf aug,0}=\bo\theta^0$,
        $\bo\theta^{\sf aug,0,(l)}=\bo\theta^{0,(l)}$, and the initialization
        leave-one-out bound $\mathfrak A_0\lesssim\delta_{\sf{loo}}^0$.
        The inductive bound \eqref{eq: auxiliary induction LOO} and
        \eqref{eq: auxiliary LOO sampling ratio} further imply
        \begin{align*}
            \mathfrak A_t
            \lesssim\sqrt{\frac{\mu r}{\pi}}
            (\delta^{\sf aug}+\tau_{\sf aug}^t).
        \end{align*}
        Condition~\ref{condition: inductive contraction for GD} and
        \eqref{eq: rowwise transfer from actual to LOO} give the explicit
        actual leave-one-out bound in
        \eqref{eq: explicit time-dependent actual LOO bound}, whereas
        \eqref{eq: time-dependent augmented LOO comparison} bounds the
        augmented coupling by the displayed bound above.
        For $l=i$, let $\mb O^{\sf{aug},t,(i)}\coloneqq\mb O(\bo\theta^{\sf{aug},t,(i)},\bo\theta^{\sf{aug},t})$. We keep the relative and global rotations explicit in the leave-one-out factors below.
        Moreover, we have the following decomposition incorporating the leave-one-out quantities: 
        \begin{align}
    & \Bignorm{\big[\big(\psi(\mc M([\bo\zeta^{\sf{aug},t},
            \mb X^{\sf{aug},t}\mb O_{\sf{aug}}^t,
            \mb Y^{\sf{aug},t}\mb O_{\sf{aug}}^t]))\\
    &\qquad - \psi(\mc M([\bo\zeta^t,
            \mb X^{\sf{aug},t}\mb O_{\sf{aug}}^t,\mb Y^t])) \big)
            \circ \bo \Omega \big]_{i,\cdot} \mb Y^t}_2 \\
    & \leq   \sigma^2 \big( \alpha_0 + \alpha_1 + \alpha_2\big) \big( \beta_0 + \beta_1\big)\\
    & \leq   C\sigma^2\sqrt{\frac{\mu \rho r}{d}}\pi\sigma_1^*\Bigg\{\kappa\kappa_\psi\Big[D_1\delta^{\sf{aug}}+\mathfrak A_t
            +\distf(\bo\theta^t,\bo\theta^{t,(i)})\Big]\sqrt{\log d} \\
            &\qquad+\frac{1}{\sqrt\pi}\Big[D_2\kappa \kappa_\psi \sqrt{\frac{\mu \rho r}{d}}(\delta^{\sf{aug}}+\tau_{\sf aug}^t)\log d
            +\mathfrak A_t+\distf(\bo\theta^t,\bo\theta^{t,(i)})\Big]\log d\Bigg\}\\
    & \leq   D_5  \sigma^2 \sqrt{\frac{\mu \rho r}{d}}\pi \sigma_1^* (\delta^{\sf{aug}}+\tau_{\sf aug}^t) \log d\\
    &\qquad + \frac{\underline c_{\sf{Hess}}}{64\underline c_\psi}
            D_2\kappa\kappa_\psi\underline c_\psi
            \sqrt{\frac{\mu\rho r}{d}}\pi\sigma_r^*
            (\delta^{\sf aug}+\tau_{\sf aug}^t)\log d \\
    & \leq   \frac{\underline c_{\sf{Hess}}}{32\underline c_\psi}
            D_2\kappa\kappa_\psi\underline c_\psi
            \sqrt{\frac{\mu\rho r}{d}}\pi\sigma_r^*
            (\delta^{\sf aug}+\tau_{\sf aug}^t)\log d
            \label{eq: decomposition of the surrogate sequence}
        \end{align}
        holds on the same event for some sufficiently large $D_2$,
        where the quantities $\alpha_0$, $\alpha_1$, $\alpha_2$, $\beta_0$, and $\beta_1$ are defined and controlled as follows 
        \begin{align*}
    & \alpha_0  \coloneqq{}  \Big\|\Big[\Big(
            \mc M([\bo\zeta^{\sf{aug},t,(i)},\mb X^{\sf{aug},t}\mb O_{\sf{aug}}^t,
            \mb Y^{\sf{aug},t,(i)}\mb O^{\sf{aug},t,(i)}\mb O_{\sf{aug}}^t]) \\
            &\quad - \mc M([\bo\zeta^{t,(i)},
            \mb X^{\sf{aug},t}\mb O_{\sf{aug}}^t,
            \mb Y^{t,(i)}\mb O^{t,(i)}])\Big) \circ \bo \Omega
            \Big]_{i,\cdot}\Big\|_2 \\ 
    & \leq{} D_6\sqrt{\frac{\mu \rho r}{d}}(\sigma_1^*)^{\frac12}
            \Bigg\{\kappa\kappa_\psi\sqrt\pi
            \Big[D_1\delta^{\sf aug}+\mathfrak A_t
            +\distf(\bo\theta^t,\bo\theta^{t,(i)})\Big]
            \sqrt{\log d} \\
            &\qquad+\Big[D_2\kappa\kappa_\psi
            \sqrt{\frac{\mu\rho r}{d}}
            (\delta^{\sf aug}+\tau_{\sf aug}^t)\log d
            +\mathfrak A_t+\distf(\bo\theta^t,\bo\theta^{t,(i)})\Big]\log d\Bigg\},\\
    & \alpha_1 \coloneqq{}  \Big\|\Big[
            \mc M([\bo\zeta^{\sf{aug},t},\mb X^{\sf{aug},t}\mb O_{\sf{aug}}^t,
            \mb Y^{\sf{aug},t}\mb O_{\sf{aug}}^t]) \\
            &\quad - \mc M([\bo\zeta^{\sf{aug},t,(i)},
            \mb X^{\sf{aug},t}\mb O_{\sf{aug}}^t,
            \mb Y^{\sf{aug},t,(i)}\mb O^{\sf{aug},t,(i)}\mb O_{\sf{aug}}^t])
            \Big]_{i,\cdot}\Big\|_2 \leq  D_6\sqrt{\frac{\mu \rho r}{d}}(\sigma_1^*)^{\frac12}\mathfrak A_t,\\
    & \alpha_2 \coloneqq{}  \Big\|\Big[
            \mc M([\bo\zeta^t,\mb X^{\sf{aug},t}\mb O_{\sf{aug}}^t,\mb Y^t]) \\
            &\quad - \mc M([\bo\zeta^{t,(i)},
            \mb X^{\sf{aug},t}\mb O_{\sf{aug}}^t,
            \mb Y^{t,(i)}\mb O^{t,(i)}])
            \Big]_{i,\cdot}\Big\|_2 \leq D_6\sqrt{\frac{\mu \rho r}{d}}(\sigma_1^*)^{\frac12}\distf(\bo\theta^t,\bo\theta^{t,(i)}),\\
    & \beta_0 \coloneqq \bignorm{\diag(\bo \Omega_{i,\cdot})\mb Y^{t,(i)}\mb O^{t,(i)}} \leq C\sqrt{\pi\sigma_1^*}+\sqrt{\log d}\Big\{\sqrt{\frac{\mu r\sigma_1^*}{p}}+\distf(\bo\theta^t,\bo\theta^{t,(i)})\Big\}\\
    &\quad \lesssim \sqrt{\pi\sigma_1^*},\\
            & \beta_1 \coloneqq \bignorm{\mb Y^t-\mb Y^{t,(i)}\mb O^{t,(i)}} \leq D_6\distf(\bo\theta^t,\bo\theta^{t,(i)}).
        \end{align*}
        To preserve conditional independence in the comparisons involving the
        two leave-one-out sequences, define their pairwise rotation
        \begin{align*}
            \widehat{\mb O}^{t,(i)}
            \coloneqq\mb O(\bo\theta^{\sf aug,t,(i)},\bo\theta^{t,(i)})
        \end{align*}
        and set
        $\mb S^{t,(i)}\coloneqq
        \mb O^{\sf aug,t,(i)}\mb O_{\sf aug}^t
        -\widehat{\mb O}^{t,(i)}\mb O^{t,(i)}$.
        Both leave-one-out iterates and $\widehat{\mb O}^{t,(i)}$ are
        measurable with respect to the held-out sigma-field. The triangle
        inequality for the four aligned factor pairs gives
        \begin{align*}
            \bignorm{\mb F^{\sf aug,t,(i)}\mb S^{t,(i)}}\fb
            &\lesssim \mathfrak A_t+D_1\delta^{\sf aug}
            +\distf(\bo\theta^t,\bo\theta^{t,(i)}),\\
            \bignorm{\mb S^{t,(i)}}
            &\lesssim\frac{\mathfrak A_t+D_1\delta^{\sf aug}
            +\distf(\bo\theta^t,\bo\theta^{t,(i)})}{\sqrt{\sigma_r^*}},
        \end{align*}
        where the second line uses the RIC lower bound on the smallest
        singular value of $\mb F^{\sf aug,t,(i)}$. In particular,
        \begin{align}
            &\mb Y^{\sf aug,t,(i)}\mb O^{\sf aug,t,(i)}\mb O_{\sf aug}^t
            -\mb Y^{t,(i)}\mb O^{t,(i)}=(\mb Y^{\sf aug,t,(i)}\widehat{\mb O}^{t,(i)}
            -\mb Y^{t,(i)})\mb O^{t,(i)}
            +\mb Y^{\sf aug,t,(i)}\mb S^{t,(i)}.
            \label{eq: pairwise LOO rotation decomposition}
        \end{align}
        Matrix Bernstein applies to each factor preceding the final
        rotation or $\mb S^{t,(i)}$ in this decomposition; the right
        multiplication by $\mb O^{t,(i)}$ does not change the norm. The
        synchronization part of a raw-noise product is therefore bounded by
        \begin{align*}
            &\frac{\sigma\sqrt{\pi\sigma_1^*r\log d}}
            {\sqrt{\sigma_r^*}}
            (\mathfrak A_t+D_1\delta^{\sf aug}
            +\distf(\bo\theta^t,\bo\theta^{t,(i)})),
        \end{align*}
        together with the corresponding tail term determined by the uniform
        summand bound. These are
        absorbed by the stated signal-strength and tail conditions. We use
        the same pairwise split in the conditional bounds for $\alpha_0$ and
        $\beta_0'$ below; thus no unattenuated static radius is multiplied by
        the raw norm of the held-out noise row.
        Here, the bound for $\beta_0$ follows by conditioning on the leave-one-out iterate, applying matrix Bernstein, and then invoking \eqref{eq: rowwise transfer from actual to LOO}. The bound for $\beta_1$ follows directly from \eqref{eq: condition -- F-norm consistency between thetat and thetatl}. The last inequality in \eqref{eq: decomposition of the surrogate sequence} then follows from \eqref{eq: epsilon GD initialization requirement}, \eqref{eq: explicit time-dependent actual LOO bound}, \eqref{eq: time-dependent augmented LOO comparison}, and the stated sampling and signal-strength conditions.

        \item Next, write the auxiliary residual relative to the population mean and insert the auxiliary leave-one-out predictor. The Lipschitz continuity of $\psi$ then gives
        \begin{align}
            &  \norm{\big[\mc P_\Omega\big(\psi(\mc M(\bo \theta^{\sf{aug},t})) - \mb R\big)\big]_{i,\cdot} (\mb Y^{\sf{aug},t}\mb O_{\sf{aug}}^t - \mb Y^t)}  \\ 
            \leq & \sigma^2(\alpha_0' + \alpha_1')(\beta_0' + \beta_1' + \beta_2') \\
            & + \norm{(\mb E\circ \bo \Omega)_{i,\cdot} (\mb Y^{\sf{aug},t,(i)}\mb O^{\sf{aug},t,(i)}\mb O_{\sf{aug}}^t - \mb Y^{t,(i)}\mb O^{t,(i)})}_2 \\ 
            & + \norm{(\mb E\circ \bo \Omega)_{i,\cdot}}_2 (\beta_1' + \beta_2') , 
            \label{eq: decomposition of the Y difference term in the auxiliary series analysis}
        \end{align}
        where the quantities above are defined as 
        \begin{align}
            \alpha_0' \coloneqq{}& \Big\|\Big[\Big(
            \mc M([\bo\zeta^{\sf{aug},t},
            \mb X^{\sf{aug},t}\mb O_{\sf{aug}}^t,
            \mb Y^{\sf{aug},t}\mb O_{\sf{aug}}^t]) \\ 
            &\quad - \mc M([\bo\zeta^{\sf{aug},t,(i)},
            \mb X^{\sf{aug},t}\mb O_{\sf{aug}}^t,
            \mb Y^{\sf{aug},t,(i)}\mb O^{\sf{aug},t,(i)}
            \mb O_{\sf{aug}}^t])\Big)\circ\bo\Omega
            \Big]_{i,\cdot}\Big\|_2, \\ 
            \alpha_1' \coloneqq{}& \Big\|\Big[\Big(
            \mc M([\bo\zeta^{\sf{aug},t,(i)},
            \mb X^{\sf{aug},t}\mb O_{\sf{aug}}^t,
            \mb Y^{\sf{aug},t,(i)}\mb O^{\sf{aug},t,(i)}
            \mb O_{\sf{aug}}^t])
            -\mc M(\bo\theta^*)\Big)\circ\bo\Omega
            \Big]_{i,\cdot}\Big\|_2, \\ 
            & \beta_0' \coloneqq \norm{\diag(\bo \Omega_{i,\cdot}) (\mb Y^{\sf{aug},t,(i)}\mb O^{\sf{aug},t,(i)}\mb O_{\sf{aug}}^t - \mb Y^{t,(i)}\mb O^{t,(i)})}, \\ 
            & \beta_1' \coloneqq \bignorm{
            \mb Y^{\sf{aug},t}\mb O_{\sf{aug}}^t
            - \mb Y^{\sf{aug},t,(i)}\mb O^{\sf{aug},t,(i)}\mb O_{\sf{aug}}^t},\\
            & \beta_2' \coloneqq \bignorm{\mb Y^{t,(i)}\mb O^{t,(i)} - \mb Y^{t}}.
        \end{align}
        These decompositions retain the sharp coupling scales: $\alpha_0'$ and
        $\beta_1'$ are controlled by $\mathfrak A_t$, $\beta_2'$ by
        $\distf(\bo\theta^t,\bo\theta^{t,(i)})$, and $\alpha_1'$ by
        $\mathfrak R_t+\mathfrak A_t+
        \distf(\bo\theta^t,\bo\theta^{t,(i)})$. The leave-one-out pair and its
        pairwise rotation are independent of the held-out row, while
        \begin{align*}
            &\distf(\bo\theta^{\sf aug,t,(i)},\bo\theta^{t,(i)})
            \lesssim D_1\delta^{\sf aug}
            +\mathfrak A_t+\distf(\bo\theta^t,\bo\theta^{t,(i)}).
        \end{align*}
        Thus \eqref{eq: pairwise LOO rotation decomposition} confines this
        enlarged bound to the conditional Bernstein term containing
        $\beta_0'$. The unconditioned factor
        $\norm{(\mb E\circ\bo\Omega)_{i,\cdot}}_2$ multiplies only
        $\beta_1'+\beta_2'\lesssim\mathfrak A_t
        +\distf(\bo\theta^t,\bo\theta^{t,(i)})$. By
        \eqref{eq: explicit time-dependent actual LOO bound} and
        \eqref{eq: time-dependent augmented LOO comparison},
        \begin{align*}
            \mathfrak A_t+\distf(\bo\theta^t,\bo\theta^{t,(i)})
            \lesssim{}&\nu^t\Big\{\kappa\kappa_\psi\sqrt{\mu\rho r}\,\delta_{\sf{row}}^0+\delta_{\sf{loo}}^0\Big\}
            +\frac{\sigma\kappa\kappa_\psi\sqrt{\kappa\mu\rho r^2\log d}}
            {\underline c_\psi\sqrt{\pi\sigma_r^*}}\\
            &+\frac{\sigma^2}{\underline c_{\sf{Hess}}}\kappa\mu r
            \sqrt{\frac{\rho\log d}{\pi(n\wedge p)}}\,\mathfrak R_t.
        \end{align*}
        The last term retains the small sampling ratio in
        \eqref{eq: auxiliary LOO sampling ratio} and is absorbed under the
        conditions used in \eqref{eq: augmented LOO one-step recursion}.
        Matrix Bernstein's inequality and the argument leading to
        \eqref{eq: decomposition of the surrogate sequence} therefore give,
        for sufficiently large $D_2$, on the joint event,
        \begin{align}
            & \norm{\big[\mc P_\Omega\big(\psi(\mc M(\bo \theta^{\sf{aug},t})) - \mb R\big)\big]_{i,\cdot} (\mb Y^{\sf{aug},t}\mb O_{\sf{aug}}^t - \mb Y^t)} \\
            \leq  &  \frac{\underline c_{\sf{Hess}}}{32\underline c_\psi} D_2 \kappa \kappa_\psi \underline c_\psi \sqrt{\frac{\mu \rho r}{d}} \pi\sigma_r^* (\delta^{\sf{aug}}+\tau_{\sf aug}^t) \log d.
        \end{align}

        \item By leveraging a similar argument to that of Lemmas~\ref{lemma: bound for Xt+1t Xt+1 - Yt+1t Yt+1},~\ref{lemma: bound for 1t Xt}, and~\ref{lemma: balancedness and orthogonality}, we may achieve the same bounds for $ \norm{\mb 1_n\t \mb X^{\sf{aug},t}}_2$ and $\norm{\mb X^{\sf{aug},t\top} \mb X^{\sf{aug},t} - \mb Y^{\sf{aug},t\top} \mb Y^{\sf{aug},t}}\fb$ as in Lemma~\ref{lemma: balancedness and orthogonality}; namely, on the joint event,
        \begin{align}
            & \norm{\mb 1_n\t \mb X^{\sf{aug},t}}_2 \lesssim  \sqrt{n} \delta^{\sf{aug}}, \label{eq: orthogonality condition for auxiliary sequence}
            \\ 
            & \bignorm{ \mb X^{\sf{aug},t\top}  \mb X^{\sf{aug},t} -  \mb Y^{\sf{aug},t\top}  \mb Y^{\sf{aug},t}}\fb \lesssim (\sigma_1^*)^{\frac12} \delta^{\sf{aug}}. \label{eq: balancedness condition for auxiliary sequence}
        \end{align}
        This combined with Lemma~\ref{lemma: balancedness and orthogonality} leads, on the joint event, to the conclusion that
        \begin{align*}
        & 2 \eta c_{\perp}\frac{\pi \sigma_r^*}{n}\norm{\mb 1_n\t \big(\mb X^{\sf{aug},t}\mb O_{\sf{aug}}^t - \mb X^t \big)}_2
        + 4 \eta c_{\sf{aug}} \pi
        \norm{\mb X^{\sf{aug},t}\mb O_{\sf{aug}}^t}\ti \\
        &\qquad \cdot \norm{(\mb X^{\sf{aug},t}\mb O_{\sf{aug}}^t)\t
        (\mb X^{\sf{aug},t}\mb O_{\sf{aug}}^t)
        -(\mb Y^{\sf{aug},t}\mb O_{\sf{aug}}^t)\t
        (\mb Y^{\sf{aug},t}\mb O_{\sf{aug}}^t)}\fb
        \\
            \leq & \eta\frac{\underline c_{\sf{Hess}}}{16\underline c_\psi}
            D_2\kappa\kappa_\psi\underline c_\psi
            \sqrt{\frac{\mu\rho r}{d}}\pi\sigma_r^*
            (\delta^{\sf aug}+\tau_{\sf aug}^t)\log d.
        \end{align*}
        \end{itemize}

        With these pieces in place, substituting them into \eqref{eq: auxiliary sequence decomposition} implies that
        \begin{align}
            & \bignorm{\mb X^{\sf{aug},t+1}_i\mb O_{\sf{aug}}^t - \mb X^{t+1}_i}_2  \leq (1 - \eta \underline c_{\sf{Hess}} \pi \sigma_r^* / 2) \bignorm{\mb X^{\sf{aug},t}_i\mb O_{\sf{aug}}^t - \mb X^{t}_i}_2 \\
            & \quad + \eta \frac{\underline c_{\sf{Hess}}}{8\underline c_\psi} D_2 \kappa \kappa_\psi \underline c_\psi \sqrt{\frac{\mu \rho r  }{d}} \pi \sigma_r^* (\delta^{\sf{aug}}+\tau_{\sf aug}^t) \log d,
        \end{align}
        holds on the joint event.
        
        It remains to replace the candidate alignment $\mb O_{\sf{aug}}^t$
        by the Procrustes rotation at time $t+1$. To retain the one-step
        factor $\eta$, we introduce a surrogate matrix as in the rowwise GD
        analysis. Set
        \begin{align*}
            \mb A^t&\coloneqq\mb F^{\sf{aug},t}\mb O_{\sf{aug}}^t,
            &\mb B^t&\coloneqq\mb F^t,\\
            \mb D_{\sf{aug}}^t&\coloneqq
            \begin{bmatrix}
                \nabla_{\mb X}L^{\sf{aug}}(\bo\theta^{\sf{aug},t})
                \mb O_{\sf{aug}}^t\\
                \nabla_{\mb Y}L^{\sf{aug}}(\bo\theta^{\sf{aug},t})
                \mb O_{\sf{aug}}^t
            \end{bmatrix},
            &\mb D^t&\coloneqq
            \begin{bmatrix}
                \nabla_{\mb X}L(\bo\theta^t)\\
                \nabla_{\mb Y}L(\bo\theta^t)
            \end{bmatrix}.
        \end{align*}
        Thus $\mb A^{t+1,\sf c}\coloneqq\mb A^t-
        \eta\mb D_{\sf{aug}}^t
        =\mb F^{\sf{aug},t+1}\mb O_{\sf{aug}}^t$ and
        $\mb B^{t+1}=\mb B^t-\eta\mb D^t=\mb F^{t+1}$. For a square
        matrix $\mb C$, write $\operatorname{skew}(\mb C)
        \coloneqq(\mb C-\mb C\t)/2$, and let
        $\mb C^{t+1}\coloneqq{\mb A^{t+1,\sf c}}\t\mb B^{t+1}$. Define
        \begin{align*}
            \mb A^{t+1,\sf{srg}}
            \coloneqq{}&\mb A^{t+1,\sf c}
            +\mb B^{t+1}
            \big({\mb B^{t+1}}\t\mb B^{t+1}\big)^{-1}
            \operatorname{skew}(\mb C^{t+1}).
        \end{align*}
        By construction, ${\mb A^{t+1,\sf{srg}}}\t\mb B^{t+1}
        =(\mb C^{t+1}+{\mb C^{t+1}}\t)/2$, which is positive definite
        under the RIC and the step-size condition. Hence the identity is the
        Procrustes alignment of $\mb A^{t+1,\sf{srg}}$ with
        $\mb B^{t+1}$.

        At time $t$, ${\mb A^t}\t\mb B^t$ is symmetric by the definition of
        $\mb O_{\sf{aug}}^t$. Rotational invariance further implies that
        ${\mb A^t}\t\mb D_{\sf{aug}}^t$ and
        ${\mb B^t}\t\mb D^t$ are symmetric. Expanding the new cross-Gram
        matrix therefore gives
        \begin{align*}
            \operatorname{skew}(\mb C^{t+1})
            ={}&-\eta\operatorname{skew}\Big(
            {\mb D_{\sf{aug}}^t}\t(\mb B^t-\mb A^t)
            +(\mb A^t-\mb B^t)\t\mb D^t\Big)
            +\eta^2\operatorname{skew}\big(
            {\mb D_{\sf{aug}}^t}\t\mb D^t\big),\\
            \bignorm{\operatorname{skew}(\mb C^{t+1})}
            \lesssim{}&\eta\bignorm{\mb A^t-\mb B^t}\fb
            \Big(\bignorm{\mb D_{\sf{aug}}^t}\fb
            +\bignorm{\mb D^t}\fb\Big)
            +\eta^2\bignorm{\mb D_{\sf{aug}}^t}\fb
            \bignorm{\mb D^t}\fb.
        \end{align*}
        The last display, the RIC bounds, and the step-size and signal-strength
        conditions ensure that
        $\bignorm{\operatorname{skew}(\mb C^{t+1})}
        \leq\sigma_r\big((\mb C^{t+1}+{\mb C^{t+1}}\t)/2\big)$.
        Since ${\mb O_{\sf{aug}}^t}\t\mb O_{\sf{aug}}^{t+1}$ is the
        polar factor of $\mb C^{t+1}$, whereas the polar factor of its
        symmetric part is $\mb I_r$, Lemma~\ref{lemma: perturbation theory of optimal rotation under F norm}
        yields
        \begin{align}
            &\norm{{\mb O_{\sf{aug}}^t}\t
            \mb O_{\sf{aug}}^{t+1}-\mb I_r}
            \lesssim\frac{\eta}{\sigma_r^*}\Big\{
            \bignorm{\mb A^t-\mb B^t}\fb
            \Big(\bignorm{\mb D_{\sf{aug}}^t}\fb
            +\bignorm{\mb D^t}\fb\Big)
            +\eta\bignorm{\mb D_{\sf{aug}}^t}\fb
            \bignorm{\mb D^t}\fb\Big\}.
            \label{eq: auxiliary one-step rotation increment}
        \end{align}
        The squared-norm argument in the proof of
        Lemma~\ref{lemma: bound on nabla Laug}, applied to the augmented
        sequence, gives $\bignorm{\mb D_{\sf{aug}}^t}\fb
        \leq\bignorm{\nabla L^{\sf{aug}}(\bo\theta^0)}_2$. Moreover,
        $\bignorm{\mb D^t}\fb$ is bounded by
        $\bignorm{\nabla L^{\sf{aug}}(\bo\theta^t)}_2
        +\bignorm{\nabla L^{\sf{diff}}(\bo\theta^t)}_2$.
        The Frobenius inductive bound \eqref{eq: auxiliary induction Frobenius},
        Lemmas~\ref{lemma: bound on nabla Laug}
        and~\ref{lemma: gradient norm at theta0}, and the balancedness bound
        therefore control every term on the right-hand side of
        \eqref{eq: auxiliary one-step rotation increment}. Consequently, the
        signal-strength and step-size conditions, after choosing $D_2$
        sufficiently large relative to $D_1$, imply
        \begin{align}
            &\bignorm{\mb F^{t+1}}\ti
            \norm{{\mb O_{\sf{aug}}^t}\t
            \mb O_{\sf{aug}}^{t+1}-\mb I_r}
            \leq \eta\frac{\underline c_{\sf{Hess}}}{8\underline c_\psi}
            D_2\kappa\kappa_\psi\underline c_\psi
            \sqrt{\frac{\mu\rho r}{d}}\pi\sigma_r^*
            \delta^{\sf{aug}}\log d.
            \label{eq: auxiliary rowwise rotation contribution}
        \end{align}

        Consequently, one has 
        \begin{align*}
            &\bignorm{\mb X^{\sf{aug},t+1}_i
            \mb O_{\sf{aug}}^{t+1}-\mb X^{t+1}_i}_2\\
            \leq &
            \bignorm{\mb X^{\sf{aug},t+1}_i\mb O_{\sf{aug}}^t
            -\mb X^{t+1}_i}_2
            +\norm{\mb X^{t+1}_i}_2
            \norm{{\mb O_{\sf{aug}}^t}\t\mb O_{\sf{aug}}^{t+1}-\mb I_r} \\ 
            \leq & \Big(1-\frac{\eta\underline c_{\sf{Hess}}\pi\sigma_r^*}{2}\Big)
            \bignorm{\mb X^{\sf{aug},t}_i\mb O_{\sf{aug}}^t - \mb X^{t}_i}_2
            +\frac{\eta\underline c_{\sf{Hess}}\pi\sigma_r^*}{4}
            D_2\kappa\kappa_\psi\sqrt{\frac{\mu\rho r}{d}}
            \delta^{\sf{aug}}\log d \\
            &\qquad+\frac{\eta\underline c_{\sf{Hess}}\pi\sigma_r^*}{8}
            D_2\kappa\kappa_\psi\sqrt{\frac{\mu\rho r}{d}}
            \tau_{\sf aug}^t\log d \\
            \leq & D_2 \kappa \kappa_\psi  \sqrt{\frac{\mu \rho r}{d}}
            (\delta^{\sf{aug}}+\tau_{\sf aug}^{t+1}) \log d.
        \end{align*}
        The transient term is absorbed by the same contraction-gap argument.

    The analysis for the intercept $\bo \zeta$ and the $\mb Y$ factor is analogous, and we omit the details. 
     Taking the
        maximum over the rows of both factors, together with the corresponding
        bound for $\bo\zeta$, verifies
        \eqref{eq: auxiliary induction rowwise} at time $t+1$.

        \bigskip

        Lastly, following the spirit of the proof of Lemma~\ref{lemma: bound for 1t Xt}, we first control the one-step candidate $\norm{\mb 1_n\t \big( \mb X^{\sf{aug},t+1}\mb O_{\sf{aug}}^t - \mb X^{t+1}\big)}_2$ via the following decomposition: 
        \begin{align}
            & \norm{\mb 1_n\t \big( \mb X^{\sf{aug},t+1}\mb O_{\sf{aug}}^t - \mb X^{t+1}\big)}_2 \\ 
    & \leq   \underbrace{(1 - 2 c_{\perp} \eta \pi \sigma_r^*) \norm{\mb 1_n\t \big( \mb X^{\sf{aug},t}\mb O_{\sf{aug}}^t - \mb X^{t}\big)}_2}_{\gamma_1}\\
    &\quad + \underbrace{\eta \norm{\mb 1_n\t \mc P_\Omega(\psi(\mc M(\bo \theta^{\sf{aug},t})) - \psi(\mc M(\bo \theta^{t}))) \mb Y^t}_2}_{\gamma_2}\\
            & + \underbrace{\eta \norm{\mb 1_n\t \mc P_\Omega(\psi(\mc M(\bo \theta^{\sf{aug},t})) - \psi(\mc M(\bo \theta^*)) -  \mb E)(\mb Y^{\sf{aug},t}\mb O_{\sf{aug}}^t - \mb Y^t)}_2}_{\gamma_3} \\
    & + \underbrace{\begin{aligned}4\eta c_{\sf{aug}}\pi
            \Big\|\mb 1_n\t\mb X^{\sf{aug},t}\mb O_{\sf{aug}}^t
            \Big[(\mb X^{\sf{aug},t}\mb O_{\sf{aug}}^t)\t
            (\mb X^{\sf{aug},t}\mb O_{\sf{aug}}^t)\\
            {}-(\mb Y^{\sf{aug},t}\mb O_{\sf{aug}}^t)\t
            (\mb Y^{\sf{aug},t}\mb O_{\sf{aug}}^t)\Big]\Big\|_2\end{aligned}}_{\gamma_4}.
        \end{align}
        \begin{itemize}
            \item It follows from \eqref{eq: auxiliary induction centering} that, on the joint event,
            \begin{align}
                & \gamma_1 \leq (1 - 2 \eta c_{\perp} \pi \sigma_r^*) D_3 \kappa \sqrt{n \mu r} \frac{\sigma^2}{c_{\perp}}  (\delta^{\sf{aug}}+\tau_{\sf aug}^t) \log d.
            \end{align}
            \item The second term $\gamma_2$ can be further written as 
            \begin{align*}
                & \gamma_2 \leq \eta  \sigma^2 \sqrt{n} \norm{\mb X^{\sf{aug},t}\mb O_{\sf{aug}}^t - \mb X^t} \max_{i\in[n]} \norm{\diag(\bo \Omega_{i,\cdot}) \mb Y^{t}}^2 \\
                & + \eta \sqrt{n} \max_{i\in[n]} \Big\|\Big[\Big(
                \psi(\mc M([\bo\zeta^{\sf{aug},t},
                \mb X^{\sf{aug},t}\mb O_{\sf{aug}}^t,
                \mb Y^{\sf{aug},t}\mb O_{\sf{aug}}^t])) \\
                &\qquad - \psi(\mc M([\bo\zeta^t,
                \mb X^{\sf{aug},t}\mb O_{\sf{aug}}^t,\mb Y^t]))\Big)
                \circ \bo \Omega \Big]_{i,\cdot}\Big\|_2 \cdot \norm{\diag(\bo \Omega_{i,\cdot})\mb Y^t}  \\
                \leq &  D_7  \eta \sqrt{n \mu r} \sigma^2\pi \sigma_1^*  (\delta^{\sf{aug}}+\tau_{\sf aug}^t) \log d
             \end{align*}
            on the joint event. Here, the first term follows from \eqref{eq: auxiliary induction rowwise}, whereas the second follows from \eqref{eq: decomposition of the surrogate sequence}.

            \item The auxiliary rowwise hypothesis
            \eqref{eq: auxiliary induction rowwise} and the actual-GD
            bounds place $\bo\theta^{\sf{aug},t}$ in the RIC, as verified
            above. Both parameter triples therefore satisfy the factor
            row bounds, and their predictor interpolants lie in
            $\mc D_\psi$. Applying
            Lemma~\ref{lemma: injectivity of P_Omega under nonlinearity}
            with endpoints $\bo\theta^{\sf{aug},t}$ and $\bo\theta^*$,
            together with the bound on
            $\distf(\bo\theta^{\sf{aug},t},\bo\theta^t)$, gives, on the
            joint event,
            \begin{align}
                & \gamma_3 \leq D_8 \eta \sqrt{n \mu r} \sigma^2\pi \sigma_1^*  \delta^{\sf{aug}} \log d.
            \end{align}

            \item Finally, rotation invariance and
            \eqref{eq: orthogonality condition for auxiliary sequence}--\eqref{eq: balancedness condition for auxiliary sequence} give
            \begin{align}
                \gamma_4
                &\leq 4\eta c_{\sf{aug}}\pi
                \bignorm{\mb 1_n\t\mb X^{\sf{aug},t}}_2
                \bignorm{{\mb X^{\sf{aug},t}}\t\mb X^{\sf{aug},t}
                -{\mb Y^{\sf{aug},t}}\t\mb Y^{\sf{aug},t}}\fb \\
                &\lesssim \eta c_{\sf{aug}}\pi
                \sqrt{n\sigma_1^*}\,(\delta^{\sf{aug}})^2
                \leq D_9\eta\sqrt{n\mu r}\,\sigma^2\pi\sigma_1^*
                \delta^{\sf{aug}}\log d.
            \end{align}
            Here, the last inequality follows from
            $c_{\sf{aug}}\asymp\underline c_\psi\leq\sigma^2$ and
            $\delta^{\sf{aug}}\ll\sqrt{\mu r\sigma_1^*}\log d$, which is
            implied by \eqref{eq: epsilon GD initialization requirement},
            \eqref{eq: definition of delta^aug}, and the stated signal-strength
            conditions.
        \end{itemize}

        Set $q_\perp\coloneqq1-2c_\perp\eta\pi\sigma_r^*$ and
        $P\coloneqq D_3\kappa\sqrt{n\mu r}\,\sigma^2
        c_\perp^{-1}\log d$. Since $q_\perp<\nu$, the same contraction-gap
        argument, with $D_3$ sufficiently large relative to $D_7$, $D_8$, and
        $D_9$, gives on the joint event that
        \begin{align}
            &\norm{\mb 1_n\t \big( \mb X^{\sf{aug},t+1}\mb O_{\sf{aug}}^t - \mb X^{t+1}\big)}_2 \\
            \leq{}&q_\perp\norm{\mb 1_n\t \big( \mb X^{\sf{aug},t}\mb O_{\sf{aug}}^t - \mb X^t\big)}_2
            +\frac{1-q_\perp}{2}P\delta^{\sf{aug}}
            +\frac{\nu-q_\perp}{2}P\tau_{\sf aug}^t \\
            \leq{}&P(\delta^{\sf{aug}}+\tau_{\sf aug}^{t+1}).
        \end{align}
        Finally, the identity
        $\mb X^{\sf{aug},t+1}\mb O_{\sf{aug}}^{t+1}
        =\mb X^{\sf{aug},t+1}\mb O_{\sf{aug}}^t
        {\mb O_{\sf{aug}}^t}\t\mb O_{\sf{aug}}^{t+1}$ gives
        \begin{align*}
            &\norm{\mb 1_n\t\big(
            \mb X^{\sf{aug},t+1}\mb O_{\sf{aug}}^{t+1}
            -\mb X^{t+1}\big)}_2\\
            \leq{}&\norm{\mb 1_n\t\big(
            \mb X^{\sf{aug},t+1}\mb O_{\sf{aug}}^t
            -\mb X^{t+1}\big)}_2
            +\norm{\mb 1_n\t\mb X^{t+1}}_2
            \norm{{\mb O_{\sf{aug}}^t}\t\mb O_{\sf{aug}}^{t+1}-\mb I_r}.
        \end{align*}
        The unused halves of the stationary and transient contraction budgets,
        together with the actual-sequence centering bound and
        \eqref{eq: auxiliary one-step rotation increment}, absorb the second
        term. This verifies \eqref{eq: auxiliary induction centering} at time
        $t+1$ with the Procrustes rotation.

        At $t=0$, \eqref{eq: auxiliary induction Frobenius}--\eqref{eq: auxiliary induction centering}
        follow because the auxiliary and actual sequences have the same
        initialization, while \eqref{eq: auxiliary induction LOO} follows from
        the initialization leave-one-out bound. It remains to verify the
        final-time assertion. The warm-start absorption in
        Theorem~\ref{thm: general GD} gives
        \begin{align*}
            &\nu^{t_0}\Big\{\kappa\kappa_\psi\sqrt{\mu\rho r}\,\delta_{\sf{row}}^0+\delta_{\sf{loo}}^0\Big\}
            \lesssim
            \frac{\sigma\kappa\kappa_\psi
            \sqrt{\kappa\mu\rho r^2\log d}}
            {\underline c_\psi\sqrt{\pi\sigma_r^*}}.
        \end{align*}
        It follows that
        \begin{align*}
            \tau_{\sf aug}^{t_0}
            &\lesssim
            \sqrt{\frac{\pi}{\mu r}}
            \frac{\sigma\kappa\kappa_\psi
            \sqrt{\kappa\mu\rho r^2\log d}}
            {\underline c_\psi\sqrt{\pi\sigma_r^*}}=\frac{\pi}{\sqrt\mu}
            \Big\{\frac{1}{\sqrt{\pi(n\wedge p)}}
            \frac{\sigma\kappa\kappa_\psi
            \sqrt{\kappa\mu r d\log d}}
            {\underline c_\psi\sqrt{\pi\sigma_r^*}}\Big\}
            \leq\delta^{\sf aug},
        \end{align*}
        where the last inequality uses $\pi\leq1$ and $\mu\geq1$. Thus all
        subsequent applications at $t=t_0$ involve only the static radius
        $\delta^{\sf aug}$.
    \end{proof}

    \subsection{Linear Approximation for Gradient Descent Iterates}
    \label{subsec: proof of inference results}
    As outlined in Section~\ref{subsubsec: key ingredients in the proof}, the linear representation requires rowwise stationarity, global control of the cross-block nuisance terms, and a rowwise score expansion. We first state the five lemmas establishing these ingredients in dependency order and then give their proofs.

    The full-space Newton heuristic would give
    \eq{
    \Big[ \int_0^1 \nabla^2 L((1 -t) \bo \theta^* + t \bo \theta^{t_0} ) \mathrm d t \Big]^{-1} \big(\nabla L(\bo \theta^{t_0}) - \nabla L(\bo \theta^*)\big) \approx -\nabla^2 L(\bo \theta^*)^{-1} \nabla L(\bo \theta^*). 
    \label{eq: heuristic of asymptotic normality}
    }
    This expression is only heuristic because $L$ has an additional scaling degeneracy absent from $L^{\sf{aug}}$, while both objectives remain rotationally invariant. We therefore avoid inverting the full Hessian of $L$. First, the auxiliary sequence under $L^{\sf{aug}}$ supplies rowwise stationarity of the final iterate. Second, a global Taylor equation for the aligned error is projected onto the range and null spaces of the singular Hessian to control the weighted nuisance terms. Finally, we return to the rowwise score decomposition, expand the resulting surrogate gradients around $\bo\theta^*$, and invert only the rowwise information matrices. This is the formal counterpart of the three-step sketch in the main text.

    We begin with the auxiliary-iterate argument for rowwise stationarity. Gradient decay under $L^{\sf{aug}}$, the comparison between the auxiliary and actual sequences, and approximate balancedness together transfer the stationary behavior of the augmented sequence to the final iterate of Algorithm~\ref{alg: GD}. The resulting bound is recorded next.
    \begin{lemma}
    \label{lemma: gradient norm control}
    Suppose that the assumptions of Theorem~\ref{thm: general GD} hold. Then, with probability at least $1-O(d^{-c})$,
    \begin{align*}
        & \max_{i\in[n]} \norm{\nabla_{\mb X_{i,\cdot}} L(\bo \theta^{t_0} )}_2 \vee \max_{j\in[p]} \bignorm{\nabla_{\mb Y^{\sf{app}}_{j,\cdot}} L(\bo \theta^{t_0} ) }_2  \lesssim \sigma^2 \kappa \kappa_\psi \sqrt{\frac{\mu \rho r}{d}}\pi \sigma_1^* \delta^{\sf{aug}} \log d.
    \end{align*}
    \end{lemma}

    We next implement the global control of the nuisance terms in the rowwise score decomposition. These terms are weighted projections of the aligned estimation error and are not sharply controlled by the rowwise GD rate alone. Introduce the augmented factors
    \begin{align*}
        \mb X^{\sf{app},t_0}&\coloneqq\big(\sqrt{\sigma_r^*/n}\mb 1_n,\mb X^{t_0}\big),
        &\mb X^{\sf{app},*}&\coloneqq\big(\sqrt{\sigma_r^*/n}\mb 1_n,\mb X^*\big),\\
        \mb Y^{\sf{app},t_0}&\coloneqq\big(\sqrt{n/\sigma_r^*}\bo\zeta^{t_0},\mb Y^{t_0}\big),
        &\mb Y^{\sf{app},*}&\coloneqq\big(\sqrt{n/\sigma_r^*}\bo\zeta^*,\mb Y^*\big),
        \qquad \mb O^{\sf{app},t_0}\coloneqq\diag(1,\mb O^{t_0}).
    \end{align*}
    The lemma below derives the required projections from the global Taylor equation for the aligned error. Its proof controls the range component through the Hessian pseudoinverse, the scaling component through approximate balancedness, and the rotational component through the Procrustes condition. The $\mb Y^{\sf{app}}$-error is projected only onto $\mb Y^*$, as required by the row-gradient nuisance term; hence no spectral-norm bound for $\mb Y^{\sf{app},*}$ is needed.

    \begin{lemma}
    \label{lemma: projection of difference onto specific directions}
        Suppose that the assumptions of Theorem~\ref{thm: general GD} hold. Then, with probability at least $1-O(d^{-c})$,
        \begin{align*}
            & \max_{j\in[p]} \norm{\big(\mb X^{\sf{app},t_0} \mb O^{\sf{app},t_0} - \mb X^{\sf{app},*}\big)\t \diag(\psi'(\mc M(\bo \theta^*))_{i,j} \ind\{(i,j) \in \Omega\})_{i\in[n]} \mb X^{\sf{app},*}} \\
            & \vee \max_{i\in[n]} \norm{\big(\mb Y^{\sf{app},t_0} \mb O^{\sf{app},t_0} - \mb Y^{\sf{app},*} \big)\t \diag(\psi'(\mc M(\bo \theta^*))_{i,j} \ind\{(i,j) \in \Omega\} )_{j\in[p]} \mb Y^*} \\
            \lesssim & \sigma\kappa\kappa_\psi\mu r^2\sqrt{\rho\log d}
            \Big\{\sqrt{\log d}+\kappa\kappa_\psi\sqrt{\mu r}\Big\}
            +\sigma^2\kappa_\psi^2\kappa^2\mu r^3\log d
            \Big\{\frac{d}{\underline c_\psi\sigma_r^*}
            +\frac{\kappa_\psi\kappa^2\sqrt{\mu}\,r\rho}{\underline c_\psi}\Big\}\\
            &\qquad+\sigma^2\pi\kappa_\psi^2\kappa^2 r
            \sqrt{\mu\rho r\sigma_1^*}\delta^{\sf{aug}}\log d.
        \end{align*}
    \end{lemma}

    With rowwise stationarity and the weighted nuisance projections available, we return to the rowwise score decomposition. Replacing the nuisance factor block by its population counterpart produces the two surrogate gradients below. For a parameter triple $[\bo\zeta,\mb X,\mb Y]$, define the augmented column-block gradient by
    \begin{align*}
        \nabla_{\mb Y^{\sf{app}}_{j,\cdot}}L([\bo\zeta,\mb X,\mb Y])
        \coloneqq\Big(\sqrt{\frac{\sigma_r^*}{n}}\nabla_{\zeta_j}L([\bo\zeta,\mb X,\mb Y]),
        \nabla_{\mb Y_{j,\cdot}}L([\bo\zeta,\mb X,\mb Y])\Big).
    \end{align*}
    This convention matches the scaling of $\mb Y^{\sf{app}}$ introduced above;
    throughout this argument, gradients with respect to matrix rows are
    represented as row vectors.
    \begin{lemma}\label{lemma: gradient replacement 1}
        Under the assumptions of Theorem~\ref{thm: general GD}, with probability at least $1- O(d^{-c})$,
        \begin{align*}
            & \max_{i\in[n]} \norm{\nabla_{\mb X_{i,\cdot}} L([\bo \zeta^*, \mb X^{t_0} \mb O^{t_0}, \mb Y^* ])}_2\\
            &\qquad\vee \max_{j\in[p]} \norm{\nabla_{\mb Y^{\sf{app}}_{j,\cdot}} L([\bo\zeta^{t_0}, \mb X^*,\mb Y^{t_0} \mb O^{t_0} ])}_2\\
            \lesssim{}&\Big[\frac{1}{\sqrt{\pi(n\wedge p)}}
            +\frac{\sigma\kappa^2\kappa_\psi\sqrt r}
            {\underline c_\psi\sqrt{\pi(n\wedge p)}}
            +\frac{\sigma\sqrt{\mu r d}}
            {\underline c_\psi\sqrt{\pi}\sigma_r^*}\Big]\sigma\kappa_\psi^2\kappa^2r^4\mu^2\rho
            \sqrt{\pi\sigma_1^*}\log d\\
            &+\sigma^2\pi\kappa_\psi^2\kappa^2\mu r^2
            \sqrt{\frac{\rho}{n\wedge p}}\sigma_1^*
            \delta^{\sf{aug}}\log d.
    \end{align*}
    \end{lemma}

    It remains to linearize the surrogate gradients around $\bo\theta^*$. Their first-order terms are the aligned row errors multiplied by the corresponding rowwise information matrices $\mb H^*_{\mb X_i}$ and $\mb H^*_{\mb Y^{\sf{app}}_{j}}$ defined in \eqref{eq: definition of fisher matrices}, while the next lemma controls the Taylor remainders. Let $\mb g_{\perp}^{t_0}\coloneqq2c_{\perp}\pi\sigma_r^*n^{-1}\mb 1_n\t\mb X^{t_0}\mb O^{t_0}$ denote the row gradient contributed by the centering penalty at the aligned iterate.
    \begin{lemma}
    \label{lemma: gradient replacement 2}
        Under the assumptions of Theorem~\ref{thm: general GD}, with probability at least $1-O(d^{-c})$,
        \begin{align*}
    & \max\Big\{\max_{i\in[n]}\Bignorm{\nabla_{\mb X_{i,\cdot}}L([\bo\zeta^*,\mb X^{t_0}\mb O^{t_0},\mb Y^*])\\
    &\qquad -\nabla_{\mb X_{i,\cdot}}L(\bo\theta^*)-\mb g_{\perp}^{t_0}-\big(\mb X^{t_0}_{i,\cdot}\mb O^{t_0}-\mb X^*_{i,\cdot}\big)\mb H^*_{\mb X_i}}_2,\\
    & \qquad\max_{j\in[p]}\Bignorm{\nabla_{\mb Y^{\sf{app}}_{j,\cdot}}L([\bo\zeta^{t_0},\mb X^*,\mb Y^{t_0}\mb O^{t_0}])\\
    &\qquad -\nabla_{\mb Y^{\sf{app}}_{j,\cdot}}L(\bo\theta^*)-\big(\mb Y^{\sf{app},t_0}_{j,\cdot}\mb O^{\sf{app},t_0}-\mb Y^{\sf{app},*}_{j,\cdot}\big)\mb H^*_{\mb Y^{\sf{app}}_{j}}}_2\Big\}\\
    & \lesssim{} \pi\sigma_1^*\Big[\frac{\sigma\kappa\kappa_\psi\sqrt{\kappa\mu\rho r^2\log d}}{\underline c_\psi\sqrt{\pi\sigma_r^*}}\Big]
            \Big[\frac{\sigma^2\sigma\kappa^2\kappa_\psi\mu r^{\frac32}\sqrt{\rho\log d}}{\underline c_\psi\sqrt{\pi(n\wedge p)}}\Big].
        \end{align*}
    \end{lemma}

    The preceding two lemmas now complete the rowwise score decomposition: Lemma~\ref{lemma: gradient replacement 1} combines rowwise stationarity with the nuisance-block replacement to control the surrogate gradients, while Lemma~\ref{lemma: gradient replacement 2} identifies the rowwise information term and its Taylor remainder. Inverting the rowwise information matrices yields the score-based linear approximation used in the main text.
    \begin{lemma}
    \label{lemma: gradient replacement 3}
        Under the assumptions of Theorem~\ref{thm: general GD}, with probability at least $1-O(d^{-c})$,
        \begin{align*}
            & \max \Big\{ \max_{i\in[n]} \norm{\mb X^{t_0}_{i,\cdot}\mb O^{t_0}-\mb X^*_{i,\cdot}-\mc P_\Omega(\mb E)_{i,\cdot}\mb Y^*(\mb H^*_{\mb X_i})^{-1}}_2,\\
            &\qquad \max_{j\in[p]} \norm{\mb Y^{\sf{app},t_0}_{j,\cdot}\mb O^{\sf{app},t_0}-\mb Y^{\sf{app},*}_{j,\cdot}-\mc P_\Omega(\mb E)_{\cdot,j}\t\mb X^{\sf{app},*}(\mb H^*_{\mb Y^{\sf{app}}_{j}})^{-1}}_2 \Big\}\\
            \lesssim{}&\Big\{\frac{\kappa^2\kappa_\psi^2\sqrt{\log d}}
            {\sqrt{\pi(n\wedge p)}}
            +\frac{\sigma\kappa^2\kappa_\psi\sqrt r}
            {\underline c_\psi\sqrt{\pi(n\wedge p)}}
            +\frac{\sigma\sqrt{\mu r d}}
            {\underline c_\psi\sqrt\pi\,\sigma_r^*}\Big\}
            \frac{\sigma\kappa_\psi^2\kappa^{\frac52}r^4\mu^2\rho\log d}
            {\underline c_\psi\sqrt{\pi\sigma_r^*}}\\
            &+\frac{\kappa^5\kappa_\psi^5\mu r^2\rho\log d}
            {c_\perp\sqrt{\pi\sigma_r^*}}
            \big\{\sigma\sqrt{\kappa\mu r\log d}
            +\xi^{\sf{gd}}+\xi^{\sf{gd},\sf{loo}}\big\}
            \eqqcolon\xi^{\sf{approx}}.
        \end{align*}
    \end{lemma}

    \begin{proof}[Proof of Theorem~\ref{thm: parameter-explicit linear approximations}]
    By Lemma~\ref{lemma: scaling equivalence}, it suffices to work in the balanced
    parametrization $\omega=1$. The assumed range of $c_\perp$, together with
    the prescribed learning rate and iteration count in
    Theorem~\ref{thm: parameter-explicit GD}, meets the calibration of
    Theorem~\ref{thm: general GD}.

    As established in the proof of
    Theorem~\ref{thm: parameter-explicit GD},
    Theorem~\ref{thm: detailed USVT}, and
    Theorem~\ref{thm: general OS}, together with
    Lemma~\ref{lemma: pipeline satisfies GD initialization}, verify the
    initialization requirements of Theorem~\ref{thm: general GD}. Consequently,
    the general gradient-descent guarantees and
    Lemmas~\ref{lemma: gradient replacement 1}--\ref{lemma: gradient replacement 3}
    apply to the output of the full pipeline.

    Lemma~\ref{lemma: gradient replacement 3}, stated in terms of the information
    matrices $\mb H^*_{\mb X_i}$ and $\mb H^*_{\mb Y^{\sf{app}}_{j}}$, gives the asserted expansion for
    each row of $\mb X^{t_0}$ and the joint expansion for each row of
    $\mb Y^{\sf{app},t_0}$. The first coordinate of the latter, multiplied by
    $\sqrt{\sigma_r^*/n}$, yields the expansion for $\zeta_j^{t_0}$, while its
    remaining $r$ coordinates yield the expansion for $\mb Y_{j,\cdot}^{t_0}$.
    To verify the required remainder, multiply the bound in
    Lemma~\ref{lemma: gradient replacement 3} by
    $\sigma\sqrt{\pi\sigma_1^*}(\log d)^2$. Substituting the concrete
    stage-two rate established in the proof of Theorem~\ref{thm: OS} into
    \eqref{eq: GD initialization rates from pipeline} gives
    \begin{align*}
        &\sigma\sqrt{\pi\sigma_1^*}(\log d)^2
        \xi^{\sf{approx}}\\
        \lesssim{}&
        \frac{\kappa^5\kappa_\psi^5\mu^2r^4\rho
        (\log d)^{\frac72}}{\sqrt{\pi(n\wedge p)}}
        +\frac{\kappa^5\kappa_\psi^5\mu^2r^{\frac92}\rho
        (\log d)^3}{\sigma\sqrt{\pi(n\wedge p)}}+\frac{\kappa^3\kappa_\psi^4\mu^{\frac52}r^{\frac92}\rho
        (\log d)^3}{\sigma\,\sigma_r^*}
        \sqrt{\frac d\pi}\\
        &+\frac{1}{c_\perp}\Bigg\{
        \underbrace{\sigma\bar\sigma\kappa^9\kappa_\psi^6\mu^{\frac32}r^3
        \rho^{\frac32}(\log d)^{\frac72}}_{\beta_{\perp,1}}
        +\underbrace{\sigma^2\kappa^{\frac{23}{2}}\kappa_\psi^7\mu^{\frac94}
        r^{\frac{13}{4}}\rho^2(\log d)^3\sqrt{\bar\sigma\sigma_1^*}
        \Big(\frac{\pi}{d}\Big)^{\frac14}\xi^{\frac52}}_{\beta_{\perp,2}}\\
        &\qquad\qquad
        +\underbrace{\frac{\sigma^3\kappa^{10}\kappa_\psi^6\mu^3r^4
        \rho^{\frac52}(\log d)^4\sigma_r^*}{\sqrt\pi\,d^{\frac32}}}_{\beta_{\perp,3}}
        \Bigg\}\lesssim1.
    \end{align*}
    It remains to verify the final bound. Since $\psi'\geq\underline c_\psi$
    on $\mc D_\psi$, the variation of $\psi$ between the two endpoints of
    $\mc D_\psi$ is at least
    $2c_{\mc D_\psi}\underline c_\psi\mu r\sigma_1^*/\sqrt{np}$ and at most
    $2\max_{x\in\mc D_\psi}|\psi(x)|$. Consequently,
    $\bar\sigma/\underline c_\psi\gtrsim
    \mu r\sigma_1^*/\sqrt{np}$. Combining this relation with the
    signal-strength condition and using $\sigma_1^*=\kappa\sigma_r^*$ gives
    the first bound below, whereas the second follows from
    \eqref{eq: effective information condition for linear approximation}:
    \begin{align}
        \pi(n\wedge p)&\gg\kappa^{24}\kappa_\psi^{10}\mu^7r^{11}\rho^3\xi^{10}(\log d)^6,\qquad
        \sigma^2\pi(n\wedge p)\gg\kappa^{24}\kappa_\psi^{10}\mu^7r^{11}\rho^3\xi^{10}(\log d)^6.
        \label{eq: effective sample sizes from linear approximation}
    \end{align}
    These two bounds control the first and second terms, respectively.
    The strengthened signal condition controls the third term directly, again
    using $\bar\sigma/\underline c_\psi\gtrsim
    \kappa_\psi/\sigma$. For the $c_\perp^{-1}$ contribution, the signal
    and sampling conditions, together with the preceding relation, give
    \begin{align*}
        \frac{\beta_{\perp,1}}{\beta_{\perp,2}}
        &\ll\frac{1}{\kappa^{\frac{17}{2}}\kappa_\psi^4\mu^2r^{\frac52}
        \rho^{\frac54}\xi^5\log d}=o(1),\\
        \frac{\beta_{\perp,3}}{\beta_{\perp,2}}
        &\lesssim\frac{\kappa^{-\frac52}\kappa_\psi^{-\frac12}\mu^{\frac14}
        r^{\frac14}\rho^{\frac14}\log d}{(\pi d)^{\frac34}\xi^{\frac52}}=o(1),
    \end{align*}
    where the last equality follows from the sampling condition. Thus
    $\beta_{\perp,2}$ dominates the two remaining terms, and
    \eqref{eq: cperp condition for linear approximation} bounds the full
    $c_\perp^{-1}$ contribution by a constant. The displayed remainder bound
    follows.
    Condition on $\bo\Omega$ and the event in
    \eqref{eq: uniform invertibility of rowwise information matrices}. After
    clipping and recentering the score entries at level $B$, the vector Bernstein
    inequality, the incoherence bounds, and the tail condition show that, with
    probability at least $1-O(d^{-c})$,
    \begin{equation*}
    \begin{aligned}
        &\max_{i\in[n]}\bignorm{\mc P_\Omega(\mb E)_{i,\cdot}\mb Y^*
        \mb H_{\mb X_i}^{*-1}}_2\vee
        \max_{j\in[p]}\bignorm{\mc P_\Omega(\mb E)_{\cdot,j}\t\mb X^{\sf{app},*}
        \mb H_{\mb Y^{\sf{app}}_{j}}^{*-1}}_2
        \lesssim\frac{\sigma\sqrt{\kappa\mu r\log d}}
        {\underline c_\psi\sqrt{\pi\sigma_r^*}}.
    \end{aligned}
    \end{equation*}
    The original and clipped score arrays
    coincide outside an event of probability $O(d^{-c-10})$, while recentering is
    negligible. Each fixed-index Bernstein bound is invoked with exceptional
    probability $O(d^{-c-2})$ before taking the union over $i$ and $j$, and the
    contribution involving $B$ is absorbed by the tail condition. Since
    $\kappa_\psi=\sigma^2/\underline c_\psi$ and
    $\kappa,\mu,r\geq1$, the ratio of the preceding linear-approximation
    remainder to the score bound is at most a constant multiple of
    $\{\kappa\kappa_\psi\sqrt{\mu r}(\log d)^{5/2}\}^{-1}=o(1)$.
    The expansions therefore imply the claimed rowwise error rate in the
    balanced parametrization.
    Finally, Lemma~\ref{lemma: scaling equivalence} multiplies the
    $\mb X$- and $\mb Y$-remainders by $\omega^{1/4}$ and $\omega^{-1/4}$,
    respectively, while leaving the intercept remainder unchanged. It applies in
    the same way to the rowwise error bound and proves the claim with probability
    at least $1-O(d^{-c})$.
    \end{proof}

    \begin{remark}
        The improvement in \eqref{eq: rowwise error from linear approximation}
        relies on the linear expansion. Without exploiting this expansion, the
        direct trajectory analysis gives only
        \begin{align*}
            &\sqrt{\frac{n}{\sigma_r^*}}
            \bignorm{\bo\zeta^{t_0}-\bo\zeta^*}_\infty
            \vee \omega^{-\frac14}
            \bignorm{\mb X^{t_0}\mb O^{t_0}-\mb X^*}\ti \vee \omega^{\frac14}
            \bignorm{\mb Y^{t_0}\mb O^{t_0}-\mb Y^*}\ti\\
    &\quad \lesssim
            \frac{\sigma\kappa\kappa_\psi
            \sqrt{\kappa\mu\rho r^2\log d}}
            {\underline c_\psi\sqrt{\pi\sigma_r^*}}.
        \end{align*}
        Thus, the linear representation removes the factor
        $\kappa\kappa_\psi\sqrt{\rho r}$ and attains the uniform concentration
        scale of its leading score terms.
    \end{remark}

    \subsubsection{Proofs of the Linear-Approximation Lemmas}
    We now prove Lemmas~\ref{lemma: gradient norm control}--\ref{lemma: gradient replacement 3} in order.

    \begin{proof}[Proof of Lemma~\ref{lemma: gradient norm control}]
        For the auxiliary sequence, the update direction is $\nabla L^{\sf{aug}}(\bo\theta^{\sf{aug},s})$. Its factor block is horizontal by rotational invariance. The RIC bounds from Lemma~\ref{lemma: auxiliary sequence analysis} allow us to apply the multiplication argument in \eqref{eq: averaged Hessian multiplication bound LOO} to this tangent direction. The squared-norm calculation in the proof of Lemma~\ref{lemma: bound on nabla Laug}, with $\eta\leq1/(2\bar c_{\sf{Hess}}c_\perp\pi\sigma_1^*)$, then gives $\bignorm{\nabla L^{\sf{aug}}(\bo\theta^{\sf{aug},s+1})}_2^2
        \leq\bigl(1-\eta\underline c_{\sf{Hess}}\pi\sigma_r^*\bigr)
        \bignorm{\nabla L^{\sf{aug}}(\bo\theta^{\sf{aug},s})}_2^2$.
        Iterating this inequality yields
        \begin{align*}
            & \norm{\nabla_{\bo \theta} L^{\sf{aug}}(\bo \theta^{\sf{aug},t_0})}_2 \leq \nu^{t_0} \norm{\nabla_{\bo \theta} L^{\sf{aug}}(\bo \theta^{\sf{aug},0})} = \nu^{t_0} \norm{\nabla_{\bo \theta} L(\bo \theta^{0})} \\
            \stackrel{\text{by Lemma~\ref{lemma: gradient norm at theta0}}}{\lesssim}  & \nu^{t_0} \big[ \sigma \sqrt{d \pi r\sigma_1^*} + \kappa_\psi \sqrt{ \kappa d \pi \sigma_1^* }\xi^{\sf{gd}} \big] \lesssim \sigma^2 \pi \sigma_1^* \delta^{\sf{aug}} / \sqrt{n}.
            \label{eq: global control on nabla Laug theta^aug}
        \end{align*}

        In the remainder of the proof, define
        $\mb O^{\sf{aug},t_0}\coloneqq
        \mb O_{\sf{aug}}^{t_0}\mb O^{t_0}$ and introduce the common-gauge
        representatives
        \begin{align*}
            \widetilde{\bo\theta}^{\sf{aug},t_0}
            &\coloneqq[\bo\zeta^{\sf{aug},t_0},
            \mb X^{\sf{aug},t_0}\mb O^{\sf{aug},t_0},
            \mb Y^{\sf{aug},t_0}\mb O^{\sf{aug},t_0}],\\
            \widetilde{\bo\theta}^{t_0}
            &\coloneqq[\bo\zeta^{t_0},
            \mb X^{t_0}\mb O^{t_0},
            \mb Y^{t_0}\mb O^{t_0}].
        \end{align*}
        Thus $\mb O^{\sf{aug},t_0}$ first aligns the auxiliary iterate with
        the actual iterate and then places both iterates in the population
        gauge. Right-multiplication by the common rotation $\mb O^{t_0}$
        preserves their Frobenius and rowwise discrepancies. We replace the
        two iterates by these representatives and suppress the tildes below.
        This leaves the linear predictors and objective values unchanged and,
        by rotational equivariance, preserves the norms of the factor-gradient
        blocks.

        Moving forward, we relate $\nabla_{\bo \theta} L^{\sf{aug}}(\bo \theta^{\sf{aug},t_0})$ to $\nabla_{\bo \theta} L(\bo \theta^{t_0})$ via Lemma~\ref{lemma: auxiliary sequence analysis} as follows:
        \begin{align}
            & \Bignorm{\nabla_{\mb X_{l,\cdot}} L(\bo \theta^{t_0} )}_2
            \leq \bignorm{\nabla_{\mb X_{l,\cdot}} L^{\sf{diff}}(\bo \theta^{\sf{aug},t_0})}_2 + \norm{\nabla_{\mb X_{l,\cdot}} L^{\sf{aug}}(\bo \theta^{\sf{aug},t_0})}_2 \\ 
            & \quad +  \norm{\mc P_{\mb X_l}\Big(\int_0^1\nabla^2_{\bo \theta} L((1 - s) \bo \theta^{\sf{aug},t_0} + s\bo \theta^{t_0}) \mathrm d s (\bo \theta^{\sf{aug},t_0} - \bo \theta^{t_0})  \Big)}_2 ,
            \label{eq: decomposition of the gradient at GD output}
        \end{align}
        where $\mc P_{\mb X_l}$ denotes the projection operator that extracts the coordinates corresponding to $\mb X_l$ from the vectorized parameter $\bo \theta \in \bb R^{p + nr + pr}$. 

        These three terms on the right-hand side are bounded separately as follows. 
        \begin{itemize}
            \item Invoking \eqref{eq: balancedness condition for auxiliary sequence} together with the row-wise consistency between $\mb X^{\sf{aug},t}$ and $\mb X^t$ from Lemma~\ref{lemma: auxiliary sequence analysis} implies that, with probability at least $1- O(d^{-c})$, 
            \begin{align}
                & \bignorm{\nabla_{\mb X_{l,\cdot}} L^{\sf{diff}}(\bo \theta^{\sf{aug},t_0})}_2 \lesssim \underline c_\psi \sqrt{\frac{\mu r}{n}} (\sigma_1^*)^{\frac12} \pi (\sigma_1^*)^{\frac12} \delta^{\sf{aug}}.
            \end{align}     
            \item Next, the relation \eqref{eq: global control on nabla Laug theta^aug} directly implies that 
            \begin{align}
                & \norm{\nabla_{\mb X_{l,\cdot}} L^{\sf{aug}}(\bo \theta^{\sf{aug},t_0})}_2 \lesssim \sigma^2 \pi \sigma_1^* \delta^{\sf{aug}} / \sqrt{n}.
            \end{align}
            \item Moving to the last term, we re-express the $\ell_2$ norm in a variational form:  
            \begin{align*}
                & \Bignorm{\mc P_{\mb X_l}\big(\int_0^1 \nabla_{\bo \theta}^2 L ((1 - s) \bo \theta^{\sf{aug}, t_0} + s \bo \theta^{t_0}) \mathrm d s \cdot  (\bo \theta^{\sf{aug}, t_0} - \bo \theta^{t_0}) \big)}_2  \\ 
                \leq & \max_{s\in[0,1]} \norm{\mc P_{\mb X_l}\big(\nabla_{\bo \theta}^2 L((1 - s) \bo \theta^{\sf{aug}, t_0} + s \bo \theta^{t_0}) \cdot  (\bo \theta^{\sf{aug}, t_0} - \bo \theta^{t_0}) \big)}_2 
                \\ 
                = & \max_{s\in[0,1]}\max_{\substack{\bo \Delta^{a} \in \bb R^{\sf{dim}(\bo \theta)},\\  \bo \Delta^a = \mc P_{\mb X_l}(\bo \Delta^a) ,
                \norm{\bo \Delta^a }_2 = 1
                }} Q(\bo \Delta^a, \bo \Delta^b,(1 - s) \bo \theta^{\sf{aug}, t_0} + s \bo \theta^{t_0}). 
            \end{align*}
            where $\bo \Delta^b \coloneqq \bo \theta^{\sf{aug}, t_0} - \bo \theta^{t_0}$ and the function $Q$ is defined as 
            \begin{align*}
                & Q(\bo \Delta^a, \bo \Delta^b, [\bo \zeta, \mb X, \mb Y]) \\ 
    & \coloneqq  \Bigip{\psi'(\mb 1_n \bo \zeta\t + \mb X \mb Y\t) \circ \mc P_{ \Omega}(\mb 1_n \bo \Delta^{a\top}_{\bo \zeta}+\bo \Delta^a_{\mb X}\mb Y\t + \mb X \bo \Delta^{a\top}_{\mb Y}),\\
    &\qquad \mc P_{ \Omega}(\mb 1_n \bo \Delta^{b\top}_{\bo \zeta}+\bo \Delta^b_{\mb X}\mb Y\t + \mb X \bo \Delta^{b\top}_{\mb Y})}\\
                & + 2 c_{\perp} \pi \frac{\sigma_r^*}{n}\bigip{\mb 1_n\t \bo \Delta_{\mb X}^a, \mb 1_n\t \bo \Delta_{\mb X}^b} + \big \langle \mc P_{\Omega}(\psi(\mb 1_n \bo \zeta\t + \mb X \mb Y\t) - \mb R), \bo \Delta^a_{\mb X}{\bo \Delta^b_{\mb Y}}\t + \bo \Delta^b_{\mb X}{\bo \Delta^a_{\mb Y}}\t \big \rangle . 
            \end{align*}
            Here, $\bo\Delta_{\bo\zeta}$ denotes the unscaled intercept component: under the bracket convention, the corresponding coordinates of $\bo\Delta=[\bo\Delta_{\bo\zeta},\bo\Delta_{\mb X},\bo\Delta_{\mb Y}]$ are $\sqrt{n/\sigma_r^*}\bo\Delta_{\bo\zeta}$. 

            For $s\in[0,1]$, write $\bo\theta_s\coloneqq
            (1-s)\bo\theta^{\sf{aug},t_0}+s\bo\theta^{t_0}$,
            $\mb X_s\coloneqq(1-s)\mb X^{\sf{aug},t_0}+s\mb X^{t_0}$,
            $\mb Y_s\coloneqq(1-s)\mb Y^{\sf{aug},t_0}+s\mb Y^{t_0}$,
            and $\mb D_{l,s}\coloneqq\diag\big(\psi'(\mc M(\bo\theta_s)_{l,j})
            \ind\{(l,j)\in\Omega\}\big)_{j\in[p]}$. Since
            $\bo\Delta^a$ is supported on the $\mb X_l$ coordinates, direct
            expansion gives
            \begin{align}
                &\max_{\substack{\bo\Delta^a=\mc P_{\mb X_l}(\bo\Delta^a),\\
                \norm{\bo\Delta^a}_2=1}}
                Q(\bo\Delta^a,\bo\Delta^b,\bo\theta_s) \\
                \leq{}&\bignorm{\mb Y_s\t\mb D_{l,s}\mb Y_s}
                \bignorm{\bo\Delta_{\mb X,l}^b}_2
                +\bignorm{\mb Y_s\t\mb D_{l,s}\bo\Delta_{\bo\zeta}^b}_2
                +\bignorm{\mb Y_s\t\mb D_{l,s}
                \bo\Delta_{\mb Y}^b\mb X_{s,l}}_2\nonumber\\
                &\quad+2c_{\perp}\frac{\pi\sigma_r^*}{n}
                \bignorm{\mb 1_n\t\bo\Delta_{\mb X}^b}_2
                +\bignorm{\Big[\mc P_\Omega\big(
                \psi(\mc M(\bo\theta_s))-\mb R\big)\Big]_{l,\cdot}
                \bo\Delta_{\mb Y}^b}_2.
                \label{eq: rowwise Hessian expansion at GD output}
            \end{align}

            We use leave-one-out conditioning for the first term in
            \eqref{eq: rowwise Hessian expansion at GD output}. Define
            \begin{align*}
                \mb O_{\sf{aug}}^{t_0,(l)}
                &\coloneqq\mb O(\bo\theta^{\sf{aug},t_0,(l)},
                \bo\theta^{t_0,(l)}),
                &\widetilde{\mb O}^{\sf{aug},t_0,(l)}
                &\coloneqq\mb O_{\sf{aug}}^{t_0,(l)}
                \widetilde{\mb O}^{t_0,(l)}.
            \end{align*}
            In the same common gauge as the full iterates, set
            \begin{align*}
                \bo\theta_s^{(l)}\coloneqq{}&(1-s)
                [\bo\zeta^{\sf{aug},t_0,(l)},
                \mb X^{\sf{aug},t_0,(l)}
                \widetilde{\mb O}^{\sf{aug},t_0,(l)},
                \mb Y^{\sf{aug},t_0,(l)}
                \widetilde{\mb O}^{\sf{aug},t_0,(l)}]\\
                &+s[\bo\zeta^{t_0,(l)},
                \mb X^{t_0,(l)}\widetilde{\mb O}^{t_0,(l)},
                \mb Y^{t_0,(l)}\widetilde{\mb O}^{t_0,(l)}],
            \end{align*}
            and let $\mb Y_s^{(l)}$ be its $\mb Y$ component. 
            
            Firstly, a leave-one-out analysis using the above quantities yields 
            $$\max_{s\in[0,1]}\bignorm{\mb Y_s\t\mb D_{l,s}\mb Y_s}
            \lesssim\sigma^2\pi\sigma_1^*. $$

            Here the sampling argument in the proof of Lemma~\ref{lemma: bounds for Hessian}, the RIC bounds \eqref{eq: rowwise consistency of parameters in geometry lemma}, and $\psi'\leq\sigma^2$ apply uniformly along the interpolation segment. For the second, third, and fifth terms in
            \eqref{eq: rowwise Hessian expansion at GD output}, we replace all
            interpolation and difference factors by their population-aligned
            leave-one-row-out counterparts. Matrix
            Bernstein controls the resulting terms, while the leave-one-out
            coupling bounds control the replacement errors. For the residual
            term, this is the decomposition used in
            \eqref{eq: decomposition of the Y difference term in the auxiliary series analysis}: the centered noise product is conditionally independent
            of row $l$, and the remaining terms follow from the Lipschitz
            continuity of $\psi$. Consequently, uniformly over $s\in[0,1]$ and
            $l\in[n]$, the right-hand side of
            \eqref{eq: rowwise Hessian expansion at GD output} is bounded by
            \begin{align*}
                &\sigma^2\pi\sigma_1^*\Big\{
                \dist\ti(\bo\theta^{\sf{aug},t_0},\bo\theta^{t_0})
                +\sqrt{\frac{\mu r}{n}}
                \distf(\bo\theta^{\sf{aug},t_0},\bo\theta^{t_0})\Big\}\\
    &\quad +c_{\perp}\frac{\pi\sigma_r^*}{n}
                \bignorm{\mb 1_n\t(\mb X^{\sf{aug},t_0}-\mb X^{t_0})}_2\\
                &\qquad\lesssim\sigma^2\kappa\kappa_\psi
                \sqrt{\frac{\mu\rho r}{d}}\pi\sigma_1^*
                \delta^{\sf{aug}}\log d.
            \end{align*}
            The first term in the last line follows specifically from the
            rowwise bound in Lemma~\ref{lemma: auxiliary sequence analysis};
            this is the source of $\kappa\kappa_\psi$. Taking a union bound
            over $l\in[n]$ gives failure probability $O(d^{-c})$.

        \end{itemize}

        Combining these bounds together with \eqref{eq: decomposition of the gradient at GD output} yields that 
        \begin{align}
            &  \max_{i\in[n]} \norm{\nabla_{\mb X_{i,\cdot}} L(\bo \theta^{t_0} )}_2 \lesssim \sigma^2 \kappa \kappa_\psi \sqrt{\frac{\mu r}{n}} \pi \sigma_1^* \delta^{\sf{aug}} \log d.
        \end{align} 
        The proof for $ \max_{l\in[p]} \bignorm{\nabla_{\mb Y^{\sf{app}}_{l,\cdot}} L(\bo \theta^{t_0} )}_2$ is similar to the one above and is hence omitted here. 

    \end{proof}

    \begin{proof}[Proof of Lemma~\ref{lemma: projection of difference onto specific directions}]
        We first exploit the Taylor expansion to obtain that 
        \begin{align*}
    & \nabla_{\bo \zeta} L(\bo \theta^{t_0}) = \nabla_{\bo \zeta} L(\bo \theta^*)  + \big( \big(\mc M(\bo \theta^{t_0}) - \mc M(\bo \theta^*) \big)\t \circ \mb H^{*\top} \big) \mb 1_n  + \big(\big(\mc M(\bo \theta^{t_0})\\
    &\quad - \mc M(\bo \theta^*) \big) \circ (\mb H^{t_0} - \mb H^*) \big)\t \mb 1_n\\
            & \nabla_{\mb X} L(\bo \theta^{t_0}) \mb O^{t_0} = \nabla_{\mb X} L(\bo \theta^*) + \big( \big(\mc M(\bo \theta^{t_0}) - \mc M(\bo \theta^*) \big) \circ \mb H^* \big) \mb Y^* \\ 
    & + \big( \big(\mc M(\bo \theta^{t_0}) - \mc M(\bo \theta^*) \big) \circ (\mb H^{t_0} - \mb H^*) \big) \mb Y^*  + \mc P_\Omega\big(\psi(\mc M(\bo \theta^{t_0})) - \mb R\big)\big(\mb Y^{t_0} \mb O^{t_0} - \mb Y^* \big)\\
    &\quad + 2 c_{\perp}\frac{\pi \sigma_r^*}{n} \mb 1_n \mb 1_n\t \mb X^{t_0} \mb O^{t_0} ,\\
            & \nabla_{\mb Y} L(\bo \theta^{t_0}) \mb O^{t_0} = \nabla_{\mb Y} L(\bo \theta^*)  + \big( \big(\mc M(\bo \theta^{t_0}) - \mc M(\bo \theta^*) \big)\t \circ \mb H^{*\top} \big) \mb X^*\\ 
            & + \big(\big(\mc M(\bo \theta^{t_0}) - \mc M(\bo \theta^*) \big) \circ (\mb H^{t_0} - \mb H^*) \big)\t \mb X^*  + \mc P_\Omega\big(\psi(\mc M(\bo \theta^{t_0})) - \mb R\big)\t \big(\mb X^{t_0} \mb O^{t_0} - \mb X^* \big), 
        \end{align*}
        where  $\mb H^{t_0} \coloneqq \mc P_\Omega\Big[\int_{0}^1 \psi'(\mb M^* + s(\mc M(\bo \theta^{t_0}) - \mb M^*)) \mathrm d s \Big]$ and $\mb H^* \coloneqq \mc P_\Omega(\psi'(\mc M(\bo \theta^*)))$. 

        Letting $\bo \Delta \coloneqq \big(\sqrt{n/\sigma_r^*}\bo \Delta_{\bo \zeta}\t, \sf{vec}(\bo \Delta_{\mb X})\t, \sf{vec}(\bo \Delta_{\mb Y})\t\big)\t$ with $\bo \Delta_{\bo \zeta} \coloneqq \bo \zeta^{t_0} - \bo \zeta^*$, $\bo \Delta_{\mb X} \coloneqq \mb X^{t_0} \mb O^{t_0} - \mb X^*$, $\bo \Delta_{\mb Y} \coloneqq  \mb Y^{t_0} \mb O^{t_0} - \mb Y^*$, $\bo\Delta_{\mb X^{\sf{app}}}\coloneqq(\mb 0,\bo\Delta_{\mb X})$, and $\bo\Delta_{\mb Y^{\sf{app}}}\coloneqq(\sqrt{n/\sigma_r^*}\bo\Delta_{\bo\zeta},\bo\Delta_{\mb Y})$, the above decompositions allow us to write
        \begin{align}
            & \nabla L(\bo \theta^*) = -\mb C \bo \Delta + \bo \Psi, \label{eq: linear equation for approximation}
        \end{align}
        where the matrix $\mb C$ collecting the effective coefficients and the residual vector $\bo \Psi$ are defined as 
        \begin{align*}
    & \mb C \coloneqq \underbrace{\sum_{(i,j) \in \Omega} \psi'(\mc M(\bo \theta^*))_{i,j} \sf{vec}\Big( \begin{matrix}
            \sqrt{\frac{\sigma_r^*}{n}}e_j^p\\ 
            e_i^n  {e^p_j}\t \mb Y^*\\ 
            e_j^p {e_i^n}\t \mb X^* 
        \end{matrix}\Big)
        \sf{vec}\Big( \begin{matrix}
            \sqrt{\frac{\sigma_r^*}{n}} e_j^p\\ 
            e_i^n {e_j^p}\t \mb Y^*\\ 
            e_j^p {e_i^n}\t \mb X^* 
        \end{matrix}\Big)\t}_{\eqqcolon \mb C_1} + \underbrace{ c_{\perp} \frac{\pi \sigma_r^*}{n} \nabla^2_{\bo \theta} \Big[\bignorm{\mb 1_n\t\mb X}_2^2 \Big]}_{\eqqcolon \mb C_2},\\
        & \bo \Psi \coloneqq \nabla L([\bo \zeta^{t_0}, \mb X^{t_0} \mb O^{t_0}, \mb Y^{t_0} \mb O^{t_0}])  - \bo \Psi_1 - \bo \Psi_2 - \bo \Psi_3 , \\ 
    & \bo \Psi_1 \coloneqq  \Big[ (\sqrt{\frac{\sigma_r^*}{n}} [\bo \Delta_{\mb X} \bo \Delta_{\mb Y}\t \circ \mb H^*]\t \mb 1_n )\t, \sf{vec}([\bo \Delta_{\mb X} {\bo \Delta_{\mb Y}\t }  \circ \mb H^*] \mb Y^*)\t,\\
    &\quad \sf{vec}([\bo \Delta_{\mb X} {\bo \Delta_{\mb Y}\t} \circ \mb H^{*}]\t \mb X^*)\t \Big]\t,\\
        & \bo \Psi_2 \coloneqq  \Big[\big(\sqrt{\frac{\sigma_r^*}{n}} \big(\big(\mc M(\bo \theta^{t_0}) - \mc M(\bo \theta^*) \big) \circ (\mb H^{t_0} - \mb H^*) \big)\t \mb 1_n \big)\t, \\ 
    & ~ \sf{vec}\big(\big( \big(\mc M(\bo \theta^{t_0}) - \mc M(\bo \theta^*) \big) \circ (\mb H^{t_0} - \mb H^*)\big) \mb Y^*\big)\t , \sf{vec}\big(\big(\big(\mc M(\bo \theta^{t_0}) - \mc M(\bo \theta^*) \big) \circ (\mb H^{t_0}\\
    &\quad - \mb H^*)\big)\t \mb X^*\big)\t   \Big]\t ,\\
    & \bo \Psi_3 \coloneqq \Big[ \mb 0, \sf{vec}\big( \big((\psi(\mc M(\bo \theta^{t_0})) - \mb R)\circ \bo \Omega\big)\big(\mb Y^{t_0} \mb O^{t_0} - \mb Y^* \big)\big)\t, \sf{vec} \big( \big((\psi(\mc M(\bo \theta^{t_0}))\\
    &\quad - \mb R)\circ \bo \Omega\big)\t \big(\mb X^{t_0} \mb O^{t_0} - \mb X^* \big)\big)\t\Big]\t.
        \end{align*}

        To proceed, we intend to control the deviation of $\mb C$ from its population counterpart $\bb E[\mb C]$. To this end, we regard $\mb C$ as a three-by-three block matrix and bound the deviation of each block separately. Specifically, we decompose
        \begin{align} 
    & \sum_{(i,j) \in \Omega} \psi'(\mc M(\bo \theta^*))_{i,j} \sf{vec}\Big( \begin{matrix}
            \sqrt{\frac{\sigma_r^*}{n}}e_j^p\\ 
            e_i^n  {e^p_j}\t \mb Y^*\\ 
            e_j^p {e_i^n}\t \mb X^* 
        \end{matrix}\Big)
        \sf{vec}\Big( \begin{matrix}
            \sqrt{\frac{\sigma_r^*}{n}} e_j^p\\ 
            e_i^n {e_j^p}\t \mb Y^*\\ 
            e_j^p {e_i^n}\t \mb X^* 
        \end{matrix}\Big)\t = \Big(\begin{matrix}
            \mb C^{(\bo \zeta, \bo \zeta)} & \mb C^{(\bo \zeta, \mb X)} & \mb C^{(\bo \zeta, \mb Y)} \\ 
            \mb C^{(\bo \zeta, \mb X)\top} & \mb C^{(\mb X, \mb X)} & \mb C^{(\mb X, \mb Y)} \\ 
            \mb C^{(\bo \zeta, \mb Y)\top} & \mb C^{(\mb X, \mb Y)\top} & \mb C^{(\mb Y, \mb Y)}
        \end{matrix} \Big).
        \end{align}
        Hence, we have
        \begin{align}
            & \norm{\mb C - \bb E[\mb C]} \leq \bignorm{\mb C^{(\bo \zeta, \bo \zeta)} - \bb E[\mb C^{(\bo \zeta, \bo \zeta)}]} + \bignorm{\mb C^{(\mb X, \mb X)} - \bb E[\mb C^{(\mb X, \mb X)}]} + \bignorm{\mb C^{(\mb Y, \mb Y)} - \bb E[\mb C^{(\mb Y, \mb Y)}]} \\ 
            & + 2\bignorm{\mb C^{(\bo \zeta, \mb X)} - \bb E[\mb C^{(\bo \zeta, \mb X)}]} + 2 \bignorm{\mb C^{(\bo \zeta, \mb Y)} - \bb E[\mb C^{(\bo \zeta, \mb Y)}]} + 2\bignorm{\mb C^{(\mb X, \mb Y)} - \bb E[\mb C^{(\mb X, \mb Y)}]}. 
        \end{align}
        In what follows, we illustrate the argument by deriving bounds for
$\bignorm{\mb C^{(\mb X,\mb X)}
-\bb E[\mb C^{(\mb X,\mb X)}]}$
and
$\bignorm{\mb C^{(\mb X,\mb Y)}
-\bb E[\mb C^{(\mb X,\mb Y)}]}$.
The remaining blocks can be handled analogously. 

\begin{itemize}
    \item Firstly, it is seen from the Bernstein inequality that, with probability at least $1 - O(d^{-c})$,  
    \begin{align}
        & \bignorm{\mb C^{(\mb X,\mb X)} -\bb E[\mb C^{(\mb X,\mb X)}]}
        \lesssim \sqrt{\pi}\sigma^2\sqrt{\frac{\mu r}{p}}\sigma_1^*\sqrt{\log d}
        +\sigma^2\frac{\mu r}{p}\sigma_1^*\log d\notag\\
        &\qquad\lesssim \sqrt{\pi}\sigma^2\sqrt{\frac{\mu r}{p}}\sigma_1^*\sqrt{\log d}.
    \end{align} 
    \item On the other hand, applying the matrix Bernstein inequality again reveals that 
    \begin{align}
        & \bignorm{\mb C^{(\mb X,\mb Y)} -\bb E[\mb C^{(\mb X,\mb Y)}]}
        \lesssim \sigma^2\sqrt{\pi r}\sqrt{\frac{\mu r}{n\wedge p}}\sigma_1^*\sqrt{\log d}
        +\sigma^2\frac{\mu r}{\sqrt{np}}\sigma_1^*\log d\notag\\
        &\qquad\lesssim \sigma^2\sqrt{\pi r}\sqrt{\frac{\mu\rho r}{d}}\sigma_1^*\sqrt{\log d}
    \end{align}
    with probability exceeding $1 - O(d^{-c})$. 
\end{itemize}

Consequently, one has with probability at least $1 - O(d^{-c})$ that
\begin{align}
    & \norm{\mb C - \bb E[\mb C]} \lesssim  \sigma^2 \sqrt{\pi r} \sqrt{\frac{\mu \rho r}{d}} \sigma_1^*\sqrt{\log d}.
    \label{eq: perturbation of C}
\end{align} 

Meanwhile, we introduce a direct product of the subspace associated with the mutual scaling ambiguity between $\mb X$ and $\mb Y$ and the one associated with rotational ambiguity, namely,
\begin{align}
& \mc V^{\sf{null}} \coloneqq \mc V^{\sf{scaling}} \oplus \mc V^{\sf{rotation}} \\ 
& \mc V^{\sf{scaling}}
\coloneqq
\big\{
\bo \Delta:
\bo \Delta_{\bo \zeta}=\mb 0,\
\bo \Delta_{\mb X}=\mb X^*\mb U,\
\bo \Delta_{\mb Y}=-\mb Y^*\mb U, \mb U = \mb U\t\big\}, \\
& \mc V^{\sf{rotation}} \coloneqq \big\{
\bo \Delta:
\bo \Delta_{\bo \zeta}=\mb 0,\
\bo \Delta_{\mb X}=\mb X^*\mb U,\
\bo \Delta_{\mb Y}=\mb Y^*\mb U, \mb U = -\mb U\t \big\}. 
\end{align} 
One can verify that $\mc V^{\sf{null}}$ is indeed contained in the null spaces of $\bb E[\mb C]$ and $\mb C$.
In addition, on the subspace $\mc V^{\sf{range}} \coloneqq (\mc V^{\sf{null}})^\perp$, the population matrix $\bb E[\mb C]$ is uniformly positive definite. Indeed, for every $\bo \Delta\in\mc V^{\sf{range}}$, one has 
\begin{align*}
    & \bo \Delta\t \bb E[\mb C] \bo \Delta \geq \pi \underline c_\psi \bignorm{\mb 1_n \bo \Delta_{\bo \zeta}\t + \mb X^*\bo \Delta_{\mb Y}\t + \bo \Delta_{\mb X}\mb Y^{*\top}}\fb^2 + c_{\perp} \frac{\pi \sigma_r^*}{n} \bignorm{\mb 1_n\t \bo \Delta_{\mb X}}_2^2\\ 
    & \geq   \pi \underline c_\psi \bignorm{ \mb X^*\bo \Delta_{\mb Y}\t + \bo \Delta_{\mb X}\mb Y^{*\top}}\fb^2 + 2 \pi \underline c_\psi \bigip{\mb 1_n \bo \Delta_{\bo \zeta}\t, \mb X^*\bo \Delta_{\mb Y}\t + \bo \Delta_{\mb X} \mb Y^{*\top}} + \pi \underline c_\psi \bignorm{\mb 1_n \bo \Delta_{\bo \zeta}\t}\fb^2\\
    &\quad +  c_{\perp} \frac{\pi \sigma_r^*}{n} \bignorm{\mb 1_n\t \bo \Delta_{\mb X}}_2^2\\
    & \geq    \pi \underline c_\psi \bignorm{ \mb X^*\bo \Delta_{\mb Y}\t + \bo \Delta_{\mb X}\mb Y^{*\top}}\fb^2 + 2 \pi \underline c_\psi \bigip{\bo \Delta_{\bo \zeta}\t \mb Y^{*}, \mb 1_n\t \bo \Delta_{\mb X} } + \pi \underline c_\psi \bignorm{\mb 1_n \bo \Delta_{\bo \zeta}\t}\fb^2\\
    &\quad +  c_{\perp} \frac{\pi \sigma_r^*}{n} \bignorm{\mb 1_n\t \bo \Delta_{\mb X}}_2^2\\
    & \geq   \pi \underline c_\psi \bignorm{ \mb X^*\bo \Delta_{\mb Y}\t + \bo \Delta_{\mb X}\mb Y^{*\top}}\fb^2 +  \pi \underline c_\psi \Bignorm{\sqrt{\frac{\underline c_\psi n}{c_\perp \sigma_r^*}}\bo \Delta_{\bo \zeta}\t \mb Y^{*} +  \sqrt{\frac{c_\perp \sigma_r^*}{\underline c_\psi n}}\mb 1_n\t \bo \Delta_{\mb X} }\fb^2\\
    & + \underbrace{\Big[ \pi \underline c_\psi n - \pi \kappa  \frac{ \underline c_\psi}{c_\perp} \underline c_\psi n\Big]}_{=c \pi \underline c_\psi n> 0} \norm{\bo \Delta_{\bo \zeta}}_2^2\\
    & \geq   c \pi\underline c_\psi \sigma_r^* \norm{\bo \Delta}_2^2
    \label{eq: lower bound on C}
\end{align*}
for some constant $c$ since ${\mb X^*}\t\bo\Delta_{\mb X}
=\bo\Delta_{\mb Y}\t\mb Y^*$ and $c_\perp \geq \kappa \underline c_\psi$. Combining \eqref{eq: perturbation of C} and \eqref{eq: lower bound on C}  yields that for every $\bo \Delta \in \mc V^{\sf{range}}$,  with probability at least $1 - O(d^{-c})$, 
\begin{align}
    & \bo \Delta\t \mb C \bo \Delta \geq \frac{1}{2} c\pi \underline c_\psi \sigma_r^* \norm{\bo \Delta}_2^2. 
\end{align}
In other words, the matrices $\mb C$ and $\bb E[\mb C]$ are both uniformly positive definite on the subspace $\mc V^{\sf{range}}$.

With these characterizations in place, we may left-multiply both sides of \eqref{eq: linear equation for approximation} by the Moore–Penrose pseudoinverse of $\mb C$ to see 
\begin{align}
    & \mc P_{\mc V^{\sf{range}}}(\bo \Delta ) = - \mb C^{\dagger} \nabla L(\bo \theta^*) + \mb C^{\dagger} \bo \Psi. 
\end{align}
For $j\in[p]$ and $k_1,k_2\in[r+1]$, define $\mb v_{\mb X,j,k_1,k_2}\in\bb R^{\operatorname{dim}(\bo\theta)}$ as the vector satisfying $\mb v_{\mb X,j,k_1,k_2}\t\bo\Delta=\big[\bo\Delta_{\mb X^{\sf{app}}}\t\operatorname{diag}\{\psi'(M^*_{i,j})\ind\{(i,j)\in\Omega\}\}_{i\in[n]}\mb X^{\sf{app},*}\big]_{k_1,k_2}$ for every $\bo\Delta$. There are at most $p(r+1)^2\lesssim d^3$ such triples. Accordingly, all fixed-triple concentration bounds below are invoked with exceptional probability $O(d^{-c-3})$ before taking a union bound. Notice that the direction vanishes when $k_1=1$, since the first column of $\bo\Delta_{\mb X^{\sf{app}}}$ is zero.

To facilitate analysis, we apply a decoupling trick to decouple the dependence between $\mb v_{\mb X,j,k_1,k_2}$ and $\bo \Delta$. Define the population-aligned leave-one-out error $\bo \Delta^{(-j)}$ from $\bo\theta^{t_0,(-j)}$ in the same way as $\bo\Delta$, using the rotation $\widetilde{\mb O}^{t_0,(-j)}$. For any two outputs $\bo\theta_1$ and $\bo\theta_2$ in the RIC, the first-order conditions of their population Procrustes rotations, together with Lemma~\ref{lemma: perturbation theory of optimal rotation under F norm}, give
\begin{align*}
    &\bignorm{[\bo\zeta_1-\bo\zeta_2,\mb X_1\mb O(\bo\theta_1,\bo\theta^*)-\mb X_2\mb O(\bo\theta_2,\bo\theta^*),\mb Y_1\mb O(\bo\theta_1,\bo\theta^*)-\mb Y_2\mb O(\bo\theta_2,\bo\theta^*)]}_2
    \lesssim \distf(\bo\theta_1,\bo\theta_2).
\end{align*}
Indeed, writing $\overline{\mb F}_k\coloneqq\mb F_k\mb O(\bo\theta_k,\bo\theta^*)$ and $\mb E_k\coloneqq\overline{\mb F}_k-\mb F^*$, population alignment makes $\overline{\mb F}_k\t\mb F^*$ symmetric. Hence $\operatorname{skew}(\overline{\mb F}_1\t\overline{\mb F}_2)=\operatorname{skew}(\mb E_1\t\mb E_2)$, whose Frobenius norm is at most $\norm{\mb E_1}\bignorm{\overline{\mb F}_1-\overline{\mb F}_2}\fb$. Lemma~\ref{lemma: perturbation theory of optimal rotation under F norm} and the RIC bounds then show that the relative rotation contributes at most $C\sqrt{\epsilon}\bignorm{\overline{\mb F}_1-\overline{\mb F}_2}\fb$, which can be absorbed for sufficiently small $\epsilon$. Invoking Theorem~\ref{thm: general GD} at the strengthened level $q=1$, Lemma~\ref{lemma: leave-one-out error} therefore reveals, simultaneously over $j\in[p]$, that with probability at least $1-O(d^{-c-1})$
\begin{align}
    & \bignorm{\bo \Delta - \bo \Delta^{(-j)}}_2 \lesssim \frac{\sigma \kappa \kappa_\psi \sqrt{\kappa \mu \rho r^2\log d}}{\underline c_\psi \sqrt{\pi \sigma_r^*}}.
\end{align}
Next, fix $j\in[p]$, replace the $j$th column of $\bo\Omega$ by an independent copy, rerun the full pipeline with the resulting mask $\bo\Omega'$, and define $\bo\Delta'$ from the primed output using its own population Procrustes rotation. Thus $\bo\Delta'$ is independent of $\bo\Omega_{\cdot,j}$ and hence of $\mb v_{\mb X,j,k_1,k_2}$. Invoke each stagewise guarantee for this primed pipeline with its concentration exponent increased by one. Coupling both outputs to the common leave-one-column-out sequence and applying the preceding population-alignment comparison then gives, for each fixed $j$, the following bound with exceptional probability $O(d^{-c-1})$:
\begin{equation}
    \bignorm{\bo \Delta - \bo \Delta'}_2 \lesssim \frac{\sigma \kappa \kappa_\psi \sqrt{\kappa \mu \rho r^2\log d}}{\underline c_\psi \sqrt{\pi \sigma_r^*}}.
\end{equation}
A union bound over $j\in[p]$ makes this coupling valid simultaneously with exceptional probability $O(d^{-c})$. This independence allows us to apply the scalar Bernstein inequality conditionally on $\bo\Delta'$. Denote $\mb v^*_{\mb X,j,k_1,k_2}\coloneqq\bb E[\mb v_{\mb X,j,k_1,k_2}]$. For each fixed $(j,k_1,k_2)$, invoke this conditional Bernstein bound with exceptional probability $O(d^{-c-3})$. Lemma~\ref{lemma: gd error controls} and a union bound over the at most $p(r+1)^2\lesssim d^3$ triples then give, uniformly over $j$, $k_1$, and $k_2$, with probability at least $1-O(d^{-c})$,
\begin{align}
    &\big|\mb v_{\mb X,j,k_1,k_2}\t\bo\Delta\big|
    \leq \big|\mb v_{\mb X,j,k_1,k_2}\t\bo\Delta'\big|
    +\norm{\mb v_{\mb X,j,k_1,k_2}}_2\norm{\bo\Delta-\bo\Delta'}_2\notag\\
    \leq{}& \big|\mb v^{*\top}_{\mb X,j,k_1,k_2}\bo\Delta'\big|
    +\sigma^2\sqrt{\frac{\mu r\sigma_1^*}{n\wedge p}}
    \Big\{\sqrt{\pi\log d}\,\bignorm{\bo\Delta'_{\mb X}}\fb
    +\log d\,\bignorm{\bo\Delta'_{\mb X}}\ti\Big\}\notag\\
    &\qquad+\norm{\mb v_{\mb X,j,k_1,k_2}}_2
    \norm{\bo\Delta-\bo\Delta'}_2\notag\\
    \lesssim{}& \big|\mb v^{*\top}_{\mb X,j,k_1,k_2}\bo\Delta'\big|
    +\sigma\kappa\kappa_\psi\mu r\sqrt{\rho\log d}
    \Big\{\sqrt{\log d}+\kappa\kappa_\psi\sqrt{\mu r}\Big\},
    \label{eq: decomposition of vXjk1k2 Delta}
\end{align}
where, on the same event,
\begin{align}
    & \norm{\mb v_{\mb X,j,k_1,k_2}}_2 \leq \sigma^2 \norm{\diag(\bo \Omega_{\cdot,j}) \mb X^{\sf{app},*}_{\cdot,k_2}}_2 \leq \sigma^2 \norm{\diag(\bo \Omega_{\cdot,j}) \mb X^{\sf{app},*}} \lesssim \sigma^2 \sqrt{\pi \sigma_1^*}.
    \label{eq: bound of vXjk1k2}
\end{align}
Here, the two conditional Bernstein contributions arise from the variance and the uniform bound on each summand, respectively. The last inequality follows from the final global and rowwise gradient-descent bounds, the independent-copy coupling bound, and the sampling condition, which absorbs the latter contribution involving $\{\pi(n\wedge p)\}^{-1/2}$.

Since $\bo \Delta$ and $\bo \Delta'$ follow the same distribution, we study the behavior of $\big|\mb v_{\mb X,j,k_1,k_2}^{*\top} \bo \Delta \big|$ for notational conciseness.
Applying the triangle inequality yields that
\begin{align}
    & \big|\mb v_{\mb X,j,k_1,k_2}^{*\top} \mc P_{\mc V^{\sf{range}}}(\bo \Delta ) \big| \leq \big|\mb v_{\mb X,j,k_1,k_2}^{*\top} \mb C^{\dagger} \nabla L(\bo \theta^*)\big|_2\\
    &\quad + \frac{2 \sigma^2 \pi (\sigma_1^*)^{\frac12} }{c \pi \underline c_\psi \sigma_r^* } \Big[ \norm{\nabla L([\bo \zeta^{t_0}, \mb X^{t_0} \mb O^{t_0}, \mb Y^{t_0} \mb O^{t_0}])}_2\\
    & + \norm{\bo \Psi_1} +  \norm{\bo \Psi_2} + \norm{\bo \Psi_3}  \Big].
    \label{eq: decomposition for the projection of Delta on vXjk1k2}
\end{align}
since $\bignorm{\mb v^*_{\mb X,j,k_1,k_2}}_2 \leq \sigma^2 \pi  \sqrt{\sigma_1^*}$.

Then it boils down to upper bounding the quantities on the right-hand side above, which is expanded in the following: 
\begin{itemize}
    \item

For each fixed $(j,k_1,k_2)$, applying the Bernstein inequality together with \eqref{eq: perturbation of C} yields the following bound with exceptional probability $O(d^{-c-3})$:
\begin{align}
    & \big|\mb v_{\mb X,j,k_1,k_2}^{*\top} \mb C^\dagger \nabla L(\bo \theta^*)\big| \\
    \leq & \big|\mb v_{\mb X,j,k_1,k_2}^{*\top} \bb E[\mb C]^\dagger \nabla L(\bo \theta^*)\big| + \bignorm{\mb v^*_{\mb X,j,k_1,k_2}}_2 \norm{\nabla L(\bo \theta^*)}_2 \bignorm{\bb E[\mb C]^{\dagger} - \mb C^{\dagger}} \\
    \lesssim & \frac{ \bignorm{\mb v^*_{\mb X,j,k_1,k_2} }_2}{ \pi \underline c_\psi \sigma_r^* } \sigma \sqrt{\pi \sigma_1^* r \log d} + \frac{ \bignorm{\mb v^*_{\mb X,j,k_1,k_2} }_2}{ \pi \underline c_\psi \sigma_r^* } B \sqrt{\frac{\mu \rho r}{d}}(\sigma_1^*)^{\frac12} \log d \\
    & + \frac{ \bignorm{\mb v^*_{\mb X,j,k_1,k_2} }_2}{ \pi \underline c_\psi \sigma_r^* } \frac{\sigma^2 \sqrt{\pi r} \sqrt{\frac{\mu \rho r}{d}} \sigma_1^*\sqrt{\log d}}{\pi \underline c_\psi \sigma_r^*} \sigma \sqrt{\pi d r\sigma_1^*} \\
     \lesssim &   \sigma \kappa_\psi^2  \kappa^2 r^{\frac32} \sqrt{\mu \rho\log d},
    \label{eq: the projection of the leave-one-out copy}
\end{align}
where we used the fact that with probability at least $1- O(d^{-c})$ 
\begin{equation}
    \norm{\nabla L(\bo \theta^*) }_2 \lesssim \sigma \sqrt{\pi d r\sigma_1^*}. 
\end{equation}
In particular, the factor-gradient blocks satisfy
$\norm{\mc P_\Omega(\mb E)\mb Y^*}\fb\leq\norm{\mc P_\Omega(\mb E)}\norm{\mb Y^*}\fb$
and its analogous $\mb X^*$-bound, where
$\norm{\mb X^*}\fb\vee\norm{\mb Y^*}\fb\lesssim\sqrt{r\sigma_1^*}$.
Taking a union bound over the at most $p(r+1)^2\lesssim d^3$ triples makes \eqref{eq: the projection of the leave-one-out copy} uniform with exceptional probability $O(d^{-c})$.

\item Moving forward, invoking Lemma~\ref{lemma: gradient norm control} implies that with probability at least $1- O(d^{-c})$
\begin{equation}
     \norm{\nabla L([\bo \zeta^{t_0}, \mb X^{t_0} \mb O^{t_0}, \mb Y^{t_0} \mb O^{t_0}])}_2 \leq \sigma^2\kappa\kappa_\psi \sqrt{\rho \mu r}\pi \sigma_1^* \delta^{\sf{aug}} \log d.  \label{eq: term 2 in the projection analysis}
\end{equation}

\item The terms $\bo\Psi_1$, $\bo\Psi_2$, and $\bo\Psi_3$ can all be controlled using the following observation. For any matrix $\mb P\in\bb R^{n\times p}$, $\mb A\in\bb R^{n\times r}$, and $\mb B,\mb C\in\bb R^{p\times r}$, we have
\begin{align*}
\norm{\mc P_\Omega\big((\mb A\mb B\t)\circ\mb P\big)\mb C}\fb
\leq
\max_{i\in[n],j\in[p]}|P_{i,j}|
\norm{\mb A}\fb
\max_{i\in[n]}\norm{\diag(\bo\Omega_{i,\cdot})\mb B}
\max_{i\in[n]}\norm{\diag(\bo\Omega_{i,\cdot})\mb C}.
\end{align*}
Indeed, let $\mb D_i\coloneqq\diag(\bo\Omega_{i,\cdot})$ and $\mb W_i\coloneqq\diag(P_{i,\cdot})$. The $i$th row on the left-hand side equals $\mb A_{i,\cdot}(\mb D_i\mb B)\t\mb W_i(\mb D_i\mb C)$, while $\norm{\mb W_i}\leq\max_{i,j}|P_{i,j}|$. Summing the squared rowwise bounds gives the display.
Quantities of the form $\norm{\diag(\bo\Omega_{i,\cdot})\mb B}$ can then be bounded through their leave-one-out counterparts. Specifically,
$$
\norm{\diag(\bo\Omega_{i,\cdot})\mb B}
\leq
\norm{\mb B-\mb B^{(i)}}
+
\norm{\diag(\bo\Omega_{i,\cdot})\mb B^{(i)}}.
$$
The second term is amenable to direct concentration because $\mb B^{(i)}$ is independent of $\bo\Omega_{i,\cdot}$. 

The final GD bounds and their leave-one-out counterparts imply
\begin{align*}
    &\bignorm{\bo\Delta_{\mb X}}\fb
    \vee\bignorm{\bo\Delta_{\mb Y}}\fb
    \lesssim\frac{\sigma\sqrt{\kappa\mu rd\log d}}
    {\underline c_\psi\sqrt{\pi\sigma_r^*}},\\
    &\max_{i\in[n]}\bignorm{\diag(\bo\Omega_{i,\cdot})\bo\Delta_{\mb Y}}\vee\max_{j\in[p]}
    \bignorm{\diag(\bo\Omega_{\cdot,j})\bo\Delta_{\mb X}}\\
    &\lesssim\frac{\sigma\kappa\kappa_\psi
    \sqrt{\kappa\mu\rho r^2\log d}}
    {\underline c_\psi\sqrt{\pi\sigma_r^*}}
    +\frac{\sigma\sqrt{\kappa\mu rd\log d}}
    {\underline c_\psi\sqrt{\sigma_r^*}}\lesssim\frac{\sigma\sqrt{\kappa\mu rd\log d}}
    {\underline c_\psi\sqrt{\sigma_r^*}},
\end{align*}
where the first term in the second bound is the actual-to-leave-one-out coupling error, while the second follows by concentration conditional on the leave-one-out factor. The last inequality follows from $\pi(n\wedge p)\gg\rho\mu^3r^2\kappa^3\kappa_\psi^4\xi^8$. The same concentration argument bounds the masked operator norms of the population factors by $\sqrt{\pi\sigma_1^*}$. Applying the preceding product inequality blockwise, and using the rank-$r$ conversion between the spectral and Frobenius norms, gives
\begin{equation}
    \norm{\bo\Psi_1}_2\lesssim
    \sigma^2 r\pi\sqrt{\sigma_1^*}
    \Big[\frac{\sigma\sqrt{\kappa\mu rd\log d}}
    {\underline c_\psi\sqrt{\pi\sigma_r^*}}\Big]^2.
    \label{eq: term 3 in the projection analysis}
\end{equation}

Next, write $\mb D\coloneqq\mc M(\bo\theta^{t_0})-\mc M(\bo\theta^*)$. For every $(i,j)\in[n]\times[p]$, its decomposition gives
\begin{align*}
    |D_{i,j}|
    \leq |\Delta_{\bo\zeta,j}|
    +\norm{\bo\Delta_{\mb X,i}}_2\norm{\mb Y_j^*}_2
    +\norm{\mb X_i^*}_2\norm{\bo\Delta_{\mb Y,j}}_2
    +\norm{\bo\Delta_{\mb X,i}}_2
    \norm{\bo\Delta_{\mb Y,j}}_2.
\end{align*}
The current rowwise bound and the incoherence of the population factors therefore yield
\begin{align*}
    \norm{\mb D}_{\max}
    &\lesssim\sqrt{\frac{\mu r\sigma_1^*}{n\wedge p}}
    \frac{\sigma\kappa\kappa_\psi
    \sqrt{\kappa\mu\rho r^2\log d}}
    {\underline c_\psi\sqrt{\pi\sigma_r^*}}
    +\Big[\frac{\sigma\kappa\kappa_\psi
    \sqrt{\kappa\mu\rho r^2\log d}}
    {\underline c_\psi\sqrt{\pi\sigma_r^*}}\Big]^2\\
    &\lesssim\frac{\sigma\kappa^2\kappa_\psi\mu r^{\frac32}
    \sqrt{\rho\log d}}
    {\underline c_\psi\sqrt{\pi(n\wedge p)}}.
\end{align*}
The intercept contribution is absorbed by the first term in the first line. The last inequality uses $\sigma_1^*/\sigma_r^*=\kappa$ and the RIC smallness of the rowwise GD error to absorb its square. Since $\sup_{x\in\mc D_\psi}|\psi''(x)|\lesssim\sigma^2$, the mean-value theorem gives
\begin{equation*}
    \norm{\mb H^{t_0}-\mb H^*}_{\max}
    \lesssim\frac{\sigma^2\sigma\kappa^2\kappa_\psi\mu r^{\frac32}
    \sqrt{\rho\log d}}
    {\underline c_\psi\sqrt{\pi(n\wedge p)}}.
\end{equation*}
Applying the preceding masked-product inequality to the four components of $\mb D$ consequently gives the following bound. Here, the first bracket is the global factor-error rate, while the two masked factor norms contribute the factor $r\pi\sigma_1^*$:
\begin{equation}
    \norm{\bo\Psi_2}_2\lesssim r\pi\sigma_1^*
    \Big[\frac{\sigma\sqrt{\kappa\mu rd\log d}}
    {\underline c_\psi\sqrt{\pi\sigma_r^*}}\Big]
    \Big[\frac{\sigma^2\sigma\kappa^2\kappa_\psi\mu r^{\frac32}
    \sqrt{\rho\log d}}
    {\underline c_\psi\sqrt{\pi(n\wedge p)}}\Big].
    \label{eq: term 4 in the projection analysis}
\end{equation}

Finally, decompose the score matrix as $\mc P_\Omega(\psi(\mc M(\bo\theta^{t_0}))-\mb R)=\mb D\circ\mb H^{t_0}-\mc P_\Omega(\mb E)$. Since $\norm{\mb H^{t_0}}_{\max}\leq\sigma^2$ and every component of $\mb D$ contains at least one estimation error, the same masked-product argument as for $\bo\Psi_1$ bounds the nonlinear contribution by the right-hand side of \eqref{eq: term 3 in the projection analysis}. For the noise contribution, response-noise concentration gives
\begin{align*}
    &\bignorm{\mc P_\Omega(\mb E)\bo\Delta_{\mb Y}}\fb
    \vee\bignorm{\mc P_\Omega(\mb E)\t\bo\Delta_{\mb X}}\fb
    \lesssim\sigma\sqrt{\pi d}
    \Big[\frac{\sigma\sqrt{\kappa\mu rd\log d}}
    {\underline c_\psi\sqrt{\pi\sigma_r^*}}\Big]\\
    &\quad \lesssim\sigma^2 r\pi\sqrt{\sigma_1^*}
    \Big[\frac{\sigma\sqrt{\kappa\mu rd\log d}}
    {\underline c_\psi\sqrt{\pi\sigma_r^*}}\Big]^2.
\end{align*}
The last inequality follows by substituting the global GD rate and using $\kappa,\kappa_\psi,\mu,r\geq1$.
Combining the nonlinear and noise contributions yields
\begin{equation}
    \norm{\bo\Psi_3}_2\lesssim
    \sigma^2 r\pi\sqrt{\sigma_1^*}
    \Big[\frac{\sigma\sqrt{\kappa\mu rd\log d}}
    {\underline c_\psi\sqrt{\pi\sigma_r^*}}\Big]^2.
    \label{eq: term 5 in the projection analysis}
\end{equation}
\end{itemize}

To sum up, plugging \eqref{eq: the projection of the leave-one-out copy}, \eqref{eq: term 2 in the projection analysis}, \eqref{eq: term 3 in the projection analysis}, \eqref{eq: term 4 in the projection analysis}, and \eqref{eq: term 5 in the projection analysis} into \eqref{eq: decomposition for the projection of Delta on vXjk1k2} finally concludes that
\begin{align}
    & \big|\mb v_{\mb X,j,k_1,k_2}^{*\top} \mc P_{\mc V^{\sf{range}}}(\bo \Delta )\big|\\
    & \lesssim   \sigma \kappa_\psi^2 \kappa^2 r^{\frac32} \sqrt{\mu \rho \log d} +  \frac{ \sigma^2 \pi \sqrt{\sigma_1^*}}{ \pi \underline c_\psi \sigma_r^* } \Big\{ \sigma^2\kappa\kappa_\psi \sqrt{\rho \mu r}\pi \sigma_1^* \delta^{\sf{aug}}\log d\\
    &\quad + \sigma^2 r \pi \sqrt{\sigma_1^*} \big[ \frac{\sigma \sqrt{\kappa \mu r d \log d}}{\underline c_\psi \sqrt{\pi \sigma_r^*}}\big]^2\\
    & +   r \pi \sigma_1^* \big[ \frac{\sigma \sqrt{\kappa \mu r d \log d}}{\underline c_\psi \sqrt{\pi \sigma_r^*}}\big]
     \big[\frac{\sigma^2\sigma\kappa^2\kappa_\psi\mu r^{\frac32}\sqrt{\rho\log d}}{\underline c_\psi\sqrt{\pi(n\wedge p)}}\big] \Big\} \\
    & \lesssim  \ \sigma\kappa_\psi^2\kappa^2r^{\frac32}
    \sqrt{\mu\rho\log d}
    +\sigma^2\kappa_\psi^2\kappa^2\mu r^2\log d
    \Big\{\frac{d}{\underline c_\psi\sigma_r^*}
    +\frac{\kappa_\psi\kappa^2\sqrt{\mu}\,r\rho}{\underline c_\psi}\Big\}\\
    &\quad +\sigma^2\pi\kappa_\psi^2\kappa^2
    \sqrt{\mu\rho r\sigma_1^*}\delta^{\sf{aug}}\log d
    \label{eq: bound for vXjk1k2 Prange Delta}
\end{align}
holds with probability at least $1- O(d^{-c})$, where we used
$\kappa_\psi=\sigma^2/\underline c_\psi$,
$\sigma_1^*=\kappa\sigma_r^*$, and $d/(n\wedge p)=\rho$.

On the other hand, we now control the components of $\bo\Delta$ in $\mc V^{\sf{scaling}}$ and $\mc V^{\sf{rotation}}$. By the optimality of the population Procrustes rotation $\mb O^{t_0}$, the matrix $\mb F^{*\top}\mb F^{t_0}\mb O^{t_0}$ is symmetric. Hence, for every skew-symmetric matrix $\mb K$,
\begin{align*}
    \bigip{\bo\Delta,[\mb 0,\mb X^*\mb K,\mb Y^*\mb K]}
    =\bigip{\mb F^{*\top}(\mb F^{t_0}\mb O^{t_0}-\mb F^*),\mb K}
    =0,
\end{align*}
{\emergencystretch=1.5em
which gives $\mc P_{\mc V^{\sf{rotation}}}(\bo\Delta)=\mb 0$.
For the scaling component, write $\mc P_{\mc V^{\sf{scaling}}}(\bo\Delta)=\allowbreak [\mb 0,\allowbreak \mb X^*\mb S,\allowbreak -\mb Y^*\mb S]$, where $\mb S=\allowbreak \mb S\t$, and let $\mb G^*\coloneqq\allowbreak {\mb X^*}\t\mb X^*=\allowbreak {\mb Y^*}\t\mb Y^*$. The normal equations for the Euclidean projection yield
$
    \mb G^*\mb S+\mb S\mb G^*
    =\allowbreak \operatorname{sym}\big({\mb X^*}\t\bo\Delta_{\mb X}
    -{\mb Y^*}\t\bo\Delta_{\mb Y}\big)
$, 
where $\operatorname{sym}(\mb A)\coloneqq\allowbreak (\mb A+\mb A\t)\allowbreak /2$.
Moreover,\par}

\begin{align*}
    &2\operatorname{sym}\big({\mb X^*}\t\bo\Delta_{\mb X}
    -{\mb Y^*}\t\bo\Delta_{\mb Y}\big)={\mb O^{t_0}}\t
    \big({\mb X^{t_0}}\t\mb X^{t_0}
    -{\mb Y^{t_0}}\t\mb Y^{t_0}\big)\mb O^{t_0}
    -{\bo\Delta_{\mb X}}\t\bo\Delta_{\mb X}\\
    &\quad +{\bo\Delta_{\mb Y}}\t\bo\Delta_{\mb Y}.
\end{align*}
Since $\lambda_{\min}(\mb G^*)=\sigma_r^*$, the Sylvester equation above implies
\begin{align*}
    \norm{\mc P_{\mc V^{\sf{scaling}}}(\bo\Delta)}_2
    \lesssim\frac{1}{\sqrt{\sigma_r^*}}
    \Big\{\bignorm{{\mb X^{t_0}}\t\mb X^{t_0}
    -{\mb Y^{t_0}}\t\mb Y^{t_0}}\fb
    +\bignorm{\bo\Delta_{\mb X}}\fb^2
    +\bignorm{\bo\Delta_{\mb Y}}\fb^2\Big\}.
\end{align*}
Lemma~\ref{lemma: balancedness and orthogonality}, \eqref{eq: final GD global error}, and the definition of $\delta^{\sf{aug}}$ therefore give, with probability at least $1-O(d^{-c})$,
\begin{align*}
    &\norm{\mc P_{\mc V^{\sf{scaling}}}(\bo\Delta)}_2
    \vee\norm{\mc P_{\mc V^{\sf{rotation}}}(\bo\Delta)}_2 \lesssim\delta^{\sf{aug}}
    +\frac{\sigma^2\kappa r d\log d}
    {\underline c_\psi^2\pi(\sigma_r^*)^{3/2}}.
\end{align*}
{\emergencystretch=1.5em
After multiplication by $\sigma^2\pi\allowbreak \sqrt{\sigma_1^*}$, the second term above is bounded by $\sigma^2\sigma^2\kappa^{3/2}r d\log d\allowbreak /(\underline c_\psi^2\sigma_r^*)$. Since $\kappa,\allowbreak \kappa_\psi,\allowbreak \mu,\allowbreak r\geq1$, this is dominated by the first component of the quadratic-error term in the display below.\par}

Using the orthogonal decomposition of $\bo\Delta$ and \eqref{eq: bound for vXjk1k2 Prange Delta} now gives
\begin{align*}
    & \big|\mb v_{\mb X,j,k_1,k_2}^{*\top} \bo \Delta\big| \leq \big|\mb v_{\mb X,j,k_1,k_2}^{*\top} \mc P_{\mc V^{\sf{range}}}(\bo \Delta )\big| + \sigma^2 \pi (\sigma_1^*)^{\frac12} \big[\norm{\mc P_{\mc V^{\sf{scaling}}}(\bo \Delta)}_2 + \norm{\mc P_{\mc V^{\sf{rotation}}}(\bo \Delta)}_2 \big]
     \\ 
    \lesssim &\ \sigma\kappa_\psi^2\kappa^2r^{\frac32}
    \sqrt{\mu\rho\log d}
    +\sigma^2\kappa_\psi^2\kappa^2\mu r^2\log d
    \Big\{\frac{d}{\underline c_\psi\sigma_r^*}
    +\frac{\kappa_\psi\kappa^2\sqrt{\mu}\,r\rho}{\underline c_\psi}\Big\}\\
    &\qquad+\sigma^2\pi\kappa_\psi^2\kappa^2
    \sqrt{\mu\rho r\sigma_1^*}\delta^{\sf{aug}}\log d.
    \label{eq: v*Xjk1k2 Delta bound}
\end{align*}

Finally, combining \eqref{eq: decomposition of vXjk1k2 Delta} and \eqref{eq: v*Xjk1k2 Delta bound} yields that with probability exceeding $1 - O(d^{-c})$
\begin{align}
    & \big|\mb v_{\mb X,j,k_1,k_2}^{\top} \bo \Delta\big| \\
    \lesssim &\ \sigma\kappa\kappa_\psi\mu r\sqrt{\rho\log d}
    \Big\{\sqrt{\log d}+\kappa\kappa_\psi\sqrt{\mu r}\Big\}
    +\sigma^2\kappa_\psi^2\kappa^2\mu r^2\log d
    \Big\{\frac{d}{\underline c_\psi\sigma_r^*}
    +\frac{\kappa_\psi\kappa^2\sqrt{\mu}\,r\rho}{\underline c_\psi}\Big\}\notag\\
    &\qquad+\sigma^2\pi\kappa_\psi^2\kappa^2
    \sqrt{\mu\rho r\sigma_1^*}\delta^{\sf{aug}}\log d.
\end{align}
Here, the second component of the leading term is
$\sigma\kappa^2\kappa_\psi^2\mu^{3/2}r^{3/2}\sqrt{\rho\log d}$,
which absorbs the first term in
\eqref{eq: v*Xjk1k2 Delta bound} because $\mu\geq1$.
Taking the maximum over $j$, $k_1$, and $k_2$ and using $\norm{\mb A}\leq(r+1)\max_{k_1,k_2}|A_{k_1,k_2}|\lesssim r\max_{k_1,k_2}|A_{k_1,k_2}|$ gives the first claimed bound.

It remains to establish the second projection. For $i\in[n]$, $k_1\in[r+1]$, and $k_2\in[r]$, define $\mb v_{\mb Y,i,k_1,k_2}\in\bb R^{\operatorname{dim}(\bo\theta)}$ by
\begin{align*}
    \mb v_{\mb Y,i,k_1,k_2}\t\bo\Delta
    \coloneqq\big[\bo\Delta_{\mb Y^{\sf{app}}}\t
    \diag\{\psi'(M^*_{i,j})\ind\{(i,j)\in\Omega\}\}_{j\in[p]}
    \mb Y^*\big]_{k_1,k_2}.
\end{align*}
For each fixed $i$, replace the $i$th row of $\bo\Omega$ by an independent copy, rerun the full pipeline, and define $\bo\Delta'$ from the resulting output using its population Procrustes rotation. The primed error is independent of $\bo\Omega_{i,\cdot}$. Moreover, the original and primed outputs can both be coupled to the same leave-one-row-out sequence. Invoking each stagewise guarantee with exceptional probability $O(d^{-c-1})$ therefore gives
\begin{align*}
    \norm{\bo\Delta-\bo\Delta'}_2
    \lesssim\frac{\sigma\kappa\kappa_\psi
    \sqrt{\kappa\mu\rho r^2\log d}}
    {\underline c_\psi\sqrt{\pi\sigma_r^*}}
\end{align*}
for each fixed $i$. A union bound makes this relation uniform over $i\in[n]$ with exceptional probability $O(d^{-c})$.

Write $\mb v^*_{\mb Y,i,k_1,k_2}\coloneqq\bb E[\mb v_{\mb Y,i,k_1,k_2}]$. Conditional on $\bo\Delta'$, scalar Bernstein applied to the $i$th-row sampling indicators yields
\begin{align*}
    &\big|\mb v_{\mb Y,i,k_1,k_2}\t\bo\Delta\big|\\
    \leq{}&\big|\mb v^{*\top}_{\mb Y,i,k_1,k_2}\bo\Delta'\big|
    +\sigma^2\sqrt{\frac{\mu r\sigma_1^*}{n\wedge p}}
    \Big\{\sqrt{\pi\log d}\,\bignorm{\bo\Delta'_{\mb Y^{\sf{app}}}}\fb
    +\log d\,\bignorm{\bo\Delta'_{\mb Y^{\sf{app}}}}\ti\Big\}\\
    &\qquad+\norm{\mb v_{\mb Y,i,k_1,k_2}}_2
    \norm{\bo\Delta-\bo\Delta'}_2\\
    \lesssim{}&\big|\mb v^{*\top}_{\mb Y,i,k_1,k_2}\bo\Delta'\big|
    +\sigma\kappa\kappa_\psi\mu r\sqrt{\rho\log d}
    \Big\{\sqrt{\log d}+\kappa\kappa_\psi\sqrt{\mu r}\Big\}.
\end{align*}
Here, each fixed-triple Bernstein bound is invoked with exceptional probability $O(d^{-c-3})$. In addition,
\begin{align*}
    \norm{\mb v_{\mb Y,i,k_1,k_2}}_2
    &\leq\sigma^2
    \norm{\diag(\bo\Omega_{i,\cdot})\mb Y^*_{\cdot,k_2}}_2
    \lesssim\sigma^2\sqrt{\pi\sigma_1^*},\\
    \norm{\mb v^*_{\mb Y,i,k_1,k_2}}_2
    &\leq\sigma^2\pi\norm{\mb Y^*_{\cdot,k_2}}_2
    \leq\sigma^2\pi\sqrt{\sigma_1^*}.
\end{align*}
Thus neither bound involves the intercept column of $\mb Y^{\sf{app},*}$.

The remaining argument follows the orthogonal decomposition already used for $\mb v_{\mb X,j,k_1,k_2}$. The only direction-specific term is $\mb v^{*\top}_{\mb Y,i,k_1,k_2}\mb C^\dagger\nabla L(\bo\theta^*)$. Applying scalar Bernstein to its row-indexed summands, using $\max_{j\in[p]}\norm{\mb Y_j^*}_2\leq\sqrt{\mu r\sigma_1^*/p}$, gives the analogue of \eqref{eq: the projection of the leave-one-out copy}, including the contribution $\sigma\kappa_\psi^2\kappa^2r^{3/2}\sqrt{\mu\rho\log d}$. This contribution is again absorbed by the second component of the leading term below. The perturbation bound for $\mb C^\dagger$, the bounds \eqref{eq: term 2 in the projection analysis}--\eqref{eq: term 5 in the projection analysis}, and the scaling and rotation controls in \eqref{eq: v*Xjk1k2 Delta bound} then yield
\begin{align*}
    \big|\mb v_{\mb Y,i,k_1,k_2}\t\bo\Delta\big|
    \lesssim{}&\sigma\kappa\kappa_\psi\mu r\sqrt{\rho\log d}
    \Big\{\sqrt{\log d}+\kappa\kappa_\psi\sqrt{\mu r}\Big\}\\
    &+\sigma^2\kappa_\psi^2\kappa^2\mu r^2\log d
    \Big\{\frac{d}{\underline c_\psi\sigma_r^*}
    +\frac{\kappa_\psi\kappa^2\sqrt{\mu}\,r\rho}{\underline c_\psi}\Big\}+\sigma^2\pi\kappa_\psi^2\kappa^2
    \sqrt{\mu\rho r\sigma_1^*}\delta^{\sf{aug}}\log d.
\end{align*}
Finally, there are at most $n(r+1)r\lesssim d^3$ triples, and
$\norm{\mb A}\leq\sqrt{r(r+1)}\max_{k_1,k_2}|A_{k_1,k_2}|
\lesssim r\max_{k_1,k_2}|A_{k_1,k_2}|$. Taking a union bound and maximizing over $i$ proves the second claimed projection with exceptional probability $O(d^{-c})$.

    \end{proof}

    \begin{proof}[Proof of Lemma~\ref{lemma: gradient replacement 1}]
        Write $\mb M_{\sf X}\coloneqq\mc M([\bo\zeta^*,\mb X^{t_0}\mb O^{t_0},\mb Y^*])$ and $\mb M_{t_0}\coloneqq\mc M(\bo\theta^{t_0})$. Rotational equivariance of the factor gradients gives the following identity for every $i\in[n]$:
        \begin{align}
            & \nabla_{\mb X_{i,\cdot}}L([\bo\zeta^*,\mb X^{t_0}\mb O^{t_0},\mb Y^*])
            -\nabla_{\mb X_{i,\cdot}}L(\bo\theta^{t_0})\mb O^{t_0} \\
            ={}&\big[\{\psi(\mb M_{\sf X})-\psi(\mb M_{t_0})\}\circ\bo\Omega\big]_{i,\cdot}\mb Y^*
            -\big[\{\psi(\mb M_{t_0})-\psi(\mc M(\bo\theta^*))-\mb E\}\circ\bo\Omega\big]_{i,\cdot}
            (\mb Y^{t_0}\mb O^{t_0}-\mb Y^*).
            \label{eq: proximity of single-row loss function}
        \end{align}
        The centering-penalty gradients cancel in this identity because the two left factors agree after alignment. Set $\bo\Delta_{\mb Y^{\sf{app}}}\coloneqq\mb Y^{\sf{app},t_0}\mb O^{\sf{app},t_0}-\mb Y^{\sf{app},*}$ and $\xi_{i,j}\coloneqq\int_0^1\psi'\big(\{(1-s)\mb M_{\sf X}+s\mb M_{t_0}\}_{i,j}\big)\mathrm ds$. Taylor's theorem and the fact that the desired $r$-vector consists of the last $r$ coordinates of the corresponding augmented vector yield
        \begin{align*}
            &\norm{\nabla_{\mb X_{i,\cdot}}L([\bo\zeta^*,\mb X^{t_0}\mb O^{t_0},\mb Y^*])
            -\nabla_{\mb X_{i,\cdot}}L(\bo\theta^{t_0})\mb O^{t_0}}_2
            \leq\alpha_1+\alpha_2+\alpha_3,
        \end{align*}
        where
        \begin{align*}
            \alpha_1&\coloneqq\norm{\mb X^{\sf{app},t_0}_{i,\cdot}\mb O^{\sf{app},t_0}\bo\Delta_{\mb Y^{\sf{app}}}\t
            \diag\big(\psi'(M^*_{i,j})\ind\{(i,j)\in\Omega\}\big)_{j\in[p]}\mb Y^*}_2,\\
            \alpha_2&\coloneqq\norm{\mb X^{\sf{app},t_0}_{i,\cdot}\mb O^{\sf{app},t_0}\bo\Delta_{\mb Y^{\sf{app}}}\t
            \diag\big(\{\xi_{i,j}-\psi'(M^*_{i,j})\}\ind\{(i,j)\in\Omega\}\big)_{j\in[p]}\mb Y^*}_2,\\
            \alpha_3&\coloneqq\norm{\big[\{\psi(\mb M_{t_0})-\psi(\mc M(\bo\theta^*))-\mb E\}\circ\bo\Omega\big]_{i,\cdot}
            (\mb Y^{t_0}\mb O^{t_0}-\mb Y^*)}_2.
        \end{align*}
        We now bound these three terms. All fixed-index concentration bounds below are invoked with exponent $c+2$, so a union bound over $i\in[n]$ and $j\in[p]$ preserves a final exceptional probability of order $O(d^{-c})$.
        \begin{itemize}
            \item Lemma~\ref{lemma: projection of difference onto specific directions} and $\bignorm{\mb X^{\sf{app},t_0}_{i,\cdot}}_2\lesssim\sqrt{\mu r\sigma_1^*/n}$ imply uniformly over $i$ that
            \begin{align*}
                \alpha_1\lesssim{}&\sqrt{\frac{\mu r\sigma_1^*}{n}}\Big[
                \sigma\kappa\kappa_\psi\mu r^2\sqrt{\rho\log d}
                \Big\{\sqrt{\log d}+\kappa\kappa_\psi\sqrt{\mu r}\Big\}\\
                &\qquad+\sigma^2\kappa_\psi^2\kappa^2\mu r^3\log d
                \Big\{\frac{d}{\underline c_\psi\sigma_r^*}
                +\frac{\kappa_\psi\kappa^2\sqrt{\mu}\,r\rho}{\underline c_\psi}\Big\}\\
                &\qquad+\sigma^2\pi\kappa_\psi^2\kappa^2 r
                \sqrt{\mu\rho r\sigma_1^*}\delta^{\sf{aug}}\log d\Big].
            \end{align*}
            The mixed component in the second line is controlled by the
            second term in the square bracket of the lemma statement.
            \item Since both endpoints of the Taylor path lie in the RIC, the mean-value theorem and the preceding max-norm analysis give
            \begin{align*}
                \max_{j\in[p]}|\xi_{i,j}-\psi'(M^*_{i,j})|
                \lesssim\frac{\sigma^2\sigma\kappa^2\kappa_\psi\mu r^{\frac32}
                \sqrt{\rho\log d}}{\underline c_\psi\sqrt{\pi(n\wedge p)}}.
            \end{align*}
            The row-leave-one-out argument used in
            \eqref{eq: term 4 in the projection analysis}, applied also to
            the scaled intercept column, gives
            \begin{align*}
                \max_{i\in[n]}\bignorm{\diag(\bo\Omega_{i,\cdot})
                \bo\Delta_{\mb Y^{\sf{app}}}}
                \lesssim\frac{\sigma\sqrt{\kappa\mu r d\log d}}
                {\underline c_\psi\sqrt{\sigma_r^*}}.
            \end{align*}
            Applying the masked-product bound used for \eqref{eq: term 4 in the projection analysis} therefore yields
            \begin{align}
                \alpha_2\lesssim \sqrt{\frac{\mu r}{n}}\pi\sigma_1^*
                \Big[\frac{\sigma\sqrt{\kappa\mu r d\log d}}
                {\underline c_\psi\sqrt{\pi\sigma_r^*}}\Big]
                \Big[\frac{\sigma^2\sigma\kappa^2\kappa_\psi\mu r^{\frac32}\sqrt{\rho\log d}}
                {\underline c_\psi\sqrt{\pi(n\wedge p)}}\Big].
            \end{align}
            Using $d=\rho(n\wedge p)$, $\sigma_1^*=\kappa\sigma_r^*$,
            and $\sigma^2=\kappa_\psi\underline c_\psi$, this becomes
            \begin{align*}
                \alpha_2\lesssim
                \frac{\sigma^2\kappa^{\frac72}\kappa_\psi^2\mu^2r^{\frac52}\rho
                \sqrt{\sigma_r^*}\log d}{\underline c_\psi\sqrt n}.
            \end{align*}

            \item For $\alpha_3$, let $\mb D\coloneqq\mb M_{t_0}-\mc M(\bo\theta^*)$. The score decomposition $\mc P_\Omega\{\psi(\mb M_{t_0})-\psi(\mc M(\bo\theta^*))-\mb E\}=\mb D\circ\mb H^{t_0}-\mc P_\Omega(\mb E)$ separates the nonlinear and noise terms. The same blockwise masked-product argument as in the proof of Lemma~\ref{lemma: projection of difference onto specific directions} bounds the nonlinear term by
            \begin{align*}
                \sigma^2 r\pi\sqrt{\sigma_1^*}\Big\{
                &\Big[\frac{\sigma\kappa\kappa_\psi\sqrt{\kappa\mu\rho r^2\log d}}
                {\underline c_\psi\sqrt{\pi\sigma_r^*}}\Big]
                \Big[\frac{\sigma\sqrt{\kappa\mu r d\log d}}
                {\underline c_\psi\sqrt{\pi\sigma_r^*}}\Big]+\sqrt{\frac{\mu r}{n}}
                \Big[\frac{\sigma\sqrt{\kappa\mu r d\log d}}
                {\underline c_\psi\sqrt{\pi\sigma_r^*}}\Big]^2\Big\}.
            \end{align*}
            To handle the noise term, use the leave-one-row-out decomposition
            \begin{align*}
                \mb Y^{t_0}\mb O^{t_0}-\mb Y^*
                ={}&\mb Y^{t_0}\mb O^{t_0}-\mb Y^{t_0,(i)}\widetilde{\mb O}^{t_0,(i)}
                +\mb Y^{t_0,(i)}\widetilde{\mb O}^{t_0,(i)}-\mb Y^*.
            \end{align*}
            The first difference is controlled by the leave-one-out coupling bound, while the second is independent of $\{E_{i,j}\Omega_{i,j}\}_{j\in[p]}$. The term arising from the uniform summand bound in the conditional Bernstein inequality is controlled by the rowwise leave-one-out error and absorbed under the sampling condition. Conditional Bernstein concentration and the final global GD bound then give
            \begin{align*}
                &\norm{\mc P_\Omega(\mb E)_{i,\cdot}(\mb Y^{t_0}\mb O^{t_0}-\mb Y^*)}_2\\
                \lesssim{}&\sigma\sqrt{\pi d}
                \Big[\frac{\sigma\kappa\kappa_\psi\sqrt{\kappa\mu\rho r^2\log d}}
                {\underline c_\psi\sqrt{\pi\sigma_r^*}}\Big]
                +\sigma\sqrt{\pi\log d}
                \Big[\frac{\sigma\sqrt{\kappa\mu r d\log d}}
                {\underline c_\psi\sqrt{\pi\sigma_r^*}}\Big].
            \end{align*}
            Consequently, $\alpha_3$ is bounded by the sum of the preceding nonlinear and noise displays.
        \end{itemize}

        Using $\sigma_1^*=\kappa\sigma_r^*$ and
        $\kappa,\kappa_\psi,\mu,r,\rho,\log d\geq1$, the terms combined in the
        lemma statement satisfy
        \begin{align*}
            &\sigma\kappa\kappa_\psi\mu r^2
            \sqrt{\frac{\mu r\rho\sigma_1^*\log d}{n\wedge p}}
            \Big\{\sqrt{\log d}+\kappa\kappa_\psi\sqrt{\mu r}\Big\}\\
            &\quad+\sqrt{\frac{\mu r\sigma_1^*}{n\wedge p}}
            \sigma^2\kappa_\psi^2\kappa^2\mu r^3\log d
            \frac{\kappa_\psi\kappa^2\sqrt\mu\,r\rho}{\underline c_\psi}
            +\frac{\sigma^2\kappa^{\frac72}\kappa_\psi^2\mu^2r^{\frac52}\rho
            \sqrt{\sigma_r^*}\log d}{\underline c_\psi\sqrt{n\wedge p}}\\
            &\lesssim\Big\{\frac{1}{\sqrt{\pi(n\wedge p)}}
            +\frac{\sigma\kappa^2\kappa_\psi\sqrt r}
            {\underline c_\psi\sqrt{\pi(n\wedge p)}}\Big\}
            \sigma\kappa_\psi^2\kappa^2r^4\mu^2\rho
            \sqrt{\pi\sigma_1^*}\log d.
        \end{align*}
        The remaining $\alpha_1$ and $\alpha_3$ contributions are controlled
        by the other terms in the lemma statement. Lemma~\ref{lemma: gradient norm control}, the orthogonality of $\mb O^{t_0}$, and the triangle inequality then establish the asserted row-gradient bound uniformly over $i\in[n]$.

{\emergencystretch=1.5em
        It remains to verify the augmented column-gradient bound. Set $\mb M_{\sf Y}\coloneqq\allowbreak \mc M([\bo\zeta^{t_0},\allowbreak \mb X^*,\allowbreak \mb Y^{t_0}\mb O^{t_0}])$ and $\bo\Delta_{\mb X^{\sf{app}}}\coloneqq\allowbreak \mb X^{\sf{app},t_0}\mb O^{\sf{app},t_0}-\mb X^{\sf{app},*}$. For every $j\in[p]$, rotational equivariance gives\par}

        \begin{align*}
            &\nabla_{\mb Y^{\sf{app}}_{j,\cdot}}L([\bo\zeta^{t_0},\mb X^*,\mb Y^{t_0}\mb O^{t_0}])
            -\nabla_{\mb Y^{\sf{app}}_{j,\cdot}}L(\bo\theta^{t_0})\mb O^{\sf{app},t_0}\\
            ={}&\big[\{\psi(\mb M_{\sf Y})-\psi(\mb M_{t_0})\}\circ\bo\Omega\big]_{\cdot,j}\t\mb X^{\sf{app},*}
            -\big[\{\psi(\mb M_{t_0})-\psi(\mc M(\bo\theta^*))-\mb E\}\circ\bo\Omega\big]_{\cdot,j}\t
            \bo\Delta_{\mb X^{\sf{app}}}.
        \end{align*}
        To make the column argument explicit, set
        $\xi_{\sf Y,i,j}\coloneqq\int_0^1\psi'\big(\{(1-s)\mb M_{\sf Y}+s\mb M_{t_0}\}_{i,j}\big)\mathrm ds$.
        Since the first column of $\bo\Delta_{\mb X^{\sf{app}}}$ is zero, Taylor's theorem gives
        \begin{align*}
            &\norm{\nabla_{\mb Y^{\sf{app}}_{j,\cdot}}L([\bo\zeta^{t_0},\mb X^*,\mb Y^{t_0}\mb O^{t_0}])
            -\nabla_{\mb Y^{\sf{app}}_{j,\cdot}}L(\bo\theta^{t_0})\mb O^{\sf{app},t_0}}_2
            \leq\beta_1+\beta_2+\beta_3,
        \end{align*}
        where
        \begin{align*}
            \beta_1&\coloneqq\norm{\big(0,\mb Y^{t_0}_{j,\cdot}\mb O^{t_0}\big)
            \bo\Delta_{\mb X^{\sf{app}}}\t
            \diag\big(\psi'(M^*_{i,j})\ind\{(i,j)\in\Omega\}\big)_{i\in[n]}\mb X^{\sf{app},*}}_2,\\
            \beta_2&\coloneqq\norm{\big(0,\mb Y^{t_0}_{j,\cdot}\mb O^{t_0}\big)
            \bo\Delta_{\mb X^{\sf{app}}}\t
            \diag\big(\{\xi_{\sf Y,i,j}-\psi'(M^*_{i,j})\}\ind\{(i,j)\in\Omega\}\big)_{i\in[n]}\mb X^{\sf{app},*}}_2,\\
            \beta_3&\coloneqq\norm{\big[\{\psi(\mb M_{t_0})-\psi(\mc M(\bo\theta^*))-\mb E\}\circ\bo\Omega\big]_{\cdot,j}\t
            \bo\Delta_{\mb X^{\sf{app}}}}_2.
        \end{align*}
        The first projection in Lemma~\ref{lemma: projection of difference onto specific directions} and
        $\norm{\mb Y^{t_0}_{j,\cdot}}_2\lesssim\sqrt{\mu r\sigma_1^*/p}$ yield
        \begin{align*}
            \beta_1\lesssim{}&\sqrt{\frac{\mu r\sigma_1^*}{p}}\Big[
            \sigma\kappa\kappa_\psi\mu r^2\sqrt{\rho\log d}
            \Big\{\sqrt{\log d}+\kappa\kappa_\psi\sqrt{\mu r}\Big\}\\
            &\qquad+\sigma^2\kappa_\psi^2\kappa^2\mu r^3\log d
            \Big\{\frac{d}{\underline c_\psi\sigma_r^*}
            +\frac{\kappa_\psi\kappa^2\sqrt{\mu}\,r\rho}{\underline c_\psi}\Big\}
            +\sigma^2\pi\kappa_\psi^2\kappa^2r
            \sqrt{\mu\rho r\sigma_1^*}\delta^{\sf{aug}}\log d\Big].
        \end{align*}
        The same max-norm estimate used for $\alpha_2$ gives
        \begin{align*}
            \beta_2\lesssim{}&\sqrt{\frac{\mu r}{p}}\pi\sigma_1^*
            \Big[\frac{\sigma\sqrt{\kappa\mu r d\log d}}
            {\underline c_\psi\sqrt{\pi\sigma_r^*}}\Big]
            \Big[\frac{\sigma^2\sigma\kappa^2\kappa_\psi\mu r^{\frac32}\sqrt{\rho\log d}}
            {\underline c_\psi\sqrt{\pi(n\wedge p)}}\Big].
        \end{align*}
        As above, this simplifies to
        \begin{align*}
            \beta_2\lesssim
            \frac{\sigma^2\kappa^{\frac72}\kappa_\psi^2\mu^2r^{\frac52}\rho
            \sqrt{\sigma_r^*}\log d}{\underline c_\psi\sqrt p}.
        \end{align*}
        Finally, the blockwise masked-product argument bounds the nonlinear part of $\beta_3$ by
        \begin{align*}
            \sigma^2 r\pi\sqrt{\sigma_1^*}\Big\{
            &\Big[\frac{\sigma\kappa\kappa_\psi\sqrt{\kappa\mu\rho r^2\log d}}
            {\underline c_\psi\sqrt{\pi\sigma_r^*}}\Big]
            \Big[\frac{\sigma\sqrt{\kappa\mu r d\log d}}
            {\underline c_\psi\sqrt{\pi\sigma_r^*}}\Big]+\sqrt{\frac{\mu r}{p}}
            \Big[\frac{\sigma\sqrt{\kappa\mu r d\log d}}
            {\underline c_\psi\sqrt{\pi\sigma_r^*}}\Big]^2\Big\}.
        \end{align*}
        For its noise part, decompose $\bo\Delta_{\mb X^{\sf{app}}}$ through the leave-one-column-out iterate. The leave-one-out component is independent of $\{E_{i,j}\Omega_{i,j}\}_{i\in[n]}$, whereas the discrepancy between the actual and leave-one-out iterates is controlled by the coupling bound. Conditional Bernstein concentration therefore gives
        \begin{align*}
            \norm{\mc P_\Omega(\mb E)_{\cdot,j}\t\bo\Delta_{\mb X^{\sf{app}}}}_2
            \lesssim{}&\sigma\sqrt{\pi d}
            \Big[\frac{\sigma\kappa\kappa_\psi\sqrt{\kappa\mu\rho r^2\log d}}
            {\underline c_\psi\sqrt{\pi\sigma_r^*}}\Big]
            +\sigma\sqrt{\pi\log d}
            \Big[\frac{\sigma\sqrt{\kappa\mu r d\log d}}
            {\underline c_\psi\sqrt{\pi\sigma_r^*}}\Big].
        \end{align*}
        Thus $\beta_3$ is bounded by the sum of the preceding two displays.
        The preceding combination inequality, with the row and column roles
        interchanged, controls the corresponding $\beta_1$ and $\beta_2$
        terms. The remaining $\beta_1$ and $\beta_3$ contributions are
        controlled by the other terms in the lemma statement.
        Combining these bounds with Lemma~\ref{lemma: gradient norm control}
        yields the asserted augmented column-gradient bound. All fixed-column concentration bounds are invoked with exceptional probability $O(d^{-c-2})$ before taking the union over $j\in[p]$.
    \end{proof}

    \begin{proof}[Proof of Lemma~\ref{lemma: gradient replacement 2}]
        Write $\bo\Delta_{\mb X}\coloneqq\mb X^{t_0}\mb O^{t_0}-\mb X^*$ and $\bo\Delta_{\mb Y^{\sf{app}}}\coloneqq\mb Y^{\sf{app},t_0}\mb O^{\sf{app},t_0}-\mb Y^{\sf{app},*}$. For $i\in[n]$ and $j\in[p]$, set
        \begin{align*}
            \xi_{\mb X,i,j}&\coloneqq\int_0^1\psi'\big(M^*_{i,j}+s(\bo\Delta_{\mb X})_{i,\cdot}{\mb Y^*_{j,\cdot}}\t\big)\,\mathrm ds,\\
            \xi_{\mb Y^{\sf{app}},i,j}&\coloneqq\int_0^1\psi'\big(M^*_{i,j}+s\mb X^{\sf{app},*}_{i,\cdot}(\bo\Delta_{\mb Y^{\sf{app}}})_{j,\cdot}\t\big)\,\mathrm ds.
        \end{align*}
        Since $\mb 1_n\t\mb X^*=\mb 0$, Taylor's theorem gives, for every $i\in[n]$,
        \begin{align*}
            &\nabla_{\mb X_{i,\cdot}}L([\bo\zeta^*,\mb X^{t_0}\mb O^{t_0},\mb Y^*])-\nabla_{\mb X_{i,\cdot}}L(\bo\theta^*)-\mb g_{\perp}^{t_0}-(\bo\Delta_{\mb X})_{i,\cdot}\mb H^*_{\mb X_i}\\
            ={}&(\bo\Delta_{\mb X})_{i,\cdot}\mb Y^{*\top}\diag\Big(\big\{\xi_{\mb X,i,j}-\psi'(M^*_{i,j})\big\}\ind\{(i,j)\in\Omega\}\Big)_{j\in[p]}\mb Y^*.
        \end{align*}
        The centering penalty does not depend on $\bo\zeta$ or $\mb Y$. Hence, the analogous expansion of the augmented column gradient gives, for every $j\in[p]$,
        \begin{align*}
            &\nabla_{\mb Y^{\sf{app}}_{j,\cdot}}L([\bo\zeta^{t_0},\mb X^*,\mb Y^{t_0}\mb O^{t_0}])-\nabla_{\mb Y^{\sf{app}}_{j,\cdot}}L(\bo\theta^*)-(\bo\Delta_{\mb Y^{\sf{app}}})_{j,\cdot}\mb H^*_{\mb Y^{\sf{app}}_{j}}\\
            ={}&(\bo\Delta_{\mb Y^{\sf{app}}})_{j,\cdot}\mb X^{\sf{app},*\top}\diag\Big(\big\{\xi_{\mb Y^{\sf{app}},i,j}-\psi'(M^*_{i,j})\big\}\ind\{(i,j)\in\Omega\}\Big)_{i\in[n]}\mb X^{\sf{app},*}.
        \end{align*}

        We work on the intersection of the uniform final-GD event, the uniform Hessian-variation event, and the uniform masked-Gram event. The first two events give
        \begin{align*}
            &\max_{i\in[n]}\norm{(\bo\Delta_{\mb X})_{i,\cdot}}_2\vee\max_{j\in[p]}\norm{(\bo\Delta_{\mb Y^{\sf{app}}})_{j,\cdot}}_2
            \lesssim\frac{\sigma\kappa\kappa_\psi\sqrt{\kappa\mu\rho r^2\log d}}{\underline c_\psi\sqrt{\pi\sigma_r^*}},\\
            &\max_{i\in[n],j\in[p]}\Big\{\big|\xi_{\mb X,i,j}-\psi'(M^*_{i,j})\big|\vee\big|\xi_{\mb Y^{\sf{app}},i,j}-\psi'(M^*_{i,j})\big|\Big\}
            \lesssim\frac{\sigma^2\sigma\kappa^2\kappa_\psi\mu r^{\frac32}\sqrt{\rho\log d}}{\underline c_\psi\sqrt{\pi(n\wedge p)}}.
        \end{align*}
        Uniform masked-factor concentration further yields
        \begin{align*}
            \max_{i\in[n]}\norm{\diag(\bo\Omega_{i,\cdot})\mb Y^*}^2
            \vee\max_{j\in[p]}\norm{\diag(\bo\Omega_{\cdot,j})\mb X^{\sf{app},*}}^2
            \lesssim\pi\sigma_1^*.
        \end{align*}
        Applying these three bounds to the two exact remainder identities proves the claimed rate. The fixed-index concentration inequalities are invoked with exponent $c+2$ before taking the unions over $i$ and $j$, so the exceptional probability remains $O(d^{-c})$.
    \end{proof}

    \begin{proof}[Proof of Lemma~\ref{lemma: gradient replacement 3}]
        We first verify the information-matrix bound used below. In the balanced parametrization,
        $\sigma_r({\mb Y^*}\t\mb Y^*)=\sigma_r^*$ and, by $\mb 1_n\t\mb X^*=\mb 0$,
        $\mb X^{\sf{app},*\top}\mb X^{\sf{app},*}=\diag(\sigma_r^*,{\mb X^*}\t\mb X^*)$.
        Moreover,
        $\max_{j\in[p]}\norm{\mb Y_j^*}_2^2\leq\mu r\sigma_1^*/p$ and
        $\max_{i\in[n]}\norm{\mb X^{\sf{app},*}_{i}}_2^2\leq\sigma_r^*/n+\mu r\sigma_1^*/n$.
        Matrix Chernoff concentration \citep{tropp2015introduction}, followed by a union bound over the $n+p$ masks, therefore gives
        \begin{align*}
            &\min_{i\in[n]}\sigma_r\big({\mb Y^*}\t\diag(\bo\Omega_{i,\cdot})\mb Y^*\big)
            \wedge\min_{j\in[p]}\sigma_{r+1}\big(\mb X^{\sf{app},*\top}
            \diag(\bo\Omega_{\cdot,j})\mb X^{\sf{app},*}\big)
            \geq\frac{1}{2}\pi\sigma_r^*
        \end{align*}
        with probability at least $1-O(d^{-c-1})$ under the stated sampling condition. Since $\psi'(M_{i,j}^*)\geq\underline c_\psi$, this event implies
        \begin{align}
            \min_{i\in[n]}\sigma_r(\mb H^*_{\mb X_i})
            \wedge\min_{j\in[p]}\sigma_{r+1}(\mb H^*_{\mb Y^{\sf{app}}_{j}})
            \geq\frac{1}{2}\underline c_\psi\pi\sigma_r^*.
            \label{eq: uniform invertibility of rowwise information matrices}
        \end{align}
        Since $\mb 1_n\t\mb X^*=\mb 0$, the centering penalty has zero gradient at $\bo\theta^*$. Consequently, the gradients at the truth satisfy the exact score identities
        $\nabla_{\mb X_{i,\cdot}}L(\bo\theta^*)=-\mc P_\Omega(\mb E)_{i,\cdot}\mb Y^*$ and
        $\nabla_{\mb Y^{\sf{app}}_{j,\cdot}}L(\bo\theta^*)=-\mc P_\Omega(\mb E)_{\cdot,j}\t\mb X^{\sf{app},*}$.

        It remains to account for the centering-gradient term in Lemma~\ref{lemma: gradient replacement 2}. The invariance of the likelihood under simultaneous intercept shifts and row shifts of $\mb X$ gives the exact identity
        \begin{align*}
            \mb 1_n\t\nabla_{\mb X}L([\bo\zeta,\mb X,\mb Y])-\nabla_{\bo\zeta}L([\bo\zeta,\mb X,\mb Y])\t\mb Y=2c_\perp\pi\sigma_r^*\mb 1_n\t\mb X.
        \end{align*}
        Consequently,
        \begin{align*}
            \mb g_\perp^{t_0}=\frac{1}{n}\Big\{\mb 1_n\t\nabla_{\mb X}L(\bo\theta^{t_0})\mb O^{t_0}-\nabla_{\bo\zeta}L(\bo\theta^{t_0})\t\mb Y^{t_0}\mb O^{t_0}\Big\}.
        \end{align*}
        The first coordinate of $\nabla_{\mb Y^{\sf{app}}_{j,\cdot}}L(\bo\theta^{t_0})$ is
        $\sqrt{\sigma_r^*/n}\,\nabla_{\zeta_j}L(\bo\theta^{t_0})$. Hence Lemma~\ref{lemma: gradient norm control} and $\norm{\mb Y^{t_0}}\lesssim\sqrt{\sigma_1^*}$ give
        \begin{align*}
            \frac{1}{n}\bignorm{\mb 1_n\t\nabla_{\mb X}L(\bo\theta^{t_0})}_2
            &\lesssim\sigma^2\kappa\kappa_\psi\sqrt{\frac{\mu\rho r}{d}}\,\pi\sigma_1^*\delta^{\sf{aug}}\log d,\\
            \frac{1}{n}\bignorm{\nabla_{\bo\zeta}L(\bo\theta^{t_0})\t\mb Y^{t_0}}_2
            &\lesssim\sqrt{\frac{\kappa p}{n}}\,\sigma^2\kappa\kappa_\psi
            \sqrt{\frac{\mu\rho r}{d}}\,\pi\sigma_1^*\delta^{\sf{aug}}\log d.
        \end{align*}
        Consequently, $p/n\leq\rho$ and \eqref{eq: uniform invertibility of rowwise information matrices} yield
        \begin{align*}
            \norm{\mb g_\perp^{t_0}}_2
            &\lesssim\sqrt{\kappa\rho}\,\sigma^2\kappa\kappa_\psi\sqrt{\frac{\mu\rho r}{d}}\,\pi\sigma_1^*\delta^{\sf{aug}}\log d,\\
            \max_{i\in[n]}\norm{\mb g_\perp^{t_0}(\mb H^*_{\mb X_i})^{-1}}_2
            &\lesssim\frac{\sigma^2}{\underline c_\psi}\kappa^{\frac52}\kappa_\psi
            \sqrt{\frac{\mu\rho r}{n\wedge p}}\,\delta^{\sf{aug}}\log d.
        \end{align*}
        We next record the algebra linking the surrogate gradients to the score expansions. Let $\mb r_{\mb X,i}$ and $\mb r_{\mb Y^{\sf{app}},j}$ denote the two residual vectors on the left-hand side of Lemma~\ref{lemma: gradient replacement 2}. Substituting the score identities into their definitions and rearranging gives the exact relations
        \begin{align*}
            &\mb X^{t_0}_{i,\cdot}\mb O^{t_0}-\mb X^*_{i,\cdot}
            -\mc P_\Omega(\mb E)_{i,\cdot}\mb Y^*(\mb H^*_{\mb X_i})^{-1}=\Big\{\nabla_{\mb X_{i,\cdot}}L([\bo\zeta^*,\mb X^{t_0}\mb O^{t_0},\mb Y^*])
            -\mb g_\perp^{t_0}\\
    &\quad -\mb r_{\mb X,i}\Big\}(\mb H^*_{\mb X_i})^{-1},\\
            &\mb Y^{\sf{app},t_0}_{j,\cdot}\mb O^{\sf{app},t_0}-\mb Y^{\sf{app},*}_{j,\cdot}
            -\mc P_\Omega(\mb E)_{\cdot,j}\t\mb X^{\sf{app},*}(\mb H^*_{\mb Y^{\sf{app}}_{j}})^{-1}=\Big\{\nabla_{\mb Y^{\sf{app}}_{j,\cdot}}L([\bo\zeta^{t_0},\mb X^*,\mb Y^{t_0}\mb O^{t_0}])\\
    &\quad -\mb r_{\mb Y^{\sf{app}},j}\Big\}(\mb H^*_{\mb Y^{\sf{app}}_{j}})^{-1}.
        \end{align*}
        By \eqref{eq: uniform invertibility of rowwise information matrices}, the maximum of the two residual norms is bounded by the sum of the right-hand sides in Lemmas~\ref{lemma: gradient replacement 1} and~\ref{lemma: gradient replacement 2}, together with $\norm{\mb g_\perp^{t_0}}_2$, divided by $\underline c_\psi\pi\sigma_r^*$.

        After this division, Lemma~\ref{lemma: gradient replacement 1} gives
        \begin{align*}
            &\Big\{\frac{1}{\sqrt{\pi(n\wedge p)}}
            +\frac{\sigma\kappa^2\kappa_\psi\sqrt r}
            {\underline c_\psi\sqrt{\pi(n\wedge p)}}
            +\frac{\sigma\sqrt{\mu r d}}
            {\underline c_\psi\sqrt\pi\,\sigma_r^*}\Big\}
            \frac{\sigma\kappa_\psi^2\kappa^{\frac52}r^4\mu^2\rho\log d}
            {\underline c_\psi\sqrt{\pi\sigma_r^*}}\\
    &\quad +\kappa_\psi^3\kappa^3\mu r^2
            \sqrt{\frac{\rho}{n\wedge p}}\,\delta^{\sf{aug}}\log d.
        \end{align*}
        The rates contributed by Lemma~\ref{lemma: gradient replacement 2}
        and $\mb g_\perp^{t_0}$ satisfy, respectively,
        \begin{align*}
            &\frac{\sigma^2\kappa^{\frac92}\kappa_\psi^3\mu^{\frac32}r^{\frac52}
            \rho\log d}{\underline c_\psi^2\pi
            \sqrt{(n\wedge p)\sigma_r^*}}
            \leq\frac{\sigma\kappa^2\kappa_\psi\sqrt r}
            {\underline c_\psi\sqrt{\pi(n\wedge p)}}
            \frac{\sigma\kappa_\psi^2\kappa^{\frac52}r^4\mu^2\rho\log d}
            {\underline c_\psi\sqrt{\pi\sigma_r^*}},\\
            &\frac{\sigma^2}{\underline c_\psi}\kappa^{\frac52}\kappa_\psi
            \sqrt{\frac{\mu\rho r}{n\wedge p}}\,\delta^{\sf{aug}}\log d
            \leq\kappa_\psi^3\kappa^3\mu r^2
            \sqrt{\frac{\rho}{n\wedge p}}\,\delta^{\sf{aug}}\log d.
        \end{align*}
        Finally, using $d=\rho(n\wedge p)$ and
        \eqref{eq: definition of delta^aug},
        \begin{align*}
            \kappa_\psi^3\kappa^3\mu r^2
            \sqrt{\frac{\rho}{n\wedge p}}\,\delta^{\sf{aug}}\log d
            ={}&\frac{\kappa^2\kappa_\psi^2\sqrt{\log d}}
            {\sqrt{\pi\mu r^3(n\wedge p)}}
            \frac{\sigma\kappa_\psi^2\kappa^{\frac52}r^4\mu^2\rho\log d}
            {\underline c_\psi\sqrt{\pi\sigma_r^*}}\\
            &+\frac{\kappa^5\kappa_\psi^5\mu r^2\rho\log d}
            {c_\perp\sqrt{\pi\sigma_r^*}}
            \big\{\sigma\sqrt{\kappa\mu r\log d}
            +\xi^{\sf{gd}}+\xi^{\sf{gd},\sf{loo}}\big\}\\
            \leq{}&\frac{\kappa^2\kappa_\psi^2\sqrt{\log d}}
            {\sqrt{\pi(n\wedge p)}}
            \frac{\sigma\kappa_\psi^2\kappa^{\frac52}r^4\mu^2\rho\log d}
            {\underline c_\psi\sqrt{\pi\sigma_r^*}}\\
            &+\frac{\kappa^5\kappa_\psi^5\mu r^2\rho\log d}
            {c_\perp\sqrt{\pi\sigma_r^*}}
            \big\{\sigma\sqrt{\kappa\mu r\log d}
            +\xi^{\sf{gd}}+\xi^{\sf{gd},\sf{loo}}\big\}.
        \end{align*}
        Combining these displays proves the stated rate.
    \end{proof}

    \subsection{Proofs of the Inferential Results}
    \label{subsec:proof-of-inference-thm}

    \begin{proof}[Proof of Theorem~\ref{thm: individual inference}]
    By Lemma~\ref{lemma: scaling equivalence}, it suffices to work in the
    balanced parametrization $\omega=1$. The information-normalized vectors
    under a general scaling differ from their balanced counterparts only by
    orthogonal transformations, which preserve the class of convex sets.

    We condition first on the sampling mask $\bo\Omega$. Let $\mc E_\Omega$ be
    the event on which \eqref{eq: uniform invertibility of rowwise information
    matrices} holds and, uniformly over all rows and columns,
    \begin{align*}
        \lambda_{\max}(\mb H^*_{\mb X_i})\vee
        \lambda_{\max}(\mb H^*_{\mb Y^{\sf{app}}_{j}})
        \lesssim \sigma^2\pi\sigma_1^*.
    \end{align*}
    Matrix Chernoff concentration gives
    $\bb P(\mc E_\Omega)\geq1-O(d^{-c})$.

    \paragraph{Step 1: Gaussian approximation of the linear terms.}
    Fix $i\in[n]$. Conditional on $\bo\Omega$, define the independent mean-zero
    vectors
    \begin{align*}
        \mb Z_{i,j}\coloneqq
        \Omega_{i,j}E_{i,j}(\mb H^*_{\mb X_i})^{-\frac12}\mb Y_j^*,
        \qquad j\in[p].
    \end{align*}
    Their conditional covariance matrices sum to $\mb I_r$. On
    $\mc E_\Omega$, moreover,
    \begin{align*}
        \max_{j:(i,j)\in\Omega}
        \bignorm{(\mb H^*_{\mb X_i})^{-\frac12}\mb Y_j^*}_2
        \lesssim\sqrt{\frac{\kappa\mu r}
        {\underline c_\psi\pi p}}.
    \end{align*}
    We implement the high-probability noise bound in
    Assumption~\ref{assumption: noise and sampling} by clipping and recentering
    each score separately at level $B$. This preserves independence. The
    original and clipped arrays can be coupled to agree except on an event of
    probability $O(d^{-c-10})$, while the resulting centering and covariance corrections
    are $o(1)$ by the exponential-family moment bounds. Suppressing these
    negligible corrections, Lemma~\ref{lemma:convex-set-gaussian-approximation}
    and $\sum_{j\in[p]}\bb E(\norm{\mb Z_{i,j}}_2^2\mid\bo\Omega)=r$ give
    \begin{align}
        &\sup_{A\in\mc C_r}\Big|
        \bb P\Big[\sum_{j\in[p]}\mb Z_{i,j}\in A\,\Big|\,\bo\Omega\Big]
        -\bb P_{\mb G\sim\mc N(\mb 0,\mb I_r)}[\mb G\in A]\Big|\notag\\
        &\quad\leq(42r^{\frac14}+16)
        \sum_{j\in[p]}\bb E\big(\norm{\mb Z_{i,j}}_2^3\mid\bo\Omega\big)
        +o(1)\notag\\
        &\quad\lesssim r^{\frac14}B
        \max_{j:(i,j)\in\Omega}
        \bignorm{(\mb H^*_{\mb X_i})^{-\frac12}\mb Y_j^*}_2
        \sum_{j\in[p]}\bb E\big(\norm{\mb Z_{i,j}}_2^2
        \mid\bo\Omega\big)+o(1)\notag\\
        &\quad\lesssim r^{\frac54}B
        \sqrt{\frac{\kappa\mu r}{\underline c_\psi\pi p}}+o(1).
        \label{eq: conditional gaussian approximation for factor score}
    \end{align}
    For a fixed $j\in[p]$, the same argument applied to the summands
    $\Omega_{i,j}E_{i,j}(\mb H^*_{\mb Y^{\sf{app}}_{j}})^{-1/2}
    \mb X^{\sf{app},*}_{i}$ gives the analogous bound with dimension $r+1$ and
    $p$ replaced by $n$. Assumption~\ref{assumption: signal strength and
    incoherence degree}(c) yields
    \begin{align*}
        (r+1)^{\frac54}B
        \sqrt{\frac{\kappa\mu r}
        {\underline c_\psi\pi(n\wedge p)}}
        \ll
        \frac{1}{\kappa^{\frac72}\kappa_\psi^{\frac72}
        \mu^{\frac12}r^{\frac14}\sqrt{\log d}}
        \lesssim(\log d)^{-\frac12}=o(1).
    \end{align*}
    Thus both conditional Gaussian-approximation errors are $o(1)$.

    \paragraph{Step 2: transfer through the linearization residuals.}
    Theorem~\ref{thm: linear approximations} and the upper information bound on
    $\mc E_\Omega$ imply that there exists a deterministic sequence
    $\varepsilon_d\lesssim(\log d)^{-2}$ such that, outside an event of
    probability $O(d^{-c})$, all information-normalized residuals are bounded
    by $\varepsilon_d$. For the fixed row $i$, write $
        \mb T_i\coloneqq(\mb H^*_{\mb X_i})^{\frac12}
        \big(\mb O^{t_0\top}\mb X_{i,\cdot}^{t_0\top}
        -\mb X_{i,\cdot}^{*\top}\big)$ and $\mb S_i\coloneqq\sum_{j\in[p]}\mb Z_{i,j}$.
    Hence $\norm{\mb T_i-\mb S_i}_2\leq\varepsilon_d$ on this event.
    For every $A\in\mc C_r$, the definition of $A^{\varepsilon_d}$ gives
    \begin{align}
        &\bb P(\mb T_i\in A\mid\bo\Omega)
        -\bb P_{\mb G\sim\mc N(\mb 0,\mb I_r)}(\mb G\in A)\notag\\
        \leq{}&\bb P\big(\norm{\mb T_i-\mb S_i}_2>
        \varepsilon_d\mid\bo\Omega\big)+\Big[\bb P(\mb S_i\in A^{\varepsilon_d}\mid\bo\Omega)
        -\bb P_{\mb G\sim\mc N(\mb 0,\mb I_r)}
        (\mb G\in A^{\varepsilon_d})\Big]\notag\\
        &+\bb P_{\mb G\sim\mc N(\mb 0,\mb I_r)}
        (\mb G\in A^{\varepsilon_d}\setminus A).
        \label{eq: upper transfer through enlarged convex set}
    \end{align}
    Since $A^{\varepsilon_d}$ is convex, the second term on the right-hand side
    is bounded by \eqref{eq: conditional gaussian approximation for factor
    score}. Lemma~\ref{lemma:convex-gaussian-boundary} bounds the last term by
    $(0.59r^{1/4}+0.21)\varepsilon_d$.

    The reverse comparison follows from the contracted set
    $A^{-\varepsilon_d}$. Indeed,
    \begin{align}
        &\bb P_{\mb G\sim\mc N(\mb 0,\mb I_r)}(\mb G\in A)
        -\bb P(\mb T_i\in A\mid\bo\Omega)\notag\\
        \leq{}&\bb P\big(\norm{\mb T_i-\mb S_i}_2>
        \varepsilon_d\mid\bo\Omega\big)+\Big[\bb P_{\mb G\sim\mc N(\mb 0,\mb I_r)}
        (\mb G\in A^{-\varepsilon_d})
        -\bb P(\mb S_i\in A^{-\varepsilon_d}\mid\bo\Omega)\Big]\notag\\
        &+\bb P_{\mb G\sim\mc N(\mb 0,\mb I_r)}
        (\mb G\in A\setminus A^{-\varepsilon_d}).
        \label{eq: lower transfer through contracted convex set}
    \end{align}
    The middle term is again bounded by \eqref{eq: conditional gaussian
    approximation for factor score}, and the last term is bounded by
    $(0.59r^{1/4}+0.21)\varepsilon_d$.

    Finally, the definition of $\xi$ and the sampling condition in
    Assumption~\ref{assumption: signal strength and incoherence degree}(a),
    together with $\pi(n\wedge p)\leq d$, imply
    $r=O\{\log d/\log\log d\}$. Consequently,
    $(r+1)^{1/4}\varepsilon_d=o(1)$. Integrating
    \eqref{eq: upper transfer through enlarged convex set} and
    \eqref{eq: lower transfer through contracted convex set} over
    $\bo\Omega$, and adding $\bb P(\mc E_\Omega^\complement)$, the
    linearization failure probability, and the clipping-coupling error, proves
    the first claim. Repeating the same comparison in dimension $r+1$ for
    $\mb Y^{\sf{app},t_0}_{j,\cdot}$ proves the second claim.
    \end{proof}
    \begin{proof}[Proof of Theorem~\ref{thm: individual asymptotic normality}]
    We work in the aligned coordinates used in the statement. Define
    \begin{align*}
        \bo\Delta_{\mb X}^{\sf est}
        &\coloneqq\hat{\mb X}^{t_0}-\mb X^*
        =\mb X^{t_0}\mb O^{t_0}-\mb X^*
        +\mb X^{t_0}\big(\hat{\mb O}^{t_0}-\mb O^{t_0}\big),\\
        \bo\Delta_{\mb Y}^{\sf est}
        &\coloneqq\hat{\mb Y}^{t_0}-\mb Y^*
        =\mb Y^{t_0}\mb O^{t_0}-\mb Y^*
        +\mb Y^{t_0}\big(\hat{\mb O}^{t_0}-\mb O^{t_0}\big).
    \end{align*}
    Since right multiplication by an orthogonal matrix preserves row norms,
    the triangle inequality, \eqref{eq: final GD rowwise errors}, and the
    signal-strength condition imply
    $\norm{\mb X^{t_0}}\ti\vee\norm{\mb Y^{t_0}}\ti
    \leq\{1+o(1)\}\big(\norm{\mb X^*}\ti\vee
    \norm{\mb Y^*}\ti\big)$. The population rotation condition in the theorem
    therefore makes the information-normalized row norms of the two
    rotation-estimation terms $o_{\bb P}(1)$.

    We next verify consistency of the plug-in information matrices. For example, adding and subtracting the information matrix formed from $\mb Y^*$ and $\psi'(\hat M_{i,j})$ gives
    \begin{align*}
        \bignorm{\hat{\mb H}_{\mb X_i}-\mb H^*_{\mb X_i}}
        \lesssim{}& \sigma^2\Big\{
        \bignorm{\diag(\bo\Omega_{i,\cdot})\bo\Delta_{\mb Y}^{\sf est}}^2
        +\sqrt{\pi\sigma_1^*}\,
        \bignorm{\diag(\bo\Omega_{i,\cdot})\bo\Delta_{\mb Y}^{\sf est}}
        \Big\}\notag\\
        &+\sigma^2\pi\sigma_1^*
        \max_{j\in[p]}|\hat M_{i,j}-M^*_{i,j}|.
    \end{align*}
    The analogous bound for $\hat{\mb H}_{\mb Y^{\sf{app}}_{j}}-\mb H^*_{\mb Y^{\sf{app}}_{j}}$ replaces $\bo\Delta_{\mb Y}^{\sf est}$ by
    $(\mb 0,\bo\Delta_{\mb X}^{\sf est})$ and the row mask by the column mask. The strengthened signal and effective-sample-size conditions used in the proof of Theorem~\ref{thm: parameter-explicit linear approximations} imply
    \begin{align*}
        &\frac{\kappa\mu r\rho^{\frac12}\sqrt{d/\pi}\,
        \sigma^2(\log d)^2}
        {\underline c_\psi^{\frac32}\sigma_r^*}=o(1),\qquad \frac{\kappa^{\frac52}\kappa_\psi^2\mu^{\frac32}r^{\frac74}
        \rho\sigma^2(\log d)^2}
        {\sqrt{\pi(n\wedge p)\underline c_\psi}}=o(1).
    \end{align*}
    The final-GD bounds and masked-Gram concentration therefore imply, uniformly over $i$ and $j$,
    \begin{align}
        &\bignorm{(\mb H^*_{\mb X_i})^{-\frac12}
        \big(\hat{\mb H}_{\mb X_i}-\mb H^*_{\mb X_i}\big)
        (\mb H^*_{\mb X_i})^{-\frac12}}=o_{\bb P}(1),\nonumber\\
        &\bignorm{(\mb H^*_{\mb Y^{\sf{app}}_{j}})^{-\frac12}
        \big(\hat{\mb H}_{\mb Y^{\sf{app}}_{j}}-\mb H^*_{\mb Y^{\sf{app}}_{j}}\big)
        (\mb H^*_{\mb Y^{\sf{app}}_{j}})^{-\frac12}}=o_{\bb P}(1).
        \label{eq: relative consistency of plug-in information}
    \end{align}
    In particular, inversion and taking principal submatrices preserve the corresponding relative $o_{\bb P}(1)$ errors.

    For the left factor, \eqref{eq: relative consistency of plug-in information} implies
    \begin{align*}
        \big(\hat{\mb X}_{i,\cdot}^{t_0}-\mb X_{i,\cdot}^*\big)
        \hat{\mb H}_{\mb X_i}
        \big(\hat{\mb X}_{i,\cdot}^{t_0}-\mb X_{i,\cdot}^*\big)\t
        =
        \bignorm{(\mb H^*_{\mb X_i})^{\frac12}
        \big(\mb O^{t_0\top}\mb X_{i,\cdot}^{t_0\top}
        -\mb X_{i,\cdot}^{*\top}\big)}_2^2+o_{\bb P}(1).
    \end{align*}
    Theorem~\ref{thm: individual inference} and Slutsky's theorem show that this converges to $\chi_r^2$. The augmented result in Theorem~\ref{thm: individual inference}, together with the second relation in \eqref{eq: relative consistency of plug-in information}, gives the standard normal limit for the first coordinate and the $\chi_r^2$ limit for the last $r$ coordinates standardized by
    $[(\hat{\mb H}_{\mb Y^{\sf{app}}_{j}}^{-1})_{2:r+1,2:r+1}]^{-1/2}$.
    This proves the claims for $\bo\zeta$ and $\mb Y$ as well as that for $\mb X$. Finally, continuity of the normal and chi-square distribution functions makes the coverage errors uniform over $\alpha\in(0,1)$.
    \end{proof}
    \begin{proof}[Proof of Theorem~\ref{thm: simultaneous inference}]
    Multiplying every coefficient in one of the two statistics by a common
    positive constant scales both that statistic and its bootstrap critical
    value by the same amount. Lemma~\ref{lemma: scaling equivalence} therefore
    reduces the proof to the balanced parametrization $\omega=1$ without
    changing either ratio of contrast norms. We prove the result
    for $\mc T_{\mb X}$ and then indicate the augmented-column analogue. There
    are three probability layers. We first condition on $\bo\Omega$, next take
    probability over the score noise $\mb E$, and finally use $\bb P_\xi$ for
    the Gaussian multipliers conditional on the observed data.

    \paragraph{Step 1: ideal score and multiplier statistics.}
    For $j\in[p]$ and $l\in[L]$, define
    \begin{align*}
        Z_{\mb X,j,l}
        \coloneqq
        \sum_{(i,k)\in\mc S_{\mb X,l}}
        a_{\mb X,l,i,k}\Omega_{i,j}E_{i,j}
        \big[\mb Y_j^*\mb H_{\mb X_i}^{*-1}\big]_k.
    \end{align*}
    Conditional on $\bo\Omega$, the vectors
    $\mb Z_{\mb X,j}\coloneqq(Z_{\mb X,j,l})_{l\in[L]}$ are independent and mean zero. For each $l\in[L]$, set
    \begin{align*}
        v_{\mb X,l}^2
        \coloneqq\sum_{j\in[p]}\bb E\big(Z_{\mb X,j,l}^2\mid\bo\Omega\big),
        \qquad
        \hat v_{\mb X,l}^2
        \coloneqq\sum_{j\in[p]}Z_{\mb X,j,l}^2.
    \end{align*}
    The ideal target and its multiplier analogue are
    \begin{align*}
        \mc T_{\mb X}^{\sharp}
        &\coloneqq\max_{l\in[L]}\sum_{j\in[p]}Z_{\mb X,j,l},\qquad 
        \mc T_{\mb X}^{\sharp[b]}
        \coloneqq\max_{l\in[L]}\sum_{j\in[p]}
        \xi_{\mb Y,j}^{[b]}Z_{\mb X,j,l}.
    \end{align*}
    This expression uses the rowwise information matrix $\mb H_{\mb X_i}^*$, not the singular Hessian of the full factorized loss.

    Let $\mc E^{\sf good}$ contain the uniform information bounds, the rowwise linear approximations, the final-GD bounds, and the rotation and weight bounds in the theorem. On this event, the conditional variance of every active contrast is bounded below and above as
    \begin{align*}
        \frac{\underline a_{\mb X}^2}{C\sigma^2\pi\sigma_1^*}
        \leq
        v_{\mb X,l}^2
        \leq
        \frac{C\bar a_{\mb X}^2}{\underline c_\psi\pi\sigma_r^*}.
    \end{align*}
    Assumption~\ref{assumption: link function} directly supplies the exponential moment
    needed for the Gaussian approximation. Indeed, the exponential-family cumulant identity and
    Assumption~\ref{assumption: link function} give, for all sufficiently small
    $|t|$,
    \begin{equation*}
        \log\bb E\{\exp(tE_{i,j})\mid M_{i,j}^*\}
        =\sum_{q=2}^{\infty}\frac{\psi^{(q-1)}(M_{i,j}^*)}{q!}t^q
        \lesssim\frac{\sigma^2 t^2}{1-C|t|}.
    \end{equation*}
    Thus the sub-exponential Orlicz norm satisfies
    $\bignorm{E_{i,j}}_{\psi_1}\lesssim1+\sigma$ uniformly.
    After standardizing the contrasts, choose a sufficiently large constant
    $C$ and define the deterministic bound
    \begin{align}
        \max_{j\in[p],\,l\in[L]}
        \bignorm{Z_{\mb X,j,l}/v_{\mb X,l}}_{\psi_1}
        \leq B_{\sf hd}
        \coloneqq C(1+\sigma)\kappa_\psi\kappa^{\frac32}
        \sqrt{\frac{\mu r}{\underline c_\psi\pi(n\wedge p)}}.
        \label{eq: simultaneous standardized envelope}
    \end{align}
    The exponential-family cumulant identity and $L_\psi=O(1)$ give
    \begin{align*}
        \bb E E_{i,j}^4
        &=\psi'''(M_{i,j}^*)+3\{\psi'(M_{i,j}^*)\}^2
        \lesssim\kappa_\psi(1+\sigma^2)\bb E E_{i,j}^2.
    \end{align*}
    Grouping terms with the same row index in $Z_{\mb X,j,l}$,
    conditional independence, the bounded support size, incoherence, and
    the information-matrix bounds therefore imply
    \begin{align*}
        \bb E\Big[\Big|\frac{Z_{\mb X,j,l}}{v_{\mb X,l}}\Big|^4
        \,\Big|\,\bo\Omega\Big]
        &\lesssim
        \frac{(1+\sigma^2)\kappa^2\kappa_\psi^2\mu r}
        {\underline c_\psi\pi(n\wedge p)}
        \bb E\Big[\frac{Z_{\mb X,j,l}^2}{v_{\mb X,l}^2}
        \,\Big|\,\bo\Omega\Big]\\
        &\lesssim B_{\sf hd}^2
        \bb E\Big[\frac{Z_{\mb X,j,l}^2}{v_{\mb X,l}^2}
        \,\Big|\,\bo\Omega\Big].
    \end{align*}
    Using
    $\sum_{j\in[p]}\bb E(Z_{\mb X,j,l}^2\mid\bo\Omega)/v_{\mb X,l}^2=1$
    then yields
    \begin{align*}
        \frac{1}{p}\sum_{j\in[p]}
        \bb E\Big[\Big|\frac{\sqrt p\,Z_{\mb X,j,l}}
        {v_{\mb X,l}}\Big|^4\,\Big|\,\bo\Omega\Big]
        &\lesssim pB_{\sf hd}^2
        \sum_{j\in[p]}\bb E\Big[
        \frac{Z_{\mb X,j,l}^2}{v_{\mb X,l}^2}
        \,\Big|\,\bo\Omega\Big]
        =pB_{\sf hd}^2.
    \end{align*}
    Hence the exponential-moment, variance, and fourth-moment conditions of
    Lemma~\ref{lemma:gaussian-multiplier-max-approximation} hold with
    $b_1=1$, $b_2=O(1)$, and $B_p=C\sqrt p\,B_{\sf hd}\geq1$, after
    enlarging the universal constant $C$.
    The same truncation and Bernstein argument used for the multiplier
    approximation gives, conditionally on $\bo\Omega$ and uniformly over
    $l\in[L]$,
    \begin{align}
        \max_{l\in[L]}
        \Big|\frac{\hat v_{\mb X,l}^2}{v_{\mb X,l}^2}-1\Big|
        \lesssim
        B_{\sf hd}\sqrt{\log(Ld)}
        +B_{\sf hd}^2\log^2(Ld)
        =o_{\bb P}(1).
        \label{eq: uniform empirical contrast variance}
    \end{align}
    In particular, with probability $1-o(1)$, all realized multiplier
    variances are bounded below by
    $\underline a_{\mb X}^2/(2C\sigma^2\pi\sigma_1^*)$.
    Moreover, for every $x\in\bb R$,
    \begin{align*}
        \{\mc T_{\mb X}^{\sharp}\leq x\}
        &\equiv
        \Big\{\Big(\sum_{j\in[p]}
        \frac{Z_{\mb X,j,l}}{v_{\mb X,l}}\Big)_{l\in[L]}
        \in\prod_{l\in[L]}(-\infty,x/v_{\mb X,l}]\Big\},
    \end{align*}
    and the same identity holds for the multiplier statistic. Hence the
    hyperrectangle Gaussian multiplier approximation used in
    Lemma~\ref{lemma:gaussian-multiplier-max-approximation} applies to the
    vectors $(\sqrt p\,Z_{\mb X,j,l}/v_{\mb X,l})_{l\in[L]}$; with the
    preceding choice, $B_p^2/p\asymp B_{\sf hd}^2$.
    Consequently, uniformly over $\alpha\in(0,1)$, the ideal multiplier critical value $\mc G_{\mb X,1-\alpha}^{\sharp}$ obeys
    \begin{align}
        \Big|\bb P\Big[\mc T_{\mb X}^{\sharp}
        \leq\mc G_{\mb X,1-\alpha}^{\sharp}\,\Big|\,\bo\Omega\Big]
        -(1-\alpha)\Big|
        \lesssim
        \Big\{\frac{(1+\sigma^2)\kappa_\psi^2\kappa^3\mu r
        \log^5(Ld)}{\underline c_\psi\pi(n\wedge p)}\Big\}^{\frac14}
        =o(1).
        \label{eq: ideal simultaneous approximation}
    \end{align}
    The last equality follows from $\log L\lesssim\log d$,
    $\kappa_\psi=\sigma^2/\underline c_\psi$, and
    \eqref{eq: effective sample sizes from linear approximation}; indeed, the
    expression inside braces is bounded by
    \begin{equation*}
        \kappa_\psi^3\kappa^3\mu r\log^5(Ld)
        \Big\{\frac{1}{\pi(n\wedge p)}
        +\frac{1}{\sigma^2\pi(n\wedge p)}\Big\}=o(1).
    \end{equation*}

    \paragraph{Step 2: comparison with the target statistic.}
    For every $i\in[n]$ and $j\in[p]$, the following rowwise decompositions hold:
    \begin{align*}
        \hat{\mb X}_{i,\cdot}^{t_0}-\mb X_{i,\cdot}^*
        &=
        \mc P_\Omega(\mb E)_{i,\cdot}\mb Y^*\mb H_{\mb X_i}^{*-1}
        +(\mb R_{\mb X}^{\sf approx})_{i,\cdot}
        +\mb X_{i,\cdot}^{t_0}\big(\hat{\mb O}^{t_0}-\mb O^{t_0}\big),\\
        \widehat{\mb Y}^{\sf{app},t_0}_{j,\cdot}-\mb Y^{\sf{app},*}_{j,\cdot}
        &=
        \mc P_\Omega(\mb E)_{\cdot,j}\t\mb X^{\sf{app},*}\mb H_{\mb Y^{\sf{app}}_{j}}^{*-1}
        +(\mb R_{\mb Y^{\sf{app}}}^{\sf approx})_{j,\cdot}
        +\big(0,\mb Y_{j,\cdot}^{t_0}(\hat{\mb O}^{t_0}-\mb O^{t_0})\big).
    \end{align*}
    In particular, the first line explicitly retains the rotation-estimation
    term that is absent from the population-aligned linear representation.
    Theorem~\ref{thm: linear approximations} gives
    \begin{equation*}
        \norm{\mb R_{\mb X}^{\sf approx}}\ti\vee
        \norm{\mb R_{\mb Y^{\sf{app}}}^{\sf approx}}\ti
        \lesssim\frac{1}
        {\sigma\sqrt{\pi\sigma_1^*}(\log d)^2}.
    \end{equation*}
    Hence, on $\mc E^{\sf good}$,
    \begin{align}
    & |\mc T_{\mb X}-\mc T_{\mb X}^{\sharp}|
        \leq Cs_0\Big\{
        \bar a_{\mb X}\norm{\mb R_{\mb X}^{\sf approx}}\ti
        +\bar a_{\mb X}\norm{\mb X^{t_0}}\ti\delta_{\sf rotation}\\
    &\quad +\delta_a\norm{\hat{\mb X}^{t_0}-\mb X^*}\ti
        \Big\}
        \eqqcolon\varepsilon_{\mb X}^{\sf target}.
        \label{eq: simultaneous target comparison}
    \end{align}
    The analogous bound for $\mc T_{\mb Y}$ replaces $\mb X^{t_0}$ by $\mb Y^{t_0}$ and the last rowwise error by
    $\bignorm{\widehat{\mb Y}^{\sf{app},t_0}-\mb Y^{\sf{app},*}}\ti$.
    The final-GD bounds and $\log L\lesssim\log d$ imply
    \begin{equation*}
        \sigma\sqrt{\pi\sigma_1^*\log L}
        \Big\{\norm{\mb X^{t_0}\mb O^{t_0}-\mb X^*}\ti
        \vee\norm{\mb Y^{\sf{app},t_0}\mb O^{\sf{app},t_0}-\mb Y^{\sf{app},*}}\ti\Big\}
        \lesssim\kappa^2\kappa_\psi^2\sqrt{\mu\rho}\,r\log d.
    \end{equation*}
    Theorem~\ref{thm: linear approximations}, the bounded contrast-norm
    ratios, and the assumptions on $\delta_{\sf rotation}$ and
    $\delta_a$ therefore give the following bound. In particular, the
    linearization residual contributes at most a constant multiple of
    $\sqrt{\log L}/(\log d)^2=o(1)$ after standardization:
    \begin{align}
        \max\Big\{
        \frac{\varepsilon_{\mb X}^{\sf target}}
        {\underline a_{\mb X}(\sigma^2\pi\sigma_1^*)^{-1/2}},
        \frac{\varepsilon_{\mb Y}^{\sf target}}
        {\underline a_{\mb Y}(\sigma^2\pi\sigma_1^*)^{-1/2}}
        \Big\}
        \sqrt{\log L}=o(1).
        \label{eq: negligible simultaneous target remainder}
    \end{align}
    Combining \eqref{eq: ideal simultaneous approximation} with Gaussian anti-concentration and \eqref{eq: negligible simultaneous target remainder} transfers the ideal approximation to the two target statistics.

    \paragraph{Step 3: comparison with the feasible bootstrap.}
    Write $\hat F_{\mb X}(x)\coloneqq
    \bb P_\xi(\mc T_{\mb X}^{[b]}\leq x\mid\mb R,\bo\Omega)$ and define
    $F_{\mb X}^{\sharp}(x)$ analogously using
    $\mc T_{\mb X}^{\sharp[b]}$. For each bootstrap draw, set
    $\mb D_{\mb Y}^{[b]}\coloneqq
    \diag(\xi_{\mb Y,j}^{[b]})_{j\in[p]}$ and introduce the population-aligned plug-in information matrix
    \begin{equation*}
        \widetilde{\mb H}_{\mb X_i}^{t_0}
        \coloneqq(\mb O^{t_0})^\top\Big\{\sum_{j\in[p],(i,j)\in\Omega}
        \psi'(\hat M_{i,j})\mb Y_j^{t_0}(\mb Y_j^{t_0})^\top\Big\}\mb O^{t_0}.
    \end{equation*}
    Since $\hat{\mb M}=\mc M(\bo\theta^{t_0})$ and
    $\hat{\mb Y}^{t_0}=\mb Y^{t_0}\hat{\mb O}^{t_0}$, rotational equivariance gives
    \begin{align*}
        \hat{\mb H}_{\mb X_i}
        &=(\hat{\mb O}^{t_0})^\top\mb O^{t_0}
        \widetilde{\mb H}_{\mb X_i}^{t_0}(\mb O^{t_0})^\top\hat{\mb O}^{t_0},\\
        (\bo\Delta_{\mb X}^{[b]})_{i,\cdot}(\hat{\mb O}^{t_0})^\top\mb O^{t_0}
        &=\hat{\mb E}_{i,\cdot}\mb D_{\mb Y}^{[b]}\mb Y^{t_0}\mb O^{t_0}
        (\widetilde{\mb H}_{\mb X_i}^{t_0})^{-1}.
    \end{align*}
    Define the ideal bootstrap row by
    $(\bo\Delta_{\mb X}^{\sharp[b]})_{i,\cdot}
    \coloneqq\mc P_\Omega(\mb E)_{i,\cdot}\mb D_{\mb Y}^{[b]}\mb Y^*
    \mb H_{\mb X_i}^{*-1}$. Adding and subtracting the population-aligned
    factors and information matrix yields the exact decomposition
    \begin{align}
        (\bo\Delta_{\mb X}^{[b]}-\bo\Delta_{\mb X}^{\sharp[b]})_{i,\cdot}
        ={}&\underbrace{(\bo\Delta_{\mb X}^{[b]})_{i,\cdot}
        \{\mb I_r-(\hat{\mb O}^{t_0})^\top\mb O^{t_0}\}}
        _{\bo\beta_{1,i}^{[b]}}\notag\\
        &+\underbrace{\{\hat{\mb E}-\mc P_\Omega(\mb E)\}_{i,\cdot}\mb D_{\mb Y}^{[b]}
        \mb Y^{t_0}\mb O^{t_0}(\widetilde{\mb H}_{\mb X_i}^{t_0})^{-1}}
        _{\bo\beta_{2,i}^{[b]}}\notag\\
        &+\underbrace{\mc P_\Omega(\mb E)_{i,\cdot}\mb D_{\mb Y}^{[b]}
        (\mb Y^{t_0}\mb O^{t_0}-\mb Y^*)
        (\widetilde{\mb H}_{\mb X_i}^{t_0})^{-1}}
        _{\bo\beta_{3,i}^{[b]}}\notag\\
        &+\underbrace{\mc P_\Omega(\mb E)_{i,\cdot}\mb D_{\mb Y}^{[b]}\mb Y^*
        \{(\widetilde{\mb H}_{\mb X_i}^{t_0})^{-1}-\mb H_{\mb X_i}^{*-1}\}}
        _{\bo\beta_{4,i}^{[b]}}.
        \label{eq: feasible multiplier decomposition}
    \end{align}
    The terms $\bo\beta_{1,i}^{[b]}$ through $\bo\beta_{4,i}^{[b]}$
    correspond, respectively, to rotation estimation, score replacement,
    factor estimation, and information-matrix inversion.

    We first control the information perturbation in the last line. The
    max-norm calculation preceding \eqref{eq: term 4 in the projection analysis},
    the masked-factor bound preceding \eqref{eq: term 3 in the projection analysis},
    and \eqref{eq: uniform invertibility of rowwise information matrices} give
    \begin{align}
        \max_{i\in[n]}\frac{\bignorm{\widetilde{\mb H}_{\mb X_i}^{t_0}-\mb H_{\mb X_i}^*}}
        {\underline c_\psi\pi\sigma_r^*}
        \lesssim{}&\frac{\sigma\kappa\kappa_\psi\sqrt{\mu r d\log d}}
        {\underline c_\psi\sqrt\pi\,\sigma_r^*}
        +\frac{\sigma\kappa^3\kappa_\psi^2\mu r^{\frac32}\sqrt{\rho\log d}}
        {\underline c_\psi\sqrt{\pi(n\wedge p)}}=o(1).
        \label{eq: bootstrap information perturbation bound}
    \end{align}
    The quadratic masked-factor error is absorbed by its linear counterpart
    under the signal-strength condition. Consequently,
    $\lambda_{\min}(\widetilde{\mb H}_{\mb X_i}^{t_0})\gtrsim
    \underline c_\psi\pi\sigma_r^*$ uniformly in $i$, and
    \begin{equation*}
        (\widetilde{\mb H}_{\mb X_i}^{t_0})^{-1}-\mb H_{\mb X_i}^{*-1}
        =(\widetilde{\mb H}_{\mb X_i}^{t_0})^{-1}
        (\mb H_{\mb X_i}^*-\widetilde{\mb H}_{\mb X_i}^{t_0})
        \mb H_{\mb X_i}^{*-1}.
    \end{equation*}

    Applying conditional Gaussian concentration uniformly over the active
    coordinates, together with the response-noise bounds and the preceding
    information estimates, shows that, on $\mc E^{\sf good}$,
    with $\bb P_\xi$-probability at least $1-O(d^{-c-2})$,
    \begin{align}
        \max_{i\in[n]}\bignorm{\bo\beta_{1,i}^{[b]}}_2
        &\lesssim\frac{\sigma\sqrt{\kappa\mu r\log d}}
        {\underline c_\psi\sqrt{\pi\sigma_r^*}}\delta_{\sf rotation},\notag\\
        \max_{i\in[n]}\bignorm{\bo\beta_{2,i}^{[b]}}_2
        &\lesssim\frac{\sigma\kappa_\psi^2\kappa^{\frac52}\mu^{\frac32}r^2
        \sqrt\rho\log d}
        {\underline c_\psi\sqrt{\pi(n\wedge p)}\sqrt{\pi\sigma_r^*}},\notag\\
        \max_{i\in[n]}\bignorm{\bo\beta_{3,i}^{[b]}}_2
        &\lesssim\frac{\sigma\sqrt{d\log d}}
        {\underline c_\psi\sqrt\pi\,\sigma_r^*}
        \frac{\sigma\kappa\kappa_\psi\sqrt{\kappa\mu\rho r^2\log d}}
        {\underline c_\psi\sqrt{\pi\sigma_r^*}},\notag\\
        \max_{i\in[n]}\bignorm{\bo\beta_{4,i}^{[b]}}_2
        &\lesssim\frac{\sigma\sqrt{\kappa\mu r\log d}}
        {\underline c_\psi\sqrt{\pi\sigma_r^*}}
        \Big\{\frac{\sigma\kappa\kappa_\psi\sqrt{\mu r d\log d}}
        {\underline c_\psi\sqrt\pi\,\sigma_r^*}
        +\frac{\sigma\kappa^3\kappa_\psi^2\mu r^{\frac32}\sqrt{\rho\log d}}
        {\underline c_\psi\sqrt{\pi(n\wedge p)}}\Big\}.
        \label{eq: bounds for feasible multiplier perturbations}
    \end{align}
    Here the score-replacement bound uses
    $\hat{\mb E}-\mc P_\Omega(\mb E)
    =\mc P_\Omega\{\psi(\mb M^*)-\psi(\hat{\mb M})\}$ and the predictor
    max-norm bound, while the factor-estimation bound uses the final-GD
    rowwise rate, including its factor $\sqrt r$. The two terms in the last
    braces are precisely the two contributions in
    \eqref{eq: bootstrap information perturbation bound}.

    The same conditional Gaussian argument gives
    \begin{equation*}
        \max_{i\in[n]}\bignorm{(\bo\Delta_{\mb X}^{[b]})_{i,\cdot}}_2
        \lesssim\frac{\sigma\sqrt{\kappa\mu r\log d}}
        {\underline c_\psi\sqrt{\pi\sigma_r^*}}.
    \end{equation*}
    Since the supports have bounded size, the elementary inequality
    $|\max_l u_l-\max_l v_l|\leq\max_l|u_l-v_l|$ gives
    \begin{equation*}
        |\mc T_{\mb X}^{[b]}-\mc T_{\mb X}^{\sharp[b]}|
        \lesssim\bar a_{\mb X}\sum_{q=1}^4
        \max_{i\in[n]}\bignorm{\bo\beta_{q,i}^{[b]}}_2
        +\delta_a\max_{i\in[n]}\bignorm{(\bo\Delta_{\mb X}^{[b]})_{i,\cdot}}_2.
    \end{equation*}
    The bounds in \eqref{eq: bounds for feasible multiplier perturbations}
    and the preceding feasible-row bound are jointly controlled by
    \begin{align*}
        \xi_{\mb X}^{\sf boot}
        \coloneqq{}&
        \frac{\sigma\kappa_\psi^2\kappa^{\frac52}\mu^{\frac32}r^2
        \sqrt\rho\log d}
        {\underline c_\psi\sqrt{\pi(n\wedge p)}\sqrt{\pi\sigma_r^*}}
        +\frac{\sigma\sqrt{\kappa\mu r\log d}}
        {\underline c_\psi\sqrt{\pi\sigma_r^*}}
        \Big(\delta_{\sf rotation}+\frac{\delta_a}{\bar a_{\mb X}}\Big)\\
        &+\frac{\sigma\sqrt{d\log d}}
        {\underline c_\psi\sqrt\pi\,\sigma_r^*}
        \frac{\sigma\kappa\kappa_\psi\sqrt{\kappa\mu^2\rho r^2\log d}}
        {\underline c_\psi\sqrt{\pi\sigma_r^*}}
        +\frac{\sigma^2\kappa_\psi^2\kappa^{\frac72}\mu^{\frac32}r^2
        \sqrt\rho\log d}
        {\underline c_\psi^2\sqrt{\pi(n\wedge p)}\sqrt{\pi\sigma_r^*}}.
    \end{align*}
    Consequently,
    \begin{align}
        \bb P_\xi\Big(
        |\mc T_{\mb X}^{[b]}-\mc T_{\mb X}^{\sharp[b]}|
        >C\bar a_{\mb X}\xi_{\mb X}^{\sf boot}\,\Big|\,\mb R,\bo\Omega\Big)
        \leq O(d^{-c-2}).
        \label{eq: feasible multiplier coupling}
    \end{align}
    The third term lifts the orders of $\mu$ and $\rho$ so that it controls
    both the factor-estimation contribution and the first information-inverse
    contribution. The last term retains the link-variation contribution from
    information-matrix inversion; this term cannot be discarded using only
    relative consistency of the plug-in information matrices.

    Using $\sigma^2=\kappa_\psi\underline c_\psi$,
    $\sigma_1^*=\kappa\sigma_r^*$, and $\log L\lesssim\log d$, we obtain
    \begin{align*}
        &\sigma\sqrt{\pi\sigma_1^*\log L}\,\xi_{\mb X}^{\sf boot}\\
        \lesssim{}&
        \frac{\kappa_\psi^3\kappa^3\mu^{\frac32}r^2\sqrt\rho
        (\log d)^{\frac32}}{\sqrt{\pi(n\wedge p)}}
        +\kappa_\psi\kappa\sqrt{\mu r}\log d
        \Big(\delta_{\sf rotation}+\frac{\delta_a}{\bar a_{\mb X}}\Big)\\
        &+\frac{\kappa^2\kappa_\psi^{\frac52}\mu r\sqrt{\rho d}
        (\log d)^{\frac32}}
        {\sqrt{\underline c_\psi\pi}\,\sigma_r^*}
        +\frac{\kappa_\psi^4\kappa^4\mu^{\frac32}r^2\sqrt\rho
        (\log d)^{\frac32}}
        {\sigma\sqrt{\pi(n\wedge p)}}=o(1).
    \end{align*}
    The first and last terms are $o(1)$ by the two bounds in
    \eqref{eq: effective sample sizes from linear approximation}, respectively;
    the third follows from the signal-strength condition. The middle term is
    $o(1)$ by the assumptions on $\delta_{\sf rotation}$ and $\delta_a$ and
    the bounded contrast-norm ratios.

    Crucially, \eqref{eq: feasible multiplier coupling} is used at the distribution level. For every $x\in\bb R$ it gives
    \begin{align}
        F_{\mb X}^{\sharp}(x-C\bar a_{\mb X}\xi_{\mb X}^{\sf boot})-O(d^{-c-2})
        \leq\hat F_{\mb X}(x)
        \leq F_{\mb X}^{\sharp}(x+C\bar a_{\mb X}\xi_{\mb X}^{\sf boot})+O(d^{-c-2}).
        \label{eq: bootstrap cdf sandwich}
    \end{align}
    Conditional on the score array, $\mc T_{\mb X}^{\sharp[b]}$ is the maximum of a Gaussian vector. By
    \eqref{eq: uniform empirical contrast variance}, its coordinate standard
    deviations are at least
    $c\underline a_{\mb X}(\sigma^2\pi\sigma_1^*)^{-1/2}$ with probability
    $1-o(1)$. Nazarov's inequality and the preceding bounds therefore imply
    \begin{align}
        \sup_{x\in\bb R}
        |\hat F_{\mb X}(x)-F_{\mb X}^{\sharp}(x)|
        \lesssim
        \frac{\bar a_{\mb X}\xi_{\mb X}^{\sf boot}}
        {\underline a_{\mb X}(\sigma^2\pi\sigma_1^*)^{-1/2}}
        \sqrt{\log L}
        +O(d^{-c-2})
        =o_{\bb P}(1).
        \label{eq: feasible bootstrap cdf consistency}
    \end{align}
    This CDF comparison, rather than a pointwise claim about two random quantiles, yields the required quantile transfer. Combining
    \eqref{eq: ideal simultaneous approximation},
    \eqref{eq: negligible simultaneous target remainder}, and
    \eqref{eq: feasible bootstrap cdf consistency} proves
    $\bb P(\mc T_{\mb X}\leq\mc G_{\mb X,1-\alpha})=1-\alpha+o(1)$.
    The conditioning on $\bo\Omega$ is removed by integrating over the mask and adding $\bb P(\mc E^{\sf good\complement})=O(d^{-c})$.

    For $\mc T_{\mb Y}$, group the independent score contributions by rows:
    \begin{align*}
        Z_{\mb Y,i,l}
        \coloneqq
        \sum_{(j,k)\in\mc S_{\mb Y,l}}
        a_{\mb Y,l,j,k}\Omega_{i,j}E_{i,j}
        \big[\mb X^{\sf{app},*}_{i}\mb H_{\mb Y^{\sf{app}}_{j}}^{*-1}\big]_k.
    \end{align*}
    Repeating the preceding argument with the $(r+1)$-dimensional augmented
    information matrices, $\bar a_{\mb Y}$, and $\underline a_{\mb Y}$ proves
    the second coverage claim. The corresponding bootstrap-comparison bound is
    obtained from $\xi_{\mb X}^{\sf boot}$ by replacing
    $\bar a_{\mb X}$ with $\bar a_{\mb Y}$ in its weight-estimation term, so
    the same calculation applies. Finally, conditional on the data, the Dvoretzky--Kiefer--Wolfowitz inequality gives a uniform empirical-CDF error
    $O_{\bb P}(N_{\sf boot}^{-1/2})$, which establishes the stated conclusion for the Monte Carlo quantiles.
    \end{proof}

    \begin{proof}[Proof of Proposition~\ref{thm:left-factor-ranking-consistency}]
    Apply Theorem~\ref{thm: simultaneous inference} to the $L=n(n-1)$
    ordered pairwise contrasts. For the contrast indexed by $(i,j)$, its only
    nonzero population weights are $s_{i,j,k}^{-1}$ at $(i,k)$ and
    $-s_{i,j,k}^{-1}$ at $(j,k)$, with the corresponding plug-in weights
    obtained by replacing $s_{i,j,k}$ with $\hat s_{i,j,k}$. Thus every
    support has size two, $\log L\lesssim\log d$, and the contrast-norm,
    rotation, and weight-estimation requirements are
    exactly those imposed in the proposition. Including both $(i,j)$ and
    $(j,i)$ turns the maximum of the signed contrasts into
    \eqref{eq: ranking stat}. Therefore,
    \begin{equation}
        \bb P\big(\mc T_k^{\sf rank}\leq\hat q_{1-\alpha}^{\sf rank}\big)=1-\alpha+o(1).
        \label{eq: simultaneous pairwise ranking event}
    \end{equation}

    On the event in \eqref{eq: simultaneous pairwise ranking event},
    $\hat X_{l,k}^{t_0}-\hat X_{i,k}^{t_0}>
    \hat q_{1-\alpha}^{\sf rank}\hat s_{l,i,k}$ implies
    $X_{l,k}^*>X_{i,k}^*$. Hence $L_i\leq\operatorname{rank}_k(i)$.
    Similarly,
    $\hat X_{l,k}^{t_0}-\hat X_{i,k}^{t_0}<
    -\hat q_{1-\alpha}^{\sf rank}\hat s_{l,i,k}$ implies
    $X_{l,k}^*<X_{i,k}^*$, and therefore
    $U_i\geq\operatorname{rank}_k(i)$. These implications hold
    simultaneously for every $i$, so
    \eqref{eq: simultaneous pairwise ranking event} yields the claimed
    coverage.
    \end{proof}

\begin{proof}[Proof of Proposition~\ref{prop:missing-entry-prediction}]
By Lemma~\ref{lemma: scaling equivalence}, it suffices to take $\omega=1$;
the claim is vacuous when $\Omega^{\complement}=\emptyset$. Set
$\widetilde{\mb X}^{t_0}\coloneqq\mb X^{t_0}\mb O^{t_0}$,
$\widetilde{\mb Y}^{t_0}\coloneqq\mb Y^{t_0}\mb O^{t_0}$,
$\widetilde{\mb X}^{\sf{app},t_0}\coloneqq
\mb X^{\sf{app},t_0}\diag(1,\mb O^{t_0})$, and
$\widetilde{\mb Y}^{\sf{app},t_0}\coloneqq
\mb Y^{\sf{app},t_0}\diag(1,\mb O^{t_0})$. The band and its bootstrap
statistics are invariant under this common rotation, so we may use the
oracle alignment throughout.

For $(i,j)\in\Omega^{\complement}$, define
\begin{align*}
    s_{\mb X,i,j}^2
    &\coloneqq(\mb Y_j^*)^\top\mb H_{\mb X_i}^{*-1}\mb Y_j^*,&
    s_{\mb Y^{\sf{app}},i,j}^2
    &\coloneqq(\mb X_i^{\sf{app},*})^\top
    \mb H_{\mb Y_j^{\sf{app}}}^{*-1}\mb X_i^{\sf{app},*}.
\end{align*}
The uniform information bounds imply, whenever $\mb Y_j^*\neq\mb 0_r$,
\begin{align}
    \sqrt{\underline c_\psi\pi\sigma_r^*}
    &\lesssim\frac{\norm{\mb Y_j^*}_2}{s_{\mb X,i,j}}
    \lesssim\sigma\sqrt{\pi\sigma_1^*},&
    \sqrt{\underline c_\psi\pi\sigma_r^*}
    &\lesssim
    \frac{\norm{\mb X_i^{\sf{app},*}}_2}
    {s_{\mb Y^{\sf{app}},i,j}}
    \lesssim\sigma\sqrt{\pi\sigma_1^*},
    \label{eq: missing-entry contrast norm bounds}
\end{align}
and the appended intercept coordinate further gives
$s_{\mb Y^{\sf{app}},i,j}\gtrsim(\sigma\sqrt{\pi n})^{-1}$.
Let
\[
    \mc J_{\mb X}\coloneqq
    \Big\{(i,j)\in\Omega^{\complement}:
    s_{\mb X,i,j}\geq
    s_{\mb Y^{\sf{app}},i,j}/(\log d)^{5/4}\Big\}.
\]
When $\mc J_{\mb X}=\emptyset$, the row-side application below is omitted
and its restricted critical value is set to zero.

We apply Theorem~\ref{thm: simultaneous inference} separately to the two
prediction components. For each $(i,j)\in\mc J_{\mb X}$ and
$u\in\{-1,1\}$, take the support $\{(i,k):k\in[r]\}$ and the oracle and
estimated weights
\[
    a_{\mb X,(i,j,u),i,k}
    \coloneqq\frac{uY_{j,k}^*}{s_{\mb X,i,j}},
    \qquad
    \hat a_{\mb X,(i,j,u),i,k}
    \coloneqq\frac{u\widetilde Y_{j,k}^{t_0}}
    {\hat s_{\mb X,i,j}}.
\]
For the column component, use the support
$\{(j,k):k\in[r+1]\}$ for every
$(i,j)\in\Omega^{\complement}$ and set
\[
    a_{\mb Y,(i,j,u),j,k}
    \coloneqq\frac{uX_{i,k}^{\sf{app},*}}
    {s_{\mb Y^{\sf{app}},i,j}},
    \qquad
    \hat a_{\mb Y,(i,j,u),j,k}
    \coloneqq\frac{u\widetilde X_{i,k}^{\sf{app},t_0}}
    {\hat s_{\mb Y^{\sf{app}},i,j}}.
\]
The two signs recover the absolute maxima. Their numbers are at most $2np$,
their supports have bounded size, and \eqref{eq: missing-entry contrast norm bounds}
verifies the contrast-norm condition in Theorem~\ref{thm: simultaneous inference}.
The active set and the oracle weights depend on $\bo\Omega$ but not on
$\mb E$, as permitted by that theorem.
Moreover, \eqref{eq: rowwise error from linear approximation},
\eqref{eq: bootstrap information perturbation bound}, and its augmented-column
counterpart give
\begin{align}
    &\max_{\substack{(i,j)\in\mc J_{\mb X},\,u\in\{-1,1\}\\k\in[r]}}
    \big|\hat a_{\mb X,(i,j,u),i,k}
    -a_{\mb X,(i,j,u),i,k}\big|
    \vee
    \max_{\substack{(i,j)\in\Omega^{\complement},\,u\in\{-1,1\}\\k\in[r+1]}}
    \big|\hat a_{\mb Y,(i,j,u),j,k}
    -a_{\mb Y,(i,j,u),j,k}\big|\notag\\
    &\qquad\lesssim_{\bb P}
    \sigma\sqrt{\pi n}(\log d)^{5/4}
    \Big\{\bignorm{\widetilde{\mb X}^{\sf{app},t_0}
    -\mb X^{\sf{app},*}}\ti
    \vee\bignorm{\widetilde{\mb Y}^{t_0}-\mb Y^*}\ti\Big\}\notag\\
    &\qquad\quad+
    \sigma\sqrt{\pi\sigma_1^*}
    \Big\{
    \max_{(i,j)\in\mc J_{\mb X}}
    \Big|\frac{\hat s_{\mb X,i,j}}{s_{\mb X,i,j}}-1\Big|
    \vee
    \max_{(i,j)\in\Omega^{\complement}}
    \Big|\frac{\hat s_{\mb Y^{\sf{app}},i,j}}
    {s_{\mb Y^{\sf{app}},i,j}}-1\Big|
    \Big\}\notag\\
    &\qquad=o_{\bb P}\Big\{
    \frac{\sqrt{\underline c_\psi\pi\sigma_r^*}}
    {\sqrt{\mu}\log d}\Big\}.
    \label{eq: missing-entry weight verification}
\end{align}
Thus all conditions of Theorem~\ref{thm: simultaneous inference} hold with
the oracle rotation, for which $\delta_{\sf rotation}=0$.

It remains only to account for the omitted row directions and the bilinear
remainder. The same rowwise expansion used in the proof of
Theorem~\ref{thm: simultaneous inference} yields
\begin{align}
    &\max_{(i,j)\in\Omega^{\complement}\setminus\mc J_{\mb X}}
    \frac{\big|(\widetilde{\mb X}_{i,\cdot}^{t_0}
    -\mb X_{i,\cdot}^*)\widetilde{\mb Y}_j^{t_0}\big|}
    {\hat s_{\mb Y^{\sf{app}},i,j}}\notag\\
    &\qquad\lesssim_{\bb P}\frac{1}{(\log d)^{5/4}}
    \max_{\substack{(i,j)\in\Omega^{\complement}\\s_{\mb X,i,j}>0}}
    \frac{\big|\mc P_\Omega(\mb E)_{i,\cdot}\mb Y^*
    \mb H_{\mb X_i}^{*-1}\mb Y_j^*\big|}{s_{\mb X,i,j}}
    +o_{\bb P}\big\{(\log d)^{-1/2}\big\}\notag\\
    &\qquad=O_{\bb P}\big\{(\log d)^{-3/4}\big\}
    +o_{\bb P}\big\{(\log d)^{-1/2}\big\}
    =o_{\bb P}\big\{(\log d)^{-1/2}\big\}.
    \label{eq: negligible inactive missing-entry directions}
\end{align}
If $s_{\mb X,i,j}=0$, then $\mb Y_j^*=\mb 0_r$, and the corresponding
term is absorbed by the factor-product remainder below. In addition,
\eqref{eq: rowwise error from linear approximation} gives
\begin{align}
    &\max_{(i,j)\in\Omega^{\complement}}
    \frac{\big|(\widetilde{\mb X}_{i,\cdot}^{t_0}-\mb X_{i,\cdot}^*)
    (\widetilde{\mb Y}_{j,\cdot}^{t_0}-\mb Y_{j,\cdot}^*)^\top\big|}
    {\hat s_{\mb X,i,j}+\hat s_{\mb Y^{\sf{app}},i,j}}
    \sqrt{\log d}\notag\\
    &\qquad\lesssim
    \frac{\sigma^3\mu\sqrt n(\log d)^{3/2}}
    {\underline c_\psi^2\sqrt\pi\,\sigma_r^*}
    =o_{\bb P}(1).
    \label{eq: negligible missing-entry product remainder}
\end{align}

Let $\hat q_{\mb X,1-\frac{\alpha}{2}}^{\sf ent}(\mc J_{\mb X})$
denote the row-side empirical bootstrap quantile with the maximum restricted
to $\mc J_{\mb X}$. The full quantile dominates this restricted quantile.
Applying Theorem~\ref{thm: simultaneous inference} at level
$1-\alpha/2$ to each collection, and using the Gaussian
anti-concentration argument in its proof, we may choose
$\varepsilon_d=o\{(\log d)^{-1/2}\}$, dominating the two remainders above
with probability $1-o(1)$, such that
\begin{align}
    \bb P\Big[
    &\max_{(i,j)\in\mc J_{\mb X}}
    \frac{\big|(\widetilde{\mb X}_{i,\cdot}^{t_0}-\mb X_{i,\cdot}^*)
    \widetilde{\mb Y}_j^{t_0}\big|}{\hat s_{\mb X,i,j}}
    \leq\hat q_{\mb X,1-\frac{\alpha}{2}}^{\sf ent}
    (\mc J_{\mb X})-4\varepsilon_d,\notag\\
    &\max_{(i,j)\in\Omega^{\complement}}
    \frac{\big|(\widetilde{\mb Y}_{j,\cdot}^{\sf{app},t_0}
    -\mb Y_{j,\cdot}^{\sf{app},*})
    \widetilde{\mb X}_i^{\sf{app},t_0}\big|}
    {\hat s_{\mb Y^{\sf{app}},i,j}}
    \leq\hat q_{\mb Y^{\sf{app}},1-\frac{\alpha}{2}}^{\sf ent}
    -4\varepsilon_d
    \Big]\geq1-\alpha-o(1).
    \label{eq: missing-entry componentwise coverage}
\end{align}
No independence between the two collections is required.

Finally, direct algebra gives
\begin{align*}
    \hat M_{i,j}-M_{i,j}^*
    ={}&(\widetilde{\mb X}_{i,\cdot}^{t_0}-\mb X_{i,\cdot}^*)
    \widetilde{\mb Y}_j^{t_0}
    +(\widetilde{\mb Y}_{j,\cdot}^{\sf{app},t_0}
    -\mb Y_{j,\cdot}^{\sf{app},*})
    \widetilde{\mb X}_i^{\sf{app},t_0}\\
    &-(\widetilde{\mb X}_{i,\cdot}^{t_0}-\mb X_{i,\cdot}^*)
    (\widetilde{\mb Y}_{j,\cdot}^{t_0}-\mb Y_{j,\cdot}^*)^\top.
\end{align*}
Combining this identity with
\eqref{eq: negligible inactive missing-entry directions}--\eqref{eq: missing-entry componentwise coverage}
yields, uniformly over $(i,j)\in\Omega^{\complement}$,
\begin{align*}
    |\hat M_{i,j}-M_{i,j}^*|
    \leq{}&
    \hat q_{\mb X,1-\frac{\alpha}{2}}^{\sf ent}
    \hat s_{\mb X,i,j}
    +\hat q_{\mb Y^{\sf{app}},1-\frac{\alpha}{2}}^{\sf ent}
    \hat s_{\mb Y^{\sf{app}},i,j}
\end{align*}
with probability at least $1-\alpha-o(1)$. Since $\psi=\Psi'$ is increasing
on the natural-parameter domain, applying $\psi$ to both endpoints completes
the proof.
\end{proof}

\subsection{Rotation Estimations based on the Gradient Descent Outputs} 
\label{subsec: rotation-free factor estimation}

We provide rigorous guarantees for the factor rotation recovery based on the gradient descent outputs and the two rotation-free estimators $\hat{\mb O}^{\sf eig}$ and $\hat{\mb O}^{\sf var}$ defined in Section~\ref{subsec: gd theory} of the main paper.
\begin{proposition}[Recovery of the factor rotation]
\label{prop:rotation-recovery}
With the notation of Section~\ref{subsec: gd theory} of the main paper, suppose that the
conditions of Theorem~\ref{thm: linear approximations} hold.

\begin{enumerate}
\item[(i)] Under the distinct-spectrum condition \eqref{eq: distinct singular value conditions} in the main paper, with probability at least
$1-O(d^{-c})$,
$$
    \min_{\mb D\in\mc D_{\pm}(r)}
    \bignorm{{\mb O^{t_0}}\t\hat{\mb O}^{\sf{eig}}-\mb D}\fb
    \lesssim
    \frac{1}{\Delta_{\sf{eig}}}
    \Big\{
    \frac{\sigma\sqrt{\log d}}
    {\underline c_\psi\sqrt{\pi}}
    +\frac{\sqrt n}{\sigma\sqrt{\pi}(\log d)^2}
    +\frac{\sigma^2 d\log d}
    {\underline c_\psi^2\pi\sigma_r^*}
    \Big\}. 
$$
Here, recall that $\mc D_{\pm}(r)$ denotes the set of diagonal sign matrices.

\item[(ii)] Suppose that $\mb X^* = \mb J_n \mb Z$, where the entries of $\mb Z$ are i.i.d. from a fixed distribution satisfying $\bb E[Z_{1,1}^4] > 3 \bb E[Z_{1,1}^2]^2$. 
Then, with probability at least
$1-O(d^{-c})$,
$$
    \min_{\mb P\in\mc P_{\pm}(r)}
    \bignorm{{\mb O^{t_0}}\t\hat{\mb O}^{\sf{var}}-\mb P}\fb
    \lesssim
    \frac{(\log d)^3}{\sqrt n}
    +\frac{\sigma\omega^{1/4}}
    {\underline c_\psi\sqrt{\pi\sigma_r^*}}. 
$$
Here, $\mc P_{\pm}(r)$ denotes the set of signed permutation matrices.
The analogous conclusions for $\mb Y^*$ follow by replacing
$(n,\omega^{1/4})$ with $(p,\omega^{-1/4})$.
\end{enumerate}

\end{proposition}

\begin{proof}
Put $(\mb L_{\mb X})_{i,\cdot}\coloneqq
\mc P_\Omega(\mb E)_{i,\cdot}\mb Y^*\mb H_{\mb X_i}^{*-1}$ and
$\mb R_{\mb X}\coloneqq
\mb X^{t_0}\mb O^{t_0}-\mb X^*-\mb L_{\mb X}$.
Theorem~\ref{thm: linear approximations} gives
$$
    \norm{\mb R_{\mb X}}\ti
    \lesssim\frac{\omega^{1/4}}
    {\sigma\sqrt{\pi\sigma_1^*}(\log d)^2}.
$$
Since $\sigma^2=\kappa_\psi\underline c_\psi$, this further implies
$$
    \norm{\mb R_{\mb X}}\ti
    =o\Big\{\frac{\sigma\sqrt r\,\omega^{1/4}}
    {\underline c_\psi\sqrt{\pi\sigma_r^*}}\Big\}.
$$
Indeed, the ratio between the rate from
Theorem~\ref{thm: linear approximations} and the displayed rate is of order
$$
    \frac{1}{\kappa_\psi\sqrt{\kappa r}(\log d)^2}=o(1).
$$
We work on the high-probability event from
Theorem~\ref{thm: linear approximations}. On this event, incoherence,
identifiability, and concentration of the observed information matrices
give
\begin{equation}
    \max_i\norm{\mb X_i^*}_2
    \lesssim
    \omega^{1/4}\sqrt{\frac{\mu r\sigma_1^*}{n}},
    \quad
    \max_j\norm{\mb Y_j^*}_2
    \lesssim
    \omega^{-1/4}\sqrt{\frac{\mu r\sigma_1^*}{p}},\quad 
    \max_i\norm{\mb H_{\mb X_i}^{*-1}}
    \lesssim
    \frac{\omega^{1/2}}
         {\underline c_\psi\pi\sigma_r^*}.
\label{P.1}
\end{equation}

We first establish the claimed improvement for the linear term. We have
$
    \omega^{-1/2}{\mb X^*}\t\mb L_{\mb X}
    =
    \sum_{i=1}^n\sum_{j=1}^p \mb S_{ij}^{\sf X}$,
    $
    \mb S_{ij}^{\sf X}
    \coloneqq
    \omega^{-1/2}\Omega_{ij}E_{ij}
    \mb X_i^*{\mb Y_j^*}\t\mb H_{\mb X_i}^{*-1}.
$
Conditional on $\Omega$, the matrices $\{\mb S_{ij}^{\sf X}\}$ are independent
and centered. By \eqref{P.1}, $\max_{i,j}\norm{\mb S_{ij}^{\sf X}}
    \lesssim
    \frac{B\mu r\kappa}
         {\underline c_\psi\pi\sqrt{np}}$,
whereas a direct calculation using
$\mb H_{\mb X_i}^*\succeq
\underline c_\psi\sum_j\Omega_{ij}\mb Y_j^*{\mb Y_j^*}\t$
gives
$$
    \max\Big\{
        \bignorm{\sum_{i,j}\bb E(\mb S_{ij}^{\sf X}{\mb S_{ij}^{\sf X}}\t\mid\Omega)},
        \bignorm{\sum_{i,j}\bb E({\mb S_{ij}^{\sf X}}\t\mb S_{ij}^{\sf X}\mid\Omega)}
    \Big\}
    \lesssim
    \frac{\sigma^2\kappa r}{\underline c_\psi^2\pi}.
$$
The matrix Bernstein inequality \citep{tropp2015introduction} therefore
yields, on its intersection with the linear-approximation event and with
probability at least $1-O(d^{-c})$,
\begin{equation}
    \omega^{-1/2}
    \bignorm{{\mb X^*}\t\mb L_{\mb X}}
    \lesssim
    \frac{\sigma\sqrt{\kappa r\log d}}
         {\underline c_\psi\sqrt{\pi}}
    +
    \frac{B\mu r\kappa\log d}
         {\underline c_\psi\pi\sqrt{np}}.
\label{P.2}
\end{equation}
The second term is absorbed by the first under
Assumption~\ref{assumption: signal strength and incoherence degree}(c).

For the distinct-spectrum case, let
$\hat{\mb G}_{\mb X}
=\omega^{-1/2}{\mb X^{t_0}}\t\mb X^{t_0}$ and define
$\mb\Delta_{\mb X}\coloneqq
\mb X^{t_0}\mb O^{t_0}-\mb X^*
=\mb L_{\mb X}+\mb R_{\mb X}$.
Then
$$
    {\mb O^{t_0}}\t\hat{\mb G}_{\mb X}\mb O^{t_0}-\bo\Sigma^*
    =
    \omega^{-1/2}\big(
        {\mb X^*}\t\mb\Delta_{\mb X}
        +{\mb\Delta_{\mb X}}\t\mb X^*
        +{\mb\Delta_{\mb X}}\t\mb\Delta_{\mb X}
    \big).
$$
Combining \eqref{P.2} with
$\norm{\mb R_{\mb X}}_F\leq\sqrt n\,
\norm{\mb R_{\mb X}}\ti$ and
$\norm{\mb X^*}=\omega^{1/4}\sqrt{\sigma_1^*}$ gives
$$
\begin{aligned}
    &\bignorm{
        {\mb O^{t_0}}\t\hat{\mb G}_{\mb X}\mb O^{t_0}
        -\bo\Sigma^*
    }\\
    &\quad\lesssim
    \frac{\sigma\sqrt{\kappa r\log d}}
         {\underline c_\psi\sqrt{\pi}}
    +\omega^{-1/4}\sqrt{n\sigma_1^*}\,
        \norm{\mb R_{\mb X}}\ti
    +\omega^{-1/2}\norm{\mb\Delta_{\mb X}}_F^2\\
    &\quad\lesssim
    \frac{\sigma\sqrt{\kappa r\log d}}
         {\underline c_\psi\sqrt{\pi}}
    +\frac{\sqrt n}{\sigma\sqrt{\pi}(\log d)^2}
    +\frac{\sigma^2\kappa d r\log d}
         {\underline c_\psi^2\pi\sigma_r^*},
\end{aligned}
$$
where the last step uses the residual bound above and
Theorem~\ref{thm: GD}. Eigenvector
perturbation \citep{stewart1990matrix}, applied separately to the $r$
simple eigenvalues, now gives
$$
    \min_{\mb D\in\mc D_\pm(r)}
    \bignorm{{\mb O^{t_0}}\t\hat{\mb O}^{\sf{eig}}-\mb D}_F
    \lesssim
    \frac{\sqrt r}{\Delta_{\sf{eig}}}
    \Big\{
    \frac{\sigma\sqrt{\kappa r\log d}}
         {\underline c_\psi\sqrt{\pi}}
    +\frac{\sqrt n}{\sigma\sqrt{\pi}(\log d)^2}
    +\frac{\sigma^2\kappa d r\log d}
         {\underline c_\psi^2\pi\sigma_r^*}
    \Big\}.
$$

It remains to consider varimax. The uncentered decomposition is
$\mb X^{t_0}\mb O^{t_0}-\mb X^*=\mb L_{\mb X}+\mb R_{\mb X}$.
By rotation equivariance, ${\mb O^{t_0}}\t\hat{\mb O}^{\sf{var}}$ maximizes
the sample criterion computed from $\mb X^*+\mb L_{\mb X}+\mb R_{\mb X}$. For
$\mb Q\in\mc O(r)$, the corresponding population criterion is
$$
    \mc V(\mb Q)
    =
    2r+\sum_{a=1}^r
    \big\{\bb E(Z_{i,a}-\bb EZ_{i,a})^4-3\big\}
    \sum_{k=1}^r Q_{ak}^4.
$$
Its maximizers are precisely the signed permutation matrices. Moreover,
for every such matrix $\mb P$ and every sufficiently small
skew-symmetric $\mb A$,
\begin{equation}
    \mc V(\mb P)-\mc V(\mb P\exp(\mb A))
    \gtrsim
    \norm{\mb A}_F^2.
\label{P.3}
\end{equation}
We next make the probability level in the empirical-moment comparison
explicit. Sub-Gaussian concentration gives
$\max_i\norm{\mb Z_i-\bb E\mb Z_i}_2\lesssim\sqrt{\log d}$ with probability
at least $1-O(d^{-c})$. On this truncation event, Bernstein's
inequality, applied to the empirical moments of degree at most four and to a
fixed-dimensional net of $\mc O(r)$, gives
\begin{align*}
    &\sup_{\mb Q\in\mc O(r)}\Big\{
        |\mc V_n(\mb Q;\mb X^*)-\mc V(\mb Q)|
        +\bignorm{\nabla\mc V_n(\mb Q;\mb X^*)-\nabla\mc V(\mb Q)}_F
        \\ &\qquad +\bignorm{\nabla^2\mc V_n(\mb Q;\mb X^*)
        -\nabla^2\mc V(\mb Q)}
    \Big\}
    \lesssim
    \frac{(\log d)^3}{\sqrt n}.
\end{align*}
Here, the derivatives are taken along the tangent spaces of $\mc O(r)$.
Expanding $\mb X^*=\mb J_n\mb Z$ in centered empirical moments shows that
the unsimplified rate is $\sqrt{\log d/n}+(\log d)^3/n$, which is bounded
by $C(\log d)^3/\sqrt n$.

Conditional on $(\mb Z,\bo\Omega)$, the rows of $\mb L_{\mb X}$ are
independent and centered. Using \eqref{P.1}, the noise-truncation event, and
Bernstein's inequality again, the terms containing exactly one copy of
$\mb L_{\mb X}$ perturb the criterion, score, and Hessian uniformly by at
most
$C\sigma\sqrt r\,\omega^{1/4}(\log d)^3/
(\underline c_\psi\sqrt{\pi\sigma_r^*n})$. Terms containing at least two
copies contribute at most
$C\sigma^2r\omega^{1/2}/(\underline c_\psi^2\pi\sigma_r^*)$; the linear
Bernstein remainder involving $B$ is absorbed by
Assumption~\ref{assumption: signal strength and incoherence degree}(c).
Finally, a deterministic expansion of the degree-four criterion and
$\norm{\mb R_{\mb X}}\ti$ bounds all terms containing $\mb R_{\mb X}$ by
$C\norm{\mb R_{\mb X}}\ti$. Under the
condition in part~(ii), the preceding three contributions are bounded by
$$
    C\Big\{\frac{(\log d)^3}{\sqrt n}
    +\frac{\sigma\sqrt r\,\omega^{1/4}}
    {\underline c_\psi\sqrt{\pi\sigma_r^*}}\Big\}.
$$
These bounds hold simultaneously with probability at least $1-O(d^{-c})$.

The resulting uniform criterion bound and \eqref{P.3}, together with the
assumed $o(1)$ condition, imply that
${\mb O^{t_0}}\t\hat{\mb O}^{\sf{var}}$ lies in an $o(1)$ neighborhood of
some $\mb P_n\in\mc P_\pm(r)$. At this signed permutation, the preceding
score bound gives
\begin{equation}
    \bignorm{\nabla\mc V_n(
        \mb P_n;\mb X^*+\mb L_{\mb X}+\mb R_{\mb X})}_F
    \lesssim
    \frac{(\log d)^3}{\sqrt n}
    +\frac{\sigma\sqrt r\,\omega^{1/4}}
    {\underline c_\psi\sqrt{\pi\sigma_r^*}}.
\label{P.4}
\end{equation}
The same uniform comparison and \eqref{P.3} show that the negative empirical
Hessian has smallest tangent-space eigenvalue bounded below by a positive
constant throughout this neighborhood. Applying
Taylor's theorem to the first-order condition and using \eqref{P.4} yields,
on the same event,
$$
    \min_{\mb P\in\mc P_\pm(r)}
    \bignorm{{\mb O^{t_0}}\t\hat{\mb O}^{\sf{var}}-\mb P}_F
    \lesssim
    \frac{(\log d)^3}{\sqrt n}
    +\frac{\sigma\sqrt r\,\omega^{1/4}}
    {\underline c_\psi\sqrt{\pi\sigma_r^*}}.
$$
The proof for $\mb Y^*$ is identical after interchanging the two factor
matrices and replacing $(n,\omega^{1/4})$ with
$(p,\omega^{-1/4})$.
\end{proof}

\section{Proof of Theorem~\ref{thm: minimax lower bound}}
\label{sec: proof of minimax lower bound}
Let $\bb P_{\bo\theta}$ denote the joint law of
$(\bo\Omega,\mc P_\Omega(\mb R))$, and let $P_u$ denote the
one-observation exponential-family law with natural parameter $u$. For any
$u,v\in\mc D_{\sf E}$, Taylor's theorem gives, for some $t\in(0,1)$,
\begin{equation}
    \mathsf{KL}(P_u\,\|\,P_v)
    =\Psi(v)-\Psi(u)-\psi(u)(v-u)
    =\frac{1}{2}\psi'\big(u+t(v-u)\big)(u-v)^2
    \lesssim
    \frac{\underline c_{\psi,{\sf E}}^2}{\sigma_{\sf E}^2}(u-v)^2,
    \label{eq: exponential-family scalar KL bound}
\end{equation}
where the last inequality follows from
$\sup_{u\in\mc D_{\sf E}}\psi'(u)=\sigma_{\sf E}^2
\asymp\underline c_{\psi,{\sf E}}$.
Since the mask law is parameter-free, for
$\bo\theta_1,\bo\theta_2\in\mc S_{\sf E}(r,s,\omega)$,
\begin{equation}
    \mathsf{KL}(\bb P_{\bo\theta_1}\,\|\,\bb P_{\bo\theta_2})
    \lesssim\pi
    \frac{\underline c_{\psi,{\sf E}}^2}{\sigma_{\sf E}^2}
    \bignorm{\mc M(\bo\theta_1)-\mc M(\bo\theta_2)}\fb^2.
    \label{eq: exponential-family matrix KL bound}
\end{equation}

For $q\geq8$, the Varshamov--Gilbert bound provides
$\mc A_q\subseteq\{-1,1\}^q$ such that
\begin{equation}
    \log|\mc A_q|\geq c_0q,
    \qquad
    d_{\sf H}(\mb a,\mb b)\geq\frac q8
    \quad\text{for all distinct }\mb a,\mb b\in\mc A_q,
    \label{eq: Varshamov code for minimax lower bound}
\end{equation}
where $d_{\sf H}$ is Hamming distance and $c_0>0$ is numerical. Define
\begin{equation*}
    \mb u_0\coloneqq\frac{1}{2\sqrt q}
    \begin{bmatrix}
        \mb 1_q\\ \mb 1_q\\-\mb 1_q\\-\mb 1_q
    \end{bmatrix},
    \qquad
    \mb w^{(\mb a)}\coloneqq\frac{1}{2\sqrt q}
    \begin{bmatrix}
        \mb a\\-\mb a\\ \mb a\\-\mb a
    \end{bmatrix},
    \quad \mb a\in\mc A_q.
\end{equation*}
They are centered unit vectors with entries of order $q^{-1/2}$ and satisfy
$\mb u_0\t\mb w^{(\mb a)}=0$. For $0<\delta\leq1/4$, set
$\mb u^{(\mb a)}\coloneqq\sqrt{1-\delta^2}\,\mb u_0
+\delta\mb w^{(\mb a)}$ for $\mb a\in\mc A_q$, write
$\mb u^{(0)}\coloneqq\mb u_0$, and set
$\mc A_q^0\coloneqq\{0\}\cup\mc A_q$. These are centered unit vectors with
entries of order $q^{-1/2}$. For distinct
$\gamma,\gamma'\in\mc A_q^0$,
\begin{equation}
    c_{\sf u}\delta^2
    \leq\bignorm{\mb u^{(\gamma)}-\mb u^{(\gamma')}}_2^2
    \leq C_{\sf u}\delta^2,
    \qquad
    {\mb u^{(\gamma)}}\t\mb u^{(\gamma')}>0,
    \label{eq: bounded centered cap packing}
\end{equation}
for numerical $c_{\sf u},C_{\sf u}>0$.

\paragraph{Intercept lower bound.}
Choose orthonormal $\mb U\in\bb R^{n\times r}$ and
$\mb V\in\bb R^{p\times r}$ such that
$\mb 1_n\t\mb U=\mb 0_r\t$,
$\bignorm{\mb U}\ti\leq C_{r_0}/\sqrt n$, and
$\bignorm{\mb V}\ti\leq C_{r_0}/\sqrt p$, and set
$\mb X_0\coloneqq\omega^{1/4}\sqrt s\,\mb U$ and
$\mb Y_0\coloneqq\omega^{-1/4}\sqrt s\,\mb V$.
These factors are centered and balanced, with all $r$ nonzero singular values
of $\mb X_0\mb Y_0\t$ equal to $s$.

Let $\mc A_p=\{\mb a^1,\ldots,\mb a^N\}$ be a Varshamov--Gilbert set
satisfying \eqref{eq: Varshamov code for minimax lower bound}. Fix a
sufficiently small $c_\delta>0$, depending only on the constants in the
theorem, and let
\begin{equation*}
    \delta_\zeta\coloneqq c_\delta
    \min\Big\{c_1,
    \frac{\sigma_{\sf E}}
    {\underline c_{\psi,{\sf E}}\sqrt{\pi n}}\Big\},
    \qquad
    \bo\zeta^{(0)}\coloneqq c_\zeta\mb 1_p,
    \qquad
    \bo\zeta^{(k)}\coloneqq
    c_\zeta\mb 1_p+\delta_\zeta\mb a^k,
    \quad k\in[N].
\end{equation*}
Set $\bo\theta_0=[\bo\zeta^{(0)},\mb X_0,\mb Y_0]$ and
$\bo\theta_k=[\bo\zeta^{(k)},\mb X_0,\mb Y_0]$ for $k\in[N]$. The choice of
$c_2$ in the theorem ensures that all these triples belong to
$\mc S_{\sf E}(r,s,\omega)$. Moreover,
\begin{align*}
    \bignorm{\mc M(\bo\theta_k)
    -\mc M(\bo\theta_0)}\fb^2
    &=np\delta_\zeta^2,\\
    \min_{0\leq k<l\leq N}
    \bignorm{\bo\zeta^{(k)}-\bo\zeta^{(l)}}_2
    &\geq c\sqrt p\,\delta_\zeta.
\end{align*}
By \eqref{eq: exponential-family matrix KL bound},
\begin{equation*}
    \max_{k\in[N]}
    \mathsf{KL}(\bb P_{\bo\theta_k}\,\|\,\bb P_{\bo\theta_0})
    \lesssim\pi
    \frac{\underline c_{\psi,{\sf E}}^2}{\sigma_{\sf E}^2}
    np\delta_\zeta^2
    \lesssim c_\delta^2p
    \leq\alpha_{\sf F}\log|\mc A_p|,
\end{equation*}
where $\alpha_{\sf F}<1/8$ is numerical.
Fano's inequality \citep[Chapter~15]{wainwright2019high} yields
\begin{equation}
    \inf_{\widehat{\bo\zeta}}
    \sup_{\bo\theta\in\mc S_{\sf E}(r,s,\omega)}
    \bb E_{\bo\theta}
    \bignorm{\widehat{\bo\zeta}-\bo\zeta}_2
    \geq c\sqrt p\Big\{c_1\wedge
    \frac{\sigma_{\sf E}}
    {\underline c_{\psi,{\sf E}}\sqrt{\pi n}}\Big\}.
    \label{eq: Fano intercept lower bound}
\end{equation}

\paragraph{Factor lower bound.}
Suppose first that $n\geq p$.
Let $m_n\coloneqq\lfloor n/8\rfloor$ and set $\delta_X\coloneqq c_\delta
\{1\wedge\sigma_{\sf E}\sqrt n/
(\underline c_{\psi,{\sf E}}\sqrt\pi\,s)\}$.
Use $q=m_n$ and $\delta=\delta_X$ in
\eqref{eq: bounded centered cap packing}, and embed the resulting vectors in
the first $4m_n$ coordinates of $\bb R^n$. On the complement, choose
$\mb U_{\sf f}\in\bb R^{(n-4m_n)\times(r-1)}$ with centered orthonormal
columns and row norm at most $C_{r_0}/\sqrt n$. For
$\gamma\in\mc A_{m_n}^0$, set
\begin{equation*}
    \mb U^{(\gamma)}\coloneqq
    \begin{bmatrix}
        \mb u^{(\gamma)}&\mb 0\\
        \mb 0&\mb U_{\sf f}
    \end{bmatrix}.
\end{equation*}
The second block is void for $r=1$. Then $\mb U^{(\gamma)}$ is centered and
orthonormal with row norm at most $C_{r_0}/\sqrt n$; retain the fixed $\mb V$
from the intercept construction.

For $\gamma\in\mc A_{m_n}^0$, set
\begin{equation*}
    \mb X^{(\gamma)}\coloneqq
    \omega^{1/4}\sqrt s\,\mb U^{(\gamma)},
    \qquad
    \mb Y_{\sf f}\coloneqq
    \omega^{-1/4}\sqrt s\,\mb V,
    \qquad
    \bo\vartheta_\gamma\coloneqq
    [c_\zeta\mb 1_p,\mb X^{(\gamma)},\mb Y_{\sf f}].
\end{equation*}
Then $\bo\vartheta_\gamma\in\mc S_{\sf E}(r,s,\omega)$, and
\begin{equation}
    \bignorm{\mc M(\bo\vartheta_\gamma)-\mc M(\bo\vartheta_0)}\fb^2
    =s^2\bignorm{\mb u^{(\gamma)}-\mb u_0}_2^2
    \leq Cs^2\delta_X^2.
    \label{eq: Fano factor matrix separation}
\end{equation}
By \eqref{eq: exponential-family matrix KL bound} and
\eqref{eq: Fano factor matrix separation},
\begin{equation*}
    \max_{\mb a\in\mc A_{m_n}}
    \mathsf{KL}(\bb P_{\bo\vartheta_{\mb a}}\,\|\,\bb P_{\bo\vartheta_0})
    \lesssim\pi
    \frac{\underline c_{\psi,{\sf E}}^2}{\sigma_{\sf E}^2}
    s^2\delta_X^2
    \lesssim c_\delta^2n
    \leq\alpha_{\sf F}\log|\mc A_{m_n}|.
\end{equation*}

For distinct $\gamma,\gamma'\in\mc A_{m_n}^0$,
${\mb U^{(\gamma)}}\t\mb U^{(\gamma')}
=\operatorname{diag}\big({\mb u^{(\gamma)}}\t
\mb u^{(\gamma')},\mb I_{r-1}\big)$ has positive diagonal entries. Hence
\begin{align*}
    &d_\omega\big((\mb X^{(\gamma)},\mb Y_{\sf f}),
    (\mb X^{(\gamma')},\mb Y_{\sf f})\big)^2\geq2s
    \Big\{r-\bignorm{{\mb U^{(\gamma)}}\t\mb U^{(\gamma')}}_*\Big\}
    =s\bignorm{\mb u^{(\gamma)}-\mb u^{(\gamma')}}_2^2
    \geq cs\delta_X^2.
\end{align*}
Fano's inequality gives
\begin{equation}
    \inf_{(\widehat{\mb X},\widehat{\mb Y})}
    \sup_{\bo\theta\in\mc S_{\sf E}(r,s,\omega)}
    \bb E_{\bo\theta}
    d_\omega\big((\widehat{\mb X},\widehat{\mb Y}),(\mb X,\mb Y)\big)
    \geq c\sqrt s\Big\{
    1\wedge\frac{\sigma_{\sf E}}{\underline c_{\psi,{\sf E}}}
    \sqrt{\frac{n}{\pi s^2}}\Big\}.
    \label{eq: Fano X-side factor lower bound}
\end{equation}

If $p>n$, let $m_p\coloneqq\lfloor p/8\rfloor$ and
$\delta_Y\coloneqq c_\delta\min\{1,
\sigma_{\sf E}\sqrt p/
(\underline c_{\psi,{\sf E}}\sqrt\pi\,s)\}$. Apply the cap construction with
$(q,\delta)=(m_p,\delta_Y)$ in the first $4m_p$ coordinates of $\bb R^p$
and denote the resulting vectors by $\mb v^{(\gamma)}$. Append the same
$r-1$ bounded-row orthonormal columns, supported on the complement, to form
$\mb V^{(\gamma)}$ with row norm at most $C_{r_0}/\sqrt p$.
Fix a centered orthonormal
$\mb U_{\sf c}\in\bb R^{n\times r}$ with row norm at most $C_{r_0}/\sqrt n$,
and set
\begin{equation*}
    \mb X_{\sf f}\coloneqq\omega^{1/4}\sqrt s\,\mb U_{\sf c},
    \qquad
    \mb Y^{(\gamma)}\coloneqq
    \omega^{-1/4}\sqrt s\,\mb V^{(\gamma)},
    \qquad
    \widetilde{\bo\vartheta}_\gamma\coloneqq
    [c_\zeta\mb 1_p,\mb X_{\sf f},\mb Y^{(\gamma)}].
\end{equation*}
These triples lie in $\mc S_{\sf E}(r,s,\omega)$, and
\begin{align*}
    \bignorm{\mc M(\widetilde{\bo\vartheta}_\gamma)
    -\mc M(\widetilde{\bo\vartheta}_0)}\fb^2
    &\leq Cs^2\delta_Y^2,\\
    \max_{\mb a\in\mc A_{m_p}}
    \mathsf{KL}(\bb P_{\widetilde{\bo\vartheta}_{\mb a}}\,\|\,
    \bb P_{\widetilde{\bo\vartheta}_0})
    &\lesssim c_\delta^2p
    \leq\alpha_{\sf F}\log|\mc A_{m_p}|.
\end{align*}
For distinct $\gamma,\gamma'\in\mc A_{m_p}^0$, the same positive
cross-Gram calculation gives
\begin{align*}
    & d_\omega\big((\mb X_{\sf f},\mb Y^{(\gamma)}),
    (\mb X_{\sf f},\mb Y^{(\gamma')})\big)^2
     \geq\omega^{1/2}\min_{\mb O\in\mc O(r)}
    \bignorm{\mb Y^{(\gamma)}\mb O-\mb Y^{(\gamma')}}\fb^2=s\bignorm{\mb v^{(\gamma)}-\mb v^{(\gamma')}}_2^2\\
    &\quad \geq c s\delta_Y^2.
\end{align*}
Fano's inequality gives
\begin{equation*}
    \inf_{(\widehat{\mb X},\widehat{\mb Y})}
    \sup_{\bo\theta\in\mc S_{\sf E}(r,s,\omega)}
    \bb E_{\bo\theta}
    d_\omega\big((\widehat{\mb X},\widehat{\mb Y}),(\mb X,\mb Y)\big)
    \geq c\sqrt s\Big\{
    1\wedge\frac{\sigma_{\sf E}}{\underline c_{\psi,{\sf E}}}
    \sqrt{\frac{p}{\pi s^2}}\Big\}.
\end{equation*}
Combining the two cases with \eqref{eq: Fano intercept lower bound} proves
Theorem~\ref{thm: minimax lower bound}. \qed

\section{Auxiliary Lemmas}

\subsection{Properties of Random Sampling Operator. }
The following lemma is an extension of \cite[Corollary~4.3]{candes2012exact}. 
\begin{lemma}
    Consider an entrywisely transformed matrix $ \psi(\mb X\mb Y\t+ \mb 1_n\bo \zeta\t) \in \bb R^{n \times p}$ with a deterministic intercept $\mb A$ and each entry independently observed with probability $\pi$. Let $\mb U \in \bb R^{n \times r}$ and $\mb V \in \bb R^{p \times r}$ be the left and right top-$r$ singular vector matrices of $\mb X\mb Y\t $, respectively, and suppose that $\psi$ is $\bar c_\psi$-Lipschitz continuous and $\sqrt{\frac{n\norm{\mb U}\ti^2}{r}} \vee \sqrt{\frac{p\norm{\mb V}\ti^2}{r}} \leq \mu$. Then there exists a constant $C_R$ such that for all $\beta >1$, $\bo \Delta_{\mb X}\in \bb R^{n \times r}$, $\bo \Delta_{\mb Y} \in \bb R^{p \times r}$, and $\bo \Delta_{\bo \zeta} \in \bb R^p$, 
    \begin{align}
        & \pi^{-1}\norm{\mc P_{\Omega}\big(\psi(\mb 1_n(\bo \zeta + \bo \Delta_{\bo \zeta})\t + (\mb X + \bo \Delta_{\mb X})(\mb Y + \bo \Delta_{\mb Y})\t) ) - \psi(\mb 1_n \bo \zeta\t+ \mb X \mb Y\t ) \big)}\fb \\ 
        \lesssim & \sqrt{\frac{1}{\pi}} \bar c_\psi\bignorm{\bo \Delta_{\mb X}\mb Y\t + \mb X{\bo \Delta_{\mb Y}}\t}\fb + \sqrt{\frac{1}{\pi}} \bar c_\psi\sqrt{n}\norm{\bo \Delta_{\bo \zeta}}_2 + \frac{1}{\pi}\bar c_\psi \bignorm{\bo \Delta_{\mb X}\bo \Delta_{\mb Y}\t}\fb
    \end{align}
    uniformly holds with probability at least $1- 3 (n\vee p)^{-\beta}$ provided that
    \begin{equation*}
    C_R\sqrt{\frac{\mu_0 (n \vee p)r ( \beta \log(n\vee p))}{\pi np }} \leq 1.
    \end{equation*}
    \label{lemma: generalized version of lemma 4.4 in candes and recht}
\end{lemma}
\begin{proof}
    It is seen from the definition of the Frobenius norm and the Lipschitz-ness of $\psi$ that 
    \begin{align*}
        & \pi^{-1}\bignorm{\mc P_{\Omega}\big(\psi(\mb 1_n (\bo \zeta + \bo \Delta_{\bo \zeta})\t + (\mb X + \bo \Delta_{\mb X})(\mb Y + \bo \Delta_{\mb Y})\t + \mb A)) - \psi(\mb 1_n \bo \zeta\t + \mb X \mb Y\t) \big)}\fb \\ 
        & \leq  \pi^{-1} \bar c_\psi \bignorm{\mc P_{\Omega}\big( \mb X \bo \Delta_{\mb Y}\t + \bo\Delta_{\mb X} {\mb Y^*}\t \big)}\fb + \pi^{-1} \bar c_\psi \bignorm{\mc P_{\Omega}(\mb 1_n \bo \zeta\t )}\fb + \pi^{-1} \bar c_\psi \bignorm{\bo \Delta_{\mb X}\bo \Delta_{\mb Y}\t}\fb \\ 
    & \stackrel{\text{\cite[Corollary~4.3]{candes2012exact}}}{\leq}  \sqrt{\frac{3}{2\pi}} \bar c_\psi \bignorm{ \mb X \bo \Delta_{\mb Y}\t + \bo\Delta_{\mb X} {\mb Y^*}\t}\fb  +  \sqrt{\frac{3}{2\pi}} \bar c_\psi \sqrt{n} \bignorm{\bo \Delta_{\bo \zeta}}_2\\
    &\quad + \pi^{-1} \bar c_\psi \bignorm{\bo \Delta_{\mb X}\bo \Delta_{\mb Y}\t}\fb
    \end{align*}
    with probability at least $1- 3 (n\vee p)^{-\beta}$. 
\end{proof}

To handle the type of quantities $\ip{ \mb S \circ \mc P_{\Omega}(\mb A \mb C\t ), \mc P_{\Omega}(\mb B \mb D\t) }$, we have the following lemma directly from \cite[Lemma~4.4]{chen2019model}. 
\begin{lemma}
    Given arbitrary matrices
$\mb A\in\mathbb R^{n\times r_1}$,
$\mb B\in\mathbb R^{n\times r_2}$,
$\mb C\in\mathbb R^{p\times r_1}$,
$\mb D\in\mathbb R^{p\times r_2}$, and $\mb S \in \bb R^{n\times p}$, it holds that 
\begin{align}
    & \big| \bigip{\mb S \circ \mc P_{\Omega}(\mb A \mb C\t ), \mc P_{\Omega}(\mb B \mb D\t) } - \pi \bigip{\mb S \circ (\mb A \mb C\t) , \mb B \mb D\t } \big| \\ 
    & \leq   \norm{\mc P_\Omega(\mb S) - \pi \mb S} \max\big\{\norm{\mb A}\fb \norm{\mb B}\ti, \norm{\mb A}\ti \norm{\mb B}\fb \big\} \max\big\{\norm{\mb C}\fb \norm{\mb D}\ti,\\
    &\quad \norm{\mb C}\ti \norm{\mb D}\fb \big\} .
\end{align}
\label{lemma: lemma 4.4 chen2019}
\end{lemma}

The following lemma associates the projection error $\norm{\mc P_\Omega(\mb M) - \pi \mb M}$ with a friendly term $\norm{\bo \Omega - \pi \mb 1_{n_1 \times n_2}}$.  
\begin{lemma}[\cite{bhojanapalli2014universal}]
    Suppose a matrix $\mb M \in \bb R^{n_1 \times n_2}$ can be decomposed as $\mb M = \mb B \mb D\t$. Then we have
    \begin{align}
        & \norm{\mc P_\Omega(\mb M) - \pi \mb M} \leq \norm{\bo \Omega - \pi \mb 1_{n_1 \times n_2}} \norm{\mb B}\ti \norm{\mb D}\ti. 
    \end{align}
    \label{lemma: bhojanapalli2014universal}
\end{lemma}

We next control the sampling deviation uniformly over admissible interpolated
predictors. Fix numerical constants $C_{\mb X},C_{\mb Y}>0$. For
$s\in\{1,2\}$, let $\bo\theta^s=[\bo\zeta^s,\mb X^s,\mb Y^s]$, where
$\bo\zeta^s\in\bb R^p$, $\mb X^s\in\bb R^{n\times r}$, and
$\mb Y^s\in\bb R^{p\times r}$, and define
\begin{equation*}
\begin{aligned}
    \mc M_\psi\coloneqq\Big\{
    &\psi'\big((1-t)\mc M(\bo\theta^1)+t\mc M(\bo\theta^2)\big):
    t\in[0,1],\\
    &\norm{\mb X^s}\ti\leq C_{\mb X}\sqrt{\frac{\mu r\sigma_1^*}{n}},
    \quad\norm{\mb Y^s}\ti\leq C_{\mb Y}\sqrt{\frac{\mu r\sigma_1^*}{p}},\\
    &\mc M(\bo\theta^s)_{i,j}\in\mc D_\psi
    \text{ for all }i\in[n],\ j\in[p],\ s\in\{1,2\}\Big\}.
\end{aligned}
\end{equation*}
Since $\mc D_\psi$ is an interval, the endpoint restriction also places
every interpolated predictor in $\mc D_\psi$.
\begin{lemma}
    Suppose that Assumption~\ref{assumption: link function} holds and that
    entries are sampled independently with probability
    $\pi\gtrsim\log d/(n\wedge p)$. Then, with probability at least
    $1-O(d^{-c-2})$, uniformly over $\mb M\in\mc M_\psi$,
    \begin{equation*}
        \norm{\mc P_\Omega(\mb M)-\pi\mb M}
        \lesssim\sigma^2\xi^4\sqrt{\pi d},
    \end{equation*}
    where $\xi$ is defined in
    Assumption~\ref{assumption: signal strength and incoherence degree},
    with $c_\xi$ sufficiently large depending on $C_{\mb X}$ and
    $C_{\mb Y}$.
    \label{lemma: injectivity of P_Omega under nonlinearity}
\end{lemma}
\begin{proof}[Proof of Lemma~\ref{lemma: injectivity of P_Omega under nonlinearity}]
Fix an admissible pair of endpoints and $t\in[0,1]$. Write
$a=C_{\mb X}\sqrt{\mu r\sigma_1^*/n}$,
$b=C_{\mb Y}\sqrt{\mu r\sigma_1^*/p}$, and
$\ell=\lceil\max\{1,8rL_\psi ab\}\rceil$. Set
\begin{equation*}
    \mb u_i=\begin{pmatrix}\mb X_i^1\\\mb X_i^2\end{pmatrix},
    \qquad
    \mb v_j=\begin{pmatrix}(1-t)\mb Y_j^1\\t\mb Y_j^2\end{pmatrix},
    \qquad z_j=(1-t)\zeta_j^1+t\zeta_j^2.
\end{equation*}
Thus $M_{i,j}=\psi'(z_j+\mb u_i\t\mb v_j)$. Partition
$[-a,a]^{2r}$ into $\ell^{2r}$ cells of side length $2a/\ell$.
For each occupied cell $E_\nu$, choose an index $i_\nu$ with
$\mb u_{i_\nu}\in E_\nu$ and use $\mb u_{i_\nu}$ as its Taylor center.
For $\mb u_i\in E_\nu$, put
$\bo\delta_i=\mb u_i-\mb u_{i_\nu}$ and
$c_{\nu,j}=z_j+\mb u_{i_\nu}\t\mb v_j$.
Both $c_{\nu,j}$ and $c_{\nu,j}+\bo\delta_i\t\mb v_j$ are entries
of the admissible interpolated predictor. Their connecting segment lies
in $\mc D_\psi$, so Assumption~\ref{assumption: link function} applies
throughout the Taylor expansion of $\psi'$.

In particular, $|\psi^{(k+1)}(x)|/k!\leq\sigma^2(k+1)L_\psi^k$
for $x\in\mc D_\psi$ and $k\geq0$. Since
$\norm{\bo\delta_i}_\infty\leq2a/\ell$ and
$\norm{\mb v_j}_\infty\leq b$, we have
$L_\psi|\bo\delta_i\t\mb v_j|\leq4rL_\psi ab/\ell\leq1/2$.
The remainder after degree $K$ is therefore at most
$\sigma^2(K+2)2^{-(K+1)}$, uniformly over all entries. It follows that
\begin{equation*}
    M_{i,j}=\sum_{k=0}^{\infty}
    \frac{\psi^{(k+1)}(c_{\nu,j})}{k!}
    (\bo\delta_i\t\mb v_j)^k,
    \qquad \mb u_i\in E_\nu.
\end{equation*}

For a multi-index $\bo\alpha\in\bb N_0^{2r}$, write
$\bo\alpha!=\prod_{h=1}^{2r}\alpha_h!$ and
$\mb w^{\bo\alpha}=\prod_{h=1}^{2r}w_h^{\alpha_h}$.
On each cell, the degree-$k$ multinomial expansion is a sum of rank-one
matrices with entries
$\ind\{\mb u_i\in E_\nu\}\bo\delta_i^{\bo\alpha}
\psi^{(k+1)}(c_{\nu,j})\mb v_j^{\bo\alpha}/\bo\alpha!$,
where $|\bo\alpha|=k$. The sum of the products of the maximum absolute
row and column coefficients is bounded by
\begin{equation*}
    \sum_{|\bo\alpha|=k}
    \frac{\sigma^2(k+1)!L_\psi^k}{\bo\alpha!}
    \Big(\frac{2ab}{\ell}\Big)^k
    =\sigma^2(k+1)\Big(\frac{4rL_\psi ab}{\ell}\Big)^k
    \leq\sigma^2(k+1)2^{-k}.
\end{equation*}
Applying Lemma~\ref{lemma: bhojanapalli2014universal} to these rank-one
terms and summing over the occupied cells and all degrees gives
\begin{align*}
    \norm{\mc P_\Omega(\mb M)-\pi\mb M}
    &\leq\sigma^2\ell^{2r}
    \norm{\bo\Omega-\pi\mb 1_n\mb 1_p\t}
    \sum_{k=0}^{\infty}(k+1)2^{-k}\\*
    &=4\sigma^2\ell^{2r}
    \norm{\bo\Omega-\pi\mb 1_n\mb 1_p\t}.
\end{align*}
This bound is deterministic and uniform over $\mc M_\psi$. Moreover,
$\ell\leq C(L_\psi\mu r^2\sigma_1^*/\sqrt{np}\vee1)$ for a fixed
$C$ depending only on $C_{\mb X}$ and $C_{\mb Y}$. Choosing $c_\xi$
sufficiently large before taking the power $2r$ gives $\ell^{2r}\leq\xi^4$.
Finally, the Bernoulli sampling condition and
\cite[Corollary~3.12 and Remark~3.13]{bandeira2016sharp}, applied to the
Hermitian dilation of $\bo\Omega-\pi\mb 1_n\mb 1_p\t$, give
$\norm{\bo\Omega-\pi\mb 1_n\mb 1_p\t}\lesssim\sqrt{\pi d}$ with
probability at least $1-O(d^{-c-2})$. This single event proves the stated
uniform bound.
\end{proof}

\subsection{Concentration Inequalities}
\paragraph{Noise matrix concentrations.}
\begin{lemma}[Spectral norms of the response and sampling noise]
    \label{lemma: noise spectral norm concentration}
    Suppose that Assumption~\ref{assumption: noise and sampling} holds, and let
    \[
        \mb E=\mb R-\mb R^*,
        \qquad
        \bar{\mb E}=\mc P_\Omega(\mb R)-\pi\mb R^*
    \]
    be the noise matrices defined in \eqref{eq: response and spectral noise matrices}. Then, with probability at least $1-O(d^{-c-2})$,
    \begin{align}
        \norm{\mc P_\Omega(\mb E)}
        &\lesssim \sigma\sqrt{\pi d}+B\sqrt{\log d}, \label{eq: response noise spectral concentration}\\
        \norm{\bar{\mb E}}
        &\lesssim \bar\sigma\sqrt{\pi d}+B\sqrt{\log d}. \label{eq: spectral noise concentration}
    \end{align}
    Under Assumption~\ref{assumption: signal strength and incoherence degree}(c), these bounds simplify to
    \[
        \norm{\mc P_\Omega(\mb E)}\lesssim \sigma\sqrt{\pi d},
        \qquad
        \norm{\bar{\mb E}}\lesssim \bar\sigma\sqrt{\pi d}.
    \]
\end{lemma}

\begin{proof}
The entries of both matrices are independent and centered. Moreover,
\[
    \bb E\!\left[\{\mc P_\Omega(\mb E)_{i,j}\}^2\right]
    \leq \pi\sigma^2,
    \qquad
    \bb E\!\left[(\bar{\mb E}_{i,j})^2\right]
    =\pi\mathrm{Var}(R_{i,j})+\pi(1-\pi)(R_{i,j}^*)^2
    \leq 2\pi\bar\sigma^2.
\]
On the high-probability event in Assumption~\ref{assumption: noise and sampling}, the entries of both matrices are bounded in magnitude by $2B$. Applying the rectangular version of \citet[Corollary~3.12 and Remark~3.13]{bandeira2016sharp} gives \eqref{eq: response noise spectral concentration} and \eqref{eq: spectral noise concentration}. The simplified bounds follow from Assumption~\ref{assumption: signal strength and incoherence degree}(c) and $d\geq n\wedge p$.
\end{proof}

\subsection{Matrix Perturbation Theory}
 
\begin{lemma}[Polar-factor perturbation; {\citealp[Lemma~23]{ma2018implicit}}]
\label{lemma: perturbation theory of optimal rotation under F norm}
Let $\mb C\in\bb R^{r\times r}$ be nonsingular, and define
\[
    \widehat{\mb H}(\mb C)
    \coloneqq
    \mb C(\mb C^\top\mb C)^{-1/2}.
\]
For any $\mb E\in\bb R^{r\times r}$ satisfying
$\norm{\mb E}\leq\sigma_r(\mb C)$ and any unitarily invariant norm
$\normiii{\cdot}$,
\[
    \normiii{\widehat{\mb H}(\mb C+\mb E)
    -\widehat{\mb H}(\mb C)}
    \leq
    \frac{2\normiii{\mb E}}
    {\sigma_{r-1}(\mb C)+\sigma_r(\mb C)}.
\]
\end{lemma}

The following lemma provides a joint perturbation bound for the two factors of a
rectangular matrix. Its proof uses the positive eigenspace of the Hermitian
dilation, which couples the left and right singular vectors and therefore
produces a common rotation for both factors.
\begin{lemma}[Joint perturbation of rectangular factors]
\label{lemma: perturbation theory for joint factors}
\label{lemma: perturbation theory for unilateral factors}
Let $\mb U^*\bo\Sigma^*\mb V^{*\top}$ be a top-$r$ SVD of
$\mb M^*\in\bb R^{n_1\times n_2}$, and let $(\mb X^*,\mb Y^*)$ be the
associated balanced factors. Thus, for some $\mb Q^*\in\mc O(r)$,
\begin{equation*}
    \mb X^*=\mb U^*(\bo\Sigma^*)^{1/2}\mb Q^*,
    \qquad
    \mb Y^*=\mb V^*(\bo\Sigma^*)^{1/2}\mb Q^*.
\end{equation*}
Let $\mb U\bo\Sigma\mb V\t$ be a top-$r$ SVD of
$\mb M\in\bb R^{n_1\times n_2}$, and define
\begin{align*}
    \mb X\coloneqq \mb U\bo\Sigma^{1/2},\quad \mb Y \coloneqq \mb V\bo\Sigma^{1/2},\quad 
    \mb W\coloneqq \frac{1}{\sqrt 2}
    \begin{pmatrix}\mb U\\ \mb V\end{pmatrix},\quad \mb W^* \coloneqq \frac{1}{\sqrt 2}
    \begin{pmatrix}\mb U^*\\ \mb V^*\end{pmatrix}.
\end{align*}
Let
\begin{equation}
    \bar{\mb O}
    \in\arg\min_{\mb O\in\mc O(r)}
    \bignorm{\mb W\mb O-\mb W^*}\fb,
    \qquad
    \mb O\coloneqq\bar{\mb O}\mb Q^*.
    \label{eq: common rotations from Hermitian dilation}
\end{equation}
Set $\sigma_{r+1}^*\coloneqq\sigma_{r+1}(\mb M^*)$. Suppose, for fixed numerical
constants $c_{\sf{gap}}>0$ and sufficiently small $c_0>0$, that
\begin{equation*}
    \norm{\mb M-\mb M^*}\leq c_0\sigma_r^*,
    \qquad
    \sigma_r^*-\sigma_{r+1}^*\geq c_{\sf{gap}}\sigma_r^*,
    \qquad
    \frac{\sigma_1^*}{\sigma_r^*}\leq\kappa.
\end{equation*}
Then, for $\diamond\in\{\mathrm{op},\mathrm F\}$, there exists a numerical
constant $C>0$ such that
\begin{equation}
    \bignorm{\mb X\mb O-\mb X^*}_{\diamond}
    \vee
    \bignorm{\mb Y\mb O-\mb Y^*}_{\diamond}
    \leq
    \frac{C\kappa}{\sqrt{\sigma_r^*}}
    \bignorm{\mb M-\mb M^*}_{\diamond},
    \label{eq: joint rectangular factor perturbation}
\end{equation}
where $\norm{\cdot}_{\mathrm{op}}$ denotes the spectral norm. 
\end{lemma}

\begin{proof}
The proof follows the perturbation argument of \citet[Lemmas~32 and
33]{ma2018implicit} after reducing the rectangular problem to a symmetric one.
For any rectangular matrix $\mb A$, define its Hermitian dilation by
\begin{equation*}
    \mc H(\mb A)
    \coloneqq
    \begin{pmatrix}
        \mb 0 & \mb A\\
        \mb A\t & \mb 0
    \end{pmatrix}.
\end{equation*}
Write $\mb H\coloneqq\mc H(\mb M)$ and
$\mb H^*\coloneqq\mc H(\mb M^*)$. The columns of $\mb W$ and $\mb W^*$
span the invariant subspaces of $\mb H$ and $\mb H^*$, respectively,
associated with their $r$ largest algebraic eigenvalues, and
\begin{equation*}
    \mb H\mb W=\mb W\bo\Sigma,
    \qquad
    \mb H^*\mb W^*=\mb W^*\bo\Sigma^*.
\end{equation*}
Moreover,
\begin{equation}
    \norm{\mb H-\mb H^*}_{\diamond}
    \leq \sqrt 2\,\norm{\mb M-\mb M^*}_{\diamond}.
    \label{eq: norm identity for Hermitian dilation}
\end{equation}
The spectrum of $\mb H^*$ consists of the singular values of $\mb M^*$,
their negative counterparts, and possibly zeros. By Weyl's inequality, if
$c_0$ is sufficiently small relative to $c_{\sf{gap}}$, the interval
$
    [\frac{\sigma_r^*+\sigma_{r+1}^*}{2},\infty)
$
contains exactly the $r$ largest algebraic eigenvalues of both $\mb H$ and
$\mb H^*$. The Davis--Kahan theorem gives 
\begin{equation}
    \bignorm{\mb W\bar{\mb O}-\mb W^*}_{\diamond}
    \lesssim
    \frac{\bignorm{\mb M-\mb M^*}_{\diamond}}{\sigma_r^*}.
    \label{eq: aligned dilation eigenspace perturbation}
\end{equation}

For convenience, set $
    \widetilde{\mb W}\coloneqq\mb W\bar{\mb O}$, $
    \widetilde{\bo\Sigma}
    \coloneqq \bar{\mb O}\t\bo\Sigma\bar{\mb O}$.
Since
$\widetilde{\bo\Sigma}
= \widetilde{\mb W}\t\mb H\widetilde{\mb W}$ and
$\bo\Sigma^*=\mb W^{*\top}\mb H^*\mb W^*$,
\eqref{eq: norm identity for Hermitian dilation} and
\eqref{eq: aligned dilation eigenspace perturbation} imply
\begin{align}
    \bignorm{\widetilde{\bo\Sigma}-\bo\Sigma^*}_{\diamond}
    &\lesssim
    \bignorm{\mb M-\mb M^*}_{\diamond}
    +\sigma_1^*
    \bignorm{\widetilde{\mb W}-\mb W^*}_{\diamond} 
    \lesssim
    \kappa\bignorm{\mb M-\mb M^*}_{\diamond}.
    \label{eq: aligned eigenvalue perturbation}
\end{align}
In addition, Weyl's inequality gives
$\lambda_{\min}(\widetilde{\bo\Sigma})\gtrsim\sigma_r^*$.
The matrix-square-root perturbation step in the proof of
\citet[Lemma~33]{ma2018implicit} therefore yields, for both
$\diamond=\mathrm{op}$ and $\diamond=\mathrm F$,
\begin{equation}
    \bignorm{\widetilde{\bo\Sigma}^{1/2}
        -(\bo\Sigma^*)^{1/2}}_{\diamond}
    \lesssim
    \frac{\kappa}{\sqrt{\sigma_r^*}}
    \bignorm{\mb M-\mb M^*}_{\diamond}.
    \label{eq: aligned square root perturbation}
\end{equation}

Finally, define the stacked factors
\begin{equation*}
    \mb F\coloneqq
    \begin{pmatrix}\mb X\\ \mb Y\end{pmatrix}
    =\sqrt 2\,\mb W\bo\Sigma^{1/2},
    \qquad
    \mb F^*\coloneqq
    \begin{pmatrix}\mb X^*\\ \mb Y^*\end{pmatrix}
    =\sqrt 2\,\mb W^*(\bo\Sigma^*)^{1/2}\mb Q^*.
\end{equation*}
Orthogonal equivariance of the matrix square root gives
\begin{equation*}
    \mb F\mb O
    =\sqrt 2\,\widetilde{\mb W}
    \widetilde{\bo\Sigma}^{1/2}\mb Q^*.
\end{equation*}
Combining this identity with
\eqref{eq: aligned dilation eigenspace perturbation} and
\eqref{eq: aligned square root perturbation} proves
\eqref{eq: joint rectangular factor perturbation}. 
\end{proof}

\subsection{Lemmas for Distributional Approximations}
\label{subsec:lemmas-distributional-approximations}

We collect several distributional comparison results that are repeatedly used to
transfer the linear representations derived above into Gaussian and
multiplier-bootstrap approximations.  Throughout this subsection, let $\mc C_d$
denote the collection of all convex Borel subsets of $\bb R^d$, and let
$\mc R_p$ denote the collection of all hyperrectangles in $\bb R^p$.  For
$C\subseteq\bb R^d$ and $\varepsilon>0$, define
$$
C^\varepsilon
=
\{x\in\bb R^d:\inf_{y\in C}\norm{x-y}_2\leq \varepsilon\},
\qquad
C^{-\varepsilon}
=
\{x\in C:\inf_{y\in C^c}\norm{x-y}_2>\varepsilon\}.
$$
For probability measures $P$ and $Q$ on the same measurable space, write
$$
\mathrm{TV}(P,Q)
=
\sup_A |P(A)-Q(A)|.
$$

We begin with a finite-dimensional Gaussian approximation over convex sets.  In
our applications, this result is useful after the target statistic has been
reduced to a fixed-dimensional linear form of independent score contributions.

\begin{lemma}[Convex-set Gaussian approximation {\citep[Theorem~1.1]{raic2019multivariate}}]
\label{lemma:convex-set-gaussian-approximation}
Let $\bo\xi_1,\ldots,\bo\xi_n$ be independent mean-zero random vectors in
$\bb R^d$, and set $\mb W=\sum_{i=1}^n\bo\xi_i$.  Suppose
$\bo\Sigma=\cov(\mb W)$ is positive definite, and let
$\mb Z\sim \mc N(\mb 0,\bo\Sigma)$.  Then
$$
\sup_{C\in\mc C_d}
\left|
\bb P(\mb W\in C)-\bb P(\mb Z\in C)
\right|
\leq
(42d^{1/4}+16)
\sum_{i=1}^n
\bb E\left[
\bignorm{\bo\Sigma^{-1/2}\bo\xi_i}_2^3
\right].
$$
\end{lemma}

For simultaneous inference, the relevant objects are often maxima of many
coordinates rather than fixed-dimensional convex-set probabilities.  The next
result gives a high-dimensional approximation for such maxima using Gaussian
multiplier critical values.

\begin{lemma}[Gaussian multiplier approximation for maxima {\citep[Theorems~2.1--2.2]{chernozhukov2022improved}}]
\label{lemma:gaussian-multiplier-max-approximation}
Let $\mb X_1,\ldots,\mb X_n$ be independent mean-zero random vectors in
$\bb R^p$, where $p\geq 2$ and $n\geq 3$.  Let $B_n\geq 1$ and let
$0<b_1\leq b_2<\infty$ be constants such that, for every
$i\in[n]$ and $j\in[p]$,
$$
\bb E\exp\{|X_{ij}|/B_n\}\leq 2,
$$
and, for every $j\in[p]$,
$$
b_1^2
\leq
\frac{1}{n}\sum_{i=1}^n \bb E X_{ij}^2,
\qquad
\frac{1}{n}\sum_{i=1}^n \bb E X_{ij}^4
\leq
B_n^2 b_2^2 .
$$
For a deterministic vector $\mb a=(a_1,\ldots,a_p)^\top$, define
$$
T_n
=
\max_{1\leq j\leq p}
\left\{
\frac{1}{\sqrt n}\sum_{i=1}^n X_{ij}
+
a_j
\right\}.
$$
Let $\bar{\mb X}=n^{-1}\sum_{i=1}^n\mb X_i$, and let
$e_1,\ldots,e_n\stackrel{\mathrm{i.i.d.}}{\sim}\mc N(0,1)$ be independent of
the data.  Define the Gaussian multiplier statistic
$$
T_n^{\sf{mb}}
=
\max_{1\leq j\leq p}
\left\{
\frac{1}{\sqrt n}\sum_{i=1}^n e_i(X_{ij}-\bar X_j)
+
a_j
\right\}.
$$
Let $c^{G}_{1- \alpha}$ be the $(1 - \alpha)$-quantile of $T_n$ with $X_{i,j}$ substituted by its Gaussian counterpart with the same mean and variance, and $c_{1-\alpha}^{\sf{mb}}$ be the conditional $(1-\alpha)$-quantile of
$T_n^{\sf{mb}}$ given $\mb X_1,\ldots,\mb X_n$.  Then there exists a constant
$C>0$, depending only on $b_1$ and $b_2$, such that
$$
\left|
\bb P\left(T_n>c_{1-\alpha}^{G}\right)-\alpha
\right| \vee \left|
\bb P\left(T_n>c_{1-\alpha}^{\sf{mb}}\right)-\alpha
\right|
\leq
C
\left\{
\frac{B_n^2\log^5(pn)}{n}
\right\}^{1/4},
\qquad 0<\alpha<1 .
$$
\end{lemma}
\begin{remark}
\label{remark:high-dim-clt}
    This result also implies that
    \begin{equation}
        \sup_{A \in \mc A}\Big|\bb P\big(\frac{1}{\sqrt{n}} \sum_{i\in[n]} \mb X_i \in A \big) - \bb P\left(\mb G \in A \right)\Big| \leq C
\left\{
\frac{B_n^2\log^5(pn)}{n}
\right\}^{1/4},
    \end{equation}
    where $\mb G \sim \mc N\big(\mb 0, \mathrm{Cov}(\frac{1}{\sqrt{n}} \sum_{i\in[n]} \mb X_i) \big)$ and $\mc A$ is the class of all hyper-rectangles in $\bb R^p$, namely, sets of the form
    $$
    A = \big\{ (w = (w_1, \cdots, w_p)\t \in \bb R^p: a_{lj} \leq w_j \leq a_{rj}\text{ for all }j\in[p] \big\}, 
    $$
    for some constants $-\infty \leq a_{lj} \leq a_{rj} \leq \infty$ with $j \in[p]$. 
\end{remark}
\begin{remark}
The empirical centering in the multiplier statistic can be omitted
without affecting its asymptotic validity. Indeed,
\begin{align*}
    &\Big|\max_{j\in[p]}
    \Big\{n^{-1/2}\sum_{i\in[n]}e_iX_{ij}+a_j\Big\}
    -T_n^{\sf mb}\Big|\leq
    \norm{\bar{\mb X}}_\infty
    \Big|n^{-1/2}\sum_{i\in[n]}e_i\Big|.
\end{align*}
{\emergencystretch=1.5em
The exponential-moment assumption and Bernstein's inequality give
\begin{equation*}
\norm{\bar{\mb X}}_\infty
=O_{\bb P}\{B_n\sqrt{\log(pn)/n}\}
\end{equation*}
whenever $B_n^2\log^5(pn)/n=o(1)$.
Since $n^{-1/2}\sum_i e_i\sim\mc N(0,1)$ conditionally on the data,
the displayed difference is
$o_{\bb P}\{(\log p)^{-1/2}\}$.
Together with a uniform empirical-variance lower bound and
Gaussian anti-concentration, this yields the same asymptotic
coverage for the uncentered multiplier statistic.\par}

\end{remark}

Finally, when a statistic is approximated by a leading Gaussian term plus a
small remainder, one needs to control how much Gaussian probability can
accumulate near the boundary of the relevant rejection region.  The following
anti-concentration bound provides this control uniformly over convex sets.

\begin{lemma}[Gaussian boundary mass for convex sets {\citep[Theorem~1.2]{raic2019multivariate}}]
\label{lemma:convex-gaussian-boundary}
Let $\mb G\sim\mc N(\mb 0,\mb I_d)$.  For each $C\in\mc C_d$, define
$$
\Gamma(C)
=
\sup_{\varepsilon>0}
\max\left\{
\frac{\bb P(\mb G\in C^\varepsilon\setminus C)}{\varepsilon},
\frac{\bb P(\mb G\in C\setminus C^{-\varepsilon})}{\varepsilon}
\right\},
\qquad
\Gamma_d
=
\sup_{C\in\mc C_d}\Gamma(C).
$$
Then
$$
\Gamma_d
\leq
0.59d^{1/4}+0.21 .
$$
\end{lemma} 

\section{Implementation Details for Numerical Experiments}
\label{sec: experimental implementation details}

We use spectral threshold
$\lambda^{\sp}=\tau\sqrt{\max(n,p)/\pi}$ with $\tau=1$, and clipping
to $[0.05,0.95]$ before applying the inverse link. No clipping is used in
the unilateral refinement step.

We run $200$ iterations of vanilla gradient descent without gradient
clipping, parameter clipping, or projection. The factor learning rate is
$\eta=m_\eta/(\pi\sigma_1)$, where $m_\eta=0.5$ in the simulations and
$m_\eta=1$ in the MetaBench analysis. Here,
$\sigma_1=\sigma_1(\mb X^*\mb Y^{*\top})$ in the estimation-error
simulations and
$\sigma_1=\sigma_1(\hat{\mb X}^{\ur}\hat{\mb Y}^{\ur\top})$ in the
coverage and real-data analyses. The blockwise step sizes are
$\eta_{\mb X}=\eta_{\mb Y}=\eta$ and
$\eta_{\bo\zeta}=\eta\,\sigma_r(\hat{\mb X}^{\ur}\hat{\mb Y}^{\ur\top})/n$.

In the MetaBench analysis, this gives $\eta=1.3688\times10^{-4}$, and
we use $1000$ bootstrap draws.

\end{appendices}

\end{document}